\documentclass[12pt]{amsart}
\usepackage[all]{xy}
\usepackage{amscd}
\usepackage{amssymb}
\usepackage{enumerate}
\usepackage{amsthm}
\usepackage[pagebackref]{hyperref}
\theoremstyle{plain}
\newtheorem{thm}{Theorem}[section]
\newtheorem{prop}[thm]{Proposition}

\newtheorem{cor}[thm]{Corollary}

\newtheorem{lem}[thm]{Lemma}
\newtheorem{lemdefn}[thm]{Lemma + definition}

\newtheorem{example}[thm]{Example}

\newtheorem{rem}[thm]{Remark}
\newtheorem{defn}[thm]{Definition}
\theoremstyle{definition}

\newtheorem{para}[subsubsection]{}

\input cyracc.def   
\define\Hom{\mathrm{Hom}\,}%

\def\j{\langle j\rangle}
\def\jj{\langle 1\rangle}

\def\MT{\cal{MT}}
\def\M{\cal M}
\def\F{\cal F}
\def\C{\cal C}

\def\B{\cal B}

\def\I{\cal I}

\def\Ltens{{\mathop\otimes\limits^{L}}}

\newcommand{\Gm}{{\mathbb{G}_m}}

\newcommand{\Ga}{{\mathbb G_a}}
\newcommand{\N}{{\mathbb{N}}}

\newcommand{\Z}{{\mathbb{Z}}}
\newcommand{\Fp}{{\mathbb F_p}}

\newcommand{\Zp}{{\mathbb Z_p}}

\newcommand{\upc}[1]{\overset {\lower 0.3ex \hbox{${\;}_{\circ}$}}{#1}}

\newcommand{\cal}{\mathcal}

\renewcommand{\theenumi}{\roman{enumi}}
\newcommand{\limi}{\displaystyle{\lim_{\longrightarrow}}}
\newcommand{\limp}{\displaystyle{\lim_{\longleftarrow}}}

\def\diagram#1{\def\normalbaselines{\baselineskip=0pt\lineskip=10pt\lineskiplimit=1pt}
  \begin{matrix}#1\end{matrix}}
\def\hfl#1#2{\smash{\mathop{\hbox to
10mm{\rightarrowfill}}\limits^{\scriptstyle#1}_{\scriptstyle#2}}}
\def\hflrev#1#2{\smash{\mathop{\hbox to
10mm{\leftarrowfill}}\limits^{\scriptstyle#1}_{\scriptstyle#2}}}
\def\hflcourte#1#2{\smash{\mathop{\hbox to
3mm{\rightarrowfill}}\limits^{\scriptstyle#1}_{\scriptstyle#2}}}
\def\hflrevcourte#1#2{\smash{\mathop{\hbox to
3mm{\leftarrowfill}}\limits^{\scriptstyle#1}_{\scriptstyle#2}}}

\def\vfl#1#2{\llap{$\scriptstyle #1$}\left\downarrow\vbox to
6mm{}\right.\rlap{$\scriptstyle #2$}}

\def\vflupcourte#1#2{\llap{$\scriptstyle #1$}\left\uparrow\vbox to
2mm{}\right.\rlap{$\scriptstyle #2$}}

\def\diagram#1{\def\normalbaselines{\baselineskip=0pt\lineskip=10pt\lineskiplimit=1pt}
\begin{matrix}#1\end{matrix}}

\def\debrom{
\makeatletter
\renewcommand{\theenumi}{(\roman{enumi})}
\renewcommand{\labelenumi}{\theenumi}
\makeatother\begin{enumerate}}
\def\finrom{\end{enumerate}}

\def\G{\Gamma}
\def\Zp{{\mathbb Z_p}}
\def\Qp{{\mathbb Q_p}}

\def\petit{\vspace{.2cm}}

\def\grand{\vspace{.6cm}}

\def\s{\sharp}
\def\O{\cal O}

\def\DC{\cal DC}

\def\S{{\cal S}}
\def\I{{\cal I}}

\begin{document}

\title{A comparison theorem for semi-abelian schemes over a smooth curve.}

\begin{abstract}
We compare flat cohomology to crystalline syntomic complexes with coefficients in two cases: 1) $p$-divisible groups over a separated $\Fp$-scheme with local finite $p$-bases, 2) semi-abelian schemes over a separated irreducible smooth curve.
\end{abstract}

\date{\today}

\author[Trihan]{Fabien Trihan}
\address{Department of Information and Communication Sciences\\
Faculty of Science and Technology, Sophia University\\
4 Yonbancho, Chiyoda-ku, Tokyo 102-0081 JAPAN.}
\email{f-trihan-52m@sophia.ac.jp}

\author[Vauclair]{David Vauclair}
\address{Laboratoire de Mathématiques Nicolas Oresme\\
Universit\'e de Caen\\Campus 2\\
14032 Caen Cedex\\
FRANCE.}
\email{david.vauclair@unicaen.fr}

\maketitle

\tableofcontents

\section{Introduction} \label{seci}


%
%
%
%
%
%

Consider an odd prime number $p$. To a scheme $X$ over $\Fp$ is naturally associated two kinds of integral $p$-adic cohomology groups. On the one hand, we have the $p$-adic \'etale cohomology groups $H^q_{et}(X,\Zp)$ and on the other hand we have the crystalline cohomology groups $H^q_{crys}(X/\Zp,\O_{X/\Zp})$. If $X$ is smooth (or more generally if $X$ is syntomic over a scheme with local finite $p$-bases) then it is well known (using e.g. \cite{FM}) that both are related by a long exact sequence \begin{eqnarray}\label{introselcst}\xymatrix{H^q_{et}(X,\Zp)\ar[r]&H^q_{crys}(X/\Zp,\O_{X/\Zp})\ar[r]^-{1-Frob}&
H^q_{crys}(X/\Zp,\O_{X/\Zp})\ar[r]^-{+1}&}
\end{eqnarray}
It is  natural to expect generalizations of this result when $\Zp$ is replaced by more general coefficients. A first (easy) step is to replace $\Zp$ by an arbitrary lisse $\Zp$-sheaf. A possible way of investigation to treat more general coefficients is to replace the \'etale topology by finer ones such as the flat or syntomic topology (note that such a change of topology does not affect cohomology for lisse $\Zp$-sheaves). Using the techniques of Fontaine-Messing again, one could presumably treat positive Tate twists $\Zp(i)$ (at least for $i\le p-1$). The Tate module of Abelian schemes have been treated in \cite{Ba}.  The purpose of this paper is to establish a suitable comparison theorem for $p$-divisible groups in the first place and  then for semi-Abelian schemes following the ideas of \cite{KT}.
\petit

Assume that $X$ is separated and has local finite $p$-bases. Consider the covariant Dieudonn\'e functor of \cite{BBM}, $D$, from the category $pdiv(X)$ of $p$-divisible groups over $X$ to the category $\cal DC(X)$ of Dieudonn\'e crystals over $X$. It is known that $D$ is fully faithful (\cite{BM}) and is in fact an equivalence (\cite{dJ1}). However, the proofs given in \cite{BM} and \cite{dJ1} do not produce explicitly a quasi-inverse functor. Instead they rely on deformation techniques in order to reduce to the case where $X$ is the spectrum of a field. In this paper we construct the \emph{syntomic complex} functor  \begin{eqnarray}\label{introSsyn}\cal S_{syn,.,X}:\cal DC(X)\to D^b(X^\N_{syn},\Z/p^.)\end{eqnarray} whose target is the bounded derived category of $\Z/p^.$-modules on the small syntomic site of $X$  (here $\Z/p^.$ denotes the pro-ring $(\Z/p^k)_{k\ge 1}$). By construction, we have a functorial distinguished triangle \begin{eqnarray}\label{introtdsyn}\xymatrix{\cal S_{syn,.,X}(D)\ar[r]&\Z/p^.\otimes_{\Z/p^{.+1}}u_{X/(\Z/p^{.+1}),*}Fil^1D_{.+1}\ar[r]^-{1-\varphi}
&u_{X/(\Z/p^.),*}D_.\ar[r]^-{+1}&}\end{eqnarray}
where $D_.$ is the restriction of $D$ to the ind crystalline topos $(X/(\Z/p^.))_{crys,syn}$, $Fil^1$ denotes the natural $mod\, p$ Hodge filtration on $D$, $u_{X/(\Z/p^.)}$ is the projection to $X_{syn}^\N$, $1$ denotes the obvious morphism and $\varphi$ is the unique $\O^{crys}_.$ semi-linear morphism such that $p\varphi$ is induced by the Frobenius of $D$. In this setting our first main result is the following.

\begin{thm} \label{thmA} (Thm. \ref{thmcomppdiv}) Let $G$ be a $p$-divisible group over $X$ and let  $G_{p^.}$ denote the projective system of $p$ power torsion subgroups of $G$ viewed as sheaves on the small syntomic site of $X$. There is a canonical isomorphism $$G_{p^.}\simeq \cal S_{syn,.,X}(D(G))$$
This isomorphism is functorial with respect to $G$ and $X$.
\end{thm}
In particular, applying $Rlimproj$ and taking cohomology in (\ref{introtdsyn}) yields a long exact sequence
\begin{eqnarray}\label{introselG}
\xymatrix{H^q_{syn}(X,T_p(G))\ar[r]&H^q_{crys}(X/\Zp,Fil^1D(G))\ar[r]^-{1-\varphi}&
H^q_{crys}(X/\Zp,D(G))\ar[r]^-{+1}&}\end{eqnarray}
which boils down to (\ref{introselcst}) for $G=\Qp/\Zp$.

Consider now a separated irreducible smooth curve $C$ over $\Fp$. Our second main result is a comparison theorem for semi-Abelian schemes over $C$. Using rigid uniformization around bad fibers, we construct a functor $D$, from the category of semi-Abelian schemes over $C$ with good reduction on a dense open $U$, to the category $\cal DC(C^\s)$ of Dieudonn\'e crystals over the logarithmic curve $(C,Z)$, where $Z$ is the reduced divisor complementary to $U$. By construction, the restriction of this functor to Abelian schemes coincides with the one of \cite{BBM}. Next we construct the \emph{twisted syntomic complex} functor \begin{eqnarray}\label{introSet}\cal S_{et,.,C^\s}(-Z):\cal DC(C^\s)\to D^b(C^\N_{et},\Z/p^.)\end{eqnarray}
whose restriction to $\cal DC(C)$ can be recovered from the functors (\ref{introSsyn}) on $X=C$ and $X=Z$ by projection from the small syntomic topos to the small \'etale one (more precisely a version of (\ref{introSsyn}) is needed above the diagram $X=(C\leftarrow Z)$, see below for explanations). Let $z:Z\to C$ denote the inclusion morphism and $\underline \G^Z:\smash{Mod(C^\N_{FL},\Z/p^.)}\to \smash{Mod(C^\N_{FL},\Z/p^.)}$ denote the functor taking $M$ to the kernel of the specialization morphism $M\to z_*z^{-1}M$. We prove the following.

\begin{thm} \label{thmB} (Thm. \ref{compss}) Consider a semi-Abelian scheme $A$ over $C$ whose restriction to $U$ is Abelian. Let $\epsilon:C_{FL}\to C_{et}$ denote the projection of the big flat topos to the small \'etale one. There is a canonical isomorphism $$R\epsilon_*R\underline \G^Z A_{p^.}\simeq \cal S_{et,.,C^\s}(-Z)(D(A))$$
This isomorphism is functorial with respect to $A$ and $Z$.
\end{thm}
Applying $Rlimproj$ and passing to cohomology yields in particular a long exact sequence
$$\xymatrix{H^{Z,q}_{FL}(C,T_p(A))\ar[r]& H^q_{crys}(C^\s/\Zp,Fil^1D(A)(-Z))\ar[r]^-{1-\varphi}&H^q_{crys}(C^\s/\Zp,D(A)(-Z))\ar[r]^-{+1}&}$$
where $H^{Z,q}_{FL}$ means \emph{flat cohomology vanishing at $Z$}. If $A$ is in fact an Abelian scheme and $Z=\emptyset$ this result follows from Thm. \ref{thmA}.

\para The results of this paper certainly won't come as a surprise to experts. In fact Thm. \ref{thmA} and Thm. \ref{thmB} are nothing more than sheafified versions of \cite{KT} Prop. 5.10 and  Prop. 5.13. These refinements will be the main ingredient for a proof of a non commutative Iwasawa main conjecture in \cite{TV}. Our original motivation was thus to check whether or not the proofs of \cite{KT} could be sheafified as well. While doing so, we discovered several problems in \emph{loc. cit.} and it finally became easier to rewrite everything than trying to fix the mistakes and/or missing arguments. It might be worth however to explain the length of this paper with regards to theirs. We encountered mainly two sources of difficulties.

The first one concerns the proof of the comparison theorem for $p$-divisible groups. The strategy of \cite{KT} Prop. 5.10 is to use the equivalence of categories between $p$-divisible groups and Dieudonn\'e crystals and to interpret cohomology as higher extension groups. Unfortunately, some confusions regarding continuous cohomology create a serious gap in the argumentation and we do not know how to fix it (some detailed explanations are given in \cite{Va to KT}). We circumvent this difficulty by following the method of \cite{Ba} instead. Namely, we show that the vanishing results for $\cal Ext$'s of \cite{Br} together with the techniques of \cite{FM} are sufficient to treat not only Abelian schemes but $p$-divisible groups as well.

The second source of difficulty concerns the proof of the comparison theorem for semi-Abelian schemes by reduction to the case of  $p$-divisible groups. Roughly speaking, the strategy of \cite{KT} Prop. 5.13 is to perform a parallel d\'evissage in flat and crystalline cohomology. Unfortunately, this d\'evissage is merely sketched in \emph{loc. cit.} and some delicate issues are left to the reader. On the side of flat cohomology, the main ingredient is to replace $A$ by a diagram of $p$-divisible groups using Raynaud's rigid uniformization around bad fibers. However a precise sheaf theoretic interpretation is not given and the cohomological consequences are left to the reader to find. Here these tasks will be achieved at the end of Sect. \ref{dipdg}.
On the side of crystalline cohomology, their first step is to define $D(A)$ by gluing \cite{BBM}'s $D(A_{|U})$ with an \emph{ad hoc} (logarithmic) Dieudonn\'e crystal on the complete neighborhood of bad points. We will check that this definition is functorial and is the only one possible by studying the structure of the category of log $1$-motives introduced in \cite{KT}. Their second step is to relate the twisted syntomic complex of $D(A)$ to the one attached to the diagram of $p$-divisible groups introduced before. In \emph{loc. cit.} explanations are only given in the local situation and globalization is left to the reader. Unfortunately the local computations explained in \emph{loc. cit.} Lem. 5.14 in presence of a lifting are not enough to perform the cohomological descent hinted in \emph{loc. cit. 5.8}. We solve this issue by using hypercoverings with divisors in the spirit of \cite{NS} (see Lem. \ref{lememb3}).

An additional new difficulty is that we have to compare the syntomic complex on the \'etale site (\ref{introSet}) (which is useful for d\'evissage) with the one on the syntomic site (\ref{introSsyn}) (appearing in Thm. \ref{thmA}). We achieve this comparison by introducing an intermediate variant (using linearized de Rham complexes). A significant part of this work is thus devoted to the construction of several variants of the syntomic complex functor and to comparing them.


\para The organization of the paper is the following.
\petit

\ref{secppi}. \emph{Preliminaries part I: category theoretical background}.

 Constructions like gluing functors, projective limits etc. are often needed at the level of derived categories. A convenient way to get rid of technical complications is to work systematically above an arbitrary diagram of (log) schemes rather than above a single (log) scheme. 
 The relevant facts are recalled in Sect. \ref{dacfc} and Sect. \ref{ftadc}. Specifically, in Sect. \ref{dacfc}, we recall some basic language and facts regarding pseudo-functors (or equivalently (co)fibered categories) and (co)lax morphisms between them. We also include a descent lemma for colax morphisms (Prop.  \ref{CBdescent2}), which will be used in Sect. \ref{paraprelimS} to prove that the syntomic complex functors do not depend on the choices of hypercoverings and Frobenius lifts. In Sect. \ref{ftadc},  we fix some terminology regarding fibered topoi (prevariable pretopologies, weakly variable topoi, see.  Def. \ref{deftopfib} \ref{deftopfibi}, \ref{deftopfibii}). We also discuss their derived categories and extension to diagrams. In addition to their elementary properties, whose use is pervasive in the text, we establish some acyclicity conditions (Lem.-Def. \ref{acycf}
 \ref{acycfiv}, \ref{acycfv}) which will be useful later.

Next, our conventions regarding the usual and crystalline (pre)topologies of a log-scheme are explained in Sect. \ref{usfls} and Sect. \ref{csfls}. In Sect. \ref{scsf}, we establish the pseudo-functoriality of the small crystalline topos. Some elementary properties of crystals are recalled in Sect. \ref{cocs} since they will be needed in Sect. \ref{laqcocs}. 
\petit

\ref{secppii}. \emph{Preliminaries part II: some properties for sheaves}.

The purpose of this chapter is to gather several technical results which are related to exactness. In Sect. \ref{pdg}
we recall some terminology (semi-abelian schemes, $1$-motives, $p$-divisible groups...) and review the well known fact that the category of 1-motives over a regular base is exact, since it will be essential for our construction of the semi-stable Dieudonné functor. Several  possible sheaf theoretic incarnations (on usual or crystalline sites) of a constant or $p$-divisible group are reviewed in Sect. \ref{topsheaves}.

In Sect. \ref{psapdozpa}, we review the basic lemma designed to divide Frobenius by $p$ (Lem. \ref{cartgen}). This involves the notion of ($L$-)normalized modules or complexes (Def. \ref{defnorm}) and the normalizing functor $\jj^*$ (Def. \ref{defj}). 

In Sect. \ref{laqcopafs} and Sect.  \ref{laqcocs}, we design a convenient framework regarding $p$-adic formal  schemes, limits and quasi-coherence. A standard feature of crystalline cohomology is that one is often led to pass from finite to $p$-adic coefficients and conversely. Inside the crystalline topos this causes no difficulty (see Lem.-Def. \ref{lemlimcrys}). When passing to realizations on the other hand (e.g. for de Rham complexes) one has to add some quasi-coherence assumptions in order to get (partial) analogous results. The relevant notions are adapted from \cite{Be5} Sect. 3.2 and discussed in Sect. \ref{paraqcohpadic} - \ref{critqcoh}.

\petit

\ref{secppiii}. \emph{Preliminaries part III: crystals and local finite $p$-bases}.

In Sect. \ref{pbale}, we review a notion of finite $p$-bases for ($p$-adic formal) log schemes which is stronger than the one defined in \cite{Ts2},  but slightly easier to handle. We also introduce the corresponding categories of (global and local) embeddings, which will be enriched later.

In Sect. \ref{mwlc} and Sect. \ref{cadrcoclf}, we assume given a closed immersion of the base log scheme into a $p$-adic formal log scheme with local finite $p$-bases and we restrict our attention to the \'etale topology. In this context, we explain how to adapt the classical theory of integrable quasi-nilpotent connections, hyper $dp$-stratifications, hyper $dp$-differential operators, linearization functors and de Rham complexes. While doing so, we pay special attention to the functoriality of the theory with respect to the base log scheme and the chosen closed immersion, since this issue is of crucial importance for the construction of the syntomic complex functors.

When dealing with crystals, one has to be careful with the notion of subobjects because the inclusion of the category of crystals inside the category of modules on the crystalline site is not left exact in general. The relevant phenomenons are discussed in Sect. \ref{epotcoc} using the previously introduced categories of embeddings (rather than a restricted crystalline site in the spirit of \cite{Be1}). In Sect. \ref{tbeld},  we discuss the relation between effective logarithmic and Cartier divisors. Putting everything together, we study to what extent a twist of a module or a crystal can be regarded as a subobject of the latter (see Lem. \ref{injtwist1} and Lem. \ref{injtwist2}, which will be used to divide Frobenius by $p$).

\petit

\ref{sectscfdc}. \emph{Twisted syntomic complexes for Dieudonn\'e crystals}.

In this chapter we define functors \begin{eqnarray}\label{introSonephisyn}\cal S^{1,\varphi}_{syn,.,X}:\DC(X)\to Mod^{1,\varphi}(X^\N_{syn},\O_.^{crys})\\
\label{introSonephiet}\cal S^{1,\varphi}_{et,.,X^\s}(-h):\DC(X^\s)\to D^+(Mod^{1,\varphi}(X^\N_{et},\O_.^{crys}))
\end{eqnarray}
The first (resp. second) one is defined when $X$ is in $\cal B_0$, the category of  locally embeddable diagrams of separated schemes whose vertices have local finite $p$-bases over $\Fp$ (resp. when $(X^\s,h)$ is in $\smash{\cal B_0^\s}$, the category of locally embeddable diagrams of separated log-schemes whose vertices have  local $p$-bases over $\Fp$, together with an effective log divisor), see Def. \ref{defB0} \ref{defB0i}. The target of these functors is the (derived) category of $(1,\varphi)$-modules, defined and briefly studied in Sect. \ref{tco1pm}. The functors (\ref{introSsyn}) and (\ref{introSet}) are respectively obtained from (\ref{introSonephisyn}) and (\ref{introSonephiet}) by forming the mapping fiber of $1-\varphi$ (see Lem. \ref{1phimod} \ref{1phimodiv}).
Aside from the definitions themselves, the main purpose of this chapter is to establish a canonical isomorphism \begin{eqnarray}\label{introRepsilon}R\epsilon_*\cal S^{1,\varphi}_{syn,.,X}D\simeq \cal S^{1,\varphi}_{et,.,X}D\end{eqnarray}
in the case $X^\s=X$ (trivial log structure). Here $\epsilon:X_{syn}\to X_{et}$ denotes the canonical morphism. This task requires several intermediate variants of the syntomic complex on the small \'etale and syntomic sites.
\petit

Let us briefly enumerate these constructions and explain the strategy leading to (\ref{introRepsilon}). In Sect. \ref{lewfacd}, we begin with some observations regarding applications of cohomological descent in crystalline cohomology and we define some categories of semi-simplicial embeddings with additional structures (Frobenius lifts, effective log divisors) which are adapted to the various constructions that we have in mind. Let us explain this roughly. A semi-simplicial global embedding $\iota:\smash{U^\s{[.]}\to Y_{[.]}^\s}$ together with a Frobenius lift is in the category $\smash{HR_F^{\s,et}}$ if $\smash{U^\s_{[.]}}$ is a hypercovering in the topos $X_{et}$. The full subcategory $\smash{HR_F^{\s,crys} \subset HR_F^{\s,et}}$ is defined by the additional condition that the logarithmic divided power envelope $\smash{T^\s_{[.]}}$  of $\iota$ is a hypercovering in the topos $\smash{(X^\s/\Zp)_{crys,et}}$. The latter is adapted to the computation of crystalline cohomology \emph{\`a la Cech} whereas the former is adapted to the computation using de Rham complexes. This leads naturally to the construction of two functors with value in a category of complexes of $(1,\varphi)$-modules: \begin{eqnarray}\label{introS1phiet}\cal S^{1,\varphi}_{et,.,T^\s_{[.]}}(-h),\cal S\Omega^{\bullet,1,\varphi}_{et,.,T^\s_{[.]}}(-h):\DC(X^\s)\to Kom^{1,\varphi}(T_{[.],.,et},\O^{crys}_.)\end{eqnarray}
(here and in the sequel the notation is abusive: these functors are not directly the ones obtained from (\ref{SetT}) applying $(-)^{loc}$, ie. by passing to diagrams (see. Sect. \ref{paraprelimS} (\ref{cons2})), but they are deduced from the latter by composition with (\ref{cons1}), ie. by restriction of Dieudonné crystals and divisors to $\smash{U^\s_{[.]}}$).
The first one (Cech) is not necessary in this text. It is nevertheless the most direct way of defining the syntomic complex and will serve as a warming up for the other constructions. The second one (de Rham) is well adapted to d\'evissage results and will be used in the next chapter.

In the case of a trivial log structure, we give a global construction of the functor (\ref{introSonephisyn}) (ie. not involving semi-simplicial embeddings). This construction has the advantage of being close to the context of \cite{FM} and \cite{Ba} and will be used in the comparison theorem for $p$-divisible groups. We also define linearized de Rham versions of the syntomic complex $\smash{\cal SL\Omega^{\bullet,1,\varphi}_{et,.,T_{[.]}}}$ and  $\smash{\cal SL\Omega^{\bullet,1,\varphi}_{syn,.,T_{[.]}}}$ on the small \'etale and syntomic site respectively. \petit

The main steps leading to (\ref{introRepsilon}) may be summarized as follows. Sect. \ref{tscotes} is concerned with the versions on the \'etale site. Using the descent lemma Prop. \ref{newglu} we prove (Lem. \ref{lemcompS} \ref{lemcompSi}, \ref{lemcompSii}) that the projection of $${\cal S\Omega^{\bullet,1,\varphi}_{et,.,T^\s_{[.]}}(-h)} \hbox{ (resp. of }\smash{\cal S^{1,\varphi}_{et,.,T^\s_{[.]}}(-h)}\hbox{, resp. of }\smash{\cal SL\Omega^{\bullet,1,\varphi}_{et,.,T^\s_{[.]}}(-h)})$$ \vspace{.1cm} to $D^+(Mod^{1,\varphi}(X_{et},\smash{\O^{crys}_.}))$ is essentially independent of the semi-simplicial embedding, as long as the latter is chosen in $\smash{HR_F^{\s,et}}$ (resp. $\smash{HR_F^{\s,crys}}$, resp. $\smash{HR_F^{\s,crys}}$). We prove furthermore (Lem \ref{lemcompS}) that the resulting three functors $\DC(X^\s)\to D^+(Mod^{1,\varphi}(X_{et},\O^{crys}_.))$ are canonically isomorphic and we define (\ref{introSonephiet}) as anyone of them (Prop. \ref{indepSet} and Def. \ref{defSX}). Sect. \ref{scotss} is concerned with the versions on the syntomic site. We prove
(Prop. \ref{compsynetloc} and Prop. \ref{compsynet} \ref{compsyneti}) that the globally defined functor (\ref{introSonephisyn}) coincides with the projection of $\smash{\cal SL\Omega^{\bullet,1,\varphi}_{syn,.,T_{[.]}}}$ to $D^+(Mod^{1,\varphi}(X_{syn},\O^{crys}_.))$. The isomorphism (\ref{introRepsilon}) is finally obtained by proving (Prop. \ref{compsynetloc}): \begin{eqnarray}R\epsilon_*\smash{\cal SL\Omega^{\bullet,1,\varphi}_{syn,.,T_{[.]}}}\simeq \smash{\cal S\Omega^{\bullet,1,\varphi}_{et,.,T_{[.]}}}\end{eqnarray}
\petit

Let us now give some explanations regarding the constructions enumerated above. 
In Sect.  \ref{trfacdfc}, we review some basics concerning the relative Frobenius and the Cartier operator. Then, we establish the crystalline version of the Cartier equivalence over $X$ (Prop. \ref{propcartier}) and deduce some elementary consequences for crystals over $X^\s$. These results will be the main ingredient of Sect. \ref{tmphfotsces}, where a crystalline subsheaf $Fil^1D$ is defined for an arbitrary Dieudonn\'e crystal $D$ and shown to satisfy a canonical isomorphism $D/Fil^1D\simeq i_*Lie(D)$ (Prop. \ref{liesurj}). Here, the tangent sheaf $Lie(D)$ is defined in terms of the Verschiebung operator modulo $p$. This filtration and a similar isomorphism are then extended to the linearized semi-simplicial crystal $\smash{L(D_{T_{[.],.}^\s})}$ and its twisted versions (Prop. \ref{complin}). In Sect. \ref{tscotes}, we check that in each three cases the Frobenius vanishes modulo $p$ on the subcomplex defined by using $Fil^1D$. Using this and the normalization functor $\jj^*$ (see Def. \ref{defj}), we are able to divide Frobenius by $p$. This yields the desired definition for the (Cech, de Rham and linearized de Rham) versions of the syntomic complex on the small \'etale site (Prop. \ref{defphi}).
In Sect. \ref{scotss}, the construction and elementary properties of the (global and linearized de Rham) versions of the syntomic complex on the syntomic site of $X$ are given using the same ingredients than in Sect. \ref{tscotes} together with some well known properties of crystalline cohomology with respect to the syntomic topology. \petit


\ref{secdofc}. \emph{D\'evissage of flat cohomology}.

Let $A/C$ be a semi-Abelian scheme with good reduction over $U=C-Z$. If $v$ is a point of $Z$ with residue field $k_v$, we denote $Z_v=Spec(k_v)$,  $C_v=Spec(O_v)$  the complete neighborhood of $v$ in $C$ and $U_v=Spec(K_v)$ the generic point of $C_v$. Consider the diagram $J^+$ with the following subschemes of $C$ as vertices: $U$ and for $v$ running in $Z$, the $Z_v$'s, the $C_v$'s and the $U_v$'s, and the inclusion morphisms between them as edges. The purpose of this chapter is to define a $p$-divisible group $H$ over $J^+$ (ie. a $p$-divisible group over each one of the subschemes in question together with base change morphisms between them) from which one can rebuild the projection of the vanishing local sections $R\underline \G^Z A_{p^.}\in D^b(C^\N_{FL},\Z/p^.)$  to the small \'etale site (see Prop. \ref{devissageFL} or (\ref{introfl4}) below for explanations). \petit

Let us explain the logical steps leading to this result. The first step is to replace the big flat site by the small \'etale one. In Sect. \ref{vc}, we begin with an elementary study of the \emph{vanishing (local) sections functor} $R\underline \G^Z$. The latter is an endofunctor of the derived category of an arbitrary topos of $C$. This study relies notably on Sect. \ref{fmfidc} (where the functorial properties of mapping fibers are discussed) and on the notions of acyclicity in fibered topoi which are discussed in Lem.-Def. \ref{acycf}. A significant drawback of the functor $R\underline \G^Z$ is  its lack of functoriality with respect to the chosen topos. For instance, if $\epsilon$ denotes the projection morphism from the big (say flat) to the small (say \'etale) topos, then $R\underline \G^Z$ does not commute to $R\epsilon_*$. This is mainly due to the following fact:
if $G$ is a group  scheme over $C$ then the small \'etale sheaf over $C$ represented by $G$ has no memory of the group scheme $Z\times_CG$ over $Z$. To circumvent this issue, we enrich the picture by replacing $C$ with the diagram $C^+:=(Z\to C)$. We introduce a functor $R\underline \G^Z(C,-)$ which is not an endofunctor, but goes from the (derived category of the chosen topos of) $C^+$ to (that of) $C$. As it turns out, this new variant of the vanishing sections functor commutes to $R\epsilon_*$. In the setting of big topoi we have moreover $R\underline \G^Z(C,(-)_{|C^+})\simeq R\underline \G^Z(-)$. Applying this to $A_{p^.}$ gives (see Cor. \ref{corsoritesvan}) \begin{eqnarray}\label{introfl1}R\epsilon_*R\underline \G^Z A_{p^.}\simeq R\underline \G^Z(C,R\epsilon_*A_{p^.,|C^+})\end{eqnarray}

In Sect. \ref{cmvfnm} we achieve the next step which is to replace $C$ (resp. $C^+$) by the diagram $J$ formed by $U$, the $U_v$'s and the $C_v$'s for $v$ in $Z$ (resp. $J^+$, obtained from $J$ by adding the $Z_v$'s). Consider the natural morphism $m:J\to C$. The functor $Rm_*$ fits into a familiar distinguished triangle of Mayer-Vietoris type (see Lem. \ref{defMV}) and is called the \emph{(complete) Mayer-Vietoris functor}. We prove that it allows to recover $R\epsilon_*A_{p^.}$ from the restriction of $A_{p^.}$ to $J$ and similarly for vanishing sections. More precisely (see Prop. \ref{MVfl}): \begin{eqnarray}\label{introfl2}R\underline \G^Z(C,R\epsilon_*A_{p^.,|C^+})\simeq Rm_*R\underline \G^{Z_J}(J,R\epsilon_*A_{p^.,|J^+})\end{eqnarray}
In addition to abstract nonsense, the proof notably relies on the fact that the direct image functor of $C_{v,et}\to C_{et}$ is exact (Lem.  \ref{iotaexact}) together with the following property of Abelian varieties which was found by Greenberg and Milne: if  $K$ denotes the function field of $C$, $v$ a point of $C$ and $K_v$ (resp. $K_v^h$) the completion (resp. henselization) of $K$ at $v$ then $H^q(K_v,A_{|K_v})$ coincides with $H^q(K_v^h,A_{|K_v^h})$ for $q\ge 1$. \petit

Let us now explain the construction of the $p$-divisible group $H$ on $J^+$ which will serve as a substitute for $A_{p^.,|J^+}$. First we have to define a $p$-divisible group $H_U$ (resp. $\smash{H_{U_v}}$, $\smash{H_{C_v}}$, $\smash{H_{Z_v}}$ for each $v$ in $Z$). When restricted to its open of good reduction $A$ gives rise directly to a $p$-divisible group. Thus we simply set $H_U:=\smash{A_{|U,p^\infty}}$ and $\smash{H_{|U_v}}:=\smash{A_{|U_v,p^\infty}}$.
Since $A$ has semi-stable reduction at $v$, we get a $p$-divisible group on $Z_v$  simply by replacing $\smash{A_{|Z_v}}$ with its connected component $\smash{A_{|Z_v}^0}$ (which is the extension of an Abelian variety $\smash{B_{k_v}}$ by a torus $\smash{T_{k_v}}$): we set $\smash{H_{Z_v}}=\smash{A_{|Z_v,p^\infty}^0}$. We know from \cite{SGA7-I} IX that $A_{|Z_v}^0$ admits a canonical lifting to a group scheme $G_v$ over $C_v$ (the \emph{Raynaud group}) which is the extension of an Abelian scheme lifting $\smash{B_{k_v}}$ by a torus lifting $\smash{T_{k_v}}$ and whose completion along the special fiber is isomorphic to the completion of $\smash{A_{C_v}^0}$. We get a $p$-divisible group on $C_v$ by setting $\smash{H_{C_v}}:=\smash{G_{v,p^\infty}}$. We make the following key observation: since $\smash{G_{v,p^k}}$ is finite over $C_v$ the previous formal isomorphism induces a morphism of group schemes $e_v:\smash{G_{v,p^k}}\to \smash{A_{|C_v,p^k}^0}$. This remark allows us to define the base change morphisms giving rise to the $p$-divisible group $H$ over $J^+$ as well as a morphism $\smash{H_{p^.}}\to \smash{A_{|J^+,p^.}^0}$ over $J^+$. This morphism is certainly not an isomorphism in general since $\smash{A_{|C_v,p^k}^0}$ might not be finite over $C_v$. In Sect. \ref{ruapoz}, we show that it is nevertheless possible to recover $\smash{A_{p^.,|J^+}^0}$ from $H_{p^.}$ using Raynaud's result concerning the generic fiber (in the sense of rigid geometry) of the morphism $e_v$. In Sect. \ref{dipdg}, we use this to prove (see Prop. \ref{riget} and the proof of Prop. \ref{devissageFL}) that \begin{eqnarray}\label{introfl3} R\underline \G^{Z_J}(J^+,R\epsilon_*A_{p^.,|J^+}) \simeq Sma\, R\underline \G^{Z_J}(J^+,R\epsilon_*H_{p^.})\end{eqnarray}
where $Sma$ is an exact functor defined over $J_{et}$, designed to neglect the generic fiber of the $C_v$ components (see Def. \ref{defsma} for more details). A key ingredient in the proof is the introduction of the small quasi-finite flat site (see Def. \ref{defqff})  which is fine enough to compute cohomology and small enough to express sheaf theoretic consequences of Raynaud's rigid uniformization. The final d\'evissage result (Prop. \ref{devissageFL}) is obtained by putting (\ref{introfl1}), (\ref{introfl2}) and (\ref{introfl3}) together:
\begin{eqnarray}\label{introfl4} R\epsilon_*R\underline \G^Z A_{p^.}\simeq Rm_*Sma\, R\underline \G^{Z_J}(J,R\epsilon_*H_{p^.})\end{eqnarray}

\petit

\ref{secdotsc}. \emph{D\'evissage of twisted syntomic complexes}.

 Our general strategy to obtain the comparison result in the case of a semi-Abelian scheme $A/C$ is to perform a d\'evissage for the twisted syntomic complex of $D_{C^\s}(A)$ which is roughly parallel to the one of Chap. \ref{secdofc} in order to reduce to the case of $p$-divisible groups. This won't be achieved before Sect. \ref{sasoc} for the following reasons. On the one hand, the semi-stable Dieudonn\'e functor $D_{C^\s}$ has yet to be defined (this will be done in the next chapter) and on the other hand, the desired reduction involves switching from $D_{C^\s}(A)$ to $D_{J}(H)$ where $H$ is the diagram of $p$-divisible groups defined in Chap. \ref{secdofc}. This will be done using a trick relying on the comparison result for $p$-divisible groups established in Sect. \ref{dopdg}. The purpose of the present chapter is nevertheless to establish two key d\'evissage results.

In Sect. \ref{cmv}, we consider the diagram of log schemes $J^\s$ which is $J$ together with the log structure coming from $C^\s$. Letting $m:J^\s\to C^\s$ denote the natural morphism, we prove (Prop. \ref{MVsyn} \ref{MVsyni}) that $m^*$ induces an equivalence between Dieudonn\'e crystals on $C^\s$ and diagrams of Dieudonn\'e crystals over $J^\s$ whose base change morphisms are invertible (we call those cartesian). By Zariski localization on $C^\s$, we reduce to the case where a lifting of $C^\s$ by a $p$-adic formal log scheme with $p$-bases over $Spf(\Zp)$ is given and then the result essentially boils down to gluing locally free modules from $J$ to $C$. Looking at de Rham complexes yields furthermore a canonical isomorphism  (Prop. \ref{MVsyn} \ref{MVsyniii}) \begin{eqnarray}\label{introisoMV}\cal S^{1,\varphi}_{et,.,C^\s}(-h)(D)\simeq Rm_*\cal S^{1,\varphi}_{et,.,J^\s}(-m^{-1}h)(m^*D)\end{eqnarray}

In Sect. \ref{alt}, we use the second variant of the local vanishing sections functor to describe the twisted syntomic complex in a specific situation. Namely consider a Dieudonn\'e crystal $D$ over $J$ and denote respectively $o^*D$ and $\rho^*D$ its pullback to $J^\s$ and $J^+$. Consider furthermore the smooth divisor $Z_J$ as an effective log divisor of $J$. Then we establish a canonical isomorphism \begin{eqnarray}\label{introisospe}\cal S^{1,\varphi}_{et,.,J^\s}(-Z_J)(o^*D)\simeq R\underline \G^{Z_J}(J,\cal S^{1,\varphi}_{et,.,J^+}(\rho^*D))\end{eqnarray}
Here the key idea is to work with semi-simplicial embeddings $\iota: U^\s_{[.]}\to Y^\s_{[.]}$ of $HR^{\s, et}_F$ having the additional property that the involved closed immersions of log schemes are exact. An easily  tractable category of such, denoted $HR_F^{\s, et,ex}$, was defined and studied in Sect. \ref{lewfacd} (see Lem. \ref{lememb3}) using blowing up as a global exactification functor. \petit

\petit

\ref{secdcfssav}. \emph{Dieudonn\'e crystals for semi-Abelian schemes}.

The goal of this chapter is to construct the Dieudonn\'e crystal $D_{C^\s}(A)\in \DC(C^\s)$ associated to a semi-Abelian scheme $A/C$. Thanks to the equivalence result of Prop. \ref{MVsyn} \ref{MVsyni}, we are reduced to define for each $v$ a $D_{C_v^\s}(A_{|C_v})\in \DC(C_v^\s)$ whose restriction to $U_v$ agrees with \cite{BBM}'s one. By Raynaud's uniformization as reviewed in Chap. \ref{secdofc}, we already know how to define a canonical \emph{log $1$-motive}  (Def. \ref{defMlog}) $M_{log}(A_{|C_v^\s})=(G_v,\G_v,f_v)$ associated to $A_{|C_v}$ (Def. \ref{defSASloc}) which is such that the $p$-divisible groups of $M_{log}(A_{|C_v^\s})_{|U_v}$ and $A_{|U_v}$ are canonically isomorphic. Our main purpose is thus to prove that there is a unique way of defining a Dieudonn\'e functor on the category $\cal M_{log}(C_v)$ of log $1$-motives over $C_v$ which is exact, compatible to finite \'etale base change and compatible with \cite{BBM}'s Dieudonn\'e functors on $p$-divisible groups over $C_v$ and $U_v$.  Here, we use the exact structure inherited from the category of Abelian sheaves on the big flat site (see Lem. \ref{1motex} and Rem. \ref{remlog1motex}).

We do this in two steps. The first one is achieved in Sect. \ref{dol1m}. By studying the structure of the category $\cal M_{log}(C'_v)$ for $C'_v$ varying among the finite \'etale $C_v$-schemes, we show (Prop. \ref{prol}) that given a stack $\cal C$ of exact categories  on the small finite \'etale site of $C_v$, extending a given exact functor $F:\cal M\to \cal C$ to $\cal M_{log}$ amounts to extending the homomorphism \begin{eqnarray}\label{introExt}Ext^1_{\cal M(C_v)}(\Z,\Z(1))\to Ext^1_{\cal C(C_v)}(F(\Z),F(\Z(1)))\end{eqnarray}
induced by $F$ to $Ext^1_{\cal M_{log}(C_v)}(\Z,\Z(1))$.
The second step will be achieved in Sect. \ref{tssdf}. It consists in proving that if $\cal C$ denotes the category of Dieudonn\'e crystals over (finite \'etale extensions of) $C_v^\s$, such an extension of (\ref{introExt})  exists and is moreover unique if one imposes the desired compatibility with the analogous homomorphism over $U_v$ (Prop. \ref{prolD1}). This verification could be achieved by brute force and explicit calculations since it is nothing more than the investigation of Kummer extensions of Dieudonn\'e crystals over $C_v$, $C^\s_v$ and $U_v$. Here we have chosen a homological approach by embedding the exact category of Dieudonn\'e crystals into the larger Abelian category of $(f,v)$-modules. This approach relies mainly on the compatibility of the Dieudonn\'e functor with Cartier duality (\cite{BBM} Chap. 5). \petit

\petit

\ref{sectct}. \emph{The comparison theorem}.

The purpose of this chapter is to prove Thm. \ref{thmA} and Thm. \ref{thmB}. The first one is proven in Sect. \ref{dopdg} and is rather independent from the rest of the text (it only uses the global construction of the syntomic complex in Sect. \ref{scotss}). The second one is then deduced by d\'evissage using the main results of the chapters 5 to 8. \petit

In Sect. \ref{dopdg}, we consider the case of $p$-divisible groups. We first consider the case of a single base scheme $X$ with local finite $p$-bases over $\Fp$. The basic idea is to reduce to the case of $\mu_{p^\infty}$ by using Cartier biduality. Following the arguments of \cite{Ba}, we begin by checking that Fontaine-Messing's exact sequence \begin{eqnarray}\label{introFMsec}\xymatrix{0\ar[r]& \mu_{p^k}\ar[rr]^{\cal L}&& \widetilde I_k^{crys}\ar[rr]^{1-{Frob\over p}}&& \O_k^{crys}\ar[r]& 0}\end{eqnarray} on $Spec(\Fp)_{syn}$ remains valid on $X_{syn}$  (Lem. \ref{lemcomp2}, \ref{lemcomp2iii}). Now if $G$ is an arbitrary $p$-divisible group on $X$ whose Cartier dual is denoted $G^*$ we show that (\ref{introFMsec}) gives rise to an exact sequence
\begin{eqnarray}\label{introFMsec2}\xymatrix{\,\,\,\,\,\,\,\,\,\,\,\, 0\ar[r]& \cal Ext^1_{X_{syn}}(G^*,\mu_{p^k})\ar[r]&  \cal Ext^1_{X_{syn}}(G^*,\widetilde I_k^{crys})\ar[r]&\cal Ext^1_{X_{syn}}(G^*,O_k^{crys})\ar[r]& 0}\end{eqnarray}   thanks to the vanishing theorem of \cite{Br}. Next we observe that the first term in (\ref{introFMsec2}) is nothing but $\smash{G_{p^k}^{**}}$, ie. $G_{p^k}$. It thus remains only to compare the second and third term to the corresponding ones in the distinguished triangle (\ref{introtdsyn}) describing $\smash{\cal S_{syn,k,X}(D_X(G))}$. For the third one, this relies directly on the fact that for certain coefficients, $\cal Ext^1$ commutes to the projection of the crystalline topos on the small syntomic topos. Next we need to check that $\smash{\cal Ext^1_{X_{syn}}(G^*,\Ga)}$ and the morphism induced by $\O_{X/(\Z/p^k)}\to i_*\Ga$ respectively coincide with $Lie(D(G))$ and the canonical morphism $D(G)\to i_*Lie(D(G))$. This is done in Lem. \ref{synsyn1} using the Cartier isomorphism on the syntomic site. Then we conclude easily from the fact that modulo the previous identifications the map induced by the Frobenius of $\O_k^{crys}$ coincides with the one induced by the Frobenius of $D_X(G)$.

Our arguments do not extend directly to the case where $X$ is a diagram (of schemes with local $p$-bases). Fortunately, the result extends nevertheless thanks to the fact (proven by the above arguments) that $\cal S_{syn,k,X}(D_X(G))$ is concentrated in degree $0$. \petit

In Sect. \ref{sasoc}, we put everything together in order to treat the case of a semi-Abelian scheme $A/C$. Namely, (\ref{introfl4}), Thm.  \ref{thmB} (applied to $H\in pdiv(J)$), (\ref{introRepsilon}) and (\ref{introisospe}) give an isomorphism \begin{eqnarray}\label{introcomp2}R\epsilon_*R\underline \G^Z A_{p^.}\simeq Rm_*Sma\cal S_{et,.,J^\s}(-Z_J)(D_J(H)_{|J^\s})\end{eqnarray}
while (\ref{introisoMV}) together with a basic property of the functor $Sma$ give \begin{eqnarray}\label{introcomp1}\cal S_{et,.,C^\s}(-Z)(D_{C^\s}(A))\simeq Rm_*Sma\cal S_{et,.,J^\s}(-Z_J)(D_{C^\s}(A)_{|J^\s})\end{eqnarray} We conclude by showing that the natural morphism $D_J(H)_{|J^\s}\to D_{C^\s}(A)_{|J^\s}$
becomes an isomorphism after applying $\cal S_{et,.,J^\s}$ and the functor $Sma$. \petit


\para A few further generalizations of such comparison results can easily be imagined. The case of finite locally free groups can be treated roughly in the same way (\cite{Pe to Va}), or alternatively as a consequence of the case of $p$-divisible groups (\cite{Va to Pe}). Other possible generalizations are log $p$-divisible  groups in the sense of \cite{Tr} or alternatively Abelian varieties with general reduction. We hope to come back to these cases in future works.

\para
Needless to say, this paper relies entirely on the previous work of K. Kato and the first author. The references to their original text \cite{KT} are so numerous that
we have given up listing them all.
We hope that the self contained approach taken here will serve as a useful reference in the future.

We are grateful to T. Tsuji for his help and suggestions, especially concerning the version of the syntomic complex involving linearized de Rham
complexes and also for the localization triangles in compactly supported log crystalline cohomology using blowing ups. We would like to thank also the anonymous referee for his careful reading and his numerous suggestions that have significantly improve the clarity of the paper.

\section{Preliminaries part I: category theoretical background} \label{secppi}

\subsection{Diagrams and (co)fibered categories} \label{dacfc} ~~ \\

We review basic facts concerning (co)fibered categories over a base category $\cal B$ and their natural extensions to the category of diagrams of $\cal B$.

\para  Recall the following definition from \cite{Gr}.

\begin{defn} \label{defdiagto} Let $\C$ denote a category.
 \debrom \item \label{defdiagtoi} A \emph{diagram $X$ of $\C$} is a functor $X:\Delta\rightarrow \C$, $\delta\mapsto X_\delta$,
$(\delta\rightarrow \delta')\mapsto (X_\delta\rightarrow X_{\delta'})$ for some small category $\Delta$ which is called \emph{the type of $X$}. We
sometimes use simplified notations such as $X/\Delta$ or $(X_\delta)$.

\item \label{defdiagtoii} Given diagrams $X$, $X'$ of respective type $\Delta$, $\Delta'$, a
\emph{morphism of diagrams $f:X\rightarrow X'$} is a couple $f=(F,\alpha)$ where $F:\Delta\rightarrow \Delta'$ is a functor and $\alpha:X\rightarrow
X'\circ F$ is a natural transformation. Using abusively the letter $f$ to designate either $f$, $F$ or $\alpha$, we sometimes write
$f_\delta:X_\delta\rightarrow \smash{X'_{f(\delta)}}$.

\item \label{defdiagtoiii} Morphisms of diagrams are composed in an obvious way and the resulting category is denoted
$Diag(\C)$. \finrom
\end{defn}

The category $Diag(\C)$  is strictly $2$-functorial with respect to the category $\C$, meaning that a functor $F:\C\rightarrow \C'$ canonically
induces a functor $Diag(\C)\rightarrow Diag(\C')$ and similarly for natural transformations between such functors, everything being strictly
compatible with composition.

We always identify $\cal C$ with a full subcategory of $Diag(\cal C)$ by sending an object to the corresponding punctual diagram. Let us emphasize that
even though diagrams of type $\Delta$ in $\cal C$ are in bijection with (and are often identified with) diagrams of type $\Delta^{op}$ in $\cal C^{op}$
there is no obvious relation between the categories $Diag(\C^{op})$ and $Diag(\C)$ or $Diag(\C)^{op}$.

\para \label{fibcof}

Let $\mathfrak Cat$ denote the $2$-category of categories and let $\cal B$ denote a category. The following result is well known.

\begin{lem} \label{eqfib} The following $2$-categories are strictly equivalent in a natural way.
\debrom
\item \label{eqfibi} The $2$-category of contravariant pseudo-functors $\cal B\to \mathfrak Cat$, pseudo morphisms and natural transformations
between them.

\item \label{eqfibii} The $2$-category of fibered categories above $\cal B$ endowed with a normalized cleavage, cartesian $\cal B$-functors and $\cal B$-natural
transformations.

\item \label{eqfibiii} The $2$-category of cofibered categories above $\cal B^{op}$ endowed with a normalized cocleavage, cocartesian
$\cal B^{op}$-functors and $\cal B^{op}$-natural transformations. \finrom
\end{lem}
\noindent Proof. This follows easily from \cite{SGA1} VI. Let us only hint the strict $2$-functors $\ref{eqfibi}\to \ref{eqfibii}$ and $\ref{eqfibi}\to \ref{eqfibiii}$ in order to fix the ideas regarding directions of arrows.  We write morphisms in $\cal B$ rather than in $\cal B^{op}$.  Start with a  contravariant pseudo-functor $\cal
F:\cal B\rightarrow \mathfrak Cat$, $X\mapsto \cal F(X)$, $f\mapsto f^*$, $(fg)^*\simeq g^*f^*$. Then:

$\ref{eqfibi}\to \ref{eqfibii}$:  The corresponding fibered category $p:\cal F_{fib}\rightarrow \cal B$ is constructed in such a way that the fiber
category $p^{-1}(X)$ above an object $X$ of $\cal B$ is $\cal F(X)$ while $Hom_{f}(\xi,\eta)=Hom_{\cal F(X)}(\xi,f^*\eta)$ if $p(\xi)=X$, $p(\eta)=Y$
and $f:X\rightarrow Y$ is a morphism in $\cal B$. The cleavage is given by the collection of tautological morphisms $f^*\eta\to \eta$.

$\ref{eqfibi}\to \ref{eqfibiii}$: The corresponding cofibered category  $q:\cal F_{cof}\rightarrow \cal B^{op}$ is constructed in such a way that the
fiber category $q^{-1}(X)$ above an object $X$ of $\cal B$ is $\cal F(X)$, while $Hom_{f}(\xi,\eta)=Hom_{\cal F(Y)}(f^*\xi,\eta)$ if $q(\xi)=X$,
$q(\eta)=Y$ and $f:X\leftarrow Y$ is a morphism in $\cal B$. The cocleavage is given by the collection of tautological morphisms $\xi\to f^*\xi$.
\begin{flushright}$\square$\end{flushright}

It is sometimes useful to notice that one can switch between the point of views \ref{eqfibii} and \ref{eqfibiii} by changing the directions of the
arrows in fibers. Namely let $\cal F:\cal B\rightarrow \mathfrak Cat$ be a contravariant pseudo-functor as in \ref{eqfibi} and consider the contravariant pseudo-functor
$\cal F^{\circ}:\cal B\rightarrow \mathfrak Cat$, defined as $X\mapsto \cal F(X)^{op}$, $f\mapsto f^*$, $g^*f^*\simeq (fg)^*$. Then we have natural
isomorphisms \begin{eqnarray}\label{chdir}(\cal F^{\circ})_{fib}/\cal B\simeq (\cal F_{cof})^{op}/\cal B&\hbox{and}&\cal (\cal F^{\circ})_{cof}/\cal
B^{op}\simeq (\cal F_{fib})^{op}/\cal B^{op}.\end{eqnarray}

Summarizing, we have found four different relative categories naturally attached to a contravariant pseudo-functor $\cal F:\cal B\to \mathfrak Cat$: two fibrations ($\cal
F_{fib}/\cal B$ and  $(\cal F^\circ)_{fib}/\cal B$) and two cofibrations ($\cal F_{cof}/\cal B^{op}$ and $(\cal F^\circ)_{cof}/\cal B^{op}$). To pick
one of these four, it is enough to specify the fiber categories and whether this is a fibration or a cofibration. We will sometimes drop the subscripts $(-)_{fib}$ and $(-)_{cof}$ from the notations without any danger of confusion.

\begin{rem} \label{remlaxmor} \debrom \item \label{remlaxmori} In practice, one sometimes has to consider $\cal B$-functors $A:\cal F_{fib}\to \cal G_{fib}$ which are not necessarily cartesian or alternatively $\cal B^{op}$-functors $B:\cal F_{cof}\to \cal G_{cof}$ which are not necessarily cocartesian. These data are no longer equivalent. More precisely:

- the data of $A$ corresponds bijectively to a \emph{lax morphism} $\cal F\to \cal G$, ie. a functor $A_X:\cal F(X)\to \cal G(X)$ for each $X$ in $\cal B$ and a natural transformation $\alpha_f:A_X f^*\to f^*A_Y$ for each $f:X\to Y$ in $\cal B$, these data being submitted to the composition constraint.

- the data of $B$ corresponds bijectively to a \emph{colax morphism} $\cal F\to \cal G$, ie. a functor $B_X:\cal F(X)\to \cal G(X)$ for each $X$ in $\cal B$ and a natural transformation $\beta_f:f^*B_Y\to B_Xf^*$ for each $f:X\to Y$ in $\cal B$, these data being submitted to the composition constraint.

\item \label{remlaxmorii} Assume that the contravariant pseudo functor $\cal F:\cal B\to \mathfrak Cat$, $X\mapsto \cal F(X)$, $f\mapsto f^*$, $(fg)^*\simeq g^*f^*$ is such that each $f^*$ admits a right adjoint, say $f_*$. Consider the contravariant pseudo-functor $\cal F':\cal B^{op}\to \mathfrak Cat$, $X\mapsto \cal F(X)$, $f\mapsto f_*$, $f_*g_*\simeq (fg)_*$ obtained by adjunction (here and in the whole remark, morphisms are considered in $\cal B$ rather than $\cal B^{op}$). We have then an isomorphism $\cal F_{cof}\simeq \cal F'_{fib}$ of $\cal B^{op}$-categories. In that situation we say that $\cal F_{cof}$ (or equivalently $\cal F'_{fib}$) is \emph{bifibered} and that the collection of adjoint pairs of functors $(f^*,f_*)$ (endowed with their adjunction morphisms $id\to f_*f^*$, $f^*f_*\to id$) together with the isomorphisms $(f^*,f_*)(g^*,g_*)\simeq ((fg)^*,(fg)_*)$ (composition of adjunction) is a \emph{bicleavage}.

\item \label{remlaxmoriii} Assume that $\cal F_{cof}$ and $\cal G_{cof}$ are bifibered. Then the set of $\cal B^{op}$-functors $B:\cal F_{cof}\to \cal G_{cof}$ is in bijection with the set of $\cal B^{op}$-functors $A:\cal F'_{fib}\to \cal G'_{fib}$. In other terms, the data of a colax morphism $(B_X:\cal F(X)\to \cal G(X),\beta_f:f^*B_Y\to B_Xf^*):\cal F\to \cal G$ corresponds bijectively to the data of a lax morphism: $(A_X:\cal F'(X)\to \cal G'(X),\alpha_f:A_Yf_*\to f_*A_X)$. To see this, set $A_X=B_X$ and use the definition of base change morphisms:  $\alpha_f$ is the composed natural transformation $A_Yf_*\to f_*f^*A_Yf_*\to f_*A_Yf^*f_*\to f_*A_Y$ where the middle arrow is induced by $\beta_f$, while $\beta_f$ is the composed natural transformation $f^*B_Y\to f^*B_Yf_*f^*\to f^*f_*B_Xf^*\to B_Xf^*$ where the middle arrow is induced by $\alpha_f$). Let us emphasize that under this bijection, the set of cocartesian $B$'s does not correspond to the set of cartesian $A$'s in general (If $A$ cartesian implies $B$ cocartesian, this might be called a \emph{base change theorem} for $A$).
\finrom
\end{rem}

\para \label{fibcofdiag}

If $X/\Delta$ is a diagram of $\cal B$ we use the simplified notations \begin{eqnarray}\label{notsectionfib}& \cal F_{fib}(X):=
\G(\smash{\Delta\times_{X,\cal B,p}\cal F_{fib}}/\Delta)\\  \label{notsectioncof} \hbox{and} &\cal F_{cof}(X):= \G(\smash{\Delta^{op}\times_{X,\cal
B^{op},q}\cal F_{cof}}/\Delta^{op})\end{eqnarray} where $\G(-)$ denotes the category of sections. These categories  identify respectively with the category of diagrams of type $\Delta$ above $X$ in
$\cal F_{fib}$ (recall that $X:\Delta\to \cal B$) and with the category of diagrams of type $\Delta^{op}$ above $X$  in $\cal F_{cof}$ (here we think of
$X:\Delta^{op}\to \cal B^{op}$). Note that for a punctual diagram $X$, we have trivially $\cal F_{fib}(X)\simeq \cal F_{cof}(X)$.


\begin{lem} \label{lemdiagcodiag} Any contravariant pseudo-functor $\cal F:\cal B\to \mathfrak Cat$ admits two natural extensions to $Diag(\cal B)$, namely:
\debrom \item \label{lemdiagcodiagi}  A contravariant pseudo-functor $\cal F^{diag}:Diag(\cal B)\to \mathfrak Cat$ satisfying $\cal F^{diag}(X)\simeq
\cal F_{fib}(X)$ for any $X$ in $Diag(\cal B)$.

\item \label{lemdiagcodiagii}  A contravariant pseudo-functor $\cal F^{codiag}:Diag(\cal B)\to \mathfrak Cat$ satisfying $\cal F^{codiag}(X)\simeq \cal F_{cof}(X)$ for any $X$ in $Diag(\cal B)$.
\finrom
\end{lem}
\noindent Proof. \ref{lemdiagcodiagi} We notice that $Diag(p):Diag(\cal F_{fib})\rightarrow Diag(\cal B)$ is automatically a fibered category and inherits a
canonical normalized cleavage from $\cal F$. Let us describe the cleavage explicitly. If $f:X/\Delta\to X'/\Delta'$ is a morphism in $Diag(\cal B)$
and $\xi'/\Delta'$ is a diagram of $\cal F_{fib}$ above $X'$ then $f^*\xi'$ is the diagram of type $\Delta$ with vertices given by the formula
$\smash{(f^*\xi')_\delta}:=\smash{f_\delta^*(\xi'_{f(\delta)})}$ and edges deduced from the edges of $\xi'$ using the universal property of the
tautological morphisms $\smash{f_\delta^*(\xi'_{f(\delta)})}\rightarrow \smash{\xi'_{f(\delta)}}$ in $\cal F_{fib}$. We define $\smash{\cal
F^{diag}}:Diag(\cal B)\rightarrow \mathfrak Cat$ as the corresponding contravariant pseudo-functor obtained by Lem. \ref{eqfib}, so that $\cal
F^{diag}(X)\simeq \cal F_{fib}(X)$ for any $X$ in $Diag(\cal B)$, as desired.

\ref{lemdiagcodiagii} Applying the previous construction to $\cal F^\circ$ instead of $\cal F$ and then reversing arrows again we find a
pseudo-functor $\cal F^{codiag}:=\cal F^{\circ,diag,\circ}:Diag(\cal B)\rightarrow \mathfrak Cat$ satisfying $$\cal F^{codiag}(X)= \cal
F^{\circ,diag}(X)^{op}\simeq (\cal F^\circ)_{fib}(X)^{op}\simeq \cal F_{cof}(X)$$ for $X$ in $Diag(\cal B)$ and we are done.

Let us describe the cocleavage of the corresponding cofibered category $(\cal F^{codiag})_{cof}/\cal B^{op}$ for the sake of clarity. If
$X/\Delta\leftarrow X'/\Delta':f$ is a morphism in $Diag(\cal B)$ and $\xi/\Delta^{op}$ is a diagram of $\cal F_{cof}$ above $X:\Delta^{op}\to \cal
B^{op}$ then $(f^*\xi)$ is the diagram of type $\Delta'^{op}$ in $\cal F_{cof}$ with vertices given by the formula
\begin{eqnarray}\label{pullbackformula}\smash{(f^*\xi)_{\delta'}}=\smash{f_{\delta'}^*(\xi_{f(\delta')})}\end{eqnarray} and edges defined in the following way: if $\delta'\leftarrow \delta'':g'$ in $\Delta'$ then the morphism \begin{eqnarray}\label{pullbackformulaedges}(f^*\xi)_{g'}:(f^*\xi)_{\delta'}\to (f^*\xi)_{\delta''} &\hbox{in $\cal
F_{cof}$}\end{eqnarray} corresponds to the composed morphism  $$\xymatrix{g'^*f_{\delta'}^*(\xi_{f(\delta')})\ar[r]^-\sim &f_{\delta''}^*f(g')^*(\xi_{f(\delta')})\ar[rr]^-{f_{\delta''}^*(\xi_{f(g')})}&&f_{\delta''}^*(\xi_{f(\delta'')})&\hbox{in $\cal F(X'_{\delta''})$}}$$ where the first isomorphism results from the commutativity of the following square $$\xymatrix{X'_{\delta'}\ar[d]_{f_{\delta'}}&X'_{\delta''}\ar[l]_-{g'}\ar[d]^{f_{\delta''}}\\ X_{f(\delta')}&\ar[l]_-{f(g')}X_{f(\delta'')}}$$
\begin{flushright}$\square$\end{flushright}

\begin{rem} \label{remdiagcodiag} Consider another contravariant pseudo-functor  $\cal G:\cal B\to \mathfrak Cat$.

\debrom \item \label{remdiagcodiagi} Every lax (resp. pseudo) morphism $\alpha:\cal F\to \cal G$ admits a canonical (resp. canonical and essentially unique) extension to a lax (resp. pseudo) morphism $\alpha:\cal F^{diag}\to \cal G^{diag}$ (use the explicit description of the cleavages on $Diag(\cal B)$).

\item \label{remdiagcodiagii} Similarly, colax (resp. pseudo) morphisms $\cal F\to \cal G$ extend canonically (resp. extend canonically and essentially uniquely) to colax (resp. pseudo) morphisms $\cal F^{codiag}\to \cal G^{codiag}$.
\finrom
\end{rem}

In the text we will use the following simplified terminology and notations. Every \emph{fibered} (resp.
\emph{cofibered}) \emph{category} will implicitly be endowed with a normalized cleavage (resp. cocleavage). In view of the Lem.
\ref{eqfib} and Lem.
\ref{lemdiagcodiag} we will often designate a fibered category $\cal G/\cal B$ or $\cal G/Diag(\cal B)$ (resp. a cofibered category $\cal G/\cal
B^{op}$ or $\cal G/Diag(\cal B)^{op}$) simply by describing the contravariant pseudo-functor $\cal F:\cal B\to \mathfrak Cat$
satisfying $\cal  F_{fib}=\cal G$ or $(\cal F^{diag})_{fib}=\cal G$ (resp. $\cal
F_{cof}=\cal G$ or $(\cal F^{codiag})_{cof}=\cal G$).


\para Let us now discuss some facts relative to change of the base category. For convenience of later reference, we look at $\cal B^{op}$-categories instead of $\cal B$-categories. Given $q_{\cal F}:\cal F\to \cal B^{op}$ and $q_{\cal G}:\cal G\to \cal B^{op}$ we denote $\smash{\cal Hom_{\cal B^{op}}}(\cal F,\cal G)$ the category where an object is a $\cal B^{op}$-functor, ie. a functor $A:\cal F\to \cal G$ such that $q_{\cal G}A=q_{\cal F}$, and where a morphism $A\to A'$ is a $\cal B^{op}$-natural transformation of functors, ie. a natural transformation $\alpha:A\to A'$ such that for all object $\xi$ of $\cal F$, $\alpha_\xi:A(\xi)\to A'(\xi)$ satisfies $q_{\cal G}(\alpha_\xi)=\smash{id_{q_{\cal F}(\xi)}}$. For any functor $\cal C\to \cal B$ and any $\cal C^{op}$-category $\cal H$ we have a natural isomorphism of categories $$\cal Hom_{\cal C^{op}}(\cal H,\cal F_{|\cal C^{op}})\simeq \cal Hom_{\cal B^{op}}(\cal H,\cal F)$$ where $\smash{\cal F_{|\cal C^{op}}}:=\smash{\cal C^{op}\times_{\cal B^{op}}\cal F}$ (strict fiber product) and where $\cal H$ is viewed as a $\cal B^{op}$-category in the obvious way on the right hand side.

\begin{lem} \label{CBdescent1} Assume that for all $X$ in $\cal B$, the fiber category $\cal C(X)$ is connected (and in particular non empty). Assume moreover that for every $x: X'\rightarrow X$ in $\cal B$ and $Y$ in $\cal C(X)$,
there exists a morphism $y:Y'\rightarrow Y$ above $x$ in $\cal C$.
Then the natural functor \begin{eqnarray}\label{bcfunct}(-)_{|\cal C^{op}}:\cal Hom_{\cal B^{op}}(\cal F,\cal G)\to \cal Hom_{\cal C^{op}}(\cal F_{|\cal C^{op}},\cal G_{|\cal C^{op}})\end{eqnarray} is fully faithful.
\end{lem}
\noindent
Proof. Consider two $\cal B^{op}$-functors $A,A':\cal F\to \cal G$. We have to show that $$\smash{Hom_{\cal B^{op}}}(A,A')\to \smash{Hom_{\cal C^{op}}(A_{|\cal C^{op}},A'_{|\cal C^{op}})}$$ is bijective. Injectivity is immediate using the fact that $\smash{\cal C^{op}\times_{\cal B^{op}}\cal F}\to \cal F$ is surjective on objects. Let us explain surjectivity. Consider  a $\cal C^{op}$-natural transformation $\alpha: \smash{A_{|\cal C^{op}}\to A'_{|\cal C^{op}}}$, ie. a collection of morphisms $$\alpha_{Y,\xi}:\smash{A_{|\cal C^{op}}(Y,\xi)=(Y,A(\xi))\to (Y,A'(\xi))=A'_{|\cal C^{op}}(Y,\xi)}$$ inducing $id_Y$ on the first component, indexed by the objects $(Y,\xi)$ of $\smash{\cal C^{op}\times_{\cal B^{op}}\cal F}$, and functorial with respect to the morphisms of the later category. By connectedness of $\cal C(q_{\cal F}(\xi))$, we note that $\alpha_{Y,\xi}$ does not depend on $Y$, but only on $\xi$. Whence a collection of morphisms $\alpha_\xi:A(\xi)\to A'(\xi)$,  indexed by the objects of $\cal F$. It remains to check that this collection is functorial with respect to $\xi$, but this, in turn, follows from the fact that given any $f:\xi'\to \xi$ in $\cal F$ there exists an object $(Y,\xi)$ of $\smash{\cal C^{op}\times_{\cal B^{op}}\cal F}$ and a morphism of the latter category extending $f$.    \begin{flushright}$\square$\end{flushright}

We say that a $\cal C^{op}$-functor $A:\smash{\cal F_{|\cal C^{op}}}\to \smash{\cal G_{|\cal C^{op}}}$
\emph{descends to a $\cal B^{op}$-functor} if it is in the essential image of $(\ref{bcfunct})$. Under the assumptions of the lemma there is then an essentially unique pair $(\overline A,\alpha)$ where $\overline A:\cal F\to \cal G$ is a $\cal B^{op}$-functor and $\alpha:\smash{\overline A_{|\cal C^{op}}}\simeq A$ is a $\cal C^{op}$-isomorphism.
We will now establish a sufficient condition for this to happen.

\begin{prop} \label{CBdescent2} Keep the assumptions of the previous lemma and assume furthermore that there exists a $\cal B$-functor $P:\cal C\times_{\cal B}\cal C\to \cal C$ together with $\cal B$-natural transformations $P\to p_i$, $i=1,2$ (here $p_i:\cal C\times_{\cal B} \cal C\to \cal C$ denote the projection functors) such that for any $X$ in $\cal B$ and any $Y$, $Y'$ in $\cal C(X)$, $P(Y,Y')$ together with its morphisms to $Y$ and $Y'$ represents the product of $Y$ and $Y'$ in $\cal C(X)$.

The  following conditions regarding a $\cal C^{op}$-functor $A:\smash{\cal F_{|\cal C^{op}}}\to \smash{\cal G_{|\cal C^{op}}}$ are equivalent:

\debrom \item \label{condCBdescenti} For every object $X$ of $\cal B$, $\xi$ of $\cal F(X)$, and every morphism $y:Y\to Y'$ of $\cal C(X)$, the image by $A$ of the morphism $(y,id_{\xi}):(Y,\xi)\to (Y,\xi)$ of $\smash{\cal C^{op}\times_{\cal B^{op}}\cal F}$ is sent to an isomorphism by the projection functor $\smash{\cal C^{op}\times_{\cal B^{op}}\cal G}\to \cal G$.

\item \label{condCBdescentii} The $\cal C^{op}$-functor $A$ descends to a $\cal B^{op}$-functor.
\finrom

\end{prop}
We begin with a preliminary lemma.

\begin{lem} \label{interCBdescent2}
\debrom \item \label{interCBdescent2i} Consider categories $\cal C$ and $\cal D$ and assume that the category $\cal C$ has finite non empty products. If $A:\cal
C^{op}\rightarrow \cal D$  is a functor which sends every arrow of $\cal C$ to an isomorphism of $\cal D$ then it is isomorphic to a constant functor.
\item \label{interCBdescent2ii}  Consider a $\cal B$-category $q_{\cal C}:\cal C\to \cal B$ satisfying the conditions of Lem. \ref{CBdescent1} and Prop. \ref{CBdescent2}. If $q_{\cal D}:\cal D\to \cal B^{op}$ is a $\cal B^{op}$-category and $A:\cal C^{op}\to \cal D$ is a $\cal B^{op}$-functor then the following conditions are equivalent:
\debrom \item \label{interCBdescent2iia} For all $X$ in $\cal B$, and all morphism $y:Y'\rightarrow Y$ in $\cal C(X)$,  the morphism $A(y):A(Y)\rightarrow A(Y')$ is
invertible in $\cal D(X)$.
\item \label{interCBdescent2iib}  There exists a section $\overline A$ of $q_{\cal D}$ and an isomorphism $\alpha:A\simeq \overline{A}\circ q_{\cal C}$.
\finrom
\finrom
\end{lem} \noindent Proof. \ref{interCBdescent2i} We have to show that up to isomorphism, the functor $A$ factors through the punctual category, or equivalently, through the category $\cal C^{rig}$ having the same set of objects than $\cal C$ but where there is exactly one arrow $Y\rightarrow Y'$ for any couple of
objects $(Y,Y')$. We will produce a functor $\overline A:\cal C^{rig,op}\to \cal D$ and an isomorphism $A\simeq \overline A\circ c$ where $c:\cal C^{op}\rightarrow\cal C^{rig,op}$ is the canonical functor.

Given $Y$ in $\cal C$, we set $\overline{A}(Y):=A(Y)$. If $Y'$ is another
object in $\cal C$ and $y$ denotes the only arrow $Y\rightarrow Y'$ in $\cal C^{rig}$, we set $\overline{A}(y):=A(p_2)^{-1}A(p_1)$, where
$$\xymatrix{Y'&\ar[l]_{p_1}Y'\times Y\ar[r]^{p_2}&Y}$$ denote the projection morphisms. Let us check that $\overline{A}$ is a functor. For
$y':Y'\rightarrow Y''$ in $\cal C^{rig}$, the equality $\overline{A}(y'y)=\overline{A}(y)\overline{A}(y')$ is ensured by the following commutative diagram of
$\cal C$:  
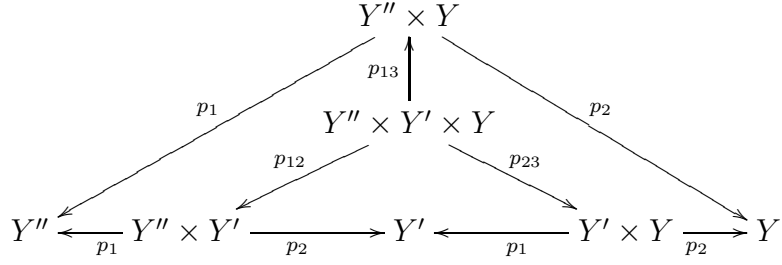
$$\xymatrix{&&Y''\times  Y\ar[lldd]_{p_{1}}\ar[rrdd]^{p_{2}}&&
\\ &&Y''\times Y'\times Y\ar[u]^{p_{13}}\ar[ld]_-{p_{12}}\ar[rd]^-{p_{23}}&&\\
Y''&Y''\times Y'\ar[l]^{p_1}\ar[r]_{p_2}&Y'&Y'\times Y\ar[l]^{p_1}\ar[r]_{p_2}&Y}$$
Next the equality $\overline{A}(id_Y)=\smash{id_{\overline A(Y)}}$ is ensured by the commutativity of the following diagram in $\cal C$: 
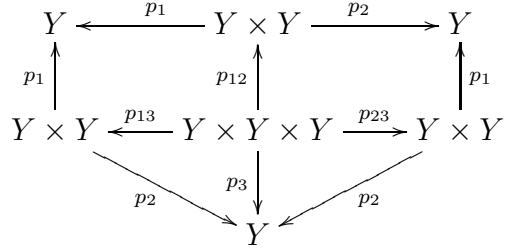
$$\xymatrix{Y&\ar[l]_{p_1}Y\times Y\ar[r]^{p_2}&Y\\
Y\times Y\ar[u]^{p_1}\ar[dr]_{p_2}&Y\times Y\times Y\ar[l]_-{p_{13}}\ar[u]^{p_{12}}\ar[r]^-{p_{23}}\ar[d]_{p_3}&Y\times Y\ar[ld]^{p_2}\ar[u]_{p_1}\\
&Y&}$$ It remains to check that $\overline{A}\circ c=A$, ie. that  $A(p_2)^{-1}A(p_1)=A(y)$
if $y:Y\rightarrow Y'$ is any arrow of $\cal C$.  Now the following commutative diagram $$
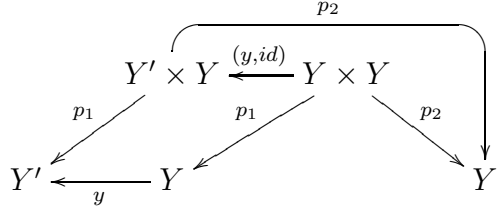
\xymatrix{&\ar[ld]_{p_1}Y'\times Y\ar @{->} @<+0pt> `u[r] `[rr]^-{p_2} [rrd]&\ar[l]_{(y,id)}\ar[ld]_{p_1}\ar[rd]^{p_2}Y\times Y&\\
Y'&\ar[l]^{y}Y&&Y}$$ of $\cal C$ reduces us to the case $Y'=Y$ and $y=id$ and thus we are done since this case already has been checked.

\ref{interCBdescent2ii} Let us prove that \ref{interCBdescent2iia} implies \ref{interCBdescent2iib}.
Thanks to \ref{interCBdescent2i} applied to the functors $A_{X}:\cal C(X)^{op}\rightarrow \cal D(X)$,  we
already know that there is a constant functor $\overline{A}_{X}$ and an isomorphism $\alpha_X:A_{X}\simeq \overline{A}_{X}$ for each object $X$ of $\cal B$. Let us denote $\overline{A}(X)$ the unique value of the functor $\overline{A}_{X}$.
We have to enrich the collection of the $\overline{A}(X)$'s into a functor $\overline{A}:\cal B^{op}\rightarrow \cal D$.
Let $x:X'\rightarrow X$ in $\cal B$ and choose a morphism $y:Y'\rightarrow Y$ in $\C$ above $x$. The following commutative diagram shows that $\overline{A}(x):=\alpha_{X',Y'}A(y)\alpha_{X,Y}^{-1}$ is well defined, ie. does not depend on the choice of $y$: 
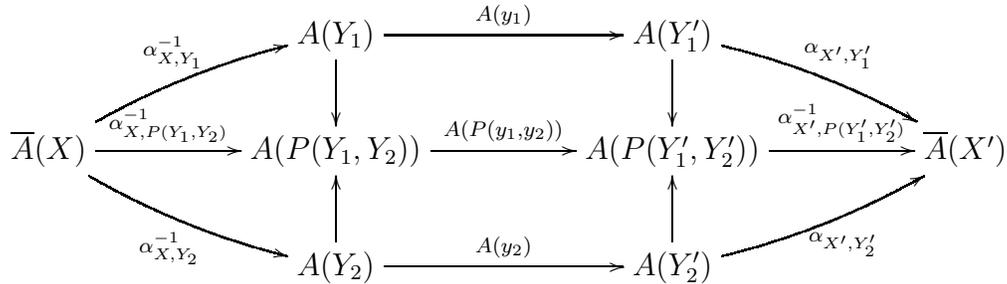
$$\xymatrix{&&A(Y_1)\ar[rr]^-{A(y_1)}\ar[d]&&A(Y_1')\ar@/^.5pc/@<-0ex>[rrd]^-{\alpha_{X',Y_1'}}\ar[d]&&\\
\overline{A}(X)\ar @/^.5pc/@<-0ex>[rru]^-{\alpha_{X,Y_1}^{-1}}\ar[rr]^-{\alpha_{X,P(Y_1,Y_2)}^{-1}}
\ar@/_.5pc/@<-0ex>[rrd]_-{\alpha_{X,Y_2}^{-1}}
&&A(P(Y_1,Y_2))\ar[rr]^-{A(P(y_1,y_2))}&&A(P(Y'_1,Y'_2))\ar[rr]^-{\alpha_{X',P(Y'_1,Y'_2)}^{-1}}&&
\overline{A}(X')
\\ &&A(Y_2)\ar[u]\ar[rr]^-{A(y_2)}&&A(Y_2')\ar@/_.5pc/@<-0ex>[rru]_-{\alpha_{X',Y_2'}}\ar[u]}$$
Since by assumption a couple of composable arrows $x:X'\rightarrow X$, $x':X''\rightarrow X'$  in $\cal B$ can always be lifted to a couple of
composable arrows $y:Y'\rightarrow Y$, $y':Y''\rightarrow Y'$ in $\cal C$ we find that $\overline{A}(xx')=\overline{A}(x') \overline{A}(x)$, so that $\overline A$ is a functor. It satisfies $q_{\cal D}\circ \overline{A}=id$ since $\alpha_X$ is above the identity and $A$ is a $\cal B^{op}$-functor. It remains to notice that, by construction of the $\overline{A}(x)$'s, the collection of isomorphisms $\smash{\alpha_{q_{\cal C}(Y),Y}}:A(Y)\simeq \overline{A}(q_{\cal C}(Y))$ is functorial with respect to $Y$, ie. gives rise to an isomorphism $A\simeq \overline{A}\circ q_{\cal C}$ as desired.
\begin{flushright}$\square$\end{flushright}

Let us now proceed with the proof of Prop. \ref{CBdescent2}. Consider the following strictly commutative diagram: $$\xymatrix{\cal Hom_{\cal B^{op}}(\cal F,\cal G)\ar[r]\ar[d]^\shortparallel & \cal Hom_{\cal C^{op}}(\cal C^{op}\times_{\cal B^{op}} \cal F,\cal C^{op}\times_{\cal B^{op}}\cal G)\ar[d]^\wr \\
\cal Hom_{\cal B^{op}}(\cal F,\cal G)\ar[r]\ar[d]^\wr& \cal Hom_{\cal B^{op}}(\cal C^{op}\times_{\cal B^{op}} \cal F,\cal G)\ar[d]^\wr \\
\cal Hom_{\cal F}(\cal F,\cal F\times_{\cal B^{op}}\cal G)\ar[r]&\cal Hom_{\cal F}(\cal C^{op}\times_{\cal B^{op}}\cal F,\cal F\times_{\cal B^{op}}\cal G)}$$
Assume that the $\cal C^{op}$-functor $A$ satisfies condition \ref{condCBdescenti}. We need to prove that $A$ is in the essential image of the top vertical arrow. Let $A_2:\smash{\cal C^{op}\times_{\cal B^{op}} \cal F\to \cal G}$  and $A'$ respectively denote the $\cal B^{op}$-functor and the $\cal F$-functor corresponding to $A$ under the right vertical isomorphisms. Then for every object $X$ of $\cal B$, $\xi$ of $\cal F(X)$, and every morphism $y:Y\to Y'$ of $\cal C(X)$, both $A_2(y,id_{\xi})$ and $A'(y,\xi)$ are isomorphisms. In particular the condition \ref{interCBdescent2iia} of Lem. \ref{interCBdescent2} is satisfied by $A'$ for $\cal C'^{op}=\smash{\cal C^{op}\times_{\cal B^{op}} \cal F}$, $\cal B'^{op}=\cal F$ and $\cal D'=\smash{\cal C^{op}\times_{\cal B^{op}}\cal G}$. It will follow that $A'$ is in the essential image of the bottom vertical arrow and the proof will be finished once checked that $\cal C'\to \cal B'$ satisfies the assumptions of Lem. \ref{interCBdescent2} \ref{interCBdescent2ii} as well. Let us check this. For $\xi$ in $\cal B'=\cal F^{op}$, the category $\cal C'(\xi)= \cal C(q_{\cal F}(\xi))\times \{\xi\}$ is connected by assumption. Next if $f:\xi'\to \xi$ in $\cal B'=\cal F^{op}$ and $(Y,\xi)$ in $\cal C'(\xi)$ it follows immediately from the assumption on $\cal C\to \cal B$ that we can find $(Y',\xi)\to (Y,\xi)$ in $\cal C'=\cal C\times_{\cal B}\cal F^{op}$ above $f$. Finally a $\cal B'$-functor $P':\smash{\cal C'\times_{\cal B'}\cal C'}\to \cal C'$ and $\cal B'$-natural transformations $P'\to p_i'$ satisfying the required conditions are obtained from $P$ and the given $\cal B$-natural transformations $P\to p_i$ by base change via $\cal F^{op}\to \cal B$.
\begin{flushright}$\square$\end{flushright}

\begin{rem} \label{remCBdescent1} Here are some technical remarks concerning Prop. \ref{CBdescent2}.

\debrom \item The condition that $P$ together with the natural transformations $P\to p_i$ represent products on the fiber categories could clearly be relaxed by asking only that the iterations of $P$ and $P\to p_i$ satisfy certain  compatibilities. We leave this generalization to the interested reader.

\item Condition \ref{condCBdescenti} can be rephrased by saying that for any morphism $y$ in any fiber category $\cal C(X)$, $A$ preserves cocartesian morphisms above $y$.
\finrom
\end{rem}

\begin{rem} \label{remCBdescent2} In practice the $\cal B^{op}$-categories $\cal F$, $\cal G$ will arise from contravariant pseudo-functors on $\cal B$ and $A$ will arise from a colax (not necessarily pseudo) morphism between contravariant pseudo-functors on $\cal C$.
\end{rem}

\subsection{Fibered topoi and derived categories} \label{ftadc}

\para \label{prelimftop}  We will use the following terminology regarding topoi. 


\begin{defn} \label{defweakmor}

\debrom

\item  \label{defweakmori} A \emph{pretopology} is a couple $(\cal C,Cov)$ where $\cal C$ is a category and $Cov$ denotes a set of families of arrows in $\cal C$ (the \emph{coverings}) defining a pretopology in the sense of \cite{SGA4-I} II, Def. 1.3.

\item  \label{defweakmorii} Consider pretopologies $(\cal C_i,Cov_i)$, $i=1,2$. A  \emph{premorphism}: $v:(\cal C_1,Cov_1)\to (\cal C_2,Cov_2)$ is a functor $V:C_1\leftarrow \C_2$ which sends $Cov_2$ into
$Cov_1$ and commutes to base change by arrows belonging to coverings of $Cov_2$.


\item  \label{defweakmoriv} Consider topoi $E_i$, $i=1,2$. A \emph{weak morphism of topoi}: $v:E_1\to E_2$ is a couple of adjoint functors $(v^{-1},v_*):E_1\to E_2$. We say that $v$ is a \emph{morphism of topoi} if $v^{-1}$ is exact. 

\finrom
\end{defn}


Morphisms of topoi are the usual notion from \cite{SGA4-I} IV. Note that premorphisms of pretopologies are called morphisms of topologies in \cite{Ar} Def. 2.4.2.  \petit

\begin{rem} \label{remweakmor} \debrom \item \label{remweakmori} Given topoi $E_i$, $i=1,2$, and a weak morphism $v:E_1\to E_2$, one can always find subcanonical sites $(\cal C_i,top_i)$ and a continuous functor $V:\cal C_2\to \cal C_1$ inducing $v$. It is moreover possible to chose $\cal C_2$ in such a way that the product (resp. fiber product) of two objects is representable; in that situation, $v^{-1}$ preserves finite products (resp. is exact) if and only if $V$ preserves the product (resp. fiber product) of two objects (this can be proven using that every object of $E_2$ is the inductive limit of objects of $\cal C_2$, resp. this is standard).
\item \label{remweakmorii} We do not know whether or not a weak morphism of topoi always arises from a premorphism of pretopologies. A morphism of topoi always does.
\finrom 
\end{rem}


\begin{lem} \label{remweakmor2} We say that a weak morphism $v:E_1\to E_2$ is ringed if one is given rings $A_i$ in the topoi $E_i$ as well as a morphism of rings $v_*A_1\leftarrow A_2$. In that situation:

\debrom \item \label{remweakmor2i} The functor $v_*:Mod(E_1,A_1)\to Mod(E_2,A_2)$ has a left adjoint, denoted $v^*$.

\item \label{remweakmor2ii} Assume that $v^{-1}$ preserves finite products (hence abelian groups, rings and modules). Then we have the familiar relation $v^*M=\smash{A_1\otimes_{A_2}v^{-1}M}$ as well as the \emph{local adjunction} isomorphism $\smash{v_*\cal Hom_{A_1}(v^*M_2,M_1)}\simeq \smash{\cal Hom_{A_2}(M_2,v_*M_1)}$.
\finrom
\end{lem}
Proof. \ref{remweakmor2i} The functor $v^*$ may be described as follows. If $V:(\cal C_2,top_2)\to (\cal C_1,top_1)$ is a continuous functor of sites inducing $v$, then $v^*M$  can be described as the  sheaf associated to the presheaf of modules $U_1\mapsto limind A_1(U_1)\otimes_{A_2(U_2)}M(U_2)$, the limit
being taken in the category of $A_1(U_1)$-modules with respect to $(U_2,U_1\to VU_2)$ running in the category $U_1/\cal C_2$.

\ref{remweakmor2ii} The claimed relation is clear. The proof of the local adjunction isomorphism given in \cite{SGA4-I} IV, Prop. 13.4 works here as well. 
\begin{flushright}$\square$\end{flushright}



\para  \label{topfib1} We will use the following variants of the notion of fibered topoi. 
 


\begin{defn} \label{deftopfib} Let $\cal B$ denote a category. 

\debrom \item \label{deftopfibi} Let $\mathfrak Pre$ denote the $2$-category of pretopologies (with premorphisms as $1$-morphisms). A \emph{prevariable pretopology $\cal P$ on $\cal B$} is a covariant pseudo-functor $\cal P:\cal B\to \mathfrak
Pre$. 

\item \label{deftopfibii} A \emph{weakly variable topos $\cal T$ on $\cal B$} is a contravariant  pseudo-functor $\cal T:\cal B\to \mathfrak Cat$ such that $\cal T(X)$ is a topos for any object $X$ in $\cal B$, and the functor $f^{-1}=\cal T(f)$ admits a right adjoint for any morphism $f$ in $\cal B$. We say that $\cal T$ is \emph{variable} if the $\cal T(f)$'s are exact. 
\finrom
\end{defn}

\begin{rem} \label{remtopfib} \debrom \item \label{remtopfibi} To a prevariable pretopology $\cal P$ on $\cal B$ is naturally associated a weakly variable topos $\cal P\tilde{}$ on $\cal B$.
\item \label{remtopfibii} If $\cal B$ is small then every variable topos on $\cal B$ arises from a prevariable pretopology on $\cal B$. \finrom
\end{rem}


\begin{lem} \label{lemfibtop} Let $\cal T:\cal B\to \mathfrak Cat$ be a weakly variable topos and consider the associated cofibered categories $\cal T_{cof}/\cal B^{op}$ and $(\cal T^{codiag})_{cof}/Diag(\cal B)^{op}$ (Lem. \ref{eqfib}, Lem. \ref{lemdiagcodiag} \ref{lemdiagcodiagii}).

\debrom
\item \label{lemfibtopi} If $X=(X_\delta)$ is a diagram of $\cal B$ then $\cal T(X):=\cal T_{cof}(X)$ is a topos. Whence a weakly variable topos $\cal T^{codiag}$ on $Diag(\cal B)$, whose $\cal T(f)$'s are described by the formulae (\ref{pullbackformula}), (\ref{pullbackformulaedges}). 
    If $\cal T$ is in fact variable then $\cal T^{codiag}$ is variable as well.

\item \label{lemfibtopii} The cofibered categories $\cal T_{cof}/\cal B^{op}$ and $\cal (\cal T^{codiag})_{cof}/Diag(\cal B)^{op}$ are in fact bifibered. The right adjoint $f_*$ of the functor $f^{-1}:=\cal T(f)$ associated to a morphism of diagrams of $\cal B$, $f:X/\Delta\rightarrow X'/\Delta'$, is explicitly described as follows. For $\xi\in \cal T(X)$ the \emph{projective limit formula},
\begin{eqnarray} \label{projformula}
(f_*\xi)_{\delta'}=\limp\ \alpha_*f_{\delta,*}\xi_{\delta}\end{eqnarray}
where the projective limit is indexed by the category $\Delta/\delta'$ (where an object is a couple $(\delta,\alpha:f(\delta)\to \delta')$), describes the vertices of $f_*\xi$. Its edges are deduced from those of $\xi$ using the covariance of the category $\Delta/\delta'$ with respect to $\delta'$.

\finrom
\end{lem}
Proof. \ref{lemfibtopi} The first statement results easily from Giraud's criterion (\cite{SGA4-I} IV, Thm. 1.2). The second statement is clear from the explicit description of the $\cal T(f)$'s.

\ref{lemfibtopii} The fact that $\cal T_{cof}/\cal B^{op}$ is bifibered follows from
the definition of a weakly variable topos. This provides $f_*$ for $f$ a morphism in $\cal B$. The reader can easily check by himself that the stated projective limit formula defines the claimed $f_*$ for $f$ a morphism of $Diag(\cal B)$
\begin{flushright} $\square$\end{flushright}


In order to make the connection with the terminology of \cite{SGA4-II} Vbis, VI, let us mention that for a variable topos $\cal T:\cal B\to \mathfrak Cat$ the category
$\cal T_{cof}/\cal B^{op}$ would be called a $\cal B^{op}$-topos in \cite{SGA4-II} Vbis, Def. 1.2.1, $\cal T_{fib}/\cal B$ would be called a fibered topos in \cite{SGA4-II} VI, Def. 7.1.1 and $\cal T_{cof}(X)$ is equivalent to the total
topos associated to the fibered topos $\Delta\times_{X,\cal B,p}\cal T_{fib}/\Delta$ (\cite{SGA4-II} VI, Prop. 7.4.7).

%
%
%

\para \label{Dfib} Let us discuss the functoriality of derived categories of modules.  We say that $(\cal T,A)$ is a ringed weakly variable topos on $\cal B$ (resp. that $(\cal P,A)$ is a ringed prevariable pretopology on $\cal B$) if $\cal T$ is as in Def. \ref{deftopfib} \ref{deftopfibii} (resp. if $\cal P$ is as in Def. \ref{deftopfib} \ref{deftopfibi}) and $A$ is a ring of the total topos $\G(\cal B^{op},\cal T)$ (resp. $\G(\cal B^{op},\cal P\tilde{})$).

\begin{lem} \label{lemDfib1} Consider a ringed weakly variable topos $(\cal T,A)$ on $\cal B$.

\debrom \item \label{lemDfib1i} There is a canonical covariant pseudo-functor $Mod(\cal T(-),A):\cal B\to \mathfrak Cat$, $X\mapsto  Mod(\cal T(X),A_X)$ (the category of modules),
$(f:X\to Y)\mapsto (f_*:Mod(\cal T(X),A_X)\to Mod(\cal T(Y),A_Y))$, $f_*g_*\simeq (fg)_*$. The associated
fibered category $Mod(\cal T(-),A)_{fib}/\cal B^{op}$ is in fact bifibered, a bicleavage being given by ``the'' natural adjunctions $(f^*,f_*)$. Similar statements hold for $X\mapsto Kom(\cal T(-),A))$ (the category of complexes of modules).  

\item \label{lemDfib1ii} If $\cal T$ is associated to a prevariable pretopology then right derivation gives rise to a covariant pseudo-functor $D^+(\cal T(-),A):\cal B\to \mathfrak Cat$, $X\mapsto D^+(\cal T(X),A_X)$ (the derived category of modules), $f\mapsto Rf_*$, $Rf_*Rg_*\simeq R(fg)_*$. If the $Rf_*$'s are of finite cohomological dimension then similar statements hold with $D$ instead of $D^+$.
\item \label{lemDfib1iii} If $\cal T$ is variable then left derivation gives rise to a contravariant pseudo-functor $D^-(\cal T(-),A)$, $X\mapsto D^-(\cal T(X),A_X)$, $f\mapsto Lf^*$, $Lg^*Lf^*\simeq L(fg)^*$. If the $Lf^*$'s are of finite homological dimension then a similar statement holds with $D$ instead of $D^-$.
\item \label{lemDfib1iv} Assume that $\cal T$ is variable and that the $Rf_*$'s are of finite cohomological dimension (resp. that the $Lf^*$'s are of finite homological dimension). Then the cofibered category $D^-(\cal T(-),A)_{cof}/\cal B^{op}$ (resp. the fibered category $D^+(\cal T(-),A)_{fib}/\cal B^{op}$)  is in fact bifibered, a bicleavage being given by natural adjunctions $(Lf^*,Rf_*)$.
\finrom
\end{lem}
Proof. This is standard. Statement \ref{lemDfib1i} is abstract nonsense using Lem. \ref{remweakmor2} \ref{remweakmor2i}. The first statement of \ref{lemDfib1ii} is proved using Cech cohomology  (see \cite{Ar} 2.4.5 or the proof of \cite{SGA4-II}, V, Prop. 4.9). The first statement of  \ref{lemDfib1iii} follows from \cite{SGA4-II}, V, Cor. 1.7.1. The last statement of \ref{lemDfib1ii} and \ref{lemDfib1iii}  is
proven using an easy truncation trick. Finally, \ref{lemDfib1iv} is a consequence of the \emph{trivial duality theorem} of \cite{SGA4-III} XVII, 4.1.3.
\begin{flushright}$\square$\end{flushright}


\begin{rem} \debrom \item 
In the above lemma, we have used the traditional boundedness conditions since it will be enough for our purpose. The reader may get rid of them by using the theory
of unbounded complexes (\cite{Sp}, \cite{KS}). 

\item In virtue of \ref{lemfibtop} \ref{lemfibtopi}, the whole lemma (including the unbounded variant) extends immediately to $Diag(\cal B)$ in the variable case.
\finrom
\end{rem}    
    
In the next lemma ,we introduce some acyclicity conditions which will later be helpful to compute derived functors in the case of diagrams (and incidentally, give partial results in the prevariable case).

\begin{lemdefn} \label{acycf} Consider a ringed weakly variable topos $(\cal T,A)$ on $\cal B$ and consider $(\cal T^{codiag},A)$, the associated ringed weakly variable topos on $Diag(\cal B)$ (see Lem. \ref{lemfibtop}).
 Let $X/\Delta$ in $Diag(\cal B)$.

\debrom

\item \label{acycfi} A module $M\in Mod(\cal T_{cof}(X),A_X)$ is flat if and only its components $M_\delta$ are flat.

\item \label{acycfii} Assume that $\cal T$ is variable.  If $f:X/\Delta\to X'/\Delta'$ is a morphism in $Diag(\cal B)$, then  $Lf^*$ can be computed componentwise, ie. $(Lf^*M)_\delta\simeq Lf_\delta^*M_{f(\delta)}$. 

\item \label{acycfiii} We say that $M\in Mod(\cal T_{cof}(X),A_X)$ has \emph{injective} (resp. \emph{flasque}) \emph{components}, or that it is \emph{componentwise injective} (resp. \emph{flasque}) if $M_\delta$ is injective (resp. flasque in the sense of  \cite{SGA4-II} V, Def. 4.1) for each $\delta$ in $\Delta$.  If now $\cal T\simeq \cal P\tilde{}$, with $\cal P$ some prevariable pretopology, we say that $M$ has \emph{$\cal P$-acyclic components}, or that it is \emph{componentwise $\cal P$-acyclic} if each $M_\delta$ is $\cal P(X_\delta)$-acyclic in the sense of \emph{loc. cit.}.

\item \label{acycfiv} We say that $M\in Mod(\cal T_{cof}(X),A_X)$ is \emph{d-injective} (resp.  \emph{d-flasque}, resp. \emph{d-$\cal P$-acyclic}, assuming $\cal T\simeq \cal P\tilde{}$ with $\cal P$ a prevariable pretopology) if it is isomorphic to $f_*M_1$ for some morphism $f:X_1/\Delta_1\to X/\Delta$ with $\Delta_1$ discrete and some module $M_1$ over $(\cal T_{cof}(X_1),A_{X_1})$ which is componentwise injective (resp. flasque, resp. $\cal P$-acyclic). The category of d-injective (resp. d-flasque, resp. d-$\cal P$-acyclic) modules is cogenerating and stable by direct images of arbitrary morphisms of diagrams.

\item \label{acycfv} Assume that $\cal T$ is variable (resp. that $\cal T$ is equivalent to $\cal P\tilde{}$ for some  prevariable pretopology $\cal P$).
If $M\in Mod(\cal T_{cof}(X),A_X)$ is d-flasque (resp. d-$\cal P$-acyclic)  then it is componentwise flasque (resp. $\cal P$-acyclic). Injectives are direct factors of d-injectives and are thus componentwise flasque (resp. $\cal P$-acyclic) as well.

\item \label{acycfvi} Assume that $\cal T$ is variable (resp. that $\cal T$ is equivalent to $\cal P\tilde{}$ for some  prevariable pretopology $\cal P$). If $f:X/\Delta\to X'/\Delta'$ is a morphism in $Diag(\cal B)$, then $Rf_*:D^+(\cal T_{cof}(X),A_X)\to D^+(\cal T_{cof}(X'),A_{X'})$ can be computed using resolutions by d-flasques (resp. d-$\cal P$-acyclics). In particular if $f':X'/\Delta'\to X''/\Delta''$ is another morphism in $Diag(\cal B)$, then $Rf'_*Rf_*\simeq R(f'f)_*$.

    In the particular case where $\Delta'=\Delta$ and $f$ induces the identity on $\Delta$ then $Rf_*$ can in fact be computed using resolutions by componentwise flasques (resp. componentwise $\cal P$-acyclics). \finrom
\end{lemdefn}

\noindent Proof. \ref{acycfi} and \ref{acycfii} are obvious from the fact that tensor products can be computed componentwise.

\ref{acycfiv} Stability by
direct images is tautological. The cogeneration statement needs only to be proven for d-injectives. If $\Delta$ itself is discrete, this results from the fact that injectives are cogenerating. In the general case, let $f:X_1/\Delta_1\to X/\Delta$ denote the inclusion of the discrete diagram underlying $X$.
Then $f_*$ preserves injectives (because $f$ is a flat morphism of ringed topoi, even though $\cal T$ is only assumed weakly variable),
 $M\to f_*f^*M$ is monomorphic (because $f^*$ is faithful) for any $M$ and the result follows.

\ref{acycfv}
Let $f:X_1/\Delta_1\to X/\Delta$ be a morphism in $Diag(\cal B)$ with $\Delta_1$ discrete. We have to prove that componentwise $\cal P$-acyclicity and componentwise
flasqueness is preserved by $f_*$. It suffices to treat the case of $\cal P$-acyclicity. Let us thus consider $M_1$ with $\cal P$-acyclic components
on $(\cal T_{cof}(X_1),A_{X_1})$. Since $\Delta_1$ is discrete, the projective limit formula of Lem. \ref{lemfibtop} \ref{lemfibtopii} describing the $\delta$'s component
of $f_*M_1$ boils down to a product of $f_{\delta_1,*}M_{\delta_1}$'s. Now each term of this product is $\cal P$-acyclic (\cite{Ar} Lem. 2.4.5 or the proof
of \cite{SGA4-II} V, Prop. 4.9)  and one may easily conclude by using the fact that for any set $I$,  the derived functors $R^q\prod_I$, $q\ge 1$ vanish on families of $\cal P$-acyclics indexed by $I$ (see \cite{SGA4-II} Vbis, Prop. 1.3.10 for details in the flasque case; the $\cal P$-acyclic case is similar).
The statement about injectives is clear since $d$-injectives are cogenerating.

\ref{acycfvi} It suffices to treat the prevariable case. Consider a diagram $X:\Delta\to \cal B$ and a functor $f_1:\Delta_1\to \Delta$. Let us denote $X_{|\Delta_1}$ the diagram $X\circ f_1$. Then we have a natural morphism $f_1:X_{|\Delta_1}\to X$ induced by the identity of the vertices $X_{f_1(\delta_1)}$'s. A morphism of this form (for some functor $f_1$) will be called a \emph{change of indices}. Note that for such $f_1$, the induced weak morphism $(\cal T(X_1),A_{X_1})\to (\cal T(X),A_X)$ is in fact a flat morphism of ringed topoi.

Let us begin with the last statement of the proposition. Since $f$ induces the identity on $\Delta$, we have  $(f_*M)_{\delta}=f_{\delta,*}M_\delta$. Now injectives have flasque components (see  \ref{acycfv}) and it follows that $Rf_*$ can be computed componentwise as well. Using this, we find  that $Rf_*$ vanishes on componentwise $\cal P$-acyclics and the statement follows.

Let us now turn to the statement for $f$ arbitrary. We begin with two facts.
\petit

\noindent \emph{Fact 1}. Consider a change of indices $f_1:X_1\to X$
where the type $\Delta_1$ of $X_1$ is discrete. If $M_1$ is a componentwise $\cal P$-acyclic module on $(\cal T_{cof}(X_1),A_{X_1})$ then $f_{1,*}M_1$ is $f_*$-acyclic.
\petit

  Since $f_1:\cal T_{cof} (X_1)\to \cal T_{cof}(X)$ is a morphism of topoi, we have $Rf_*Rf_{1,*}M_1\simeq R(ff_1)_*M_1$. Now the arguments of \ref{acycfv} show that $R^qf_{1,*}M_1$ and $R^q(ff_1)_*M_1$ vanish for $q\ge 1$ and the fact follows.
 \petit


\noindent \emph{Fact 2}. If $\smash{X\mathop\rightarrow\limits^{a}X'_{|\Delta} \mathop\rightarrow\limits^{b} X'}$ is the natural factorization of the morphism $f:X/\Delta\to X'/\Delta'$ then  $Rf_*\simeq Rb_*Ra_*$. \petit

Let $\Delta_1$ be the discrete category underlying $\Delta$ and consider the following commutative diagram \begin{eqnarray}\label{diagacycfvi}
\xymatrix{X_{|\Delta_1}\ar[r]^-{a_1}\ar[d]^-{f_1}&(X'_{|\Delta})_{|\Delta_1}\ar[d]^{f_1'}\\
X\ar[r]^-a&X'_{|\Delta}\ar[r]^-b&X'}\end{eqnarray}
Since the category of modules of the form $f_{1,*}M_1$ with $M_1$ injective is cogenerating in $Mod(\cal T_{cof}(X),A_X)$ (see the proof of \ref{acycfiv}), it is sufficient to check that $Rf_*f_{1,*}M_1\simeq Rb_*Ra_*f_{1,*}M_1$ for $M_1$ injective. This, in turn, is proven by the following series of isomorphisms: \begin{eqnarray}\label{ac1}Rb_*Ra_*f_{1,*}M_1&\simeq &Rb_*a_*f_{1,*}M_1\\ \label{ac2}&\simeq &Rb_*f_{1,*}'a_{1,*}M_1\\ \label{ac3}&\simeq & b_*f_{1,*}'a_{1,*}M_1\\ \label{ac4}&\simeq & f_*f_{1,*}M_1\\ \label{ac5} &\simeq & Rf_*f_{1,*}M_1\end{eqnarray} The isomorphisms (\ref{ac1}) and (\ref{ac5})  are due to the fact that $f_{1,*}$ preserves injectives (because $f_1:(\cal T(X_1),A_{X_1})\to (\cal T(X),A_X)$ is a flat morphism of topoi). The isomorphisms (\ref{ac2}) and (\ref{ac4}) are by commutativity of the diagram (\ref{diagacycfvi}). The isomorphism (\ref{ac3}) proceeds from fact 1 since $f_1'$ is a change of indices and $a'_{1,*}M_1$ has flasque components. \petit

Let us now prove that d-$\cal P$-acyclics are $f_*$-acyclic. Consider a diagram $X_1$ whose type $\Delta_1$ is discrete and let $f_1:X_1\to X$ denote an arbitrary morphism. As before we have a natural diagram
\begin{eqnarray}\label{diagacycfvibis}
\xymatrix{X_1\ar[r]^-{a_1}\ar[d]^-{f_1}&X'_{|\Delta_1}\ar[d]^{f_1'}\\
X\ar[r]^-a&X'_{|\Delta}\ar[r]^-b&X'}\end{eqnarray}
Let $M_1$ be a componentwise $\cal P$-acyclic module on $(\cal T_{cof}(X_1),A_{X_1})$. The following series of isomorphisms prove the $f_*$-acyclicity of $f_{1,*}M_1$:
\begin{eqnarray}\label{ac6} Rf_*f_{1,*}M_1&\simeq & Rb_*Ra_*f_{1,*}M_1\\ \label{ac7}&\simeq & Rb_*a_*f_{1,*}M_1\\ \label{ac8}&\simeq & Rb_*f'_{1,*}a_{1,*}M_1\\ \label{ac9}&\simeq & b_*f'_{1,*}a_{1,*}M_1\end{eqnarray}
The isomorphism (\ref{ac6}) follows from fact 2. The isomorphism (\ref{ac7}) follows from the fact that $f_{1,*}M_1$ has $\cal P$-acyclic components and $a$ induces the identity on $\Delta$. The isomorphism (\ref{ac8}) is obvious and (\ref{ac9}) follows from fact 1 since $f_1'$ is a change of indices while $a_{1,*}M_1$ has $\cal P$-acyclic components.
\begin{flushright}$\square$\end{flushright}

\para Let us now discuss the functoriality of $Mod(\cal T^{codiag}(-),A)$ and derived categories with respect to  $(\cal T,A)$. 

\begin{defn} \label{defmorT}
\debrom \item \label{defmorTi} A \emph{premorphism} $g:\cal P\to \cal P'$ between prevariable pretopologies on $\cal B$ is a collection of
premorphisms $g_X:\cal P(X)\to \cal P'(X)$ together with compatibility isomorphisms $g_{X'}\cal P(f)\simeq \cal P'(f)g_X$ satisfying the composition
constraint.

\item \label{defmorTii}
A \emph{weak morphism} $g:\cal T\to \cal T'$ between weakly variable topoi on $\cal B$ is a collection of functors $g_X^{-1}:\cal T(X)\leftarrow \cal
T'(X)$ having right adjoints together with compatibility isomorphisms $f^{-1}g_{X'}^{-1}\simeq g_X^{-1}f^{-1}$ satisfying the composition constraint.
If the $g_X^{-1}$'s are exact, we say that $g$ is a \emph{morphism}. \finrom
\end{defn}

If $(\cal P,A)$ and $(\cal P',A')$ are ringed prevariable pretopologies we define a premorphism $(\cal P,A)\to (\cal P',A')$ as the data of a
premorphism $g:\cal P\to \cal P'$ together with a ring morphism $A'\to g_*A$. We use a similar terminology for weak morphisms between weakly ringed topoi.

\begin{lem} \label{funcT}
Any weak morphism $g:(\cal T,A)\to (\cal T',A')$ of ringed weakly variable topoi on $\cal B$ naturally extends to a weak morphism $\cal (T^{codiag},A)\to \cal (T'^{codiag},A')$ of ringed weakly variable topoi on $Diag(\cal B)$ satisfying the following properties. 

\debrom
\item \label{funcTi} For any $X/\Delta$ in $Diag(\cal B)$, the induced couple of adjoint functors $(g_X^*,g_{X,*}):Mod(\cal T_{cof}(X),A_X)\to Mod(\cal T'_{cof}(X),A'_{X})$ satisfies  $(g_X^*M')_\delta=g_{X_\delta}^*M'_\delta$ and $(g_{X,*}M)_{\delta}=g_{X_\delta,*}M_{\delta}$.

\item \label{funcTii} If $f:X/\Delta\to X'/\Delta'$ is a morphism in $Diag(\cal B)$, there are natural isomorphisms $g_{X',*} f_*\simeq f_*g_{X,*}$ and  $g_X^*f^*\simeq f^*g_{X'}^*$. These isomorphisms satisfy the composition constraint for composable $f$'s. 
\finrom 






\end{lem} Proof. The obvious extension of $g$ (hinted in Rem. \ref{remdiagcodiag}) is defined componentwise, ie. satisfies $(g_X^{-1}\xi)_\delta=\smash{g_{X_\delta}^{-1}\xi_\delta}$ for any $\delta\in \Delta$, $\xi\in \cal T_{cof}(X)$ and $X/\Delta$ in $Diag(\cal B)$. That it defines a weak morphism of weakly variable topoi means that there are natural isomorphisms $g_X^{-1}f^{-1}\simeq f^{-1}g_{X'}^{-1}$ (which follow immediately from  the pullback formula (\ref{pullbackformula})) satisfying the composition constraint. Statements \ref{funcTi} and \ref{funcTii} follow formally using Lem. \ref{remweakmor2} \ref{remweakmor2i}. 
\begin{flushright}$\square$\end{flushright}

Let us turn to derived categories.

\begin{lem} \label{DfuncT}
Consider a weak morphism of ringed weakly variable topoi $g:(\cal T,A)\to (\cal T',A')$ and assume either $(A)$ that it is associated to a  premorphism of ringed pretopologies $(\cal P,A)\to (\cal P',A')$, or $(B)$ that $\cal T$, $\cal T'$ are variable and that $g$ is a morphism. 

\debrom \item  \label{DfuncTi} Let $X/\Delta \in Diag(\cal B)$.

- In case $(A)$, the functor $Rg_{X,*}:D^+(\cal T_{cof}(X),A_X)\to D^+(\cal T_{cof}(X),A'_X)$ can be computed componentwise ie.  $(Rg_{X,*}M)_\delta \simeq
Rg_{X_\delta,*}M_\delta$.

- In case $(B)$, the functor $Lg_X^*:D^-(\cal T_{cof}'(X),A'_X)\to D^-(\cal T_{cof}(X),A_X)$ can be computed componentwise ie. $(Lg_X^*M)_\delta\simeq
Lg_{X_\delta}^*M_\delta$.

\item \label{DfuncTii}  Let $f:X/\Delta\to X'/\Delta'$ be a morphism in $Diag(\cal B)$.

    - In case $(A)$, there is a canonical isomorphism $Rg_{X',*}Rf_*\simeq Rf_*Rg_{X,*}$. These isomorphisms satisfy the composition constraint for composable $f$'s, ie. give rise to a pseudo-morphism $Rg_*:D^+(\cal T^{codiag}(-),A)\to D^+(\cal T'^{codiag}(-),A')$ between covariant pseudo-functors on $Diag(\cal B)$.

    - In case $(B)$, there is a canonical isomorphism $Lg^*_{X}Lf^*\simeq Lf^*Lg_{X'}^*$. These isomorphisms satisfy the composition constraint for composable $f$'s, ie. give rise to a pseudo-morphism $Lg^*:D^-(\cal T'^{codiag}(-),A')\to D^-(\cal T^{codiag}(-),A)$  between contravariant pseudo-functors on $Diag(\cal B)$.

\item \label{DfuncTiii} Consider another morphism $g':(\cal T',A')\to (\cal T'',A'')$ between ringed variable topoi on $\cal B$.

    - In case $(A)$, there is a canonical isomorphism $Rg'_*Rg_*\simeq R(g'g)_*$. When $g$ and $g'$ vary, the resulting collection of isomorphisms satisfy the associativity constraint.

    - In case $(B)$, there is a canonical isomorphism $L(g'g)^*\simeq Lg^*Lg'^*$. When $g$ and $g'$ vary, the resulting collection of isomorphisms satisfy the associativity constraint.
\finrom

\end{lem}
Proof. \ref{DfuncTi} The statement about $Lg_X^*$ follows from Lem.-Def. \ref{acycf} \ref{acycfi}. The statement about $Rg_{X,*}$ follows from the fact that injectives are componentwise $\cal P$-acyclic (Lem.-Def. \ref{acycf} \ref{acycfv}). The statements of \ref{DfuncTiii} follow.

\ref{DfuncTii} The statement about derived inverse images is clear using flat resolutions. The statement about derived direct images follows from the fact that $d$-$\cal P$-acyclics are preserved by $f_*$ and $g_{X,*}$, as well as $f_*$-acyclic by Lem. + Def. \ref{acycf} \ref{acycfvi} and componentwise $\cal P$-acyclic (hence $g_{X,*}$-acyclic) by Lem. + Def. \ref{acycf} \ref{acycfv}.
\begin{flushright}$\square$\end{flushright}





%
%

\para \label{projtop} Let us discuss projective systems.

\begin{defn} \label{defprojtop} Let $\N$ denote the set of positive integers viewed as a category with exactly one arrow $k\to k'$ if $k\ge k'$. Given a weakly variable topos $\cal T$ over $\cal B$,  we define $\cal T^\N$ as follows. For $X$ in $\cal B$, $\cal T^\N(X)$ is the category  of projective systems $k\mapsto \xi_k$ (we also use the notation $(\xi_k)$ or $\xi_.$) indexed by $\N$. For $f:X\to Y$, $f^{-1}:\cal T^\N(Y)\to \cal T^\N(X)$ is defined componentwise: $(f^{-1}\xi_.)_k=f^{-1}\xi_k$.
\end{defn}

Note that direct images can be computed componentwise  as well: $(f_*\xi_.)_k=f_*\xi_k$. Here are some immediate properties of $\cal T^\N$.

\begin{lem} \label{lemprojtop} Consider a weakly variable topos $\cal T$ on $\cal B$.  \debrom
\item \label{lemprojtopi} There is a natural isomorphism $(\cal T^\N)^{codiag}\simeq (\cal T^{codiag})^\N$ of weakly variable topoi on $Diag(\cal B)$.
\item \label{lemprojtopii} If $\cal T$ is associated to a prevariable pretopology then so is $\cal T^\N$.
\item \label{lemprojtopiii} The functor $\xi_.\mapsto \xi_k$ defines a morphism of weakly variable topoi $\iota_k:\cal T\to \cal T^\N$.
\item \label{lemprojtopiv} The functor $\xi\mapsto (k\mapsto \xi)$ defines a morphism of weakly variable topoi $l:\cal T^\N\to \cal T$.
\finrom
\end{lem}
Proof. Only \ref{lemprojtopii} deserves an explanation. If $\cal T(X)\simeq \cal P(X)\tilde{}$ then $\cal T(X)^\N$ is naturally equivalent to the category of sheaves on $\cal P(X)\times \mathbb N^{op}$ endowed with the pretopology for which $Cov(U,k)$ is the set of families $((U_i,k)\to (U,k))_i$ with $(U_i\to U)_i$ in $Cov(U)$. The statement results from this.
\begin{flushright}$\square$\end{flushright}

The following lemma will often be used implicitly.
\begin{lem} \label{lemprojcomp} Let $g:\cal T\to \cal T'$ denote a weak morphism of weakly variable topoi on $\cal B$. Consider rings $A_.$, $A'_.$ in the respective total topos of $\cal T^\N$, $\cal T'^\N$ and a morphism of rings $A'_.\to g_*A_.$.

\debrom \item \label{lemprojcompi}  Consider a morphism $f:X/\Delta\to X'/\Delta'$ in $Diag(\cal B)$. The functor $Rf_*:D^+(\cal T^\N_{cof}(X),A_{.,X})\to D^+(\cal T^\N_{cof}(X'),A_{.,X'})$ can be computed $k$ by $k$: $\iota_k^{-1}Rf_*M_.\simeq Rf_*M_k$.

\item \label{lemprojcompii} Consider $X/\Delta$ in $Diag(\cal B)$. The functor $Rg_{X,*}:D^+(\cal T^\N_{cof}(X),A_{.,X})\to D^+(\cal T'^\N_{cof}(X),A'_{.,X})$ can be computed $k$ by $k$: $\iota_k^{-1}Rg_{X,*}M_.\simeq Rg_{X,*}M_k$.
\finrom
\end{lem}
Proof. We may view $\cal T(X)^\N$ as the total topos of the constant variable topos $k\mapsto \cal T(X)$ on $\N^{op}$. With this in mind, \ref{lemprojcompi} (resp.  \ref{lemprojcompii}) would be a formal consequence of Lem. \ref{DfuncT} \ref{DfuncTi} (resp. Lem. \ref{DfuncT} \ref{DfuncTii}) under the additional assumption that $g$ is induced by a premorphism of pretopologies. The result still holds without this assumption however, because the category $\cal I$ of modules $M_.$ with each $M_k$ injective has the following properties: 1) it is cogenerating (use the monomorphism $M_.\hookrightarrow (k\mapsto \prod_{k'\le k}M_{k'})$) and 2) bounded below complexes with objects in $\cal I$ which are acyclic are sent by $f_*$ and $g_{X,*}$ to acyclic complexes.
\begin{flushright}$\square$\end{flushright}


\subsection{Usual sites for log schemes} \label{usfls}

\para  \label{usualtop} Let $top$ denote a property for morphisms of schemes which is a stable under isomorphism, base change and composition. We say that \emph{$top$ has fiber products} if $X\times_Y Z$ is $top$ over $S$ whenever $X/S$, $Y/S$ and $Z/S$ are $top$.

We use the notation $\cal Sch^\s$ for the category of fine log schemes as defined in  \cite{Ka2} and we look at $\cal Sch$ as a full subcategory of $\cal Sch^\s$ by endowing any scheme with the trivial log structure. For $X^\s$ in $\cal Sch^\s$, we denote $X$ the underlying scheme and $M_X\to \cal O$ the log structure. Finally we say that a morphism of $\cal Sch^\s$ is $top$ if the underlying morphism of $\cal Sch$ is $top$.

\begin{defn} \label{defusualtop} Consider $X^\s$ in $Sch^\s$.
\debrom
\item \label{defusualtopi} The \emph{$\s$-big $top$ pretopology} $TOP^\s(X^\s)$ is the category $Sch^\s/X^\s$ together with surjective families of strict $top$  morphisms as coverings. The associated topos  $X^\s_{TOP^\s}$ is called the \emph{$\s$-big topos of $X^\s$}.

\item \label{defusualtopii} The \emph{big $top$ pretopology} $TOP(X^\s)$ is the full subcategory of $Sch^\s/X^\s$ formed by the strict $Y^\s/X^\s$'s     together with surjective families of strict $top$  morphisms as coverings. The associated topos $\smash{X^\s_{TOP}}$ is called the \emph{big $top$ topos of $X^\s$}.

\item \label{defusualtopiii} The \emph{small $top$ pretopology} $top(X^\s)$ is the full subcategory of $Sch^\s/X^\s$ formed by the strict $top$ $Y^\s/X^\s$'s     together with surjective families of strict $top$  morphisms as coverings. The associated topos $\smash{X^\s_{top}}$ is called the \emph{small $top$ topos of $X^\s$}.
\finrom
\end{defn}

The $\s$-big, big and  small $top$ pretopologies are pseudo-functorial with respect to $X^\s$ in the sense that they give rise to prevariable pretopologies, and thus to weakly variable topoi  $\cal T=\smash{(-)}_{TOP^\s}$, $(-)_{TOP}$ and $(-)_{top}$  on $\cal B=\cal
Sch^\s$. The pullback functors are exact in the $\s$-big and big case. In the small case, they are known to preserve finite products in general (Rem. \ref{remweakmor} \ref{remweakmori}), and they are exact if $top$ has fiber products.

The inclusion functors may be viewed as premorphisms of pretopologies and thus induce weak morphisms of weakly variable topoi
\begin{eqnarray}\label{defp} \xymatrix{(-)_{TOP^\s}\ar[r]^\pi&(-)_{TOP}\ar[r]^\pi&(-)_{top}}\end{eqnarray}
The pullback functors of the first $\pi$ are exact while those of the second one are only known to  preserve finite products (unless $top$ has fiber products). For a fixed $X^\s$, both $\smash{\pi_{X^\s}}$'s admit a right section $\smash{r_{X^\s}}$ (ie.
$\smash{r_{X^\s}^{-1}}=\smash{\pi_{X^\s,*}}$). It follows in particular that the functors $\smash{\pi_{X^\s}^{-1}}$ are fully faithful. The morphisms $r_{X^\s}$ are not pseudo-functorial with respect to $X^\s$ however.

If $top_1$ is finer than or equal to $top_2$ then we write $top_1\preceq top_2$ and we have a natural morphism
\begin{eqnarray}\label{defepsilon}\xymatrix{(-)_{TOP^\s_1}\ar[r]^-\epsilon & (-)_{TOP^\s_2}}\end{eqnarray}
The same holds for $TOP_i$ and $top_i$, except that in the latter case the functor $\epsilon^{-1}$ is only known to preserve finite products (unless $top_2$ has  fiber products). 
The weak morphisms $\epsilon$ and $\pi$ canonically pseudo-commute to each other.

\para \label{usualtopusual}
The properties used most frequently as $top$ in this paper are the following: \petit

- $fl$ means \emph{flat, locally of finite presentation},

- $syn$ means \emph{of complete intersection} (see Def. 19.3.1, Def. 19.3.6 and Prop. 19.3.7 of \cite{EGA4-IV}),

- $et$ means \emph{\'etale},

- $zar$ means \emph{open immersion}.

\petit \noindent Note that $X^\s_{FL}$ is thus equivalent to the usual \emph{fppf} topos of the underlying scheme as defined e.g. in \cite{SGA3-I} IV, Sect. 6.3.

We have
morphisms of variable topoi \begin{eqnarray}\label{usualepsilon}\xymatrix{(-)_{FL^\s}\ar[r]^\epsilon&
(-)_{SYN^\s}\ar[r]^\epsilon&(-)_{ET^\s}\ar[r]^\epsilon&(-)_{ZAR^\s}}\end{eqnarray} and similarly for big strict or small topoi (except that
$(-)_{fl}\to (-)_{syn}$ is only a weak morphism of weakly variable topoi).

\para \label{usualO} If $fl\preceq top$ then $U^\s\mapsto \G(U,\O_U)$ naturally defines a sheaf of rings on $TOP^\s(Spec(\Z))$ (here $\O_U$ denotes the structural ring of the scheme $U$).  We use the following notation:

\begin{defn} \label{defusualO} Let $X^\s$ in $\cal Sch^\s$. The ring of $X^\s_{TOP^\s}$, $X^\s_{TOP}$ or $X^\s_{top}$ obtained by restriction of the above sheaf of rings is called the \emph{structural ring} and is denoted $\O$.
\end{defn}

The weakly variable topoi defined in Sect. \ref{usualtop} are naturally ringed by $\O$ and the weak morphisms $\pi$, $\epsilon$ are naturally ringed as well.

%

\subsection{Crystalline sites for log schemes} \label{csfls}

\label{sectioncrys}


\para We review crystalline sites in the absolute case. Let $\Sigma_k$ denote $Spec(\Z/p^k)$ if $1\le k<\infty$ and $Spf(\Zp)$ if $k=\infty$. We denote $\cal Sch^\s/\Sigma_\infty$ or $\smash{\cal Sch^\s_{p,nil}}$ the full subcategory of $\cal Sch^\s$ formed by the fine log schemes on which $p$ is locally  nilpotent.  All divided powers and divided power envelopes considered are implicitly compatible with the divided power structure of $(\Zp,(p))$.

\begin{defn} \label{defcrystop} Let $X^\s$ in $\cal Sch^\s/\Sigma_1$, $\Sigma=\Sigma_k$ with $1\le k\le \infty$. Consider $top$ as in Sect. \ref{usualtop} and assume $fl \preceq top \preceq zar$.

\debrom

\item \label{defcrystopi} A \emph{dp-thickening} in $\cal Sch^\s/\Sigma$ is a quadruple $(U^\s,T^\s,\iota,\gamma)$ where $U^\s$, $T^\s$ are in $\cal Sch^\s/\Sigma$, $\iota:U^\s\to T^\s$ is an exact closed immersion and $\gamma$ is a divided power structure on the closed immersion $U\to T$ underlying $\iota$. We often use simplified notations such as $(U^\s,T^\s)$ or even $T^\s$. A morphism of dp-thickenings is a couple $f=(f_U,f_T)$ where $f_U$, $f_T$ are compatible morphisms of log schemes and $f_T$ is compatible with divided power structures.

\item \label{defcrystopii} We say that a morphism of dp-thickenings $f:(U'^\s,T'^\s)\to (U^\s,T^\s)$ is \emph{cartesian} if the underlying commutative square $$\xymatrix{U'^\s\ar[r]\ar[d]_{f_U}&T'^\s\ar[d]^{f_T}\\
    U^\s\ar[r]&T^\s}$$
    is cartesian in $\cal Sch^\s$. We say that $f$ is \emph{$top$ cartesian} if $f_T$ is moreover strict and  $top$.

\item \label{defcrystopiii} The $\s$-big crystalline pretopology $\smash{CRYS^\s_{top}(X^\s/\Sigma)}$ is defined as follows. An object $(U^\s/X^\s,T^\s,\gamma,\iota)$ of the underlying category  is a dp-thickening in $\cal Sch^\s/\Sigma$ together with a morphism $U^\s\to X^\s$. A morphism is a morphism $(f_U,f_T)$ of dp-thickenings such that $f_U$ is an $X^\s$-morphism.    A covering  is a surjective family of $top$ cartesian morphisms. The associated topos $\smash{(X^\s/\Sigma)_{CRYS^\s,top}}$ is called the $\s$-big crystalline $top$ topos.

\item \label{defcrystopiv} The big crystalline pretopology $\smash{CRYS_{top}(X^\s/\Sigma)}$  is the full subcategory of \break $\smash{CRYS^\s_{top}(X^\s/\Sigma)}$ formed by the $(U^\s/X^\s,T^\s,\gamma,\iota)$'s with $U^\s/X^\s$ strict together with surjective families of $top$ cartesian morphisms as coverings. The associated topos $\smash{(X^\s/\Sigma)_{CRYS,top}}$ is called the big crystalline $top$ topos.

\item \label{defcrystopv} The small crystalline pretopology  $\smash{crys_{top}(X^\s/\Sigma)}$ is the full subcategory of \break $\smash{CRYS^\s_{top}(X^\s/\Sigma)}$ formed by the $(U^\s/X^\s,T^\s,\gamma,\iota)$'s with $U^\s/X^\s$ strict $top$ together with  surjective families of $top$ cartesian morphisms as coverings. The associated topos $\smash{(X^\s/\Sigma)_{crys,top}}$ is called the small crystalline $top$ topos.
    \finrom
\end{defn}

As in the case of usual topoi, the inclusion functors induce weak morphisms \begin{eqnarray}\label{defpcrys}\xymatrix{(X^\s/\Sigma)_{CRYS^\s,top}\ar[r]^{\pi}&(X^\s/\Sigma)_{CRYS,top}\ar[r]^{\pi}&
(X^\s/\Sigma)_{crys,top}}\end{eqnarray} The pullback functor of the first (resp. second) $\smash{\pi}$ is exact (resp. preserves finite products). Both $\pi$'s admit a non functorial right section $\smash{r}$.

For $1\le k\le k'\le \infty$ the inclusion of $\smash{CRYS^\s_{top}}(X^\s/\Sigma_k)$ into  $\smash{CRYS^\s_{top}(X^\s/\Sigma_{k'})}$ is cocontinuous and
thus induces a morphism
\begin{eqnarray}\label{defiotacrys}\xymatrix{(X^\s/\Sigma_k)_{CRYS^\s,top}\ar[r]^{\iota_{k,k'}}
&(X^\s/\Sigma_{k'})_{CRYS^\s,top}}
\end{eqnarray}
We will also use the simplified notation $\iota_{k}:=\iota_{k,\infty}$. The same holds for $CRYS$ or $crys$.

If $top_1$ is finer than $top_2$ then  we have a natural morphism \begin{eqnarray}\label{defepsiloncrys}
\xymatrix{(X/\Sigma)_{CRYS^\s,top_1}\ar[r]^{\epsilon}&(X/\Sigma)_{CRYS^\s,top_2}}\end{eqnarray}
The same holds for $CRYS$ or $crys$ except that in the latter case $\epsilon^{-1}$ is only known to preserve finite products (unless $top_2$ has fiber products). 
The weak morphisms $\epsilon$ and $\pi$ canonically pseudo-commute to each other.


\para
The forgetful functor $CRYS^\s_{top}(X^\s/\Sigma)\to TOP^\s(X^\s)$, $(U^\s/X^\s,T^\s,\iota,\gamma)\mapsto U^\s/X^\s$ is a premorphism of pretopologies
and induces a morphism of topoi $i$ as follows:

\begin{eqnarray}\label{defiu} \xymatrix{X^\s_{TOP^\s}\ar[rr]_-i&&(X^\s/\Sigma)_{CRYS^\s,top}\ar @{-->} @/_2pc/ [ll]_-u}\end{eqnarray}

We do not know whether or not the forgetful functor is cocontinuous as well for an arbitrary $top$. When this is the case it also induces a morphism $u$ as
above which is a right retraction for $i$ (ie. $i_*=u^{-1}$, $ui\simeq id$). The situation is similar for big and small topoi.

\begin{lem} \label{lemiucocont}  The forgetful functor is cocontinuous (and thus induces a morphism $u$) in the cases $top=zar$, $et$ or $syn$.  \end{lem}
Proof. We only treat the case of $\s$-big sites. The other cases are similar. Consider the following property: \petit

$(lift)$ For every $(U^\s,T^\s)$ in $CRYS^\s_{top}(X^\s/\Sigma)$ and every strict $top$ morphism $V^\s\rightarrow U^\s$ there exists a surjective
family of $top$ morphisms $(V_i\rightarrow V)$ such that each $V_i/U$ admits a lifting to a $top$ morphisms $T_i/T$. \petit

Remark that if we set $T_i^\s:=T^\s\times_TT_i$ and $V_i^\s:=V^\s\times_VV_i$ and if we endow the closed immersions $\iota_i:V_i\to T_i$ with the
divided powers extending those of $V\to T$ (using  \cite{BO} Cor. 3.22 and flatness of  $T_i/T$) then the family $(\smash{(V_i^\s,T_i^\s)}\to (U^\s,T^\s))$ is a covering in the sense of Def.
\ref{defcrystop} \ref{defcrystopiii}. The property $(lift)$ would thus imply the desired cocontinuity:  any covering of $TOP(X^\s)$ can be refined by the image of
a covering in $\smash{CRYS_{top}^\s(X^\s/\Sigma)}$ under the forgetful functor.

Let us check $(lift)$ in each case. If $top$ is $zar$ then $V/U$ itself admits a lifting to $T$ namely the open subscheme of $T$ which has the same
underlying topological space as $V$. The case $top=et$ is similar using the fact that the categories of \'etale schemes over $U$ and $T$ are
naturally equivalent (\cite{SGA4-II} VIII, Thm. 1.1).
In the case $top=syn$, recall that syntomic $T$-schemes are characterized by the property that any (one is enough) closed embedding into a
smooth $T$-scheme is transversally regular relatively to $T$ (\cite{EGA4-IV} Prop. 19.3.7). Replacing $T$, $U$ and $V$ by open coverings we may assume that these
are affine schemes and that  $V$ is a closed subscheme into an open of $\mathbb A^{d}_{U}$ defined by a sequence $\underline
x=(x_1,\dots, x_n)$ which is transversally regular relatively to $U$. Then  any lift of the sequence $\underline x$ to the corresponding open of $\mathbb A^d_T$ is transversally regular relatively to $T$ (use \cite{EGA4-IV} Prop. 19.8.2 to replace $U$ and $T$ by a Noetherian affine schemes and then apply Prop. 2.5 and Rem. 2.6 (d) of \cite{Mi1} I) and thus gives rise
to the desired syntomic lift of $V/U$.
\begin{flushright}$\square$\end{flushright}


\begin{lem} \label{iotacrysex} Let $top=zar$, $et$ or $syn$ then the functor $\iota_{k,k',*}$ in (\ref{defiotacrys}) admits a right adjoint and is thus exact.
\end{lem}
Proof. It suffices to prove that the continuous functor $T^\s\mapsto T^\s_k$ inducing $\iota_{k,k'}$ is cocontinuous as well. The proof is similar to Lem. \ref{lemiucocont}.
\begin{flushright}$\square$\end{flushright}

\para \label{paralambdareal}
Keep the notations and assumptions of Def. \ref{defcrystop}. If  $T^\s=(U^\s,T^\s,\iota,\gamma)$ is an object of $CRYS^\s_{top}(X^\s/\Sigma)$, we have weak morphisms \begin{eqnarray}\label{deflambda}\xymatrix{(X^\s/\Sigma)_{CRYS^\s,top}&&
\ar[ll]_-{f_{T^\s}}(X^\s/\Sigma)_{CRYS^\s,top}/T^\s\ar[rr]^-{\lambda_{T^\s}} &&\smash{T^\s_{top}}}\end{eqnarray} The left one is sometimes denoted
$\smash{f_{T^\s/X^\s}}$ if a reference to $X^\s$ is useful. It is the usual localization morphism. Recall that the topos on the middle naturally identifies with the topos of sheaves on $\smash{CRYS^\s_{top}(X^\s/\Sigma)/T^\s}$ (induced topology). Modulo this identification, the right weak morphism is induced
by the premorphism of pretopologies $T'^\s/T^\s$ $\mapsto$ $(\smash{U^\s\times_{T^\s}T'^\s},T'^\s,\iota',\gamma')$ where $\iota'$ is the base change of $\iota$ by the flat morphism $T'^\s/T^\s$ and $\gamma'$ is the unique divided power structure on $\iota'$ extending $\gamma$ (\cite{BO}
Prop. 3.21). Its pullback is only known to preserve finite products (unless $top$ has fiber products).  
 Both $\smash{f_{T^\s}}$ and $\smash{\lambda_{T^\s}}$ are pseudo-functorial with respect to $T^\s$.

In the case $CRYS$ or $crys$, we have similarly a morphism $\smash{f_{T^\s}}$ and a weak morphism $\smash{\lambda_{T^\s}}$. Both are pseudo-compatible with $\pi$ in the obvious way.


\begin{defn} \label{defreal}For $T^\s$ as above and $F$ in $(X^\s/\Sigma)_{CRYS^\s,top}$ (resp.   $(X^\s/\Sigma)_{CRYS,top}$, resp. $(X^\s/\Sigma)_{crys,top}$) we use the following notations: \petit

- $F_{|T^\s}:=f_{T^\s}^{-1}F$ is the \emph{restriction of $F$ to $T^\s$}.

- $F_{T^\s}:=\lambda_{T^\s,*}F_{|T^\s}$ is the \emph{realization of $F$ on $T^\s$}.
\end{defn}

Some properties of these functors and their relation to the morphisms (\ref{defiu}) will be explained in Sect. \ref{cocs}.

\para \label{crysO} The rules $(U^\s,T^\s)\mapsto \G(T,\O_T)$ and $(U^\s,T^\s)\mapsto \G(U,\O_U)$ naturally define two sheaves of rings on $CRYS^\s(\Sigma_1/\Sigma_\infty)$.

\begin{defn} \label{defcrysO} Let $X^\s\in\cal Sch^\s$, $1\le k\le \infty$.
\debrom

\item \label{defcrysOi} The ring of $\smash{(X^\s/\Sigma_k)_{CRYS^\s,top}}$, $\smash{(X^\s/\Sigma_k)_{CRYS,top}}$ or $\smash{(X^\s/\Sigma_k)_{crys,top}}$  obtained by restriction of the first sheaf of rings above is called the \emph{structural ring} and is denoted $\O$ or  $\smash{\O_{X^\s/\Sigma_k}}$ depending on the context.

\item \label{defcrysOii} The ring of $\smash{(X^\s/\Sigma_k)_{CRYS^\s,top}}$, $\smash{(X^\s/\Sigma_k)_{CRYS,top}}$ or $\smash{(X^\s/\Sigma_k)_{crys,top}}$ obtained by restriction of the second sheaf of rings above is denoted $\Ga$.
\finrom
\end{defn}

The  topoi $\smash{((X^\s/\Sigma_k)_{CRYS^\s,top}}$, $\smash{((X^\s/\Sigma_k)_{CRYS,top}}$, $\smash{((X^\s/\Sigma_k)_{crys,top}}$ and the weak morphisms $\pi$, $\smash{\iota_{k,k'}}$, $\epsilon$, $i$,  $\smash{f_{T^\s}}$ and $\smash{\lambda_{T^\s}}$ are naturally ringed by $\O$ as well.

\begin{defn} \label{decOcrys} Assume that the morphism $u:\smash{(X^\s/\Sigma_k)_{CRYS^\s,top}\to X^\s_{TOP^\s}}$ exists. We set $$\O_k^{crys}:=u_*\O$$
We use similar notations in the setting of big or small topoi. In the case $k=\infty$ we also use the notation $\O^{crys}$ instead of $\smash{\O_\infty^{crys}}$.
\end{defn}



\subsection{Small crystalline sites, functoriality} \label{scsf} ~~ \\



We use the definitions and notations introduced in \ref{sectioncrys}. In particular $\Sigma=\Sigma_k$ with $1\le k\le\infty$ and $fl \preceq top \preceq zar$.

\para 
We explain pseudo-functoriality of the $\s$-big and small crystalline topoi $(X^\s/\Sigma)$ with respect to $X^\s$ in $\cal Sch^\s/\Sigma_1$. A similar discussion holds when either the $\s$-big or small topos is replaced by the big one but we omit it to
avoid lengthy statements.


\begin{defn} \label{deffunctcrys} Let $f:X'^\s\to X^\s$ be a morphism in $\cal Sch^\s/\Sigma_1$.
\debrom
\item \label{deffunctcrysi} We define  $f_{CRYS^\s}:(X'^\s/\Sigma)_{CRYS^\s,top}\to (X^\s/\Sigma)_{CRYS^\s,top}$ as the morphism induced by the forgetful cocontinuous functor $f_0:\smash{CRYS^\s_{top}(X'^\s/\Sigma)}\to \smash{CRYS^\s_{top}(X^\s/\Sigma)}$.

\item \label{deffunctcrysii}  We define a weak morphism  $f_{crys}:(X'^\s/\Sigma)_{crys,top}\to (X^\s/\Sigma)_{crys,top}$ by the formula  $$f_{crys}:=\pi f_{CRYS^\s}r$$
    where $\pi:\smash{(X'^\s/\Sigma)_{CRYS^\s,top}}\to \smash{(X'^\s/\Sigma)_{crys,top}}$ and $r:\smash{(X^\s/\Sigma)_{crys,top}}\to \smash{(X^\s/\Sigma)_{CRYS^\s,top}}$ are as in
    (\ref{defpcrys}).
    \finrom
\end{defn}

Later in the text, we will simply write $f$ instead of $\smash{f_{CRYS^\s}}$ or $f_{crys}$. For the following discussion however it is convenient to keep the indices.

It follows immediately from the definition that there are natural isomorphisms $$\smash{(ff')_{CRYS^\s}\simeq f_{CRYS^\s}f'_{CRYS^\s}}$$ satisfying the composition constraint  (look at cocontinuous functors). In other words, $\smash{(-/\Sigma)_{CRYS^\s,top}}$ is a variable topos on $\smash{\cal Sch^\s/\Sigma_1}$. The case of small crystalline topoi requires more work.

\begin{prop} Assume that $top$ is finer than or equal to $et$.  \label{propA1}

\debrom
\item \label{propA1i} There are natural isomorphisms $\smash{(fg)_{crys}}\rightarrow \smash{f_{crys}g_{crys}}$ satisfying the composition constraint. In other words $\smash{(-/\Sigma)_{crys,top}}$ is a weakly variable topos on $\cal Sch^\s/\Sigma_1$.
\item \label{propA1ii} The pullback functor of $f_{crys}$ preserves finite products and is exact if $top$ has fiber products. 
\finrom
\end{prop}
Proof. \ref{propA1i}  Since $r$ is a right section for $\pi$ there is a natural morphism $id\to r\pi$. Using this for each fine $\Sigma_1$-log scheme, we find a family of morphisms \begin{eqnarray}\label{morfcrys}(ff')_{crys}\to f_{crys}f'_{crys}\end{eqnarray} satisfying the composition constraint. We will now prove that  (\ref{morfcrys}) is an isomorphism in three steps. \petit

\emph{Step 1}. Consider  $f:X'^\s\to X^\s$ and $f':X''^\s\to X'^\s$ some morphisms in $\cal Sch^\s/\Sigma_1$. The morphism (\ref{morfcrys}) is invertible in the following cases:

 \debrom \item \label{claimA1i} $f'$ is strict $top$, or

 \item \label{claimA1ii} $f$ is strict \'etale.
\finrom
\petit

In case \ref{claimA1i}, a straightforward verification shows that $f'_{crys}$ is in fact induced by the cocontinuous functor $f'_0:\smash{crys_{top}(X''^\s/\Sigma_1)\to crys_{top}(X'^\s/\Sigma_1)}$ (compute inverse images) and we may conclude from the resulting isomorphism  $\smash{rf'_{crys}}\simeq \smash{f'_{CRYS^\s}r}$. In case \ref{claimA1ii}, we need the following fact which is a consequence of \cite{SGA4-II} VIII, Thm. 1.1 together with \cite{BO} Cor. 3.22. \petit

 \emph{Fact}. The category of \'etale cartesian dp-thickenings over a given dp-thickening $(U^\s,T^\s)$ is naturally equivalent to the category of strict \'etale log schemes over $U^\s$.
\petit

Using this fact,
we find that the cocontinuous functors underlying $\smash{f_{CRYS^\s}}$ and $f_{crys}$ have compatible right adjoints. We may then conclude from the resulting isomorphism $\smash{f_{CRYS^\s}\pi\simeq \pi f_{crys}}$.


\petit

\emph{Step 2}. We may always assume that the schemes $X$, $X'$ and $X''$ are affine and that $f':X''^\s\to X'^\s$ has a chart \begin{eqnarray}\label{chartstep2}\xymatrix{X''^\s\ar[r]\ar[d]^{c''}&X'^\s\ar[d]^{c'}\\
(Spec(\Z[P'']),P'')\ar[r]^{ch}&(Spec(\Z[P']),P')}\end{eqnarray} where $P'$ and $P''$ are finitely generated integral monoids.
\petit

Let us explain this. Given arbitrary $f$ and $f'$, we can always find a family of commutative diagrams $$\xymatrix{X''^\s_{\lambda}\ar[r]^-{f'_{\lambda}}\ar[d]^{h''_{\lambda}}&X'^\s_{\lambda}
\ar[r]^-{f_\lambda}\ar[d]^-{h'_\lambda}&X^\s_\lambda\ar[d]^-{h_\lambda}\\
X''^\s\ar[r]^-{f'}&X'^\s\ar[r]^-f&X^\s}$$
where $h_\lambda$, $h'_\lambda$ and $h''_\lambda$ are strict \'etale,
$X_\lambda$, $X'_\lambda$, $X''_\lambda$ and $f'_\lambda$ satisfy the assumptions of Step 2 and the family of the $h''_\lambda$'s is surjective. These diagrams induce squares $$\xymatrix{h''^{-1}_{\lambda,crys}f'^{-1}_{crys}f^{-1}_{crys}\ar[d]^\wr\ar[r]&h''^{-1}_{\lambda,crys}(ff')^{-1}_{crys}\ar[d]^\wr\\
f'^{-1}_{\lambda,crys}f^{-1}_{\lambda,crys}h^{-1}_{\lambda,crys}\ar[r]&(f_\lambda f'_\lambda)_{crys}h^{-1}_{\lambda,crys}}$$
where the vertical maps are isomorphisms by Step 1 and which are commutative thanks to the composition constraint. It remains to notice that the family of the $\smash{h''^{-1}_{\lambda,crys}}$'s is conservative (indeed the essential images of the underlying cocontinuous functors are generating thanks to the Fact used in the proof of Step 1).
\petit

\emph{Step 3}. The morphism (\ref{morfcrys}) is invertible under the assumptions of Step 2.  \petit

The reader may easily establish the formula \begin{eqnarray}\label{formulafcrys}\smash{f_{crys,*}F(U^\s,T^\s)\simeq \limp_{\cal C_{f,(U^\s,T^\s)}}
F(U'^\s,T'^\s)}\end{eqnarray}  where the projective limit is indexed on the category $\cal C_{f,(U^\s,T^\s)}$ of couples $(\smash{(U'^\s/X'^\s,T'^\s)},\break \tilde f:\smash{f_0(U'^\s/X'^\s,T'^\s)\rightarrow (U^\s/X^\s,T^\s)})$ where $\smash{(U'^\s/X'^\s,T'^\s)}$ is an object of  $crys_{top}(\smash{X'^\s/\Sigma})$, $\smash{f_0(U'^\s/X'^\s,T'^\s)}$ is the object of $\smash{CRYS^\s_{top}(X^\s/\Sigma)}$ obtained by composition of $\smash{U'^\s\to X'^\s}$ with $\smash{f:X'^\s\to X^\s}$ and $\tilde f$ is a morphism in $\smash{CRYS^\s_{top}(X^\s/\Sigma)}$.  The projective limit remains the same if the category $\cal C_{f,(U^\s,T^\s)}$ is replaced by the full subcategory $\smash{\cal C^*_{f,(U^\s,T^\s)}}$ defined by the conditions that the scheme underlying $T'^\s$ is affine and that the composed morphism $\smash{U'^\s\to X'^\s\mathop\to\limits^{c'} Spec(\Z[P'],P')}$ extends to $T'^\s$ (this happens \'etale locally on $T'$ thanks to \cite{Ka2} Lem. 2.10). Let us denote $h=ff'$. The reader may check that the evaluation of the direct image functors underlying (\ref{morfcrys}) at a given $F$ and $(U^\s,T^\s)$ is \begin{eqnarray}\label{morfcrysUT}\limp_{\cal C_{h,(U^\s,T^\s)}^*} F(U''^\s,T''^\s)\to  \limp_{\cal C_{f,f',(U^\s,T^\s)}^*}F(U''^\s,T''^\s)\end{eqnarray} where $\smash{\cal C_{h,(U^\s,T^\s)}^*}$ and ${\cal C_{f',(U'^\s,T'^\s)}^*}$ are defined in the same way as  $\smash{\cal C^*_{f,(U^\s,T^\s)}}$ while ${\cal C^*_{f,f',(U^\s,T^\s)}}$ is the category where an object is a couple $(\smash{((U''^\s/X''^\s,T''^\s),\tilde f')},\smash{((U'^\s/X'^\s,T'^\s),\tilde f)})$ whose first (resp. second) argument is in $\smash{\cal C^*_{f',(U'^\s,T'^\s)}}$ (resp. $\smash{\cal C^*_{f,(U^\s,T^\s)}}$) and  where morphisms are defined in the natural way.

In order to prove that (\ref{morfcrys}) is an isomorphism, it is sufficient to prove that this is the case when $T$ is affine. In that case we will prove that
the natural ``composition'' functor \begin{eqnarray}\label{compoCf}\smash{\cal C^*_{f,f',(U^\s,T^\s)}}\to \smash{\cal C^*_{h,(U^\s,T^\s)}} \end{eqnarray} is cofinal, ie. that the category $\cal C:=(\smash{(U''^\s/X''^\s,T''^\s),\tilde h)}/\smash{\cal C^*_{f,f',(U^\s,T^\s)}}$ is connected for any $(\smash{(U''^\s/X''^\s,T''^\s),\tilde h)}$ in $\smash{\cal C^*_{h,(U^\s,T^\s)}}$. \petit

Let us first prove that $\cal C$ is non empty, ie. that there exist $\smash{((U'^\s_1/X'^\s,T_1'^\s),\tilde f_1)}$ in  $\smash{\cal C^*_{f,(U^\s,T^\s)}}$, $\smash{((U''^\s_1/X''^\s,T_1''^\s),\tilde f'_1)}$ in  $\smash{\cal C^*_{f',(U'^\s_1,T'^\s_1)}}$ and a morphism $e_1$ in $crys_{top}(X''^\s/\Sigma)$ rendering the following diagram commutative in  $\smash{CRYS^\s_{top}(X^\s/\Sigma)}$: \begin{eqnarray}\label{diaghCff}
\xymatrix{h_0(U''^\s/X''^\s,T''^\s)\ar[dd]^-{\tilde h}\ar[r]^-{h_0(e_1)}&f_0f'_0(U''^\s_1/X''^\s,T''^\s_1)\ar[d]^-{f_0(\tilde f'_1)}\\ &
f_0(U'^\s_1/X'^\s,T'^\s_1)\ar[dl]^-{\tilde f_1}\\
(U^\s/X^\s,T^\s)}
 \end{eqnarray}

\noindent
Consider the affine schemes $\smash{U':=X'\times_XU}$ and choose polynomial schemes $\smash{Y''=Spec(\Zp[\underline x''])}$, $\smash{Y'=Spec(\Zp[\underline x'])}$ as well as closed immersions $\smash{\iota''_1}$, $\smash{\iota'_1}$ and morphisms $y''$, $y'$ fitting into a commutative diagram of the form $$\xymatrix{&U''\ar[r]\ar[dl]\ar[d]^{\iota''_1}&U'\ar[d]^{\iota'_1}\\ T''\ar[r]^-{y''}&Y''\ar[r]^-{y'}&Y'}$$
Putting this together with a chart for $f'$ as in (\ref{chartstep2}) together with the choice of a lifting to $\smash{T''^\s}$ of the morphism $\smash{U''^\s\to X''^\s\mathop\to\limits^{c''}Spec(\Z[P''],P'')}$ easily produces a commutative diagram
of fine log schemes $$\xymatrix{&X''^\s\ar[r]&X'^\s\ar[r]&X^\s\\
&U''^\s\ar[r]\ar[d]\ar[dl]\ar[u]&U'^\s\ar[r]\ar[d]\ar[u]&U^\s\ar[d]\ar[u]\\ T''^\s\ar[r]&Y''\times Spec(\Z[P''],P'')\times T^\s\ar[r]& Y'\times Spec(\Z[P'],P')\times T^\s\ar[r]& T^\s}$$
The desired object of $\cal C$ is obtained by forming logarithmic divided power envelopes of the bottom vertical closed immersions with respect to the divided power of $(U^\s,T^\s)$.

It remains to prove that two objects in $\cal C$ are always related by a chain of arrows. Let us simply denote $\smash{((\tilde f_1,\tilde f_1'),e_1)}$ the object of $\cal C$ just constructed and consider an other arbitrary object $\smash{((\tilde f_2,\tilde f_2'),e_2)}$. We will now hint the construction of a chain $$\smash{((\tilde f_1,\tilde f_1'),e_1)}\leftarrow \smash{((\tilde f_3,\tilde f_3'),e_3)} \rightarrow \smash{((\tilde f_4,\tilde f_4'),e_4)}\leftarrow \smash{((\tilde f_1,\tilde f_1'),e_1)}$$ and the proof will be finished. The reader may easily guess how to define a commutative diagram of closed immersions 
{\small$$\hspace{-.5cm}\xymatrix{ &(U''^\s,Y''\times Spec(\Z[P''],P'')\times T^\s)\ar[r]& (U'^\s,Y'\times Spec(\Z[P'],P')\times T^\s)\ar[rd]&  \\  (U''^\s,T''^\s)\ar[ru]\ar[rd]\ar[r]&(U''^\s,Y''\times Spec(\Z[P''],P'')\times T''^\s_2)\ar[r]\ar[u]\ar[d]& (U'^\s_2, Y'\times Spec(\Z[P'],P')\times T'^\s_2)\ar[r]\ar[u]\ar[d]& (U^\s,T^\s) \\ &(U''^\s_2,T''^\s_2)\ar[r]& (U'^\s_2,T'^\s_2)\ar[ru]&}$$}

\noindent whose first terms are above $$\xymatrix{X''^\s\ar@{=}[rrr]&&&X''^\s\ar[rrrrr]&&&&&X'^\s\ar[rrr]&&&X^\s}$$
and whose bottom (resp. top) line is obtained from $\smash{((\tilde f_2,\tilde f_2'),e_2)}$ by forgetting divided powers (resp. was used in the definition of $\smash{((\tilde f_1,\tilde f_1'),e_1)}$). We define $\smash{((\tilde f_3,\tilde f_3'),e_3)}$ (resp.  $\smash{((\tilde f_4,\tilde f_4'),e_4)}$) as the object of $\cal C$ obtained from the second line
by forming logarithmic divided powers envelopes of the middle (resp. bottom) line with respect to the divided powers of $(U^\s,T^\s)$. This gives a chain
$\smash{((\tilde f_1,\tilde f_1'),e_1)}\leftarrow \smash{((\tilde f_3,\tilde f_3'),e_3)} \rightarrow \smash{((\tilde f_4,\tilde f_4'),e_4)}$. We conclude by noticing that the universal property of logarithmic divided power envelopes with respect to $(U^\s,T^\s)$ gives a morphism $\smash{((\tilde f_2,\tilde f_2'),e_2)}\to \smash{((\tilde f_4,\tilde f_4'),e_4)}$.

\ref{propA1ii} The first statement follows from the fact that the product of two elements is representable in the small crystalline site (hence $p^{-1}$ preserves finite products). The second statement follows from a straightforward adaptation of  \cite{Be1} Chap. 3,
Cor. 2.1.4, which shows that the small crystalline site has fiber products as soon as $top$ does.
\begin{flushright}$\square$\end{flushright}



\para Let us now investigate the behaviour of some usual (weak) morphisms with respect to crystalline functoriality.

\begin{prop} \label{compatcrys1} Let $fl\preceq top'\preceq top\preceq et$.

\debrom \item \label{compatcrys1i} The weak morphisms (\ref{defpcrys}), (\ref{defepsiloncrys}) and (\ref{defiu}) naturally define weak morphisms of weakly variable topoi on $\cal Sch^\s/\Sigma_1$ fitting into a canonically pseudo-commutative cube $$\xymatrix{&(-)_{top'}\ar[rr]^(.25)i\ar'[d]^(.5)\epsilon[dd]&&(-/\Sigma)_{crys,top'}\ar[dd]^(.25)\epsilon\\
(-)_{TOP'^\s}\ar[dd]^(.25)\epsilon\ar[rr]^(.25)i\ar[ru]^(.35)\pi&&(-/\Sigma)_{CRYS^\s,top'}\ar[ru]^(.35)\pi
\ar[dd]^(.25)\epsilon&\\
&(-)_{top}\ar'[r]^(.5)i[rr]&&(-/\Sigma)_{crys,top}\\
(-)_{TOP^\s}\ar[rr]^(.25)i\ar[ru]^(.35)\pi&&(-/\Sigma)_{CRYS^\s,top}\ar[ru]^(.35)\pi&}$$

\item \label{compatcrys1ii} If $top$ and $top'$ satisfy the property $(lift)$ considered in Lem. \ref{lemiucocont}, then (\ref{defiu}) naturally defines a morphism of weakly variable topoi as well and there is a canonically pseudo-commutative cube
$$\xymatrix{&(-)_{top'}\ar'[d]^(.5)\epsilon[dd]&&(-/\Sigma)_{crys,top'}\ar[ll]_(.25)u\ar[dd]^(.25)\epsilon\\
(-)_{TOP'^\s}\ar[dd]^(.25)\epsilon\ar[ru]^(.35)\pi&&(-/\Sigma)_{CRYS^\s,top'}\ar[ll]_(.35)u\ar[ru]^(.35)\pi
\ar[dd]^(.25)\epsilon&\\
&(-)_{top}&&(-/\Sigma)_{crys,top}\ar'[l]_(.5)u[ll]\\
(-)_{TOP^\s}\ar[ru]^(.35)\pi&&(-/\Sigma)_{CRYS^\s,top}\ar[ll]_(.35)u\ar[ru]^(.35)\pi&}$$

The isomorphism $ui\simeq id$ is moreover functorial and compatible with $\pi$ and $\epsilon$.
\finrom
\end{prop}
Proof. \ref{compatcrys1i} A similar pseudo-commutative cube above a fixed $X^\s$ in $\cal Sch^\s/\Sigma_1$ is obvious by looking at continuous functors. Let $\smash{(v_{X^\s})_{X^\s\in \cal Sch^\s/\Sigma_1}}$ denote one of the twelve collections of weak morphisms involved. We need to achieve the following tasks: \petit

\noindent (A) Enrich $\smash{(v_{X^\s})}$ into a pseudo-morphism $v$, ie. define a functoriality isomorphism $fv_{X^\s}\simeq v_{X'^\s}f$ for each $f:X'^\s\to X^\s$, this family of isomorphisms being submitted to the composition constraint. \petit

\noindent (B) Check that the isomorphisms expressing the pseudo-commutativity of the cube above a fixed $X^\s$ are functorial with respect to $X^\s$, ie. compatible with the functoriality isomorphisms of (A).
\petit

We first explain how to achieve task (A).
Let us begin with $v=\pi$. The case of $\smash{\pi_{X^\s}^{-1}}:\smash{X^\s_{top}}\to \smash{X^\s_{TOP^\s}}$ is obvious by looking at continuous functors. Let us thus consider the case of $\smash{\pi_{X^\s}^{-1}}:\smash{(X^\s/\Sigma)_{crys,top}}\to \smash{(X^\s/\Sigma)_{CRYS^\s,top}}$. The morphisms $\smash{\pi_{X'^\s}^{-1}r_{X'^\s}^{-1}}\to id$ for $f:X'^\s\to X^\s$ give a family of morphisms
\begin{eqnarray}\label{plax}\smash{\pi_{X'^\s}^{-1}f_{crys}^{-1}\to f_{CRYS^\s}^{-1}\pi_{X^\s}^{-1}}\end{eqnarray}
satisfying the composition constraint. If now $f'':X''^\s\to X'^\s$ is a morphism in $\cal Sch^\s/\Sigma_1$, then (\ref{plax}) becomes an isomorphism after applying $r_{X''^\s}^{-1}f'_{CRYS^\s}$ (Prop. \ref{propA1} \ref{propA1i}). It follows that (\ref{plax}) is in fact an isomorphism and task (A) is achieved in the case $v=\pi$.



Next we consider the case $v=\epsilon$. Here again the case of usual topoi is obvious by looking at continuous functors. The case of big crystalline topoi is clear as well (sheafification pseudo-commutes to restriction). The case of small crystalline topoi can be deduced as follows. We define an isomorphism $\epsilon_{X^\s,*}f_{crys,*}\simeq f_{crys,*} \epsilon_{X'^\s,*}$ by composition: \begin{eqnarray}\label{ef1}f_{crys,*}\epsilon_{X'^\s,*}&\simeq \label{ef2}&f_{crys,*}\epsilon_{X'^\s,*}\pi_{X'^\s,*}\pi_{X'^\s}^{-1}\\
\label{ef3}&\simeq &f_{crys,*}\pi_{X'^\s,*}\epsilon_{X'^\s,*}\pi_{X'^\s}^{-1}\\
\label{ef4}&\simeq &\pi_{X^\s,*}f_{CRYS^\s,*}\epsilon_{X'^\s,*}\pi_{X'^\s}^{-1}\\
\label{ef5}&\simeq &\pi_{X^\s,*}\epsilon_{X^\s,*}f_{CRYS^\s,*}\pi_{X'^\s}^{-1}\\
\label{ef6}&\simeq &\epsilon_{X^\s,*}\pi_{X^\s,*}f_{CRYS^\s,*}\pi_{X'^\s}^{-1}\\
\label{ef7}&\simeq &\epsilon_{X^\s,*}f_{crys,*}\pi_{X'^\s,*}\pi_{X'^\s}^{-1}\\
\label{ef8}&\simeq &\epsilon_{X^\s,*}f_{crys,*}\end{eqnarray}
The composition constraint for this family of isomorphism is a formal consequence (exercise!) of the previously established composition constraint for the following weak morphisms of variable topoi:
$$\begin{array}{l}\pi:(-/\Sigma)_{CRYS^\s,top'}\to (-/\Sigma)_{crys,top'}\\
\pi:(-/\Sigma)_{CRYS^\s,top}\to (-/\Sigma)_{crys,top}\\
\epsilon: (-/\Sigma)_{CRYS^\s,top'}\to (-/\Sigma)_{CRYS^\s,top}\end{array}$$

Let us finally consider the case $v=i$. The case of big topoi is obvious by looking at cocontinuous functors. The case of small topoi follows formally, as in the case $v=\epsilon$.

Let us now explain task (B). The left hand face is clear by looking at continuous functors. The right, top and bottom faces are tautological from the definition of functoriality isomorphisms for $\epsilon$ and $i$ in the setting of small crystalline topoi. Using that $\pi$ admits a right retraction, we see that it only remains to consider the front face ie. to check that the following diagram commutes: $$\diagram{f_{CRYS^\s}i_{X'^\s}\epsilon_{X'^\s}&\simeq &i_{X^\s}f_{CRYS^\s}\epsilon_{X'^\s}&\simeq & i_{X^\s}\epsilon_{X^\s}f_{CRYS^\s} \cr \wr | &&&&\wr| \cr
f_{CRYS^\s}\epsilon_{X'^\s}i_{X'^\s}&\simeq &\epsilon_{X^\s}f_{CRYS^\s}i_{X'^\s}&\simeq &\epsilon_{X^\s}i_{X^\s}f_{CRYS^\s}}$$
This causes no difficulty when looking at inverse image functors.

\ref{compatcrys1ii} Task (A) for $v=u$ in the setting of big topoi is obvious by looking at cocontinuous functors. The case of small topoi follows formally, as in the case $v=\epsilon$. Regarding task (B), it is sufficient to consider the front face. By cocontinuity of the functor $\smash{CRYS^\s_{top'}(X^\s/\Sigma)}\to TOP'^\s(X^\s)$,  $(U^\s,T^\s)\to U^\s$ we find that the $top'$ sheafification of $(U^\s,T^\s)\mapsto F(U)$ and $U^\s\mapsto F(U^\s)$ coincide. This defines an isomorphism $$\epsilon^{-1}_{X^\s}u_{X^\s}^{-1}\simeq u_{X^\s}^{-1}\epsilon_{X^\s}^{-1}$$ whose compatibility with $f_{CRYS^\s}^{-1}$ is immediate.

The functoriality of the isomorphism $\smash{id\simeq u_{X^\s}i_{X^\s}}$ is clear by looking at cocontinuous functors. Compatibility with $\pi_{X^\s}$ and $\epsilon_{X^\s}$ causes no difficulty either looking at direct or inverse images.

\begin{flushright}$\square$\end{flushright}

\subsection{Crystals on crystalline sites} \label{cocs} ~~ \\

We review the definition and basic properties of crystals.

\para Let us begin with some properties of the realization functors $\smash{(-)_T^\s}$ defined in Def. \ref{defreal}.

\begin{lem} \label{lemreal} \debrom
\item  \label{lemreali}
Consider $(U^\s,T^\s)$ in $\smash{CRYS^\s_{top}(X/\Sigma)}$. The functor $top(T^\s)\to$ \break $\smash{CRYS^\s_{TOP}(X/\Sigma)}/(U^\s,T^\s)$,
$T'^\s/T^\s\mapsto \smash{(U^\s\times_{T^\s}T'^\s,T'^\s)/(U^\s,T^\s)}$  is fully faithful continuous and cocontinuous. The topology of $top(T^\s)$ is induced by the topology of $\smash{CRYS^\s_{top}(X/\Sigma)/(U^\s,T^\s)}$ via this functor.

\item \label{lemrealii} The realization functor $\smash{(-)_{T^\s}}:\smash{(X^\s/\Sigma)_{CRYS^\s,top}\to T^\s_{top}}$ is the inverse image functor of a morphism of topoi. The corresponding functor between categories of $\O$-modules commutes in particular to arbitrary limits and commutes to tensor products.

\item \label{lemrealiii} Sending a sheaf $F$ on $\smash{CRYS^\s_{top}(X/\Sigma)}$ to $\xi=((U^\s,T^\s)\mapsto \smash{F_{T^\s}},(f:(U'^\s,T'^\s)\to (U^\s,T^\s))\mapsto (f^{-1}\smash{F_{T^\s}}\to \smash{F_{T'^\s}}))$ induces a natural equivalence  $$\smash{(X^\s/\Sigma)_{CRYS^\s,top}}\simeq\smash{Real(CRYS^\s_{top}(X/\Sigma))}$$ where the right category denotes the full subcategory of $\smash{\Gamma(CRYS^\s_{top}(X/\Sigma),(-)_{top})}$ whose objects are the sections $\xi:(U^\s,T^\s)\mapsto \xi(U^\s,T^\s)$, $(f:(U'^\s,T'^\s)\to (U^\s,T^\s))$ $\mapsto (\xi(f):f^{-1} \xi(U^\s,T^\s)\to \xi(U'^\s,T'^\s))$ satisfying the condition that $\xi(f)$ is an isomorphism if $f$ is $top$ cartesian.

\item \label{lemrealiv} Consider $top'$ finer than $top$. A sheaf $F$ on $\smash{CRYS^\s_{top}(X/\Sigma)}$ is a sheaf for $top'$ if and only if the corresponding $\xi$ satisfies the following descent condition: if $\smash{f_i:(U_i^\s,T_i^\s)\to (U^\s,T^\s)}$ is a covering in $\smash{CRYS^\s_{top}(X/\Sigma)}$ then $$\xi(U^\s,T^\s)\simeq Ker (\prod_if_{i,*}\xi(U_i^\s,T_i^\s)\rightrightarrows \prod_{j,k}f_{jk,*}\xi(U_{jk}^\s,T_{jk}^\s))$$
    where $f_{jk}:(U_{jk}^\s,T_{jk}^\s)\to (U^\s,T^\s)$ is the product of $f_j$ and $f_k$ computed in the category $\smash{CRYS^\s(X^\s/\Sigma)/(U^\s,T^\s)}$. The full subcategory of $\smash{Real(CRYS^\s_{top}(X^\s/\Sigma))}$ formed by such sections will be denoted $\smash{Real_{top'}(CRYS^\s_{top}(X^\s/\Sigma))}$.
\finrom
\end{lem}
Proof. Everything is straightforward from Def. \ref{defcrystop}. \begin{flushright}$\square$\end{flushright}

Let us emphasize that the morphism of topoi mentioned in \ref{lemrealii} is not pseudo-functorial with respect to $(U^\s,T^\s)$.

\begin{rem} \label{remreal} Here are some complementary remarks.

\debrom
\item \label{remreali} Let $T^\s=(U^\s,T^\s)$ in $\smash{CRYS^\s_{top}(X^\s/\Sigma)}$. An obvious local variant of  Lem. \ref{lemreal} \ref{lemrealiii}, \ref{lemrealiv} provides equivalences \begin{eqnarray}\label{remrealeq1}(X^\s/\Sigma)_{CRYS^\s_{top}}/T^\s&\simeq & Real(CRYS^\s_{top}(X^\s/\Sigma)/T^\s)\\
    (X^\s/\Sigma)_{CRYS^\s_{top'}}/T^\s&\simeq & Real_{top'}(CRYS^\s_{top}(X^\s/\Sigma)/T^\s)\end{eqnarray}

\item  \label{remrealii} Consider $T^\s$, $T'^\s$, a morphism $h:T'^\s\to T^\s$ in $\smash{CRYS^\s_{top}(X^\s/\Sigma)}$, $f_{T^\s}$, the structural morphism of $T^\s$ viewed as an object of the topos and $f:X'^\s\to X^\s$, a morphism in $\cal Sch^\s/\Sigma_1$. The functors $$\begin{array}{l}\lambda_{T^\s}^{-1}:T^\s_{top}\to (X^\s/\Sigma)_{CRYS^\s_{top}}/T^\s\\ f_{T^\s}^{-1}:(X^\s/\Sigma)_{CRYS^\s_{top}}\to (X^\s/\Sigma)_{CRYS^\s_{top}}/T^\s\\
    h^{-1}:(X^\s/\Sigma)_{CRYS^\s_{top}}/T^\s\to  (X^\s/\Sigma)_{CRYS^\s_{top}}/T'^\s\\
    f^{-1}:(X^\s/\Sigma)_{CRYS^\s_{top}}\to (X'^\s/\Sigma)_{CRYS^\s_{top}}\end{array}$$ have the following convenient translation in terms of the corresponding categories of realizations: $(\smash{\lambda_{T^\s}^{-1}}F):(\smash{h_1:T_1^\s\to T^\s})\mapsto \smash{h^{-1}_1F}$, $f_T^{-1}\xi:(T_1^\s/T^\s)\mapsto \xi(T_1^\s)$, $h^{-1}\xi:\smash{(T_1^\s/T'^\s)}\mapsto \xi(T_1^\s/T^\s)$, $f^{-1}\xi:T_1^\s\mapsto \xi(T_1^\s)$.

\item \label{remrealiii} In Lem. \ref{lemreal} \ref{lemrealiv}, the descent condition needs only to be checked for families $f_i:T_i^\s\to T^\s$ where the morphisms of schemes underlying the $f_i$'s are affine (see \cite{SGA4-I} II, Prop. 2.3 and \cite{EGA2} Prop. 1.5.1).

\item \label{remrealiv} Statements Lem. \ref{lemreal} \ref{lemrealiii}, \ref{lemrealiv}, Rem. \ref{remreal}  \ref{remreali}, \ref{remrealii},\ref{remrealiii}  have obvious counterparts for the categories of modules.

\item \label{remrealv} Statements Lem. \ref{lemreal} \ref{lemreali}, \ref{lemrealii}, \ref{lemrealiii} and the  Rem. \ref{remreal} \ref{remreali} (\ref{remrealeq1}), \ref{remrealiii}, as well as their counterparts for modules hold verbatim if $CRYS^\s$ is replaced by $crys$. The same is true for Rem. \ref{remreal} \ref{remrealii} if one assumes that $f$ is strict and $top$.

\finrom
\end{rem}

Let us write down some compatibilities which will be useful later on.

\begin{lem} \label{lemreal2} Consider $X^\s$ in $\cal Sch^\s/\Sigma_1$ and let $T^\s=(U^\s/X^\s,T^\s,\iota,\gamma)$ in $\smash{CRYS^\s_{top}(X^\s/\Sigma)}$ (resp. $\smash{crys_{top}(X^\s/\Sigma)}$). Let $\smash{f_{T^\s}}$ and $\smash{\lambda_{T^\s}}$ as in (\ref{deflambda}) and let $\pi$ denote the weak morphism of projection from $\s$-big to small crystalline or usual $top$ topoi.

\debrom  \item \label{lemreal2i} There are canonical isomorphisms 
    $$\begin{array}{rrl}&(i_*F)_{(U^\s,T^\s)}\simeq \pi_*\iota_*(F_{|U^\s}) & \hbox{(resp.} (i_*F)_{(U^\s,T^\s)}\simeq \iota_*F_{|U^\s}\hbox{ )} \\ \label{iotalambdai2} \hbox{and} &   \pi_*((i^{-1}G)_{|U^\s}) \simeq G_{(U^\s,U^\s)}  &\hbox{(resp.} (i^{-1}G)_{|U^\s}\simeq G_{(U^\s,U^\s)}  \hbox{ )}\end{array}$$
which are functorial with respect to $F$ in the $\s$-big (resp. small) $top$ topos and $G$ in the $\s$-big (resp. small) crystalline $top$ topos.

\item \label{lemreal2ii} If $U^\s=X^\s$ then we have a canonical morphism
$$\begin{array}{rl}\pi\iota uf_{T^\s}\rightarrow \lambda_{T^\s} & \hbox{(resp. } \iota uf_{T^\s}\rightarrow \lambda_{T^\s} \hbox{ )}\end{array}$$
This is in fact an isomorphism if $top=zar$ or $et$.


 \finrom
\end{lem}
Proof. The isomorphisms of \ref{lemreal2i} are immediate from the definitions. Let us prove \ref{lemreal2ii} in the case $CRYS^\s$ using the local realization functors $$(-)_{T'^\s/T^\s}:(X^\s/\Sigma)_{CRYS^\s_{top}}/T^\s\to T'^\s_{top}$$ underlying the equivalence (\ref{remrealeq1}). Given $F$ in $T^\s_{top}$ and $h:(U'^\s,T'^\s)\to (U^\s,T^\s)$ in \break $CRYS^\s_{top}(X^\s/\Sigma)/T^\s$, we have (using  \ref{lemreal2i}, Rem. \ref{remreal} and full faithfulness of $\pi^{-1}$)  $$\begin{array}{rrcl}& (\lambda_{T^\s}^{-1}F)_{T'^\s/T^\s}&\simeq &h^{-1}F\\
\hbox{and}&(f_{T^\s}^{-1}u^{-1}\iota^{-1}\pi^{-1}F)_{T'^\s/T^\s}&\simeq &(i_*\iota^{-1}\pi^{-1}F)_{T'^\s}\\
&&\simeq & \pi_*\iota_*(\iota^{-1}\pi^{-1}F)_{|U'^\s}\\ &&\simeq &\pi_*\iota_*\iota^{-1}\pi^{-1}h^{-1}F\\
&&\simeq &\iota_*\iota^{-1}h^{-1}F
\end{array}$$
The adjunction morphism $id\to \iota_*\iota^{-1}$ for varying $T'^\s/T^\s$'s and $F$'s gives rise to the desired morphism $\smash{\pi\iota uf_{T^\s}\to \lambda_{T^\s}}$. If now $top=et$ (resp. $zar$) then $\iota:\smash{U'^\s_{top}\to T'^\s_{top}}$ is an equivalence (resp. an isomorphism) and this morphism is thus invertible as claimed. The case of small crystalline topoi is similar.

%

\begin{flushright}$\square$\end{flushright}

\para \label{crystalconditions} Let us come to the definition of crystals.

\begin{lem} \label{lemcrystal1} (The crystal condition) Consider $M$ in  $\smash{Mod((X^\s/\Sigma)_{CRYS^\s,top},\O)}$ and  $T^\s$ in $\smash{CRYS^\s_{top}}(X^\s/\Sigma)$.

\begin{enumerate}

\item \label{lemcrystal1-tae} The following conditions are equivalent.  \debrom

\item \label{lemcrystal1i} The adjunction morphism $\smash{\lambda_{T^\s}^*M_{T^\s}}\rightarrow \smash{M_{|T^\s}}$ is invertible.

\item \label{lemcrystal1ii} The base change morphism $\smash{h^*M_{T'^\s}\rightarrow M_{T^\s}}$ is invertible for all  $h:T'^\s\to T^\s$ in $\smash{CRYS^\s_{top}}(X^\s/\Sigma)/T^\s$.
\finrom

\item \label{lemcrystal1iii} Let $\smash{(T^\s_i)}$ denote a family of objects in $\smash{CRYS^\s_{top}(X^\s/\Sigma)}$ which covers the final object of the topos. If the conditions of \ref{lemcrystal1-tae} hold with $T^\s=\smash{T_i^\s}$ for each $i$ then it holds for any $T^\s$.
\end{enumerate}

The same holds verbatim with $crys$ instead of $CRYS^\s$.
\end{lem}
Proof.  \ref{lemcrystal1-tae} This results from the description of $\smash{\lambda_{T^\s}^*}$ in terms of realizations (Rem. \ref{remreal} \ref{remrealii} and \ref{remrealiv}). 

\ref{lemcrystal1iii} if $T^\s$ is arbitrary, the conditions of  \ref{lemcrystal1-tae} for $T^\s$ are $top$ local on $T^\s$ and it is thus sufficient to check them under the assumption that $\smash{Hom(T^\s,T_i^\s)}\ne \emptyset$. In that case the verification is straightforward. The same remark holds with $crys$ instead of $CRYS^\s$.
\begin{flushright}$\square$\end{flushright}

\begin{defn} \label{defcrystal} We say that $M$ is a \emph{crystal} of $\smash{((X^\s/\Sigma)_{CRYS^\s,top},\O)}$ if the equivalent conditions of Lem. \ref{lemcrystal1} are satisfied for all $T^\s$'s. The full subcategory of $\smash{Mod((X^\s/\Sigma)_{CRYS^\s,top},\O)}$ formed by crystals is denoted $\smash{\cal Crys ((X^\s/\Sigma)_{CRYS^\s,top},\O)}$. We use similar definitions and notations in the context of small crystalline topoi.
\end{defn}

The category of crystals enjoys the following useful properties.

\begin{lem} \label{lemcrystal2} Assume that $top$ is finer than or equal to $et$.
 \debrom \item \label{lemcrystal2i} The adjunction $(\pi_*,\pi^*)$ for $\O$-modules on crystalline sites induces an equivalence $$\xymatrix{\smash{\cal Crys ((X^\s/\Sigma)_{CRYS^\s,top},\O)}\ar[rr]^-{(\pi^*,\pi_*)}&&\smash{\cal Crys ((X^\s/\Sigma)_{crys,top},\O)}}$$
\item \label{lemcrystal2ii}  Consider a morphism $g:X'^\s\to X^\s$ in $\cal Sch/\Sigma_1$. Let $T'^\s$ in $\smash{CRYS^\s_{top}(X'^\s/\Sigma)}$, $T^\s$ in $\smash{CRYS^\s_{top}(X^\s/\Sigma)}$ and $h:T'^\s\to T^\s$ a morphism of dp-thickenings which is compatible with $g$. There is a natural morphism $$h^*(M_{T^\s})\to(g^*M)_{T'^\s}$$ which is functorial with respect to $M$ in $\smash{Mod(\smash{(X^\s/\Sigma)_{CRYS^\s,top}},\O)}$. This is an isomorphism if $M$ is a crystal.
    The category of crystals is in particular stable by $g^*$. The same is true with $crys$ instead of $CRYS^\s$ with no restriction on $g$.
\item \label{lemcrystal2iii}  If $top$ is moreover coarser than or equal to $syn$, then the contravariant pseudo-functor $\smash{\cal Crys(\smash{(-/\Sigma)_{CRYS^\s,top}},\O)}:\cal Sch^\s/\Sigma_1\to \mathfrak Cat$ is  a stack for $top$ (ie. satisfies $top$ descent). The same statement holds verbatim with $crys$ instead of $CRYS^\s$.
\finrom
\end{lem}
Proof. Let us make a preliminary observation. Assume given a surjective family of $top$ morphisms  $(f_i:\smash{U_i^\s}\to X^\s)$ where the $\smash{U_i^\s}$'s are affine and possess fine charts. Let $\smash{T_{i,l}^\s}$ denote the logarithmic divided power envelope of $\smash{U_i^\s}$ inside a log scheme of the form $(Spec(\Z/p^l[(x_\lambda)_{\lambda\in \Lambda}][\N^e]),\N^e)$. If $k<\infty$ (resp. $k=\infty$) then the family $(\smash{U_i^\s},\smash{T_{i,k}^\s})_i$ (resp. $(\smash{U_i^\s},\smash{T_{i,l}^\s})_{i,l}$)  lies in the small crystalline $top$ site. If the $f_i$'s are in fact syntomic then it follows easily from the lifting property of syntomic schemes along dp-thickenings that the family $(\smash{U_i^\s},\smash{T_{i,k}^\s})_i$ (resp. $(\smash{U_i^\s},\smash{T_{i,l}^\s})_{i,l}$) covers the final object in the $\s$-big crystalline $top$ topos of $(X/\Sigma_k)$ as well as in the small one. If $top$ is finer than or equal to $et$ then we may always find $f_i$'s which are \'etale. In that case, we find in particular that the final object of the $\s$-big crystalline $top$ topos can be covered by objects of the small crystalline $top$ (and even $et$) site.

\ref{lemcrystal2i} If $T^\s$ is in $crys_{top}(X^\s/\Sigma)$, we have the following pseudo-commutative diagram of ringed topoi and weak morphisms
$$\xymatrix{((X^\s/\Sigma)_{CRYS^\s,top},\O)\ar[d]_-\pi&((X^\s/\Sigma)_{CRYS^\s,top}/T^\s,\O)\ar[rd]^-{\lambda_{T^\s}}
\ar[l]_-{f_{T^\s}}\ar[d]_-{\pi_{|T^\s}}&\\
((X^\s/\Sigma)_{crys,top},\O)&((X^\s/\Sigma)_{crys,top}/T^\s,\O)
\ar[l]^-{f_{T^\s}}\ar[r]_-{\lambda_{T^\s}}&(T^\s_{top},\O)}$$
and the following functors are fully faithful: $$\begin{array}{l}\pi^*:Mod((X^\s/\Sigma)_{crys,top},\O)\to Mod((X^\s/\Sigma)_{CRYS^\s,top},\O)\\
\pi_{|T^\s}^*:Mod((X^\s/\Sigma)_{crys,top}/T^\s,\O)\to ((X^\s/\Sigma)_{CRYS^\s,top}/T^\s,\O)\\
\lambda_{T^\s}^*:Mod(T^\s_{top},\O)\to Mod((X^\s/\Sigma)_{crys,top}/T^\s,\O)\\
\lambda_{T^\s}^*:Mod(T^\s_{top},\O)\to Mod((X^\s/\Sigma)_{CRYS^\s,top}/T^\s,\O)\end{array}$$
According to Def. \ref{defcrystal} (resp. Def. \ref{defcrystal} and Lem. \ref{lemcrystal1}\ref{lemcrystal1iii} together with the above preliminary observation), a module $M$ of $\smash{(X^\s/\Sigma)_{crys,top},\O)}$ (resp. $\smash{(X^\s/\Sigma)_{CRYS^\s,top},\O)}$) is a crystal if and only if $\smash{f_{T^\s}^*M}$ is in the essential image of $\smash{\lambda_{T^\s}^*}$ for all $T^\s$'s in the small $top$ crystalline site. With the above diagram in mind, one sees immediately that the condition of being a crystal is preserved by the functors $\pi^*$ and $\pi_*$. It remains to check that $\pi_*$ is fully faithful on the category $\smash{\cal Crys((X^\s/\Sigma)_{CRYS^\s,top},\O)}$. If $M$ is in the latter category and $T^\s$ is in the small crystalline site then $\smash{f_{T^\s}^*M}$ is in the essential image of $\smash{\lambda_{T^\s}^*}$, hence of $\smash{\pi_{|T^\s}^*}$, ie. $$\pi_{|T^\s}^*\pi_{|T^\s,*}f_{T^\s}^*M\simeq f_{T^\s}^*M$$ Using the preliminary observation again, it follows that $\pi^*\pi_*M\simeq M$ as desired.

\ref{lemcrystal2ii} Let us use the notations $\smash{f_{T^\s/X^\s}}=f_{T^\s}$ and $\smash{\lambda_{T^\s/X^\s}}=\lambda_{T^\s}$ in order to emphasize the dependance in $X^\s$.
We have the following pseudo-commutative diagram of ringed topoi and weak morphisms:
$$\xymatrix{((X'^\s/\Sigma)_{CRYS^\s,top},\O)\ar[d]_-g&&((X'^\s/\Sigma)_{CRYS^\s,top}/T'^\s,\O)
\ar[ll]_-{f_{T'^\s/X'^\s}}\ar[rrd]^-{\lambda_{T'^\s/X'^\s}}\ar[d]^-\wr&&
\\ ((X^\s/\Sigma)_{CRYS^\s,top},\O)&&((X^\s/\Sigma)_{CRYS^\s,top}/T'^\s,\O)\ar[d]^-h\ar[ll]_-{f_{T'^\s/X^\s}}
\ar[rr]_-{\lambda_{T'^\s/X^\s}}&&(T'^\s_{top},\O)\ar[d]^-h\\
&&((X^\s/\Sigma)_{CRYS^\s,top}/T^\s,\O)\ar[llu]^-{f_{T^\s/X^\s}}
\ar[rr]_-{\lambda_{T^\s/X^\s}}&&(T^\s_{top},\O)}$$
If $M$ is a module of $\smash{((X^\s/\Sigma)_{CRYS^\s,top},\O)}$ the desired morphism is defined by composition as follows: $$h^*\lambda_{T^\s/X^\s,*}f_{T^\s/X^\s}^*M\to \lambda_{T'^\s/X^\s,*}h^*f_{T^\s/X^\s}^*M\simeq \lambda_{T'^\s/X^\s,*}f_{T'^\s/X^\s}^*M\simeq \lambda_{T'^\s/X'^\s,*}f_{T'^\s/X'^\s}^*g^*M$$
If $M$ is a crystal then $\smash{f_{T^\s/X^\s}^*M}$ is in the essential image of  $\smash{\lambda_{T^\s/X^\s}^*}$ and the first arrow is thus invertible (by full faithfulness of $\smash{\lambda_{T^\s/X^\s}^*}$ and $\smash{\lambda_{T'^\s/X^\s}^*}$).


\ref{lemcrystal2iii} Thanks to Lem. \ref{lemcrystal1}\ref{lemcrystal1iii} and the preliminary observation we find that the property of being a crystal is local for $top$ as long as $top$ is finer than or equal to $et$ and coarser than or equal to $syn$.
\begin{flushright}$\square$\end{flushright}

\begin{rem} If $top$ is coarser than $et$ (e.g. if $top=zar$), then Lem. \ref{lemcrystal2}, \ref{lemcrystal2i}, \ref{lemcrystal2ii}  are still valid if we impose that $X^\s$ has charts $top$ locally. In this setting the analogue of \ref{lemcrystal2iii} is that we get a stack for $top$ over the full subcategory of $\cal Sch/\Sigma_1$ having $top$ local charts.
\end{rem}

\section{Preliminaries part II: some properties for sheaves}\label{secppii}

\subsection{$1$-motives and $p$-divisible groups} \label{pdg} ~~ \\

We review some well known definitions and properties concerning finite locally free groups, $p$-divisible groups, Abelian schemes and $1$-motives.

\para \label{deffullyexact}  Part of this reminder is devoted to put an exact structure on the relevant categories. The reader may consult \cite{Bu} for a reminder of basic facts concerning exact categories. Let us introduce a terminology which is adapted both for our present purpose and for later ones as well (Sect. \ref{epotcoc}). Consider a functor $F$ between exact categories $\C_1$, $\C_2$ whose exact structures are both denoted $e$ for simplicity. We say that $F$ is \emph{$e$-exact} if it sends short $e$-exact sequences of $\cal C_1$ to short $e$-exact sequences of $\cal C_2$. We say that it \emph{reflects $e$-exactness} if a short sequence of $\cal C_1$ whose image by $F$ is  $e$-exact in $\cal C_2$ is already $e$-exact. A \emph{fully $e$-exact subcategory} of $\cal C_1$ is a full additive subcategory which is closed by extensions and equipped with the exact structure induced by $e$.

If $\cal C_1$ and $\cal C_2$ are Abelian and $e$ denotes their natural exact structure, then $F$ is $e$-exact if and only if it is exact in the usual sense (ie. if it commutes to finite limits). An e-exact functor $F$ reflects $e$-exactness if and only if it is moreover faithful. A \emph{fully Abelian subcategory} of an Abelian category is a full subcategory which is Abelian, closed by extensions and which is such that the inclusion functor is exact. In other words it is a fully $e$-exact subcategory which is stable by kernels and cokernels.


\para Let us fix some terminology.

\begin{defn} \label{defflfgppdiv}
Let $X$ be a scheme.

\debrom \item \label{defflfgppdivi}  A \emph{finite locally free group} over $X$ is an Abelian group of $X_{FL}$ which is representable by a finite locally free group scheme.

\item \label{defflfgppdivii} A \emph{$p$-divisible group $G$} over $X$ is an Abelian group of $X_{FL}$ which is such that $p:G\to G$ is epimorphic, the kernel $G_{p^k}$ of $p^k:G\to G$ is a finite locally free group for any $k$ and  $limind G_{p^k}\simeq G$. \finrom
\end{defn}

\begin{rem} \debrom \item  Finite locally free is equivalent to finite flat and locally of finite presentation.

\item Finite implies affine (\cite{EGA2} Prop. 6.1.4). \finrom
\end{rem}

The \emph{Cartier dual} of a finite locally free (resp. $p$-divisible) group $G$ is defined as follows: \begin{eqnarray}\label{defcartierduality1} &G^*:=\cal Hom(G,\Gm)&\\ \label{defcartierduality2} \hbox{(resp.} &G^*:=\limi_k (G_{p^k})^*&\hbox{)}\end{eqnarray}

Using that $\cal Ext^1(G,\Gm)$ vanishes for any finite locally free group $G$ (\cite{SGA7-I} VIII, Prop. 3.3.1) we find an exact sequence
$$\xymatrix{0\ar[r]&(G_{p^k})^*\ar[r]^{(p^l)^*}&(G_{p^{k+l}})^*\ar[r]^{(p^k)^*}&(G_{p^{k+l}})^*\ar[r]&
(G_{p^k})^*\ar[r]&0}$$ for a $p$-divisible group $G$. It follows in particular that $(G_{p^k})^*$ is naturally isomorphic to $(G^*)_{p^k}$ and that
$G^*$ is a $p$-divisible group. This implies that the obvious biduality isomorphism $G\simeq G^{**}$ for finite locally free groups naturally extends
to $p$-divisible groups.

\begin{lem} \label{gpflfexact} \debrom \item \label{gpflfexacti} Finite locally free groups form a fully $e$-exact subcategory $flf(X)$ of the category of Abelian groups of $X_{FL}$. This is in fact a fully Abelian subcategory if $X=Spec(K)$ for some field $K$.

\item \label{gpflfexactii} $p$-divisible groups form a fully $e$-exact subcategory $pdiv(X)$ of the category of Abelian groups of $X_{FL}$.
    \finrom
\end{lem}
Proof. \ref{gpflfexacti} Let us prove the first assertion. Consider an exact sequence $0\to G\to G'\to G''\to 0$ where $G$ and $G''$ are finite
locally free groups. Looking at $G'$ as an object of $\smash{G''_{FL}}\simeq X_{FL}/G''$ by \cite{SGA4-I} III, Prop. 5.4, we find that it is a torsor under
the pullback of $G$ to $G''$.
It is thus represented by a finite locally free group scheme over $G''$ in virtue of Cor. 5.7 and Cor. 7.9 of \cite{SGA1} VIII.
We conclude by composition with the finite locally free morphism $G''\to X$. The second assertion follows from \cite{Gr3} 212-17, Cor. 7.3. Statement \ref{gpflfexactii} follows from \ref{gpflfexacti}.
\begin{flushright}$\square$\end{flushright}

\begin{lem} \label{flfsyn} A finite locally free group scheme is syntomic over its base scheme $X$. \end{lem}
Proof. Since finite locally free groups are flat and locally of finite presentation we may assume (\cite{EGA4-IV}, Def. 19.3.6) that $X=Spec(K)$ for some
field $K$. By descent, we may furthermore assume that $K$ is algebraically closed (\cite{EGA4-IV} Prop. 19.3.9 (ii)). In that case, the category of finite
locally free groups is fully Abelian in  $Ab(X_{FL})$ (Lem. \ref{gpflfexact} \ref{gpflfexacti}). Each one of its objects admits moreover a composition
series with simple quotients isomorphic to either one of the groups $\mu_p$, $\alpha_p$ or $\Z/l$ ($l$ any prime number) by \cite{Oo} Chap. I.2. We may thus conclude by d\'evissage, since the property of being syntomic is preserved by forming extensions (use the same arguments as in the
proof of Lem. \ref{gpflfexact}).
\begin{flushright}$\square$\end{flushright}

\para The following definitions are standard. 

\begin{defn} \label{defmotabtor} Let $X$ be a scheme.
\debrom \item  \label{defmotabtori} An \emph{Abelian scheme} (or \emph{Abelian variety} if $X$ is the spectrum of a field)  is an Abelian group of $X_{FL}$ which is representable by a proper smooth $X$-group scheme whose fibers are geometrically connected.

\item  \label{defmotabtoriii} A \emph{torus} is an Abelian group of $X_{FL}$ which is representable and locally isomorphic to a finite direct sum of $\Gm$'s.

\item \label{defmotabtorii} A \emph{semi-Abelian scheme} is an Abelian group of $X_{FL}$ which is representable by a smooth $X$-group scheme whose fibers are successive extensions of an \'etale finite locally free group, an Abelian variety and a torus.

\item  \label{defmotabtoriv} A \emph{twisted constant group}  is an Abelian group of $X_{FL}$ which is representable and locally isomorphic to a finite direct sum of $\Z$'s.

\item  \label{defmotabtorv} A \emph{$1$-motive} is a complex $[\G\to G]$ of Abelian groups of $X_{FL}$ placed in degrees $[0,1]$ where $\G$ is a twisted constant group and $G$ is the extension of an Abelian scheme by a torus.
    \finrom
\end{defn}

\begin{lem} \label{lemG} \debrom \item \label{lemGi} A torus (resp. twisted constant groups) becomes isomorphic to a finite direct sum of $\Gm$'s (resp. $\Z$'s) over a suitable \'etale covering.

\item \label{lemGii} If $T$ is a torus and $A$ is an Abelian scheme then every morphism $T\to A$ or $A\to T$ is trivial.

\item \label{lemGiii} If $G$ is an extension of an Abelian scheme $B$ by a torus $T$ then $G$ is automatically representable by a smooth group scheme. Moreover, $G$ admits a unique maximal subtorus in the sense of \cite{SGA3-II} XII, Def. 1.3 (see also Rem. 1.4, a' of loc. cit.) which is precisely $T$ and varies functorially with respect to $G$.
\finrom
\end{lem}

Proof. \ref{lemGi} This follows from 
 \cite{SGA3-II} X, Cor. 4.5 (resp. \cite{SGA3-II} X, Prop. 5.3). 

\ref{lemGii} Use \cite{SGA3-II} XV, Lem. 8.3 and the Thm. 3.1 and Cor. 3.8 of \cite{Mi3} in the first case and  \cite{SGA7-I} VIII, Sect. 3.2 in the second case. 

\ref{lemGiii} The first assertion results easily from \cite{SGA1} VIII, Cor. 7.9. 

The existence of the maximal torus results from \cite{SGA3-II} XV, Cor. 8.17], 
     and its unicity follows from \cite{SGA3-II} XII, Thm. 7.1, b. 
The remaining statements follows from  \ref{lemGii}. 
\begin{flushright}$\square$\end{flushright}

\begin{lem} \label{1motex} Assume that $X$ is regular. The category  $\cal M(X)$ of $1$-motives is fully $e$-exact in the category of complexes of Abelian groups of $X_{FL}$.
\end{lem}
Proof. We have to show that $\cal M(X)$ is stable by extensions. The category of twisted constant groups is stable by extensions thanks to \cite{SGA3-II} X, Lem. 5.4. 
The only point deserving explanations is thus the following \petit




\emph{Claim}. Consider an exact sequence $0\to G'\to G\to G''\to 0$ of Abelian groups in $X_{FL}$. If $G'$, $G''$ are respectively an extension of an Abelian scheme $B'$, $B''$ by a torus $T'$, $T''$ then $G$ is an extension of an Abelian scheme $B$ by a torus $T$ as well. We have moreover $dim_XT=dim_XT'+dim_XT''$ and $dim_XB=dim_XB'+dim_XB''$. \petit

First, we note that $G$ is representable by a smooth group scheme (use \cite{SGA1} VIII, Cor. 7.9 to reduce to the case where $G'$ is an Abelian scheme and then conclude by \cite{Ra1} Prop. XIII, 2.6). \petit

Let us now prove the claim assuming that $G$ has a maximal subtorus, say $T$. In virtue of \cite{SGA3-II} XV, Lem. 8.3, 
the kernel $H'$ of $T\to G''$ is a multiplicative subgroup of finite type in $G'$ (\cite{SGA3-II} IX, Def. 1.1, 1.2).
It clearly contains $T'$ and the quotient $H'/T'$ is a multiplicative subgroup of finite type in $B'$ (\cite{SGA3-II} IX, Prop. 2.7).
It follows easily from Lem. \ref{lemG} \ref{lemGii} that it is in fact finite locally free. The image of $T\to G''$ on the other hand, is a subtorus of $G''$ (\cite{SGA3-II} XV, Lem. 8.3) and is in fact equal to $T''$ (\cite{SGA3-II} XII, Thm. 7.1, e)). 
Forming cohomology of the vertical complexes in the following commutative diagram with exact lines: \begin{eqnarray}\label{exG1}\xymatrix{0\ar[r]&T'\ar[r]\ar[d]&G'\ar[r]\ar[d]&B'\ar[r]\ar[d]&0\\
0\ar[r]&T\ar[r]\ar[d]&G\ar[r]\ar[d]&G/T\ar[r]\ar[d]&0\\
0\ar[r]&T''\ar[r]&G''\ar[r]&B''\ar[r]&0}\end{eqnarray}
produces an exact sequence \begin{eqnarray}\label{secAQSD}\xymatrix{0\ar[r]&H'/T'\ar[r]&B'\ar[r]&G/T\ar[r]&B''\ar[r]&0}\end{eqnarray}
Since $X$ is normal, a result of Grothendieck (\cite{Ra1} Sect. XI, 1) ensures that $B'$ is projective. It thus follows from \cite{Ra2} 5. Thm. 1, applications (a) (iii) that the quotient of $B'$ by the finite locally free group $H'/T'$ is an Abelian scheme. As a result $B:=G/T$ is the extension of two Abelian schemes and is thus an Abelian scheme as well (using \cite{Ra1} Prop. XIII, 2.6 again). The assertion about dimensions is clear from the exact sequence (\ref{secAQSD}), since $dim_XH'/T'=0$. \petit

The above proof applies in particular to the case where $X$ is the spectrum of a field. If now $X$ is regular we have to show that $G$ has a maximal subtorus. This follows from \cite{SGA3-II} XV, Cor. 8.17, thanks to the dimension formulae in the fibers.

\begin{flushright}$\square$\end{flushright}

\begin{lem} \label{pdivfunct} \debrom \item \label{pdivfuncti} Consider Abelian schemes $A$, $A'$ and an isogeny $f:A\to A'$ (ie. $f$ is a finite faithfully flat morphism of group schemes (\cite{SGA3-III} XXII, Def. 4.2.9). The kernel of $f$ is a finite locally free group.

\item \label{pdivfunctii} If $M=[\G\to G]$ is a $1$-motive then $M_{p^\infty}:=\Qp/\Zp\otimes^LM$ is concentrated in degree $0$ and is a $p$-divisible group. If $X$ is regular, the resulting functor $$(-)_{p^\infty}:\M(X)\to pdiv(X)$$ is $e$-exact in the sense of Lem. \ref{1motex}.

\item \label{pdivfunctiii} Consider a semi-Abelian scheme $A$. Its connected component $A^0$  (see \cite{SGA3-I} VIB Def. 3.1) is semi-Abelian as well. Multiplication by $p$ on $A$ is flat locally of finite presentation and locally quasi-finite. It is moreover surjective if $A=A^0$.
\finrom \end{lem}
Proof. \ref{pdivfuncti} Since $f$ is finite and flat, the same is true for $ker\, f$. The latter group scheme is furthermore locally of finite presentation since the unit section $X\to A'$ is a closed immersion (thanks to $A'/X$ being separated) whose ideal is finitely generated.

\ref{pdivfunctii} The case where $M=[0\to A]$ with $A$ an Abelian scheme follows from \ref{pdivfuncti}, thanks to the well known fact that  $p:A\to A$ is an isogeny  (\cite{Gr2}). The general case proceeds from there by d\'evissage using Lem. \ref{gpflfexact}.

\ref{pdivfunctiii} According to \cite{SGA3-I} VIB, Thm. 3.10,
$A^0$ is representable by an open subgroup scheme. The first statement follows from this and the definition. The other statements follow easily from \ref{pdivfunctii} applied to the connected components of the fibers (use \cite{EGA4-III} Thm. 11.3.10 for flatness).
\begin{flushright}$\square$\end{flushright}



%



\begin{rem} In fact a result of Raynaud ensures that every finite locally free group arises as the kernel of an isogeny between Abelian schemes after suitable Zariski localization on $X$ (\cite{BBM} Thm. 3.1.1). \end{rem}

\vspace{.5cm}

Let us finally discuss the group of components when the base is a curve.

\begin{lem} \label{components} Assume that $X$ is locally Noetherian of dimension $1$. Let $z:Z\to X$ denote the inclusion morphism of a reduced closed subcheme and let $U$ denote its complement. Consider a semi-Abelian scheme $A/X$ whose restriction to $U$ is Abelian. Its \emph{component group}  $\Phi:=A/A^0$ satisfies the following.

\debrom \item \label{componentsi} The sheaf $z^{-1}\Phi$ is representable by an \'etale group scheme.

\item \label{componentsii} The adjunction morphism $\Phi\to z_*z^{-1}\Phi$ is invertible.
\finrom
\end{lem}
Proof. Under our assumptions, $Z$ is a discrete scheme. Statement \ref{componentsi} thus follows simply from the fact that the formation of $A^0$ is compatible to base change (see \cite{SGA3-I} VIB Prop. 3.3).

\ref{componentsii} Consider the following commutative diagram with exact lines in $Ab(X_{FL})$:
$$\xymatrix{0\ar[r]&A^0\ar[r]\ar[d]&A\ar[r]\ar[d]&\Phi\ar[r]\ar[d]&0\\
0\ar[r]&z_*z^{-1}A^0\ar[r]&z_*z^{-1}A\ar[r]&z_*z^{-1}\Phi}$$
It follows immediately from the definition of $A^0$ (see \cite{SGA3-I} VIB Def. 3.1) that the kernel of the middle vertical morphism is contained in $A^0$. It thus remains only to observe that the right vertical morphism is epimorphic, thanks to \ref{componentsi} and \cite{SGA1}, I, Prop. 8.1.
\begin{flushright}$\square$\end{flushright}

\subsection{Technical remarks on certain $top$ sheaves}\label{topsheaves} ~~ \\

We review basic facts about constant sheaves and $p$-divisible groups viewed as sheaves on usual or crystalline sites for certain  $top$.

\para Let us begin with constant sheaves.

\begin{lem} \label{lemcst} Assume that $fl \preceq top \preceq zar$.
\debrom \item \label{lemcsti} Let $X^\s$ in $\cal Sch^\s$ and let us (temporarily) denote $\smash{\underline M}$ in $\smash{X^\s_{TOP^\s}}$  the constant sheaf associated to a set $M$. Then $\underline M$ is representable by the $X^\s$-scheme  $X^{\s,(M)}$ consisting of a direct sum of copies of $X^\s$ indexed by $M$. The same is true in the setting of big or small topoi.

\item \label{lemcstii} Consider the morphism $i:X^\s_{TOP^\s}\to (X^\s/\Sigma)_{CRYS^\s,top}$. Then $i_*X^{\s,(M)}$ is the constant sheaf  associated to $M$. The same is true in the setting of big or small topoi.

\item \label{lemcstiii} The formation of constant sheaves is compatible with:

 - the restriction functors $\pi_*$ between the $\s$-big, big and small crystalline or usual $top$ topoi (note that compatibility with $\pi^{-1}$ is tautological).

 - the change of topology functors $\epsilon_*$ (compatibility with $\epsilon^{-1}$ is tautological as well).

 - the functor $i_*$ from usual to crystalline topoi (compatibility with $i^{-1}$ is tautological as well).
\finrom
\end{lem}
Proof. We only prove \ref{lemcsti} since \ref{lemcstii} is similar and \ref{lemcstiii} follows from \ref{lemcsti} and \ref{lemcstii}. Let $\check{H}{}^0(\cal U,-)$ denote the Cech functor associated to a covering $\cal U=(\smash{U^\s_\lambda\to U^\s})_{\lambda\in \Lambda}$ of some $U^\s$ in $\cal Sch^\s$ and let $\check{ H}{}^0(U,-)$ denote the direct limit of the previous functors for $\cal U$ running in the category of coverings of $U^\s$. If $M_0$ denotes the constant presheaf associated to $M$ then $M_1:=U^\s\mapsto \check{H}{}^0(U^\s,M_0)$ is the separated presheaf associated to $M$ and $U^\s\mapsto \check{H}{}^0(U^\s,M_1)$ is the sheaf $\underline M$ associated to $M$. These may be computed as follows. The value of $M_1$ on $U^\s$ is the final set if $U$ is the empty scheme and $M$ otherwise. Next $\check{H}{}^0(\cal U,M_1)$ is a product of copies of $M$ indexed by the quotient $\overline \Lambda$ of $\Lambda$ by the equivalence relation $\lambda\simeq \lambda'$ $\Leftrightarrow$ $\exists n\ge 1$, $\exists \lambda_0,\dots, \lambda_n$, $\lambda_0=\lambda$, $\lambda_n=\lambda'$,  $\forall 0\le k\le n-1$, $\smash{U_{\lambda_k}\times_U U_{\lambda_{k+1}}}\ne \emptyset$. Consider $\lambda\in \Lambda$ with equivalence class $\overline \lambda$ in $\overline \Lambda$. Define $\smash{U_{\overline \lambda}^\s}$ as the open sub log scheme of $U^\s$ whose underlying topological space is the union of the images of the $U_{\lambda'}$'s with $\lambda'\simeq \lambda$. This defines an open partition $\cal P_{\cal U}$ of $U^\s$. A morphism $f:U^\s\to X^{\s,(M)}$ on the other hand, defines an open partition $\cal P_f$ of $U^\s$ as well. With this in mind, we find that the natural morphism $\check{H}{}^0(\cal U,M_1)\to Hom(U^\s,X^{\s,(M)})$ is injective and that $f$ is in its image if and only if $\cal P_f$ is refined by $\cal P_{\cal U}$. The result follows, since every open partition is of the form $\cal P_{\cal U}$ for an adequate covering $\cal U$.
\begin{flushright}$\square$\end{flushright}

In view of this lemma the constant sheaf on a usual or crystalline site represented by a set or group $M$ will simply be denoted $M$ from now on.



\para It is sometimes useful to view locally free groups or $p$-divisible groups as sheaves on big or small $top$ (crystalline or usual) sites rather than on the big $fl$ usual site. The following lemma will be helpful.


\begin{lem} \label{lemlimind} Assume that $fl\preceq top \preceq top'\preceq zar$ and let $X^\s$ in $\cal Sch^\s$. The functors $\epsilon_*:\smash{X^\s_{TOP^\s}}\to \smash{X^\s_{TOP'^\s}}$, $\epsilon_*:\smash{(X^\s/\Sigma)_{CRYS^\s,top}}\to \smash{(X^\s/\Sigma)_{CRYS^\s,top'}}$ and $i_*:\smash{X^\s_{TOP^\s}}\to \smash{(X^\s/\Sigma)_{CRYS^\s,top}}$ commute to filtrant inductive limits. The same holds in the context of big or small sites.
\end{lem}
Proof. Consider a filtrant inductive system $(F_i)$ in $\smash{X^\s_{TOP^\s}}$ (resp. $\smash{(X^\s/\Sigma)_{CRYS^\s,top}}$). Let $F$ denote the direct limit computed in the category of presheaves and let $G$ denote the associated sheaf for the $zar$ topology. \petit

\emph{Fact}. The presheaf $G$ on $TOP^\s(X^\s)$ (resp. $CRYS^\s_{top}(X^\s/\Sigma)$) is in fact a sheaf (for $top$). It is in particular the direct limit computed in $\smash{X^\s_{TOP^\s}}$ (resp. $\smash{(X^\s/\Sigma)_{CRYS^\s,top}}$).
\petit

We may use the characterization of $top$ sheaves given in Prop. 6.2.1 and Cor. 6.2.3 of \cite{SGA3-I} IV in the following situation (with the notations of \emph{loc. cit.}): $\cal C$ is the category underlying $TOP^\s(X^\s)$ (resp. $CRYS^\s_{top}(X^\s/\Sigma)$), $\cal C'$ is the full subcategory of $\cal C$ formed by the $U^\s/X^\s$'s where $U$ is an affine scheme (resp. the $(U^\s/X^\s,T^\s,\iota,\gamma)$'s where $T$ is an affine scheme), $P$ is the set of open coverings and $P'$ is the set of finite surjective families of $top$ strict morphisms (resp. of cartesian $top$ strict morphisms) of finite presentation in $\cal C'$. We must check that $G$ satisfies the descent condition for the families which belong to  $P$ or $P'$. To start with, we notice that $F$ satisfies the descent condition for the families which are elements of $P'$. Indeed in the category of sets, filtrant inductive limits commute to finite projective limits. Fact will thus be proven if we check that $F$ and $G$ coincide on the objects of $\cal C'$. Now $F$ is clearly a separated presheaf for $zar$ and $G$ can thus be computed from $F$ by applying the Cech functor only once. The Cech computation in question can be done using only finite families thanks to quasi-compacity and Fact follows.

It follows immediately from Fact that $\epsilon_*$ commutes to filtrant inductive limits. The same is true for  $i_*$, as seen by reduction to the case $top=zar$ (where $i_*$ commutes to arbitrary limits) using this together with the full faithfulness of $\epsilon_*:\smash{(X^\s/\Sigma)_{CRYS^\s,top}}\to \smash{(X^\s/\Sigma)_{CRYS^\s,zar}}$: $$\begin{array}{rcl}\limi i_*F_i&\simeq &\limi \epsilon^{-1}\epsilon_*i_*F_i \\ &\simeq &
\epsilon^{-1}\limi \epsilon_*i_*F_i\\ & \simeq & \epsilon^{-1}\limi i_*\epsilon_* F_i \\ &\simeq & \epsilon^{-1} i_*\limi \epsilon_* F_i \\ & \simeq & \epsilon^{-1} i_*\epsilon_* \limi F_i \\ & \simeq & \epsilon^{-1}\epsilon_*i_*\limi F_i \\ & \simeq &
i_*\limi F_i\end{array}$$

The case where $TOP^\s$ and $CRYS^\s$ are replaced respectively by $TOP$ and $CRYS$  or $top$ and $crys$ follows formally using that $p$ pseudo-commutes to $\epsilon$ and $i$, $p^{-1}$ is fully faithful and $p_*$ commutes to inductive limits.
\begin{flushright}$\square$\end{flushright}

\begin{cor} \label{corpdivtop} Assume that $fl\preceq top \preceq syn$. Consider the functor $\epsilon_*$ (resp. $\pi_*\epsilon_*$, resp. $i_*\epsilon_*$, resp. $i_*\pi_*\epsilon_*$)
from $Ab(X_{FL})$ to $Ab(X_{TOP})$ (resp. $Ab(X_{top})$, resp. $Ab((X/\Sigma)_{CRYS,top})$, resp. $Ab((X/\Sigma)_{crys,top})$). We denote  $flf_{TOP}(X)$ and $pdiv_{TOP}(X)$ (resp. $flf_{top}(X)$ and $pdiv_{top}(X)$, resp. $flf_{CRYS,top}(X/\Sigma)$ and $pdiv_{CRYS,top}(X/\Sigma)$, resp.  $flf_{crys,top}(X/\Sigma)$ and $pdiv_{crys,top}(X/\Sigma)$) the essential image of $flf(X)$ and $pdiv(X)$.

\debrom
\item \label{corpdivtopi} Let $G$ in $pdiv(X)$. Multiplication by $p$ is epimorphic on $\epsilon_*G$ (resp. $\pi_*\epsilon_*G$, resp. $i_*\epsilon_*G$, resp. $i_*\pi_*\epsilon_*G$).

\item \label{corpdivtopii} The restriction of the functor $\epsilon_*$ (resp. $\pi_*\epsilon_*$, resp.  $i_*\epsilon_*$, resp. $i_*\pi_*\epsilon_*$) to  $flf(X)$ or $pdiv(X)$ is fully faithful and $e$-exact.


\item \label{corpdivtopiii} Assume now that $fl\preceq top'\preceq top\preceq syn$. The adjunctions $(\epsilon^{-1}, \epsilon_*)$, $(\pi^{-1},\pi_*)$ and  $(i^{-1},i_*)$ induce equivalences between the vertices of the following cube:
$$\xymatrix{&flf_{top'}(X)\ar[rr]^(.25){(i^{-1},i_*)}\ar'[d]^(.5){(\epsilon^{-1}, \epsilon_*)}[dd]&&flf_{crys,top'}(X/\Sigma)\ar[dd]^(.25){(\epsilon^{-1}, \epsilon_*)}\\
flf_{TOP'}(X)\ar[dd]^(.25){(\epsilon^{-1}, \epsilon_*)}\ar[rr]^(.25){(i^{-1},i_*)}\ar[ru]^(.35){(\pi^{-1},\pi_*)}&&flf_{CRYS,top'}(X/\Sigma)\ar[ru]^(.35){(\pi^{-1},\pi_*)}
\ar[dd]^(.25){(\epsilon^{-1}, \epsilon_*)}&\\
&flf_{top}(X)\ar'[r]^(.5){(i^{-1},i_*)}[rr]&&flf_{crys,top}(X/\Sigma)\\
flf_{TOP}(X)\ar[rr]^(.25){(i^{-1},i_*)}\ar[ru]^(.35){(\pi^{-1},\pi_*)}&&flf_{CRYS,top}(X/\Sigma)\ar[ru]^(.35){(\pi^{-1},\pi_*)}&}$$
and similarly for $pdiv$.

\finrom
\end{cor}
Proof. \ref{corpdivtopi} Since $\pi_*$ preserves epimorphisms, it is sufficient to consider the case of $\epsilon_*G$ and $i_*\epsilon_*G$. We may furthermore assume that $top=syn$. In that case, $i_*$ preserves epimorphisms as well and it suffices to look at the case of $\epsilon_*G$. Now the result follows immediately from the fact that the morphism of group schemes representing $p:\epsilon_*G_{p^{k+1}}\to \epsilon_*G_{p^k}$ is surjective and syntomic (indeed $G_{p^{k+1}}\to G_{p^k}$ is a torsor of $(G_{p^k})_{FL}$ under the syntomic group scheme $G_p$ (Lem. \ref{flfsyn})).

\ref{corpdivtopii} Let $v$ denote either one of the weak morphisms of topoi $\epsilon$, $\pi\epsilon$, $i\epsilon$, $i\pi\epsilon$. Let us explain full faithfulness. It is sufficient to show that the adjunction morphism $v^{-1}v_*G\to G$ is invertible whenever $G$ is in $flf(X)$ or $pdiv(X)$. Thanks to Lem. \ref{lemlimind}, it suffices to consider the case of $flf(X)$. Now the functors $\epsilon_*:Ab(X_{FL})\to Ab(X_{TOP})$, $i_*:Ab(X_{TOP})\to Ab((X/\Sigma)_{CRYS,top})$
and $i_*:Ab(X_{top})\to Ab((X/\Sigma)_{crys,top})$ being fully faithful themselves, we are reduced to the case of $\pi_*:flf_{TOP}(X)\to flf_{top}(X)$. The result then follows from Lem. \ref{flfsyn}.

Let us explain $e$-exactness. Consider a short $e$-exact sequence $0\to G'\to G\to G''\to 0$ of $flf(X)$ or $pdiv(X)$. We have to show that the image of $G\to G''$ by $v_*$ is  epimorphic. Passing to $p^k$-torsion points in the case of $pdiv(X)$, we are reduced to the case of $flf(X)$. Since $\pi_*:Ab(X_{TOP})\to Ab(X_{top})$ and $\pi_*:Ab((X/\Sigma)_{CRYS,top}\to Ab((X/\Sigma)_{crys,top}$ are exact, we may restrict our attention to big topoi. We may then assume that $top=syn$. The case $v=\epsilon$ results from Lem. \ref{flfsyn} and the case $v=i\epsilon$ follows by exactness of $i_*:Ab(X_{SYN})\to Ab((X/\Sigma)_{CRYS,syn})$.

\ref{corpdivtopiii} It remains to show that the functors $\epsilon^{-1}:Ab(X_{TOP})\to Ab(X_{FL})$, $\pi^{-1}:Ab(X_{top})\to Ab(X_{TOP})$, $\pi^{-1}:Ab((X/\Sigma)_{crys,top})\to Ab((X/\Sigma)_{CRYS,top})$,
$i^{-1}:Ab((X/\Sigma)_{CRYS,top})\to Ab(X_{TOP})$ and $i^{-1}:Ab((X/\Sigma)_{crys,top})\to Ab(X_{top})$ preserve the corresponding categories of finite locally free groups and $p$-divisible groups. Let $v$ denote one of the five weak morphisms $\epsilon$, $\pi$, $i$ involved. We have to check that the adjunction morphism $v^{-1}v_*G\to G$ is invertible when $G$ is a locally free group or a $p$-divisible group. By Lem.  \ref{lemlimind}, we may assume that $G$ is a finite locally free group. Aside from the cases already treated in \ref{corpdivtopii}, it remains to explain the case
$v=\pi:(X/\Sigma)_{CRYS,top}\to (X/\Sigma)_{crys,top}$. For $G$ in $flf_{TOP}(X)$ we will show that $\pi^{-1}\pi_*i_*G\simeq i_*G$. We may always assume that $\Sigma=\Sigma_k$ with $1\le k<\infty$ and that $X$ (hence $G$) is an affine scheme. Choose a closed immersion $G\hookrightarrow Y$ where $Y=Spec(\Z/p^k[(x_i)_{i\in I}])$ ($I$ possibly infinite) and let $(G,T)$, $\smash{(G,T^{(1)})}$ respectively denote the divided power envelope of $G$ into $Y$, $Y\times Y$. The result follows from the observation that $i_*G$ (resp. $\pi_*i_*G$) may be described as the cokernel of the couple of morphisms $\smash{(G,T^{(1)})}\rightrightarrows (G,T)$ (induced by projections) in the topos $(X/\Sigma)_{CRYS,top}$ (resp. $(X/\Sigma)_{crys,top}$).
\begin{flushright}$\square$\end{flushright}

\subsection{Projective systems and $p$-divisibility over a $\Z/p^.$-algebra} \label{psapdozpa}


\para \label{ARmod}
Consider a ringed variable topos on $\N$: $k\mapsto E_k$, $(k'\ge k)\mapsto (\iota_{k,k'}:(E_k,A_k)\to \smash{(E_{k'},A_{k'})})$ and let $(E_.,A_.)$ denote the associated ringed total topos.
An object $F_.$  of $E_.$ (resp. a module $M_.$ of $(E_.,A_.)$) is thus a projective system (in the sense of the associated cofibered category) whose $k^{th}$ term $F_k$ (resp. $M_k$) is in $E_k$ (resp. $Mod(E_{k},A_{k})$). Concretely speaking, it thus consists of the data of the $F_k$'s (resp. $M_k$'s) together with transition morphisms $\iota_{k,k'}^{-1}F_{k'}\rightarrow F_k$ (resp. $\iota_{k,k'}^*M_{k'}\rightarrow M_k$) for $k'\ge k$ satisfying the obvious composition constraint.

\begin{defn} \label{defj} If $j\ge 1$, we define an endomorphism $\j$ of the ringed topos $(E_.,A_.)$ by setting $(\j^{-1} F_.)_k=\iota_{k,k+j}^{-1}F_{k+j}$ for $k\ge 1$, $(\j_* F_.)_k$ a final object for $1\le k\le j$ and  $(\j_* F_.)_k=\iota_{k,k+j,*}F_{k-j}$  for $k\ge j+1$), while $\j^{-1}A_.\to A_.$ is the natural morphism.
\end{defn}

We note that there is a natural isomorphism $\j\simeq \jj^j$ as well as a natural morphism $id\to \j$ (given by $(\j^{-1} F_.)_k=\iota_{k,k+j}^{-1}F_{k+j}\to F_k$). In the following definition, we make use of the derived functor $L\j^*$ for unbounded complexes (\cite{KS} 18.6.9).


\begin{defn} \label{defnorm} \debrom \item \label{defnormi} We say that $M_.$ in the category $Mod(E_{.},A_{.})$ of modules  is normalized if $\jj^*M_.\simeq M_.$ (ie. if  $\iota_{k,k+1}^*M_{k+1}\simeq M_{k}$ for all $k\ge 1$).
\item \label{defnormii} We say that $M_.$ in the derived category $D(E_{.},A_{.})$ is $L$-normalized if $L\jj^*M_.\simeq M_.$ (ie. if  $L\iota_{k,k+1}^*M_{k+1}\simeq M_{k}$ for all $k\ge 1$).
\finrom \end{defn}

The following lemma will be used to ``divide the Frobenius'' by $p$ in the definition of syntomic complexes. For $i\ge 0$, we use the notation $(E_{.+i},A_{.+i})$ for the ringed total topos associated to $k\mapsto (E_{k+i},A_{k+i})$, $(k\le k')\mapsto \smash{\iota_{k+i,k'+i}}$. The morphisms $(\smash{E_{k},A_{k}})\to (\smash{E_{k+i},A_{k+i}})$ for varying $k$ and $i$ give rise to natural morphisms denoted as follows: $$\xymatrix{(\smash{E_{k},A_{k}}) \ar @{->} @<+0pt> `d[r] `[rrrr]_-{\iota_{.+i,.+i+j}} [rrrr] \ar[rr]^{\iota_{k,.+i}}&&(\smash{E_{.+i},A_{.+i}})\ar[rr]^{\iota_{.,.+i+j}}
&&(\smash{E_{.+i+j},A_{.+i+j}})}$$

\begin{lem} \label{cartgen} Assume that each one of the morphisms of topoi $\iota_{k,k+1}:E_k\to E_{k+1}$ is an equivalence and that $A_.$ is a flat normalized $\Z/p^{.}$-algebra. Consider modules $L_.$, $M_.$, $N_.$ over $(E_.,A_.)$.

\debrom
\item \label{cartgeni} The module $M_.$ is $L$-normalized if and only if it is normalized and $\Z/p^.$-flat.
\item  \label{cartgenii} For  $j\ge 1$, there is a functorial exact sequence $$\xymatrix{\iota_{.,.+j,*}\j^*M_.\ar[r]^-{{p^j}}& M_{.+j}\ar[r]& \iota_{j,.+j,*}\iota_{j,.+j}^*M_{.+j}\ar[r]& 0}$$
of modules on $(E_{.+j},A_{.+j})$. If $M_.$ is  $\Z/p^{.}$-flat the exact sequence can be extended with a $0$ on the left. If $M_.$ is normalized then
the last term can be replaced by $\iota_{j,.+j,*}M_{j}$. \petit


\item \label{cartgeniii} If $N_.$ is normalized and $M_.$ is $L$-normalized then the group $\smash{Hom(N_.,M_.)}$ of $A_{.}$-linear morphisms is $p$-torsion free. \petit


\item \label{cartgeniv} Consider an exact sequence $\smash{0\to N_.\mathop\rightarrow\limits M_.\to L_.\to 0}$ of modules on $(E_{.},A_{.})$ where $L_.$ is normalized and $M_.$ is $L$-normalized. If $L_.$ is killed by $p^j$ then $\j^*N_.$ and $\tau_{\ge -1}L\j^*L_.$ are $L$-normalized and fit into a canonical distinguished triangle of $D(E_{.},A_{.})$: \begin{eqnarray}\label{dtj1}\xymatrix{\j^*N_.\ar[r] &\j^*M_.\ar[r]& \tau_{\ge -1}L\j^*L_.\ar[r]&\j^*N_.[1]}\end{eqnarray}

\noindent The associated long exact sequence of cohomology is
\begin{eqnarray}\label{exj1}\xymatrix{0\ar[r]&L_.'\ar[r]&\j^*N_.\ar[r]&M_.\ar[r]&L_.\ar[r]&0}\end{eqnarray} where $L'_.$ denote a projective
system where $L'_k=L_k$ if $k\ge j$ and $\smash{\iota_{k,k'}^*L'_{k'}\to L_k}$ is zero if $k'\ge k+j$.
\finrom
\end{lem}
Proof.  The morphisms $\iota_{1,k}:E_1\to E_k$ for varying $k$ induce an equivalence $E_1^\N\simeq E_.$. We may thus assume that $(E_.,A_.)$ is of the form $(E^\N,A_.)$.

\ref{cartgeni} Since $A_.$ is flat and normalized over $\Z/p^.$, $\smash{\cal Tor_q^{A_{k+j}}(A_k,-)}$ and   $\smash{\cal Tor_q^{\Z/p^{k+j}}(\Z/p^k,-)}$ coincide. Using this, we find that $M$ is $L$-normalized if it is  $\Z/p^.$-flat and normalized. The other implication follows from the fact that $M_k$ is $\Z/p^k$-flat if and only if $\smash{\cal Tor_q^{\Z/p^k}(\Z/p,M_k)}$ vanishes for $q\ge 1$.

\ref{cartgenii} The claimed exact sequence in $(E^\N,A_{.+j})$ boils down to $$\xymatrix{M_{.+j}/p^{.}\ar[r]^-{{p^j}}&M_{.+j}\ar[r]&M_{.+j}/p^j\ar[r]&0}$$
which is straightforward from $0\to \Z/p^.\to \Z/p^{.+j}\to \Z/p^j\to 0$.

\ref{cartgeniii} Recall that $(\j^*\j_*M_.)_k\simeq M_k/p^{k-j}$ if $k\ge j+1$ and $0$ otherwise. Multiplication by $p^j$ on $M_.$ thus factors through $\j^*\j_*M_.$. Since $M_.$ is $\Z/p^.$-flat, the resulting morphism $\j^*\j_*M_.\to M_.$ is a monomorphism. We conclude by the following commutative square where the top arrow is an isomorphism because $M_.$ is normalized: $$\xymatrix{Hom(N_.,M_.)\ar[d]^{p^j}\ar[rr]^-{\j^*}&& Hom(\j^*N_.,\j^*M_.)\ar[d]^\wr \\Hom(N_.,M_.)&&\ar[ll]_-{Hom(N_.,p^j)}Hom(N_.,\j_*\j^*M_.)}$$

\ref{cartgeniv} Applying $L\j^*$ to the given exact sequence gives a distinguished triangle $$\xymatrix{L\j^*N_.\ar[r]&L\j^*M_.\ar[r]&L\j^*L_.\ar[r]&\j^*N_.[1]}$$ The claimed distinguished triangle follows by truncation, since the second term is concentrated in degree $0$. Let us prove that $\j^*N_.$ and $\tau_{\ge -1}L\j^*L_.$ are $L$-normalized. By the distinguished triangle, it suffices to prove this for $\tau_{\ge -1}L\j^*L_.$. Since $A_.$ is flat and normalized over $\Z/p^.$, we may furthermore assume that $A_.=\Z/p^.$. Since $L_.$ is normalized and killed by $p^j$, it is of the form $\Z/p^.\otimes_\Z L$ for some $L$ killed by $p^j$. But then \begin{eqnarray}\label{eq1cartgeniv}\tau_{\ge -1}L\j^*L_.&\simeq &\tau_{\ge -1}\Z/p^.\Ltens_{\Z/p^{.+j}}\Z/p^{.+j}\otimes_\Z L\\
 \label{eq2cartgeniv}&\simeq &\tau_{\ge -1}\Z/p^.\Ltens_{\Z/p^{.+j}}L\end{eqnarray} in $D(E^\N,\Z/p^.)$. Consider the natural morphism \begin{eqnarray} \label{eq3cartgeniv}\Z/p^.\Ltens_{\Z}L\to \Z/p^.\Ltens_{\Z/p^{.+j}}L\end{eqnarray}
in $D(E^\N,\Z/p^.)$. We claim that this morphism becomes invertible after $\tau_{\ge -1}$. To check this, it is sufficient to look at the $k$-th component and restrict scalars to $\Z$. Then (\ref{eq3cartgeniv}) gives a morphism of $D(E,\Z)$ which is represented by the following morphism of complexes placed in degrees $]-\infty,0]$: \begin{eqnarray}\label{eq4cartgeniv}[L\mathop\rightarrow\limits^{p^k} L]\to [\dots \to L\mathop\rightarrow\limits^{p^j}  L \mathop\rightarrow\limits^{p^k}L\mathop\rightarrow\limits^{p^j}  L \mathop\rightarrow\limits^{p^k} L]\end{eqnarray}
Our claim thus results from the fact that $p^j$ is zero on $L$. The right hand side in (\ref{eq2cartgeniv}) is thus isomorphic to the left hand side of (\ref{eq3cartgeniv}).  It is in particular $L$-normalized.
The description of $L'_.$ in the long exact sequence of cohomology is obtained by letting $k$ vary in the righthand side of (\ref{eq4cartgeniv}): in  $D(E_.,\Z)$ we have an isomorphism $${\tau_{\ge -1}L\j^*L_.}\simeq [L'_.\mathop\rightarrow\limits^{p^.} L_.]$$
Let us emphasize however that such an isomorphism certainly does not hold in general at the level of $D(E_.,\Z/p^.)$.
\begin{flushright}$\square$\end{flushright}

Assume now given a ringed topos $(E,A)$ together with a morphism of variable ringed topoi on $\N^{op}$: $\iota_k:(E_k,A_k)\to (E,A)$, $k\ge 1$. In this situation we have a morphism $$l:(E_.,A_.)\to (E,A)$$
and the typical example of normalized module (resp. an $L$-normalized complex) is $l^*M$ (resp. $Ll^*M$).

\begin{lem} \label{lemek} Assume that each $\iota_k:E_k\to E$ is an equivalence. Assume that $A$ is $a$ flat $\Zp$-algebra and $A_.\simeq \Z/p^.\otimes l^{-1}A$. The restriction of the functor $$Rl_*:D^+(E_.,A_.)\to D^+(E,A)$$ to the full subcategory of $L$-normalized complexes is fully faithful.
\end{lem}
Proof. Since $Ll^*$ preserves $D^+$ it suffices to prove that the adjunction morphism $Ll^*Rl_*M_.\to M_.$ is invertible for $M_.$ $L$-normalized. As in the proof of Lem. \ref{cartgen}, we may assume that $E_.=E^\N$ and $A=\Zp$. Then we invoke \cite{Ek} Lem. 2.2 (i).
\begin{flushright}$\square$\end{flushright}

\begin{example} \label{exRl} Let $E_.=X_{syn}^\N$, $E=X_{syn}$, $A=\Zp$ and $G\in pdiv(X)$.

\debrom \item The projective system $G_{p^.}$ is $L$-normalized.

\item The isomorphism $Ll^*Rl_*G_{p^.}\simeq G_{p^.}$ means that the sheaves $R^ql_*G_{p^.}$ are uniquely $p$-divisible (but not necessarily trivial) for $q\ge 1$. For instance if $X=\Sigma_1$ then the restriction of $R^1l_*\Z/p^.$ to the small \'etale site is $\Qp$.
\finrom
\end{example}

\subsection{Limits and quasi-coherence on $p$-adic formal schemes} \label{laqcopafs}


\para \label{defl} We adopt the following conventions regarding $p$-adic schemes, $p$-adic fine log schemes and their associated sites. As usual $\Sigma_k$ denotes $Spec(\Z/p^k)$ (see Sect. \ref{sectioncrys}).
\begin{defn} \label{defpadic} \debrom

\item \label{defpadici}   A \emph{$p$-adic scheme} $X$ is an ind object of the category of schemes which is   isomorphic to the inductive system of schemes $\smash{X_.=(X_k)_{k\ge 1}}$ where $X_k$ is the reduction mod $p^k$ of a $p$-adic formal scheme in the sense of \cite{EGA1}. The category of $p$-adic schemes is denoted $\cal Sch_p$. The category $\cal Sch_{p,nil}$ of schemes where $p$ is locally nilpotent is identified with a full subcategory of $\cal Sch_p$ in the natural way.


\item \label{defpadicii} A \emph{$p$-adic log scheme} $X^\s$ is an ind object of the category of log schemes which is isomorphic to an inductive system of log schemes $\smash{X^\s_.=(X^\s_k)_{k\ge 1}}$ satisfying $\smash{X^\s_k\simeq  \Sigma_k\times X^\s_{k+1}}$ and such that the underlying ind-object of the category of schemes is in $\cal Sch_p$. The category of $p$-adic log schemes is denoted $\smash{\cal Sch_p^\s}$. The category $\cal Sch^\s$ of fine log schemes and $\cal Sch_p$ are identified with full subcategories of $\smash{\cal Sch_p^\s}$ in the natural way.
\finrom
\end{defn}

With the notations of \ref{defpadicii}, $X^\s$ is  $\smash{indlim\, X^\s_k}$, ie. the inductive limit of the $\smash{X_k^\s}$'s computed in the category of ind-objects of the category of log schemes, or equivalently, in $\smash{\cal Sch^\s_p}$. It follows moreover from \cite{SGA4-I} I, Prop. 8.9.1 that $\smash{X^\s_k\simeq \Sigma_k\times X^\s}$ where the product is computed in $\smash{\cal Sch^\s_p}$. Given any $X^\s$ in $\smash{\cal Sch_p^\s}$, we will thus use the notation $X^\s_k:=\Sigma_k\times X^\s$ without any danger of confusion.



In virtue of Prop. 10.6.2 and Cor. 10.6.4 of \cite{EGA1} Chap. 0, the category $\cal Sch_p$ is equivalent to the usual category of $p$-adic formal schemes (without any
Noetherian assumption). The category $\cal Sch_p^\s$, on the other hand, is probably too large to be a good category  (see the remark after Def. \ref{defOM} below) but this won't matter for our purpose.

\begin{defn} \label{usualtoppadic} Let $top$ be as in Sect. \ref{usualtop} and $X^\s$ in $\smash{\cal Sch_p^\s}$.
\debrom
\item \label{usualtoppadici} The pretopology $\smash{TOP^\s(X^\s)}$ is $\smash{\cal Sch^\s_p}/X^\s$ endowed with the pretopology for which a covering $\smash{(U_i^\s\to U)}$ is a family whose reduction mod $p^k$ is a covering in $TOP^\s(X^\s)$ for all $k\ge 1$. The associated topos is denoted $\smash{X_{TOP^\s}^\s}$.
    The pretopology $\smash{TOP(X^\s)}$ (resp. $\smash{top(X^\s)}$) is the full subcategory of  $\smash{TOP^\s(X^\s)}$ formed by the $U^\s/X^\s$'s whose reduction mod $p^k$ is in  $\smash{TOP(X_k^\s)}$ (resp. $\smash{top(X_k^\s)}$) for all $k\ge 1$ and where coverings are defined as in $\smash{TOP^\s(X^\s)}$. The associated topos is denoted $\smash{X^\s_{TOP}}$ (resp. $\smash{X_{top}^\s}$).
\item \label{usualtoppadicii} The topos  $\smash{X_{.,TOP^\s}^\s}$ is the $\s$-big $top$ topos of $X_.$ viewed as a diagram of type $\N^{op}$. Explicitly an object of $\smash{X_{.,TOP^\s}^\s}$ is thus a collection $(F_k)$, $F_k\in \smash{X_{k,TOP^\s}^\s}$ together with transition morphisms $F_{k+1}\to \smash{\iota_{k,k+1,*}F_k}$ (here $\iota_{k,k+1}$ denotes the inclusion $\smash{X_k^\s\to X_{k+1}^\s}$). The topoi $\smash{X^\s_{.,TOP}}$ and $\smash{X_{.,top}^\s}$ are defined similarly.
    \finrom
\end{defn}

We have moreover a pseudo-commutative diagram of topoi  \begin{eqnarray}\label{iotalTOP}\xymatrix{X^\s_{k,TOP^\s}\ar[rr]_-{\iota_{k,.}}\ar
@<+2pt> `u[r] `[rrrr]^-{\iota_k} [rrrr]
&&X^\s_{.,TOP^\s}\ar[rr]_-{l}&&X^\s_{TOP^\s}}\end{eqnarray} where the morphisms $\iota_{k,.}$ and $l$ have the following explicit description (as usual the description of transition morphisms is implicitly left to the reader).

\noindent - The functor $\iota_{k,.}^{-1}$ takes $F_.$ to $F_k$. The functor $\iota_{k,.,*}$ takes $F$ to $F'_.$ where $F'_{k'}=\iota_{k,k'_*}F$ if $k'\ge k$ and $F'_{k'}$ is a final object of $\smash{X^\s_{k,TOP^\s}}$   otherwise.
%

\noindent - The functor $l^{-1}$ takes $F$ to $F'_.$ where $F'_k=\iota_{k}^{-1}F$.  The functor $l_*$ takes $F_.$ to $limproj_k \iota_{k*}F_k$.
%
\petit

When $X^\s$ varies, (\ref{iotalTOP}) is a diagram of morphisms between variable topoi on $\smash{\cal Sch^\s_p}$. In the case of big and small topoi, it has a weak analogue. These diagrams are
pseudo-compatible with the weak morphisms of projection between $\s$-big, big and small topoi and with the weak morphisms involving a change of topology.

\begin{defn} \label{defOM} Let $top$ be as in Sect. \ref{usualtop} and $X^\s$ in $\smash{\cal Sch_p^\s}$.
\debrom \item \label{defOMi} The structural ring of $\smash{X^\s_{.,TOP}}$ is defined as $\O:=\O_.$ where $\O_k$
denotes the structural ring of $\smash{X_{k,TOP}^\s}$. The structural ring $\O$ of $\smash{X_{TOP}^\s}$ is defined as $\O:=l_*\O_.$. The structural ring of the small topoi are defined by restriction.
\item \label{defOMii} Consider the morphism of monoids $\smash{M_{X_k}\to \O_k}$ of $\smash{X_{k,et}}$ defining the log structure of $\smash{X^\s_k}$. Letting $k$ vary defines a monoid $\smash{M_{X,.}}$ in $\smash{X_{.,et}^\s}$ and a morphism $\smash{M_{X,.}\to \O}$. The monoid $M_X$ of $\smash{X_{et}^\s}$ is defined as $M_X:=\smash{l_*M_{X,.}}$.
\finrom
\end{defn}

It does not seem clear to us whether or not $M_X\to \O$ is always a fine log structure in the sense of \cite{Sh}. This won't be a problem for us since in practice our $p$-adic log schemes will always have nice charts (see Sect. \ref{parapb}).

\para \label{paraqcohsch} Before discussing quasi-coherence over $p$-adic schemes, let us  recall a possible definition for quasi-coherent modules on the usual sites of a fine log scheme. Log structure play essentially no role in the following discussion. In particular $\s$-big topoi could be replaced by big topoi in Def. \ref{defqcohsch}, Lem.-Def. \ref{lemqcohsch} and Lem. \ref{acycqcohsch} below.

\begin{defn} \label{defqcohsch} Consider $X^\s$ in $\cal Sch^\s$. Let $prop\in \{qcoh,lf,lfft\}$. Let $Mod_{prop}(\smash{X^\s_{zar}},\O)$ denote the category of quasi-coherent modules (resp. locally free, resp. locally free of finite type) on the scheme $X$ if $prop=qcoh$ (resp. $lf$, resp. $lfft$). Consider the projection morphism $\pi$ from the $\s$-big to the small topoi. We define $Mod_{prop}(\smash{X^\s_{ZAR^\s}},\O)$ as the essential image of the functor $$\pi^*:Mod_{prop}(\smash{X^\s_{zar}},\O)\to Mod(\smash{X^\s_{ZAR^\s}},\O)$$
\end{defn}

A reasonable definition of quasi-coherence on the other usual sites requires the following lemma.

\begin{lemdefn} \label{lemqcohsch} Let $top$ be as in Sect. \ref{usualtop} and let $prop$ be as in Def. \ref{defqcohsch}.
\debrom \item \label{lemqcohschi} If $M$ is in $Mod_{qcoh}(\smash{X^\s_{ZAR^\s}},\O)$ then it is in fact a sheaf on $TOP^\s(X^\s)$. We may thus define $$Mod_{prop}(\smash{X^\s_{TOP^\s}},\O):=Mod_{prop}(\smash{X^\s_{ZAR^\s}},\O)$$
\item \label{lemqcohschii} Let $\epsilon$ denote the morphism of change of topology. The essential images of the following functors coincide: $$\begin{array}{l} \pi_*:Mod_{prop}(\smash{X^\s_{TOP^\s}},\O)\to Mod(\smash{X^\s_{top}},\O)\\ \epsilon^*: Mod_{prop}(\smash{X^\s_{zar}},\O)\to Mod(\smash{X^\s_{top}},\O)\end{array}$$
    We define $Mod_{prop}(\smash{X^\s_{top}},\O)$ as the essential image in question.
\item \label{lemqcohschiii} Consider the following pairs of adjoint functors $$\begin{array}{l}(\pi^*,\pi_*):Mod(\smash{X^\s_{TOP^\s}},\O)\to Mod(\smash{X^\s_{top}},\O)\\ (\epsilon^*,\epsilon_*):
    Mod(\smash{X^\s_{top}},\O)\to Mod(\smash{X^\s_{zar}},\O)\end{array}$$
    Each one of these four functors preserves $prop$. The induced adjunctions on $Mod_{prop}$ are equivalences.
\item \label{lemqcohschiv} The contravariant pseudo-functors $Mod_{prop}((-)_{TOP^\s},\O)$ and $Mod_{prop}((-)_{top},\O):\cal Sch^\s\to \mathfrak Cat$ are stacks for $fl$ (ie. satisfy $fl$ descent).
\item \label{lemqcohschv} A module $M$ of $\smash{(X^\s_{TOP^\s},\O)}$ or $\smash{(X^\s_{top},\O)}$ satisfies $prop$ if and only if its restrictions to the elements of a surjective family of $top$ morphisms do.
\finrom
\end{lemdefn}
Proof. Recall the following fact from sheaf theory. Let $\C$ denote a full subcategory of $\cal Sch^\s/X^\s$, such that $U^\s\in \cal C$ and $V^\s/U^\s\in top(U^\s)$ imply $V^\s\in \cal C$. Consider on the one hand the ringed site $((\C,zar),\O)$ (resp. $((\C,top),\O)$) obtained by endowing $\C$ with the $zar$ (resp. $top$) topology and the structural ring. Consider on the other hand $Mod((-)_{zar},\O)$ as a bifibered category over $\cal C^{op}$.
Sending a module $F$ on $((\cal C,zar),\O)$ to the collection $\xi$ of its restrictions $\xi(U^\s):=\smash{F_{|zar(U^\s)}}$ together with the natural base change morphisms $\xi(f):f^*\smash{F_{|zar(U^\s)}} \to \smash{F_{|zar(U'^\s)}}$ for $f:U'^\s\to U^\s$ realizes an equivalence of categories between the modules over $((\cal C,zar),\O)$ and a full subcategory, say $\cal F_{zar}$, of the category of sections of the cofibered category $Mod((-)_{zar},\O)/\cal C^{op}$. Explicitly a section $\xi$ of the latter cofibered category is in $\cal F_{zar}$ if and only if  $\xi(f):f^*\xi(U^\s)\to \xi(U'^\s)$ is invertible for every open immersion  $f$. Under this equivalence the category of sheaves of modules over $((\cal C,top),\O)$ corresponds to a full category, say  $\cal F_{zar,top}$, of $\cal F_{zar}$. Explicitly, a section $\xi$ in $\cal F_{zar}$ is in fact in $\cal F_{zar,top}$ if and only if it satisfies the following property ($top$ descent): if $\smash{(f_i:U_i^\s\to U^\s)}$ is a $top$ surjective family then $$\xi(U^\s)\simeq Ker(\prod_i f_{i,*}\xi(U_i^\s)\rightrightarrows \prod_{j,k}(f_j\times_{U^\s} f_k)_*\xi(U_j^\s\times_{U^\s}U_k^\s))$$

Apply this to the case where $\C$ is either $\cal Sch^\s/X^\s$ or the category underlying $top(X^\s)$. Starting from a module $M$ of $\smash{(X^\s_{zar},\O)}$ the pullback of $M$ to $((\cal C,zar),\O)$ corresponds to the collection $\xi$ of pullbacks $\xi(U^\s)=f^*M$, $f:U^\s\to X^\s$ endowed with the obvious base change (iso)morphisms. If $M$ is quasi-coherent then it follows from descent theory that $\xi$ automatically satisfies the property of $top$ descent. In other terms, the pullback of $M$ to $((\cal C,zar),\O)$ is already a sheaf for $top$, and thus coincides with the pullback of $M$ to $((\cal C,top),\O)$. All statements are straightforward from this remark.
\begin{flushright}$\square$\end{flushright}

\begin{rem} \debrom \item The last part of the proof is related to \cite{SGA4-II} VII, Rem. 2.1, c. The statement in \emph{loc. cit.} (where no log structure appear) would imply in particular that the collection $\xi$ of pullbacks of $M$ satisfies $top$ descent whenever $top$ is coarser than $fl$, for any module $M$ of $\smash{(X^\s_{zar},\O)}$. This is incorrect without the assumption that $M$ is quasi-coherent. The reader can already find counterexamples when $top=et$ and $\C=et(X^\s)$.

\item In the same spirit, let us emphasize that in the previous lemma, \ref{lemqcohschv} is not a formal consequence of \ref{lemqcohschiv} but also relies on the (trivial) fact that $Mod_{prop}((-)_{TOP^\s},\O)$ and $Mod_{prop}((-)_{top},\O)$ are stacks for $top$.\finrom
\end{rem}

\begin{rem} \label{remqcohsch1} \debrom \item \label{remqcohsch1i} The category $Mod_{qcoh}(\smash{X^\s_{top}},\O)$ is Abelian and the inclusion functor into $Mod(\smash{X^\s_{top}},\O)$ is exact. In the case $top=zar$, this follows from Cor. 1.3.9 and Thm. 1.4.1 of \cite{EGA1} Chap. 1. In the general case, this follows from the flatness of the morphism $\epsilon:(\smash{X^\s_{top}},\O)\to (\smash{X^\s_{zar}},\O)$ (which, in turn, results from the interpretation of $\epsilon^*$ in terms of collections of small Zariski sheaves as in the above proof).
\item \label{remqcohsch1ii} Consider a subsheaf $F$ of a module $M$ in $(\smash{X^\s_{top}},\O)$. We may form the submodule $<F>$ of $M$ generated by $F$ (ie. the image in $M$ of the free module with basis $F$). It is not always the case that $<F>$ is quasi-coherent even if $M$ is. This is the case however if $M$ is quasi-coherent and $F$ is a locally constant sheaf, as seen using Lem.-Def. \ref{lemqcohsch} \ref{lemqcohschv}.
\item \label{remqcohsch1iii} Consider a family of sections $(f_i)_{i\in I}\subset \Gamma(X,\O)$. By \ref{remqcohsch1ii}, the $f_i$'s generate a quasi-coherent ideal of $\O$ and thus a closed subscheme of $X$ that we will denote $V_X((f_i)_{i\in I})$ (\cite{EGA1}, Chap. 1, Def. 4.1.3).
\finrom
\end{rem}

\begin{lem} \label{acycqcohsch} Consider a quasi-coherent module $M$ on $(\smash{X^\s_{TOP^\s}},\O)$ (resp. $(\smash{X^\s_{top}},\O)$).
    \debrom \item \label{acycqcohschi} Consider $top'$ between $zar$ and $top$ and let  $\epsilon$ denote the morphism $\smash{X^\s_{TOP^\s}\to X^\s_{TOP'^\s}}$ (resp. $\smash{X^\s_{top}\to X^\s_{top'}}$). Then $M$ is acyclic for $\epsilon_*$.
\item \label{acycqcohschii} Let $f:X^\s\to X'^\s$ denote a morphism in $\cal Sch^\s$ whose underlying morphism of schemes is affine. Then $M$ is acyclic for $f_*$. The module $f_*M$ is moreover quasi-coherent.
\finrom
\end{lem}
Proof. We only explain the case of small sites since the case of $\s$-big sites is similar (or alternatively follows formally). Let $\cal C({X^\s},top)$ denote the full subcategory of $top(X^\s)$ formed by the $U^\s/X^\s$'s with $U$ affine. Endow $\cal C({X^\s},top)$ with the pretopology of surjective families of $top$ morphisms.
The associated topos is thus equivalent to $\smash{X^\s_{top}}$ (see \cite{Ar} Lem. 3.1.3). It follows from faithfully flat descent theory that the Cech cohomology of a quasi-coherent module for an arbitrary covering in $\cal C({X^\s},top)$ is trivial in degree $\ge 1$. In other terms, quasi-coherent modules are $\cal C({X^\s},top)$-acyclic. Statement \ref{acycqcohschi} (resp. the first statement of \ref{acycqcohschii}) thus follows from the fact that $\epsilon$ (resp. $f$) arises from a premorphism $\cal C({X^\s},top)\to \cal C({X^\s},top')$ (resp. $\cal C({X^\s},top)\to \cal C({X'^\s},top')$). The second statement of \ref{acycqcohschii} follows from Prop. 1.2.6 and Cor. 1.5.2 of \cite{EGA2},  together with the characterization explained in the proof of Lem.-Def. \ref{lemqcohsch} for quasi-coherent modules as cocartesian sections of the cofibered category of quasi-coherent modules on the small Zariski site of a variable base.
\begin{flushright}$\square$\end{flushright}



\para \label{paraqcohpadic} We discuss briefly a first candidate for the notions of quasi-coherence on a $p$-adic log scheme $X^\s$. Logarithmic structures are left aside from the discussion since $\smash{X^\s_{top}}\simeq X_{top}$ and $\smash{X_{.,top}^\s\simeq X_{.,top}}$. We restrict furthermore to the \'etale topology which has the advantage that the morphisms $\iota_{k,k+1}:X_{k,et}\to X_{k+1,et}$ and $\iota_k:X_{k,et}\to X_{et}$ are equivalences.

\begin{defn} Consider $X$ in $\cal Sch_p$.
\debrom \label{defqcohpadic}
\item \label{defqcohpadici} We say that a module $M_.$ of $(X_{.,et},\O)$ is \emph{quasi-coherent} if each $M_k$ is quasi-coherent in the sense of Lem.-Def. \ref{lemqcohsch}, \ref{lemqcohschii}.
\item \label{defqcohpadicii} Let $M$ denote a module on $(X_{et},\O)$. We say that $M$ is \emph{quasi-coherent} if $l^*M$ is quasi-coherent in the sense of \ref{defqcohpadici} and $M\to l_*l^*M$ is an isomorphism.
\finrom
\end{defn}

\begin{rem} \label{remqcohpadic} By definition, a $p$-adic scheme admits a covering by open sub $p$-adic schemes which are affine. Since quasi-coherent modules are acyclic on affine schemes, the Mittag-Leffler criterion of \cite{BO} Lem. 7.20 ensures the following:
\debrom \item \label{remqcohpadici} The natural morphism $\iota_k^{-1}\O/p^k\to \O_k$ is invertible.
\item \label{remqcohpadicii} More generally if $M_.$ is a normalized quasi-coherent module of $(X_{.,et},\O)$ then $M_k\simeq \smash{\iota_k^*l_*M_.}\simeq \smash{\iota_k^{-1}l_*M_./p^k}$.
\item \label{remqcohpadiciii} The functor $l_*$ is in particular fully faithful on the category of normalized quasi-coherent modules. Whence a tautological equivalence $$\xymatrix{Mod_{norm,qcoh}(X_{.,et},\O)\ar[rr]_-\sim^-{(l^*,l_*)}&&Mod_{qcoh}(X_{et},\O)}$$
\finrom
\end{rem}

We say that $X$ is flat over $\Sigma_\infty$  if each $X_k$ is flat over $\Sigma_k$. We say that $X\to Y$ is a closed immersion if each $X_k\to Y_k$ is a closed immersion.

\begin{lem} \label{lemflatpadic} Consider $X$ in $\cal Sch_p$.
\debrom \item \label{lemflatpadici} The structural ring $\O$ of $X_{et}$ is $\Zp$-flat if and only if $X$ is flat over $\Sigma_\infty$.
\item \label{lemflatpadicii} Assume that $X$ is flat over $\Sigma_\infty$ and consider a closed immersion $Y\to X$. Then $Y$ is flat over $\Sigma_\infty$ if and only if the corresponding quasi-coherent ideal $I_.$ of $(X_{.,et},\O)$ is normalized.
\finrom
\end{lem}
Proof. \ref{lemflatpadici} follows from \cite{MW} Lem. 2.1 and \ref{lemflatpadicii} is straightforward.
\begin{flushright}$\square$\end{flushright}

\para \label{Lqcohpadic} A drawback of the category of quasi-coherent modules in $(X_{et},\O)$ as defined in Def. \ref{defqcohpadic} is that it fails to be Abelian. This may justify the following alternative definition.

\begin{defn} \label{defLqcohpadic} Consider a flat $p$-adic scheme $X$ over $\Sigma_\infty$.
\debrom \item  \label{defLqcohpadici} Consider $M_.$ in $D^+(X_{.,et},\O)$. We say that $M_.$ is \emph{quasi-coherent} if its cohomology modules are quasi-coherent in the sense of Def. \ref{defqcohpadic} \ref{defqcohpadici}.
\item \label{defLqcohpadicii} Consider $M$ in $D^+(X_{et},\O)$. We say that $M$ is \emph{$L$-quasi-coherent} if $Ll^*M$ is quasi-coherent in the sense of \ref{defLqcohpadici} and $M\to Rl_*Ll^*M$ is an isomorphism.
\finrom
\end{defn}

\begin{rem} \label{remLqcohpadic} In virtue of Rem. \ref{remqcohpadic} \ref{remqcohpadici}, the morphism $l:(X_{.,et},\O)\to (X_{et},\O)$ satisfies the assumptions of Lem. \ref{lemek}. Whence a tautological equivalence $$\xymatrix{D^+_{Lnorm,qcoh}(X_{.,et},\O)\ar[rr]_-\sim^-{(Ll^*,Rl_*)}&&D^+_{Lqcoh}(X_{et},\O)}$$
where the subscripts on the left (resp. the subscript on the right) derived category refer to  Def. \ref{defnorm} \ref{defnormii}  and Def. \ref{defLqcohpadic} \ref{defLqcohpadici} (resp. Def. \ref{defLqcohpadic} \ref{defLqcohpadicii}).
\end{rem}

A drawback of the category $D^+_{Lqcoh}(X_{et},\O)$ is that it is not stable by truncation. This issue might probably be resolved as in \cite{Ek}. We do not pursue this however since in practice the result of Lem. \ref{cartgen} will be enough for our purpose.

\para \label{critqcoh} The notions of quasi-coherence (Def. \ref{defqcohpadic}) and $L$-quasi-coherence (Def. \ref{defLqcohpadic}) do not seem to be simply related. We will content ourselves with the following lemma.

\begin{lem} \label{lemqcohLqcoh} Assume that $X$ is flat over $\Sigma_\infty$. Consider a module $M$ of $(X_{et},\O)$ which is such that each $\iota_k^*M$ is
quasi-coherent.

\debrom \item \label{lemqcohLqcohi} Assume that $M$ is $\Zp$-flat. Then $M$ is quasi-coherent  if and only it is $L$-quasi-coherent.

\item \label{lemqcohLqcohii} If $M$ is killed by $p^k$ for some $k\ge 0$ then it is quasi-coherent and $L$-quasi-coherent.

\finrom
\end{lem}
Proof. Let us denote $M_.=l^*M$.

\ref{lemqcohLqcohi} Since $M$ is $\Zp$-flat, the module $M_.$ is both normalized and $L$-normalized. Thus it suffices to observe that $l_*M_.\simeq Rl_*M_.$ by Mittag-Leffler (\cite{BO} Lem. 7.20).

\ref{lemqcohLqcohii} We have to prove that the morphisms $M\to l_*l^*M$ and $M\to Rl_*Ll^*M$ are invertible. This is clear for the first one. To deduce it for the second one, it suffices to check that $l_*M_.\simeq Rl_*M$ and $Rl_*\tau_{\le -1}Ll^*M=0$. The first condition follows from the fact that $M_.$ is constant from rank $k$ (ie. $\smash{\iota_{k,k',*}M_k\simeq M_{k'}}$ for $k'\ge k$). The second condition follows by \cite{Ek} Lem. 1.1, thanks to the fact that $\tau_{\le -1}Ll^*M=N_.[1]$ where $N_.$ is essentially zero (ie. $\smash{\iota_{k',k''}^*N_{k''}\to N_{k'}}$ is zero for $k''>>k'$) as seen using the natural $\Zp$-flat resolution of $\Z/p^.$.
 \begin{flushright}$\square$\end{flushright}

\subsection{Limits and quasi-coherence on crystalline sites} \label{laqcocs}

\para \label{paraprocrystop} Consider $X^\s$ in $\cal Sch^\s/\Sigma_1$. Letting $k$ vary in Def. \ref{defcrystop} defines a variable topos on $\N$: $k\mapsto (X^\s/\Sigma_k)_{CRYS^\s,top}$,  $(k'\ge k)\mapsto \iota_{k,k'}:(X^\s/\Sigma_k)_{CRYS^\s,top}\to (X^\s/\Sigma_{k'})_{CRYS^\s,top}$. The same is true for big and small crystalline $top$ topoi.

\begin{defn} \label{defprocrystop} The topos $(X^\s/\Sigma_.)_{CRYS^\s,top}$ is the total topos associated to the above variable topos. Explicitly an object is thus a collection $(F_k)$, $F_k\in  \smash{(X^\s/\Sigma_.)_{CRYS^\s,top}}$, together with transition morphisms $F_{k+1}\to \iota_{k,k+1,*}F_k$. The topoi
$\smash{(X^\s/\Sigma_.)_{CRYS,top}}$ and $\smash{(X^\s/\Sigma_.)_{crys,top}}$ are defined similarly.
\end{defn}

The functoriality properties of this topos with respect to $X^\s$ and $top$ are similar to the case where $k$ is fixed. We have moreover a pseudo-commutative diagram of topoi

\begin{eqnarray}\label{diagiotalcrys} &&{\xymatrix{(X^\s/\Sigma_k)_{CRYS^\s,top}\ar[rr]_-{\iota_{k,.}}\ar @<+2pt> `u[r] `[rrrr]^-{\iota_k} [rrrr]
&&(X^\s/\Sigma_.)_{CRYS^\s,top}\ar[rr]_-{l}&&(X^\s/\Sigma_\infty)_{CRYS^\s,top}}}\end{eqnarray}
The formulae describing  (\ref{iotalTOP}) describe (\ref{diagiotalcrys}) as well.
When $X^\s$ varies, (\ref{diagiotalcrys}) is a diagram of morphisms between variable topoi on $\cal Sch^\s$. In the case of big and small topoi it has a weak analogue. These diagrams are pseudo-compatible with the weak morphisms of projection between $\s$-big, big and small topoi and with the weak morphisms of changing the topology. \petit


The morphisms $i$, $u$ of (\ref{defiu}) naturally induce morphisms
\begin{eqnarray}\label{defproiu} \xymatrix{X^{\s,\N}_{TOP^\s}\ar[rr]_-i&&(X^\s/\Sigma_.)_{CRYS^\s,top}\ar @{-->} @/_2pc/ [ll]_-u}\end{eqnarray} and similarly for big and small topoi. The functoriality properties of these morphisms are similar to the case where $k$ is fixed. \petit

Assume now given $T^\s$ in $\smash{\cal Sch^\s_p}$ as well as an exact divided power immersion
$(U^\s,T^\s,\iota,\gamma)$.
In $\smash{(X^\s/\Sigma_k)_{CRYS^\s,top}}$ we have natural isomorphisms $\smash{(U^\s,T^\s_k)}\rightarrow \smash{\iota_{k,k+1}^{-1} (U^\s, T^\s_{k+1})}$. Let us simply denote $\smash{T^\s_.}$ the object $\smash{(X^\s/\Sigma_.)_{CRYS^\s,top}}$ obtained by inverting  these isomorphisms. Then  (\ref{deflambda}) naturally induces a morphism $\smash{f_{T_.^\s}}$ and a weak morphism ${\lambda_{T_.^\s}}$ as follows:
\begin{eqnarray}\label{defprolambda} \xymatrix{(X^\s/\Sigma_.)_{CRYS^\s,top}&&
\ar[ll]_-{f_{T_.^\s}}(X^\s/\Sigma_.)_{CRYS^\s,top}/T_.^\s\ar[rr]^-{\lambda_{T_.^\s}} &&\smash{T^\s_{.,top}}} \end{eqnarray}
Generalizing Def. \ref{defreal}, we define the \emph{restriction $\smash{F_{.,|T^\s_.}:=f_{T_.^\s}^{-1}F_.}$ of $F_.$ to $\smash{T^\s_.}$} and the \emph{realization   $\smash{F_{.,T^\s_.}:=\lambda_{T_.^\s,*}F_{|T^\s_.}}$ of $F_.$ on $\smash{T^\s_.}$}. Similar observations and notations hold with $CRYS$ or $crys$ instead of $CRYS^\s$.

\para \label{paralimcrys} We discuss limits in crystalline topoi. The relevant morphism of ringed topoi is $$l:(\smash{(X^\s/\Sigma_.)_{CRYS^\s,top}},\O)\to (\smash{(X^\s/\Sigma_\infty)_{CRYS^\s,top}},\O)$$
where both $\O$'s denote the structural ring, ie. $\O:=l^{-1}\O$. The following lemma shows that the situation is much simpler here than in Sect. \ref{paraqcohpadic} and Sect. \ref{critqcoh}.

\begin{lemdefn} Consider $X^\s$ in $\cal Sch^\s$. \label{lemlimcrys}
\debrom \item \label{lemlimcrysi} Let us say that $M_.$ in $\smash{Mod((X^\s/\Sigma_.)_{CRYS^\s,top},\O)}$ is \emph{normalized} if $\smash{\iota_{k,k+1}^{-1}M_{k+1}\simeq M_k}$ for all $k\ge 1$. The morphism $l$ induces an equivalence of categories
$$\xymatrix{Mod_{norm}((X^\s/\Sigma_.)_{CRYS^\s,top},\O)\ar[rr]_-\sim^-{(l^{-1},l_*)}&&
Mod((X^\s/\Sigma_\infty)_{CRYS^\s,top},\O)}$$
\item \label{lemlimcrysii} If $M_.$ is a normalized module of  $((X^\s/\Sigma_.)_{CRYS^\s,top},\O)$ then it is acyclic for $l_*$.
\item \label{lemlimcrysiii} Statements \ref{lemlimcrysi} and \ref{lemlimcrysii} hold verbatim if $CRYS^\s$ is replaced by $CRYS$ or $crys$.
\finrom
\end{lemdefn}
Proof. \ref{lemlimcrysi}  For $M_.$ in $Mod((X^\s/\Sigma_.)_{CRYS^\s,top},\O)$ and $(U^\s,T^\s)$ in $CRYS^\s(X^\s/\Sigma_\infty)$ we have $$l_*M_.(U^\s,T^\s)\simeq \limp_{k}M_k(U^\s,T_k^\s)$$
Consider a module $M$ of $((X^\s/\Sigma_\infty)_{CRYS^\s,top},\O)$. The evaluation of the adjunction morphism $M\to l_*l^{-1}M$ at $(U^\s,T^\s)$ thus reads $M(U^\s,T^\s)\to limproj_k \smash{M(U^\s,T_k^\s)}$. It is clearly an isomorphism since $p$ is locally (and hence may be assumed) nilpotent on $T^\s$. Consider now a module $M_.$ of $((X^\s/\Sigma_.)_{CRYS^\s,top},\O)$. Consider a fixed integer $k_0$. The  $k_0$-th component (ie. the image by $\smash{\iota_{k_0,.}^{-1}}$) of the adjunction morphism $l^{-1}l_*M_.\to M_.$ evaluated at some $(U^\s,T^\s)$ in $\smash{CRYS^\s(X^\s/\Sigma_{k_0})}$ reads $limproj_k M_k(U^\s,T^\s_k)\to M_{k_0}(U^\s,T^\s)$. If $M_.$ is normalized then it is invertible as well since the projective system on the left stabilizes for $k\ge k_0$.


\ref{lemlimcrysii} We have to show that $\smash{(R^ql_*M_.)_{T^\s}}=0$ for every $q\ge 1$ and  $T^\s=(U^\s,T^\s)$ in
$CRYS_{top}(X^\s/\Sigma_\infty)$. Localizing if necessary, we can assume that $p$ is nilpotent on $T$, say $\smash{p^{k_0}}=0$.
First we notice that
$\smash{(-)_{T^\s}} \circ Rl_*$ coincides with the derived functor of $\smash{(-)_{T^\s}} \circ l_*$. Now $\smash{(-)_{T^\s}} \circ l_*\simeq l_*\circ \smash{(-)_{T_.^\s}}$ where $\smash{(-)_{T_.^\s}}=\smash{\lambda_{T_{.}^\s,*}}\smash{f_{T_.}^{-1}}$ as in (\ref{defprolambda}). Next we
claim that $R(l_*\circ \smash{(-)_{T_.^\s}})\simeq Rl_*\circ R\smash{(-)_{T_.^\s}}$. This claim simply results from Lem. \ref{DfuncT} \ref{DfuncTii} applied to the following prevariable pretopologies on $\N^{op}\cup \{\infty\}$: $\cal P:k\mapsto \smash{CRYS^\s_{top}(X^\s/\Sigma_k)}$, $\cal P':k\mapsto \smash{top(T^\s_k)}$, the premorphism of pretopologies $g:\cal P\to \cal P'$ given by $T'^\s/T_k^\s\mapsto (U^\s\times_{T^\s_k}T'^\s,T'^\s)$ (note that $g_{\N^{op},*}=\smash{(-)_{T_.^\s}}$) and the morphism of diagrams $f:\N^{op}\to \{\infty\}$ (note that $f_*=l_*$). Now $\smash{(-)_{T_.^\s}}$ is exact (Lem. \ref{lemreal}\ref{lemrealii}) and we are thus finally reduced to check that $R^ql_*\smash{M_{T_.^\s}}$ vanishes for $q\ge 1$. We conclude by observing that the transition morphisms $$M_{T^\s_k+1}\to \iota_{k,k+1,*}M_{T^\s_k}$$ are invertible for $k\ge k_0$.
\ref{lemlimcrysiii} The proofs of \ref{lemlimcrysi} and \ref{lemlimcrysii} work as well for $CRYS$ or $crys$.
\begin{flushright}$\square$\end{flushright}

\begin{rem} \label{remlcrys}  The statements of Lem.-Def. \ref{lemlimcrys} would hold verbatim for modules over an arbitrary ring.
\end{rem}

\begin{rem}\label{lemlcrys}

\debrom

\item  \label{lemlcrysi} The functors $l_*$ and $\smash{\iota_k^{-1}}$ clearly preserve the crystal condition of Lem. \ref{lemcrystal1}. The equivalence of Lem.-Def.  \ref{lemlimcrys} \ref{lemlimcrysi} thus induces an equivalence (with obvious notations) $$\xymatrix{\cal Crys_{norm}((X^\s/\Sigma_.)_{CRYS^\s,top},\O)\ar[rr]_-\sim^-{(l^{-1},l_*)}&&
\cal Crys((X^\s/\Sigma_\infty)_{CRYS^\s,top},\O)}$$
\item  \label{lemlcrysii} Let $prop$ be as in Lem.-Def. \ref{lemqcohsch}. We say that a module $M$ has $prop$ realizations if each $M_{T^\s}$ is $prop$. This condition is clearly preserved by $l_*$ and $l^{-1}$. Whence (with obvious notations):  $$\xymatrix{Mod_{norm,prop}((X^\s/\Sigma_.)_{CRYS^\s,top},\O)\ar[rr]_-\sim^-{(l^{-1},l_*)}&&
    Mod_{prop}((X^\s/\Sigma_\infty)_{CRYS^\s,top},\O)}$$
\item \label{lemlcrysiii} The above statements \ref{lemlcrysi} and \ref{lemlcrysii} hold as well if $CRYS^\s$ is replaced by $CRYS$ or $crys$. 

\item \label{lemlcrysiv} Let $\Sigma=\Sigma_.$ or $\Sigma_k$, $k\le \infty$ and let $\pi$ denote the
    weak morphism of projection from the $\s$-big crystalline topos to the small crystalline topos.  The adjunction  $(\pi^*,\pi_*)$ for $\O$-modules induces an equivalence $$\xymatrix{\cal Crys((X^\s/\Sigma)_{CRYS^\s,top},\O)\ar[rr]_-\sim^-{(\pi^*,\pi_*)}&&
\cal Crys((X^\s/\Sigma)_{crys,top},\O)}$$
This equivalence clearly preserves the property of $prop$ realizations (we will sometimes use the terminology \emph{$prop$ crystal} instead of \emph{crystal with $prop$ realizations}). In the case $\Sigma=\Sigma_.$, the property $norm$ is preserved as well. A similar statement holds with $CRYS^\s$ or $crys$ replaced by $CRYS$.
\finrom
\end{rem}


\begin{lem} \label{crysepsilon} Let $\Sigma$ denote either  $\Sigma_.$ or  $\Sigma_k$, $k\le \infty$  and consider the usual $$\epsilon:((X^\s/\Sigma)_{CRYS^\s,top'},\O)\to ((X^\s/\Sigma)_{CRYS^\s,top},\O)$$ for some $top'\preceq top$.
\debrom \item \label{crysepsiloni} The adjunction $(\epsilon^*,\epsilon_*)$ for $\O$-modules induces an equivalence     $$\xymatrix{\cal Crys_{qcoh}((X^\s/\Sigma)_{CRYS^\s,top'},\O)\ar[rr]_-\sim^-{(\epsilon^*,\epsilon_*)}&&
\cal Crys_{qcoh}((X^\s/\Sigma)_{CRYS^\s,top},\O)}$$
\item \label{crysepsilonii} Modules with quasi-coherent realizations are acyclic for $\epsilon_*$.
\item \label{crysepsiloniii} Statements \ref{crysepsiloni} and \ref{crysepsilonii} hold verbatim if $CRYS^\s$ is replaced by $CRYS$ or $crys$.
\finrom
\end{lem}
Proof. \ref{crysepsiloni} We may assume that $\Sigma=\Sigma_k$. Recall that $\epsilon_*$ is fully faithful. We need to check the following: a. (resp. b.) the functor $\epsilon_*$ (resp. $\epsilon^*$) preserves the condition of being a crystal with quasi-coherent realizations and c. if $M$ is a such a crystal then $M\simeq \epsilon_*\epsilon^*M$. By abstract nonsense, it is in fact sufficient to prove a. for arbitrary $top$ and $top'$ and b., c. for $top=zar$.

Let us prove a. Consider $M$ in $\cal Crys_{qcoh}((X^\s/\Sigma)_{CRYS^\s,top'})$ and let us check that $\epsilon_*M$ is a crystal with quasi-coherent realizations. Since $\epsilon_*$ pseudo-commutes to the realization functors, we only have to check the crystal condition, ie. that for all $h:\smash{T_1^\s\to T_2^\s}$, $\smash{h^*(\epsilon_*M)_{T_1^\s}}\to  \smash{(\epsilon_*M)_{T_2^\s}}$ is invertible. Since both the source and the target are quasi-coherent, this is equivalent to $\smash{\epsilon^*h^*(\epsilon_*M)_{T_1^\s}}\to  \smash{\epsilon^*(\epsilon_*M)_{T_2^\s}}$ being invertible. Now the latter morphism identifies with $\smash{h^*M_{T_1^\s}\to M_{T_2^\s}}$ (recall that by quasi-coherence $\epsilon^*\epsilon_*M_{T_i^\s}\simeq M_{T_i^\s}$). It is thus invertible indeed by the crystal condition for $M$.

Let us prove c. in the case $top=zar$. Consider  $M$ in $\cal Crys_{qcoh}((X^\s/\Sigma)_{CRYS^\s,zar})$. The corresponding section $\xi$ of $Mod((-)_{top'},\O)$ over $CRYS^\s_{zar}(X^\s/\Sigma)$ satisfies that $\xi(f)$ is invertible for any $top$ cartesian $f$ (and in fact for any $f$). Now it follows easily from this and \cite{SGA1} VIII, Thm. 1.1, that the variant of the descent condition of Lem. \ref{lemreal} \ref{lemrealiv} as suggested in Rem. \ref{remreal} \ref{remrealiv} automatically holds, ie. that $M$ is in fact a sheaf for $top'$.
In other terms, $M\to \epsilon_*\epsilon^*M$ is invertible.

Let us now prove b. in the case $top=zar$. Consider  $M$ in $\cal Crys_{qcoh}((X^\s/\Sigma)_{CRYS^\s,zar})$ and let us check that $\smash{h^*(\epsilon^*M)_{T_1^\s}}\to  \smash{(\epsilon^*M)_{T_2^\s}}$ is invertible. This will result from the crystal condition for $M$ if we show that $\smash{\epsilon^*(M_{T_i^\s})\to (\epsilon^*M)_{T_i^\s}}$ is invertible. Let us interpret the latter morphism in the category $\cal F_{zar,top'}$ defined during the proof of Lem.-Def. \ref{lemqcohsch} for $\cal C=top'(\smash{T_i^\s})$. By c. above, the source of the latter morphism identifies with $\smash{(f:T'^\s\to T_i^\s)}\mapsto \smash{M_{T'^\s}}$ whereas the target identifies with $\smash{(f:T'^\s\to T_i^\s)}\mapsto f^*\smash{M_{T_i^\s}}$ by quasi-coherence of $\smash{M_{T_i^\s}}$. We may thus conclude by the crystal condition for $M$.

The cases $CRYS$ and $crys$ follow formally from the case $CRYS^\s$ by Rem. \ref{lemlcrys} \ref{lemlcrysiv}.

\ref{crysepsilonii} Thanks to Lem. \ref{DfuncT} \ref{DfuncTi},  it suffices to treat the case $\Sigma=\Sigma_k$ here as well. Now since $\smash{(-)_{T^\s}}$ and $\epsilon_*$ are induced by commuting continuous functors we have the following isomorphisms: \begin{eqnarray}(R\epsilon_*M)_{T^\s}&\simeq& R((-)_{T^\s}\circ \epsilon_*) M\\ &\simeq &  R(\epsilon_*\circ (-)_{T^\s}) M\\
&\simeq &  R\epsilon_*\circ M_{T^\s}\end{eqnarray}
We conclude by Lem. \ref{acycqcohsch} \ref{acycqcohschi}. The case  $CRYS$ or $crys$ is similar (or alternatively follows from the case $CRYS^\s$).
\begin{flushright}$\square$\end{flushright}

\section{Preliminaries part III: crystals and local finite $p$-bases}\label{secppiii}

\subsection{Local $p$-bases and local embeddings} \label{pbale}

\para \label{parapb} We recall the notion of $p$-bases for morphisms of schemes introduced by Kato in \cite{Ka3} as well as a simplified notion of $p$-bases for morphisms of ($p$-adic) log schemes.

\begin{defn} \label{defpb} \debrom

\item \label{defpbi} A morphism $f:X\to Y$ between schemes of characteristic $p$ is \emph{relatively perfect} if the relative Frobenius morphism $$F^{(X/Y)}:X\to X^{(p/Y)}$$ is invertible.

\item \label{defpbii} A morphism $f:X\to Y$ between $p$-adic  schemes is \emph{relatively perfect} if each $f_k:X_k\to Y_k$ is formally \'etale and if $f_1:X_1\to Y_1$ is relatively perfect in the sense of \ref{defpbi}.

\item \label{defpbiii} A \emph{finite $p$-basis} (of cardinal $d$) for a morphism $f:X\to Y$ of $p$-adic  schemes is a $d$-tuple $\underline s\in \G(X,\O_X)^d$ such that the induced morphism $$(\underline s,f):X_k\to {\mathbb A}_{\Sigma_k}^d\times Y_k$$ to the affine space (of dimension $d$) over $Y$ is relatively perfect for each $k$.

\item \label{defpbiv} A \emph{finite $p$-basis} (of cardinal $(d,e)$) for a morphism of $p$-adic log schemes $f:X^\s\to Y^\s$ is a couple $(\underline s,\underline t)$ where $\underline s\in \G(X,\O_X)^d$, $\underline t\in \G(X,M_X)^e$ such that the induced morphism $$(\underline s,\underline t,f):\smash{X_k^\s}\to {\mathbb A}_{\Sigma_k}^{d}\times ({\mathbb A}_{\Sigma_k}^e,\N^e)\times Y^\s_k$$ to the affine space of dimension $d+e$ over $Y_k$, with log structure induced by $M_{Y_k}$ and the canonical one on $\smash{{\mathbb A}^e_{\Sigma_k}}$, is strict and relatively perfect for each $k$.

\item \label{defpbv} A morphism of $p$-adic log schemes $f:X^\s\to Y^\s$ has \emph{local finite $p$-bases} if there exist strict \'etale coverings $\smash{(Y_{i}^\s\to Y^\s)}$, $\smash{(X_{ij}^\s\to X^\s\times_{X^\s}Y_i^\s)}$  such that each $\smash{X_{ij}^\s}$ has finite  $p$-bases over $Y_i^\s$.
\finrom
\end{defn}


\begin{rem} In \cite{Ts2} Def. 1.4, the author defines a $p$-basis for a morphism $f:X^\s\to Y^\s$  of fine $\Sigma_k$-log schemes as a set of elements $(b_\lambda)_\lambda$ in $\G(X,\O_X)$ together with a chart $ch$ for $f$ having certain properties. If  $(\underline s,\underline t)$ is a finite $p$-basis in the sense of Def. \ref{defpb} \ref{defpbiv}, then $\underline t$ can be viewed as a chart $ch$ for the morphism $f$ and the couple $((s_i)_{i},ch)$ is a $p$-basis in the sense of \emph{loc. cit}. Our definition is more restrictive however since we only consider log structures of the type $\N^e$. The reader may consult \cite{CV} for less restrictive notions of local finite $p$-bases. 
\end{rem}

Let us gather some known properties which will be used freely in the text. We begin with facts about relative perfectness.

\begin{lem} \label{lemrp}
\debrom
\item \label{lemrpi} If a morphism of $p$-adic schemes is \'etale (ie. if its reductions $mod\, p^k$ are \'etale) then it is relatively perfect.

\item \label{lemrpii} Relative perfectness as defined in Def. \ref{defpb} \ref{defpbi},  \ref{defpbii} is a Zariski local notion on the source. It is moreover stable by base change and composition.


\item \label{lemrpiii} Consider a relatively perfect morphism $X\to Y$ over $\Sigma_1$. If $Spec(A)$ is an affine open subscheme of $X$, $I$ is a finitely generated ideal in $A$ and $\smash{\hat A}$ is the $I$-adic completion of $A$, then $Spec(\smash{\hat A})$ is relatively  perfect over $Y$.

\item \label{lemrpiv} Consider a commutative square in $\cal Sch_p$: $$\xymatrix{U\ar[r]\ar[d]_i&X\ar[d]^f\\ T\ar @{-->}[ru]^h\ar[r]&Y}$$ If $f$ is relatively perfect and if $i$ is a nilimmersion of order $n$ (ie. if each $U_k\to T_k$ is a closed immersion defined by an ideal whose sections on every open are nilpotent of order $n$) then there exists a unique arrow $h$ making the full diagram commutative.

\item \label{lemrpv} Let $X$ in $\cal Sch_p$. The category of relatively perfect $p$-adic schemes over $X$ and $X_1$ are naturally equivalent. If $Y/X$ is relatively perfect then $Y_k/X_k$ is flat if and only if $Y_1/X_1$ is flat.

\item \label{lemrpvi} Consider a relatively perfect morphism $f:X\to Y$ over $\Sigma_1$. If $Y$ is regular then $f$ is flat. If $Y$ is regular and $X$ is locally Noetherian then $X$ is regular.
\finrom
\end{lem}

Proof. Statement \ref{lemrpi} is proven in \cite{Ka1} Lem. 1.3 and \ref{lemrpii} is elementary.  Statement \ref{lemrpiii} is a particular case of \cite{dJ1} Lem. 1.1.3. Let us briefly indicate the proof of \ref{lemrpiv}. By formal \'etaleness of $X_k/Y_k$, we are easily reduced to the case where $Y=Y_1$. But
then, the result follows from the fact that for $p^r\ge n$, $U^{(p^r/T)}\simeq T$ while $X\simeq X^{(p^r/Y)}$. In \ref{lemrpv} both statements follow from \cite{Ka} Lem. 1 together with \ref{lemrpiv}. The first statement of \ref{lemrpvi} is a particular case of \cite{Ka4} Prop. 1.5 (whose proof refers to \cite{Ka1} Prop. 5.2). The second statement is a consequence of formal smoothness over $\Sigma_1$ (\cite{EGA4-I} Chap. 0, Thm. 22.5.8).
\begin{flushright}$\square$\end{flushright}

Now, some facts about finite $p$-bases.

\begin{lem} \label{lempb} \debrom
\item \label{lempbi} Consider a log scheme of the form $X^\s=(X,Z)$ where $X$ is smooth over $\Sigma_1$ and $Z$ is a normal crossing divisor  (\cite{Ka2} (1.5) (1), (2.5)). Then both $X$ and $X^\s$ have local finite $p$-bases over $\Sigma_1$.

\item \label{lempbii} Morphisms of $p$-adic log schemes having finite $p$-bases or local finite $p$-bases are stable by composition and base change.

\item \label{lempbiii} Consider a commutative square in $\cal Sch^\s_p$
 $$\xymatrix{U^\s\ar[r]\ar[d]_i&X^\s\ar[d]^f\\ T^\s\ar @{-->}[ru]^h\ar[r]&Y^\s}$$
 If $f$ has local finite $p$-bases and if $i$ is an exact nilimmersion of order $n$ (ie. its reductions are exact nilimmersions of order $n$) then after replacing $T^\s$ by a strict \'etale covering one can find $h$ rendering the diagram commutative.

\item \label{lempbiv} Consider a $p$-adic log scheme $Y^\s$ and a morphism $X^\s\to Y^\s_1$ with $p$-basis $(\underline s,\underline t)$. Up to canonical isomorphism, there exists a unique triple $(\smash{\tilde X^\s},\tilde f:\smash{\tilde X^\s}\to Y^\s, (\tilde{\underline s},\tilde{\underline t}))$ where $\tilde f:\smash{\tilde X^\s}\to Y^\s$ is a lifting of $f$ and $(\tilde {\underline s},\tilde{\underline t})$ is a $p$-basis for $\tilde f$ lifting $(\underline s,\underline t)$.

\item \label{lempbv} If $(\underline s,\underline t)$ is a $p$-basis for a morphism of $p$-adic log schemes $f:X^\s\to Y^\s$ then the module of logarithmic differentials is free: \begin{eqnarray}\label{isopbdiff}\Omega_{X^\s/Y^\s}:=l_*\Omega_{X^\s_./Y^\s_.}=\left(\oplus_{i=1}^d \O ds_i\right)\oplus \left( \oplus_{i=1}^e \O d\log t_i\right)\end{eqnarray}

\item \label{lempbvi} Consider a $p$-adic log scheme $X^\s$ and $(\underline s,\underline t)\in \G(X,\O_X)^d\times \G(X,M_X)^e$ such that $(\underline s_1,\underline t_1)$ is a $p$-basis for $X^\s_1/\Sigma_1$. Then $(\underline s,\underline t)$ is a $p$-basis for $X^\s/\Sigma_k$ (resp. $X^\s/\Sigma_\infty$) if and only if $X$ is flat over $\Sigma_k$ (resp. $\Sigma_\infty$). When this is the case, the morphism $X_k\to \mathbb A^{d+e}_{\Sigma_k}$ induced by $(\underline s,\underline t)$ is flat (resp. for any $k\ge 1$).
\finrom
\end{lem} Proof. \ref{lempbi} Up to \'etale localization we may assume that $X$ is Noetherian and that  $Z$ is a strict normal crossing divisor (\cite{dJ2} Sect. 2.4), ie. $Z=\cup_{i=1}^e Z_i$ (reduced scheme structure) with each $Z_J:=\cap_{i\in J}Z_i$ smooth over $\Sigma_1$ and of codimension $\# J$ in $X$ for each $J\subset \{1,\dots, e\}$. By Zariski localization and induction on $e$, it suffices to assume that $Z_{\{1,\dots, e\}}$ is non empty and to find $p$-bases at the neighborhood of each point of $Z_{\{1,\dots, e\}}$ (by convention  $Z_{\{1,\dots, e\}}=X$ if $e=0$). Consider a closed point $x$ of $Z_{\{1,\dots, d\}}$ and choose for each $i$ a generator $t_i$ of the ideal of $\smash{Spec(\O_{X,x})}\times_XZ_i$ in $Spec(\smash{\O_{X,x}})$.
Since $Spec(\O_{X,x})\times_XZ_{\{1,\dots, e\}}$ is regular of codimension $e$ in $Spec(\O_{X,x})$,  it follows from \cite{Se} Chap. III Prop. 15 and Chap. IV Prop. 22, that $(t_1,\dots,t_e)$ can be completed into a system of parameters generating the maximal ideal $\mathfrak m_{X,x}$ of $\O_{X,x}$, say $(t_1,\dots, t_e,s_1,\dots, s_d)$. The resulting morphism $(\underline s,\underline t):Spec(\O_{X,x})\to \smash{\mathbb A^{d+e}_{\Sigma_1}}$ extends to a morphism $U\to \smash{\mathbb A^{d+e}_{\Sigma_1}}$ which is \'etale at $x$ (\cite{EGA4-IV} Prop. 17.5.3) for some open neighborhood $U$ of $x$. Shrinking $U$ if necessary, we may assume that this morphism is \'etale everywhere (\cite{EGA4-IV} Cor. 17.11.4) and that $Z\cap U=\cup_{i=1}^e V(t_i)$. Then $(\underline s,\underline t)$ is a $p$-basis for $U^\s:=U\times_XX^\s$ and the proof is finished.

Statement \ref{lempbii} is immediate from the definition and  \ref{lempbiii} and \ref{lempbiv} follow respectively from Lem. \ref{lemrp} \ref{lemrpiv} and \ref{lemrpv}. Statement
\ref{lempbv} follows from \cite{EGA4-I} Chap. 0, Cor. 20.7.7 and \cite{Ka2} Prop. 3.12.  In \ref{lempbvi}, the equivalence  follows from \cite{Ka4}  Prop. 1.4 (whose proof relies on \cite{Ka1} Prop. 5.2 and \cite{Ka} Lem. 1) and the last statement follows from Lem. \ref{lemrp} \ref{lemrpv}, \ref{lemrpvi}.
together with the fact that
$\iota_{1,k,*}(\smash{M_{X_1}/\Gm})\simeq \smash{M_{X_k}/\Gm}$.
\begin{flushright}$\square$\end{flushright}

\begin{rem} \label{examplesofpb} Consider a perfect field $k$. Here are some basic examples of a $p$-adic log scheme $X^\s=X_.^\s$ with finite $p$-bases over $\Sigma_\infty$.
\debrom \item  \label{examplesofpbi} $Y_.=Spec(W_.(k)[S_1,\dots, S_d,T_1,\dots,T_e])$ or $Spec(W_.(k)[[S_1,\dots, S_d,T_1,\dots,T_e]])$ and  $Y_.^\s=(Y_.,(\N^e\to \O, e_i\mapsto
T_i))$. A  $p$-basis is given by $(\underline s,\underline t)$ where $s_i=S_i$ and $t_i=e_i$.

\item \label{examplesofpbii} $Y_.=Spec(W_.(k)[S_1,\dots, S_d,T_1,\dots,T_e][(S_1-\alpha_1)^{-1},\dots, (S_d-\alpha_d)^{-1}])$ with $\alpha_i\in W(k)^\times$ or $Spec(W_.(k))$ $[[S_1,\dots, S_d, T_1,\dots,T_e]])$,
$Y_.^\s=(Y_.,(\N^e\to \O, e_i\mapsto T_i))$. A $p$-basis is given by $(\emptyset,\underline t)$ where $t_i=S_i-\alpha_i$, $i\le d$ and $t_{i}=e_{i-d}$,
$i\ge d+1$.

\item \label{examplesofpbiii} Example \ref{examplesofpbi} can be covered \'etale locally by Example \ref{examplesofpbii}: enlarge $k$ if necessary, so that its cardinal becomes at least $d+2$; choosing $\alpha_1=\dots =\alpha_d=\alpha$ for $d+1$ values of $\alpha$ whose reduction mod $p$ are distinct and non zero gives the desired covering. As a result, if $Y^\s$ is an arbitrary $p$-adic log scheme with local finite $p$-bases, then it has local finite finite $p$-bases of the form $(\emptyset,\underline t)$.
\finrom
\end{rem}

\para \label{paraemb} We define categories of embeddings using the notion of finite $p$-bases introduced in  Sect. \ref{parapb}.

\begin{defn}  \debrom \label{defemb1} \item \label{defemb1i} We define the category $Emb^\s$ and a full subcategory $Emb^{\s,glob}\subset Emb^\s$ as follows:

- An object of $Emb^\s$ (resp. $Emb^{\s,glob}$), called a \emph{local} (resp. \emph{global}) \emph{embedding}, is a triple $(U^\s/X^\s,Y^\s,\iota)$ where
$U^\s$ and $X^\s$ are fine separated $\Sigma_1$-log schemes, $Y^\s$ is a $p$-adic log scheme with local finite $p$-bases over $\Sigma_\infty$,
$U^\s/X^\s$ is strict \'etale surjective (resp. $U^\s=X^\s$) and  $\iota:U^\s\rightarrow Y^\s$ is a closed immersion (ie. for all $k$, $U^\s\to \smash{Y^\s_k}$ is a closed immersion in
the sense of log schemes).

- A morphism  $f:(U^\s/X^\s,Y^\s,\iota)\rightarrow  (U'^\s/X'^\s,Y'^\s,\iota')$ consists of a triple of compatible morphisms $(f_X,f_U,f_Y)$.

\item \label{defemb1ii} The category $Emb^\s$ (resp.  $Emb^{\s,glob}$) is viewed as a $Sch^\s/\Sigma_1$-category via the forgetful functor
$(U^\s/X^\s,Y^\s,\iota)\rightarrow X^\s$. The fiber above $X^\s$ is denoted $Emb^\s(X^\s)$ (resp. $Emb^{\s,glob}(X^\s)$) and its objects are called
\emph{local} (resp. \emph{global}) \emph{embeddings for $X^\s$}. We say that $X^\s$ is \emph{locally} (resp. \emph{globally}) \emph{embeddable} if
$Emb^\s(X^\s)$ (resp.  $Emb^{\s,glob}(X^\s)$) is non empty.

\item \label{defemb1iii} A morphism $f=(f_X,f_U,f_Y)$ in $Emb^\s$ (resp. $Emb^{\s,glob}$) is said to be \emph{above $f_X$}, or to \emph{extend $f_X$ locally} (resp. \emph{globally}). We say that $f$ \emph{lifts
$f_X$ locally} (resp. \emph{globally}) if the commutative square $$\xymatrix{U^\s\ar[d]_{f_U}\ar[r]^{\iota}&Y^\s\ar[d]^{f_Y}\\  U'^\s\ar[r]^{\iota'}&
Y'^\s}$$ is moreover cartesian. \finrom
\end{defn}

We often use simplified notations, such as $(U^\s,Y^\s)$ or even $Y^\s$ to designate an object of $Emb^\s$ or $Emb^{\s,glob}$. In practice it will also be useful to consider the following smaller categories.

\begin{defn} \label{defemb1var} \debrom \item \label{defemb1vari} Let $\cal Sch^{\s,slfpb}/\Sigma_1$ (resp. $\cal Sch^{slfpb}/\Sigma_1$) denote the category of separated log schemes  with local finite $p$-bases over $\Sigma_1$ (resp. and whose log structure is trivial).

\item \label{defemb1varii} Let $Emb^{\s,glob,lfpb}$ and  $Emb^{\s,lfpb}$ denote the respective full subcategories of $Emb^{\s,glob}$ and $Emb^\s$ defined by the condition that $X^\s$ (hence $U^\s$) has local finite $p$-bases over $\Sigma_1$.

\item \label{defemb1variii} Let $Emb^{glob}$, $Emb^{glob,lfpb}$, $Emb$ and $Emb^{lfpb}$ denote the respective full subcategories of $Emb^{\s, glob}$, $Emb^{\s,glob,lfpb}$, $Emb^\s$ and $Emb^{\s,lfpb}$ defined by the condition that $X^\s$ (hence $U^\s$) and $Y^\s$ have trivial log structures.
\finrom
\end{defn}

\begin{rem} If $X=X^\s$ then $X$ is locally (resp. globally) embeddable if and only if $Emb(X)$ (resp. $Emb^{glob}(X)$) is non empty (just forget log structures on $Y^\s$).
\end{rem}

\begin{lem} \label{lememb1}
Let $X^\s$ be a separated fine log scheme over $\Sigma_1$.
\debrom
\item \label{lememb1i} The categories $Emb^{\s,glob}(X^\s)$ and $Emb^\s(X^\s)$ have finite non empty products.

\item \label{lememb1ii} Consider a morphism $f_X:X'^\s\rightarrow X^\s$ and assume that $Emb^\s(X'^\s)$ is non empty. For any $Y^\s\in Emb^\s(X^\s)$, there exists $Y'^\s \in Emb^\s(X'^\s)$ and $f:Y'^\s\to Y^\s$ in $Emb^\s$ above $f_X$.

\item \label{lememb1iii} If $X^\s/\Sigma_1$ has finite $p$-bases, then the category $Emb^{\s,glob}(X^\s)$ contains \emph{liftings} (ie. objects satisfying $X^\s\simeq \smash{Y^\s_1}$). Closed sub log schemes of $X^\s$ are in particular globally embeddable.

\item \label{lememb1iv} The category $Emb^\s(X^\s)$ is non empty in the following cases:

\noindent 1) $X^\s/\Sigma_1$ has local finite $p$-bases.

\noindent 2)  $X/\Sigma_1$ is of finite type.

\noindent 3) $X=Spec(A)$ where $A$ is the completion of a $\Fp$-algebra of finite type.


In the latter two cases, the category $Emb^\s(X^\s)$ contains objects $(U^\s/X^\s,Y^\s)$ with $U$ affine and $Y^\s=(Spf(\Zp\{\N^e\}),\N^e)$ for some $e\ge 0$.

\item \label{lememb1v} If $X^\s$ has local finite $p$-bases and $Y^\s\in Emb^\s(X^\s)$ then its logarithmic divided power envelope $T^\s:=D(X^\s,Y^\s)$ is flat over $\Sigma_\infty$.
    
\item \label{lememb1vi} Statements \ref{lememb1i} - \ref{lememb1v} hold without $\s$'s.  
\finrom
\end{lem}
Proof. \ref{lememb1i} The product $(U^\s/X^\s, Y^\s,\iota)\times (U'^\s/X^\s, Y'^\s,\iota')$ is represented by $$(U^\s\times_{X^\s}U'^\s/X^\s, Y^\s
\times Y'^\s,\iota'':U^\s\times_{X^\s}U'^\s\rightarrow Y^\s\times Y'^\s)$$ where $\iota''$ is the natural morphism deduced from $\iota\times\iota'$. It
is indeed a closed immersion thanks to $X^\s$ being separated (\cite{EGA1} Chap. 1, Prop. 5.4.2).

\ref{lememb1ii} Let us begin by choosing an arbitrary $(U'^\s, Y'^\s,\iota')$ in $Emb^\s(X'^\s)$. Observe that $U'^\s\times_{X^\s}U^\s\simeq
U'^\s\times_{X'^\s}(X'^\s\times_{X^\s}U^\s)$ is a strict \'etale surjective $X'^\s$-scheme via the first projection and that the morphism
$\iota'':U'^\s\times_{X^\s}U^\s\rightarrow Y'^\s\times Y^\s$ is a closed immersion. Projecting on second factors, we get the desired morphism above
$f_X$: $$f:(U'^\s\times_{S^\s}U^\s,Y'^\s\times Y^\s,\iota'')\rightarrow (U^\s, Y^\s,\iota)$$

\ref{lememb1iii} The first statement follows from Lem. \ref{lempb} \ref{lempbiv} and the second statement follows from the first one.

\ref{lememb1iv} Case 1 follows from \ref{lememb1iii}. Case 2: Note that if a local embedding of the claimed type exists for each $\smash{X_i^\s}$ in a finite family, then it also exists for the disjoint union of the $\smash{X_i^\s}$'s.  Replacing $X^\s$ with a strict \'etale covering if necessary, we may thus assume given a
fine chart $P\rightarrow M_X$ . Choosing a surjection $\N^e\rightarrow P$ and a closed immersion $X\rightarrow Spec(\Fp[\N^d])$ gives a closed
immersion $\iota:X^\s\rightarrow Spec(\Fp[\N^{d+e}],\N^e)$. Replacing $X^\s$ with an open covering and composing $\iota$ with an appropriate translation, we may assume that this closed immersion factors through
$Spec(\Fp[\N^{d+e}],\N^{d+e})$ and the result follows. Case 3 is deduced from Case 2 by completion (Lem. \ref{lemrp} \ref{lemrpiii}).

\ref{lememb1v}  Since the question is strict-\'etale local on $T^\s$ and since  $et(T^\s)$ and $et(X^\s)$ are naturally equivalent we may assume that $X^\s$ is separated and has finite $p$-bases. But then we can find a lifting $\smash{\tilde X^\s}$ in $Emb^{\s,glob}(X^\s)$ and it follows from \cite{Ts2} Prop. 1.8 that
$\smash{T'^\s_k}:=D(X^\s,\smash{Y_k\times \tilde X^\s_k})$ is simultaneously an algebra of divided power polynomials above $\smash{\tilde X^\s_k}$
and $\smash{T^\s_k}$. Since $\smash{\tilde X^\s_k}$ is flat over $\Sigma_k$, the result follows by faithfully flat descent along
$\smash{T'^\s_k}/\smash{T^\s_k}$.
\begin{flushright}$\square$\end{flushright}

\subsection{Modules with logarithmic connection} \label{mwlc}

\label{appnablamod}

\para \label{intronabla} We review the classical relation between crystals and modules with connection on the small \'etale site. Throughout this section, we let $1\le k<\infty$ and we consider a closed immersion $X^\s\to Y^\s$ and its logarithmic divided power envelope (resp. of order $n$) $T^\s=D(X^\s,Y^\s)$ (resp. $T^{\s n}=D^n(X^\s,Y^\s)$),  where $X^\s$ is over $\Sigma_1$ and $Y^\s$ has local finite finite $p$-bases over $\Sigma_k$.  This means in particular that $Y^\s$ has a strict \'etale covering by $\smash{Y_\lambda^\s}$'s having $p$-bases over $\Sigma_k$ of the form $(\emptyset,\underline t)$ (see. Def. \ref{defpb} \ref{defpbiv} and Rem. \ref{examplesofpb} \ref{examplesofpbiii}; this condition could be relaxed as in \cite{CV}). For the purpose of local descriptions, we fix such a strict covering and we use furthermore the following notations: $\smash{X_\lambda^\s/X^\s}$ denotes the base change of $\smash{Y^\s_\lambda/Y^\s}$ and $\smash{T_\lambda:=D(X^\s_\lambda,Y^\s_\lambda)}$.
We also fix a $p$-basis $(\emptyset,\underline t)$ for each $\smash{Y_\lambda^\s}$ ($e$ and $\underline t=(t_1,\dots,t_e)$ thus depend
on $\lambda$ even though we do not write it in order to keep notations reasonable).

All tensor products and inner homomorphisms are taken with respect to the ring $\cal O$ unless mentioned otherwise.


\para  \label{theringP} Let us denote $\smash  Y^{\s(i)}$ the $(i+1)$-th fold product of $Y^\s$ with itself over $\Sigma_k$  and $(X^\s,T^{\s(i)})$ (resp. $(X^\s,T^{\s(i),n})$) the object of $Crys(X^\s/\Sigma_k)$ which is the  logarithmic divided power envelope (resp. of order $n$) $D(X^\s,Y^{\s(i)})$ (resp. $D^n(X^\s,Y^{\s(i)})$) of the diagonal morphism $X^\s\rightarrow Y^{\s(i)}$. Let $${\cal P}^{\s(i)}_{T,Y}\hbox{ (resp. } {\cal P}^{\s(i),n}_{T,Y} \hbox{ )}$$  denote the ring of $Y_{zar}$ which is obtained from the structural ring  of $T^{\s(i)}$ (resp. $T^{\s(i),n}$) by pullback and pushforward along the obvious morphisms $$Y_{zar}\mathop\leftarrow\limits
T_{zar}^{(0)}\mathop\rightarrow\limits^\sim T_{zar}^{(i)}\hbox{ (resp. }Y_{zar}\mathop\leftarrow\limits T_{zar}^{(0),n}\mathop\rightarrow\limits^\sim
T_{zar}^{(i),n}\hbox{ )}$$ Let $\smash{d_0,\dots, d_i:\cal P^{(0)}_{T,Y}\rightarrow\cal P^{\s(i)}_{T,Y}}$ or $\smash{\cal
P^{(0),n}_{T,Y}\rightarrow\cal P^{\s(i),n}_{T,Y}}$ denote the ring morphisms corresponding to the projections. Unless explicitly mentioned, these
rings will be viewed as $\O_Y$-algebras via $d_0$ and the natural $\O_Y$-algebra structure of $\smash{\cal P^{\s(0)}_{T,Y}}$ (resp. $\smash{\cal
P^{\s(0),n}_{T,Y}}$).

Consider now $f:Y'^\s\rightarrow Y^\s$ and denote $X':=X^\s\times_{Y^\s} Y'^\s$. Consider the logarithmic divided
power envelope (resp. of order $n$)  $T'^{\s(i)}$ (resp. $T'^{\s(i),n}$) of $X'^\s\to Y'^{\s,i+1}$ and denote similarly $\smash{\cal
P^{\s(i)}_{T',Y'}}$ and $\smash{\cal P^{\s(i),n}_{T',Y'}}$ the corresponding $\O_{Y'}$-algebras of $Y'_{zar}$.

\begin{lem} \label{lemdefP1}
If $f$ is strict \'etale then the natural morphism $$\smash{f^*\cal P^{\s(i)}_{T,Y}\rightarrow \cal P^{\s(i)}_{T',Y'}}\hbox{ (resp. }\smash{f^*\cal P^{\s(i),n}_{T,Y}\rightarrow \cal P^{\s(i),n}_{T',Y'}\hbox{ )}}$$ of $Mod(Y'_{zar},\cal O_{Y'})$ is invertible. Here $f^*$ denotes the module pullback functor of the morphism $f:(Y'_{zar},\O_{Y'})\to (Y_{zar},\O_Y)$.
\end{lem}
\noindent Proof. We have to show that the squares of the following commutative diagram are cartesian
$$\xymatrix{X'^\s \ar[r]\ar[d] & T'^{\s(i),n}\ar[r]\ar[d]&  T'^{\s(i)}\ar[r]\ar[d] \ar[r]^{p_0}\ar[d] & Y'^\s\ar[d]\\
X^\s \ar[r] & T^{\s(i),n}\ar[r]&  T^{\s(i)}\ar[r] \ar[r]^{p_0} & Y^\s}$$
Here the arrows denoted $p_0$ are induced by the first projection $\smash{Y^{\s(i)}\to Y^\s}$ and $\smash{Y'^{\s(i)}\to Y'^\s}$. Since logarithmic
divided power envelopes  and logarithmic divided power envelopes of order $n$ commute to strict \'etale base change, this, in turn, is equivalent to
the natural morphism $T'^{\s(i),n}\to D(X'^\s,Y'^\s\times Y^{\s(i-1)})$ being invertible. It suffices to prove that the latter morphism is strict \'etale and that it lifts the identity of $X'^\s$. The result thus follows from the fact that it admits a natural decomposition as $$D(X'^\s,Y'^{\s(i)})\to D(X'^{\s,i+1/X^\s},Y'^{\s(i)})\to
D(X'^\s,Y'^\s\times Y^{\s(i-1)})$$ where the fist arrow is the open immersion induced by the diagonal open immersion $X'^\s\to X'^{\s,i+1/X^\s}$ and the second is the base change of the strict \'etale morphism $Y'^{\s(i)}\to Y^\s\times Y^{\s(i-1)}$ (note that $\smash{X'^\s\times_{Y'^\s\times Y^{\s(i-1)}}Y'^{\s(i)}}\simeq \smash{X'^{\s,i+1/X^\s}}$).
\begin{flushright}$\square$\end{flushright}

This implies in particular that the quasi-coherent modules $\smash{\cal P^{\s(i)}_{T',Y'}}$ and  $\smash{\cal P^{\s(i),n}_{T',Y'}}$ for varying strict \'etale $Y^\s$-log schemes $Y'^\s$
satisfy descent and thus come from algebras of $(Y_{et},\O)$ (see the proof of Lem.-Def. \ref{lemqcohsch}). The latter will simply be denoted $\smash{\cal P^{\s(i)}_{T,Y}}$ and  $\smash{\cal P^{\s(i),n}_{T,Y}}$ as well.
\begin{lem} \label{lemdefP2} Let $\iota:\smash{(T_{et}^\s,\O)\rightarrow (Y^\s_{et},\O)}$ denote the tautological morphism.
\debrom
\item \label{lemdefP2i} There are canonical algebras $\cal P^{\s(i)}_T,\, \cal P^{\s(i),n}_T$ in $(T^\s_{et},\O)$ such that $$\iota_*\cal P^{\s(i)}_T\simeq \cal P^{\s(i)}_{T,Y}\hspace{1cm} \hbox{and} \hspace{1cm} \iota_*\cal P^{\s(i),n}_T\simeq \cal P^{\s(i),n}_{T,Y}$$
\item \label{lemdefP2ii} The algebra $\smash{\cal P_T^{\s(0)}}$ is $\O$ itself. The algebra $\smash{\cal P_T^{\s(i)}}$ has the following explicit local description:
    $$(\cal P^{\s(i)}_{T})_{|T^\s_{\lambda}}\simeq \cal
(\cal P^{\s(0)}_T)_{|T^\s_{\lambda}}<u^{0,\mu}_a-1>_{1\le \mu\le i, 1\le a\le d}$$ with $u^{0,\mu}_{a}$ the unique section of  $\smash{\cal P^{\s(i)}_T}$ over $\smash{T^\s_{\lambda}}$
satisfying $u^{0,\mu}_{a}d_0t_a=d_\mu t_a$. We use the following standard notations in the case $i=1$:  $\tau^\s_a:=u^{0,1}_a-1$ and $dlog\, t_a$ is
its image in $\smash{\cal P^{\s(1),1}(T^\s_{\lambda})}$.
\item \label{lemdefP2iii} Let $\smash{\cal P^{\s(i)}_Y}$, $\smash{\cal P^{(i),n}_Y}$ denote the algebras of $(Y_{et},\O)$ whose construction is similar to $\smash{\cal P^{\s(i)}_T}$, $\smash{\cal P^{(i),n}_T}$ but where the closed immersion $X^\s\to Y^\s$ is replaced from the start by $Y^\s_1\to Y^\s$ (recall that $Y_1:=\smash{\Sigma_1\times_{\Sigma_k} Y}$ and that $Y=D(Y_1,Y)$). Then we have canonical isomorphisms
    $$\iota^*\cal P^{\s(i)}_Y\simeq \cal P^{\s(i)}_T  \hspace{1cm} \hbox{and} \hspace{1cm} \iota^*\cal P^{\s(i),n}_Y\simeq \cal P^{\s(i),n}_T$$
\finrom
\end{lem}
\noindent Proof. \ref{lemdefP2i} It is clear from the construction that $\smash{\cal P^{\s(i)}_{T,Y}}$ and $\smash{\cal P^{\s(i),n}_{T,Y}}$ are supported on
the closed subtopos $X_{et}$ of $Y_{et}$. The result follows. The local description \ref{lemdefP2ii} is a straightforward consequence of the one
given in \cite{Ts2} Prop. 1.8. This description shows that the natural morphisms $\smash{\iota^*\cal P^{\s(i)}_Y\to \cal P^{\s(i)}_T}$ and $\smash{\iota^*\cal P^{\s(i),n}_Y\to \cal P^{\s(i),n}_T}$
are invertible as claimed in \ref{lemdefP2iii}.
\begin{flushright}$\square$\end{flushright}

It follows from \ref{lemdefP2iii} that one has a canonical exact sequence of $\smash{\cal P_T^{\s(1)}}$-modules of $\smash{T^\s_{et}}$.
\begin{eqnarray}\label{defs}\diagram{0&\hflcourte{}{}& \Omega_{T^\s}&\hfl{}{}& \cal P^{\s(1),1}_T&\hfl{s}{}&\O&\hflcourte{}{}&0}\end{eqnarray} where
we denote abusively \begin{eqnarray}\label{defdiffT}\Omega_{T^\s}:=\iota^*\Omega_{Y^\s}\end{eqnarray} and
 $s$ is the morphism induced by the
diagonal immersion. After scalar restriction to $\O$ via $d_0$ (resp. $d_1$) the morphism $s$ has a canonical splitting given by $d_0$ (resp.
$d_1$). We will further denote \begin{eqnarray}\label{cander}d:=d_1-d_0:\cal O\rightarrow\Omega_{T^\s}\subset \smash{\cal
P^{\s(1),1}_T}\end{eqnarray} the canonical $\O$-derivation (it is not universal in general because $\Omega_{T^\s}$ is \emph{not} the module of K\"ahler differentials, see (\ref{defdiffT})).

\para  \label{formulesD}
Since for any strict \'etale $Y'^\s/Y^\s$ as above, the diagonal closed immersion $X'^\s\rightarrow T'^{\s(1)}\times_{T'^\s}T'^{\s(1)}$ is exact
(because $p_0:T'^{\s(1)}\to T'^\s$ is strict) and has divided powers  (as in \cite{Be3} Prop. 2.1.3,  use flatness of the projections $p_0$, $p_1$ to extend the divided powers in two different compatible ways), the morphism $Y'^{\s 2}\times_{Y'}Y'^{\s 2}\rightarrow Y'^{\s 2}$,
$(a,b,b,c)\mapsto (a,c)$ induces $\smash{T'^{\s(1)}\times_{T'^\s}T'^{\s(1)}}\rightarrow T'^{\s(1)}$. The resulting  morphisms for varying $Y'^\s$ are
compatible with each other (use that $T'^{\s(1)}\simeq Y'^\s\smash{\times_{Y^\s}}T^{\s(1)}$) and gives rise to morphisms as follows in
$Mod(\smash{T^\s_{et}},\O)$: \begin{eqnarray}\label{deltaP}\begin{array}{rcl}\delta&:&\cal P^{\s(1)}_T\rightarrow \cal
P^{\s(1)}_T{_{d_1}\otimes_{d_0}} \cal P^{\s(1)}_T\\ \delta^{n,n'}&:&\cal P^{\s(1),n+n'}_T\rightarrow \cal P^{\s(1),n'}_T{_{d_1}\otimes_{d_0}} \cal
P^{\s(1),n}_T.\end{array}\end{eqnarray} These morphisms are compatible with their counterparts for $\cal P^{\s(1)}_{Y}$. Let us write an explicit
local formula. Over $\smash{T^\s_\lambda}$ we have
$$\delta(\tau^\s_a+1)=(\tau^\s_a+1)\otimes (\tau^\s_a+1).$$





\begin{defn} \label{defDiff}
We define $\cal D$, the \emph{$\O$-algebra of $dp$-differential operators of $X^\s$ inside $Y^\s$}, as the following module of
$(\smash{T^\s_{et},\O})$
$$\begin{array}{rcl}
\cal D^{\s}:=\limi_n \cal D^{\s n} &\hbox{ where }& \cal D^{\s n}:=\cal Hom_{\cal O}(\cal P_T^{\s(1),n},\O)
\end{array}$$
together with the composition rule  $\partial \partial':=\partial\circ (id\otimes \partial')\circ \delta$.
\end{defn}
We refer to \cite{Mo} Sect. 2.3.2 or \cite{BO} Def. 4.4 for the verification that this composition rule defines indeed an $\O$-algebra structure.
Let us only write some local formulae. Recall that $\smash{(\tau^{\s[\underline n]})_{\underline n}}$  is a basis of  the restriction of $\smash{\cal
P_T^{\s(1)}}$ to $\smash{T^\s_{\lambda,et}}$. Let $\smash{(\partial^\s_{[\underline n]})_{\underline n}}$ denote the dual basis of the restriction of
$\cal D^\s$ to $\smash{T^\s_{\lambda,et}}$. For $f\in \G(\smash{T^\s_\lambda,\O})$ and $\partial^\s_a:=\smash{\partial_{e_a}}\in
\G(\smash{T^\s_\lambda,\cal D^\s})$ (here $e_a=(0,\dots,0,1,0,\dots, 0)$ is the $d$-tuple with $1$ in $a^{th}$ position) we have:

\begin{eqnarray}\label{formulD} \begin{array}{rcl} \hbox{in $\G(T^\s_{\lambda},\iota^*\Omega_{Y^\s})$}&: & d(f)=\sum_a \partial^\s_a(f)dlog\,  t_a \hbox{,}\\
\hbox{in $\G(T^\s_{\lambda},\cal D^\s)$}&:& \partial^\s_a f-\partial^\s_a(f)=f\partial^\s_a \hbox{,}\\  \hbox{in $\G(T^\s_{\lambda},\cal D^\s)$}&:& \partial_a^\s\partial_b^\s=\partial^\s_b\partial^\s_a\hbox{ and }  \smash{\partial^\s_{[\underline n]}}={\Pi_{a=1}^d\Pi_{i=0}^{n_a-1}(\partial^\s_{a}-i)} \hbox{.}
\end{array}\end{eqnarray}

\noindent In other terms, the restriction of $\smash{\cal D^{\s}}$ to $\smash{T^\s_{\lambda}}$
is a ring of non commutative polynomials with the $\partial^\s_a$'s as commuting variables.



\para \label{etlognabla} Following closely \cite{BO}, we discuss briefly integrable quasi-nilpotent connections and hyper $dp$-stratifications.

\begin{lem} \label{defnabla} Let $M$ in $Mod(\smash{T^\s_{et},\O})$. The following data \ref{defnablai}, \ref{defnablaii}, \ref{defnablaiii} are equivalent.
\debrom

\item \label{defnablai} A morphism  $\nabla:M\rightarrow M\otimes\smash{\Omega_{T^\s}}$ of Abelian groups in $\smash{T^\s_{et}}$ satisfying
$\nabla(f x)=f \nabla(x)+x\otimes d(f)$ for any strict \'etale $\smash{T'^\s/T}$, $f\in \G(T'^\s,\O)$ and $x\in \G(T'^\s,M)$.
\item \label{defnablaii} A morphism $\theta_1:M\rightarrow M\otimes_{d_0} {\cal P^{\s(1),1}_T}$ in $Mod(\smash{T^\s_{et}},\O)$ where the target is endowed with its right
$\O$-module structure (ie. the one coming from the $\O$-algebra structure of $\smash{{\cal P^{\s(1),1}_T}}$ given by $d_1$) which gives back the
identity of $M$ when composed with $\smash{M\otimes_{d_0}{\cal P^{\s(1),1}_T}}\rightarrow M_T$.
\item \label{defnablaiii} A morphism $\nabla:{\cal D^{\s(1),1}}{_{d_1}\otimes} M\rightarrow M$ in $Mod(\smash{T^\s_{et}},\O)$ where the source is endowed with its left
$\O$-module-structure (ie. the one coming from the $\O$-algebra structure of $\smash{{\cal P^{\s(1),1}_T}}$ given by $d_0$) which gives back the
identity of $M$ when composed with $M\rightarrow \smash{{\cal D^{\s(1),1}}{_{d_1}\otimes} M}$.  \petit \finrom Such data is called a
\emph{(logarithmic) connection} on $M$.
\end{lem}
\noindent Proof. $\ref{defnablai}\Leftrightarrow \ref{defnablaii}$ is by setting $\theta_1(x)=(x\otimes 1)+\nabla(x)$ and $\ref{defnablaii}\Leftrightarrow
\ref{defnablaiii}$ is by adjunction using that $\smash{{\cal P^{\s(1),1}_T}}$ is locally free of finite type.
\begin{flushright}$\square$\end{flushright}

\begin{lem} \label{defnablaint} Let $(M,\nabla)$ be a module with connection as in Lem. \ref{defnabla}.  The following conditions \ref{defnablainti}, \ref{defnablaintii}, \ref{defnablaintiii} are equivalent.
\debrom
\item \label{defnablainti} The \emph{curvature} morphism $$K:M\rightarrow M\otimes \wedge^2\smash{\Omega_{T^\s}}$$  is zero. Here $K$ is the morphism of Abelian groups of $\smash{T^\s_{et}}$ defined as
$K(x)=\nabla^1(\nabla(x))$ where $\nabla^1(x\otimes \omega)=\nabla(x)\wedge \omega+ x\otimes d^1(\omega)$ and $d^1(y\, dlog\,  z)=dy\wedge dlog\,
z$.
\item \label{defnablaintii} There exists a morphism $\theta_2$ rendering the following square of $Mod(\smash{T^\s_{et}},\O)$ commutative:
    $$\xymatrix{M\ar[r]^-{\theta_1}\ar[d]_{\theta_2}&M\otimes_{d_0}\cal P^{\s(1),1}\ar[d]^{\theta_1\otimes 1}\\ M\otimes_{d_0}\smash{{\cal P^{\s(1),2}_T}}\ar[r]_-{1\otimes \delta^{1,1}}& M\otimes_{d_0}\smash{{\cal P^{\s(1),1}_T}}{_{d_1}\otimes_{d_0}}\smash{{\cal P^{\s(1),1}_T}}}$$
    Here the tensor products modules are viewed as $\O$-modules via $d_1$ on the last factor.
\item \label{defnablaintiii} There exists a morphism $\nabla_2$ rendering the following square of $Mod(\smash{T_{et},\O})$ commutative:
    $$\xymatrix{{\cal D^{\s 1}}{_{d_1}\otimes_{d_0}}{\cal D^{\s 1}}{_{d_1}\otimes} M\ar[d]_{\nabla_2}\ar[r]^-{1\otimes \nabla}&
    {\cal D^{\s 1}}{_{d_1}\otimes} M\ar[d]^-\nabla\\
    {\cal D^{\s 2}}{_{d_1}\otimes} M\ar[r]&M}$$
    Here the bottom arrow is induced by the composition in $\cal D^\s$ and the tensor product is viewed as $\O$-modules, via left multiplication on the first factor. \finrom
When these conditions are satisfied we say that the connection $\nabla$ is \emph{integrable}.
\end{lem}
\noindent Proof. A straightforward computation with our chosen $p$-basis of $\smash{Y_\lambda^\s}$ shows that
$$\hspace{-.5cm}\begin{array}{rcl}(\theta_1\otimes 1)(\theta_1(x))&=&x\otimes 1\otimes 1+\sum_{a} (\nabla (\partial^\s_a\otimes x))\otimes
\delta^{1,1}(\tau^\s_a) \\&&+\sum_a (\nabla( \partial^\s_a\otimes (\nabla (\partial^\s_a-1)\otimes(x))\otimes \delta^{1,1}(\tau_a^{\s[2]})+\sum_{a\le
b} (\nabla\partial^\s_a)(\nabla\partial^\s_b)(x)(\delta^{1,1}(\tau_a^\s\tau_b^\s))\\ &&+\tilde K(x)\\ \hbox{where }\tilde K(x)&=&\sum_{a<b}
((\nabla(\partial^\s_b\otimes (\nabla(\partial^\s_a\otimes x))))-(\nabla(\partial^\s_b\otimes (\nabla(\partial^\s_a\otimes
x)))))\otimes(\tau^\s_b\otimes \tau^\s_a)\end{array}$$ and that $\tilde K(x)$ is sent to $K(x)$ via the canonical morphism
$$M\otimes_{d_0}\smash{{\cal P^{\s(1),1}}}{_{d_1}\otimes_{d_0}}\smash{{\cal P^{\s(1),1}}}\rightarrow M\otimes_{d_0} Coker(\delta^{1,1})\simeq
M\otimes \wedge^2\smash{{\cal P^{\s(1),1}}}$$
This shows that conditions \ref{defnablainti} and \ref{defnablaintii} are equivalent. Conditions \ref{defnablaintii} and \ref{defnablaintiii} on the
other hand are clearly equivalent by adjunction.
\begin{flushright}$\square$\end{flushright}
If $\nabla$ is an integrable connection on $M$, then it follows from the explicit description of $\cal D^\s$ that there is a unique structure of $\cal
D^\s$-module on $M$ extending $\nabla$ (ie. such that $\partial x:=\nabla(\partial\otimes x)$ for any $\partial$ in $\cal D^{\s 1}$ and $x$ in $M$). By
adjunction, one deduces a right $\O$-linear morphism \begin{eqnarray}\label{deftheta}\theta:M\rightarrow \limp_n(M\otimes_{d_0}\smash{{\cal
P_T^{\s(1),n}}})\end{eqnarray} lifting $\theta_1$ and satisfying the cocycle condition $(\theta\otimes 1)(\theta(x))=\delta(\theta(x))$. The explicit local
description of $\theta$ is given by the following \emph{Taylor formula}:
\begin{eqnarray}\label{devtaylor}\theta(x)=(\sum_{|\underline n|\le n}\partial^\s_{[\underline n]}x\otimes \tau^{\s[\underline n]})_n\end{eqnarray}

\begin{lem} \label{defqnil} Consider a module with integrable connection $(M,\nabla)$. The following conditions \ref{defqnili}, \ref{defqnilii} are equivalent.
\debrom \item  \label{defqnili} For all strict \'etale $\smash{T'^\s}$ over $T^\s$ and $x\in \G(\smash{T'^\s,M})$,  all but a finite number of the $\smash{\partial^\s_{[\underline
n]}x}$'s vanish in $\G(T'^\s,M)$.
\item \label{defqnilii} The morphism (\ref{deftheta}) factors through a morphism
$$\theta:M\to M\otimes_{d_0} \cal P_T^{\s(1)}$$
\finrom When these conditions are satisfied, we say that the integrable connection $\nabla$ is \emph{quasi- nilpotent}.
\end{lem}
\noindent Proof. This is straightforward from the formula (\ref{devtaylor}).
\begin{flushright}$\square$\end{flushright}



\begin{prop} \label{hyperstrat} Let $M$ in $Mod(\smash{T^\s_{et}},\O)$. The data of an integrable quasi-nilpotent connection $\nabla$ on $M$ is equivalent to a \emph{hyper $dp$-stratification} ie. a $\cal P^{\s,(1)}_T$-linear isomorphism $$\varepsilon:\cal
P_T^{\s(1)}{_{d_1}\otimes} M\simeq M\otimes_{d_0}\cal P_T^{\s(1)}$$  satisfying the cocycle condition $$(\cal P_T^{\s(2)}\smash{_{d_{0,1}}\otimes_{\cal
P_T^{\s(1)}}} \varepsilon)\circ(\cal P_T^{\s(2)}\smash{_{d_{1,2}}\otimes_{\cal P_T^{\s(1)}}} \varepsilon)=(\cal P_T^{\s(2)}\smash{_{d_{0,2}}\otimes_{\cal
P_T^{\s(1)}}} \varepsilon)$$
 \end{prop}
Proof. One deduces $\varepsilon$ from $\theta$ by scalar extension via $d_1$. Then, condition Lem. \ref{defnablaint} \ref{defnablaintii} may be translated into
the cocycle condition using that $\cal P^{\s(2)}_T\simeq \cal P^{\s(1)}_T{_{d_1}\otimes_{d_0}}\cal P^{\s(1)}_T$ and that via this identification,
$\delta: \cal P^{\s(1)}_T\rightarrow \cal P^{\s(1)}_T{_{d_1}\otimes_{d_0}}\cal P^{\s(1)}_T$ translates into $d_{0,2}$.
\begin{flushright}$\square$\end{flushright}




\para  \label{crystalss}

We use the following notation.

\begin{defn} \label{defnablamod} The category $\nabla$-$Mod^{pd}(X^\s,Y^\s)$, also abusively denoted $\nabla$-$Mod(T^\s)$, is defined as follows. An object is a module with quasi-nilpotent integrable connection. A morphism $(M,\nabla)\to (M',\nabla')$ is an $\O$-linear morphism which is compatible with the given connections. If $*$ is either $qcoh$, $lf$ or $lfft$ then we denote furthermore $\nabla$-$Mod_*(T^\s)$ the full subcategory formed by the $(M,\nabla)$ with $M$ satisfying $*$.
\end{defn}
\noindent Note that, as in the case of the algebra of $dp$-differential operators, the reference to $T^\s$ only is abusive since the category depends in fact \emph{a
priori} on the immersion $X^\s\to Y^\s$ rather than $T^\s$ itself.

We have the following analogue of \cite{BO} Thm. 6.6.

\begin{prop}  \label{crystals}  There is a canonical equivalence
$$Crys_*((X^\s/\Sigma_k)_{crys,et},\O)\simeq \nabla\hbox{-}Mod_*(T^\s).$$
\end{prop}
Proof. The local lifting property for log schemes with local finite $p$-bases together with the universal property of logarithmic divided power envelopes
ensure that $T$ covers the final object of $(X^\s/\Sigma_k)_{crys,et}$. As a result,  the category of modules of
$((X^\s/\Sigma_k)_{crys,et},\O)$ is equivalent to the category of modules of $((X^\s/\Sigma_k)_{crys,et}/T^\s,\O)$ endowed with a descent datum
satisfying the cocycle condition. It is thus equivalent to the category of modules on $(T_{et},\O)$
endowed with the extra data coming from the descent datum (Lem. \ref{lemcrystal1} \ref{lemcrystal1i} and  \ref{lemcrystal1iii}). It only remains to notice that a descent datum $\varepsilon:p_1^*M\simeq p_0^*M$ on a module
$M$ over $(X^\s/\Sigma_k)_{crys,et}/T^{\s}$ exactly translates into a descent datum
$\varepsilon:p_1^*\smash{M_{T^\s}}\simeq p_0^*\smash{M_{T^\s}}$ on the realization $\smash{M_{T^\s}}$ ie. (using the diagonal equivalence
$T_{et}\rightarrow \smash{T_{et}^{(1)}}$) into  a $\smash{\cal P^{\s(1)}_T}$-linear isomorphism $\varepsilon:\smash{\cal P^{\s(1)}_T{_{d_1}\otimes}
M\simeq M\otimes_{d_0}\cal P^{\s(1)}_T}$.
We may conclude by Prop. \ref{hyperstrat} in the case $*=\emptyset$. The other cases follow (see Rem. \ref{lemlcrys} \ref{lemlcrysii}, \ref{lemlcrysiv} for the meaning of the category on the left).
  \begin{flushright}$\square$\end{flushright}

%

%

\para \label{functcrys} We explain inverse images for modules with connection.

\begin{lem} \label{functnabla} Let $X\to Y$ and $X'\to Y'$ be as in Sect. \ref{intronabla} and assume given a commutative square $$\xymatrix{X\ar[r]\ar[d]_{f_X}&Y\ar[d]^{f_Y}\\ X'\ar[r]&Y'}$$ Let $f_T:T^\s\to T'^\s$ denote the morphism obtained by forming logarithmic divided power envelopes. Let $M'$ in $\cal Crys((X'^\s/\Sigma_{k,crys,et}),\O)$  and consider its pullback  $M:=f_X^*M'$ in $\cal Crys((X^\s/\Sigma_k)_{crys,et},\O)$. If $M'$ corresponds to $(N',\nabla')$ in $\nabla'$-$Mod(T'^\s)$, then $M$ corresponds to $(N,\nabla)$ in $\nabla$-$Mod(T^\s)$, where $N\simeq f_T^*N'$ and $\nabla$ has the following alternative characterizations:
\debrom
\item \label{functnablai} if $\nabla'$ corresponds to $\varepsilon':\smash{\cal P^{\s(1)}_{T'}\otimes N'\simeq N'\otimes\cal P^{\s(1)}_{T'}}$ on $\smash{(T'^{\s}_{et},\cal
P^{\s(1)}_{T'})}$ then $\nabla$ corresponds to the morphism $\varepsilon:\smash{\cal P_T^{\s(1)}}\otimes N\simeq N\otimes\smash{\cal P_{T}^{\s(1)}}$
deduced from $\varepsilon'$ by pullback via $\smash{(T_{et},\cal P^{\s(1)}_T)}\rightarrow \smash{(T'_{et},\cal P^{\s(1)}_{T'})}$.

\item \label{functnablaii} if $\nabla'$ corresponds to $\theta': N'\rightarrow  N'\otimes\smash{\cal P^{\s(1)}_{T'}}$ on $\smash{(T'^\s_{et},\O)}$ (the $\O$-module structure
on the target is via $d_1$) then $\nabla$ corresponds to the morphism $\smash{\theta:N\rightarrow  N\otimes\cal P^{\s(1)}}$ deduced from $\theta'$ by
pullback via $f_T:(\smash{T^\s_{et}},\O)\rightarrow (\smash{T'^\s_{et}},\O)$ and the natural ``base change morphism'' $f_T^* (\smash{N'\otimes\cal
P^{\s(1)}_{T'}})\rightarrow  \smash{N\otimes\cal P^{\s(1)}_T}$.

\item \label{functnablaiii} $\nabla$ is the unique connection on $M$ rendering the following diagram of Abelian groups of $T_{et}$ $$\xymatrix{N\ar[r]^-\nabla&N\otimes
\Omega_{T^\s}
\\ f_T^{-1}N'\ar[u]\ar[r]^-{f_T^{-1}\nabla'}&\ar[u]f_T^{-1}(N'\otimes \Omega_{T'^\s})}$$
commutative, the right vertical arrow being induced by $\smash{f_T^{-1}N'}\rightarrow \smash{N}$, $\smash{f_T^{-1}\O}\rightarrow \smash{\O}$ and $\smash{f_T^{-1}\Omega_{T'^\s}}\rightarrow \smash{\Omega_{T^\s}}$. \finrom
\end{lem}
Proof. The first statement is Lem. \ref{lemcrystal2} \ref{lemcrystal2ii}  since $N'=\smash{M'_{T'^\s}}$ and $N=\smash{M_{T^\s}}$.
Characterization \ref{functnablai} is also an easy consequence of Lem. \ref{lemcrystal2} \ref{lemcrystal2ii} by looking at the commutative following commutative diagram of $CRYS^\s_{et}(X'^\s/\Sigma_1)$ $$\xymatrix{(X,T^\s)\ar[d]&\ar[l]_-{p_0}(X^\s,T^{\s(1)})\ar[d]\ar[r]^-{p_1}&(X,T^\s)\ar[d]\\
(X'\to T'^\s)&\ar[l]_-{p_0}(X'^\s,T'^{\s(1)})\ar[r]^-{p_1}&(X'^\s,T'^\s)}$$ Characterization \ref{functnablaii} follows immediately and
\ref{functnablaiii} as well (note that by linearity, $\theta$ is characterized by its composition with the $f_T^{-1}\O$-linear map
$f_T^{-1}N'\rightarrow N$).
\begin{flushright}$\square$\end{flushright}

\subsection{Crystalline and de Rham cohomology of crystals (linearization functors)} \label{cadrcoclf} ~~ \\
\label{crysdRcoh}

We discuss the de Rham interpretation of the \'etale crystalline cohomology of a crystal in the local case following the exposition of \cite{BO}. We keep
the notations and assumptions of Sect. \ref{appnablamod}. Unless mentioned otherwise, tensor products and inner homomorphisms are taken with respect to the
ring $\O$.

\para Let us begin with the definition of the category of hyper $dp$-differential operators.
\begin{defn} \label{defHdp} The category $Hdp^{dp}(X^\s,Y^\s)$, also abusively denoted $Hdp(T^\s)$, is defined as follows:

- an object is a module of $(\smash{T^\s_{et}},\O)$,

- the set $Hdp(M,N)$ of morphisms from $M$ to $N$ is $Hom(\smash{\cal P_T^{\s(1)}}_{d_1}\otimes M,N)$ (the tensor product is viewed as an $\O$-module via
$d_0$ on $\smash{\cal P^{\s(1)}_T}$) and composition is defined by the formula $\smash{f\circ_{Hdp}g}:=f\circ (1\otimes g)\circ \delta$.
\end{defn}

We note that if $\O^{d=0}$ denotes the kernel of the canonical derivation $d:\O\to \smash{\Omega_{T^\s}}$ (\ref{cander}) we have the forgetful
functor \begin{eqnarray}\label{forHdpMod}Hdp(T^\s)\to Mod(T^\s_{et},\O^{d=0})\end{eqnarray}

.

\begin{lem} \label{tensnablaHdp} There is a functor $$\otimes:\nabla\hbox{-}Mod(\smash{T^\s})\times Hdp(T^\s)\rightarrow Hdp(T^\s)$$ which sends:

- a couple of objects $((M,\nabla),N)$ to $M\otimes N$,

- a couple of morphisms $(f:(M_1,\nabla_1)\to (M_2,\nabla_2),g:\smash{\cal P^{\s(1)}_T}\otimes N_{1}\rightarrow N_{2})$ to the composed morphism
$(f\otimes g)\circ (\varepsilon_1\otimes 1):\smash{\cal P^{\s(1)}_T}\otimes M_1\otimes N_{1}\rightarrow M_2\otimes  N_{2}$ where $\varepsilon_1$ is the hyper
$dp$-stratification corresponding to $\nabla_1$.
\end{lem}
Proof. Use the cocycle condition for $\varepsilon_1$ in order to check the compatibility with respect to composition in $Hdp(T^\s)$.
\begin{flushright}$\square$\end{flushright}

\para We review the de Rham complex of a module with logarithmic connection.

\begin{lem} \label{dR} \debrom \item \label{dRi}
There exists a unique complex of $Mod(T_{et},\O^{d=0})$ $$\diagram{\Omega^\bullet_{T^\s}&:=&[\O&\hflcourte{d^0}{}&
\Omega_{T^\s}&\hflcourte{d^1}{}& \dots &\hflcourte{d^{e-1}}{}& \wedge^e\Omega_{T^\s}]}$$ called the \emph{de Rham complex}, whose differentials are characterized by the
following formulae:

- $d^0=d:=d_1-d_0:\O\rightarrow \smash{\Omega^1_{T^\s_.}\subset \cal P_T^{\s(1),1}}$,

- $d^{q+1}(\omega\wedge dlog\, m)=(d^q\omega)\wedge dlog\, m)$ (or equivalently $d^{q+q'}\omega\wedge
\omega'=d^q\omega\wedge\omega'+(-1)^q\omega\wedge d^{q'}\omega'$)

\item \label{dRii} There is a canonical complex in $Hdp(T^\s)$ which is sent to $\smash{\Omega_{T^\s}^\bullet}$ under the forgetful functor (\ref{forHdpMod}). Its differentials $\tilde d^q$ are given by the formulae:

-  $\tilde d^0:=1-d_1s$ (with $s$ as in (\ref{defs})),

- $\smash{\tilde d^q(\lambda\otimes \omega):=(\tilde d^0\lambda)\wedge \omega +s(\lambda)d^q\omega}$. \finrom
\end{lem}
Proof. \ref{dRi} A simplicial construction of the differentials complex can be carried out by adapting \cite{Il} Chap. VIII, Prop. 1.2.8 to our case
(denote $I=Ker\, s:\smash{\cal P_T^{\s(1)}}\rightarrow \O$ and check directly that the well defined maps $d^0=d_1-d_0:\O\rightarrow I$,
$d^{q+1}:I\otimes I^{\otimes q}\rightarrow I^{\otimes q+2}$, $a\otimes b\mapsto (1\otimes a-\delta a+ a\otimes 1)-a\otimes d^{q}b$ induce
$d^q:\smash{\wedge^q(I/I^{[2]})}\rightarrow \smash{\wedge^{q+1} (I/I^{[2]})}$ satisfying the desired formulae).

\ref{dRii} The reader may check that $\tilde d^{q+1}\tilde d^q=0$.
\begin{flushright}$\square$\end{flushright}

\begin{defn} \label{dRM} For $M=(M,\nabla)$ in $\nabla$-$Mod(T^\s)$ we define $$\Omega^\bullet_{T^\s}(M):=(M,\nabla)\otimes \Omega^\bullet_{T^\s}$$
in $Hdp(T^\s)$. Here the de Rham complex on the right is meant in the sense of Lem. \ref{dR} \ref{dRii} and the tensor product is Lem. \ref{tensnablaHdp}.
\end{defn}

The differential of degree $q$ of this complex is denoted $\tilde \nabla^q$ and can be described explicitly  by the following formula:  $$\smash{\tilde
\nabla^q}(\lambda\otimes m\otimes \omega)=(m\otimes \smash{\tilde d^0}(\lambda)+s(\lambda)\nabla(m))\wedge\omega+s(\lambda)m\otimes d^q\omega$$

\para We discuss the linearization functor. If $M$ is an $\O$-module, we usually view $\smash{\cal P^{\s(1)}_T\otimes M}$ as an $\O$-module via $d_0$ and $\smash{M\otimes \cal P^{\s(1)}_T}$ as an $\O$-module via $d_1$. In order to keep notations simple, we will use the obvious isomorphism of $\O$-modules \begin{eqnarray}\label{defexch}exch:\smash{\cal P^{\s(1)}_T\otimes M}\simeq \smash{M\otimes \cal P^{\s(1)}_T}\end{eqnarray}
(this is not a hyper $dp$-stratification, only a triviality) obtained by ``exchanging factors in $\smash{\cal P^{\s(1)}_T}$'' (in terms of local coordinates
$exch(\tau^\s_a\otimes m)=m\otimes (1+\tau^\s_a)^{-1}-1$).

\begin{lem} \label{defLT} There is a functor $$L_{T^\s}:Hdp(T^\s)\to \nabla\hbox{-}Mod(T^\s)$$ which sends:

- an $\O$-module $M$ to $\smash{L_{T^\s}}(M):=(\smash{\cal P^{\s(1)}_T\otimes M},\nabla)$ where $\nabla$ is ``derivation on the first factor'', ie.
corresponds to the $d_1$-linear morphism $$\xymatrix{ \smash{M\otimes \cal P^{\s(1)}_T}\ar[r]_-{exch}^\sim&\smash{M\otimes \cal
P^{\s(1)}_T}\ar[r]_-{id\otimes \delta}&M\otimes \cal P^{\s(1)}_T\otimes \cal P^{\s(1)}_T\ar[r]^\sim_-{exch^{-1}}&(\cal P^{\s(1)}_T\otimes M)\otimes
\cal P^{\s(1)}_T \ar @{<-} @<+0pt> `u[l] `[lll]_-{\theta} [lll] }$$

- an $\O$-linear morphism $f:\cal P^{\s(1)}_T\otimes M\to N$ to the composed $\O$-linear morphism
$$\xymatrix{\cal P^{\s(1)}_T \otimes M\ar[r]_-{(\delta\otimes id)}& \cal P^{\s(1)}_T \otimes \cal P^{\s(1)}_T \otimes M\ar[r]_-{id\otimes f}&\cal P^{\s(1)}_T \otimes N\ar @{<-} @<+0pt> `u[l] `[ll]_-{\smash{L_{T^\s}}(f)} [ll]}$$
\end{lem}
Proof. Once noticed that $\smash{L_{T^\s}}(f)=exch^{-1}\circ (id\otimes \delta)\circ (f\otimes id)\circ exch$ one may check that everything boils
down to the coassociativity of $\delta$.
\begin{flushright}$\square$\end{flushright}

Let us write a local formula for the connection $\nabla$ of $\smash{L_{T^\s}}(M)$: $$\nabla({\sum_{\underline
n}\prod_a((1+\tau^\s_a)^{-1}-1)^{[n_a]}\otimes x_{\underline n}})= \smash{\sum_{\underline n}\sum_b
(1+\tau_b^\s)^{-1}\prod_a((1+\tau^\s_a)^{-1}-1)^{[n_a]}\otimes x_{\underline n+e_b}\otimes d\log t_b}$$

\begin{lem} \label{isoLtens}
There is a canonical isomorphism in $\nabla\hbox{-}Mod(\smash{T^\s})$
\begin{eqnarray} M\otimes \smash{L_{T^\s}}(N)\simeq \smash{L_{T^\s}}(M\otimes N)\end{eqnarray}
which is functorial with respect to  $(M,N)$ in $\nabla\hbox{-}Mod(\smash{T^\s})\times Hdp(T^\s)$. Here, the first tensor product is meant in the
sense of Lem. \ref{tensnablaHdp} while the second is the usual one in $\nabla\hbox{-}Mod(\smash{T^\s})$ (the connection on $M''=M\otimes M'$ is the one
corresponding to $\varepsilon''=(\varepsilon'\otimes 1)\circ(1\otimes \varepsilon)$ ie.  $\nabla''=\nabla\otimes 1+1\otimes \nabla'$).
\end{lem}
Proof. The claimed isomorphism is
 $$\xymatrix{M\otimes \cal P_T^{\s(1)}\otimes N\ar[r]^-\sim_-{exch_N}&N\otimes \cal P_T^{\s(1)}\otimes M\ar[r]_-\sim^-{id\otimes \varepsilon_M}& N\otimes
M\otimes \cal P_T^{\s(1)}\ar[r]^\sim_-{exch_{M\otimes N}}&\cal P^{\s(1)}_T\otimes M\otimes N}$$
Compatibility with connections and bifunctoriality is easily checked using the cocycle condition for $\varepsilon_M$.
\begin{flushright}$\square$\end{flushright}

\para We have now all the ingredients needed to state and prove the logarithmic version of the Poincar\'e lemma.

\begin{lem} \label{defaug}  The morphism $d_0:\O\to \smash{\cal P^{\s(1)}_T}$ induces a quasi-isomorphism
$$\xymatrix{\O\ar[r]^-{aug}&\smash{L_{T^\s}}(\Omega^\bullet_{T^\s})}$$ of complexes in the Abelian category  $\nabla\hbox{-}Mod(T^\s)$.
\end{lem}
Proof. The differentials of the complex of $\O$-modules underlying $\smash{\smash{L_{T^\s}}(\Omega_{T^\s}^\bullet)}$ are given by the following
formulae, which are the logarithmic variant of \cite{BO} Lem. 6.11.  Recall that $\smash{\smash{L_{T^\s}}(\Omega^q_{T^\s})}=\smash{\cal P_T^{\s(1)}\otimes
\Omega^q_{T^\s}}$ is a locally free $\O$-module with local basis formed by the elements $\smash{\tau^{\s[\underline n]}\otimes \tau^\s_{i_1}\wedge
\tau^\s_{i_2} \dots \wedge \tau^\s_{i_q}}$ indexed by $\underline n\in \N^e$, $1\le i_1<i_2 <i_q\le e$. We have then $$d^q\smash{(\tau^{\s[\underline
n]}\otimes \tau^\s_{i_1}\wedge \tau^\s_{i_2} \dots \wedge \tau^\s_{i_q})}=\sum_{a=1^e}(1+\tau^\s_a)\tau^{\s[\underline n-e_a]}\otimes
\tau^\s_a\wedge\tau^\s_{i_1}\wedge \tau^\s_{i_2} \dots \wedge \tau^\s_{i_q}$$ with the convention that $\tau^{\s[\underline n]}=\prod
\tau_a^{\s[n_a]}=0$ as soon as one of the $n_a$'s is $-1$. Now we recognize a Koszul complex which is a resolution of $\O$. Namely, this is the
$e^{th}$ tensor power of the length $1$ complex $[\O<\tau^\s> \rightarrow \O<\tau^\s>]$, $\tau^{[n]}\mapsto (\tau^\s+1)\tau^{\s[n-1]}$, which is
itself a resolution of $\O$ (use that the series $log(1+\tau^\s)$ is convergent).
\begin{flushright}$\square$\end{flushright}

\para In order to achieve the computation of crystalline cohomology using the de Rham complex, we need an interpretation of the linearization functor in terms of the crystalline site. We only treat the case of $\s$-big crystalline topoi for simplicity but the exact same results and proofs apply with $CRYS$ or $crys$ instead of $CRYS^\s$.

\begin{lem} \label{lemft}
If $M$ is a \emph{local crystal of $\smash{((X^\s/\Sigma)_{CRYS^\s,et},\O)}$} ie. is in the essential image of the fully faithful functor $$\lambda_{T^\s}^*:Mod(T^\s_{et},\O)\to Mod((X^\s/\Sigma_k)_{CRYS^\s,et}/T^\s,\O)$$
then $f_{T^\s,*}M$ is a crystal.
\end{lem}
\noindent Proof. 
We have to check condition Lem. \ref{lemcrystal1} \ref{lemcrystal1ii}
for $T^\s$ (see Lem. \ref{lemcrystal1}\ref{lemcrystal1iii}) ie. that the natural morphism

\begin{eqnarray}\label{ftcrys1}h^*(f_{T^\s,*}M)_{T^\s}\to (f_{T^\s,*}M)_{T'^\s}\end{eqnarray}
is invertible for each $h:(U'^\s,T'^\s,\gamma')\to (X^\s,T^\s,\gamma)$ in $CRYS^\s_{et}(X^\s/\Sigma_k)$ (see Rem. \ref{remreal} \ref{remrealii}). By \cite{Be1} I Cor. 1.7.2, there is a unique divided power structure $\gamma_1'$ producing the upper left hand corner in the following cartesian square of $CRYS^\s_{et}(X^\s/\Sigma_k)$. $$\xymatrix{(U'^\s,T'^\s\times_{T^\s,p_0}T^{\s(1)},\gamma_1')\ar[r]^-{h}\ar[d]^{p_0}&(X^\s,T^{\s(1)},\gamma_1)\ar[d]^{p_0} \\ (U'^\s,T'^\s,\gamma')\ar[r]^-h& (X^\s,T^\s,\gamma)}$$
We may thus compute the source of (\ref{ftcrys1}) as follows $$\begin{array}{rcl} h^*(f_{T^\s,*}M)_{T^\s}&= &h^*\lambda_{T^\s,*}f_{T^\s}^{-1}f_{T^\s,*}M\\ &\simeq &
h^*\lambda_{T^\s,*}p_{0,*}p_1^*M\\ &\simeq &h^*p_{0,*}\lambda_{T^{\s(1)},*}p_1^*M \\ & \simeq & h^* p_{0,*} M_{T^{\s(1)}} \hspace{2cm} \hbox{(see Rem. \ref{remreal} \ref{remrealii}, \ref{remrealiv})} \\ &\simeq & p_{0,*} h^*M_{T^{\s(1)}} \end{array}$$
Here the last base change isomorphism is due to the fact that the vertical arrows denoted $p_0$ are affine on the underlying schemes and induce equivalences $\smash{(T'^\s\times_{T^\s}T^{\s(1)})_{et}}\simeq \smash{T'^\s_{et}}$, $\smash{T^{\s(1)}_{et}\simeq T^\s_{et}}$ (the main point is the isomorphism  $h^*p_{0,*}\O\simeq p_{0,*}h^*\O$ which follows from \cite{EGA2} Cor. 1.5.2 and the last statement of Lem. \ref{acycqcohsch} \ref{acycqcohschii}).
The target of (\ref{ftcrys1}) on the other hand may be computed as follows.
$$\begin{array}{rcl} (f_{T^\s,*}M)_{T'^\s}&\simeq & \lambda_{T'^\s,*}f_{T'^\s}^*f_{T^\s,*}M \\
&\simeq & \lambda_{T'^\s,*}p_{0,*}p_1^*M\\
&\simeq & \smash{p_{0,*} \lambda_{T'^\s\times_{T^\s,p_0}T^{\s(1)},*}p_1^*M}\\
&\simeq & \smash{p_{0,*} M_{T'^\s\times_{T^\s,p_0} T^{\s(1)}}} \end{array}$$
The morphism (\ref{ftcrys1}) then identifies with the one obtained from the natural one $$h^*M_{T^{\s(1)}} \to M_{T'^\s\times_{T^\s,p_0}T^{\s(1)}}$$ and we may conclude using the assumption on $M$ (see Rem. \ref{remreal} \ref{remrealii}, \ref{remrealiv}).
\begin{flushright}$\square$\end{flushright}

\begin{lem} \label{compLft}
\debrom \item \label{compLfti} The following diagram is canonically pseudo-commutative: $$\xymatrix{Mod(T^\s_{et},\O)\ar[d]_{nat}\ar[rr]_-{\ref{lemft}}^-{f_{T^\s,*}\lambda_{T^\s}^*}&&\cal Crys((X^\s/\Sigma_k)_{CRYS^\s,et},\O)\ar[d]^{(-)_{T^\s}}_{\ref{crystals}} \\
Hdp(T^\s)\ar[rr]^-{\smash{L_{T^\s}}}&&\nabla\hbox{-}Mod(T^\s)}$$
\item \label{compLftii} If $M$ is a crystal on $\smash{((X^\s/\Sigma_k)_{CRYS^\s,et},\O)}$ then $\smash{f_{T^\s}^*M}$ satisfies the condition of  Lem. \ref{lemft}. Consider the adjunction morphism \begin{eqnarray}\label{adjft}\xymatrix{M\ar[r]&f_{T^\s,*}f_{T^\s}^*M}\end{eqnarray} as a morphism in $\cal Crys\smash{((X^\s/\Sigma_k)_{CRYS^\s,et},\O)}$. Via the equivalence of Prop. \ref{crystals}, (\ref{adjft}) translates into a morphism of $\nabla\hbox{-}Mod(T^\s)$ \begin{eqnarray}\label{defaugT}\xymatrix{M_{T^\s}\ar[r]^-{aug}&\smash{L_{T^\s}}(M_{T^\s})}\end{eqnarray} which may be described as the following composed morphism in $Mod(\smash{T^\s_{et},\O})$:
    $$\xymatrix{M_{T^\s}\ar[r]^-{id\otimes d_0}&M_{T^\s}\otimes \cal P^{\s(1)}_{T}\ar[r]_{\sim}^-{\varepsilon_{M_{T^\s}}^{-1}}&\cal P^{\s(1)}_{T}\otimes M_{T^\s}}$$
\finrom
\end{lem}
Proof. \ref{compLfti} One checks that for a module $M$ over $(T^\s_{et},\O)$, the $\O$-modules $\smash{L_{T^\s}}(M)$ and
$\smash{\lambda_{T^\s,*}f_{T^\s}^*f_{T^\s,*}\lambda_{T^\s}^*M}$ both identify naturally with $\smash{p_{0,*}p_1^*M}$ and that the resulting
isomorphism is compatible with connections. The proof of \ref{compLftii} is left to the reader.
\begin{flushright}$\square$\end{flushright}


\begin{lem} \label{uL} \debrom \item \label{uLi} Consider the morphism of topoi $\iota:\smash{X^\s_{et}}\to \smash{T^\s_{et}}$. There is a natural isomorphism of rings in $\smash{T^\s_{et}}$:
$$\iota_*\O^{crys}\simeq \O^{d=0}.$$
\item \label{uLii} The following diagram is canonically pseudo-commutative: $$\xymatrix{Hdp(T^\s)\ar[r]^-{\smash{L_{T^\s}}}\ar[d]^-{(\ref{forHdpMod})}&\nabla\hbox{-}Mod(T^\s)\ar[r]_-\sim^-{(\ref{crystals})}&
\cal Crys((X^\s/\Sigma_k)_{CRYS^\s,et},\O)\ar[d]^{u_*}\\
Mod(T^\s_{et},\O^{d=0})\ar[rr]^{\iota^{-1}}&&Mod(X^\s_{et},\O^{crys}_k)}$$ \finrom
\end{lem}
Proof. This follows from Lem. \ref{lemreal2}
\ref{lemreal2ii} and Lem. \ref{compLft} \ref{compLfti}. \begin{flushright}$\square$\end{flushright}

\para \label{compdR}
We are now in a position to prove the expected analogue of Thm. 6.12 and Thm. 7.1 of \cite{BO} . Here, we
restrict ourselves to the small \'etale crystalline site. Let
$L$ denote the following composed functor:

\begin{eqnarray}\label{defLcrys}\xymatrix{Hdp(T^\s)\ar[r]_-{L_{T^\s}}^-{\ref{defLT}}&\nabla\hbox{-}Mod(T^\s)
\ar[r]^-{\ref{crystals}}_-\sim &\cal Crys((X^\s/\Sigma_k)_{crys,et},\O)\ar @{<-} @<+0pt> `u[l] `[ll]_-{L} [ll]}\end{eqnarray}

\begin{prop} \label{dRlocal} Consider a crystal $M$ on $\smash{(X^\s/\Sigma_k)_{crys,et}}$.
\debrom \item  \label{dRlocali} The morphism  (\ref{defaugT})
induces a quasi-isomorphism $$\xymatrix{M\ar[r]^-{aug}&L(\Omega_{T^\s}^\bullet(M_{T^\s}))}$$ in the category of complexes of modules of
$(\smash{(X^\s/\Sigma_k)_{crys,et}},\O)$.

\item \label{dRlocalii} The quasi-isomorphism \ref{dRlocali} induces an isomorphism
$$\xymatrix{Ru_*M\ar[r]^-\sim&\Omega^\bullet_{T^\s}(M_{T^\s})}$$
    in the derived category of modules on $(\smash{X^\s_{et},\O_k^{crys}})$.
\finrom
\end{prop}
Proof. \ref{dRlocali} We begin with the case $M=\O$. We already know from Lem. \ref{defaug} that $\smash{L(\Omega_{T^\s})}$ is a resolution of $\O$ via $d_0$
in the Abelian category $\smash{\cal Crys((X^\s/\Sigma_k)_{crys,et},\O)}$. To see that it is in fact a resolution in
$\smash{Mod((X^\s/\Sigma_k)_{crys,et},\O)}$, it suffices to notice that $$\xymatrix{{\O_{T^\s}}\ar[r]&L_{T^\s}(\Omega^\bullet_{T^\s})}$$ remains a
quasi-isomorphism when pulled back to $(T'_{et},\O)$ via an arbitrary morphism $g:\smash{T'^\s}\rightarrow \smash{T^\s}$ of the crystalline site (indeed the modules
$\O$ and $\smash{L_{T^\s}\Omega^q_{T^\s}}$ are flat). The case of an arbitrary crystal $M$ follows from the isomorphism $$M\otimes
L(\Omega^\bullet_{T^\s})\simeq L(\Omega_{T^\s}^\bullet(M_{T^\s}))\hbox{ in $Mod((X^\s/\Sigma_k)_{crys,et},\O)$}$$ (use Lem. \ref{isoLtens}) by flatness of
the modules $\smash{L(\Omega^q_{T^\s})}$.

\ref{dRlocalii} According to Lem. \ref{uL} \ref{uLii}, we have a canonical isomorphism $$u_*L(\Omega_{T^\s}^\bullet(M_{T^\s}))\simeq
\Omega_{T^\s}^\bullet(M_{T^\s})$$ in the category of complexes of modules on $(X^\s_{et},\O_k^{crys})$. The result will thus follow from the
following \petit

\noindent \emph{Claim.} The module $L(N)$ of $\smash{((X^\s/\Sigma_k)_{crys,et},\O)}$ is $u_*$-acyclic for any $N$ in $Mod(\smash{T^\s_{et}},\O)$.
\petit

Let us prove the claim. Let $T'^\s$ denote an arbitrary object in $\smash{crys_{et}(X^\s/\Sigma_k)}$ and form the product $T''^\s=T'^\s\times^{crys}T^\s$ computed in the crystalline site. If $p_0:T''^\s\to T'^\s$ and $p_1:T''^\s\to T^\s$ denote the canonical projections, we have the following natural isomorphisms: $$\begin{array}{rcl}(-)_{T'^\s}\circ f_{T^\s,*}&\simeq & \lambda_{T'^\s,*}f_{T'^\s}^*f_{T,*}\\
& \simeq &\lambda_{T'^\s,*}p_{0,*}p_1^* \\
&\simeq & p_{0,*}\lambda_{T''^\s,*}p_1^*.\end{array}$$ Since $p_0:\smash{T''^\s_{et}}\to \smash{T'^\s_{et}}$ is an equivalence, these isomorphisms show that the functor $\smash{f_{T^\s,*}}$ is exact. Whence isomorphisms $$\begin{array}{rcll}Ru_{*}L(N)&\simeq &Ru_{*}(f_{T^\s,*}\lambda_{T^\s}^*N)\\
&\simeq &R(u_*f_{T^\s,*})\lambda_{T^\s}^*N\end{array}$$ and the claim follows since $\smash{u_*f_{T^\s,*}}$ is isomorphic to
$\iota^{-1}\smash{\lambda_{T^\s,*}}$ (Lem. \ref{lemreal2} \ref{lemreal2ii}) hence exact.
\begin{flushright}$\square$\end{flushright}

\para \label{chL} Letting $k$ vary in Sect. \ref{appnablamod} and \ref{crysdRcoh} causes no difficulty and we leave this to the reader. Funtoriality with respect to the closed immersion $X^\s\to Y^\s$ fixed at the beginning of Sect. \ref{appnablamod} is less straightforward and so we briefly give some explanations.

Recall from Lem. \ref{functnabla} that the category $\nabla\hbox{-}Mod(-)$ is naturally a contravariant pseudo-functor in the argument $X^\s\to Y^\s$. The category $Hdp(-)$ on the other hand is naturally a covariant pseudo-functor with respect to $X^\s\to Y^\s$, as explained in the following lemma.
\begin{lem} \label{functHdp} Consider a morphism
$f=(f_X,f_Y):(X^\s\rightarrow Y^\s)\rightarrow (X'^\s\rightarrow Y'^\s)$ between closed immersions as in Sect. \ref{intronabla} and let
$f_T:T^\s\rightarrow T'^\s$ denote the morphism obtained by forming logarithmic divided power envelopes. There is a natural functor
$$f_*:Hdp(T^\s)\to Hdp(T'^\s)$$ sending a module $N$ to $f_{T,*}N$. \end{lem} Proof. We have to specify the effect of $f_*$ on morphisms. We define
the image of a morphism $g:N_1\to N_2$ in $Hdp(T^\s)$ (ie. $g:\smash{\cal P^{\s(1)}_T}\otimes N_{1}\rightarrow N_{2}$ in $Mod\smash{(T^\s_{et},\O)}$)
as the composed morphism $\smash{\cal P^{\s(1)}_{T'}}\otimes \smash{f_{T,*}}N_{1}\rightarrow \smash{f_{T,*}}(\smash{\cal P^{\s(1)}_T}\otimes
N_{1})\rightarrow \smash{f_{T,*}}N_{2}$ (the first arrow is the base change morphism). One checks without difficulty that this defines indeed a functor.
\begin{flushright}$\square$\end{flushright}

\begin{lem} \label{functpoincare} Let $f=(f_X,f_Y)$ be as in Lem. \ref{functHdp}.
\debrom \item \label{functpoincarei} 
There is a canonical \emph{base change morphism} $$ch_f:f^*L_{T'^\s}(\Omega^\bullet_{T'^\s}(M'))\to L_{T^\s}(\Omega^\bullet_{T^\s}(f^*M))$$ in
$\nabla$-$Mod(T^\s)$ which is functorial with respect to $M'$ in $\nabla$-$Mod(T'^\s)$. The collection of the $ch_f$'s satisfies the composition
constraint $ch_g\circ g^*ch_f=ch_{fg}$.

\item \label{functpoincareii} The base change morphism renders the following square of $\nabla$-$Mod(T^\s)$ commutative: \begin{eqnarray}\label{functLdR}\xymatrix{f^*M'\ar[r]^-{aug}\ar[rd]_-{f^*aug}& L_{T^\s}(\Omega^\bullet_{T^\s}(f^*M))\\
    & \ar[u]_{ch_f}f^*L_{T'^\s}(\Omega^\bullet_{T'^\s}(M))}\end{eqnarray}
\finrom
\end{lem}
Proof. Let us only explain the construction of $ch_f$. First, we notice that there is a natural base change morphism
\begin{eqnarray}\label{functlinearization}f^*L_{T'^\s}(f_*N)\rightarrow L_{T^\s}(N)\end{eqnarray} in $\nabla$-$Mod(T^\s)$ for $N$ in $Hdp(T^\s)$. Next,
it follows from the description of $f^*$ in Lem. \ref{functnabla}, that there is a natural base change morphism
\begin{eqnarray}\smash{f_{T}^*\Omega_{T'^\s}}(M')\rightarrow \smash{\Omega_{T^\s}}(f^*M')\end{eqnarray} for $M'$ in $\nabla$-$Mod(T'^\s)$.  From
there, one easily deduces a morphism   \begin{eqnarray}\label{functdRcomplex}\Omega^\bullet_{T'^\s}(M')\rightarrow
f_*\Omega^\bullet_{T^\s}(f^*M')\end{eqnarray} of complexes in $Hdp(T'^\s)$. The morphism $ch_f$ is obtained by combining (\ref{functlinearization})
and (\ref{functdRcomplex}).
\begin{flushright}$\square$\end{flushright}





\subsection{Exactness for crystals} \label{epotcoc} ~~ \\

The purpose of this section is to discuss some exactness properties of the categories of crystals over $(X^\s/\Sigma_.)_{crys,et}$ or $(X^\s/\Sigma_\infty)_{crys,et}$ and of the realization functors attached to local embeddings.

\para We use the terminology of Sect. \ref{deffullyexact}.


\begin{lem} \label{mexactcrys} Let $X^\s$ in $\cal Sch^\s/\Sigma_1$, $\Sigma=\Sigma_.$ or $\Sigma_\infty$ and $prop\subset \{qcoh,lf, lfft,norm\}$. The full subcategory $\cal Crys_{prop}((X^\s/\Sigma)_{crys,et},\O)$
of $Mod((X^\s/\Sigma)_{crys,et},\O)$ is closed by extensions. The induced exact structure will be denoted $e_M$ (\emph{exactness of modules}). Concretely speaking, a short sequence of crystals is thus $e_M$-exact if and only if each one of its realizations is exact.
\end{lem}
Proof. This is clear from the exactness of the realization functors (Lem. \ref{lemreal} \ref{lemrealii}).
\begin{flushright}$\square$\end{flushright}

In practice the exact structure $e_M$ is not very useful. In the subsequent paragraphs we will review another (weaker) canonical exact structure.

\para We begin with a preliminary result concerning the exactness of the pullback functor for modules with connections.

\begin{defn} \label{defflatlift} Consider a morphism $f_X:X^\s\to X'^\s$ between fine separated $\Sigma_1$-log schemes. We say that \emph{$f_X$ admits flat liftings to local embeddings} if there exist $(U^\s/X^\s,Y^\s,\iota)$ in $Emb^\s(X^\s)$, $(U'^\s/X'^\s,Y'^\s,\iota')$ in $Emb^\s(X'^\s)$ and a morphism $f=(f_X,f_U,f_Y)$ between them which lifts $f_X$ locally (Def. \ref{defemb1} \ref{defemb1iii}) and which is such that each one of the morphisms $Y_k\to Y'_k$ underlying $f_Y$ is flat.
\end{defn}

\begin{lem} \label{lemflatlift1} If $X'^\s$ is locally embeddable and $f_X:X^\s\to X'^\s$ is a strict \'etale morphism between separated $\Sigma_1$-log schemes then $f_X$ admits flat liftings to local embeddings.
 \end{lem}

Proof. Replacing $X'^\s$ by an \'etale covering if necessary, we may assume given a global embedding $(X'^\s,Y'^\s,\iota')$. We conclude using \cite{SGA1} I, Prop. 8.1.
\begin{flushright}$\square$\end{flushright}

Consider a global embedding $Y^\s$ of $X^\s$ with logarithmic divided power envelope denoted $T^\s$. In this situation, we may define the category of module with quasi-nilpotent integrable connection over $T_k^\s$ as in Def. \ref{defnablamod}. Let us denote  $\nabla\hbox{-}Mod^{pd}_.(X^\s,Y^\s)$ the category obtained by letting $k$ vary (ie. the category of sections of the cofibered category over $\N$ corresponding to the pseudo-functor defined by $k\mapsto \nabla\hbox{-}Mod^{pd}(X^\s,Y^\s_k)$ and Lem. \ref{functnabla}). This category is clearly Abelian.

\begin{lem}
\label{lemflatlift2} Consider a morphism $f=(f_X,f_Y):(X^\s,Y^\s,\iota)\to (X'^\s,Y'^\s,\iota')$ in $Emb^{\s,glob}$. If $f_X$ admits flat liftings to local embeddings then the pullback functor $$f^*:\nabla\hbox{-}Mod^{pd}_.(X'^\s,Y'^\s) \to \nabla\hbox{-}Mod^{pd}_.(X^\s,Y^\s)$$ described in Lem. \ref{functnabla} is exact.

\end{lem}
Proof. Denote $T^\s$ and $T'^\s$ the respective logarithmic divided power envelope of $\iota$ and $\iota'$. Since $et(X^\s)\simeq et(T^\s)$ and $et(X'^\s)\simeq et(T'^\s)$ we find that the problem is \'etale local on $X'$ and $X$. We may thus assume given a flat global lifting $f_1$ of $f_X$, ie. a morphism of $Emb^{\s,glob}$ of the form $f_1=(f_{X},f_{1,Y}):(X^\s,\smash{Y_1^\s})\rightarrow (X'^\s, \smash{Y'^\s_1})$ such that $f_{1,Y}$ lifts $f_X$ and is flat. Letting $\smash{Y_2^\s}$ (resp. $\smash{Y'^\s_2}$) denote the product  $\smash{Y^\s\times \smash{Y_1^\s}}$ (resp. $\smash{Y'^\s\times \smash{Y'^\s_1}}$)  computed in $Emb^{\s,glob}(X^\s)$ (resp. $Emb^{\s,glob}(X'^\s)$) and $f_2=f\times f_1$, we get from  Prop. \ref{crystals} a pseudo commutative diagram $$\xymatrix{\nabla\hbox{-}Mod_.^{pd}(X^\s,Y^\s)\ar[r]^{p_1^*}& \nabla\hbox{-}Mod_.^{pd}(X^\s,Y_2^\s)&\nabla\hbox{-}Mod_.^{pd}(X^\s,Y_1^\s)\ar[l]_{p_2^*}\\
\nabla\hbox{-}Mod_.^{pd}(X'^\s,Y'^\s)\ar[r]^{p_1^*}\ar[u]^{f^*}&
\nabla\hbox{-}Mod_.^{pd}(X'^\s,Y'^\s_2)\ar[u]^{f_2^*}&\nabla\hbox{-}Mod_.^{pd}(X'^\s,Y'^\s_1)
\ar[l]_{p_2^*}\ar[u]^{f_1^*}}$$ where horizontal arrows are
equivalence. We are thus reduced to the case where $f=f_1$, but then, the statement follows from   Lem. \ref{functnabla}, since the logarithmic divided power envelope $T_1^\s$ of
$Y_1^\s$ is flat over the logarithmic divided power envelope $T_1'^\s$ of $Y_1'^\s$.
\begin{flushright}$\square$\end{flushright}

\para \label{paraexactcrysk}  We are now in a position to discuss the case of crystals over $(X^\s/\Sigma_.)_{crys,et}$.
\begin{prop} \label{exactcrysk} Consider $X^\s$ in $\cal Sch^\s/\Sigma_1$ which is locally embeddable and let $T^\s$ denote the logarithmic divided power envelope of some local embedding $Y^\s$ for $X^\s$.
\debrom \item \label{exactcryski} The category $\cal Crys((X^\s/\Sigma_.)_{crys,et},\O)$ is Abelian and $$(-)_{T_.^\s}:\cal Crys((X^\s/\Sigma_.)_{crys,et},\O)\to Mod(T^\s_{.,et},\O)$$
is an exact faithful functor.
\item \label{exactcryskii} Crystals with quasi-coherent realizations form a fully Abelian subcategory \break $\cal Crys_{qcoh}((X^\s/\Sigma_.)_{crys,et},\O)$ of  $\cal Crys((X^\s/\Sigma_.)_{crys,et},\O)$.
\item \label{exactcryskiii} If $prop\subset \{norm, qcoh, lf, lfft\}$ then     $\cal Crys_{prop}((X^\s/\Sigma_.)_{crys,et},\O)$ is closed by extensions in the Abelian category $\cal Crys((X^\s/\Sigma_.)_{crys,et},\O)$. The induced exact structure will be denoted $e$ (\emph{exactness of crystals}).
\finrom
\end{prop}
Proof. \ref{exactcryski} Let $\Delta$ denote the category whose objects are $[\nu]$, $0\le\nu \le 2$ and whose arrows are the strictly increasing maps. Consider the diagram of type $\Delta^{op}$ whose vertex $\smash{Y_{[\nu]}^\s}=(\smash{U_{[\nu]}^\s},\smash{Y_{[\nu]}^\s})$ is the $\nu+1$-th fold product of $Y^\s$ computed in the category $Emb^\s(X^\s)$. It follows from Prop. \ref{crystals} and Lem. \ref{lemcrystal2} \ref{lemcrystal2iii} that $\cal Crys((X^\s/\Sigma)_{crys,et},\O)$ is equivalent to the category $\Gamma_{cart}(\cal F/\Delta)$ of cartesian sections of the cofibered category $\cal F$ over $\Delta$ corresponding to $[\nu]\mapsto \cal F([\nu]):= \smash{\nabla\hbox{-}Mod_.^{pd}(Y^\s_{[\nu]})}$, $(f:[\nu]\to [\nu'])\mapsto f^*$. The result now follows, by exactness of the functors $f^*$ (Lem. \ref{lemflatlift1} and Lem. \ref{lemflatlift2}).

\ref{exactcryskii} Using the crystal condition, we find easily that the conditions of having quasi-coherent realizations can be tested by realization on the given $T^\s_.$. The result follows by exactness of the realization functor $\smash{(-)_{T^\s_.}}$ of \ref{exactcryski}. The proof of \ref{exactcryskiii} is similar.
\begin{flushright}$\square$\end{flushright}

\begin{rem} \label{remexactcrysk} \debrom
\item \label{remexactcryski}
Concretely, a short sequence of crystals $\cal E_.:0\to M_{1,.}\to M_{2,.}\to M_{3,.}\to 0$ is $e_M$-exact (resp. $e$-exact) if all realizations of the $\cal E_k$'s are exact (resp. if the realization of $\cal E_k$ at $\smash{T^\s_k}$ is exact for each $k$, $T^\s$ being fixed as in  Prop. \ref{exactcrysk}).
The exact structure $e_M$ is thus stronger than $e$. This may be reformulated by saying that $e$ is not induced by the inclusion of the category of crystals into the category of modules over $((X^\s/\Sigma_.)_{crys,et},\O)$. Instead it is induced by the inclusion into the category of modules over a suitable variant of the restricted crystalline site of \cite{Be1} IV, 2.
\item \label{remexactcryskii} If a cofibered category has a cleavage by functors commuting to direct limits, then its category of sections has the same exactness properties than its fibers. It thus follows from Prop. \ref{exactcrysk} \ref{exactcryskii} that the category  $\cal Crys((X^\s/\Sigma_.)_{crys,et},\O)$ is Abelian as well when $X^\s$ is a diagram whose vertices are locally embeddable (no condition is needed on the edges here).
\finrom
\end{rem}

\para \label{paraexactcrysinfty} Let us turn to the case $k=\infty$. As a preliminary observation, let us notice that the equivalence of Rem. \ref{remqcohpadic} \ref{remqcohpadiciii} endows $Mod_{qcoh}(T^\s_{et},\O)$ with a canonical exact structure which will be denoted $e$.

\begin{prop} \label{exactcrysinfty} Let $X^\s$, $Y^\s$, $T^\s$ as in Prop. \ref{exactcrysk}.
\debrom

\item \label{exactcrysinftyi} The fully faithful functor  $$l^{-1}:\cal Crys((X^\s/\Sigma_\infty)_{crys,et},\O)\to \cal Crys((X^\s/\Sigma_.)_{crys,et},\O)$$ induces an exact structure on its source which will be denoted $e$ (\emph{exactness of crystals}).

\item \label{exactcrysinftyii} If $prop\in \{qcoh, lf, lfft\}$ then $\smash{\cal Crys_{prop}((X^\s/\Sigma_\infty)_{crys,et},\O)}$ is a fully $e$-exact subcategory of $\smash{\cal Crys((X^\s/\Sigma_\infty)_{crys,et},\O)}$.

\item \label{exactcrysinftyiii} Consider the \emph{realization functor} $\smash{(-)_{T^\s}}:Mod(\smash{(X^\s/\Sigma_\infty)_{crys,et}},\O)\to Mod(\smash{T^\s_{et}},\O)$ defined by the formula $\smash{M_{T^\s}}=\smash{l_*(l^{-1}M)_{T^\s_.}}$. The induced functor $$(-)_{T^\s}:\cal Crys_{qcoh}((X^\s/\Sigma_{\infty})_{crys,et},\O)\to Mod_{qcoh}(T^\s_{et},\O)$$ 
is faithful, $e$-exact and reflects $e$-exactness.
\finrom
\end{prop}
Proof. \ref{exactcrysinftyi} This follows from the $crys$ variant of the equivalence of Rem. \ref{lemlcrys} \ref{lemlcrysi} together with the case $prop=norm$ in Prop. \ref{exactcrysk} \ref{exactcryskiii}.

\ref{exactcrysinftyii} This follows from Prop. \ref{exactcrysk} \ref{exactcryski}
thanks to the fact that the functor  $$l^{-1}:\cal Crys_{prop}((X^\s/\Sigma_\infty)_{crys,et},\O)\to \cal Crys_{prop,norm}((X^\s/\Sigma_.)_{crys,et},\O)$$ is an equivalence (Rem. \ref{lemlcrys} \ref{lemlcrysi}, \ref{lemlcrysii}, \ref{lemlcrysiii}).

\ref{exactcrysinftyiii} We have a pseudo-commutative square as follows:
$$\xymatrix{\cal Crys_{qcoh}((X^\s/\Sigma_{\infty})_{crys,et},\O)\ar[d]^-{(-)_{T^\s}}\ar[r]^-{l^{-1}}&\cal Crys_{qcoh}((X^\s/\Sigma_.)_{crys,et},\O)\ar[d]^-{(-)_{T^\s_.}}\\
Mod_{qcoh}(T^\s_{et},\O)\ar[r]^-{l^*}&Mod(T^\s_{.,et},\O)}$$
Indeed, if $M$ is a crystal with quasi-coherent realizations in $((X^\s/\Sigma_{\infty})_{crys,et},\O)$ then $\smash{(l^{-1}M)_{T^\s_.}}$ is normalized and quasi-coherent, hence satisfies $\smash{l^*l_*(l^{-1}M)_{T^\s_.}\simeq (l^{-1}M)_{T^\s_.}}$. Now, the statement follows formally from the fact that the horizontal arrows and the right vertical one are faithful $e$-exact and reflect $e$-exactness (Prop. \ref{exactcrysk} \ref{exactcryski}).
\begin{flushright}$\square$\end{flushright}

\para \label{paraexactcryslf} If one is interested only in crystals with locally free realizations the situation is simpler.

\begin{prop} \label{exactcryslf} The exact structures $e_M$ and $e$ coincide on $\cal Crys_{lf}((X^\s/\Sigma_{\infty})_{crys,et},\O)$ and $\cal Crys_{lf}((X^\s/\Sigma_.)_{crys,et},\O)$.
\end{prop}
Proof. In both cases, one is reduced to prove that a sequence of locally free crystals $0\to M'\to M\to M''\to 0$ in $(\smash{(X^\s/\Sigma_k)_{crys,et},\O})$ which is $e$-exact is $e_M$-exact as well. Let $T^\s$ be as in Prop. \ref{exactcrysk} and consider an arbitrary $T'^\s$ in $crys_{et}(X^\s/\Sigma_k)$. Replacing $T'^\s$ with a covering if necessary, we may always assume that there exists a morphism $\smash{g:T'^\s\to T^\s_k}$. Consider the following commutative diagram $$\xymatrix{0\ar[r]&M'_{T'^\s}\ar[r]&M_{T'^\s}\ar[r]&M''_{T'^\s}\ar[r]&0\\
0\ar[r]&g^*M'_{T_k^\s}\ar[u]\ar[r]&g^*M_{T_k^\s}\ar[u]\ar[r]&g^*M''_{T_k^\s}\ar[u]\ar[r]&0}$$
The vertical arrows are invertible by the crystal condition and the bottom line is exact by local freeness of $M''$. The top horizontal line is thus exact as well and we are done by the last statement of Lem. \ref{mexactcrys}.
\begin{flushright}$\square$\end{flushright}

%
%
%
%
%
%


\begin{rem} \label{remexactcryskiii} Statements Lem. \ref{lemflatlift2}, Prop. \ref{exactcrysk} and Prop. \ref{exactcryslf} have obvious counterparts over $\Sigma_k$ if $k<\infty$. \end{rem}

\subsection{Twists by effective logarithmic divisors} \label{tbeld}

\para \label{twlogtop} We recall the basics of twisting in a general setting. Let $(E,A)$ be a ringed topos and assume given an \emph{integral log structure} $\alpha:M_E\to A$ (ie. $M_E$ is an integral monoid of $E$, $\alpha$ is a morphism of monoids, $M_E^\times:=\alpha^{-1}A^\times$ is sent isomorphically to $A^\times$ by $\alpha$ and the natural morphism of $M_E$ into its fraction group $M_E^{gp}$ is monomorphic).

\begin{defn} \label{deftwlogtop} Let $(E,A,M_E,\alpha)$ as above.

\debrom

\item \label{deftwlogtopi} An \emph{effective log divisor} $h$ is a section of the monoid $M_E/M_E^\times$.

\item \label{deftwlogtopii} If $h$ is an effective log divisor, the associated line bundle $$A(-h):=P_h\wedge^{M_E^\times} A$$ is defined as follows: $P_h$ is the $\smash{M_E^\times}$-torsor determined by $h$ (ie. $P_h:=q^{-1}\{h\}$  where $q:M_E\to \smash{M_E/M_E^\times}$)  and $A(-h)$ is the quotient of $P_h\times A$ by the action of $\smash{M_E^\times}$ given by $m.(p,a)=(pm^{-1},\alpha(m)a)$ endowed with the addition $(p,a)+(pm,b)=(p,a+\alpha(m)b)$ and the action of $A$ given by $a.(p,b)=(p,ab)$.

\item \label{deftwlogtopiii} If $M$ is a module of $(E,A)$ and $h$ an effective log divisor, we set $$M(-h):=M\otimes_A A(-h)$$
\finrom
\end{defn}

We make some remarks about these definitions:

- if $\tilde h\in \G(U,M_E)$ is a preimage of $h$ for some $U$ in $E$ then it induces a trivialization $A(-h)_{|U}\simeq A_{|U}$,  $(m\tilde h,a)\mapsto ma$.

- the \emph{twisting functor} $(-h):Mod(E,A)\to Mod(E,A)$ is exact since $A(-h)$ is locally isomorphic to $A$.

- the \emph{twisted module} $M(-h)$ can as be described as a contracted product as well: $M(-h)\simeq P_h\wedge^{M_E^\times} M$.

\begin{lem} \label{lemtwlogtop}
\debrom
\item \label{lemtwlogtopi} Let $h\in \G(E,M_E/M_E^\times)$ and assume given a factorization $h=h_1+h_2$. Then we have canonical isomorphism  and morphism of $A$-modules as follows: $$A(-h_1)\otimes_A A(-h_2)\buildrel{\sim}\over\rightarrow A(-h)\hbox{\hspace{3cm}}nat:A(-h)\rightarrow A(-h_1)$$
    \item \label{lemtwlogtopii} If  $f:(E',M_{E'}\to A')\to (E,M_E\to A)$ is a morphism of integral log ringed topoi, $h\in \G(E,M_E/M_E^\times)$ and $h'=f^{-1}h$, there are canonical isomorphisms
    $$f^*(M(-h))\simeq (f^*M)(-h')\hspace{1cm}\hbox{and} \hspace{1.5cm}
(f_*M')(-h)\simeq f_*(M'(-h'))$$ which are functorial with respect to $M\in Mod(E,A)$ and $M'\in Mod(E',A')$.
\finrom
\end{lem}
Proof. \ref{lemtwlogtopi} The isomorphism is given by the formula $(m_1,a_1)\otimes (m_2,a_2)\mapsto (m_1m_2,a_1a_2)$. The second morphism is induced
by the morphism $$A(-h_2)=P_{h_2}\wedge^{M_E^\times} A\to A\wedge^{A^\times} A\to A$$ where the first arrow is induced by $\alpha$ and the second is
multiplication in  $A$.

\ref{lemtwlogtopii} Since $A(-h)$ is locally isomorphic to $A$ both isomorphisms follow formally from the more obvious isomorphism $f^*(A(-h))\simeq
A'(-h')$.
\begin{flushright}$\square$\end{flushright}

\para We explain briefly the relation between effective logarithmic and Cartier divisors.

\begin{defn} \label{defcartierdiv} Consider a scheme $X$ and let $top$ be a topology between $zar$ and $fl$.  \debrom

\item \label{defcartierdivi}  Let $U/X$ in $top$  and $f\in \G(U,\O)$. We say that $f$ is a \emph{non zero divisor} if multiplication by $f$  on the structural ring of $U_{top}$ is monomorphic. Let us denote $S_{top}(U)\subset \G(U,\O)$ the monoid of non zero divisors so that $S:U/X\mapsto S_{top}(U)$ is a sheaf of monoids on $top(X)$.

\item \label{defcartierdivii} The \emph{$top$ sheaf of effective Cartier divisors on $X$} is the quotient sheaf $S_{top}/\Gm$. A \emph{$top$ effective Cartier divisor on $X$} is a global section of this sheaf.
\finrom
\end{defn}

For $top=zar$ this definition agrees with the one explained in \cite{Kl}.

\begin{lem} \label{epsiloncardiv} Let $\epsilon:X_{top}\to X_{zar}$ denote the natural morphism. There is a natural isomorphism $\epsilon_*(S_{top}/\Gm)\simeq S_{zar}/\Gm$. The notion of $top$ effective Cartier divisors is in particular independent of $top$.
\end{lem}
Proof. Note first that by flatness of the morphism  $\epsilon:(X_{top},\O)\to (X_{zar},\O)$, we have $\epsilon_*S_{top}=S_{zar}$. The point is thus to prove that $\epsilon_*S_{top}/\epsilon_*\Gm\simeq \epsilon_*(S_{top}/\Gm)$. If $S_{top}$ was a group, then we could conclude directly from the fact that $\Gm$ is $\epsilon_*$-acyclic. We may still conclude however thanks to the following fact (whose proof is straightforward and left to the reader).

\emph{Claim.} Consider a pretopology $(\C,Cov)$. If $G$ is a sheaf of groups acting freely on a sheaf of sets $F$, then for any $U$ in $\C$ there is a monomorphism $F(U)/G(U)\to F/G(U)$ whose image coincides with the preimage of zero under the boundary map $F/G(U)\to \check{H}^1(U,G)$.
\begin{flushright}$\square$\end{flushright}


\begin{defn} \label{flatcharts} Consider a fine log scheme $X^\s$ over $\Sigma_k$.
\debrom

\item \label{flatchartsi} A \emph{flat fine chart} for $X^\s$ over $\Sigma_k$ is a strict flat morphism \break $X^\s\rightarrow (Spec(\Z/p^k[P]),P)$ with $P$ a finitely generated integral monoid.

\item \label{flatchartsii} We say that $X^\s/\Sigma_k$ has \emph{local flat fine charts} if there is a surjective strict \'etale family $\smash{U_\lambda^\s\to X^\s}$ where each $\smash{U_\lambda^\s/\Sigma_k}$ has a flat fine chart.
    \finrom
\end{defn}

Let us point out that in virtue of Lem. \ref{lemrp} \ref{lemrpvi}, a finite $p$-basis over $\Sigma_k$ automatically induces a flat fine chart over $\Sigma_k$.

\begin{lem} \label{lemflatcharts1} Assume that $X^\s$ has local flat fine charts over $\Sigma_k$.

\debrom
\item \label{lemflatcharts1i} The logarithmic structure $\alpha:M_X\to \O$ on $(X_{et},\O)$ induces monomorphisms $M_X\hookrightarrow S_{et}$ and $M_X/\Gm\hookrightarrow S_{et}/\Gm$.

\item \label{lemflatcharts1ii} Assume that $X^\s$ has a $p$-basis $(\underline s,\underline t)$. Then $M_X$ identifies with the submonoid sheaf of $S_{et}$ generated by $\Gm$ and the $t_i$'s. In particular if $z_i:V_X(t_i)\to X$ denote the inclusion morphism of closed subscheme defined by $t_i$ (Rem. \ref{remqcohsch1} \ref{remqcohsch1iii}) then $M_X/\Gm\simeq \prod_i z_{i,*}\N$.
\finrom
\end{lem}
Proof. The proof of \ref{lemflatcharts1i} and the first statement of \ref{lemflatcharts1ii} is given in \cite{Ka2} 1.5 (3), (1.10). For the proof of the last statement we may assume that $X$ is an affine space as in Rem. \ref{examplesofpb} \ref{examplesofpbi}. We have to describe the sub monoid sheaf of $S_{et}/\Gm$ generated by the $t_i$'s. Consider an \'etale $X$-scheme $U$ and denote $z_{i,U}:V_U(t_i)\to U$. Since $U$ is regular, the Zariski sheaf of its Cartier divisors identifies with its sheaf of Weil divisors and we find that the sub monoid Zariski sheaf generated by the $t_i$'s identifies with $\prod_i z_{i,U,*}\N$. The resulting collection of Zariski sheaves for variable $U$ clearly satisfies the \'etale descent condition and corresponds to the \'etale sheaf $\prod_i z_{i,*}\N$ as claimed.
\begin{flushright}$\square$\end{flushright}

\begin{lemdefn}\label{suppcenter} Consider an arbitrary fine log scheme $X^\s$. The set of points where the log structure is non trivial is closed. The corresponding reduced closed subscheme will be denoted $Supp(M_X)$. The ideal generated by the image of $M_X$ into $\O$ defines a closed subscheme which we denote $Center(M_X)$.
\end{lemdefn}

Proof. The fact that the non trivial locus of the log structure is closed results e.g. from  \cite{Sh} Prop. 2.3.1. Finally, one can use Rem. \ref{remqcohsch1} \ref{remqcohsch1ii} to show that the ideal generated by the image of $M_X\to \O$ is quasi-coherent. 
\begin{flushright}$\square$\end{flushright}

\begin{rem} \label{indepZzpb}
 If $X^\s$ has a $p$-basis $(\underline s,\underline t)$ then Lem. \ref{lemflatcharts1} implies that the closed subschemes defined in Lem.-Def. \ref{suppcenter} have the following explicit description: $Supp(M_X)=V_X(\prod t_i)$ and $Center(M_X)=V_X(t_1,\dots, t_e)$ (see Rem. \ref{remqcohsch1} \ref{remqcohsch1ii}). Note in particular that the latter is reduced too.
\end{rem}

\para \label{twistTOP}
Let $X^\s$ in $\cal Sch_{p}$. The collection of log structures $\smash{M_{X_k}\to \O}$ defines an integral log structure  $\smash{M_{X_.}\to \O}$ on the ringed topos $(\smash{X_{.,et}^\s},\O)$. Note that by exactness of the direct image of $\iota_{1,k}:X_1^\s\to X_k^\s$ the effective  log divisors on $X_k$ are in bijection with those of $X_1$. In particular $\G(X_{.,et},M_{X_.}/M_{X_.}^\times)\simeq \G(X_1,M_{X_1}/\Gm)$.

\begin{lem} \label{injtwist1} Let $h_.\in \G(X_{.,et},M_{X_.}/M_{X_.}^\times)$.
\debrom

\item \label{injtwist1i} The twisting functor can be computed componentwise: $\iota_k^{-1}(M_.(-h_.))\simeq M_k(-h_k)$.
\item \label{injtwist1ii} The image of the canonical morphism $nat_.:\O(-h_.)\rightarrow \O$ is the ideal generated by $\alpha(h_.)$. If $X^\s_k$ has flat fine charts over $\Sigma_k$ then $nat_k$ is a monomorphism in the category
$Mod(\smash{X^\s_{k,et}},\O)$.

\item \label{injtwist1iii}  The twisting functor $(-h_.):Mod(\smash{X^\s_{.,et}},\O)\to Mod(\smash{X^\s_{.,et}},\O)$ is exact and preserves the properties $norm$, $qcoh$, $lf$, $lfft$.
\finrom
\end{lem} Proof. Statement \ref{injtwist1i} is a special case of Lem. \ref{lemtwlogtop} \ref{lemtwlogtopii}.  In \ref{injtwist1ii}, the first assertion follows from the fact that \'etale locally, $nat_k$ can be described as multiplication on $\O$ by $\smash{\alpha(\tilde h_k)}$ with $\tilde h_k$ a preimage of $h$ giving a trivialization as explained after Def. \ref{deftwlogtop}. The second assertion follows from Lem. \ref{lemflatcharts1} \ref{lemflatcharts1i}. Statement \ref{injtwist1iii} is clear from the fact that the properties in question are \'etale local.
\begin{flushright}$\square$\end{flushright}

Consider $X^\s$ in $\cal Sch/\Sigma_1$. The collection of log structures  $M_T\rightarrow \O$ for $T^\s=(U^\s,T^\s,\iota)$ in $Crys_{et}(X^\s/\Sigma_k)$ defines an integral log structure $\smash{M_{X/\Sigma}\rightarrow \O_{X/\Sigma}}$ on the ringed topos $\smash{((X^\s/\Sigma)_{crys,et},\O)}$ ($\Sigma=\Sigma_.$ or $\Sigma_k$, $1\le k\le \infty$). Note that an effective log divisor $h$ of $X^\s$ induces canonically an effective log divisor $h_{T^\s}$ of $T^\s$ since $\iota$ is an exact closed immersion. Whence in particular a bijection $\G(\smash{(X^\s/\Sigma)_{crys,et}},\smash{M_{X/\Sigma}/\O_{X/\Sigma}^\times})\simeq\G(X,M_X/\Gm)$.


\begin{lem} \label{injtwist2} Let $X^\s$, $T^\s$ and $h$ be as above. \debrom

\item \label{injtwist2i} The twisting functor $(-h):Mod((X/\Sigma)_{crys,et},\O)\to Mod((X/\Sigma)_{crys,et},\O)$ is pseudo-compatible with the realization functor $\smash{(-)_{T^\s_.}}$, ie. $\smash{M(-h)_{T^\s_.}}\simeq \smash{M_{T^\s_.}(-h_{T_.})}$. In particular it preserves crystals.

\item \label{injtwist2ii} If $X^\s$ is locally embeddable then $nat:\O(-h)\to \O$ is a monomorphism in the category of crystals.
\finrom
\end{lem}
Proof. Statement \ref{injtwist2i} follows from Lem. \ref{lemtwlogtop} \ref{lemtwlogtopii}. Statement \ref{injtwist2ii} is a consequence of Prop. \ref{exactcrysk} \ref{exactcryski} and Lem. \ref{injtwist1} \ref{injtwist1ii}.
\begin{flushright}$\square$\end{flushright}



\section{Twisted syntomic complexes for Dieudonn\'e crystals} \label{sectscfdc} ~~ \\

The purpose of this chapter is to define and relate several variants of the syntomic complex functors announced in (\ref{introSsyn}) and (\ref{introSet}). On the \'etale site we carry on three parallel constructions of a twisted syntomic complex (Def. \ref{deffil} and Def.  \ref{defSonephi}) which are shown to be isomorphic (Lem. \ref{lemcompS}). Then we give two constructions on the syntomic site in the case of trivial log structures and we relate them to the previous ones (Def. \ref{deffilsynmod}, Def. \ref{defSonephisyn} and Prop. \ref{compsynet}).

\subsection{The category of $(1,\varphi)$-modules} \label{tco1pm} ~~ \\

Roughly speaking (each variant of) the syntomic functor will be defined as the mapping fiber of $1-\varphi$ where $1$, $\varphi$ are morphisms to be defined. The purpose of this section is to summarize some basic facts about  the category keeping track of $1$ and $\varphi$ themselves (rather than their mapping fiber).

\para \label{paraonephi} It will be convenient to look at the category of $(1,\varphi)$-modules as a category of modules of a suitable ringed topos (namely $(E^{1,\varphi},(A,A,id_A,F))$, see Lem. \ref{1phimod} \ref{1phimodii} below) arising from a ringed topos $(E,A)$ equipped with an endomorphism $F$ of $A$. The following constructions will be used e.g. when $\smash{E=X_{et}^\N}$, $\smash{T^\s_{[.],.,et}}$, $\smash{X_{syn}^\N}$ etc. in Sect. \ref{tscotes} and Sect. \ref{scotss}.

\begin{lem} \label{1phitop} Let $E$ denote a topos and consider the category $E^{1,\varphi}$ of $E$-valued presheaves on the category
$1\leftleftarrows 0$. Explicitly:

- an object is a quadruple $X=(X_1,X_0,1_X,\varphi_X)$ where $X_1$ and $X_0$ are objects of $E$ and $1_X$ and  $\varphi_X$ are morphisms $X_1\rightarrow X_0$ in
$E$.

- a morphism $X\rightarrow Y$ is a couple of morphisms $X_1\rightarrow Y_1$, $X_0\rightarrow Y_0$ which are compatible with $1_X,$ and $1_Y$ as well as with
$\varphi_X$ and $\varphi_Y$.

\debrom
\item \label{1phitopi} The category $E^{1,\varphi}$ is a topos.

\item \label{1phitopii} There are two morphisms of topoi  $0, 1:E\to E^{1,\varphi}$ such that $1^{-1}(X_1,X_0,1,\varphi)=X_1$, $1_*X=(X,e_E,1,1)$, $0^{-1}(X_1,X_0,1,\varphi)=X_0$ and $0_*X=(X\times X,X,p_1,p_2)$. Here $e_E$ denotes a final object in $E$.

\item \label{1phitopiii}  There is a morphism of topoi $\varpi:E^{1,\varphi}\rightarrow E$ such that $\varpi^{-1}X=(X,X,id_X,id_X)$ and
$\varpi_*(X_1,X_0,1,\varphi)=Ker(1,\varphi)$. One has in particular  $\varpi \circ 1\simeq \varpi \circ 0\simeq id_E$ canonically.

\item \label{1phitopiv} The topos $E^{1,\varphi}$ and the morphisms $1$, $0$ and  $\varpi$ are canonically pseudo-functorial with respect to $E$.

\item \label{1phitopv} The functors $1^{-1}$, $0^{-1}$ and $\varpi^{-1}$ have left adjoints, $1_0$, $0_0$ and $\varpi_0$, satisfying $1_0(X)=(X,X\sqcup X,i_1,i_2)$, $0_0(X)=(\emptyset_E,X,0,0)$ and $\varpi_0(X_1,X_0,1,\varphi)=Coker(1,\varphi)$. Here $\emptyset_E$ denotes an initial object in $E$.
\finrom
\end{lem}
Proof. \ref{1phitopi} is clear from Giraud's criterion and the remaining statements are immediate.  \begin{flushright}$\square$\end{flushright}

The following lemma explains how the topos $E^{1,\varphi}$ gives rise to $(1,\varphi)$-modules.

\begin{lem} \label{1phimod} Assume that $(E,A)$ is a ringed topos, and that $A$ is equipped with a given endomorphism $F$. Let $(F)$ denote the endomorphism of the ringed topos $(E,A)$ whose underlying morphism of topoi is the identity of $E$ and whose underlying morphism of rings is $F$. Let us consider the category $Mod^{1,\varphi}(E,(A,F))$ (or
simply $Mod^{1,\varphi}(E,A)$ if there is no ambiguity about $F$) where:

- an object is a quadruple $M=(M_1,M_0,1_M,\varphi_M)$ where $M_1$ and $M_0$ are modules of $(E,A)$ and $1_M:M_1\rightarrow M_0$ and $\varphi_M:M_1\rightarrow
(F)_*M_0$ are morphisms in $Mod(E,A)$.

- a morphism $M\rightarrow N$ is a couple of morphisms $M_1\rightarrow N_1$, $M_0\rightarrow N_0$ in $Mod(E,A)$ which are compatible with $1_M$ and $1_N$ as well as with $\varphi_M$ and $\varphi_N$.

\debrom \item \label{1phimodi} The category $Mod^{1,\varphi}(E,A)$ is canonically isomorphic to the category of modules of the ringed topos $(E^{1,\varphi},(A,A,id_A,F))$.

\item
\label{1phimodii} The morphisms $1$, $0$ and $\varpi$ induce tautological morphisms of ringed topoi
$$\xymatrix{(E,A)\ar@<+.7ex>[r]^-1\ar@<-.7ex>[r]_-0&(E^{1,\varphi},(A,A,id_A,F))\ar[r]^-\varpi&(E,A^{F=1})}$$
where $A^{F=1}:=\varpi_*A$ denotes the subring of fixed points of the endomorphism $F$.

\item
\label{1phimodiii}  The pullback functors for modules $1^*$, $0^*$ and $\varpi^*$ have left adjoints, $1_!$, $0_!$ and $\varpi_!$, satisfying $1_!(M)=(M,M\oplus
(F)^*M,i_1,i_2)$, $0_!(M)=(0,M,0,0)$ and $\varpi_!(M_1,M_0,1,\varphi)=Coker(1-\varphi)$. The second one is exact and the other two are left derivable. There
is moreover an exact sequence in  $Mod^{1,\varphi}(E,A)$ $$\xymatrix{0\ar[r]&0_!M_0\ar[r]& M\ar[r]& 1_*M_1\ar[r]&0}$$ which is functorial with respect
to $M=(M_1,M_0,1_M,\varphi_M)$ in $Mod^{1,\varphi}(E,A)$.


\item
\label{1phimodiv}  In $D^+(E,A^{F=1})$ one has a canonical distinguished triangle
$$\xymatrix{R\varpi_*M\ar[r]&M_1\ar[rr]^{1_M-\varphi_M}&&M_0\ar[r]&R\varpi_*M[1]}$$ which is functorial
with respect to  $M=(M_1,M_0,1_M,\varphi_M)$ in $D^+(Mod^{1,\varphi}(E,A))$.

\finrom

\end{lem}
Proof.  In \ref{1phimodi}, it suffices to observe that a module action $$(A,A,id_A,F)\times (M_1,M_0,1_M,\varphi_M)\to (M_1,M_0,1_M,\varphi_M)$$ is equivalent to the data of module actions $A\times M_i\to M_i$ for which $1_M$ is linear and $\varphi_M$ is $F$-linear. The proof of \ref{1phimodii} and  \ref{1phimodiii} is routine. Let us explain \ref{1phimodiv}. By scalar restriction to $(A^{F=1},A^{F=1},\smash{id_{A^{F=1}}},\smash{id_{A^{F=1}}})$, we can assume that $A^{F=1}=A$ (recall that $R\varpi_*$ commutes to restricting scalars (\cite{SGA4-II} V, Prop. 5.1, (2)). But then, the result will follow from Lem. \ref{lemMF2} applied to the case where the functors $F_1$,
$F_2$, $F$ and $F_0$ respectively send $M=(M_1,M_0,1_M,\varphi_M)$ to $M_1$, $M_0$, $1_M-\varphi_M$ and $\varpi_*M$, once checked that $1_M-\varphi_M$ is
surjective as soon as $M$ is an injective object of $Mod^{1,\varphi}(E,A)$. Now, the adjunction map between the forgetful functor $M\mapsto (M_1,M_0)$ and its right adjoint provides a monomorphism $M\hookrightarrow M'$ where $M'=(M_1\times M_0\times M_0,M_0,p_2,p_3)$ and $p_2$, $p_3$ denote the second and third projections. This monomorphism is split when $M$ is injective and it is thus sufficient to notice that $\smash{1_{M'}-\varphi_{M'}}=p_2-p_3$ is surjective.
%
\begin{flushright}$\square$\end{flushright}


\subsection{Local embeddings with Frobenius and cohomological descent} \label{lewfacd} ~~ \\

In this section, we define several categories of semi-simplicial local embeddings with additional structures which are adapted to the constructions we have in mind. We also prove that they satisfy the assumptions of Lem. \ref{CBdescent2}, which provides us a canonical way of gluing local constructions.

\para \label{prelimemb}
As a motivation, let us begin with two applications of cohomological descent in the small \'etale crystalline topos. Let $X^\s$ in $\cal Sch^\s/\Sigma_1$  and consider a semi-simplicial $p$-adic $dp$-thickening  $\smash{T^\s_{[.]}}=(\smash{U^\s_{[.]}/X^\s,T^\s_{[.]}})$. Denoting again 
$\smash{T^\s_{[.]}}=limind_k\, \smash{T^\s_{[.],k}}$ the associated semi-simplicial crystalline sheaf, we have the following pseudo-commutative diagram (Lem. \ref{lemreal2} \ref{lemreal2ii}):

\begin{eqnarray}\label{descentdiaget} \xymatrix{((X^\s/\Sigma_.)_{crys,et},\O)\ar[d]_u& ((U^\s_{[.]}/\Sigma_.)_{crys,et},\O)\ar[l]_-{f_{U^\s_{[.],.}}} \ar[d]_u
&((X^\s/\Sigma_.)_{crys,et},\O)/T^\s_{[.],.}\ar[l]_-{f_{T^\s_{[.],.}/U^\s_{[.],.}}}\ar[d]^{\lambda_{T^\s_{[.],.}}}  \ar @{->} @<+0pt> `u[l] `[ll]_-{f_{T^\s_{[.],.}}} [ll]\\
(X^{\s,\N}_{et},\O^{crys}_.)&
\ar[l]_-{f_{U^\s_{[.],.}}}(U^{\s,\N}_{[.],et},\O^{crys}_.)\ar[r]^-{\iota}_\sim &(T^{\s}_{[.],.,et},\O^{crys}_.)
\ar @{->} @<+0pt> `d[l] `[ll]^-{f_{T^\s_{[.],.}}} [ll]
}\end{eqnarray}

\noindent Here, we have used the simplified notation $\O^{crys}_.$ instead of $\iota_*\O^{crys}_.$, where $\iota:\smash{U^\s_{[.]}\to T^\s_{[.]}}$. In the following, we refer to \cite{SGA4-II} V,
Sect. 7.3.1.1 for the notion of \emph{hypercovering} in a topos.

%

\begin{prop} \label{cdOmega} Let $X^\s$ in $\cal Sch/\Sigma_1$.

\debrom \item \label{cdOmegai} Assume that $\smash{T^\s_{[.]}}$ is a hypercovering of  $(X^\s/\Sigma_\infty)_{crys,et}$. For any $M_.$ in  \break
$\smash{D^+((X^\s/\Sigma_.)_{crys,et},\O)}$, the adjunction morphism $M_.\rightarrow \smash{Rf_{T^\s_{[.],.},*}f_{T^\s_{[.],.}}^{-1}M_.}$ is invertible and
induces an isomorphism
\begin{eqnarray}\label{cohdesccrys} Ru_*M_.\simeq Rf_{T^\s_{[.],.},*}\lambda_{T_{[.],.},*}f_{T^\s_{[.],.}}^{-1} M_.&\hspace{.5cm} \hbox{in
$D^+(X^{\s,\N}_{et},\O_.^{crys})$}\end{eqnarray}

\item \label{cdOmegaii} Assume that $\smash{U_{[.]}^\s}$ is a hypercovering of $\smash{X^\s_{et}}$ and that $\smash{T_{[.]}^\s}$ is the logarithmic divided power
envelope of some semi-simplicial object  $\smash{(U_{[.]}^\s,Y_{[.]}^\s)}$ of $Emb^{\s,glob}$. For any module $M_.$ on
$\smash{((X^\s/\Sigma_.)_{crys,et},\O)}$ the adjunction morphism $M_.\rightarrow \smash{Rf_{U^\s_{[.],.},*}f_{U^\s_{[.],.}}^{-1}M_.}$ is invertible. If
$M_.$ is a crystal this isomorphism together with the Poincar\'e resolution (see Prop. \ref{dRlocal} \ref{dRlocali} and Lem. \ref{functpoincare}) induce
\begin{eqnarray}\label{cohdescetdR} Ru_*M_.\simeq Rf_{T^\s_{[.],.},*}\Omega^\bullet_{T^\s_{[.],.}}(M_.)&\hspace{.5cm} \hbox{in
$D^+(X^{\s,\N}_{et},\O_.^{crys})$}\end{eqnarray}
\finrom
\end{prop}
\noindent Proof. The first statement of \ref{cdOmegai} (resp. \ref{cdOmegaii}) is an easy consequence of \cite{SGA4-II} V, Thm. 7.3.2 applied to the
hypercovering $\smash{T_{[.],k}^\s}$ (resp. $\smash{i_*U^\s_{[.]}}$) of the topos $(X^\s/\Sigma_k)_{crys,et}$. The second statement of \ref{cdOmegai}
follows by pseudo-commutativity of (\ref{descentdiaget}). Let us now explain the second statement of \ref{cdOmegaii}. Using Lem. \ref{functpoincare}, we get a
semi-simplicial version of the Poincar\'e Lemma (Prop. \ref{dRlocal} \ref{dRlocali}): $$f_{U^\s_{[.]}}^{-1}M_.\simeq
L(\Omega^\bullet_{T^\s_{[.],.}}(M_{.,|U^\s_{[.]}}))\hbox{ in $D^+((U^\s_{[.]}/\Sigma_.)_{crys,et},\O)$}$$
Using Lem. \ref{uL} \ref{uLii} as in the proof of Prop. \ref{dRlocal} \ref{dRlocali} we get
$$Ru_*f_{U^\s_{[.]}}^{-1}M_.\simeq \iota^{-1}\Omega^\bullet_{T^\s_{[.]}}(M_{.,|U^\s_{[.]}})\hbox{ in $D^+((U^{\s,\N}_{[.],et},\O_.^{crys}))$}$$ and the result
follows by (\ref{descentdiaget}).
\begin{flushright}$\square$\end{flushright}

In the statement of the
proposition, we have used the simplified notation $$\smash{\Omega^\bullet_{T^\s_{[.],.}}(M_.)}:=\Omega^\bullet_{T^\s_{[.],.}}(M_{.,|U^\s_{[.]}})$$
and we will continue to do so in the rest of the text.


%

\begin{rem} \label{remcomp}  Assume that $\smash{T^\s_{[.]}}$ satisfies simultaneously the conditions of Prop. \ref{cdOmega} \ref{cdOmegai} and \ref{cdOmegaii} and let $M_.$ denote a crystal. Consider the following morphisms in $D^+(T^\s_{[.],.,et},\O_.^{crys})$:
\begin{eqnarray}\label{maps123}\xymatrix{(M_.)_{T^\s_{[.],.}}\ar[r]_-{\sim}^-{\ref{dRlocal} \ref{dRlocali}}&
(L\Omega^\bullet_{T^\s_{[.],.}}(M_.))_{T^\s_{[.],.}}&\ar[l]\iota_*u_*
L\Omega^\bullet_{T^\s_{[.],.}}(M_.)&\ar[l]_-{\ref{uL} \ref{uLii}}^-{\sim} \Omega^\bullet_{T^\s_{[.],.}}(M_.)}\end{eqnarray}
where the middle arrow is deduced from  Lem. \ref{lemreal2} \ref{lemreal2ii} as follows: \begin{eqnarray}\label{middle123}\iota_*u_*\to
\iota_*u_*f_{T^\s_{[.],.},*}f_{T^\s_{[.],.}}^*\to \lambda_{T^\s_{[.],.},*}f_{T^\s_{[.],.}}^*.\end{eqnarray} On the one hand it follows easily from (\ref{descentdiaget}) and the proof of Prop. \ref{cdOmega} \ref{cdOmegai} and \ref{cdOmegaii} that the middle arrow of (\ref{maps123}) in question becomes
invertible after applying $$\smash{Rf_{T^\s_{[.],.},*}}:D^+(\smash{T^\s_{[.],.,et},\O_.^{crys}})\rightarrow
D^+(\smash{(X^{\s,\N}_{et},\O^{crys}_.)}).$$ Stupid truncation on the other hand provides  a morphism from the last
term of (\ref{maps123}) to the first one. As explained in \cite{BD}, one can prove that the image of the latter morphism under
$\smash{Rf_{T^\s,[.],.,*}}$ is also invertible. We have not tried to compare these two isomorphisms.
\end{rem}

In view of constructing twisted syntomic complexes above fine log schemes or diagrams of such it will be convenient to introduce the following
variants of the categories of embeddings defined in Def. \ref{defemb1} and Def. \ref{defemb1var}.



%
%
%

\begin{defn} \label{defemb2} Consider $*\subset \{\s,lfpb\}$.

\debrom \item \label{defemb2i} Let $Emb^*_{F}$ (resp. $\smash{Emb_{F}^{*,glob}}$) denote the following category. An object
$(U^\s/X^\s,Y^\s,\iota,\smash{\tilde F})$ is an object of $Emb^*$ (resp. $Emb^{*,glob}$) together with a Frobenius lift $\smash{\tilde F}$ on $Y^\s$. A morphism is a morphism
of $Emb^*$ which is compatible with Frobenius lifts.


\item \label{defemb2ii} Let $\smash{HR_{F}^{*,et}}$ (resp. $\smash{HR_{F}^{*,crys}}$) denote the category of semi-simplicial objects of the category $\smash{Emb^*_{F}}$ of
the form  $\smash{Y^\s_{[.]}}=\smash{(U^\s_{[.]}/X^\s, Y^\s_{[.]},\iota_{[.]},\tilde F_{[.]})}$ for some constant semi-simplicial object $X^\s$ of
$\smash{\cal Sch^*}/\Sigma_1$ and such that $\smash{U^\s_{[.]}}$ (resp. the logarithmic divided power envelope $\smash{T_{[.]}^\s}$ of $\iota_{[.]}$) is
a hypercovering of the topos $\smash{X^\s_{et}}$ (resp. $(X^\s/\Sigma_\infty)_{crys,et}$).

\item \label{defemb2iii} Let $\smash{\cal Sch_{div}^\s}/\Sigma_1$ denote the category of couples $(X^\s,h)$ where $X^\s=(X,M_X)$ is in $\cal Sch^\s/\Sigma_1$ and $h\in
\G(X,M_X/\Gm)$ is an effective log divisor on $X^\s$. A morphism $(X^\s,h)\rightarrow (X'^\s,h')$ is a morphism $f:X^\s\rightarrow X'^\s$ which is \emph{compatible with
log divisors} in the sense that $h$ divides $f^*h'$ in $\G(X,M_X/\Gm)$. If $\cal C=Emb^*_F$ or $\smash{Emb_F^{*,glob}}$ (resp. $\smash{HR_F^{*,et}}$ or
$\smash{HR_F^{*,crys}}$), we view it as a category above $\smash{\cal Sch^\s}/\Sigma_1$ via the forgetful functor  $Y^\s\mapsto X^\s$ (resp.
$\smash{Y^\s_{[.]}}\mapsto X^\s$). We further denote $\cal C_{div}$ the category of couples $(Y^\s,h)$ (resp. $(\smash{Y^\s_{[.]}},h)$) where $h$ is an effective log
divisor on $X^\s$ and we view it as a category above $\smash{\cal Sch_{div}^\s}/\Sigma_1$ via $(Y^\s,h)\mapsto (X^\s,h)$ (resp.
$(\smash{Y^\s_{[.]}},h)\mapsto (X^\s,h)$).

\item \label{defemb2iv} Let $X^\s$ (resp. $(X^\s,h)$) denote a diagram in the category $\cal Sch^\s/\Sigma_1$ (resp. $\cal Sch^\s_{div}/\Sigma_1$). If $\cal C=\smash{Emb^*_F}$,
$\smash{Emb_F^{*,glob}}$, $\smash{HR_F^{*,et}}$ or $\smash{HR_F^{*,crys}}$ we denote simply $\cal C(X^\s)$  the fiber at $X^\s$ of
$\smash{Diag(\C)}/\smash{Diag(\cal Sch^*/\Sigma_1)}$ and $\cal C(X^\s,h)\simeq \cal C(X^\s)$ the fiber at $(X^\s,h)$ of
$\smash{Diag(\C_{div})}/\smash{Diag(\cal Sch_{div}^*}/\Sigma_1)$. \finrom
\end{defn}

The following lemma is a complement to Lem. \ref{lememb1}.

\begin{lem} \label{lememb2} Let $X^\s$, $X'^\s$ denote diagrams of $\cal Sch^\s/\Sigma_1$ whose vertices are separated.

\debrom
\item \label{lememb2i} 
The product of two objects is representable in the categories $\smash{Emb^\s_{F}(X^\s)}$, $Emb_F^{\s,glob}(X^\s)$,  $\smash{HR_{F}^{\s,et}}(X^\s)$ and $\smash{HR_{F}^{\s,crys}}(X^\s)$. This product is moreover functorial with respect to $X^\s$.

\item \label{lememb2ii} The category $\smash{HR_{F}^{\s,et}}$ contains $\smash{HR_{F}^{\s,crys}}$. These categories have a non empty fiber above $X^\s$ if and only if
$\smash{Emb^\s_{F}}(X^\s)$ is non empty.

\item \label{lememb2iii} Consider a morphism $f_X:X^\s\rightarrow X'^\s$  and assume that $\smash{Emb^\s_F(X^\s)}$ is non empty. Any given object above $X'^\s$ in the category $\smash{Diag(Emb^\s_{F})}$ (resp. $Diag(\smash{HR_{F}^{\s,et}})$, resp. $Diag(\smash{HR_{F}^{\s,crys}})$)  is the target of a
morphism above $f_X$.

\item \label{lememb2iv} Assume that the type $\Delta$ of $X^\s$ is such that $\delta/\Delta$ (the category of arrows with source $\delta$) is finite for any $\delta$ in $\Delta$ (this is the case if $\Delta$ is finite or if $\Delta=\N$). The category $\smash{Emb^\s_{F}(X^\s)}$ is non empty if $\smash{Emb^\s(X^\s_\delta)}$ is non empty for each $\delta\in \Delta$ (this is in particular the case if the $\smash{X^\s_\delta}$'s have local finite $p$-bases). 
\item Statements \ref{lememb2i} - \ref{lememb2iv} hold without $\s$'s.     
\finrom \end{lem}

Proof. \ref{lememb2i} As in Lem. \ref{lememb1}, we can form products componentwise (note
that the property of being a hypercovering is stable by products by \cite{SGA4-II} V, Lem. 7.3.4).

\ref{lememb2ii} Since the property of being a hypercovering is expressed in terms of fiber products,  it is  preserved by the inverse image functor of any morphism of topoi. The first claim thus follows from the fact that $\smash{U_{[.]}^\s}=i^{-1}\smash{T_{[.]}^\s}$. To get the second
one, it suffices to observe that the coskeleton construction makes sense in $\smash{Emb_{F}^\s(X^\s)}$ thanks to \ref{lememb2i} and gives a right adjoint of the functor
$\smash{HR_{F}^{\s,crys}(X^\s)}\to \smash{Emb^\s_{F}(X^\s)}$, $\smash{Y^\s_{[.]}\mapsto Y^\s_{[0]}}$.

\ref{lememb2iii} The proof of Lem. \ref{lememb1} \ref{lememb1ii} works here as well.

\ref{lememb2iv} Assume given local embeddings $\smash{(U_\delta^\s/X_\delta^\s,Y_\delta^\s,\iota_\delta )}$ for each $\delta$ in $\Delta$. Replacing
$\smash{Y_\delta^\s}$ by a strict \'etale  $\smash{Y_\delta^\s}$-$p$-adic log scheme if necessary,  and similarly with  $\smash{U_\delta^\s}$, we may
assume that each $\smash{Y_\delta^\s}$ has $p$-bases and also that $Y_\delta$ (and thus $U_\delta$) is separated. Choosing  $p$-bases  induces
Frobenius lifts $\tilde F_\delta$. The collection $\smash{(U_\delta^\s/X_\delta^\s,Y_\delta^\s,\iota_\delta,\tilde F_\delta)}$ for $\delta\in \Delta$
can be seen as a diagram of type $|\Delta|$ (the discrete subcategory underlying $\Delta$). Now we observe that the forgetful functor
$for:\smash{(\cal Sch_{p}^\s)^\Delta}\rightarrow \smash{(\cal Sch_{p}^\s)^{|\Delta|}}$ has a right adjoint, say $cofor$. Explicitly $cofor$ sends a
collection $\smash{(Z^\s_\delta)}$ to the diagram $\delta\mapsto \smash{\prod Z^\s_{\delta'}}$ where the product is indexed on the finite set of objects
$\delta\rightarrow \delta'$ in $\delta/\Delta$.
 Note that, thanks to the finiteness assumption on the $\delta/\Delta$'s, the vertices of $cofor(Z_\delta)$ have local finite $p$-bases if and only each $Z_\delta$ has.
Let $|X|$ denote the diagram of type $|\Delta|$ underlying $X$. The forgetful functor $\smash{for_{X^\s}:(\cal Sch^\s)^\Delta/X^\s\to (\cal Sch^\s)^{|\Delta|}/|X|}$ admits similarly a right adjoint $\smash{cofor_{X^\s}}$ sending $\smash{(Z^\s_\delta)}$ to $\delta\mapsto \smash{\prod_{X^\s_\delta} (X^\s_{\delta}\times_{X^\s_{\delta'}} Z_{\delta'})}$. By separatedness of the $X_\delta$'s, we find that the natural morphism $\smash{cofor_{X^\s} (Z_\delta)}\to cofor (Z_\delta)$ is a closed immersion (\cite{EGA1} Chap. 1, Prop. 5.4.2). The vertices of $\smash{cofor_{X^\s}(Z_\delta)}$ are moreover separated over $\Sigma_1$ (resp. have local finite $p$-bases over $\Sigma_1$, resp. are \'etale above the vertices of $X$) if and only the same is true for each $Z_\delta$.
Since the Frobenius lifts $\smash{\smash{\tilde F_\delta}}$ induce a Frobenius lift on $cofor (Y^\s_\delta)$ by
functoriality, we have all the ingredients to build an object $\smash{Emb^\s_{F}(X^\s)}$.
\begin{flushright}$\square$\end{flushright}

For some purpose, it will be useful to work with exact closed immersions rather than closed immersions. The reader may consult \cite{NS} for more elaborate variants of the following lemma.

\begin{lem}\label{lememb3}  Let us denote $\smash{Emb^{\s,ex}_F}$ the full subcategory of $\smash{Emb^{\s}_F}$ formed by the objects $(U^\s/X^\s,Y^\s,\iota,\tilde F)$
where $Y^\s$ has local finite $p$-bases of the form $(\underline s,t)$ (and thus has local charts of type $\N$) and $\iota$ is an \emph{exact} closed immersion. Also let
$\smash{HR^{\s,et,ex}_{F}}$ and $\smash{HR^{\s,crys,ex}_F}$ denote the respective full subcategories of $\smash{HR^{\s,et}_{F}}$ and $\smash{HR^{\s,crys}_F}$
formed by the $\smash{Y^\s_{[.]}}$'s for which each $\smash{Y^\s_{[\nu]}}$ is in $\smash{Emb^{\s,ex}_F}$. Consider a diagram $X^\s$ of $\cal Sch^\s/\Sigma_1$ whose vertices are separated.

\debrom
\item \label{lememb3i} The categories $\smash{Emb^{\s,ex}_F(X^\s)}$, $\smash{HR^{\s,et,ex}_{F}(X^\s)}$ and $\smash{HR^{\s,crys,ex}_F(X^\s)}$ have finite non empty products.

\item \label{lememb3ii} Assume that $X_\delta^\s$ has local charts of type $\N$ and that  $\delta/\Delta$ is finite for all $\delta$ in $\Delta$. The categories  $\smash{Emb^{\s,ex}_F(X^\s)}$, $\smash{HR^{\s,et,ex}_{F}(X^\s)}$ and $\smash{HR^{\s,crys,ex}_F(X^\s)}$ are non empty if $\smash{Emb^\s(X^\s_\delta)}$ is non empty for each $\delta\in \Delta$ (this is in particular the case if the $\smash{X^\s_\delta}$'s have local finite $p$-bases).
\finrom
\end{lem}
Proof. \ref{lememb3i} Consider the full subcategory $\smash{Emb^{\s,diagonal}_F}$  of $\smash{Emb^{\s}_F}$ formed by the objects $(U^\s/X^\s,Y^\s,\iota,\tilde F)$ such that \'etale locally $Y^\s$ admits $p$-bases $(\underline s,\underline t)$ satisfying that each one of the $\iota^*t_i$'s is a chart for $U^\s$ (ie. $\N\to M_U$, $1\mapsto \iota^*t_i$ is a chart). On the one hand, $\smash{Emb^{\s,diagonal}_F}$ clearly contains $\smash{Emb^{\s,ex}_F}$. On the other hand,  the category
$\smash{Emb^{\s,diagonal}_F}(X^\s)$ is a full subcategory of $\smash{Emb^{\s}_F(X^\s)}$ which  is stable by finite non empty products, as the reader will check easily. To prove that
$\smash{Emb^{\s,ex}_F}(X^\s)$ has finite non empty products as claimed, it is thus enough to show that the inclusion functor $inc:\smash{Emb^{\s,ex}_F}\rightarrow
\smash{Emb^{\s,diagonal}_F}$ has a right adjoint $ex$ which is a $\smash{\cal Sch^\s/\Sigma_1}$-functor.
For $(U^\s/X^\s,Y^\s,\iota,\smash{F_{Y^\s}})$ in $Emb^{\s,diagonal}_F$, we set $ex(U^\s/X^\s,Y^\s,\iota,\smash{F_{Y^\s}}):=(U^\s/X^\s,\smash{\tilde Y^\s},\smash{\tilde
\iota},\smash{F_{\tilde Y^\s}})$ where $\smash{\tilde Y^\s}:=(\smash{\tilde Y^\s_k})$, $\tilde\iota:X^\s\rightarrow \smash{\tilde Y^\s}$ and $\smash{F_{\tilde Y^\s}}$ are to be described now.

Let $Z_k=Supp(\smash{M_{Y_k}})$ and $z_k=Center(\smash{M_{Y_k}})$ (see Lem.-Def. \ref{suppcenter}). Consider the blowing up $\smash{\overline Y_k}$ of $Y_k$ centered at $z_k$. We define $\smash{\tilde Y_k}$
as the open complement of the strict transform of $Z_k$ in $\smash{\overline Y_k}$. Etale locally on $Y_k$ we thus have explicitly \begin{eqnarray}\tilde
Y_k\simeq Spec_{Y_k}(\O_{Y_k}\otimes_{\Z[t_1,\dots,t_e]}\Z[t_1,y_{1,2}^{\pm 1},\dots,y_{1,e}^{\pm 1}])\end{eqnarray} where the tensor product is taken
with respect to $t_1\mapsto t_1$ and $\smash{t_i\mapsto y_{1,i}t_1}$, $i\ge 2$. The universal property of blowing ups gives a morphism $\smash{\overline
\iota:U\rightarrow \overline Y_k}$ above $\iota:U\to Y_k$.
Since each $\iota^*t_i$ is a chart for $U^\s$, there exists, for each
$i\ge 2$, a unique $u_{1,i}\in \Gm(U_i)$ satisfying $\smash{\iota^*t_i=u_{1,i}\iota^*t_1}$.
In particular, $\smash{\overline \iota}$ factors through $\smash{\tilde Y_k}$. Endowing
$\smash{\tilde Y_k}$ with the log structure induced by $\smash{Y_k^\s}$ and letting $k$ vary, we find  an exact closed immersion $\smash{\tilde \iota}:U\to \smash{\tilde
Y^\s}$, where $\smash{\tilde Y^\s}$  admits local finite $p$-bases of the form $((s_1,\dots,s_d,\smash{y_{1,2}},\dots,\smash{y_{1,e}}),t_1)$ and is the desired $p$-adic log scheme. The construction of $(U^\s/X^\s,\smash{\tilde Y^\s},\tilde
\iota)$ is clearly functorial with respect to $(U^\s/X^\s,Y^\s,\iota)$. The required Frobenius lift $\smash{F_{\tilde Y^\s}}$ may thus be defined by
functoriality and this ends the construction of the functor $ex$. We leave it to the reader to check that it is canonically right adjoint to $inc$
(one adjunction morphism is given by the structural morphism  $\smash{\tilde Y^\s}\rightarrow Y^\s$ of the blowing up and the other one is the identity).

Let us emphasize that the functor $ex$ does not affect logarithmic divided power envelopes since the blowing up is log \'etale (\cite{Ka2} Prop. 5.3, Def. 5.4).
Forming logarithmic divided power envelopes commutes in particular to finite non empty products computed in the category $\smash{Emb^{\s,ex}_F}(X^\s)$.
It follows immediately that  $\smash{HR^{\s,et,ex}_{F}(X^\s)}$ and $\smash{HR^{\s,crys,ex}_F(X^\s)}$ also have finite non empty products.

\ref{lememb3ii} Since $\smash{Emb^{\s,ex}_F(X^\s)}$ has finite non empty products which are compatible with logarithmic divided power envelopes, the coskeleton construction provides a right adjoint to $\smash{HR^{\s,crys,ex}_F(X^\s)}\to \smash{Emb^{\s,ex}_F(X^\s)}$, $\smash{Y^\s_{[.]}}\mapsto \smash{Y^\s_{[0]}}$. It is thus sufficient to prove that $Emb^{\s,ex}_F(X^\s)$ is non empty if the $X_\delta^\s$'s are locally embeddable. Start from a local embedding $(\smash{U^\s_\delta/X^\s_\delta},\smash{Y^\s_\delta},\iota_\delta)$ for each $\delta$. Applying e.g. Lem. \ref{lememb2}  \ref{lememb2iii} to an appropriate strict \'etale  surjective $\smash{X'^\s_\delta\to X_\delta^\s}$,
we may always assume that $\smash{U^\s_\delta}$ has charts of type $\N$, say $f:\smash{U_\delta^\s}\to (Spec(\Z[t]),(0))$. We may assume furthermore given a family $\smash{(Y^\s_{\delta,\lambda})_\lambda}$ and compatible isomorphisms $\smash{Y^\s_{\delta,k}}\simeq \sqcup_\lambda \smash{Y^\s_{\delta,\lambda,k}}$ such that each $\smash{Y^\s_{\delta,\lambda}}$ has $p$-bases of the form $(\underline s,\underline t)=\smash{((s_{1},\dots,s_{d_\lambda}),(t_1,\dots, t_{e_\lambda}))}$. Let $\smash{U_{\delta,\lambda}^\s:=U^\s\times_{Y_\delta^\s}Y^\s_{\delta,\lambda}}$. There is at least one index $i$ such that the image of $t_i$ in $\G(U_{\delta,\lambda},M_{U_{\delta,\lambda}}/\Gm)$ coincides with the image of $f^*t$. Consider the log structure on $\smash{Y_\lambda}$ such that  $\smash{((s_{1},\dots,s_{d_\lambda},t_1,\dots,t_{i-1},t_{i+1},\dots,t_{e_\lambda}),t_i)}$ is a $p$-basis. The resulting log structure on $Y_\delta$ defines an object of $\smash{Emb_F^{\s,ex}(X^\s_\delta)}$. We may now apply the $cofor$ construction, as in the proof of Lem. \ref{lememb2} \ref{lememb2iv}, in order to produce an object in $\smash{Emb^{\s,diagonal}_F(X^\s)}$. We conclude  using the functor $ex:\smash{Emb^{\s,diagonal}_F(X^\s)}\to \smash{Emb^{\s,ex}_F(X^\s)}$ defined in the proof of \ref{lememb3i}.
\begin{flushright}$\square$\end{flushright}

\begin{rem} As already observed in Lem. \ref{lememb2} \ref{lememb2ii}, a semi-simplicial local embedding $\smash{Y^\s_{[.]}}$ which is in $\smash{HR_F^{\s,crys}}$ is automatically in $\smash{HR_F^{\s,et}}$. This implication is obviously strict since the property of being in $\smash{HR_F^{\s,et}}$ only depends on $\smash{U^\s_{[.]}}$. In particular, it does not see log structures. To illustrate this, let us start with $X^\s$ as in Lem. \ref{lememb3} \ref{lememb3ii} and some $Y^\s$ in $\smash{Emb_F^{\s,diagonal}(X^\s)}$. The coskeleton construction provides some $\smash{Y^\s_{[.]}}$ in $\smash{HR_F^{\s,crys}(X^\s)}$. The underlying $\smash{Y_{[.]}}$ obtained by forgetting log structures is clearly in $\smash{HR^{\s,crys}(X)}$ as well. Apply now the functor $ex$ defined in the above proof. The semi-simplicial local embedding $\smash{\tilde Y_{[.]}}$ of $X$ obtained by forgetting log structures on $\smash{\tilde Y^\s_{[.]}:=ex(Y^\s_{[.]})}$ is in $\smash{HR_F^{et}(X)}$ (because $\smash{\tilde Y^\s_{[.]}}$ is in $\smash{HR_F^{\s,crys}(X^\s)}\subset \smash{HR_F^{\s,et}(X^\s)}$) but certainly not in $\smash{HR_F^{crys}(X)}$ in general (unless $X=X^\s$). Indeed, let $\smash{\tilde T^\s_{[.]}}$ denote the logarithmic divided power envelope of $\smash{\tilde Y^\s_{[.]}}$. Since the closed immersions $\smash{U_{[\nu]}^\s\to \tilde Y^\s_{[\nu]}}$ are exact, the underlying $\smash{\tilde T_{[.]}}$ coincides with the divided power envelope of $\smash{\tilde Y_{[.]}}$. Now $$\smash{Rf_{\tilde T^\s_{[.],.},*}}\O\simeq \smash{Rf_{\tilde T_{[.],.},*}}\O$$ computes the crystalline cohomology of $(X^\s/\Sigma_.)$, but certainly not that of $(X/\Sigma_.)$.
\end{rem}

\para \label{paraprelimS}
We explain how to reduce the construction of the syntomic complex to the case where a global embedding (with effective log divisor and Frobenius lift) is given. The first step is to pass from global embeddings to hypercoverings. The second step is to show the essential independence with respect to the choice of hypercoverings (see Lem. \ref{newglu} below).

\begin{defn} \label{defB0} 
\debrom \item \label{defB0i}
Let $\cal B_0^\s$ denote the full subcategory of $\smash{Diag(\cal Sch^{\s,slfpb}_{div})}$ formed by the $(X^\s,h)$'s such that $X^\s$ is locally embeddable (this is the case if the type $\Delta$ of $X^\s$ satisfies the finiteness condition of \ref{lememb2} \ref{lememb2iv}). We also denote $\cal B_0$ the full category formed by the $X$'s (with no log structure). 

\item \label{defB0ii} Let $\cal C_0^\s$ denote either one of the categories $Diag(HR^{\s,et,lfpb}_{F,div})$ or $Diag(HR^{\s,crys,lfpb}_{F,div})$.
\finrom
\end{defn}

Consider a contravariant pseudo-functor $\cal F:\smash{\cal Sch^{\s}}\to \mathfrak Cat$, as well as a ringed weakly variable topos $(\cal T,A)$ on $\smash{\cal Sch^\s}$, which is either variable or associated to a prevariable pretopology. Let us denote $\cal F^{emb}$ (resp. $\cal K^{emb}$) the contravariant pseudo-functor obtained by composing $\cal F$ (resp. $(Kom^+(\cal T(-),A),f\mapsto f^*)$) with the forgetful functor $\smash{Emb^{\s,glob,lfpb}_{F,div}}\to \cal Sch^\s$, $Y^\s=(X^\s\hookrightarrow Y^\s,F,h)\mapsto X^\s$.

Our starting point is a given  colax morphism \begin{eqnarray}\label{Semb}\cal S^{emb}:\cal F^{emb}\to \cal K^{emb}\end{eqnarray} between contravariant pseudo-functors on $\smash{Emb^{\s,glob,lfpb}_{F,div}}$. In other terms, we assume given a collection of functors $\cal S^{emb}_{Y^\s}:\cal F(X^\s)\to Kom^+(\cal T(X^\s),A)$ indexed by the $Y^\s$'s in $\smash{Emb^{\s,glob,lfpb}_{F,div}}$, as well as a collection of natural transformations $f^*\smash{\cal S^{emb}_{Y'^\s}}\to \smash{\cal S^{emb}_{Y^\s}}f^*$, indexed by the morphisms $f:Y^\s\to Y'^\s$ in $\smash{Emb^{\s,glob,lfpb}_{F,div}}$, and subject to the composition constraint. The data of $\cal S^{emb}$ is equivalent to that of a $(\smash{Emb^{\s,glob,lfpb}_{F,div}})^{op}$-functor \begin{eqnarray}\label{Sembcof}\cal S^{emb}:(\cal F^{emb})_{cof}\to (\cal K^{emb})_{cof}\end{eqnarray}

From this data, we are going to construct a $\smash{\cal B^{\s,op}_0}$-functor \begin{eqnarray}\label{Sglob}\cal S^{glob}:(\cal F^{glob})_{cof}\to (\cal D^{glob}{}')_{fib}\end{eqnarray}
where $\cal F^{glob}$ (resp. $\cal D^{glob}{}'$) is the contravariant (resp. covariant) pseudo-functor on $\smash{\cal B^{\s}_0}$ obtained by composing $\cal F^{codiag}$ (resp. $(D^+(\cal T^{codiag}(-),A),f\mapsto Rf_*)$) with the forgetful functor $\smash{\cal B^\s_0}\to Diag(\cal Sch^\s)$, $(X^\s,h)\mapsto X^\s$. This, in turn, amounts to a collection of functors $\smash{\cal S^{glob}_{X^\s}}:\cal F^{codiag}(X^\s)\to D^+(\cal T^{codiag}(X^\s),A)$, together with a collection of base change morphisms $\smash{\cal S^{glob}_{X^\s}}\to Rf_* \smash{\cal S^{glob}_{X^\s}}f^*$ indexed by the morphisms of $\smash{\cal B^\s_0}$ and subject to the composition constraint.

We need to introduce several intermediate pseudo-functors: $\cal F^{loc}$ (resp. $\cal K^{loc}$) is the contravariant pseudo-functor on $\smash{\cal C^{\s}_0}$ obtained from $\cal F^{codiag}$ (resp. $(Kom^+(\cal T^{codiag}(-),A),f\mapsto f^*)$) by composition with the functor $tot:\smash{\cal C^{\s}_0}\to Diag(\cal Sch^\s)$ sending $\smash{(U^\s_{[.]}/X^\s, Y^\s_{[.]},\iota_{[.]},\tilde F_{[.]},h)}$ to the underlying diagram $([\nu],\delta)\mapsto \smash{U^\s_{[\nu],\delta}}$ of type $Simp^{op}\times \Delta$, where $Simp$ denotes the category of ordered sets $[\nu]=\{0 \le \dots \le \nu\}$. Next $\cal K^{loc}{}'$ (resp. $\cal D^{loc}{}'$) is the covariant pseudo-functor on  $\smash{\cal C^{\s}_0}$ obtained from  $(Kom^+(\cal T^{codiag}(-),A),f\mapsto f_*)$ (resp.  $(D^+(\cal T^{codiag}(-),A),f\mapsto Rf_*)$) by composition with $tot$.  We make the following observations.

- The collection of the functors $\smash{f_{U^\s_{[.]}}^*}:\cal F^{codiag}(X^\s)\to \cal F^{codiag}(\smash{U^\s_{[.]}})$, for $\smash{Y^\s_{[.]}}=(\smash{f_{U^\s_{[.]}}}:\smash{U^\s_{[.]}}\to X^\s, \smash{Y^\s_{[.]}},\iota_{[.]},\tilde F_{[.]},h)$, running in $\smash{\cal C^\s_0}$ naturally gives rise to a pseudo-morphism \begin{eqnarray} \label{cons1} \cal F^{glob}_{|\cal C^{\s}_0}\to \cal F^{loc}\end{eqnarray} between contravariant pseudo-functors on $\smash{\cal C^{\s}_0}$.

- Applying $(-)^{codiag}$ (Lem. \ref{lemdiagcodiag} \ref{lemdiagcodiagii} and Rem. \ref{remdiagcodiag} \ref{remdiagcodiagii}) to (\ref{Semb}), then restricting via the functor $\smash{\cal C^{\s}_0}\to Diag(\smash{Emb^{\s,glob,lfpb}_{F,div}})$, sending $\smash{Y^\s_{[.]}/\Delta}$ to the underlying diagram of $\smash{Emb^{\s,glob,lfpb}_{F,div}}$ of type $Simp^{op}\times \Delta$, we get a colax morphism \begin{eqnarray} \label{cons2}\cal S^{loc}:\cal F^{loc}\to \cal K^{loc}\end{eqnarray} between contravariant functors on $\smash{\cal C^{\s}_0}$.

- There is an isomorphism of $\smash{\cal C^{\s,op}_0}$-categories (Rem. \ref{remlaxmor} \ref{remlaxmorii}) \begin{eqnarray}\label{cons3}(\cal K^{loc})_{cof}\simeq (\cal K^{loc}{}'_{fib})\end{eqnarray}

-  The compatibility of right derivation with composition (Lem.-Def. \ref{acycf} \ref{acycfvi}) yields a lax morphism \begin{eqnarray}\label{cons4}\cal K^{loc}{}'\to \cal D^{loc}{}'\end{eqnarray} between covariant pseudo-functors on $\smash{\cal C^{\s}_0}$.

- The collection of the functors $\smash{Rf_{U^\s_{[.]},*}}:D^+(\cal T^{codiag}(U^\s_{[.]}),A)\to D^+(\cal T^{codiag}(X^\s),A)$, for $\smash{Y^\s_{[.]}}$ running in $\smash{\cal C^\s_0}$, naturally gives rise to a pseudo-morphism \begin{eqnarray} \label{cons5} \cal D^{loc}{}'\to \cal D^{glob}{}'_{|\cal C^{\s}_0}\end{eqnarray} between covariant pseudo-functors on $\smash{\cal C^{\s}_0}$.

Applying $(-)_{cof}$ to (\ref{cons1}) and (\ref{cons2}), $(-)_{fib}$ to (\ref{cons3}) and (\ref{cons4}) and putting everything together, we get a $\smash{\cal C^{\s,op}_0}$-functor \begin{eqnarray} \label{cons6} \cal S_{\cal C^{\s,op}_0}: ((\cal F^{glob})_{cof})_{|\cal C^{\s,op}_0}\to  (\cal D^{glob}{}')_{fib})_{|\cal C^{\s,op}_0}\end{eqnarray}
In other terms, we have obtained a collection of functors \begin{eqnarray}\label{SglobC}\begin{array}{rcccl}\smash{\cal S_{Y^\s_{[.]}}}&:&\cal F^{codiag}(X^\s)&\to & D^+(\cal T^{codiag}(X^\s),A)\\ && \xi &\mapsto & \smash{Rf_{U^\s_{[.]},*}} \cal S^{loc}(\smash{f_{U^\s_{[.]}}^*\xi})\end{array}\end{eqnarray}
(here $\cal S^{loc}(\smash{f_{U^\s_{[.]}}^*\xi})$ is the object of $Kom^+(\cal T^{codiag}(\smash{U^\s_{[.]}},A))$ collecting the $\cal S^{emb}({f_{U^\s_{[\nu],\delta}}^*\xi})$'s and their base change morphisms for varying $([\nu],\delta)$'s), indexed by the objects $\smash{Y^\s_{[.]}}$ of $\smash{\cal C^{\s,op}_0}$, as well as a collection of base change morphisms $\smash{\cal S_{Y^\s_{[.]}}}\to Rf_*\smash{\cal S_{Y'^\s_{[.]}}}f^*$, indexed by the morphisms $f:\smash{Y'^\s_{[.]}}\to \smash{Y^\s_{[.]}}$ in $\smash{\cal C^{\s}_0}$, and subject to the composition constraint.

\begin{lem} \label{newglu} The functor (\ref{cons6}) descends to a $\smash{\cal B^{\s,op}_0}$-functor (\ref{Sglob})  if and only the base change morphism $\smash{\cal S_{Y^\s_{[.]}}}\to \smash{\cal S_{Y'^\s_{[.]}}}$ is invertible for each $X^\s$ in $\smash{\cal B^\s_0}$ and each morphism $f:\smash{Y'^\s_{[.]}}\to \smash{Y^\s_{[.]}}$ in $\cal C^\s_0(X^\s)$.
\end{lem} 
Proof. This follows from \ref{CBdescent2}, thanks to \ref{lememb2} \ref{lememb2i},  \ref{lememb2ii},  \ref{lememb2iii}. \begin{flushright}$\square$\end{flushright}

\begin{rem} If one is given a global embedding $Y^\s\in \smash{Emb^{\s,glob,lfpb}_{F,div}}$, one easily checks that there is a canonical isomorphism $\smash{\cal S^{emb}_{Y^\s}}\simeq  \smash{\cal S^{glob}_{X^\s}}$ (look at $Y^\s$ as a constant  hypercovering of $X^\s$).
\end{rem}

\grand

\subsection{The relative Frobenius and Cartier's descent for crystals} ~~ \\
\label{subsectioncartier} \label{trfacdfc}

The purpose of this section is to review the crystalline interpretation of Cartier's descent for crystals with trivial $p$-curvature in the case of local finite $p$-bases (see Prop. \ref{propcartier}) and to gather the exactness results needed for the definition of the mod $p$ Hodge filtration of Dieudonn\'e crystals. Logarithmic structures play essentially no role until Sect. \ref{cartierfactor} where we discuss elementary consequences of the Cartier equivalence on $X$ for the category of crystals on $(X^\s/\Sigma_1)$.


\para \label{sectintrocartier} Let us begin with a review of the relative Frobenius and of the inverse Cartier operator in the language of topoi.

We denote $F$ or $\smash{F_{X^\s}}$ the absolute Frobenius endomorphism of $X^\s$ in $\cal Sch^\s/\Sigma_1$. If $U^\s$ is in $TOP^\s(X^\s)$, we denote $F_0U^\s$ the object deduced from $U^\s$  by composing the structural morphism  $U^\s\to X^\s$ with $\smash{F_{X^\s}}$. The absolute Frobenius of $U^\s$ is thus a morphism $$F_{U^\s}:F_0U^\s\to U^\s \hbox{ in $TOP^\s(X^\s)$.}$$ Similarly, if $(U^\s,T^\s)$ is an object of $\smash{CRYS^\s_{top}(X^\s/\Sigma_1)}$, we denote $F_0(U^\s,T^\s)$ the object deduced from $(U^\s,T^\s)$ by  composing the structural morphism  $U^\s\to X^\s$ with  $\smash{F_{X^\s}}$. The absolute Frobenius of $U^\s$ and $T^\s$ is then a morphism $$F_{(U^\s,T^\s)}:F_0(U^\s,T^\s)\to (U^\s,T^\s)\hbox{ in  $\smash{CRYS^\s_{top}(X^\s/\Sigma_1)}$.}$$ In both case we have obtained a natural transformation of functors \begin{eqnarray} \label{defFcoc}F:F_0\to id\end{eqnarray}

\begin{defn} \label{defFrel}
\debrom \item \label{defFreli} Let $F:\smash{X^\s_{TOP^\s}}\to \smash{X^\s_{TOP^\s}}$ or $\smash{(X^\s/\Sigma_1)_{CRYS^\s,top}}\to \smash{(X^\s/\Sigma_1)_{CRYS^\s,top}}$ denote the morphism induced by $F:X^\s\to X^\s$. We define the \emph{relative Frobenius} \begin{eqnarray}\label{eqdefFreli}F^{(-/X^\s)}:F\to id\end{eqnarray} as the morphism induced by (\ref{defFcoc}) viewed as a natural transformation between cocontinuous functors.

\item \label{defFrelii} In the case $TOP$ or $top$ (resp. $CRYS$ or $crys$) we deduce \emph{relative Frobenius} morphism from \ref{defFreli} as follows \begin{eqnarray}\label{eqdefFrelii}F^{(-/X^\s)}:F=\pi Fr\to \pi r\simeq id\end{eqnarray}
    Here $\pi$ and $r$ are as in $(\ref{defpcrys})$ and the comment below.
\finrom
\end{defn}
The morphism $\smash{F^{(-/X^\s)}}$ can equivalently be seen as a natural transformation $F_*\to id$, or $id\to F^{-1}$. Note that in the case of $\smash{X^\s_{TOP}}$, the latter natural transformation (or equivalently its restriction to $TOP(X^\s)$, ie. the natural transformation from the identity to the functor ``base change by $\smash{F_{X^\s}}$) is what is more commonly called the relative Frobenius.

The relative Frobenius on the usual and $top$ crystalline topoi are compatible with each other via the functoriality isomorphism $Fi\simeq iF$ and $Fu\simeq uF$ (when $u$ exists, e.g. if $top\in \{zar,et,syn\}$ by Lem. \ref{lemiucocont}).
Recall that $ui\simeq id$, canonically, and that there is a natural morphism $nat:iu\to id$ (provided by $i_*\simeq u^{-1}$).

\begin{lemdefn} \label{defCinv} Let $top=et$ or $syn$.
\debrom \item \label{defCinvi} Consider the Frobenius morphism $F:\smash{(X^\s/\Sigma_1)_{CRYS^\s,top}}\to \smash{(X^\s/\Sigma_1)_{CRYS^\s,top}}$. There exists a unique morphism $C^{-1}$, called the \emph{inverse Cartier morphism}, factorizing $\smash{F^{(-/X^\s)}}$ as follows: $$\xymatrix{F\ar[dr]_{C^{-1}}\ar[rr]^-{F^{(-/X^\s)}}&&id\\ &iu\ar[ur]_{nat}&}$$

\item \label{defCinvii} In the case $CRYS$ or $crys$ the \emph{inverse Cartier morphism} $C^{-1}$ is defined from  \ref{defCinvi} as follows: $$\xymatrix{F=\pi Fr\ar[dr]_{C^{-1}}\ar[rr]^{F^{(-/X^\s)}}&&\pi r=id\\ &\pi iur=iu\ar[ur]_{nat}&}$$
\finrom
\end{lemdefn}
Proof. \ref{defCinvi} Recall that the morphisms $F$, $id$ and $iu$ are respectively induced by the cocontinuous functor sending $(U^\s,T^\s)$ to $F_0(U^\s,T^\s)$, $(U^\s,T^\s)$ and $(U^\s,U^\s)$.
Since these functors are continuous as well, it is not difficult to check that the claimed factorization of $\smash{F^{(-/X^\s)}}$ is equivalent to a functorial factorization $$\xymatrix{F_0(U^\s,T^\s)\ar[dr]_-{C^{-1}}\ar[rr]^-{F}&&(U^\s,T^\s)\\ &(U^\s,U^\s)\ar[ru]_{nat}&}$$ Uniqueness is clear and existence results from the divided power structure on the ideal $I_{U/T}$ of the closed immersion $U\to T$, using $x^p=p!x^{[p]}=0$.
\begin{flushright}$\square$\end{flushright}

\noindent As in the case of the relative Frobenius morphism, the inverse Cartier morphism $C^{-1}$ can equivalently be seen as a natural transformation $F_*\to i_*u_*$ or $u^{-1}i^{-1}\to F^{-1}$.
Of particular importance in Sect. \ref{tmphfotsces} will be the morphism \begin{eqnarray}\label{defCinvbis}C^{-1}:i^{-1}\to u_*F^{-1}\end{eqnarray}
which is deduced from the latter by adjunction.

The following compatibilities will be used in Sect. \ref{dopdg}.

\begin{lem}

\label{compFext} \debrom \item \label{compFexti} Consider Abelian groups $A$, $B$ in $(X^\s/\Sigma_1)_{CRYS^\s,top}$ or  $X^\s_{TOP^\s}$. The following square is commutative: $$\xymatrix{R\cal Hom(A,B)\ar[r]^-{F^{(B/X^\s)}}\ar[d]_{F^{(R\cal Hom(A,B)/X^\s)}}&R\cal Hom(A,F^{-1}B)\\ F^{-1}R\cal Hom(A,B)\ar[r]_-\sim& R\cal Hom(F^{-1}A,F^{-1}B)\ar[u]_{F^{(A/X^\s)}}}$$

\item \label{compFextii} Let $top$ be as in Lem.-Def. \ref{defCinv} and consider Abelian groups $A$, $B$ in $(X^\s/\Sigma_1)_{CRYS^\s,top}$. The following diagram is commutative:
$$\hspace{-.2cm}\xymatrix{R\cal Hom(A,B)\ar[d]_{C^{-1}}\ar[r]^-{C^{-1}}& R\cal Hom(A,i_*Ru_*F^{-1}B)\ar[r]^-{adj}_-\sim &i_*Ru_*R\cal Hom(u^{-1}i^{-1}A,F^{-1}B)\ar[d]^-{nat}\\ i_*Ru_*F^{-1}R\cal
Hom(A,B)\ar[r]^-\sim&  i_*Ru_*R\cal Hom(F^{-1}A,F^{-1}B)\ar[r]^-{F^{(A/X)}}&i_*Ru_*R\cal Hom(A,F^{-1}B)}$$

\item \label{compFextiii} If $X^\s=X$ (resp. and if the absolute Frobenius morphism of $X$ is $top$) then \ref{compFexti} and \ref{compFextii} hold \emph{verbatim} with big (resp. small) sites instead of $\s$-big sites.
\finrom
\end{lem}
Proof. It is sufficient to prove analogous compatibilities before deriving. We begin with a general fact whose  proof is left to the reader. If $f,g:E\to E'$ are any two morphisms of topoi and $\alpha:f\to g$ is a morphism between them, then the following natural square of $Ab(E')$ is commutative for any $A$ in $Ab(E')$, $B$ in $Ab(E)$:
\begin{eqnarray}\label{sqfg}\xymatrix{\cal Hom(A,f_*B)\ar[d]^\wr\ar[r]&\ar[d]^\wr\cal Hom(A,g_*B)\\
f_*\cal Hom(f^{-1}A,B)\ar[r]&g_*\cal Hom(g^{-1}A,B)}\end{eqnarray}
Let us now prove the lemma.

\ref{compFexti} Recall that the natural transformation $id\to F^{-1}$ can be deduced from $F_*\to id$ by composition as follows: $$\xymatrix{id\ar[r]& F_*F^{-1}\ar[r]&F^{-1}}$$ We thus have to prove that the exterior square of the following diagram is commutative:
$$\xymatrix{\cal Hom(A,B)\ar[r]\ar[d]&\cal Hom(A,F_*F^{-1}B)\ar[d]^\wr\ar[r]&\cal Hom(A, F^{-1}B)\ar[d]_{\parallel}\\
F_*F^{-1}\cal Hom(A,B)\ar[d]\ar[r]^-\sim &F_*\cal Hom(F^{-1}A,F^{-1}B)\ar[d]\ar[r]&\cal Hom(A,F^{-1}B)\ar @{<-} @<+2pt> `d[d] [dl]^-{F^{(A/X^{\s})}} \\
F^{-1}\cal Hom(A,B)\ar[r]^-\sim &\cal Hom(F^{-1}A,F^{-1}B)
&}$$
It suffices to prove that the interior squares are commutative. The left ones cause no difficulty and the top right one follows from (\ref{sqfg}) applied  to the morphism $\smash{F^{(-/X^\s)}}:F\to id$ (note that the bottom right triangle is tautologically commutative by definition of the bottom horizontal arrow in (\ref{sqfg})).

\ref{compFextii} Here the arrows denoted $C^{-1}$ are meant as the ones induced by the morphism $C^{-1}:id\to i_*u_*F^{-1}$. The latter can be decomposed as $$\xymatrix{id\ar[r]^{}&F_*F^{-1}\ar[rr]^{F_*\to i_*u_*}&& i_*u_*F^{-1}}$$
We thus have to prove that the exterior square of the following diagram is commutative:
$$\xymatrix{\cal Hom(A,B)\ar[r]^-{}\ar[d]_{}&\cal Hom(A,F_*F^{-1}B)\ar[d]_\wr\ar[r]^-{F_*\to i_*u_*}&\cal Hom(A,i_*u_*F^{-1}B)\ar[d]_{adj}\\
F_*F^{-1}\cal Hom(A,B)\ar[d]_{F_*\to i_*u_*}\ar[r]^-\sim&F_*\cal Hom(F^{-1}A,F^{-1}B)\ar[r]\ar[d]_{F_*\to i_*u_*}&i_*u_*\cal Hom(u^{-1}i^{-1}A,F^{-1}B)\ar[d]_{nat}\\
i_*u_*F^{-1}\cal Hom(A,B)\ar[r]_-\sim&i_*u_*\cal Hom(F^{-1}A,F^{-1}B)\ar[r]_-{F^{(A/X^\s)}}&i_*u_*\cal Hom(A,F^{-1}B)}$$
It suffices to prove that the interior squares are commutative. The left ones cause no difficulty and the top right one  is  (\ref{sqfg}) applied  to the morphism $C^{-1}:F\to iu$. Commutativity of the bottom right square results from the commutativity of the triangle $$\xymatrix{&u^{-1}i^{-1}A\ar[dl]_-{C^{-1}}\\F^{-1}A&A\ar[l]^-{F^{(A/X^\s)}}\ar[u]_-{nat}}$$

\ref{compFextiii} The previous proofs remain valid as long as $F$ is a localization morphism (so that $\cal Hom$ commutes to $F^{-1}$).
\begin{flushright}$\square$\end{flushright}

\para Let us now discuss the linear variants of $C^{-1}$ and $F^{(-/X^\s)}$. The following definitions are taken from \cite{BBM} (4.3.4.2).

\begin{defn} \label{defphitop} Let $top=et$ or $syn$. We use the following notations.
\debrom \item \label{defphitopi} We denote respectively $(C^{-1})$, $(nat)$, $(F)$, $(F)$, the morphism of ringed topoi which is the identity of $\smash{(X^\s/\Sigma_1)_{CRYS^\s,top}}$ (or $X^\s_{TOP^\s}$) together with the morphism of rings  $C^{-1}:\Ga=u^{-1}i^{-1}\O\to F^{-1}\O\simeq \O$, $nat:\O\to i_*i^{-1}\O=\Ga$, $F:\O\to \O$, $x\mapsto x^p$, $F:\Ga\to \Ga$, $x\mapsto x^p$ ($C^{-1}:\O=i^{-1}\O\to u_*F^{-1}\O\simeq \smash{\O_1^{crys}}$, $nat:\smash{\O_1^{crys}=u_*\O\to u_*i_*i^{-1}\O=\O}$, $F:\smash{\O_1^{crys}}\to \smash{\O_1^{crys}}$, $x\mapsto x^p$, $F:\O\to \O$, $x\mapsto x^p$).

\item \label{defphitopii} We denote
 $$\xymatrix{((X^\s/\Sigma_1)_{CRYS^\s,top},\O)\ar[r]^-\phi&(X^\s_{TOP^\s},\O)}$$ the morphism of ringed topoi defined by the morphism of topoi $u$ together with the morphism of rings $C^{-1}:\O\rightarrow \smash{\O_1^{crys}}$. 

\item \label{defphitopiii} We use similar notations in the context of big or small topoi instead of $\s$-big ones (the natural functoriality morphisms $F^{-1}\O\to \O$ might not be invertible in this context but still can be used in the same way for the above definitions).
\finrom
\end{defn}

Let us gather some compatibilities into a lemma.

\begin{lem} \label{compatcartier} Let $top=et$ or $syn$.
\debrom \item \label{compatcartieri} There is a canonically pseudo-commutative diagram of ringed topoi
$$\xymatrix{(({X^\s}/\Sigma_1)_{CRYS^\s,top},\Ga)\ar[r]^{(nat)}&(({X^\s}/\Sigma_1)_{CRYS^\s,top},\O)
\ar[r]^{(C^{-1})}\ar[d]_u\ar[rd]^\phi&
(({X^\s}/\Sigma_1)_{CRYS^\s,top},\Ga)\ar[d]_{u}\\
(X^\s_{TOP^\s},\O)\ar[u]^{i}\ar[ru]^i\ar[r]^{(nat)}&(X^\s_{TOP^\s},\O^{cris}_1)\ar[r]^{(C^{-1})}&
(X^\s_{TOP^\s},\O)}$$
This diagram is pseudo-functorial with respect to $X^\s$. The same is true in the context of big or small topoi and the resulting diagrams are
pseudo-compatible via the projection weak morphisms $\pi$.

\item \label{compatcartierii} We have  canonical isomorphisms
$$\begin{array}{rcl}
i_*\simeq u^{-1}&:&Mod(S^\s_{TOP^\s},\O)\rightarrow Mod(({S^\s}/\Sigma_1)_{CRYS^\s,top},\Ga)\\
(C^{-1})^*\simeq(F)^*(nat)_*&:&Mod(({S^\s}/\Sigma_1)_{CRYS^\s,top},\Ga)\rightarrow (({S^\s}/\Sigma_1)_{CRYS^\s,top},\O) \\
\phi^*\simeq (F)^*i_*&:&Mod(S^\s_{TOP^\s},\O)\rightarrow Mod((S^\s/\Sigma_1)_{CRYS^\s,top},\O))\end{array}$$
and similarly in the context of big or small topoi. \finrom
\end{lem}
Proof. \ref{compatcartieri} results easily from the definitions.

\ref{compatcartierii} The first isomorphism is clear. The second one results from the fact that $nat:\O\rightarrow \Ga$ is epimorphic and the third
one follows formally.
\begin{flushright}$\square$\end{flushright}

\begin{lem} \label{cartierlin} \debrom \item \label{cartierlini} The morphism $F^{(-/X^\s)}$ of Def. \ref{defFrel} (in either the $\s$-big, big, small, crystalline or usual topoi setting) induces a morphism between the following (weak) endomorphisms of ringed topoi: \begin{eqnarray}\label{defFrellin}\xymatrix{F\ar[rr]^{F^{(-/X^\s)}}&& (F)}\end{eqnarray}
This morphism is functorial in $X^\s$ and  compatible with $\pi$, $\epsilon$, $\lambda$, $i$, $\phi$.
\item \label{cartierlinii} The morphism $C^{-1}$ of Lem.-Def. \ref{defCinv} (in either the $\s$-big, big or small crystalline topoi setting) induces a morphism between the following (weak) endomorphisms of ringed topoi: \begin{eqnarray} \label{defcartierlin}\xymatrix{F\ar[rr]^{C^{-1}}&& i\phi}\end{eqnarray}
This morphism is functorial in $X^\s$ and  compatible with $\pi$ and $\epsilon$.
\finrom
\end{lem}
Proof. Left to the reader.
\begin{flushright}$\square$\end{flushright}


\begin{lem}  \label{freliso} Assume that $top$ is coarser than or equal to $et$.
\debrom \item \label{frelisoi} On $(X^\s_{top},\O)$, the relative Frobenius $\smash{F^{(-/X^\s)}:F\to (F)}$ is an isomorphism.

\item \label{frelisoii} On  $\cal Crys((X^\s/\Sigma_1)_{CRYS^\s,top},\O)$, the relative Frobenius $\smash{F^{(-/X^\s)}}:(F)^*\to F^*$ is an isomorphism. In particular, the functor $(F)^*$ preserves crystals of $((X^\s/\Sigma_1)_{CRYS^\s,top},\O)$. The same is true with $CRYS$ or $crys$ instead of $CRYS^\s$.
\finrom
\end{lem}
Proof.  \ref{frelisoi} Since $F^{(U^\s/X^\s)}:U^\s\to F^{-1}U^\s$ is an isomorphism for any $U^\s$ in $top(X^\s)$, we find that $F_*M\to (F)_*M$ is invertible for any module $M$.

\ref{frelisoii} It suffices to prove that $(F)^*M_{|T^\s}\to F^*M_{|T^\s}$ is invertible for any $T^\s$ in $CRYS^\s_{top}(X^\s/\Sigma_1)$. This follows from \ref{frelisoi} by compatibility of the relative Frobenius with $\lambda_{T^\s}$ since \break $\smash{M_{|T^\s}\simeq \lambda_{T^\s}^*M_{T^\s}}$.
\begin{flushright}$\square$\end{flushright}

\begin{rem} \label{remfreliso} The statements of Lem. \ref{freliso} have the following consequences for $zar\preceq top \preceq fl$ (using Lem.-Def. \ref{lemqcohsch} and  Lem.  \ref{crysepsilon}).

\debrom \item \label{remfrelisoi} If $M$ is a quasi-coherent module of the $\s$-big, big or small $top$ topos of $X^\s$ then $\smash{F^{(M/X^\s)}:(F)^*M\to F^*M}$ is an isomorphism.

\item \label{remfrelisoii} If $M$ is a quasi-coherent crystal of the $\s$-big, big or small crystalline $top$ topos of $(X^\s/\Sigma_1)$ then  $\smash{F^{(M/X^\s)}:(F)^*M\to F^*M}$ is an isomorphism.
\finrom
\end{rem}

\para \label{cartiergen}

Consider a separated log scheme $X^\s$ with local finite $p$-bases over $\Sigma_1$ and let $Y^\s=T^\s=X^\s$. Let $X_\lambda$ and $(\emptyset,\underline t)$ be as in section \ref{appnablamod}, so that the $\O$-algebras $\smash{\cal
P_{X}^{\s(1)}}$ and $\smash{\cal D^{\s}}$ of $\smash{X^\s_{et}}\simeq X_{et}$ have the following local description:
$$\begin{array}{ccc}(\cal P_X^{\s(1)})_{|X_\lambda}\simeq \cal O_{|X_\lambda}<\underline \tau^\s >&\hbox{\hspace{1cm} and \hspace{.5cm}}&\cal D^{\s}_{|X_\lambda}\simeq \cal O_{|X_\lambda}[\underline
\partial^\s]\end{array}$$

Let us furthermore denote $\cal P_X^{(1)}$ and $\cal D$ the rings of $X_{et}$ corresponding to the underlying scheme without log structures $X=Y=T$, so that explicitly $$\begin{array}{ccc}(\cal P_X^{(1)})_{|X_\lambda}\simeq \cal O_{|X_\lambda}<\underline \tau>&\hbox{\hspace{1cm} and \hspace{.5cm}}&\cal D_{|X_\lambda}\simeq \cal
O_{|X_\lambda}[\underline \partial]\end{array}$$ where $\tau_a\mapsto t_a\tau_a^\s$ (resp. $\partial^\s_a\mapsto t_a\partial_a$) under the natural ring
homomorphism $$\begin{array}{ccc} \cal P_X^{(1)}\rightarrow \cal P_X^{\s(1)} &\hbox{\hspace{1cm} and \hspace{.5cm}}&\cal D^\s\rightarrow \cal
D\end{array}$$
evaluated at $X_\lambda$.

\begin{defn} \label{defKprings} We use the following notations.
\debrom \item \label{defKpringsi} Let $\smash{\cal P^{(1),F}_X}$ denote the structural ring of $(X{\times_{F,X,F}}X)_{et}$  viewed as an object of
$X_{et}$  via the diagonal immersion $X\rightarrow X{\times_{F,X,F}}X$ and endowed with the structure of $\O$-algebra coming from the first projection.

\item \label{defKpringsii} Let $\cal D^{F}$ denote the module $\cal Hom(\smash{\cal
P^{(1),F}_X},\O)$ of $(X_{et},\O)$. \finrom
\end{defn}

\begin{lem} \label{lemKprings} \debrom \item \label{lemKpringsi} The $\O$-module  $\smash{\cal P^{(1),F}_X}$ naturally identifies with a direct summand of the $\O$-module  $\smash{\cal P^{(1)}_X}$, which is moreover stable by multiplication and comultiplication. Explicitly: $$(\cal P_X^{(1),F})_{|X_\lambda}\simeq \cal O_{|X_\lambda}[\tau_1,\dots, \tau_d]/(\tau_1^p, \dots, \tau_d^p).$$

 \item \label{lemKpringsii} The $\O$-module $\cal D^F$ has a natural $\O$-algebra structure for which it is a quotient of $\cal D$. Explicitly: $$(\cal D^F)_{|X_\lambda}\simeq \cal O_{|X_\lambda}[\partial_1,\dots, \partial_d]/(\partial_1^p,\dots,\partial_d^p)$$
\finrom
\end{lem}
Proof. Since the diagonal closed immersion $X\to X\times_{F,X,F}X$ has divided powers, we have a canonical morphism $inc:X\times_{F,X,F}X\to T^{(1)}$. Since on the other hand the ideal of the closed immersion $X\times_{F,X,F}X\to X\times X$ is generated by $p^{th}$ powers, the morphism $inc$ has a canonical retraction $ret:T^{(1)}\to X\times_{F,X,F}X$. All statements follow easily from the fact that $inc$ and $ret$ are compatible with the projection morphisms $p_0$, $p_1$.
\begin{flushright}$\square$\end{flushright}

%

Recall the following equivalences of categories (Prop. \ref{crystals}) $$\begin{array}{lrcl}&\cal Crys((X^\s/\Sigma_1)_{crys,et},\O)&\simeq& \nabla\hbox{-}Mod(X^\s)\\
\hbox{and}&\cal Crys((X/\Sigma_1)_{crys,et},\O)&\simeq &\nabla\hbox{-}Mod(X)\end{array}$$


\begin{lemdefn} \label{defKp}
Consider $M$ in $\cal Crys((X/\Sigma_1)_{crys,et},\O)$ corresponding to $(M_X,\nabla)$ in $\nabla\hbox{-}Mod(X)$. The following conditions are equivalent

\debrom
\item \label{defKpi} The corresponding morphism $\theta:M_X\rightarrow M_X\otimes \smash{\cal P_X^{(1)}}$ takes its values in $M_X\otimes \cal P_X^{(1),F}$ viewed as a
submodule of $M_X\otimes \smash{\cal P_X^{(1)}}$.

\item \label{defKpia} The corresponding hyper $dp$-stratification  $\varepsilon:\smash{\cal P^{(1)}_X}\otimes M_X\simeq M_X\otimes \smash{\cal P^{(1)}_X}$ sends $ \cal P_X^{(1),F}\otimes M_X$ into $M_X\otimes \cal P_X^{(1),F}$.

\item \label{defKpii} The action of $\cal D$ on $M_X$ factors through the quotient ring $\cal D^F$.
\finrom

When they hold we say that $M$ (or equivalently $(M_X,\nabla)$) has \emph{trivial $p$-curvature}. The category $\cal Crys^F((X/\Sigma_1)_{crys,et},\O)$ of crystals with trivial $p$-curvature is a full subcategory of $\cal Crys((X/\Sigma_1)_{crys,et},\O)$ which is stable by subobjects and quotient objects. In particular it is Abelian and the inclusion functor is exact.
\end{lemdefn}
Proof. This is straightforward using the correspondences explained in Lem. \ref{defnabla} and Prop. \ref{hyperstrat}.
\begin{flushright}$\square$\end{flushright}

\begin{rem} \label{remstratcartier}
If $M$ is in $\cal Crys^F((X/\Sigma_1)_{crys,et},\O)$ then the resulting isomorphism $$\varepsilon^F:\smash{\cal P^{(1),F}_X}\otimes M_X\simeq M_X\otimes \smash{\cal P^{(1),F}_X}$$ satisfies the cocycle condition,  since $\varepsilon$ does. We also note that $\varepsilon$ can be recovered from $\varepsilon^F$ by scalar extension via $\smash{\cal P^{(1),F}_X}\to \smash{\cal P^{(1)}_X}$.
\end{rem}

The following result gives an interpretation of Cartier equivalence using the morphism $$\phi:((X/\Sigma_1)_{crys,et},\O)\rightarrow (X_{et},\O)$$ We will follow
mainly the arguments of \cite{Be4} Thm. 2.3.6. This interpretation is peculiar to small \'etale sites and the reader is referred to \cite{FM} II, Sect. 1.6 for a good understanding of the
situation in the setting of small syntomic sites.

\begin{prop} \label{propcartier} Recall that $X$ has local finite $p$-bases over $\Sigma_1$.
The functor $\phi^*:Mod(X_{et},\O)\rightarrow Mod((X/\Sigma_1)_{crys,et},\O)$ induces an equivalence of categories $$Mod(X_{et},\O)\simeq
Crys^{F}((X/\Sigma_1)_{crys,et},\O)$$
The properties $qcoh$, $lf$, $lfft$ are moreover preserved under this equivalence and a quasi-inverse is induced by $\phi_*$.
\end{prop}
Proof. We begin with a lemma.

\begin{lem} \label{lemcartier} \debrom \item \label{lemcartieri} Consider a cartesian square of schemes $$\diagram{U_1&\hfl{\iota}{}&T_1\cr \vfl{h}{}&\square&\vfl{h}{}\cr U_2&\hfl{\iota}{}&T_2}$$ where the  horizontal arrows are finite morphisms. The base change morphism $h^*\iota_*M\rightarrow \iota_*h^*M$ is invertible for any $M$ in $Mod(U_{2,et},\O)$.

\item \label{lemcartierii} Let $f:Y\rightarrow X$ denote a finite surjective morphism and set $Y':=Y\times_XY$, $Y'':=Y\times_XY\times_XY$. Let $prop\in \{\emptyset, qcoh,lf,lfft\}$. The category $Mod_{prop}(X_{et},\O)$ is equivalent to that of
$\O$-modules on $Y_{et}$ satisfying $prop$ with descent datum relatively to $f$.

\item \label{lemcartieriii} Consider a commutative triangle  $$\xymatrix{U'\ar[d]_h\ar[r]^{\iota'}&T\\U\ar[ru]_{\iota}&}$$
    in $\cal Sch/\Sigma_1$ and assume
that $\iota$ and $\iota'$ are nilimmersions  of order $p$. The adjunction morphism $id\to h_*h^*$ induces an isomorphism $(F)^*\iota_*M\simeq (F)^*\iota'_*h^*M$ for any $M$ in $Mod(U_{et},\O)$.
\finrom
\end{lem}
Proof. \ref{lemcartieri} and \ref{lemcartierii} are respective variants of Prop. 5.5 and Thm. 9.4 of \cite{SGA4-II} VIII. The proofs given there for Abelian sheaves can easily be adapted.

\ref{lemcartieriii} Since $\iota$ and $\iota'$ are nilimmersions of order $p$, there exist unique $\overline F$ and $\overline F'$ factorizing the absolute Frobenius as follows: $$\xymatrix{T\ar[rd]_{\overline F}\ar[r]^{\overline F'}&U'\ar[d]_h\ar[r]^{\iota'}&T\\&U\ar[ru]_{\iota}&}$$
Now $\iota_*:Mod(U_{et},\O)\to Mod(T_{et},\O)$ is a fully faithful functor. Hence, $$F^*\iota_*M\simeq \overline F^*\iota^*\iota_*M\simeq \overline F^*M$$
Similarly $F^*\iota'_*h^*M\simeq \overline F'^*h^*M$. The natural morphism $$F^*\iota_*M\to F^*\iota'_*h^*M$$ thus identifies with the transitivity isomorphism $$\overline F^*M\simeq \overline F'^*h^*M$$
The result follows by Lem. \ref{freliso} \ref{frelisoi}.
\begin{flushright}$\square$\end{flushright}



\noindent
We can now proceed with the proof of Prop. \ref{propcartier}.  The first thing to check is that, for any $M$ in $Mod(X_{et},\O)$,  $\phi^*M$ is a crystal. Consider $(h_{T},h_U):(U_1,T_1,\iota_1,\gamma_1)\rightarrow (U_2,T_2,\iota_2,\gamma_2)$  in $Crys_{et}(X/\Sigma_1)$.
We have the following compatible isomorphisms in $Mod(T_{1,et},\O)$:

$$\diagram{(\phi^*M)_{T_1}&\mathop\simeq\limits^{\ref{compatcartier} \ref{compatcartierii}} &
((F)^*i_*M)_{T_1}&\mathop\simeq\limits^{\ref{lemreal} \ref{lemrealii}}
&(F)^*(i_*M)_{T_1}&\mathop\simeq\limits^{\ref{lemreal2} \ref{lemreal2i}} &
(F)^*\iota_{1,*}M_{|U_1}\cr
\vflupcourte{}{}&&\vflupcourte{}{}&&\vflupcourte{}{}&&\vflupcourte{ch}{}\cr
h_T^*(\phi^*M)_{T_2}&\mathop\simeq\limits^{\ref{compatcartier} \ref{compatcartierii}} &
h_T^*((F)^*i_*M)_{T_2}&\mathop\simeq\limits^{\ref{lemreal} \ref{lemrealii}}
&(F)^*h_T^*(i_*M)_{T_2}&\mathop\simeq\limits^{\ref{lemreal2} \ref{lemreal2i}} &
(F)^*h_T^*\iota_{2,*}M_{|U_2}\cr
}$$
It thus suffices to prove that the arrow denoted $ch$ is invertible. This, in turn, follows from Lem. \ref{lemcartier}
\ref{lemcartieri} and \ref{lemcartieriii} applied respectively to the square and triangle of the commutative diagram:
$$\xymatrix{U_1
\ar[rd]_-{\iota_1}
\ar[r]\ar @<+2pt> `u[r] `[rr]^-{h_U} [rr] &U_1\times_{T_2}T_1\ar[r]
\ar[d]
&U_2
\\&T_1\ar[r]^-{h_T}&T_2\ar@{<-}[u]_-{\iota_2}
}$$

We have thus checked that $\phi^*M$ is a crystal. Let $(M_X,\varepsilon)$ denote the corresponding hyper $dp$-stratification. Using the isomorphism $\phi^*M\simeq (F)^*i_*M$ (Lem. \ref{compatcartier} \ref{compatcartierii}) one checks easily that $\varepsilon$ satisfies the condition of Lem.-Def. \ref{defKp} \ref{defKpia}, ie. has trivial $p$-curvature. Using the diagonal equivalence $\iota:X_{et}\simeq (X\times_{F,X,F}X)_{et}$ and Lem. \ref{freliso} \ref{frelisoi}, we may translate the corresponding $\varepsilon^F$ on $(\phi^*M)_X\simeq (F)^*M$ as a descent datum on $F^*M$ along $F$. This descent datum is the canonical one and we may thus conclude by
Lem. \ref{lemcartier} \ref{lemcartierii} applied to $f=F$ (note that $F$ is finite since $X^\s$ has local finite $p$-bases).

It remains to check that $\phi_*$ preserves $prop$ but this also follows from Lem. \ref{lemcartier} \ref{lemcartierii} which ensures that $M$ is $prop$ if and only if
$\phi^*M$ is.
%
%
\begin{flushright}$\square$\end{flushright}

\para \label{cartierfactor}
We explain some exactness properties of the category of crystals with respect to the following morphisms: \begin{eqnarray}\label{defo}
\xymatrix{((X^\s/\Sigma_1)_{crys,et},\O)\ar[r]^-o &((X/\Sigma_1)_{crys,et},\O)\ar[r]^-{\phi_X}&(X_{et},\O)\ar @{<-} @<+0pt> `u[l] `[ll]_-{\phi_{X^\s}}
[ll]}\end{eqnarray} (here as in the previous paragraphs we have identified $(X^\s_{et},\O)$ and $(X_{et},\O)$). For simplicity we assume that $X^\s$ has a fixed global finite $p$-basis of the form $(\emptyset,\underline t)$.

%
%

\begin{defn} \label{deftfr} \debrom \item \label{deftfri} We say that a module $M$ on $(X_{et},\O)$ is \emph{$\underline t$-torsion free} if multiplication by $\smash{t_i}$ on $M$ is monomorphic for all $i$. The fully $e$-exact subcategory of $Mod(X_{et},\O)$ formed by such modules is denoted $\smash{Mod_{\underline t\hbox{-}fr}(X_{et},\O)}$.
\item \label{deftfrii} We say that a crystal on $((X/\Sigma_1)_{crys,et},\O)$ is \emph{$\underline t$-torsion free} if its realization at $X$ is $\underline t$-torsion free. The fully $e$-exact subcategory of $\cal Crys((X/\Sigma_1)_{crys,et},\O)$ formed by such crystals is denoted $\smash{\cal Crys_{\underline t\hbox{-}fr}}((X/\Sigma_1)_{crys,et},\O)$. We use similar notations with $X^\s$ instead of $X$.
\finrom
\end{defn}

\begin{rem} \label{remtfr} \debrom \item \label{remtfri} Def. \ref{deftfr} \ref{deftfri}  does not depend on the chosen $p$-basis $(\emptyset, \underline t)$. Indeed if $M$ is {$\underline t$-torsion free} then multiplication on $M$ by any local section of the monoid of $X^\s$ is monomorphic as well.

\item \label{remtfrii} A similar remark holds for Def. \ref{deftfr} \ref{deftfrii}. Moreover if a crystal $M$ on $(X^\s/\Sigma_1)$ is {$\underline t$-torsion free} and if $T^\s$ is the logarithmic divided power envelope of a local embedding of $X^\s$, then multiplication on the realization of $M$ at $T^\s$ by any local section of the monoid of $T^\s$ is monomorphic as well (Lem. \ref{lemflatlift1} and Lem. \ref{lemflatlift2}).

\item \label{remtfriii} Locally free modules on $(X_{et},\O)$ and locally free crystals of $((X/\Sigma_1)_{crys,et},\O)$ or $((X^\s/\Sigma_1)_{crys,et},\O)$ are \emph{$\underline t$-torsion free} in virtue of Lem. \ref{lemflatcharts1} \ref{lemflatcharts1i}. We thus have the following inclusions of fully $e$-exact subcategories (Sect. \ref{deffullyexact}):
    $$Mod_{lf}(X_{et},\O)\subset Mod_{\underline t\hbox{-}fr}(X_{et},\O)\subset Mod(X_{et},\O)$$ and similarly for the categories of crystals of $((X/\Sigma_1)_{crys,et},\O)$ or $((X^\s/\Sigma_1)_{crys,et},\O)$.
\finrom
\end{rem}

The category of \emph{$\underline t$-torsion free}  crystals is well behaved with respect to $o:X^\s\to X$. Using this fact, the next lemma gathers some consequence of Prop. \ref{propcartier}.

\begin{lem} \label{lemcartierex2}
\debrom \item \label{lemcartierex2i} The functors \begin{eqnarray}\label{ocartier}o^*&:&\cal Crys((X/\Sigma_1)_{crys,et},\O)\to \cal Crys((X^\s/\Sigma_1)_{crys,et},\O)\\
\label{phicartier}\phi_X^*&:&Mod(X_{et},\O)\to \cal Crys((X/\Sigma_1)_{crys,et},\O)\\
\label{phiscartier}\phi_{X^\s}^*&:&Mod(X_{et},\O)\to \cal Crys((X^\s/\Sigma_1)_{crys,et},\O)\end{eqnarray}
preserve $\underline t$-torsion freeness and are conservative for this property. The respective right adjoints $\smash{\phi_{X,*}}$, $\smash{\phi_{X^\s,*}}$ of the latter two preserve $\underline t$-torsion freeness as well.

\item \label{lemcartierex2ii} The functors (\ref{ocartier}), (\ref{phicartier}) and (\ref{phiscartier}) are exact. The functor (\ref{phicartier}) is fully faithful. The functor (\ref{ocartier}) (resp. (\ref{phiscartier})) is faithful and its restriction to $\underline t$-torsion free crystals (resp. modules) is fully faithful.
\item \label{lemcartierex2iii}
Consider a short exact sequence  $\cal E:0\to M_1\to M_2\to M_3\to 0$ of the category  $\smash{\cal Crys((X/\Sigma_1)_{crys,et},\O)}$ or $\smash{\cal Crys((X^\s/\Sigma_1)_{crys,et},\O)}$ accordingly. Let  $prop\in\{{\underline t\hbox{-}fr},qcoh,lf, lfft\}$.

-  If $M_2\simeq \smash{\phi_X^*}M'_2$ for some $M'_2$ in $Mod(X_{et},\O)$ then $\cal E\simeq \smash{\phi_X^*}\cal E'$ for some short exact sequence $\cal E':0\to M_1'\to M_2'\to M_3'\to 0$ of $Mod(X_{et},\O)$. Moreover $M_i'$ satisfies $prop$ if and only if $M_i$ does.

- If $M_1$, $M_3$ are in $\smash{\cal Crys_{{\underline t\hbox{-}fr}}((X^\s/\Sigma_1)_{crys,et},\O)}$ and $M_2\simeq o^*M'_2$ for some $M'_2$ in $\smash{\cal Crys((X/\Sigma_1)_{crys,et},\O)}$ (resp. $M_2\simeq \smash{\phi_{X^\s}^*}M'_2$ for some $M'_2$ in $Mod(X_{et},\O)$) then $\cal E\simeq o^*\cal E'$ for some short exact sequence $\cal E':0\to M_1'\to M_2'\to M_3'\to 0$ of $\smash{\cal Crys((X/\Sigma_1)_{crys,et},\O)}$ (resp. $Mod(X_{et},\O)$).
 Moreover $M_i'$ satisfies $prop$ if and only if $M_i$ does.

\finrom
\end{lem}
Proof. \ref{lemcartierex2i}  That $\smash{\phi_X^*}$ (resp. $\smash{\phi_{X,*}}$, resp. $\smash{\phi_{X^\s,*}}$) preserves $\underline t$-torsion freeness follows easily from the third isomorphism in Lem. \ref{compatcartier} \ref{compatcartierii}, using the flatness of $F:X\to X$
(resp. follows from the fact that $\smash{\phi_{X,*}}M$ is a submodule of $(F)_*M_X$, resp. follows from the fact that $\smash{\phi_{X^\s,*}}M$ is a submodule of $(F)_*M_{X^\s}$). For a crystal $M$, the isomorphism $(o^*M)_{X^\s}\simeq M_X$ shows that $o^*$ does not affect the property of being $\underline t$-torsion free. Since $\smash{\phi_{X^\s}}=\smash{\phi_{X^\s}} o$, it only remains to prove that $\smash{\phi_{X^\s}^*}$ is conservative for  $\underline t$-torsion freeness. Now, this follows from Prop. \ref{propcartier}, thanks to the fact that $\underline t$-torsion freeness is preserved by $\smash{\phi_{X^\s,*}}$.

\ref{lemcartierex2ii} In terms of modules with connection, $o^*$ sends $(M,\nabla)$ to
$(M,\nabla')$ where $\nabla'$ is deduced from $\nabla$ via $\smash{\Omega_{X/\Sigma_1}}\rightarrow
\smash{\Omega_{X^\s/\Sigma_1}}$. This shows that  (\ref{ocartier}) is exact and faithful. Recall that $\Omega_{X/\Sigma_1}\simeq \oplus_i \O dt_i$, $\Omega_{X^\s/\Sigma_1}\simeq \oplus_i \O dlog(t_i)$  and $dt_i\mapsto t_idlog(t_i)$; in particular $\smash{M\otimes \Omega_{X^\s/\Sigma_1}}\to M  \otimes \smash{\Omega_{X/\Sigma_1}}$ is monomorphic as long as $M$ is $\underline t$-torsion free. This remark shows that the restriction of $o^*$ to $\underline t$-torsion free crystals is fully faithful. The functor (\ref{phicartier}) is exact and fully faithful in virtue of Lem.-Def. \ref{defKp} and Prop. \ref{propcartier}. The statement about (\ref{phiscartier}) follows by composition.

\ref{lemcartierex2iii} The case of $\smash{\phi_X^*}$ is clear by Prop. \ref{propcartier}, Lem.-Def. \ref{defKp} and \ref{lemcartierex2i}. In the case of $o^*$ we have to check that the connection of $\smash{M_{1,X^\s}}$ and $\smash{M_{3,X^\s}}$ have no logarithmic poles, or equivalently that $\smash{M_{1,X^\s}}$ is stable by the $\partial_i$'s. This is true since $\partial^\s_i=t_i\partial_i$ and $\smash{M_{3,X^\s}}$ is $\underline t$-torsion free.
The case of $\smash{\phi_{X^\s}^*}$ follows by the isomorphism $\smash{\phi_{X^\s}^*}\simeq o^*\smash{\phi_X^{*}}$.
\begin{flushright}$\square$\end{flushright}


\begin{lem} \label{compatcartierex}  The essential image of the functor $$(F)^*:\cal Crys_{\underline t\hbox{-}fr}((X^\s/\Sigma_1)_{crys,et},\O)\to \cal Crys_{\underline t\hbox{-}fr}((X^\s/\Sigma_1)_{crys,et},\O)$$  (see Lem. \ref{freliso} \ref{frelisoii})  is contained in the essential image of \begin{eqnarray}\label{defphit}\phi_{X^\s}^{*}:Mod_{\underline t\hbox{-}fr}(X_{et},\O)\to \cal Crys_{\underline t\hbox{-}fr}((X^\s/\Sigma_1)_{crys,et},\O)\end{eqnarray}
\end{lem}

Proof. It suffices to prove that the following diagram is pseudo-commutative:
$$\xymatrix{\cal Crys((X^\s/\Sigma_1)_{crys,et},\O)&\\ \cal Crys((X^\s/\Sigma_1)_{crys,et},\O)\ar[u]^{(F)^*}\ar[r]^-{i^{-1}}&Mod(X_{et},\O)\ar[lu]_-{\phi_{X^\s}^{*}}}$$
Let $M'$ in $\cal Crys((X^\s/\Sigma_1)_{crys,et},\O)$ and let $(M'_X,\nabla')$ denote the corresponding module with connection. Using Lem. \ref{freliso}
\ref{frelisoii} and the description of $F^*$ recalled in Lem. \ref{functnabla}, we find that $M:=(F)^*M'$ corresponds to $(M_X,\nabla)$, where
$M_X=(F)^*M'_X$ and $\nabla$ is the trivial connection $\lambda\otimes x\mapsto (1\otimes x)\otimes d\lambda$. The claimed pseudo-commutativity follows easily using the isomorphism $\smash{\phi_{X^\s}^{*}}\simeq (F)^*i_*$ (Lem. \ref{compatcartier} \ref{compatcartierii}).
\begin{flushright}$\square$\end{flushright}

\subsection{The mod $p$ Hodge filtration on the small crystalline \'etale site}
\label{tmphfotsces} ~~ \\

In this section we define the tangeant sheaf and the mod $p$ Hodge filtration of a twisted Dieudonn\'e crystal over $X^\s$.  Then, we deduce a filtration on the corresponding linearized crystal and de Rham complexes. Finally, we use the functor $\jj^*$ to normalize $Fil^1$ and get the morphism $\varphi$ occurring in the definition of syntomic complexes. Here we have to assume that $X^\s$ has local finite $p$-bases over $\Sigma_1$ since the construction relies on Cartier's descent.

\para
Recall the following definitions (see e.g. \cite{dJ1} Rem. 2.4.10 for \ref{defDCii}).

\begin{defn} \label{defDC} Let $X^\s$ in $\cal Sch^\s/\Sigma_1$.
\debrom \item \label{defDCi} A \emph{Dieudonn\'e crystal}  on $X^\s$ is a triple $(D,f,v)$ where $D$ is a crystal of locally free modules of finite type
of $(X^\s/\Sigma_\infty)_{crys,et}$ and $f:F^*D\rightarrow D$ and $v:D\rightarrow F^*D$ satisfy $fv=p$ and $vf=p$. A morphism of Dieudonn\'e crystals is a morphism of crystals which is compatible with the operators $f$ and $v$. The category of Dieudonn\'e crystals on $X^\s$
is denoted ${\cal DC}(X^\s)$.

\item \label{defDCii} Assume that $X^\s$ is locally embeddable. A \emph{truncated Dieudonn\'e crystal of level $1$} on $X^\s$ is a triple $(D,f,v)$ where $D$ is a crystal of locally free modules of finite type of
$({X^\s}/\Sigma_1)_{crys,et}$, and  $f:F^*D\rightarrow D$ and  $v:D\rightarrow F^*D$ fit into an exact sequence $D\rightarrow F^*D\rightarrow D\rightarrow F^*D$ of $\cal Crys((X^\s/\Sigma_1)_{crys,et},\O)$.
The category of  truncated Dieudonn\'e crystals of level $1$  on $X^\s$ is
denoted $\DC_1(X^\s)$.
\finrom
\end{defn}

\begin{lem} \label{exDC} Assume that $X^\s$ is locally embeddable.
\debrom \item \label{exDCi} There is an exact structure $e$ on $\cal DC(X^\s)$ (resp. $\cal DC_1(X^\s)$) such that the forgetful functor to  $\cal Crys_{lfft}((X^\s/\Sigma_\infty)_{crys,et},\O)$ (resp.
$\cal Crys_{lfft}((X^\s/\Sigma_1)_{crys,et},\O)$)  is $e$-exact and reflects exactness.

\item \label{exDCii} The natural morphism $\iota_1:((X^\s/\Sigma_1)_{crys,et},\O)\to ((X^\s/\Sigma_\infty)_{crys,et},\O)$ induces an $e$-exact functor $$\iota_1^{-1}:\cal DC(X^\s)\to \cal DC_1(X^\s)$$
We use the notation $(\overline D,\overline f,\overline v):=\iota_1^{-1}(D,f,v)$.
\finrom
\end{lem}
Proof. \ref{exDCi} We say that  $0\to (D,f,v)\to (D',f',v')\to (D'',f'',v'')\to 0$ of $\cal DC(X^\s)$ (resp. $\cal DC_1(X^\s)$) is exact if the underlying sequence $0\to D\to D'\to D''\to 0$ of $\cal Crys((X^\s/\Sigma_\infty)_{crys,et},\O)$ (resp. $\cal Crys((X^\s/\Sigma_1)_{crys,et},\O)$) is exact. The reader may check directly that this defines an exact structure, ie. satisfies the axioms of \cite{Bu} Def. 2.1 using that $F^*$ is $e$-exact on $\cal Crys_{lfft}((X^\s/\Sigma_\infty)_{crys,et},\O)$ (resp. $\cal Crys_{lfft}((X^\s/\Sigma_1)_{crys,et},\O)$).

\ref{exDCii}  Let us check that the claimed functor is well defined. Consider $(D,f,v)$ in ${\cal DC}(X^\s)$ and let $D_k$, $f_k:F^*D_k\to D_k$ and $v_k:D_k\to F^*D_k$ respectively denote the image of $D$, $f$ and $v$ under the
functor $\iota_k^{-1}:\cal Crys(X^\s/\Sigma_\infty,\O)\to \cal Crys(X^\s/\Sigma_k,\O)$. We need to show that the sequence
$$\xymatrix{D_1\ar[r]^-{v_1}&F^*D_1\ar[r]^-{f_1}&D_1\ar[r]^-{v_1}&F^*D_1}$$ is exact in the category of crystals over $(X^\s/\Sigma_1)$. Since
$\iota_{1,2,*}:\cal Crys(X^\s/\Sigma_1,\O)\to \cal Crys(X^\s/\Sigma_k,\O)$  is exact faithful and $\iota_{1,2,*}D_1\simeq D_2/p$,
it is equivalent to show the exactness of $$\xymatrix{D_2/p\ar[r]^-{v_2}&F^*D_2/p\ar[r]^-{f_2}&D_2/p\ar[r]^-{v_2}&F^*D_2/p}$$
Since $p=f_2v_2$ we have
$$F^*D_2/(Ker\,
f_2+Im\, v_2)\mathop{\hookrightarrow}\limits^{f_2} D_2/ Im\, p$$
Now  $Ker\, p$ contains $Ker\, f_2$ and coincides with $Im\, p$  (check
this using e.g. Lem. \ref{lememb1} \ref{lememb1v}, or Lem. \ref{lempb} \ref{lempbiv}, \ref{lempbvi}). Exactness at the second term follows. Exactness at the third term is similar. The functor $\iota_1^*$ is thus well defined. Its exactness is clear by local freeness.
\begin{flushright}$\square$\end{flushright}


\begin{rem} \label{eqDC} Replacing the small \'etale crystalline site with the $\s$-big, big or small $top$ crystalline site, $fl\preceq top\preceq et$, we may define similarly categories of Dieudonn\'e crystals $\smash{\DC_{CRYS^\s,top}(X^\s)}$, $\smash{\DC_{CRYS,top}(X^\s)}$, $\smash{\DC_{crys,top}(X^\s)}$.
\debrom
\item \label{eqDCi} The various categories of Dieudonn\'e crystals obtained in that way are all naturally equivalent by Rem. \ref{lemlcrys} and Lem. \ref{crysepsilon} (note that for a quasi-coherent crystal $M$, we have $\smash{F^*\epsilon_*M\simeq \epsilon_*F^*M}$ and $\smash{F^*\pi_*M\simeq \pi_*F^*M}$) and thus naturally endowed with an exact structure $e$.

\item \label{eqDCii} Using Prop. \ref{exactcryslf} we find that the forgetful functor from the category of Dieudonn\'e crystals on the $\s$-big, big or small $top$ crystalline site to the corresponding category of $\O$-modules is $e$-exact (because locally free crystals are flat, hence acyclic for $\epsilon^*$ and $\pi^*$) and reflects exactness (because quasi-coherent crystals are acyclic for $\epsilon_*$ and $\pi_*$, see Lem. \ref{crysepsilon} \ref{crysepsilonii}).
\finrom
\end{rem}

\begin{defn} \label{deflie} Assume that $X^\s/\Sigma_1$ has local finite $p$-bases.

\debrom \item \label{defliei} Let $D=(D,f,v)$ in $\DC_1(X^\s)$.

\noindent - Recall the morphism $\phi$ from Def. \ref{defphitop} \ref{defphitopii}. We define $Lie(D)$ as the following module of $(\smash{X_{et}},\O)$:
 $$Lie(D):=\phi_*(Coker\, v)$$

\noindent - The \emph{Hodge filtration} on the module of $((X^\s/\Sigma_1)_{crys,et},\O)$ underlying $D$ is defined as follows: $Fil^2D:=0$,
$Fil^0D:=D$ and  $$Fil^1D:=Ker(can_D:D\rightarrow i_{*}i^{-1}D\rightarrow i_{*}Lie(D))$$ where the
first one is the adjunction morphism and the second one is the inverse Cartier operator
(\ref{defcartierlin}): $i^*D\rightarrow \phi_*F^*D$. \petit

\noindent \item \label{deflieii} Let $D=(D,f,v)$ in $\DC(X^\s)$.

\noindent - We define $Lie(D)$ as $Lie(\overline D)$.

\noindent - The \emph{$mod\, p$ Hodge filtration} on the module of $((X^\s/\Sigma_\infty)_{crys,et},\O)$ underlying $D$ is defined as follows: $Fil^hD$ is the
inverse image of $\iota_{1,*}Fil^h \overline D$ by the adjunction morphism $D\rightarrow \iota_{1,*}\overline D$. \finrom
\end{defn}




\begin{prop} \label{liesurj} Assume that $X^\s/\Sigma_1$ has local finite  $p$-bases. Consider $D$ in $\DC(X^\s)$.
\debrom \item \label{liesurji}
 The morphism $\smash{can_{\overline D}}:\overline D\rightarrow i_{*}Lie(\overline D)$ occurring in the definition of $Fil^1\overline D$
is an epimorphism. It induces an exact sequence in $Mod((X^\s/\Sigma_\infty)_{crys,et},\O)$: $$\xymatrix{0\ar[r] &Fil^1D\ar[r]&D\ar[r]^-{can_D}& i_*Lie(D)\ar[r]&0}$$
\item \label{liesurjii} Consider a morphism $f:X'^\s\to X^\s$ where $X'^\s$ has local finite $p$-bases as well. There is a natural base change isomorphism $ch_f$ rendering the following square commutative: $$\xymatrix{i^{-1}f^*D\ar[rr]^-{can_{f^*D}}&&Lie(f^*D)\\
     f^*i^{-1}D\ar[u]^\wr\ar[rr]^-{f^*(can_D)}&&f^*Lie(D)\ar[u]^\wr_{ch_f}}$$ The family of the $ch_f$'s satisfies the composition constraint.
\item \label{liesurjiii} The functor $Lie:\DC(X^\s)\to Mod(\smash{X^\s_{et}},\O)$ is $e$-exact.
\finrom
\end{prop}

Proof. \ref{liesurji} Up to \'etale localization, we may assume given a $p$-basis of the form $(\emptyset,\underline t)$ as in Sect. \ref{cartiergen}, \ref{cartierfactor} (see Rem. \ref{examplesofpb} \ref{examplesofpbiii}). Since $\overline D$ is a crystal, $\overline D\rightarrow i_{*}i^{-1}\overline D$ is epimorphic in
$Mod((X^\s/\Sigma_1)_{crys},\O)$ (by Lem. \ref{lemreal2} \ref{lemreal2i} and the crystal condition). We want to prove that
$i_{*}i^{-1}\overline D\rightarrow i_*Lie(\overline D)$ is epimorphic, ie. that
\begin{eqnarray}\label{comparrrow}\xymatrix{i^{-1}\overline D\ar[r]^{C^{-1}}& \phi_*F^*\overline D\ar[r]^{\phi_*(nat)}& \phi_*Coker\overline v}\end{eqnarray}  is epimorphic in $Mod(X_{et},\O)$. Here $nat$ denotes the tautological morphism $F^*\overline D\to Coker\, \overline v$.
Let us begin with the following \petit

\noindent $\emph{Claim}$: The following adjunction morphisms are isomorphisms: \begin{eqnarray}\label{isoliesurj1}i^{-1}\overline D&\to& \phi_*\phi^*i^{-1}\overline D\\
\label{isoliesurj2}\phi^*\phi_*F^*\overline D&\to& F^*\overline D\\
\label{isoliesurj3}\phi^*\phi_*Coker\,\overline v&\to& Coker\, \overline v.\end{eqnarray}
 \petit

Let us prove the claim. The isomorphism   $i^{-1}\overline D\simeq \smash{\overline D_{X^\s}}$ (Lem. \ref{lemreal2} \ref{lemreal2i}) shows that
$i^{-1}\overline D$ is  $\underline t$-torsion free. The first isomorphism of the claim follows by full faithfulness of the restriction of $\phi^*$ to
$\underline t$-torsion free modules (see Lem. \ref{lemcartierex2} \ref{lemcartierex2ii}, where $\phi^*$ was denoted $\phi^{\s,*}$). Thanks to Lem. \ref{lemcartierex2} \ref{lemcartierex2i}, the second and third isomorphism will follow if we show that $F^*\overline D$ and $Coker\, \overline v$ are in the essential image of the fully faithful functor
(\ref{defphit}). In the case of $F^*\overline D$, this results directly from Lem. \ref{compatcartierex} and Lem. \ref{freliso} \ref{frelisoii}. Next, we notice that the crystal $Coker\,
\overline v$ is $\underline t\hbox{-}$torsion free since it is a subcrystal of $\overline D$.
We may thus conclude that $Coker\,
\overline v$ is in the essential image of $\phi^*$ as well by Lem. \ref{lemcartierex2} \ref{lemcartierex2iii}, and this ends the proof of the claim.

 We will now prove that the composed arrow (\ref{comparrrow}) is epimorphic. As the reader may check, we have a commutative diagram \begin{eqnarray}\label{cartieriso}\xymatrix{i^{-1}\overline D\ar[r]^{C^{-1}}\ar[d]_{adj}& \phi_*F^*\overline D \\
\phi_*\phi^*i^{-1}\overline D\ar[r]_\sim&\phi_*(F)^*\overline D\ar[u]_{\phi_*F^{\overline D/X^\s}}}\end{eqnarray} where the lower isomorphism is by
pseudo-commutativity of the triangle in the proof of Lem. \ref{compatcartierex}. The left (resp. right) vertical arrow is invertible as well by the above claim (resp. Lem. \ref{freliso} \ref{frelisoii}). The first arrow in  (\ref{comparrrow}) is thus invertible.

Let us now investigate the second arrow in (\ref{comparrrow}). Consider the following commutative square of $\cal Crys(X^\s/\Sigma_1,\O)$: $$\xymatrix{F^*\overline D\ar[rr]^-{nat}&&Coker\, \overline v\\
\phi^*\phi_*F^*\overline D\ar[u]_{adj}^\wr\ar[rr]^-{\phi^*\phi_*(nat)}&&\phi^*\phi_*Coker\,\overline v\ar[u]_{adj}^\wr}$$ Vertical arrows are
isomorphisms by the claim. The arrow $nat$ is epimorphic and thus $\phi^*\phi_*(nat)$ also. We conclude that  $\phi_*(nat)$ is epimorphic as well since $$\phi^*:Mod(X_{et},\O)\to \cal Crys((X^\s/\Sigma_1)_{crys,et},\O)$$ is conservative for epimorphisms (it is a faithful exact functor between Abelian categories by Lem. \ref{lemcartierex2} \ref{lemcartierex2i}, \ref{lemcartierex2ii} and  Prop. \ref{exactcrysk} \ref{exactcryski}).  We have thus proven that the composed arrow (\ref{comparrrow}) is epimorphic  as desired.

We already know that the bottom sequence of the following tautologically commutative diagram of $Mod((X^\s/\Sigma_\infty)_{crys},\O)$ is exact
(recall that $\iota_{1,*}$ is exact by Lem. \ref{iotacrysex}): $$\diagram{0&\hflcourte{}{}&Fil^1D&\hfl{}{}&D&\hfl{}{}&i_*Lie(D)&\hflcourte{}{}&0\cr
&&\vfl{}{}&&\vfl{}{}&&\vfl{\wr}{}&\cr 0&\hflcourte{}{}&\iota_{1,*}Fil^1\overline D&\hfl{}{}&\iota_{1,*}\overline
D&\hfl{}{}&\iota_{1,*}i_{*}Lie(\overline D)&\hflcourte{}{}&0}$$ Exactness of the top sequence follows formally since the left square
is cartesian (by definition of $Fil^1D$) and the middle vertical arrow is epimorphic (since $D$ is a crystal).

\ref{liesurjii} The definition of the morphism $ch_f$ is purely formal from the canonical isomorphisms $\iota_1 f\simeq f\iota_1$, $Ff\simeq Ff$, $\phi f\simeq f\phi$  and $if\simeq fi$. The composition constraint for the $ch_f$'s moreover follows from the compatibility of these isomorphisms with respect to composition of different $f$'s. It remains to check that $ch_f$ is in fact an isomorphism. We can always assume that both $X^\s$ and $X'^\s$ have global $p$-bases as in the proof of \ref{liesurji}. There is the following commutative diagram on $(X'_{et},\O)$:

 $$\xymatrix{
Lie(f^*\overline D)\ar[r]^-1_-\sim& \phi_*f^*Coker\, \overline v& \phi_*f^*\phi^*Lie(\overline D)\ar[l]^-\sim_-{(\ref{isoliesurj3})}\ar[r]_-\sim & \phi_*\phi^*f^*Lie(\overline
D)
\\
 f^*Lie(\overline D)\ar[r]_-=\ar[u]^{ch_{Lie}} &
f^*\phi_*Coker\, \overline v\ar[u]^{ch} & f^*\phi_*\phi^*Lie(\overline D) \ar[u]^{ch}\ar[l]^-\sim_-{(\ref{isoliesurj3})}& f^*Lie(\overline D)\ar[l]_-2\ar[u]^2}$$
where the isomorphism denoted $1$ is defined by identifying $f^*Coker\, \overline v$
and $Coker\, f^*\overline v$. Let us finally check that the arrows denoted $2$ are isomorphisms. It suffices to check that $Lie(\overline D)$ and $f^*Lie(\overline D)$ are $\underline t$-torsion free (Lem. \ref{lemcartierex2} \ref{lemcartierex2i}, \ref{lemcartierex2ii}). This is indeed the case by  Lem. \ref{lemcartierex2} \ref{lemcartierex2i} since $\phi^*Lie(\overline D)\simeq Coker\, \overline v$ and $\phi^*f^*Lie(\overline D)\simeq f^*Coker \overline v\simeq Coker f^*\overline v$ are $\underline t$-torsion free as already explained.

\ref{liesurjiii} The question is local so we may again assume given a $p$-basis as in \ref{liesurji}. The functor $$\begin{array}{rcl}\DC(X^\s)&\to &\smash{\cal Crys((X^\s/\Sigma_1)_{crys,et},\O)}\\ D&\mapsto& Coker\, \overline v\end{array}$$ is $e$-exact since $\iota_1^{-1}:\cal DC(X^\s)\to \cal DC_1(X^\s)$ is $e$-exact (Lem. \ref{exDC} \ref{exDCii}) hence $e_M$-exact (Prop. \ref{exactcryslf} and Rem. \ref{remexactcryskiii}) and there is an isomorphism $Coker\, \overline v\simeq Ker\, \overline v$ (induced by $\overline f:F^*\overline D\to D$). Next, we observe that $\phi_*:\smash{\cal Crys((X^\s/\Sigma_1)_{crys,et},\O)}\to \smash{Mod(X_{et},\O)}$ has the following property thanks to Lem. \ref{lemcartierex2} \ref{lemcartierex2ii}, \ref{lemcartierex2iii}: if $0\to M_1\to M_2\to M_3\to 0$ is $e$-exact and the $M_i$'s are $\underline t$-torsion free and in the essential image of $\phi^*$ then  $0\to \phi_*M_1\to \phi_*M_2\to \phi_*M_3\to 0$ is exact.
\begin{flushright}$\square$\end{flushright}

\begin{lem} \label{lemlielfft} The module $Lie(D)$ is quasi-coherent. It is locally free of finite type if $X$ is locally 
Noetherian . \end{lem}

Proof. The realization of $Coker\, \overline v$ at $X^\s$ has a right resolution as follows  $$(Coker\, \overline v)_{X^\s}\mathop\rightarrow\limits^{\overline f} [D_{X^\s}\mathop\rightarrow\limits^{\overline v} F^*D_{X^\s}\mathop\rightarrow\limits^{\overline f} D_{X^\s}\mathop\rightarrow\limits^{\overline v} \dots]$$
It is thus quasi-coherent
and even locally free of finite type under the additional assumption that $X$ is locally
noetherian (hence  regular by Lem. \ref{lemrp} \ref{lemrpvi}). The same is true for $Lie(D)$ by Lem. \ref{lemcartierex2} \ref{lemcartierex2iii}.
\begin{flushright}$\square$\end{flushright}

\para \label{parafilL} Let us consider a separated log scheme $X^\s/\Sigma_1$ with local finite $p$-bases and assume moreover given a global embedding $X^\s\to Y^\s$ whose logarithmic divided power envelope is denoted $\iota:X^\s\to T^\s$. We will use the following notations.

- As mentionned in Sect. \ref{chL}, letting $k$ vary in (\ref{defLcrys}) causes no difficulty and defines a functor $L_.$ below. We define furthermore a functor $L$ by commutativity of the following
triangle:
\begin{eqnarray}\label{functorLsanspoint}\xymatrix{Hdp_{norm}(T^\s_.,\O)\ar[rd]_-{L}\ar[r]^-{L_.}&\cal
Crys_{norm}((X^\s/\Sigma_.)_{crys,et},\O)\ar[d]^-{l_*}\\
&Crys((X^\s/\Sigma_\infty)_{crys,et},\O)} \end{eqnarray}

- If $M$ is a module on $\smash{((X^\s/\Sigma_\infty)_{crys,et},\O)}$ we write $M_.:=l^{-1}M$ (resp. $\smash{M_{T_.^\s}}:=\smash{M_{.,T_.^\s}}$) the
associated normalized module of $\smash{((X^\s/\Sigma_.)_{crys,et},\O)}$ (resp. its realization on $(\smash{T^\s_{.,et},\O})$).

- If $M$ is a crystal of $\smash{((X^\s/\Sigma_\infty)_{crys,et},\O)}$, we write $\smash{\Omega^\bullet_{T^\s_.}(M)}:=\smash{\Omega^\bullet_{T^\s_.}(M_.)}$ the associated de Rham complex in $(\smash{T^\s_{.,et}},\smash{\O^{crys}_.})$.

\begin{prop} \label{complin}  Consider $D\in \cal DC(X^\s)$.


\debrom \item \label{complini} There exists a canonical morphism $can_L$ rendering the right square below commutative. We let  $Fil^1(L(D_{T_.^\s}))$
denote its kernel. Whence a tautological morphism of exact sequences on $((X^\s/\Sigma_\infty)_{crys,et},\O)$: $$\xymatrix{0\ar[r]
&Fil^1D\ar[d]\ar[r]&D\ar[d]_{aug}^{(\ref{defaugT})}\ar[r]^-{can_D}& i_*Lie(D)\ar[d]^\parallel\ar[r]&0\\ 0\ar[r]
&Fil^1(L(D_{T_.^\s}))\ar[r]&L(D_{T_.^\s})\ar[r]^-{can_L}& i_*Lie(D)\ar[r]&0}$$

\item \label{complinii} There is a natural isomorphism of exact sequences on $(T^\s_{.,et},\O^{crys}_.)$ for any $h\in \G(X,M_X/\Gm)=\G((X/\Sigma_\infty)_{crys,et},M_{X/\Sigma_\infty}/\O_{X/\Sigma_\infty}^\times)$
$$\xymatrix{0\ar[r] &\iota_*u_*l^{-1}(Fil^1(L(D_{T_.^\s}))(-h))\ar[d]^\wr\ar[r]&\iota_*u_*l^{-1}(L(D_{T_.^\s})(-h))\ar[d]^\wr_{\ref{uL} \ref{uLii}}
\ar[r]^-{can_L}& \iota_*l^{-1}(Lie(D)(-h))\ar[d]^\wr_{\ref{lemreal2} \ref{lemreal2i}}\ar[r]&0\\ 0\ar[r]& (Fil^1D)(-h)_{T_.^\s}\ar[r]&D(-h)_{T_.^\s}\ar[r]^-{can_D}& (i_*Lie(D)(-h))_{T_.^\s}\ar[r]&0}$$
where the top (resp. bottom) one is deduced from the bottom  (resp. top) one in \ref{complini} by twisting and applying $\iota_*u_*l^{-1}$ (resp. $\smash{(l^{-1}(-))_{T^\s_.}}$).

\item \label{compliniii} The diagrams of \ref{complini} and \ref{complinii} are naturally lax functorial with respect to the chosen global embedding and effective log divisor. This means that there is a natural base change morphism $ch_f:f^*\ref{complini}\to \ref{complini}'$  in $Mod((X'^\s/\Sigma_\infty)_{crys,et},\O)$ (resp. $ch_f:f^*\ref{complinii} \to \ref{complinii}'$ in $Mod(\smash{(T'^\s_{.,et},\O_.^{crys})})$) for each morphism $f:(X'^\s\to Y'^\s)\to (X^\s\to Y^\s)$ in $\smash{Emb^{\s,glob}}$ (resp. for each morphism $f:(X'^\s\to Y'^\s)\to (X^\s\to Y^\s)$ in $\smash{Emb^{\s,glob}}$ and each $h'$ dividing $f^*h$) and that the family of the $ch_f$'s satisfies the composition constraint.
\finrom
\end{prop}
Proof. \ref{complini} The only point requiring explanations is the definition of $can_L$. Recall from Def. \ref{deflie} \ref{defliei}, \ref{deflieii} that $can_D$ factorizes via the adjunction morphisms $D\to i_*i^{-1}D$. By Lem.-Def. \ref{lemlimcrys} \ref{lemlimcrysi},  it is thus sufficient to provide a canonical factorization $\alpha$ of $D_.\to
i_*i^*D_.$ via $D_.\to \smash{f_{T_.^\s,*}f_{T_.^\s}^{-1}D_.\simeq L_.(D_{T_.^\s})}$. We do this by observing that the right vertical arrow is
invertible in the following commutative diagram of adjunction morphisms:
\begin{eqnarray}\label{defalpha}\xymatrix{D_.\ar[r]\ar[d]&i_*i^{-1}D_.\ar[d]^\wr\\ f_{T_.^\s,*}f_{T_.^\s}^{-1}D_.\ar[r]\ar  @{-->} [ru]^\alpha&f_{T_.^\s,*}f_{T_.^\s}^{-1}i_*i^{-1}D_.}\end{eqnarray}

\ref{complinii} It suffices to check that the right hand square is commutative. Using Lem.-Def. \ref{lemlimcrys} \ref{lemlimcrysi}, this will follow from the commutativity of the following diagram $$\xymatrix{\iota_*u_*((f_{T_.^\s,*}f_{T_.^\s}^{-1}D_.)(-h))\ar[d]^\wr_{\ref{lemreal2} \ref{lemreal2ii}}\ar[r]^-\alpha & \iota_*u_*((i_*i^{-1}D_.)(-h))\ar[r]^\sim&\iota_*i^{-1}(D_.(-h))\ar[d]^\wr_{\ref{lemreal2} \ref{lemreal2i}}\\
\lambda_{T^\s_.,*}f_{T_.^\s}^{-1}(D_.(-h))\ar[rr]&&\lambda_{T^\s_.,*}f_{T_.^\s}^{-1}i^{-1}i_*(D_.(-h))}$$
where the definition of each arrow uses Lem. \ref{lemtwlogtop} \ref{lemtwlogtopii} except the bottom horizontal one. We are thus reduced to noticing the commutativity of the following diagram where $1:\smash{\iota_*u_*}\to \smash{\lambda_{T^\s_.,*}f_{T_.^\s}^{-1}}$ is the natural morphism deduced from Lem. \ref{lemreal2} \ref{lemreal2ii}, $2$ expresses the isomorphism $ui=id$ and the unlabeled arrows are adjunction morphisms.
$$\xymatrix{\iota_*u_*f_{T_.^\s,*}f_{T_.^\s}^{-1}\ar[r]\ar[d]^1 \ar @{->} @<+0pt> `u[r] `[rr]^-{\alpha} [rr]
\ar @{->} @<+0pt> `l[d] `[dd]_-{\ref{lemreal2} \ref{lemreal2ii}} [dd]&
\iota_*u_*f_{T_.^\s,*}f_{T_.^\s}^{-1}i_*i^{-1}\ar[d]^1&\iota_*u_*i_*i^{-1}\ar[d]_\wr^2\ar[l]^-\sim\\
\lambda_{T^\s_.,*}f_{T_.^\s}^{-1}f_{T_.^\s,*}f_{T_.^\s}^{-1}\ar[r]\ar[d]& \lambda_{T^\s_.,*}f_{T_.^\s}^{-1}f_{T_.^\s,*}f_{T_.^\s}^{-1}i_*i^{-1}\ar[d]
& \iota_*i^{-1}\ar[d]^-{\ref{lemreal2} \ref{lemreal2i}}\\
\lambda_{T^\s_.,*}f_{T_.^\s}^{-1}\ar[r] & \lambda_{T^\s_.,*}f_{T_.^\s}^{-1}i_*i^{-1}&\ar[l]^-=
\lambda_{T^\s_.,*}f_{T_.^\s}^{-1}i_*i^{-1}\ar[lu]^-\sim\ar @{<-} @<+0pt> `r[u] `[uu]_-{1} [uu]}$$

\ref{compliniii} Let us explain the morphism $f^*\ref{complini}\to \ref{complini}'$. It is sufficient to define compatible morphisms for the vertices of the right hand square. Using obvious notations we define:

\noindent - $ch_f:f^*D\to f^*D$ as the identity;

\noindent - $ch_f:f^*i_*Lie(D)\to Lie(f^*D)$ using $if\simeq fi$ and the base change morphism for $Lie$ (Prop. \ref{liesurj} \ref{liesurjii});

\noindent - $ch_f:f^*L(\smash{D_{T_.^\s}})\to L'((f^*D)_{T'^\s})$ as the natural morphism $f^*l_*\smash{f_{T_.^\s,*}f_{T_.^\s}}\to l_*\smash{f_{T'^\s_.,*}f_{T'^\s_.}}$ resulting from $lf\simeq fl$ and $f\smash{f_{T'^\s_.}}\simeq \smash{f_{T^\s_.}}f$.

In each case the composition constraint for the $ch_f$'s results from the fact that each one of the morphisms used to build $ch_f$ satisfies the composition constraint. Compatibility with $can_D$ (resp. $can_L$) follows from Prop. \ref{liesurj} \ref{liesurjii} (resp. Prop. \ref{liesurj} \ref{liesurjii} and (\ref{defalpha})). Compatibility with $aug$ follows  from (\ref{functLdR}).

Let us explain the morphism $f^*\ref{complinii}\to \ref{complinii}'$. We define:

\noindent - $ch_f:\smash{f^*(D(-h)_{T^\s_.})\to (f^*D(-h'))_{T'^\s_.}}$ using Lem. \ref{lemcrystal2} \ref{lemcrystal2ii} and Lem. \ref{lemtwlogtop} \ref{lemtwlogtopii}, \ref{lemtwlogtopi}.

\noindent - $ch_f:\smash{f^*\iota_*u_*(L_.(D_{T^\s_.})(-h))\to \iota'_*u_*(L'_.((f^*D)_{T'^\s_.})(-h'))}$, using the morphisms  $\iota f\simeq f\iota'$, $uf\simeq fu$, $\smash{L_.(D_{T^\s_.})(-h)}\simeq \smash{L_.(D(-h)_{T^\s_.})}$ (recall that $\smash{(-)_{T^\s_.}=\lambda_{T^\s_.,*}f_{T^\s_.}^{-1}}$, $\smash{L_.\simeq f_{T^\s_.,*}\lambda_{T^\s}^*}$ and apply Lem. \ref{lemtwlogtop} \ref{lemtwlogtopii}), $\smash{f^*L_.((-)_{T^\s_.})}\to \smash{L'_.((f^*(-))_{T'^\s_.})}$ (see above), $f^*(D(-h))\to f^*D(-h')$ (Lem. \ref{lemtwlogtop} \ref{lemtwlogtopi}, \ref{lemtwlogtopii}), $\smash{L'_.(f^*D(-h')_{T'^\s_.})}\simeq \smash{L'_.(f^*D_{T'^\s_.})(-h')}$ (Lem. \ref{lemtwlogtop} \ref{lemtwlogtopii}).

\noindent - $ch_f:f^*\iota_*l^{-1}(Lie(D)(-h))\to \iota'_*l^{-1}(Lie(f^*D))$, using $\iota f\simeq f\iota'$, $lf\simeq fl$ and the base change morphism for $Lie$;

\noindent - $ch_f:\smash{f^*(i_*Lie(D)(-h))_{T^\s_.}}\to i_*Lie(f^*D)(-h)_{T'^\s_.}$, using Lem. \ref{lemcrystal2} \ref{lemcrystal2ii}, $lf\simeq fl$ and $if\simeq fi$.

Here again the composition constraint causes no difficulty. Let us explain the compatibility with respect to the boundaries of the right hand square in Prop. \ref{complin} \ref{complinii}. Compatibility with the horizontal arrows denoted $can_L$ and $can_D$ just follows from the corresponding compatibility in \ref{complini}. Compatibility with the left vertical arrow (denoted Lem. \ref{uL} \ref{uLii}) easily reduces to the commutativity of the following diagram $$\xymatrix{\iota uf_{T^\s_.}f\ar[r]\ar[d]^{\ref{lemreal2}\ref{lemreal2ii}}& \iota uff_{T'^\s_.}\ar[r]&\iota f uf_{T'^\s_.}\ar[r]&f\iota'uf_{T'^\s_.}\ar[d]^{\ref{lemreal2}\ref{lemreal2ii}}\\ \lambda_{T^\s_.} f\ar[rrr]&&& f\lambda_{T'^\s_.}}$$
which, in turn, causes no difficulty. Compatibility with the right vertical arrow follows formally (alternatively, it would follows from the  functoriality of Lem. \ref{lemreal2} \ref{lemreal2ii} in a sense that the reader can imagine).


\begin{flushright}$\square$\end{flushright}

\begin{defn} \label{deffil} Fix an effective log divisor $h\in \G(X,M_X/\Gm)$. We define three functors with value in the category $\smash{Fil^{0,1}Kom(T_{.,et}^\s,\O_.^{cris})}$ of complexes of modules of $\smash{(T_{.,et}^\s,\O_.^{cris})}$ endowed with a one step filtration $Fil^1\subset Fil^0$:
$${Fil}^._{.,T^\s}(-h),\, {Fil}^.\Omega^\bullet_{.,T^\s}(-h),\, {Fil}^.L\Omega^\bullet_{.,T^\s}(-h): \DC(X^\s)\to Fil^{0,1}Kom(T_{.,et}^\s,\O_.^{cris})$$
Their description is the following (in \ref{deffiliii}, $Fil^1$ is viewed as a subcomplex of $Fil^0$ via the isomorphism $\smash{(L(D_{T^\s_.}))(-h)\simeq L(D(-h)_{T_.^\s})}$ resulting from Lem. \ref{lemtwlogtop} \ref{lemtwlogtopii}):

\debrom \item  \label{deffili} $\smash{{Fil}^0_{.,T^\s}(-h)(D)}$ is ${D(-h)_{T^\s_{.}}}$ placed in degree $0$ and $\smash{{Fil}^1_{.,T^\s}(-h)(D)}:=\smash{(Fil^1 D)(-h)_{T^\s_{.}}}$.

\item \label{deffilii} $\smash{Fil^0\Omega^\bullet_{.,T^\s}(-h)(D)}$ is $\smash{\Omega^\bullet_{T^\s_.}(D(-h))}$ and  $\smash{Fil^1\Omega^\bullet_{.,T^\s}(-h)(D)}$ is the following subcomplex: $$\begin{array}{l}
{Fil}^1\Omega^0_{.,T^\s}(-h)(D):=(Fil^1 D)(-h)_{T^\s_{.}}
\\ \hbox{and }{Fil}^1\Omega^q_{.,T^\s}(-h)(D):={Fil}^0\Omega^q_{.,T^\s}(-h)(D) \hbox{ for $q\ge 1$}\end{array}$$

\item \label{deffiliii} $\smash{Fil^0L\Omega^\bullet_{.,T^\s}(-h)(D)}$ is $\smash{(L(\Omega^\bullet_{T^\s_.}(D(-h))))_{T^\s_.}}$ and $\smash{Fil^1L\Omega^\bullet_{.,T^\s}(-h)(D)}$ is the following subcomplex: $$\begin{array}{l}
Fil^1L\Omega^0_{.,T^\s}(-h)(D):= (Fil^1(L(D_{T^\s_.})))(-h)_{T^\s_.}
\\ \hbox{and } Fil^1L\Omega^q_{.,T^\s}(-h)(D):= Fil^0L\Omega^q_{.,T^\s}(-h)(D) \hbox{ for $q\ge 1$}\end{array}$$
\finrom
\end{defn}

\begin{rem} \label{remdeffil} Using Prop. \ref{complin} \ref{compliniii} and Lem. \ref{functpoincare} we find that the functors of Def. \ref{deffil} are subject to natural base change morphisms turning them into colax morphisms between contravariant pseudo functors on $\smash{Emb^{\s,glob,lfpb}_{div}}$:

$${Fil}^._{.,(-)^{dp}}(-),\, {Fil}^.\Omega^\bullet_{.,(-)^{dp}}(-),\, {Fil}^.L\Omega^\bullet_{.,(-)^{dp}}(-): \DC(-)\to Fil^{0,1}Kom((-)^{dp}_{.,et},\O_.^{cris})$$
\end{rem}

\para
We will now use the functor $\jj^*$ defined in Def. \ref{defj}  in order to get a normalized version of the filtered complexes defined in Def. \ref{deffil}. To begin with, we replace the ring $\smash{\cal O^{crys}_.}$ (which is not $\Z/p^.$-normalized for the \'etale topology) by the following.

\begin{defn} \label{defocrysmod}
Recall the morphism of ringed topoi $l:\smash{T^\s_{.,et}}\to \smash{T^\s_{et}}$. We set $$\widetilde\O^{crys}_.:=\Z/p^.\otimes l^{-1}\O^{crys}$$
\end{defn}

\begin{defn} \label{deffilmod} Consider the endomorphism $\jj$ of the ringed topos $\smash{(T^\s_{.,et},\widetilde\O^{crys}_.)}$ defined as in Def. \ref{defj}.
Let $Fun^._.$ denote one of the three functors defined in Def. \ref{deffil}. We define a functor with values in the category of arrows of complexes of modules of $\smash{(T^\s_{.,et},\widetilde O^{crys}_.)}$, $$\widetilde{Fun}{}^i_.:\DC(X^\s)\to Kom(T^\s_{.,et},\widetilde O^{crys}_.)^{[1]}$$
by setting $\smash{\widetilde{Fun}{}^i_.(D):=\jj^*{Fun}{}^i_.(D)}$.
\end{defn}

\begin{rem} \label{remdeffilmod} The morphism $\jj$ is functorial with respect to $T^\s$ in the obvious way. The three functors defined in Def. \ref{deffilmod} naturally extend to colax morphisms between contravariant pseudo-functors on $Emb^{\s,glob,lfpb}_{div}$:
$$\widetilde{Fil}{}^._{.,(-)^{dp}}(-),\, \widetilde{Fil}{}^.\Omega^\bullet_{.,(-)^{dp}}(-),\, \widetilde{Fil}{}^.L\Omega^\bullet_{.,(-)^{dp}}(-): \DC(-)\to Kom((-)^{dp}_{.,et},\O_.^{cris})^{[1]}$$
as in Rem. \ref{remdeffil}.
\end{rem}

The following lemma explains the difference between $\widetilde{Fun}{}^i_.$ and $\smash{Fun^i_.}$.

\begin{lem} \label{lemfil} Let $(T^\s,h)$ denote the logarithmic divided power envelope of some $(Y^\s,h)$ in $\smash{Emb^{\s,glob,lfpb}_{div}}$ and let $\smash{Fun^._.}$ (resp. $\smash{\widetilde{Fun}{}^._.}$), denote one of the three  functors defined in Def. \ref{deffil} and  Rem. \ref{remdeffil} (resp. Def. \ref{deffilmod} and Rem. \ref{remdeffilmod}).
 \debrom \item \label{lemfili} There is a natural morphism of functorial distinguished triangles in $\smash{D^b(T^\s_{.,et},\widetilde \O^{crys}_.)}$
$$\xymatrix{\widetilde{Fun}{}^1_.(D)\ar[r]^-1\ar[d]&\widetilde{Fun}{}^0_.(D)\ar[r]^-{can}\ar[d]& \widetilde{Lie}{}_{.,T^\s}(-h)(D)\ar[r]^-{+1}\ar[d]&\\
Fun^1_.(D)\ar[r]^-1&Fun^0_.(D)\ar[r]^-{can}&Lie_{.,T^\s}(-h)(D)\ar[r]^-{+1}&}$$
where $\smash{Lie_{.,T^\s}(-h)(D)}$ is the direct image via $\smash{(X^\N_{et},\O)}\to \smash{(T^\s_{.,et},\widetilde{\O}{}_.^{crys})}$ of the constant projective system $Lie(D)(-h)$ and $\smash{\widetilde{Lie}_{.,T^\s}(-h)(D)}:=\smash{\tau_{\ge_{-1}} L\jj^*Lie_{.,T^\s}(-h)(D)}$. The cohomology of the latter complex is described explicitly as follows: $H^0$ is $\smash{Lie_{.,T^\s}(-h)(D)}$; the $k$-th component of $H^{-1}$ is $\smash{Lie_{k,T^\s}(-h)(D)}$ as well but the transition morphisms of the projective system are zero.

\item \label{lemfilii} The objects of the complexes $\smash{{Fun}{}^i_.(D)}$ (resp.  $\smash{\widetilde{Fun}{}^i_.(D)}$)
come from normalized (resp. $L$-normalized) quasi-coherent modules on $\smash{T_{.,et}}$ by scalar restriction via $\smash{{\O}{}^{crys}_.}\to \O$ (resp.  $\smash{\widetilde{\O}{}^{crys}_.}\to \O$). They are in particular normalized (resp. $L$-normalized) and $l_*$-acyclic. The complexes $\smash{\widetilde{Fun}{}^i_.(D)}$ are $L$-normalized and $\smash{\widetilde{Lie}_{.,T^\s}(-h)(D)}$ as well.
\finrom
\end{lem}

Proof. Everything is straightforward from Lem. \ref{cartgen} \ref{cartgeniv}, once observed that $T$ is flat over $\Sigma_\infty$ (Lem. \ref{lememb1} \ref{lememb1v}) and that the tautologically normalized $\Z/p^.$-algebra  $\smash{\widetilde{\O}{}^{crys}_.}$ is consequently flat as well (note that by Lem. \ref{lemflatpadic} \ref{lemflatpadici} the structural ring of $T_{et}$ is $\Zp$-flat, hence $\O_{crys}$ as well).
\begin{flushright}$\square$\end{flushright}

\para \label{fillim} The difference between $\smash{\widetilde{Fun}{}^i_.}$ and $\smash{Fun^i_.}$ disappears when taking limits. Let us discuss this briefly.

\begin{defn} \label{deffillim}
Let $\smash{Fun^._.}$ denote one of the three functors defined in Def. \ref{deffil}.
 We define a functor $$Fun^.:\DC(X^\s)\to Fil^{0,1}Kom(T^\s_{et},\O^{crys})$$
by setting $\smash{Fun^.(D)}:=l_*\smash{Fun^._.(D)}$.
\end{defn}

\begin{lem} \label{lemdeffillim} \debrom

\item \label{lemdeffillimi} Consider the natural morphism $l:(T^\s_{.,et},\smash{\widetilde{\cal O}{}^{crys}_.})\to (T^\s_{et},\cal O^{crys})$. The images of the vertical arrows in Lem. \ref{lemfil} \ref{lemfili} by $Rl_*$ are isomorphisms. The distinguished triangle obtained by applying $Rl_*$ to either one of the horizontal lines boils down to an exact sequence $$\xymatrix{0\ar[r]&Fun^1(D)\ar[r]&Fun^0(D)\ar[r]&Lie_{T^\s}(-h)(D)\ar[r]&0}$$ where $\smash{Lie_{T^\s}(-h)(D)}$ is the direct image of $Lie(D)(-h)$ via $\smash{(X^\s_{et},\cal O)}\to \smash{(T^\s_{et},\O^{crys})}$. Conversely,  the top line of Lem. \ref{lemfil} \ref{lemfili} can be obtained by applying $Ll^*$ to this exact sequence.

\item \label{lemdeffillimii} The objects of the complexes $\smash{Fil^i_{T^\s}(-h)(D)}$, $\smash{Fil^i\Omega^\bullet_{T^\s}(-h)(D)}$, $\smash{Fil^iL\Omega^\bullet_{T^\s}(-h)(D)}$ and $Lie_{T^\s}(-h)(D)$ come from quasi-coherent and $L$-quasi-coherent modules on $T^\s_{et}$ by scalar restriction via $\smash{\O^{crys}\to \O}$.
    \finrom
\end{lem}

Proof. \ref{lemdeffillimi} To prove the first assertion, it suffices to notice that the image of $\smash{L^1\jj^*Lie_{.,T^\s}(-h)(D)}$ by $Rl_*$ vanishes (the transition morphisms are zero).
The remaining assertions of \ref{lemdeffillimi}, \ref{lemdeffillimii} are straightforward from
 Lem. \ref{lemfil} \ref{lemfilii}.
\begin{flushright}$\square$\end{flushright}



\subsection{Twisted syntomic complexes on the \'etale site} ~~ \\

\label{tscotes} In this section we complete the construction of the twisted syntomic complexes (on the \'etale site) for Dieudonn\'e crystals over an $X^\s$ having local finite $p$-bases over $\Sigma_1$. We first define the morphism $\varphi$ in presence of a global embedding with Frobenius lift (Prop. \ref{defphi}) and we conclude using the construction of Sect. \ref{paraprelimS}.

\para \label{froblift} We begin with some basic facts about liftings of the relative Frobenius on small \'etale site.

\begin{lem} \label{lemfroblift} Consider a $p$-adic log scheme $\smash{T^\s}$ endowed  with a Frobenius lift $\smash{\tilde F_{T^\s}}$.

\debrom \item \label{lemfroblifti} The relative Frobenius $\smash{F^{(-/T_1^\s)}:id\rightarrow F^{-1}}$ on $\smash{T^\s_{1,et}}$ uniquely extends to
a natural transformation of endofunctors of $\smash{T^\s_{.,et}}$:
$$
\label{frelT}
\tilde F^{(-/T^\s_.)}:id\rightarrow \tilde F_{T_.^\s}^{-1}$$
called the \emph{lifted relative Frobenius} (attached to $\smash{\tilde F_{T^\s}}$).

\item \label{lemfrobliftii} Consider the ringed topos $(\smash{T^\s_{.,et}},\O)$. There exists a unique endomorphism $\tilde F$ of the ring $\O$ which simultaneously extends:

- the endomorphism of the structural ring of $\smash{T^\s_{.,zar}}$ defining $\smash{\tilde F_{T_.^\s}}$, and

- the Frobenius endomorphism $F:x\mapsto x^p$ of the structural ring of $\smash{T^\s_{1,et}}$.

\noindent Explicitly this endomorphism can be obtained by composing the lifted relative Frobenius $\smash{\tilde F^{(\O/T^\s_.)}}:\O\rightarrow \smash{\tilde F_{T_.^\s}^{-1}\O}$ with the functoriality morphism along $\smash{\tilde F_{T_.^\s}}$,  $\smash{\tilde F_{T_.^\s}^*}:\smash{\tilde F_{T_.^\s}^{-1}\O}\to \O$.
\finrom
\end{lem}
Proof. Everything follows easily from the fact that $\smash{T_{1,et}^\s}\to \smash{T_{k,et}^\s}$ is an equivalence.
\begin{flushright}$\square$\end{flushright}

\begin{rem} \label{remfroblift} Consider $X^\s$ in $\cal Sch^\s/\Sigma_1$ and let $\smash{\cal O^{crys}_.:=u_*\O_{X/\Sigma_.}}$ in $\smash{X^{\s,\N}_{et}}$ as usual.  \debrom

\item \label{remfroblifti} The ring $\smash{\cal O^{crys}_.}$ is endowed with a canonical \emph{Frobenius endomorphism} extending $F:\smash{\cal O_1^{crys}}\to \smash{\cal O_1^{crys}}$, $x\mapsto x^p$. It is defined by composing the relative Frobenius $\smash{F^{(\cal O^{crys}_./X^\s)}}:\smash{\cal O^{crys}_.}\to \smash{F_{X^\s}^{-1}\cal O^{crys}_.}$ with the functoriality morphism along $F_X$, $\smash{F_X^*}:\smash{F_{X^\s}^{-1}\cal O^{crys}_.}\to \smash{\cal O^{crys}_.}$.

\item \label{remfrobliftii} Assume now given $T^\s$ and $\smash{\tilde F_{T^\s}}$ as in Lem. \ref{lemfroblift} and a nilimmersion $\iota:X^\s\to T^\s$, with $X^\s$ in $\smash{\cal Sch^\s/\Sigma_1}$ and $T^\s$ in $\smash{\cal Sch^\s_p}$. With the usual abuse of notation, we have a natural ring homomorphism $\smash{\O^{crys}_.\to \O}$ in $\smash{T^\s_{.,et}}$. The Frobenius $F$ of \ref{remfroblifti} and the endomorphism $\smash{\tilde F}$ of Lem. \ref{lemfroblift} \ref{lemfrobliftii} are compatible via this homomorphism.

\item \label{remfrobliftiii} As in (\ref{defFrellin}), the lifted relative Frobenius of Lem. \ref{lemfroblift} \ref{lemfroblifti} has a linear version. More precisely it induces a morphism $\smash{\tilde F^{(-/T_.^\s)}}:\smash{(\tilde F)^*}\to \smash{\tilde F_{T_.^\s}^*}$ between endomorphisms of $(\smash{T^\s_{.,et}},\O)$ and a morphism $\smash{\tilde F^{(-/T_.^\s)}}:\smash{(F)^*\to \smash{\tilde F_{T_.^\s}^*}}$ between endomorphisms of $(\smash{T^\s_{.,et}},\smash{\O^{crys}_.})$. Both morphisms are compatible via $\smash{\O^{crys}_.\to \O}$.
\finrom
\end{rem}

\para \label{phiet} We will now show that ``Frobenius is uniquely divisible by $p$ on $\smash{\widetilde{Fil}{}^1_.}$''.

\begin{defn} \label{defFr} Let $\smash{(T^\s,\tilde F_{T^\s})}$ denote the logarithmic divided power envelope of some  $\smash{(X^\s,Y^\s,\tilde F_{Y^\s})}$ in $\smash{Emb^{\s,glob}_{F}}$. If $h\in \G(X,M_X/\Gm)$ and $(D,f,v)$ in $\DC(X^\s)$ we use the following notations. \petit

\noindent - We let $f:F^*(D(-h))\to D(-h)$ in  $\cal Crys((X^\s/\Sigma_\infty)_{crys,et},\O)$ denote $$\xymatrix{F^*(D(-h))\ar[r]_-\sim^-{\ref{lemtwlogtop}
\ref{lemtwlogtopii}}& (F^*D)(-ph)\ar[r]^-{\ref{lemtwlogtop} \ref{lemtwlogtopi}} & (F^*D)(-h)\ar[r]^-{f}& D(-h)}$$


\noindent - We let $Fr:\smash{D(-h)_{T^\s_.}\to (\smash{\tilde F})_*D(-h)_{T^\s_.}}$ in $Mod(\smash{T^\s_{.,et}},\O)$ (or
$Mod(\smash{T^\s_{.,et}},\smash{\O^{crys}_.})$) denote
$$\xymatrix{D(-h)_{T^\s_.}\ar[r]^-{\tilde F^{(-/T^\s_.)}}&
(\tilde F)_*\tilde F_{T_.^\s}^*(D(-h)_{T^\s_.})\ar[r]_-\sim^{\ref{lemcrystal2} \ref{lemcrystal2ii}}& (\tilde F)_*(F^*(D(-h)))_{T^\s_.}\ar[r]^-f& (\tilde F)_*D(-h)_{T^\s_.}}$$
%

\noindent - We let  $Fr:\smash{\Omega^\bullet_{T^\s_.}(D(-h))}\to (\smash{\tilde F})_*\smash{\Omega^\bullet_{T^\s_.}(D(-h))}$ in
$Kom(\smash{T^\s_{.,et},\O^{crys}_.})$ denote
$$\xymatrix{\Omega^\bullet_{T^\s_.}(D(-h))
\ar[r]^-{{\tilde F^{(-/T_.^\s)}}}&
(F)_*\tilde F^{*}_{T_.^\s}\Omega^\bullet_{T^\s_.}(D(-h))\ar[r]^-{(\ref{functdRcomplex})}&
(F)_*\Omega^\bullet_{T^\s_.}(F^*(D(-h)))\ar[r]^-f&(F)_*\Omega^\bullet_{T^\s_.}(D(-h))}$$

\noindent - We let  $Fr:\smash{(L(\Omega^\bullet_{T^\s_.}(D(-h))))_{T^\s_.}}\to (\tilde F)_*\smash{(L(\Omega^\bullet_{T^\s_.}(D(-h))))_{T^\s_.}}$ in
$Kom(\smash{T^\s_{.,et},\O^{crys}_.})$ denote
$$\xymatrix{(L(\Omega^\bullet_{T^\s_.}(D(-h))))_{T^\s_.}\ar[r]^-{{\tilde F^{(-/T_.^\s)}}}&
(F)_*\tilde F^{*}_{T_.^\s}(L(\Omega^\bullet_{T^\s_.}(D(-h))))_{T^\s_.}\ar[r]^-{(\ref{functLdR})}&
(F)_*(L(\Omega^\bullet_{T^\s_.}(F^*(D(-h)))))_{T^\s_.}\ar[d]^-f\\
&&(F)_*(L(\Omega^\bullet_{T^\s_.}(D(-h))))_{T^\s_.}}
$$
\end{defn}

\begin{prop} \label{defphi} Let $\smash{(X^\s,Y^\s,\iota, \tilde F,h)}$ in $\smash{Emb^{\s,glob,lfpb}_{F,div}}$ and denote $T^\s$ the logarithmic divided power of $\iota$. Let  $\smash{\widetilde{Fun}{}^._.}$ denote one of the three functors defined in Def. \ref{deffilmod} and let $D\in \DC(X^\s)$. There exists a unique morphism $\varphi$ in $\smash{Kom(\smash{T^\s_{.,et},\widetilde \O_.^{crys}})}$ rendering the following square commutative: $$\xymatrix{{\widetilde{Fun}{}^1_.(D)}\ar[r]^-\varphi\ar[d]^1&(\tilde F)_*{\widetilde{Fun}{}^0_.(D)}\ar[d]^p\\
{\widetilde{Fun}{}^0_.(D)}\ar[r]^-{Fr}&(\tilde F)_*{\widetilde{Fun}{}^0_.(D)}}$$ The morphism $\varphi$ is functorial with respect to $D$ and
compatible with the base change morphisms arising from morphisms in $\smash{Emb^{\s,glob,lfpb}_{F,div}}$ (Rem. \ref{remdeffilmod}).
\end{prop}
Proof. Unicity, functoriality and compatibility to base change follows from Lem. \ref{cartgen} \ref{cartgeniii} and Lem. \ref{lemfil} \ref{lemfilii}.  Let us prove existence. By Lem. \ref{cartgen} \ref{cartgenii}, we have the bottom exact sequence in the following commutative diagram:
$$\xymatrix{&&Fun^1_{.+1}(D)\ar[r]\ar[d]^-{Fr\circ 1}\ar@{-->}_-?[dl]&
\iota_{1,.+1,*}Fun^1_1(D)\ar[d]^-{Fr\circ 1}\\
0\ar[r]&(\tilde F)_*\iota_{.,.+1,*}\widetilde{Fun}{}^0_.(D)\ar[r]^-{p}&(\tilde F)_*Fun^0_{.+1}(D)\ar[r]&
(\tilde F)_*\iota_{1,.+1,*}Fun^0_1(D)\ar[r]&0}$$
We want to prove the existence of the dotted arrow marked $?$, since the desired $\varphi$ will follow by applying $\smash{\iota_{.,.+1}^*}$ to it (note that $\smash{\iota_{.,.+1,*}}$ is fully faithful and commutes to $\smash{(\tilde F)_*}$). It is thus sufficient to prove that the morphism \begin{eqnarray}\label{Frmodp}Fr\circ
1:Fun^1_1(D)\to (F)_*Fun^0_1(D)\hbox{ in $Kom(T^\s_{1,et},\widetilde\O^{crys}_1)$}\end{eqnarray} vanishes. Let us examine the three cases of
Def. \ref{defFr}. Using the definition of the base change morphisms (\ref{functdRcomplex}) and (\ref{functLdR}), it is clear that (\ref{Frmodp})
vanishes in degree $\ge 1$ in all cases.
It thus remains to prove that \begin{eqnarray}\label{Fr0D}Fr_1:\smash{D(-h)_{T_1^\s}\to (F)_*D(-h)_{T_1^\s}}\hbox{ and}\\
\label{Fr0L} Fr_1:\smash{(L(D(-h)_{T^\s_.}))_{T_1^\s}}\to \smash{ (F)_* (L(D(-h)_{T^\s_.}))_{T_1^\s}}\end{eqnarray} respectively vanish on
$\smash{(Fil^1D)(-h)_{T_1^\s}}$ and  $\smash{(Fil^1L(D_{T^\s_.}))(-h)_{T^\s_1}}$. Thanks to  Lem. \ref{injtwist2} \ref{injtwist2ii}, we see that
$$M(-h)_{T_1^\s}\hookrightarrow M_{T_1^\s}$$ if $M$ is  either one of the locally free crystals $D$ or $L(D_{T^\s_.})$. This monomorphism is
compatible with $Fr_1$ and we may thus assume that $h=1$.

Consider the following commutative diagram on $((X^\s/\Sigma_1)_{crys},\O)$: $$\xymatrix@1{\ar @<+2pt> `u[r] `[rrrr]^{\overline Fr} [rrrr]\ar@<-4pt> `l[d]
`[dd] [rrdd]^-{can_D}\overline D\ar[d]^-{adj}\ar[rr]^{F^{(\overline D/X^\s)}}&& (F)_*F^*\overline D\ar[rr]^(.4){(F)_*\overline f} &&(F)_*\overline D
\\ i_*i^*\overline D\ar[rr]^-{i_*i^*(C^{-1})} &&i_*\phi_*F^*\overline D\ar[u]_{nat}
\ar[rr]^-{i_*\phi_*\overline f}\ar[d]&&i_*\phi_*\overline D\ar[u]_{nat}\\
&& i_*\phi_*Coker\overline v\ar[rru]_-{i_*\phi_*\overline f}}$$ where $nat:i_*\phi_*\to (F)_*$ is deduced from the isomorphism $\phi^*\simeq
(F)^*i_*$ (Lem. \ref{compatcartier} \ref{compatcartierii}). Since $Fil^1\overline D=Ker\, can_D$, it must be in the kernel of  $\overline Fr$. This ends the
proof in the case $M=D$ since $Fr_1$ coincides with the realization of $\overline Fr$ at $\smash{T_1^\s}$.

The case $M=L\smash{(D_{T^\s_.})}$ is proven similarly using the following commutative diagram of $Mod((X^\s/\Sigma_1)_{crys,et},\O)$:

$$\xymatrix@1{\ar @<+2pt> `u[r] `[rrrrrrr]^{\overline Fr} [rrrrrrr]\ar@<-4pt> `l[d] `[ddd] [rrrrddd]^-{can_L}f_{T_1^\s,*}f_{T_1^\s}^*\overline D\ar[dd]^{(\ref{defalpha})}\ar[rrrr]_{ch\circ F^{(
f_{T_1^\s,*}f_{T_1^\s}^*\overline D/X^\s)}}&&&~~\hspace{1.5cm}~~& (F)_*f_{T_1^\s,*}f_{T_1^\s}^*F^*\overline D\ar[rrr]_{(F)_*f_{T_1^\s,*}f_{T_1^\s}^*\overline f} &&~~\hspace{1cm}~~&(F)_*f_{T_1^\s,*}f_{T_1^\s}^*\overline D\\
&&&&(F)_*F^*\overline D\ar[rrr]^{(F)_*\overline f}\ar[u]_{adj}&&&(F)_*\overline D\ar[u]_{adj}
\\ i_*i^*\overline D\ar[rrrr]^-{i_*i^*(C^{-1})} &&&&i_*\phi_*F^*\overline D\ar[u]_{nat}
\ar[rrr]^-{i_*\phi_*\overline f}\ar[d]&&&i_*\phi_*\overline D\ar[u]_{nat}\\
&&&& i_*\phi_*Coker\overline v\ar[rrru]_-{i_*\phi_*\overline f}}$$
\begin{flushright}$\square$\end{flushright}

\begin{defn} \label{defSonephi} As in Prop. \ref{defphi}, let $T^\s$ denote the logarithmic divided power envelope of some $\smash{(X^\s,Y^\s,\iota, \tilde F,h)}$ in $\smash{Emb^{\s,glob,lfpb}_{F,div}}$. For $D\in \DC(X^\s)$, we define three complexes of $(1,\varphi)$-modules of $\smash{(T^\s_{.,et},\widetilde \O^{crys}_.)}$ as follows.  $$\begin{array}{rcl}\cal S^{1,\varphi}_{et,.,T^\s}(-h)(D)&:=&(\widetilde{Fil}{}^1_{.,T^\s}(-h)(D),\widetilde{Fil}{}^0_{.,T^\s}(-h)(D),1,\varphi)\\
\cal S\Omega^{\bullet,1,\varphi}_{et,.,T^\s}(-h)(D)&:=&(\widetilde{Fil}{}^1\Omega^\bullet_{.,T^\s}(-h)(D),\widetilde{Fil}{}^0\Omega^\bullet_{.,T^\s}(-h)(D),1,\varphi)\\
\cal SL\Omega^{\bullet,1,\varphi}_{et,.,T^\s}(-h)(D)&:=&(\widetilde{Fil}{}^1L\Omega^\bullet_{.,T^\s}(-h)(D),
\widetilde{Fil}{}^0L\Omega^\bullet_{.,T^\s}(-h)(D),1,\varphi)\end{array}$$
where the $\varphi$'s in the right hand terms denote the ones defined in Prop. \ref{defphi}.

By functoriality with respect to $D$ and $\smash{(X^\s,Y^\s,\iota, \tilde F,h)}$, this defines three colax morphisms between contravariant pseudo functors on $\smash{Emb^{\s,glob,lfpb}_{F,div}}$:
 \begin{eqnarray}\label{SetT}\cal S^{1,\varphi}_{et,.},\, \cal S\Omega^{\bullet,1,\varphi}_{et,.},\, \cal SL\Omega^{\bullet,1,\varphi}_{et,.}:\DC(-)\to Kom^{1,\varphi}((-)^{dp}_{.,et},\widetilde \O^{crys}_.)\end{eqnarray}
Here $\smash{Kom^{1,\varphi}(T^\s_{[.],.,et},\widetilde \O^{crys}_.)}$ denotes the category of complexes of the category of $(1,\varphi)$-modules defined in Lem. \ref{1phimod} for the ringed topos $\smash{(T^\s_{[.],.,et},\widetilde \O^{crys}_.)}$ and the endomorphism of $\smash{\widetilde \O^{crys}_.}$ induced by the Frobenius endomorphism of $\smash{\O^{crys}_.}$ (Rem. \ref{remfroblift} \ref{remfroblifti}).
\end{defn}


We will now apply the globalizing construction explained in Sect. \ref{paraprelimS} with $\cal F=\cal DC(-)$ and $(\cal T,A)=((-)^{1,\varphi,\N}_{et},\smash{\widetilde \O^{crys}_.})$. Note that by the topological invariance of the \'etale site we have, with the notations of \emph{loc. cit.}, an equivalence $\cal K^{emb}\simeq   \smash{Kom^{1,\varphi}((-)^{dp}_{.,et},\widetilde \O^{crys}_.)}$ between contravariant pseudo-functors on
$\smash{Emb^{\s,glob,lfpb}_{F,div}}$. Letting $\cal S$ denote one of the three colax morphisms (\ref{SetT}) we denote $\cal S^{emb}:\cal F^{emb}\to \cal K^{emb}$ (see. (\ref{Semb})) the colax morphism induced by this equivalence, and $\cal S^{loc}:\cal F^{loc}\to \cal K^{loc}$ the colax morphism obtained in (\ref{cons2}), whose data is equivalent to a collection of functors
\begin{eqnarray}\label{SHR}\cal S(-h_{|U^\s_{[.]}}):\DC(U^\s_{[.]})\to Kom^{1,\varphi}(T^\s_{[.],.,et},\widetilde \O^{crys}_.)\end{eqnarray} indexed by the objects $\smash{Y^\s_{[.]}}=(\smash{U^\s_{[.]}/X^\s,Y^\s_{[.]},\iota_{[.]},\tilde F_{[.]}})$ of $\smash{Diag(HR^{\s,et,lfpb}_{F})}$ equipped with an effective logarithmic divisor $h$ on $X^\s$, together with a collection of base change morphisms, indexed by the morphisms of $\smash{Diag(HR^{\s,et,lfpb}_{F,div})}$, and satisfying the composition constraint. 

%


\begin{lem} \label{lemcompS} \debrom \item \label{lemcompSi} The three colax morphisms $\cal S^{loc}$ just defined are canonically related by natural transformations above $Diag(HR^{\s,et,lfpb}_{F,div})$ as follows:
$$\xymatrix{(\cal S^{1,\varphi}_{et,.})^{loc}\ar[r]^-{\alpha'}&  (\cal SL\Omega^{\bullet,1,\varphi}_{et,.})^{loc}&\ar[l]_-{\alpha''} (\cal S\Omega^{\bullet,1,\varphi}_{et,.})^{loc}}$$
Fix a locally embeddable diagram $(X^\s,h)$ of $\smash{\cal Sch^{\s,lfpb}_{div}/\Sigma_1}$. Let $\smash{T^\s_{[.]}}$ denotes the logarithmic divided power envelope of some object $\smash{Y^\s_{[.]}}$ in  $HR^{\s,*}_{F}(X^\s)$ (Lem. \ref{lememb2} \ref{lememb2ii}, \ref{lememb2iv}).

- If $*=et$ and $\smash{D_{[.]}}\in \DC(\smash{U^\s_{[.]}})$, the morphism $\alpha'$ induces an isomorphism  $$\cal
S^{1,\varphi}_{et,.,T^\s_{[.]}}(-h_{|U^\s_{[.]}})(D_{[.]})\simeq \cal SL\Omega^{\bullet,1,\varphi}_{et,.,T^\s_{[.]}}(-h_{|U^\s_{[.]}})(D_{[.]})\hbox{ in
$D(Mod^{1,\varphi}(T_{[.],.,et}^{\s},\widetilde\O^{crys}_.))$}$$

- If $*=crys$ and $D\in \DC(X^\s)$, the morphism $\alpha''$ induces an isomorphism $$R f_{T^\s_{[.],.},*}\cal
S\Omega^{\bullet,1,\varphi}_{et,.,T^\s_{[.]}}(-h_{|U^\s_{[.]}})(D_{|U^\s_{[.]}})\simeq R f_{T^\s_{[.],.},*}\cal
SL\Omega^{\bullet,1,\varphi}_{et,.,T^\s_{[.]}}(-h_{|U^\s_{[.]}})(D_{|U^\s_{[.]}}) \hbox{ in $D(Mod^{1,\varphi}(X^{\s,\N},\widetilde\O^{crys}_.))$}$$

\item \label{lemcompSii} Let $\cal S$ denote one of the three functors (\ref{SetT}) and fix a diagram $(X^\s,h)$ of $\smash{\cal Sch^{\s,slfpb}_{div}/\Sigma_1}$ as well as $D\in \DC(X^\s)$.
If  $\smash{T_{[.]}'^\s\to T_{[.]}^\s}$ is the logarithmic divided power envelope of some morphism  $g$ of $HR^{\s,*,flpb}_{F,div}$ above $f:(X'^\s,h')\to (X^\s,h)$, the base change morphism for (\ref{SglobC}) reads:
$$Rf_{T^\s_{[.],.},*}\cal S_{T^\s_{[.]}}(-h_{|U^\s_{[.]}})(D_{|U^\s_{[.]}})\to Rf_*R f_{T^\s_{[.],.},*}\cal S_{T'^\s_{[.]}}(-h'_{|U'^\s_{[.]}})((f^*D)_{|U'^\s_{1,[.]}})\hbox{ in $D(Mod^{1,\varphi}(X^{\s,\N},\widetilde\O^{crys}_.))$}$$

It is an isomorphism if $f$  (but not necessarily $g$) is the identity and one of the following conditions hold:

- $\cal S$ is either $\cal S^{1,\varphi}_{et,.}$ or $\cal SL\Omega^{\bullet,1,\varphi}_{et,.}$  and $*=crys$.

- $\cal S$ is $\cal S\Omega^{\bullet,1,\varphi}_{et,.}$ and $*=et$.
\finrom \end{lem}

\noindent Proof. \ref{lemcompSi} Putting together Lem. \ref{lemfil} and Prop. \ref{cdOmega} gives morphisms $$\xymatrix{\widetilde{Fil}{}^i_{.,T^\s_{[.]}}(-h_{|U_{[.]}})&\ar[l] \widetilde{Fil}{}^iL\Omega^\bullet_{.,T^\s_{[.]}}(-h_{|U_{[.]}})\ar[r]&
\widetilde{Fil}{}^iL\Omega^\bullet_{.,T^\s_{[.]}}(-h_{|U_{[.]}})}$$ which are compatible with $1$ and $Fr$. They are compatible with $\varphi$ as well by Lem. \ref{cartgen} \ref{cartgeniii}. Whence the claimed morphisms $\alpha'$ and $\alpha''$. An immediate d\'evissage (using Lem. \ref{1phimod} \ref{1phimodiii} and Lem. \ref{lemdeffillim} \ref{lemdeffillimi}) together with Lem. \ref{cdOmega} and Rem. \ref{remcomp} show that $\alpha'$ is always an isomorphism and that $\alpha''$ becomes an isomorphism once we apply $\smash{Rf_{T^\s_{[.],.},*}}$ to it.
Similar arguments apply in \ref{lemcompSii}.
\begin{flushright}$\square$\end{flushright}


The above lemma shows that the assumptions of Lem. \ref{newglu} are satisfied. From this we get immediately the following 
result which roughly says that:

- ${Rf_{T^\s_{[.],.},*}\cal S_{T^\s_{[.]}}(-h_{|U^\s_{[.]}})(D_{|U^\s_{[.]}})}$ is essentially independant of the chosen $\smash{Y^\s_{[.]}}$, 

- the base change morphism in Lem. \ref{lemcompS} \ref{lemcompSii} does not depend on the choice of $g$, and still exists even if $g$ does not.

\begin{prop} \label{indepSet} \debrom \item \label{indepSeti} Up to a unique isomorphism there exists a unique couple $(\cal S^{1,\varphi}_{et,.},\alpha)$ such that the following conditions are satisfied:

- $\cal S^{1,\varphi}_{et,.}$ is a $\smash{\cal B_0^\s}$-functor $\DC(-)_{cof}\to \smash D(Mod^{1,\varphi}((-)^{\N},\widetilde\O^{crys}_.))'_{fib}$ above the category  defined in Def. \ref{defB0}. In other terms,  it is a collection of functors \begin{eqnarray}\label{SX}\cal S^{1,\varphi}_{et,.,X^\s}(-h):\DC(X^\s)\to
D(Mod^{1,\varphi}(X^{\s,\N},\widetilde\O^{crys}_.))\end{eqnarray} indexed by the $(X^\s,h)$'s of $\smash{\cal B^\s_0}$ together
with a canonical family of base change morphisms $$\cal S^{1,\varphi}_{et,.,X_2^\s}(-h_2)\to Rf_*\cal S^{1,\varphi}_{et,.,X_1^\s}(-h_1)f^*$$
indexed by the morphisms $f:(X_1^\s,h_1)\to (X_2^\s,h_2)$ of $\smash{\cal B^\s_0}$  and satisfying the cocycle condition.

- $\alpha$ is a collection of functorial isomorphisms indexed by the diagrams of $\smash{HR^{\s,et,lfpb}_{F,div}}$: $$\alpha_{Y^\s}:\cal
S^{1,\varphi}_{et,.,X^\s}(-h)(D)\to Rf_{T^\s_{[.],.},*}\cal S\Omega^{\bullet,1,\varphi}_{et,.,T^\s_{[.]}}(-h_{|U^\s_{[.]}})(D_{|U^\s_{[.]}})$$ A morphism
$f:(X_1^\s,h_1)\to (X_2^\s,h_2)$ induces a base change morphism on the left side. If $g:\smash{Y_{1,[.]}^\s\to Y^\s_{2,[.]}}$ is furthermore a morphism of
$Diag(HR^{\s,et}_{F})$ above $X^\s_1\to X^\s_2$ then it induces a base change morphism on the right side. Both base change morphisms are compatible via $\alpha$.

\item \label{indepSetii} If $(X^\s,h)$ is in $\smash{\cal B^\s_0}$ and $\smash{Y_{[.]}^\s}$ is in $HR^{\s,crys}_{F}(X^\s)$ then \ref{indepSeti} and Lem. \ref{lemcompS} \ref{lemcompSi} induce canonical isomorphisms as follows $$\begin{array}{c}\alpha'\alpha:\cal S^{1,\varphi}_{et,.,X^\s}(-h)(D)\mathop\to\limits^\sim Rf_{T^\s_{[.],.},*}\cal SL\Omega^{\bullet,1,\varphi}_{et,.,T^\s_{[.]}}(-h_{|U^\s_{[.]}})(D_{|U^\s_{[.]}})\\
\alpha''^{-1}\alpha'\alpha:\cal S^{1,\varphi}_{et,.,X^\s}(-h)(D)\mathop\to\limits^\sim Rf_{T^\s_{[.],.},*}\cal
S^{1,\varphi}_{et,.,T^\s_{[.]}}(-h_{|U^\s_{[.]}})(D_{|U^\s_{[.]}})\end{array}$$  Those morphisms are compatible with the base change morphism attached to
morphisms of $Diag(HR^{\s,crys,lfpb}_{F,div})$ in the same sense as above. \finrom
\end{prop}
\begin{flushright}$\square$\end{flushright}

\begin{defn} \label{defSX} If $(X^\s,h)$ is in $\smash{\cal B^\s_0}$ and $D$ is in $\DC(X^\s)$, we define the \emph{syntomic complex of $D$ twisted by $(-h)$ on the \'etale site}  as follows: $$\cal S_{et,.,X^\s}(-h)(D):=R\varpi_*\cal S^{1,\varphi}_{et,.,X^\s}(-h)(D)\hbox{ in $D(\smash{X^{\s,\N,1}},\smash{\widetilde\O^{crys,F=1}_.})$}$$ where $\smash{\cal S^{1,\varphi}_{et,.,X^\s}(-h)}$ is the functor (\ref{SX}) and $\varpi:(\smash{X^{\s,\N,1,\varphi}_{et}},\smash{\widetilde\O^{crys}_.})\to (\smash{X^{\s,\N}_{et}},\smash{\widetilde\O^{crys,F=1}_.})$ denotes the projection morphism of Lem. \ref{1phimod} \ref{1phimodii}. Finally, we set $$\cal S_{et,X^\s}(-h)(D):=Rl_*\cal S_{et,.,X^\s}(-h)(D)\hbox{ in $D(\smash{X^{\s}},\smash{\widetilde\O^{crys,F=1}})$}$$
where $\smash{\widetilde\O^{crys,F=1}}:=l_*\smash{\widetilde\O_.^{crys,F=1}}$.
\end{defn}

\begin{lem} \label{Sisexact} The functor (\ref{SX}) sends short exact sequences of Dieudonn\'e crystals to distinguished triangles.
\end{lem}

Proof. This follows immediately from Prop. \ref{liesurj} \ref{liesurjiii}. \begin{flushright}$\square$\end{flushright}

\begin{prop} \label{dtet} Consider $(X^\s,h)$ in $\smash{\cal B^\s_0}$.

\debrom
\item \label{dteti}  If $\smash{T^\s_{[.]}}$ is the logarithmic divided power envelope of some $\smash{Y^\s_{[.]}}$ in $\smash{HR^{\s,et}_{F}(X^\s)}$, there are  canonical distinguished triangles in $D(\smash{X^{\s,\N,1}},\smash{\widetilde\O^{crys,F=1}_.})$
{\small$$\begin{array}{c} Rf_{T^\s_{[.],.},*}\widetilde{Fil}{}^1\Omega^{\bullet}_{.,T^\s_{[.]}}(-h_{|U^\s_{[.]}})(D_{|U^\s_{[.]}})
\mathop\to\limits^{1} Rf_{T^\s_{[.],.},*}\widetilde{Fil}{}^0\Omega^{\bullet}_{.,T^\s_{[.]}}(-h_{|U^\s_{[.]}})(D_{|U^\s_{[.]}})
\to Rf_{T^\s_{[.],.},*}\widetilde{Lie}_{.,T^\s_{[.]}}(D_{|U^\s_{[.]}})(-h_{|T^\s_{[.]}})\mathop\to\limits^{+1}\\
\cal S_{et,.,X^\s}(-h)(D)\to Rf_{T^\s_{[.],.},*}\widetilde{Fil}{}^1\Omega^{\bullet}_{.,T^\s_{[.]}}(-h_{|U^\s_{[.]}})(D_{|U^\s_{[.]}})
\mathop\to\limits^{1-\varphi} Rf_{T^\s_{[.],.},*}\widetilde{Fil}{}^0\Omega^{\bullet}_{.,T^\s_{[.]}}(-h_{|U^\s_{[.]}})(D_{|U^\s_{[.]}})
\mathop\to\limits^{+1}\end{array}$$ }

\item \label{dtetii} If  $\smash{Y^\s_{[.]}}$ is in fact in $\smash{HR^{\s,crys}_{F}}$, there are analogous distinguished triangles with ``$\smash{\widetilde{Fil}{}^iL\Omega^{\bullet}}$'' or ``$\smash{\widetilde{Fil}{}^i}$'' instead of ``$\smash{\widetilde{Fil}{}^i\Omega^{\bullet}}$''. The resulting couples of triangles are moreover canonically isomorphic.
\finrom 
\end{prop}
Proof. The first distinguished triangle is provided by Lem. \ref{lemfil} \ref{lemfili} and the second by  Lem. \ref{1phimod} \ref{1phimodiv}.  \begin{flushright}$\square$\end{flushright}

\begin{rem} Here again the situation becomes clearer after passing to the limit. Indeed, from Prop. \ref{cdOmega} \ref{cdOmegai}, we find that the distinguished triangles of Prop. \ref{dtet} \ref{dtetii} in the case $``\smash{\widetilde{Fil}{}^i}''$ become
$$\begin{array}{c}\xymatrix{Ru_*((Fil^1D)(-h))\ar[r]^-{1}&Ru_*(D(-h))\ar[r]&Lie(D)(-h)\ar[r]^-{+1}&}\\
\xymatrix{\cal S_{et,X^\s}(-h)(D)\ar[r]&Ru_*((Fil^1D)(-h))\ar[r]^-{1-\varphi}&Ru_*(D(-h))\ar[r]^-{+1}&}\end{array}$$
in $D(X^\s_{et},\O^{crys,F=1})$. Conversely, it follows from Lem. \ref{lemfil} \ref{lemfilii}, that the triangles of Prop. \ref{dtet} can be retrieved from these
ones by applying $Ll^*$ (note that $\smash{\widetilde \O^{crys,F=1}_.}$ is isomorphic to $\smash{\Z/p^.\otimes \O^{crys,F=1}\simeq \Z/p^.}$, as will be recalled in Lem. \ref{lemcomp2} \ref{lemcomp2iii}).
\end{rem}

\subsection{Syntomic complexes on the syntomic site}
\label{scotss} ~~ \\

%
%

We explain two constructions of the syntomic complex using the syntomic topology for schemes without log structures. The first one is global, ie.
does not involve local embeddings or cohomological descent. The second one is local and involves the linearized de Rham complex. It will serve as a
bridge to the previous constructions. We will consider the following situations in parallel until Def. \ref{defSonephisyn}: \petit

- The \emph{global situation}, where $X$ is a separated scheme with local finite $p$-bases over $\Sigma_1$.

- The \emph{local situation}, where  $(\iota:X\to Y,F)$ is an object of $Emb^{glob,lfpb}_F$ and $T$ denotes the divided power envelope of $\iota$.

\petit

%
%
%
%


\petit

Let us gather some technical facts in a lemma.

\begin{lem} \label{lemsyn} Consider the global situation.

\debrom \item \label{lemsyni} Let $\epsilon:(\smash{X^\N_{syn},\O}) \to (\smash{X^\N_{et},\O})$ denote the natural morphism. If $M_.$ is a module on $(X^\N_{et},\O)$ (resp. a quasi-coherent local crystal on  
$((X/\Sigma_.)_{crys,et}/T_{.},\O)$)  then the following morphisms are invertible $$\begin{array}{rcl}  \epsilon^*i_*M_.&\to &i_*\epsilon^*M_.\\
\epsilon^*f_{T_{.},*}M_.&\to& f_{T_{.},*}\epsilon^*M_.\end{array}$$
where the first (resp. second) $\smash{f_{T_.}}$ denotes the localisation morphism attached to $T$ as in (\ref{descentdiaget}).

\item \label{lemsynii} Let $\epsilon:((X/\Sigma_.)_{crys,syn},\O)\to ((X/\Sigma_.)_{crys,et},\O)$ denote the natural morphism.  If $M_.$ is a crystal of $((X/\Sigma_.)_{crys,et},\O)$ then  $\epsilon^*M_.$ is a crystal of $((X/\Sigma_.)_{crys,syn},\O)$. If $M_.$ is furthermore quasi-coherent, locally free or locally free of finite type then the same is true for $\epsilon^*M$. In that case, $M_.\to R\epsilon_*\epsilon^*M_.$ is moreover an isomorphism.

\item \label{lemsyniii} If $M_.$ is a quasi-coherent crystal of  $((X/\Sigma_.)_{crys,syn},\O)$, then it is acyclic for the functor $u_*$. The following natural morphism (induced by $id\to i_*i^*$) is moreover an epimorphism in $Mod(X_{syn},\O)$: $$u_*M_.\to i^{*}M_.$$

\item \label{lemsyniv} The ring $\O^{crys}_k$ of $X^\N_{syn}$ is flat over $\Z/p^k$ and the natural morphism $\smash{\O^{crys}_{k+1}/p^k}\to \smash{\O^{crys}_k}$ is invertible. More generally, if $M_.$ is a locally free crystal of \break $((X/\Sigma_.)_{crys,syn},\O)$ then $u_*M_.$ is flat over $\Z/p^.$ and $$(u_{*}M_{.+1})/p^.\simeq u_*M_.$$
    \finrom
\end{lem}
Proof. \ref{lemsyni} Let us explain the first isomorphism. First, we note that we have an isomorphism $\epsilon^{-1}i_*\simeq i_*\epsilon^{-1}$ for
Abelian sheaves expressing the compatibility of $u$ and $\epsilon$. The claimed analogous morphism for modules will follow formally if we prove that
$\epsilon^*\Ga\to \Ga$ is an isomorphism, ie. that $\epsilon^*I\twoheadrightarrow I$ (here $I$ denotes the canonical ideal in the structural ring of
the crystalline sites). This, in turn, can be checked easily, using that each affine object $(U',T')$ of $crys_{syn}(X/\Sigma_\infty)$ admits a morphism to some
affine $(U,T)$ in $crys_{et}(X/\Sigma_\infty)$, where $T$ is the divided power envelope of $U$ inside a polynomial algebra of the form
$\Z/p^k[x_\alpha,y_\beta]$ and where the $x_{\alpha}$  (resp. $y_\beta$) are sent to generators of the algebra of $U$ (resp. the ideal of $U'$ inside
$T'$).


The second isomorphism is a formal consequence of Lem. \ref{lemft} since $\epsilon_*$ is fully faithful on the category of quasi-coherent crystals (Lem. \ref{crysepsilon} \ref{crysepsiloni}) and of quasi-coherent local crystals (easy variant of Lem. \ref{crysepsilon} \ref{crysepsiloni}).

\ref{lemsynii} This is a repetition of Lem. \ref{crysepsilon} \ref{crysepsiloni}, \ref{crysepsilonii}.


\ref{lemsyniii} To prove the first statement it is sufficient to show that the syntomic sheaf associated to the presheaf $U\mapsto
H^q((U/\Sigma_k)_{crys,syn},M_{k|U})$ vanishes for all $k\ge 1$ and $q\ge 1$. We have isomorphisms
$$\begin{array}{rcll}H^q((U/\Sigma_k)_{crys,syn},M_{k|U})&\simeq &H^q((U/\Sigma_k)_{crys,et},\epsilon_* (M_{k|U}))&\hbox{by \ref{lemsynii}} \\ &\simeq &
H^q(U_{et},\Omega^\bullet_{D(U,Y_k)/\Sigma_k}(\epsilon_* (M_{k|U})))&\hbox{if $Y\in Emb^{glob}(U)$}\end{array}$$ by Prop. \ref{dRlocal} \ref{dRlocalii} and the second
isomorphism is functorial with respect to morphisms in $Emb^{glob}$. Now if $U$ and $Y$ are affine, the latter is a subquotient of
$$\G(U_{et},\Omega^q_{D(U,Y_k)/\Sigma_k}(\epsilon_* (M_{k|U})))\simeq M_k(U,D(U,Y_k)) \otimes_{\O(Y_k)}\Omega^q_{Y_k/\Sigma_k}(Y_k)$$ It is thus
sufficient to notice that any global section $\omega=fdg$  of $\smash{\Omega^q_{Y_k/\Sigma_k}}$ vanishes when restricted to
$\smash{Y'_k:=Spec_{Y_k}(\O[x]/(x^{p^k}-g))}$, which is a syntomic covering of $Y_k$. The second statement is proven similarly.

\ref{lemsyniv} Let us prove that $u_*M$ is flat over $\Z/p^k$ as soon as $M$ is a locally free crystal (not necessarily of finite type) on
$((X/\Sigma_k)_{crys,syn},\O)$. We will show that $\smash{u_*M\otimes_{\Z/p^k}(-)}$ preserves monomorphisms $N\hookrightarrow N'$ of $\Z/p^k$ modules
on $X_{crys,syn}$. Since the question is local, we can always assume that $X$ affine and choose a lifting $\smash{\widetilde X}$ with $p$-bases over
$\Z/p^k$. Let $N'':U\mapsto N'(U)/N(U)$ denote the cokernel presheaf. It is sufficient to show that the sheaf associated to the following presheaf
vanishes: $$K:U\mapsto Tor^{\Z/p^k}_1(u_*M(U),N''(U))$$
We will prove that  for any $U$ affine and $s\in K(U)$, there is a syntomic covering $U'\to U$ killing $s$. Let $U/X$ syntomic with $U$ affine and
choose a transversally regular  immersion $U\to \mathbb A^n_X$ above $X$. Let us choose a transversally regular sequence $(\underline x)$ above $X$
defining this immersion. Then, any lift $\smash{\tilde{\underline x}}$ of $\underline x$ to $Y:=\smash{\mathbb A^n_{\widetilde X}}$ is transversally
regular (reduce to the case where $\smash{Y}$ is noetherian by  \cite{EGA4-IV} Prop. 19.8.2 and then use \cite{Mi1} I, Rem. 2.6 (d)). Let us choose one such lift and denote
$\smash{\widetilde U}$ the corresponding syntomic $\smash{\widetilde X}$-scheme.
By \cite{Be3} Prop. 1.5.3 (i) with $m=0$, the divided power envelope
$D(\smash{\widetilde U},Y)$ is flat over $\smash{\widetilde X}$. Let us emphasize that $D(\smash{\widetilde U},Y)$ coincides with $D(U,Y)$ since all
divided powers are intended to be compatible with $p$  (see \emph{loc. cit.} 1.3.1). Applying Prop. \ref{dRlocal} \ref{dRlocalii}, forming derived global
sections and truncating we get a distinguished triangle of $\Z/p^k$-modules: $$\xymatrix{u_*M(U)\ar[r]&
\Gamma(U_{et},\Omega_{D(U,Y)/\Sigma_k}^\bullet(\epsilon_*(M_{|U})))\ar[r]&\tau_{\ge
1}\Gamma(U_{et},\Omega_{D(U,Y)/\Sigma_k}^\bullet(\epsilon_*(M_{|U})))\ar[r]^-{+1}&}$$ Since the non zero objects of the middle term are flat over
$\Z/p^k$ and placed in positive degrees this induces  $$Tor^{\Z/p^k}_1(u_*M(U),N''(U))\simeq Tor^{\Z/p^k}_2(\tau_{\ge
1}\Gamma(U_{et},\Omega_{D(U,Y)/\Sigma_k}^\bullet(\epsilon_*(M_{|U}))),N''(U))$$ This Abelian group thus admits a finite (here we use that $Y$ has a finite
$p$-basis) filtration whose graduations are subquotients of $$Tor^{\Z/p^k}_{2+q}(H^q(U_{et},\Omega_{D(U,Y)/\Sigma_k}^\bullet(\epsilon_*(M_{|U}))),N''(U))$$ with
$q\ge 1$.  Applying repeatedly the argument of \ref{lemsyniii}, produces a syntomic covering $Y'/Y$ such that the image of $s$ vanishes in
$$Tor^{\Z/p^k}_2(\tau_{\ge 1}\Gamma(U'_{et},\Omega_{D(U',Y')/\Sigma_k}^\bullet(M_{|U'})),N''(U'))$$ where $U':=U\times_YY'$. The result follows by
functoriality of the above distinguished triangle with respect to $(U\to Y)$.
\begin{flushright}$\square$\end{flushright}



\begin{defn} \label{deffilsyn} Consider the global (resp. local) situation and let $D\in \DC(X)$.

\debrom \item \label{deffilsyni} We define a quasi-coherent module on $(X_{syn},\O)$  as follows
$$\begin{array}{rcl}Lie^{syn}(D)&:=&\epsilon^*Lie(D)\end{array}$$

\item \label{deffilsynii} We define a quasi-coherent crystal on  $\smash{((X/\Sigma_\infty)_{crys,syn},\O)}$ as follows
$$\begin{array}{crclc}&Fil^0D^{syn}=D^{syn}&:=&\epsilon^*D&\\
\hbox{(resp. }&Fil^0L^{syn}\Omega^{\bullet}_{T_.}(D)=L^{syn}\Omega^{\bullet}_{T_.}(D)&:=
&\epsilon^* L\Omega^{\bullet}_{T_.}(D)&\hbox{ )}\end{array}$$

\item \label{deffilsyniii} We define a sub-module of the crystal just defined above as follows $$\begin{array}{crclc}&Fil^1D^{syn}&:=&Ker(can^{syn}_D:D^{syn}\to i_*Lie^{syn}(D))&\\
\hbox{(resp. }&Fil^1L^{syn}\Omega^{\bullet}_{T_{.}}(D)&:= &Ker(can^{syn}_L:L^{syn}\Omega^{\bullet}_{T_{.}}(D)\to i_*Lie^{syn}(D))&\hbox{
)}\end{array}$$ where $can^{syn}_D$ (resp. $can^{syn}_L$) is the arrow deduced from the corresponding one in Prop. \ref{complin} \ref{complini} using the base change
morphism of Lem. \ref{lemsyn} \ref{lemsyni}. \finrom
\end{defn}

Next we project to $X_{syn}^\N$ and modify the resulting $\O^{crys}_.$-modules (resp. complexes) exactly as in Def. \ref{deffilmod}. Here we use the
functor $\jj^*$ for $(E_.,A_.)=(\smash{X^\N_{syn}},\smash{\O^{crys}_.})$.  Concretely the functor $\jj^*$ is
thus simply $M_.\mapsto \smash{\O^{crys}_.\otimes_{\O^{crys}_{.+1}}M_{.+1}}$.
A pleasant feature of the present setting is that the ring $\O^{crys}_.$ need not to be modified, thanks to Lem. \ref{lemsyn} \ref{lemsyniv}.

\begin{defn} \label{deffilsynmod} Consider the global (resp. local) situation and let $D\in \DC(X)$.
\debrom \item \label{deffilsynmodi} We define a filtered module (resp. complex)  as follows
$$\begin{array}{crcllc}&Fil^{i,crys}_.D&:=&u_{*}l^{-1}Fil^iD^{syn}&\hbox{in   $Mod(X_{syn}^\N,\O^{crys}_.)$, $i=0,1$}&\\ \hbox{(resp. }&
Fil^{i}L^{crys}\Omega^\bullet_{.,T}(D)&:=&u_{*}l^{-1}Fil^iL^{syn}\Omega^{\bullet}_{T_{.}}(D) &\hbox{in  $Kom(X_{syn}^\N,\O^{crys}_.)$, $i=0,1$}&\hbox{
)}\end{array}$$

\item \label{deffilsynmodii} We define a couple of modules (resp. complexes) as follows
$$\begin{array}{crcllc}& \widetilde{Fil}{}^{i,crys}_.D&:=&\jj^*Fil^{i,crys}_{.}D&\hbox{in   $Mod(X_{syn}^\N,\O^{crys}_.)$, $i=0,1$}&\\
\hbox{(resp. }&\widetilde{Fil}{}^{i}L^{crys}\Omega^\bullet_{.,T}(D)&:=&\jj^*Fil^{i}
L^{crys}\Omega^\bullet_{.,T}(D)&\hbox{in $Kom(X_{syn}^\N,\O^{crys}_.)$, $i=0,1$}&\hbox{ )}\end{array}$$
\finrom
\end{defn}
Note that in virtue of Lem. \ref{lemsyn} \ref{lemsyniv}, $Fil^{i,crys}$ and $\smash{\widetilde{Fil}{}^{i,crys}}$ coincide for $i=0$.

\begin{prop} \label{tdfilonesyn} Consider the global (resp. local) situation and consider $D\in \DC(X)$. Let $\smash{Fun^._.}$ (resp. $\smash{{\widetilde{Fun}}{}^._.}$) denote one of the two couples of functors defined in Def. \ref{deffilsynmod} \ref{deffilsynmodi} (resp. \ref{deffilsynmodii}).
\debrom

\item \label{tdfilonesyni} There is a natural morphism of functorial distinguished triangles in \smash{$D(X^\N_{syn},\O_.^{crys})$}:
    $$\xymatrix{\widetilde{Fun}{}^1_.(D)\ar[r]^-1\ar[d]&\widetilde{Fun}{}^0_.(D)\ar[r]^-{can}\ar[d]& \widetilde{Lie}^{syn}{}_.(D)\ar[r]^-{+1}\ar[d]&\\
Fun^1_.(D)\ar[r]^-1&Fun^0_.(D)\ar[r]^-{can}&Lie^{syn}_.(D)\ar[r]^-{+1}&}$$
where $\smash{Lie^{syn}_.(D)}$ is obtained by scalar restriction to $\smash{\O_.^{crys}}$ from $Lie^{syn}(D)$ viewed as a constant projective system of $\O$-modules while $\smash{\widetilde{Lie}^{syn}_.(D)}:=\smash{\tau_{\ge_{-1}} L\jj^*Lie^{syn}_.(D)}$.

\item \label{tdfilonesynii} The objects of the complexes $\smash{\widetilde{Fun}{}^i_.}$ are $L$-normalized $\smash{\O^{crys}_.}$-modules. The objects of the complexes $\smash{\widetilde{Fil}{}^iL^{crys}\Omega^\bullet_{.,T}(D)}$ are moreover $\epsilon_*$-acyclic.

\item \label{tdfilonesyniii} In the case $Fun^._.=\smash{Fil^.L^{crys}\Omega^\bullet_{.,T}}$, the image of the diagram in \ref{tdfilonesyni} by $R\epsilon_*$ is naturally isomorphic to the image by $\iota^{-1}$ of the diagram in Lem. \ref{lemfil} \ref{lemfili} in the case $Fun^._.=\smash{Fil^.\Omega^\bullet_{.,T}}$. More precisely, there is a canonical commutative cube as follows in $\smash{Kom(X_{et}^\N,\widetilde \O^{crys}_.)}$: $$\xymatrix{\iota^{-1}\widetilde{Fil}{}^1\Omega^\bullet_{.,T}(D)\ar[dd]\ar[rr]^(.68)\sim\ar[rd]^-1&&\epsilon_*\widetilde{Fil}{}^1L^{crys}\Omega^\bullet_{.,T}(D)\ar[rd]^-1\ar'[d][dd]\\
    &\iota^{-1}\widetilde{Fil}{}^0\Omega^\bullet_{.,T}(D)\ar[dd]\ar[rr]^(.6)\sim&&\epsilon_*\widetilde{Fil}{}^0L^{crys}\Omega^\bullet_{.,T}(D)\ar[dd]\\
    \iota^{-1}{Fil}{}^1\Omega^\bullet_{.,T}(D)\ar[rd]^-1\ar'[r][rr]^(.35)\sim&&\epsilon_*{Fil}{}^1L^{crys}\Omega^\bullet_{.,T}(D)\ar[rd]^-1\\
    &\iota^{-1}{Fil}{}^0\Omega^\bullet_{.,T}(D)\ar[rr]^(.6)\sim&&\epsilon_*{Fil}{}^0L^{crys}\Omega^\bullet_{.,T}(D)}$$
\finrom
\end{prop}
Proof. \ref{tdfilonesyni} Using Prop. \ref{complin} \ref{complini}, Def. \ref{deffilsyn} \ref{deffilsyniii}, Lem. \ref{lemsyn} \ref{lemsyni} and the isomorphism
$l^{-1}i_*\simeq i_*l^{-1}$ we find the following exact sequences on $((X/\Sigma_.)_{crys,syn},\O)$
$$\begin{array}{c}\xymatrix{0\ar[r]&l^{-1}Fil^1D^{syn}\ar[r]&l^{-1}Fil^0D^{syn}\ar[r]&i_*l^{-1}Lie^{syn}(D)\ar[r]&0} \\
\xymatrix{0\ar[r]&l^{-1}Fil^1L^{syn}\Omega_{T_.}^\bullet(D)\ar[r]&l^{-1}Fil^0L^{syn}\Omega_{T_.}^\bullet(D)
\ar[r]&i_*l^{-1}Lie^{syn}(D)\ar[r]&0}\end{array}$$
These sequences stay exact after applying $u_*$ thanks to the second statement of Lem. \ref{lemsyn} \ref{lemsyniii}. Whence the bottom distinguished triangle. The top one follows by Lem. \ref{cartgen} \ref{cartgeniv}.

\ref{lemsynii} The first statement is part of Lem. \ref{cartgen} \ref{cartgeniv} as well. Let us prove the acyclicity statements. First, we note that $\smash{Fil^0L^{crys}\Omega^q_{.,T}}(D)$ is $\epsilon_*$-acyclic by the conjunction of Lem. \ref{lemsyn} \ref{lemsynii} and \ref{lemsyniii}. The same is true for $Lie^{syn}_.(D)$ and $\smash{L^1\jj^*Lie^{syn}_.(D)}$ by Lem. \ref{acycqcohsch} \ref{acycqcohschi}. Given the isomorphisms $$\widetilde Fil^1L^{crys}\Omega^q_{.,T}\simeq Fil^0L^{crys}\Omega^q_{.,T}$$ for $q\ge 1$ and the exact sequence $$\xymatrix{0\ar[r]&L^1\jj^* Lie^{syn}_.(D)\ar[r]&\widetilde {Fil}{}^1L^{crys}\Omega^0_{.,T}(D)\ar[r]&Fil^0L^{crys}\Omega^q_{.,T}(D)\ar[r]^-{can}&Lie^{syn}_.(D)\ar[r]&0}$$ it is thus sufficient to check that the image of $can$ by $\epsilon_*$ remains epimorphic. This is indeed the case by Prop. \ref{liesurj} \ref{liesurji}, thanks to the following compatible isomorphisms \begin{eqnarray}\label{carrelemsyn}\xymatrix{\epsilon_*Fil^0L^{crys}\Omega^0_{.,T}(D)\ar[r]\ar[d]^\wr&\epsilon_*Lie^{syn}_.(D)\ar[d]^\wr\\
\iota^{-1}Fil^0\Omega^0_{.,T}(D)\ar[r]&\iota^{-1}Lie_{.,T}(D)}\end{eqnarray} resulting from Lem. \ref{lemreal2} \ref{lemreal2ii} together with the full faithfulness of $\epsilon^*$ on quasi-coherent modules and quasi-coherent crystals.

\ref{lemsyniii} The arguments giving (\ref{carrelemsyn}) provide the bottom commutative square and the rest follows formally.
\begin{flushright}$\square$\end{flushright}

\para \label{phisyn}
From here, the construction is similar to the case of the \'etale topology.

\begin{defn} \label{defFrsyn} Let $(D,f,v)\in \cal DC(X)$. We use the following notations. \petit

\noindent - In the global situation, we let $Fr:\smash{Fil^{0,crys}_.D}\to (\smash{F})_*Fil_.^{0,crys}D$ in $\smash{Mod(X^\N_{syn},\O^{crys}_.)}$
denote
$$\xymatrix{u_*l^*\epsilon^*D\ar[rr]^-{F^{(-/X^\N)}}&&
(F)_*F^{*}u_*l^*\epsilon^*D\ar[r]& (F)_*u_*l^*\epsilon^*F^*D\ar[r]^-f& (F)_*u_*l^*\epsilon^*D}$$

\noindent - In the local situation, we let  $Fr:\smash{Fil^{0}L^{crys}\Omega^\bullet_{.,T}(D)}\to (F)_*\smash{Fil^{0}L^{crys}\Omega^\bullet_{.,T}(D)}$
in \break $Kom(\smash{X_{syn}^\N,\O^{crys}_.})$ denote
$$\xymatrix{u_*l^*\epsilon^*L(\Omega^\bullet_{T_.}(D))\ar[rr]^-{{F^{(-/X^\N)}}}&&
(F)_*F^{*}u_*l^*\epsilon^*L(\Omega^\bullet_{T_{.}}(D))\ar[r]&
(\smash{F})_*u_*l^*\epsilon^*F^*L(\Omega^\bullet_{T_{.}}(D))
\ar[d]^-{(\ref{functLdR})}\\
&&(\smash{F})_*u_*l^*\epsilon^*L(\Omega^\bullet_{T_{.}}(D))&\ar[l]_-f(\smash{F})_*
u_*l^*\epsilon^*L(\Omega^\bullet_{T_{.}}(F^*(D)))}
$$
\end{defn}

\begin{prop} \label{defphisyn} Consider the global (resp. local situation) and let $D\in \DC(X)$. There exist unique morphisms $\varphi$ in  $\smash{Kom(X^\N_{syn},\O^{crys}_.)}$ rendering
the following squares commutative: $$\begin{array}{cc}\xymatrix{{\widetilde{Fil}{}^{1,crys}_.D}\ar[r]^-\varphi\ar[d]^1&(F)_*
{\widetilde{Fil}{}^{0,crys}_.D}\ar[d]^p\\
{\widetilde{Fil}{}^{0,crys}_.D}\ar[r]^-{Fr}&(F)_*{\widetilde{Fil}{}^{0,crys}_.D}}
&\xymatrix{Fil^{1}L^{crys}\Omega^\bullet_{.,T}(D)\ar[r]^-\varphi\ar[d]^1&(F)_*
Fil^{0}L^{crys}\Omega^\bullet_{.,T}(D)\ar[d]^p\\
Fil^{0}L^{crys}\Omega^\bullet_{.,T}(D)\ar[r]^-{Fr}& (F)_*Fil^{0}L^{crys}\Omega^\bullet_{.,T}(D)}\end{array}
$$
These morphisms are functorial with respect to $D$ and $X$ or $Y$ accordingly.
\end{prop}
Proof. Using Lem. \ref{cartgen} \ref{cartgenii}, \ref{cartgeniii} as in the proof of Prop. \ref{defphi}, we are reduced to prove that the first component
\begin{eqnarray}
\label{Fr0Dsyn}Fr_1:Fil^{0,crys}_1D\to (F)_*Fil^{0,{crys}}_1D \hbox{ and}\\
\label{Fr0Lsyn} Fr_1:\smash{Fil^{0}L^{crys}\Omega^q_{1,T}}(D)\to (F)_* \smash{Fil^{0}L^{crys}\Omega^q_{1,T}}(D)\end{eqnarray}
of the morphisms $Fr$ vanish respectively on $\smash{Fil^{1,crys}_1D}$ and $\smash{Fil^{1}L^{crys}\Omega^q_{1,T}(D)}$. For $q\ge 1$, the arrow
(\ref{Fr0Lsyn}) is zero since the base change morphism (\ref{functLdR}) occurring in its definition (see Def. \ref{defFrsyn}) is zero. Next, let
$M$ denote either the module $D$ or the complex $\smash{L(\Omega^q_{T_{.}}(D))}$ on $\smash{((X/\Sigma_\infty)_{crys,et},\O)}$ and $f:F^*M\to M$ the morphism given by the
Frobenius of the Dieudonn\'e crystal (using  (\ref{functLdR}) in the second case). The reader may check that  $Fr_1:u_*\iota_1^{-1}\epsilon^*M\to (F)_*u_*\iota_1^{-1}\epsilon^*M$ identifies with the composed morphism

\begin{eqnarray}\label{FR0syn}\xymatrix{u_*\epsilon^*\iota_1^{-1}M\ar[rr]^-{u_*\epsilon^*\overline Fr}&&u_*\epsilon^*(F)_*M\ar[r]&(F)_*u_*\epsilon^*\iota_1^{-1}M}\end{eqnarray}

\noindent where $\iota_1:\smash{((X/\Sigma_1)_{crys,et},\O)}\to \smash{((X/\Sigma_\infty)_{crys,et},\O)}$ and $\overline Fr$ is the morphism introduced during the proof of Prop. \ref{defphi}. Consider now the commutative diagram
\begin{eqnarray}\label{diagFrsyn}\xymatrix{&\epsilon^*\iota_1^{-1}Fil^1 D\ar[r]\ar@{->>}[d]&\epsilon^*\iota_1^{-1}D\ar[rr]^-{\epsilon^*\iota_1^{-1}can_D}\ar[d]^\wr&&
\epsilon^*\iota_1^{-1}i_*Lie(D)\ar[r]\ar[d]_{\ref{lemsyn} \ref{lemsyni}}^\wr&0\\
0\ar[r]&\iota_1^{-1}Fil^1D^{syn}\ar[r]&\iota_1^{-1}\epsilon^*D\ar[rr]^{\iota_1^{-1}can_D^{syn}}&&
\iota_1^{-1}i_*\epsilon^*Lie(D)&}\end{eqnarray}
where the first line is exact and is obtained from Prop. \ref{complin} \ref{complini} by applying the right exact functor $\epsilon^*\iota_1^{-1}$ while the second one is exact as well and is deduced from  Def. \ref{deffilsyn} \ref{deffilsyniii} by the exact functor $\iota_1^{-1}$. It has been shown in Prop. \ref{defphi} that $\overline Fr$ vanishes on the kernel of
$\iota_1^{-1}can_D$. It follows in particular that $\epsilon^* \overline Fr$ vanishes on the   image of the left vertical arrow in (\ref{diagFrsyn}).
Applying $u_*$, we find that (\ref{Fr0Dsyn}) vanishes on $\smash{Fil^{1,crys}_1D}:=\smash{u_*\iota_1^{-1}Fil^1D^{syn}}$ as desired. The case  of
(\ref{Fr0Lsyn}) is similar.
\begin{flushright}$\square$\end{flushright}

\begin{defn}\label{defSonephisyn}
\debrom \item \label{defSonephisyni} In the global situation, we define a functor $\smash{\cal S^{1,\varphi}_{syn,.,X}}:\DC(X)\to
Mod^{1,\varphi}(X^\N_{syn},\O^{crys}_.)$ by setting $$\begin{array}{rcl}\cal
S^{1,\varphi}_{syn,.,X}(D)&:=&(\widetilde{Fil}{}^{1,crys}_{.}(D),\widetilde{Fil}{}^{0,crys}_{.}(D),1,\varphi)\end{array}$$
By functoriality with respect to $D$ and $X$, this defines a colax morphism $$\cal S^{1,\varphi}_{syn,.}:\cal DC(-)\to
Mod^{1,\varphi}((-)^\N_{syn},\O^{crys}_.)$$ between contravariant pseudo-functors on $\smash{\cal Sch^{sflpb}/\Sigma_1}$ or $\smash{Diag(\cal Sch^{sflpb}/\Sigma_1)}$ (Lem. \ref{lemdiagcodiag} \ref{lemdiagcodiagii}, Rem. \ref{remdiagcodiag} \ref{remdiagcodiagii}).

\item \label{defSonephisynii}In the local situation, we define a functor  $\smash{\cal SL\Omega^{\bullet,1,\varphi}_{syn,.,T}}:\DC(X)\to Kom^{1,\varphi}(X^\N_{syn},\O^{crys}_.)$ by setting
$$\begin{array}{rcl}\cal SL\Omega^{\bullet,1,\varphi}_{syn,.,T}(D)&:=&(\widetilde{Fil}{}^{1}L^{crys}\Omega^\bullet_{.,T}(D),
\widetilde{Fil}{}^{0}L^{crys}\Omega^\bullet_{.,T}(D),1,\varphi)\end{array}$$

By functoriality with respect to $D$ and $(X,Y,\iota,\tilde F)$, this defines a colax morphism $$\cal SL\Omega^{\bullet,1,\varphi}_{syn,.}:\DC(-)\to Kom^{1,\varphi}((-)^\N_{syn},\O^{crys}_.)$$
between contravariant pseudo-functors on $\smash{Emb^{glob,lfpb}_F}$ or $Diag(\smash{Emb^{glob,lfpb}_F})$.
\finrom
\end{defn}

Here, as in the case of the \'etale topology, one could globalize the functor $\smash{\cal SL\Omega^{\bullet,1,\varphi}_{syn,.}}$ using Sect. \ref{paraprelimS}. This won't be necessary for our purpose.

\begin{prop} \label{compsynetloc} Let $X$ denote a diagram of $\cal Sch^{slfpb}/\Sigma_1$ and assume given $Y$ in $Emb^{glob}_F(X)$. For $D\in DC(X)$, there are canonical isomorphisms
$$\begin{array}{rcll}  \cal S^{1,\varphi}_{syn,.,X}(D)&\simeq &\cal SL\Omega^{\bullet,1,\varphi}_{syn,.,T}(D)&\hbox{in $D(Mod^{1,\varphi}(X^\N_{syn},\O^{crys}_.))$}\\
\iota^{-1}\cal S\Omega^{\bullet,1,\varphi}_{et,.,T}(D)&\simeq &R\epsilon_*\cal SL\Omega^{\bullet,1,\varphi}_{syn,.,T}(D)&\hbox{in $D(Mod^{1,\varphi}(X^\N_{et},\widetilde\O^{crys}_.)).$}
 \end{array}$$
These isomorphisms are functorial with respect to $D$ and $Y$ in the obvious sense. \end{prop}

\noindent Proof. Recall that there is a natural
quasi-isomorphism \begin{eqnarray}\label{qisaz}D\to L\Omega^\bullet_{T_.}(D)\end{eqnarray} in $Kom((X/\Sigma_\infty)_{crys,et},\O)$ (Prop. \ref{dRlocal} \ref{dRlocali}). Applying $\epsilon^*$ and
taking Def. \ref{deffilsyn} into account, we get compatible morphisms $\smash{Fil^iD^{syn}}\to \smash{Fil^iL^{syn}\Omega_{T_.}^\bullet(D)}$, $i=0,1$. Apply
then $l^*$, $u_*$ and $\jj^*$ to get a couple of morphisms \begin{eqnarray}\label{moraz}\widetilde{Fil}{}_.^{i,crys}D\to
\widetilde{Fil}{}^iL^{crys}\Omega^\bullet_{.,T}\hbox{, $i=0,1$}\end{eqnarray} which are compatible with the morphism $1$. These morphisms are clearly
compatible with $Fr$ and thus with $\varphi$ also, thanks to Lem. \ref{cartgen} \ref{cartgeniii}. Up to now, we have thus obtained a canonical morphism $$\cal
S^{1,\varphi}_{syn,.,X}(D)\to \cal SL\Omega^{\bullet,1,\varphi}_{syn,.,T}(D)$$ in $Kom^{1,\varphi}(X^\N_{syn},\O^{crys}_.)$. It thus remains to check that this
is a quasi-isomorphism. By Lem. \ref{1phimod} \ref{1phimodiii} and
Prop. \ref{tdfilonesyn} \ref{tdfilonesyni}, we see that it suffices to check that
(\ref{moraz}) is a quasi-isomorphism for $i=0$, ie. that the image of  (\ref{qisaz}) under $u_*l^{-1}\epsilon^*$ remains a quasi-isomorphism. This is indeed the case by  Lem. \ref{lemsyn} \ref{lemsyniii} (note also that local freeness implies the acyclicity of $\epsilon^*$). Similarly, the second isomorphism easily reduces to Lem. \ref{lemreal2} \ref{lemreal2ii} using  Lem. \ref{1phimod} \ref{1phimodiii}, Prop. \ref{tdfilonesyn} \ref{tdfilonesyniii} and Lem. \ref{lemsyn} \ref{lemsyni}, \ref{lemsynii}.
\begin{flushright}$\square$\end{flushright}

%
%
%

Let us turn to the global situation. For a semi-simplicial $\smash{T_{[.]}}$ as in Sect. \ref{prelimemb}, we have the following partial counterpart for (\ref{descentdiaget}):
$$ \xymatrix{((X/\Sigma_.)_{crys,syn},\O)\ar[d]_u& ((U_{[.]}/\Sigma_.)_{crys,syn},\O)\ar[l]_-{f_{U_{[.],.}}} \ar[d]_u
&((X/\Sigma_.)_{crys,syn},\O)/T_{[.],.}\ar[l]_-{f_{T_{[.],.}/U_{[.],.}}} \ar @{->} @<+0pt> `u[l] `[ll]_-{f_{T_{[.],.}}} [ll]\\
(X^{\N}_{syn},\O^{crys}_.)&\ar[l]_-{f_{U_{[.],.}}}
(U^{\N}_{[.],syn},\O^{crys}_.)
}$$
This diagram and (\ref{descentdiaget}) are compatible with each other in the obvious way via the morphisms
$\epsilon:\smash{((X^\s/\Sigma_.)_{crys,syn},\O)}\to \smash{((X^\s/\Sigma_.)_{crys,et},\O)}$, $\smash{\epsilon:(X^{\N}_{syn},\O^{crys}_.)}\to
\smash{(X^{\N}_{et},\O^{crys}_.)}$  and their localizations.

\begin{prop} \label{compsynet} Consider $X$ in $\cal B_0$ and $\smash{Y_{[.]}}$ in $HR^{et}_F(X)$ with divided power envelope $\smash{T_{[.]}}$.
\debrom \item \label{compsyneti} The adjunction morphisms induce the following isomorphisms: $$\begin{array}{rcll}\cal S^{1,\varphi}_{syn,.,X}(D)&\simeq&Rf_{U_{[.],.},*}\cal S^{1,\varphi}_{syn,.,U_{[.]}}(D)&\hbox{in $D(Mod^{1,\varphi}(X^\N_{syn},\O^{crys}_.))$}\\
\cal S^{1,\varphi}_{et,.,X}(D)&\simeq&Rf_{U_{[.],.},*}\cal S^{1,\varphi}_{et,.,U_{[.]}}(D)&\hbox{in $D(Mod^{1,\varphi}(X^\N_{et},\widetilde\O{}^{crys}_.))$}\\
\end{array}$$

\item \label{compsynetii} The isomorphisms of \ref{compsyneti} and Prop. \ref{compsynetloc} induce $$\begin{array}{rcll}R\epsilon_*\cal S^{1,\varphi}_{syn,.,X}(D)&\simeq&\cal
S^{1,\varphi}_{et,.,X}(D)&\hbox{in $D(Mod^{1,\varphi}(X^\N_{et},\widetilde\O{}^{crys}_.))$}\end{array}$$ This isomorphism is functorial with respect to $D$ and
independent of the choice of $\smash{Y_{[.]}}$. \finrom \end{prop}
Proof. \ref{compsyneti} The first (resp. second) morphism is invertible by cohomological descent in $X_{syn}$ (resp. $X_{et}$) and d\'evissage using Lem. \ref{1phimod} \ref{1phimodiii} and Prop. \ref{tdfilonesyn} \ref{tdfilonesyni} (resp. Lem. \ref{1phimod} \ref{1phimodiii} and Lem. \ref{lemfil} \ref{lemfili}).

\ref{compsynetii} The claimed isomorphism follows from \ref{compsyneti} by Prop. \ref{compsynetloc} applied to $U_{[.]}$ viewed as a diagram of $\cal Sch^{slfpb}/\Sigma$ and $Y_{[.]}$ as an object of $Emb^{glob}_F$. The independence with respect to the choice of $Y_{[.]}$ follows from the connectedness of the category
$\smash{HR^{et}_F(X)}$.
\begin{flushright}$\square$\end{flushright}

\begin{defn} \label{defSsyn} If $X$ is a diagram of separated schemes with local finite $p$-bases and $D$ is in $\cal DC(X)$, we define the \emph{syntomic complex of $D$ on the syntomic site} as follows: $$\cal S_{syn,.,X}(D):=R\varpi_*\cal S^{1,\varphi}_{syn,.,X}(D)$$
\end{defn}

\begin{rem}  \label{tdSsyn} In $D(X^\N_{syn},\O_.^{crys,F=1})$, one has canonical distinguished triangles $$\xymatrix{\widetilde{Fil}^{1,crys}_.(D)\ar[r]^-{1}&\widetilde{Fil}^{0,crys}_.(D)\ar[r]&\widetilde{Lie}_.^{syn}(D)\ar[r]^-{+1}&}$$
$$\xymatrix{\cal S_{syn,.,X}(D)\ar[r]&\widetilde{Fil}^{1,crys}_.(D)\ar[r]^-{1-\varphi}&\widetilde{Fil}^{0,crys}_.(D)\ar[r]^-{+1}&}$$
where $\widetilde{Fil}^{0,crys}_.(D)=u_*l^{-1}D^{syn}$ and
 $\widetilde{Fil}^{1,crys}_.(D)=\jj^*u_*l^{-1}Fil^1D^{syn}$ are concentrated in degree $0$. Let us also point out that $D^{syn}$ and $Fil^1D^{syn}$ are $u_*$-acyclic as already noticed.
\end{rem}


\section{D\'evissage of flat cohomology} \label{secdofc}

The purpose of this section is to prove that the flat cohomology (vanishing at $Z$) of the N\'eron model of an Abelian variety (with good reduction outside $Z$) can be recovered from a suitable diagram
of $p$-divisible groups. This will be achieved in Prop. \ref{devissageFL}.

\subsection{Functorial mapping fibers in derived categories} \label{fmfidc} ~~ \\

We explain a basic construction, designed to avoid the difficulties caused by the lack of functoriality of mapping cones and mapping fibers with
respect to morphisms in the derived category. We use the notations $Kom(\cal A)$ and $D(\cal A)$ to designate respectively the category of complexes
and the derived category of a given Abelian category $\cal A$.

\para \label{prelimMF1} Let $[1]$ denote the category $\{0\to 1\}$. If $\cal C$ is a category,  $\cal C^{[1]}$ is thus the category of arrows of $\cal C$.
Given an Abelian category $\cal A$, we have an obvious functor \begin{eqnarray}\label{MF11}\xymatrix{D(\cal A^{[1]})\ar[r]& D(\cal
A)^{[1]}}\end{eqnarray} This functor is not an equivalence in general. The objects of the source category will be called \emph{true arrows}. Consider
the \emph{mapping fiber} functor \begin{eqnarray}\label{MF12}\xymatrix{Kom(\cal A)^{[1]}\ar[r]&Kom(\cal A)}\end{eqnarray} taking a morphism of
complexes $f:A\to B$ to the simple complex associated to $f$ viewed as a double complex placed in degrees $[0,1]\times ]-\infty,+\infty[$. The
functor (\ref{MF12}) does not induce a functor $D(\cal A)^{[1]}\to D(\cal A)$. We can nevertheless make the following definition.

\begin{defn} \label{functoMF} The \emph{functorial mapping fiber}
$$MF:D(\cal A^{[1]})\to D(\cal A)$$ is the functor deduced from (\ref{MF12}) by the universal property of derived categories.
\end{defn}
By definition we thus have a tautological distinguished triangle \begin{eqnarray}\label{MF} \xymatrix{MF(f)\ar[r]& A\ar[r]&
B\ar[r]^{+1}&}\end{eqnarray} for any true arrow $f:A\to B$.  Assume now that $\cal A$ has enough injectives. Consider left exact functors $F_i:\cal
A\to \cal A'$, $i=1,2$ and a natural transformation $t:F_1\rightarrow F_2$. The triple $F=(F_1,F_2,t)$ can be viewed as a left exact functor $\cal
A\to \smash{\cal A'^{[1]}}$. Deriving, we get $RF:D^+(\cal A)\rightarrow D^+(\smash{\cal A'^{[1]}})$. Now (\ref{MF}) yields a distinguished triangle
\begin{eqnarray}\label{MFt}\xymatrix{MF(RF(M))\ar[r]& RF_1(M)\ar[r]& RF_2(M)\ar[r]^-{+1}&}\end{eqnarray} which is functorial with respect to $M$ in
$D^+(\cal A)$.

\begin{lem} \label{lemMF1} Let $\cal A$ and $F$  be as above and consider moreover left
exact functors $F_3:\cal A'\to \cal A''$ and $F_4:\cal A'''\to \cal A$, where $\cal A'$ and $\cal A'''$ have enough injectives.

\debrom
\item \label{lemMF1i} The functor
$\smash{F_3^{[1]}}:\smash{\cal A'^{[1]}}\to \smash{\cal A''^{[1]}}$ induced by $F_3$ is right derivable. The derived functor $\smash{RF_3^{[1]}}$ can
be computed componentwise. There is a canonical isomorphism $$RF_3(MF(f'))\simeq MF(RF_3^{[1]}(f'))$$
    in $\smash{D^+(\cal A'')}$ which is functorial with respect to $f'$ in $\smash{D^+(\cal A'^{[1]})}$.

\item \label{lemMF1ii} Assume that $F_1$ and $F_2$ send injectives of $\cal A$ to $F_3$-acyclic objects.
There is a canonical isomorphism of functors $D^+(\cal A)\to D^+(\cal A''^{[1]})$: $$RF_3^{[1]}\circ RF\simeq R(F_3^{[1]}\circ F)$$

\item \label{lemMF1iii} Assume that $F_4$ sends injectives of $\cal A'''$ to $F_i$-acyclic objects $i=1,2$.
There is a canonical isomorphism of functors $D^+(\cal A''')\to D^+(\cal A^{[1]})$: $$RF\circ RF_4\simeq R(F\circ F_4)$$

\finrom
\end{lem}
Proof. \ref{lemMF1i} The first statement follows from the fact that there are enough injectives in $\smash{\cal A'^{[1]}}$. The second and third
statements follow from the fact that if $f':A'\to B'$ is an injective object of the category $\smash{\cal A'^{[1]}}$, then both $A'$ and $B'$ are
injective as well.

\ref{lemMF1ii} This follows from the fact that $RF_3^{[1]}$ can be computed using true arrows whose source and target are acyclic for
$RF_3$.

\ref{lemMF1iii} This follows from the fact that $F_i$-acyclicity for $i=1,2$ is equivalent to $F$-acyclicity.
\begin{flushright}$\square$\end{flushright}


\para \label{prelimMF2}
Let us discuss an alternative point of view on the distinguished triangle $(\ref{MFt})$.

\begin{lem} \label{lemMF2} Let $\cal A$ and $F=(F_1,F_2,t)$ be as in Lem. \ref{lemMF1} and set
$F_0(M):=Ker(F_1(M)\rightarrow F_2(M))$. The functor $F_0$ is right derivable and there is a natural transformation $$RF_0\to MF \circ RF$$ This
natural transformation is an isomorphism if the following holds: \petit

\emph{(Condition)} If $I$ is an injective object of $\cal A$ then $F_1(I)\to F_2(I)$ is epimorphic.
\end{lem}
Proof. The first statement is clear, since $F_0$ is left exact and $\cal A$ has enough injectives. The claimed arrow is induced by the obvious isomorphisms $RF_0\simeq MF\circ R(F_0\to 0)$  and the natural transformation $R(F_0\to 0)\to R(F_1\to F_2)$. If the condition holds then $0\rightarrow F_0(I)\rightarrow F_1(I)\rightarrow
F_2(I)\rightarrow 0$ is exact for $I$ injective and we may conclude.
\begin{flushright}$\square$\end{flushright}


\subsection{Some useful diagrams and functors} \label{sudaf}

\para The following notations will be used from now on.
\begin{eqnarray}\label{diagcomplet}
\diagram{Z_v&\hfl{z_v}{}&C_v&\hflrev{o_v}{}&C_v^\s&\hflrev{j_v}{}&U_v\cr  \vfl{{\iota_{Z_v}}}{}&&\vfl{{\iota_{C_v}}}{}&&\vfl{{\iota_{C^\s_v}}}{}&&\vfl{{\iota_{U_v}}}{}\cr
Z&\hfl{z}{}&C&\hflrev{o}{}&C^\s&\hflrev{j}{}&U\cr}
\end{eqnarray}
- $C$ is an irreducible  smooth curve over $\Sigma_1=Spec(\Fp)$. \\
- $Z$ is a smooth effective divisor on $C$. In other terms, $Z$ is a disjoint union of $Spec(k_v)$'s, with $v$ running through a finite number of closed points of $C$. \\
- $C^\s$ is the log scheme $(C,Z)$ and $U$ is the complementary open subscheme of $Z$ in $C$. \\
- For each point $v$ of $Z$, we let $Z_{v}=Spec(k_v)$ denote the corresponding reduced closed subscheme, as well as $C_v=Spec(O_v)$, $C_v^\s=(C_v,Z_v)$ and $U_v=Spec(K_v)$, where $O_v$ is the completion of the local ring $O_{C,v}$ and $K_v$ is the fraction field of $O_v$.\\
- The arrows $z,j,o,z_v,j_v,o_v,{\iota_{Z_v}},{\iota_{C_v}},{\iota_{C_v^\s}},{\iota_{U_v}}$ are the obvious ones. \petit

\para Let us introduce notations for some variants of the diagram (\ref{diagcomplet}).
We consider star shaped diagrams $J$ and $J^\s$ in which $U$ is the central vertex while the branches are indexed by $v$ running through the
points of $Z$ and  have the following form (we only draw one branch for convenience): \begin{eqnarray}\label{defdiag}\begin{array}{rcl} 
J&:=&(\xymatrix{C_v&U_v\ar[l]_{j_v}\ar[r]^{\iota_{U_v}}&U})\\
J^\s&:=&(\xymatrix{C_v^\s&U_v\ar[l]_{j_v}\ar[r]^{\iota_{U_v}}&U})\end{array}\end{eqnarray}

Let us consider moreover the discrete diagram \begin{eqnarray}\label{defdiagZ}\begin{array}{rcl} Z_J&:=&(Z_v)\end{array}\end{eqnarray} whose vertices
are indexed by the points of $Z$ (this discrete diagram of schemes should not be confused with the scheme $Z$ itself). With these notations, we have a
natural commutative diagram as follows in the category $Diag(\cal Sch^\s/\Sigma_1)$.
\begin{eqnarray}\label{diagdiag}\xymatrix{Z_J\ar[r]^-{z_J}
\ar[d]_{m_Z}&J\ar[d]^m&\ar[l]_-{o_J}J^\s\ar[d]^{m^\s}\\ Z
\ar[r]^-z&C&\ar[l]_-oC^\s}
\end{eqnarray}


\para \label{diagopi}
Let us discuss Mayer-Vietoris functors and triangles which are relevant to completion at points of $Z$.
Recall that an object of the topos $J_{et}$ is a diagram of the form $$\xymatrix{F_{C_v}\ar[r]^{F_{j_v}}&F_{U_v}&F_U\ar[l]_{F_{\iota_{U_v}}}}$$ in
$(-)_{et,cof}$ above $\smash{\Delta_J^{op}}$. In other terms, $F_U$ is an object of $U_{et}$ and for each $v$ in $Z$, $F_{C_v}$ (resp. $F_{U_v}$) is an object of $C_{v,et}$ (resp. $U_{v,et}$) while
$F_{j_v}:j_v^{-1}F_{C_v}\to F_{U_v}$ and $F_{\iota_{U_v}}:\iota_{U_v}^{-1}F_{U}\to F_{U_v}$ are morphisms of $U_{v,et}$.

\begin{lem} \label{defMV} \emph{(complete Mayer-Vietoris functors)} Consider the natural morphism  $$m:(J_{et}^\N,\Z/p^.)\to (C_{et}^\N,\Z/p^.)$$

\debrom \item \label{defMVi} For $M=(M_{C_v}\to M_{U_v}\leftarrow M_U)$ in $\smash{Mod(J_{et}^\N,\Z/p^.)}$, we have canonically $$m_*M\simeq Ker\left( (\prod_{v\in
Z}\iota_{C_v,*}M_{C_v})\times j_*M_U\to \prod_{v\in Z}j_*\iota_{U_v,*}M_{U_v}\right)$$

\item \label{defMVii} For $M=(M_{C_v}\to M_{U_v}\leftarrow M_U)$ in $\smash{D^+(J_{et}^\N,\Z/p^.)}$, we have a canonical isomorphism  $$Rm_*M\simeq MF(R\left( (\prod_{v\in Z}\iota_{C_v,*}(-)_{C_v})\times j_*(-)_U
    \to \prod_{v\in Z}j_*\iota_{U_v,*}(-)_{U_v}\right)(M))$$
    where $MF$ is the functorial mapping fiber of Def. \ref{functoMF}. We get in particular, a canonical distinguished triangle $$\xymatrix{Rm_*M\ar[r]& (\prod_{v\in Z}R\iota_{C_v,*}M_{C_v})\times Rj_*M_U
    \ar[r]& \prod_{v\in Z}Rj_*R\iota_{U_v,*}M_{U_v}\ar[r]^-{+1}&}$$
\finrom
\end{lem}
Proof. \ref{defMVi} This is nothing but the projective limit formula of Lem. \ref{lemfibtop} \ref{lemfibtopii}.

\ref{defMVii} The claimed isomorphism will follow from Lem. \ref{lemMF2} once proven that $$\xymatrix{(\prod_{v\in Z}\iota_{C_v,*}M_{C_v})\times
j_*M_U\ar[rrrrrr]^-{(\prod_v\iota_{C_v,*}M(j_v)')\circ p_1-(j_*M(\iota_{U_v})')_v \circ p_2}&&&&&& \prod_{v\in Z}j_*\iota_{U_v,*}M_{U_v}}$$ is
epimorphic when $M$ is injective. Here, $M(j_v)':M_{C_v}\to j_{v,*}M_{U_v}$ denotes the morphism deduced from
$M(j_v):j_v^{-1}M_{C_v}\to M_{U_v}$ by adjunction and similarly for $M(\iota_{U_v})'$. Let $|J|$ denote the discrete diagram underlying $J$ and
$f:|J|\to J$ denote the obvious morphism. By injectivity of $M$, the natural monomorphism $M\to f_*f^{-1}M$  splits and it is thus sufficient to
prove that $(\prod_v\iota_{C_v,*}N(j_v)')\circ p_1-(j_*N(\iota_{U_v})')_v \circ p_2$ is epimorphic for $N=f_*f^{-1}M$. Now we see on the following explicit formula  $$f_*f^{-1}M=(M_{C_v}\times
j_{v,*}M_{U_v}\mathop{\longrightarrow}\limits^{p_2}M_{U_v}\mathop{\longleftarrow}\limits^{p_v\circ p_{1}}(\prod_{w\in Z}\iota_{U_w,*}M_{U_w})\times M_U)_v$$
that $N(j_v)'$ is a split epimorphism (it is the second projection $M_{C_v}\times
j_{v,*}M_{U_v}\to j_{v,*}M_{U_v}$) and we are done.
The distinguished triangle follows by (\ref{MFt}), since the component functors $\smash{(-)_{C_v}}$, $(-)_U$ and $\smash{(-)_{U_v}}$  send injectives to flasque modules  by Lem-Def. \ref{acycf} \ref{acycfv}.
\begin{flushright}$\square$\end{flushright}

\begin{rem} \debrom \item  Lem. \ref{defMV} has been stated for $\Z/p^.$-modules for convenience of later reference. It would clearly work for arbitrary rings.
\item Restricting Lem. \ref{defMV} \ref{defMVii} to the $k$-th component yields a similar isomorphism and a distinguished triangle in $D^+(C_{et},\Z/p^k)$.
    \finrom
\end{rem}

%
%
%
%
%

\para In practice the following lemma will be useful to compute $Rm_*$.

\begin{lem} \label{iotaexact} The following functors are exact: $$\begin{array}{rcccc}\smash{\iota_{U_v,*}}&:&Ab(U_{v,et})&\to& Ab(U_{et})\\ \smash{\iota_{C_v,*}}&:&Ab(C_{v,et})&\to& Ab(C_{et})\end{array}$$
\end{lem}
Proof. Let us explain the case of $\smash{\iota_{C_v,*}}$ since the case of $\smash{\iota_{U_v,*}}$ is easier (or alternatively, follows formally).  We will check that for any geometric point $\overline x$ of $C$, the functor $F\mapsto \smash{(\iota_{C_v,*}F)_{\overline x}}$ is exact. Obviously we can assume that $C$ is affine, say $C=Spec(A)$. With our usual notations, $K$ is thus the fraction field of $A$.

Consider a closed point $w$ of $Spec(A)$. Let $A_w$ denote the completion of $A$ at $w$ and $K_w$ its fraction field. Let us fix a separable closure $K^{sep}$ (resp. $K_w^{sep}$) of $K$ (resp. $K_w$) as well as an embedding $\iota_{w}^{sep}:K^{sep}\to K_w^{sep}$ for each $w$. We denote $G$ the Galois group of $K^{sep}/K$ and $D_{w^{sep}}$ (resp. $I_{w^{sep}}$) the decomposition (resp. inertia) group of $w^{sep}$ inside $G$. By Krasner's Lemma (\cite{La} II, 1¤7 Prop. 4), $D_{w^{sep}}$ thus identifies with $Gal(K_w^{sep}/K_w)$.

Let $\cal C$ denote the ordered set of subalgebras  $A'$ of $K^{sep}$ of finite type over $A$ such that the ring $A'$ is integrally closed.
For $A'$ in $\cal C$, we identify the set of  closed points of $Spec(A')$ with a subset of the set of places of the fraction field $K'$ of $A'$. We use the notation $w'$ to designate either the place of $K'$ induced by $\iota_w^{sep}$ or, when it exists, the corresponding point of $Spec(A')$. The completion of $A'$ (resp. $K'$) with respect to $w'$ is denoted $A'_{w'}$ (resp. $K'_{w'}$). We use respectively the notation $G'$ (resp. $D'_{w'}$, resp. $I_{w'}$) for the Galois group of $K'/K$ (resp. the image of $D_{w^{sep}}$, resp. of $I_{w^{sep}}$ inside it). If $H$ is a subgroup of $G'$ we have a pseudo-commutative diagram of topoi {\scriptsize$$\hspace{-0cm}\xymatrix{&&Spec(K_v)_{et}&\\
B_{K_v,D'_{v'}}/Spec(K'_{v'})\ar[rru]^-\star\ar[r]^-\star\ar[d]& B_{K_v,G'}/Spec(K_v\otimes_K K')\ar[ru]^-\star \ar[d]&\ar[l] B_{K_v,H}/Spec(K_v\otimes_K K')\ar[u]\ar[r]^-\star\ar[d]& Spec(K_v\otimes_K K'^H)_{et}\ar[lu]\ar[d]\\
B_{D'_{v'}}\ar[r]&B_{G'}&\ar[l]B_H\ar[r]&Set}$$}
which is moreover pseudo-functorial with respect to $Spec(K')$ (or $Spec(A')$) in the obvious way.
Here, we have used the notation $B_{H}$ (resp. $B_{R,H}$) for the topos of left $H$-sets (resp. left $H$-objects in $Spec(R)_{et}$) and the schemes $Spec(K_v\otimes_K K')$ are endowed with the left action of $G$ obtained by inverting the natural action of $G$ on $K'$. All the arrows of the diagram are simply obtained by functoriality of classifying topoi (\cite{SGA4-I} IV, Sect. 4.5) and localization. Let us only mention that in the bottom line, the direct images functors are respectively induction from $D'_{v'}$ to $G'$, from $H$ to $G'$ as described e.g. in \cite{Se2} and fixed points $(-)^H$. The squares of this diagram are subject to obvious base change isomorphisms expressing that the inverse image functors of the horizontal arrows pseudo-commute to the direct image functors of the vertical arrows.

If now, $K'/K$ is Galois, then $Spec(K_v\otimes_K K')$ is a $G'$-torsor in $Spec(K_v)_{et}$. In this situation the arrows indicated with a $\star$ are equivalences. Via these equivalences, the direct image functor of the first (resp. second) vertical arrow simply sends $F\in Spec(K_{v,et})$ to the set $\G(K'_{v'},F)$ (resp. $\G(K_v\otimes_K K',F)$) endowed with the left action of  $D'_{v'}$ (resp. $G'$) induced by the natural action of $D'_{v'}$ on $K'_{v'}$ (resp. $G'$ on $K'$). \petit


We are now in a position to prove the  exactness of the functor  $F\mapsto \smash{(\iota_{C_v,*}F)_{\overline x}}$. We consider the following cases.

\emph{Case 1}. The geometric point $\overline x$ lies above the generic point $\eta$ of $C$: $\overline x:Spec(K^{sep})\to C$. Consider the following filtrant subsets: $\cal C_{\overline x}\subset \cal C$ contains only the $A'$'s which are \'etale over $A$, $\cal C_{\overline x}'\subset \cal C_{\overline x}$ contains only the $A'$'s whose fractions field $K'$ is Galois over $K$ and whose spectrum has no point above $v$. We have a series of functorial isomorphisms
\begin{eqnarray}(\iota_{C_v,*}F)_{\overline x}&\simeq & \limi_{A'\in \cal C'_{\overline x}} \G(A_v\otimes_A A',F)\\
&\simeq & \limi_{A'\in \cal C'_{\overline x}} \G(K_v\otimes_K K',F)\\ &\simeq& \limi_{A'\in \cal C'_{\overline x}} \G(Ind_{G'}^{D'_{v'}}K'_{v'},F)\\ &\simeq& \limi_{A'\in \cal C'_{\overline x}} Ind_{G'}^{D'_{v'}}\G(K'_{v'},F)\\
&\simeq & Ind_G^{D_{v^{sep}}}F_{\overline \eta_v}\end{eqnarray}
where $\overline \eta_v:Spec(K_v^{sep})\to Spec(A_v)$ is the geometric point satisfying $\iota^{sep}\overline \eta_v=\overline x$ and $Ind_G^{D_{v^{sep}}}$ means discrete induction.
The first isomorphism comes from the fact that $\cal C_{\overline x}'$ is cofinal inside $\cal C_{\overline x}$. The second one occurs since $Spec(A')$ has no points above $v$. The third one is induced by the isomorphism \begin{eqnarray}\label{indK}K_v\otimes_K K'&\simeq & Ind_{G'}^{D'_{v'}}K'_{v'}\\a\otimes b&\mapsto &(g\mapsto a\iota^{sep}(gb))\end{eqnarray} The fourth one is by pseudo-commutativity of the first square in the above diagram of classifying topoi. Since discrete induction commutes to filtrant
inductive limits, the fifth one comes from the fact that the functor $\cal C'_{\overline x}\to \cal C_{\overline \eta_v}$, $A'\mapsto K'_{v'}$, where $\cal C_{\overline \eta_v}$ denotes the set of \'etale sub algebras of $K_v^{sep}/K_v$, is cofinal. The desired exactness follows from the fact that $Ind_G^{D_{v^{sep}}}$ is exact.
\petit

\emph{Case 2}. The geometric point $\overline x$ lies over a closed point $w\ne v$.  Consider the following filtrant sets: $\cal C_{\overline x}\subset \cal C$ contains only the \'etale $A'/A$'s whose spectrum contains $w'$, $\cal C_{\overline x}'\subset \cal C$ contains only the $A'$'s satisfying 1) $A'/A$ is \'etale outside $w$, 2) $K'/K$ is Galois and 3) $A'$ is the integral closure of $A[1/s]$ in $K'$ for some $s$ satisfying $w(s)=0$, $v(s)>0$. We define a functor $(-)'':\cal C'_{\overline x} \to \cal C_{\overline x}$, $A'\mapsto A'':=A'^{I'_{w'}}[1/t]$ for some $t$ whose valuation is zero at $w'$ and at places which are not above $w$ and is non zero  at the other places above $w$. We have a series of functorial isomorphisms \begin{eqnarray}(\iota_{C_v,*}F)_{\overline x}&\simeq & \limi_{A'\in \cal C'_{\overline x}}\G(A_v\otimes_A A'',F)\\ &\simeq &\limi_{A'\in \cal C'_{\overline x}}\G((Ind_{G'}^{D'_{v'}}K'_{v'})^{I'_{w'}},F)\\
&\simeq &  \limi_{A'\in \cal C'_{\overline x}} (Ind_{G'}^{D'_{v'}}\G(K'_{v'},F))^{I'_{w'}}\\
&\simeq & (Ind_G^{D_{v^{sep}}} F_{\overline \eta_{v}})^{I_w^{sep}}\end{eqnarray}
where $\overline \eta_{v}$ is as in Case 1. The first isomorphism holds because the functor $(-)''$ is cofinal. The second one follows from (\ref{indK}) since $A_v\otimes_A A''=A_v\otimes_AA'^{I'_{w'}}$ and $Spec(A')$ has no points above $v$. The third isomorphism follows from the diagram of classifying topoi, by pseudo-commutativity of the first and third squares, together with the base change isomorphism in the second square. The last isomorphism is similar to Case 1. Exactness follows using that $D_{v^{sep}}$ intersects trivially the conjugacy class of $I_{w^{sep}}$ by Krasner's lemma and the approximation lemma of \cite{Se1} I, p.23.
\petit

\emph{Case 3}. The geometric point $\overline x$ lies above $v$. Consider the following filtrant sets: $\cal C_{\overline x}\subset \cal C$ contains only the \'etale $A'/A$'s whose spectrum contains $v'$, $\cal C_{\overline x}'\subset \cal C$ contains only the $A'$'s satisfying 1) $A'$ is \'etale outside $v$, 2) $K'/K$ is Galois and 3) $A'$ is the integral closure of $A[1/s]$ in $K'$ for some $s$ satisfying $v(s)=0$. We define a functor  $(-)'':\cal C'_{\underline x} \to \cal C_{\underline x}$, $A'\mapsto A'':=A'^{I'_{v'}}[1/t]$ for some $t$ whose valuation is zero at $v'$ and at places which are not above $v$ and is non zero  at the other places above $v$. We have the following series of functorial isomorphisms
\begin{eqnarray}(\iota_{C_v,*}F)_{\overline x}&\simeq & \limi_{A'\in \cal C'_{\overline x}}\G(A_v\otimes_A A'',F)\\ &\simeq &\limi_{A'\in \cal C'_{\overline x}}\G((Ind_{G'}^{D'_{v'}}K'_{v'})\times_{K'_{v'}} A'^{I'_{v'}}_{v'},F)\\
&\simeq &  \limi_{A'\in \cal C'_{\overline x}} (Ind_{G'}^{D'_{v'}}\G(K'_{v'},F))\times_{\G(K'_{v'},F)} \G(A'^{I'_{v'}}_{v'},F) \\
&\simeq & (Ind_G^{D_{v^{sep}}} F_{\overline \eta_v})\times_{F_{\overline \eta_v}} F_{\overline x_v}\end{eqnarray}
where $\overline \eta_v$ is as in Case 1 and $\overline x_v$ is the geometric point of $Spec(A_v)$ inducing $\overline x$. The first isomorphism holds because $(-)''$ is cofinal. The second isomorphism is deduced from (\ref{indK}) by taking fixed points under $I'_{v'}$ and imposing integrality at $v'$ (note that the fiber product is taken with respect to the adjunction map, explicitly described as $f\mapsto f(1)$ on the left and the inclusion map on the right).
The third isomorphism is obtained by descent along the surjective \'etale morphism $$Spec(Ind_{G'}^{D'_{v'}}K'_{v'})\sqcup Spec(A'^{I'_{v'}}_{v'})\to Spec((Ind_{G'}^{D'_{v'}}K'_{v'})\times_{K'_{v'}} A'^{I'_{v'}}_{v'})$$ The last isomorphism holds because $\cal C'_{\overline x}\to \cal C_{\overline \eta_v}$, $A'\mapsto A'_{v'_1}$ and $\cal C'_{\overline x}\to \cal C_{\overline x_v}$, $A'\mapsto A''_{v''}$ are cofinal if $v'_1\ne v'$ is a place of $K'$ above $v$ and $\cal C_{\overline x_v}\subset \cal C$ contains only the \'etale $A'/A$'s whose spectrum contains $v'$. Exactness then follows using the isomorphism $$\begin{array}{rcl}(Ind_G^{D_{v^{sep}}} F_{\overline \eta_v})\times_{F_{\overline \eta_v}} F_{\overline x_v}&\simeq &Cont(D_{v^{sep}}\backslash (G-D_{v^{sep}}),F_{\overline \eta_v})\times F_{\overline x_v}\\ (f,a)&\mapsto& (f\circ s,a)\end{array}$$ provided by a continuous section $s$ of $G\to D_{v^{sep}}\backslash G$ (\cite{Se2} 1.2, Prop. 1).
 \begin{flushright}$\square$\end{flushright}

\para \label{parasmash} We now introduce some \emph{ad hoc} functors designed to neglect the generic fiber of $C_v$. Recall the well known equivalence between $C_{v,et}$ and triples $(F_1,F_2,f)$
where $F_1\in Z_{v,et}$,  $F_2\in U_{v,et}$ and  $f:F_1\to z_v^{-1}j_{v,*}F_2$ (\cite{Mi1} II, Thm. 3.10).


\begin{defn} \label{defsma} \emph{(Smashing functors)}

\debrom \item \label{defsmai} The functor $$Sma_v: (C_v\leftarrow U_v)_{et}\to (C_v\leftarrow U_v)_{et}$$ is defined by sending $F_1\to F_2$ to
$F_3\to F_2$ where $F_3$ corresponds to the triple $(z_v^{-1}F_1,F_2,z_v^{-1}F_1\to z_v^{-1}j_{v,*}F_2)$.

\item \label{defsmaii} The functor $$Sma:J_{et}\to J_{et}$$ is defined by sending $(F_{C_v}\to F_{U_v}\leftarrow F_U)$ to $(G_{C_v}\to F_{U_v}\leftarrow
F_U)$ where $(G_{C_v}\to F_{U_v})=Sma_v(F_{C_v}\to F_{U_v})$. \finrom
\end{defn}

These functors have the following properties.

\begin{lem} \label{lemsma}
\debrom
\item \label{lemsmai}
The functor $Sma_v$ (resp. $Sma$) is exact. It induces in particular an endofunctor of the categories  $\smash{Mod((C_v\leftarrow U_v)_{et}^\N,\Z/p^.)}$ and $\smash{D((C_v\leftarrow U_v)_{et}^\N,\Z/p^.)}$ (resp. $\smash{Mod(J_{et}^\N,\Z/p^.)}$ and $\smash{D(J_{et}^\N,\Z/p^.)})$.

\item \label{lemsmaii} The morphism $F(j_v):j_v^{-1}F_{C_v}\to F_{U_v}$  for variable $F$ induces natural transformations $Sma_v\to Id$ (resp. $Sma\to Id$). If $F(j_v)$ (resp. each $F(j_v)$) is invertible, then the previous natural transformation induces an isomorphism: $$\begin{array}{r}Sma_v(F)\simeq F\\ \hbox{(resp. $Sma(F)\simeq F$)}\end{array}$$

\item \label{lemsmaiii} Consider an object $M=\smash{(M_{C_v}\to M_{U_v})}$ (resp. $M=\smash{(M_{C_V}\to M_{U_v}\leftarrow M_{U})}$) of $\smash{D((C_v\leftarrow U_v)_{et}^\N,\Z/p^.)}$ (resp. $\smash{D(J_{et}^\N,\Z/p^.)}$). Then $$\begin{array}{r}Sma_v(M)=0\\ \hbox{(resp. $Sma(M)=0$)} \end{array}$$ if and only if  $\smash{M_{U_v}}=0$ and $\smash{z^{-1}_vM_{C_v}}=0$  (resp. $M_{U}=0$, $M_{U_v}=0$ and $\smash{z^{-1}_vM_{C_v}}=0$).
\finrom \end{lem}
Proof. Everything follows directly from the definition of $Sma$ and $Sma_v$.
\begin{flushright}$\square$\end{flushright}

\begin{cor} \label{corsma} Consider a morphism $f:M\to N$ in $\smash{D(J_{et}^\N,\Z/p^.)}$. The induced morphism $Sma(M)\to Sma(N)$ is invertible if and only if the following morphisms are invertible: $f_U:M_U\to N_U$, $\smash{f_{U_v}}:\smash{M_{U_v}}\to \smash{N_{U_v}}$ and $\smash{z_v^{-1}f_{C_v}}:\smash{z_v^{-1}M_{C_v}}\to \smash{z_v^{-1}N_{C_v}}$. \end{cor}
Proof. Apply Lem. \ref{lemsma} \ref{lemsmaiii} to a cone of the morphism $f$.
\begin{flushright}$\square$\end{flushright}

\subsection{Vanishing cohomology} \label{vc} ~~\\

We define vanishing cohomology in a general setting and study its behaviour when changing the topos.


\para We begin with simple diagram theoretic constructions regarding a morphism of diagrams $i:Y\to X$. In the applications, $i$ will be taken either as the closed immersion $z:Z\to C$, the closed immersion $z_J:Z_J\to J$ or slight variants of these.

\begin{defn} \label{defitype}
\debrom
\item \label{defitypei} Consider small categories $\Delta$ and $\Delta'$. We say that a functor $i:\Delta'\to \Delta$ is \emph{extremal} if it is injective on objects, fully faithful, and if $\Delta'/\delta$ is empty  for $\delta\notin i(\Delta')$ (ie. $Hom(i(\delta'),\delta)=\emptyset$ for $\delta'\in \Delta'$, $\delta\notin i(\Delta')$).

\item \label{defitypeii} Let $i:\Delta'\to \Delta$ be an extremal functor. We denote $\Delta^+$ the full subcategory of $\Delta\times [1]$ formed by the $(\delta,0)$'s with $\delta\in i(\Delta')$ and all $(\delta,1)$'s. Here, $[1]$ denotes the category $\{0\to 1\}$. We have natural functors $i^+:\Delta'\to \Delta^+$, $\delta'\mapsto (i(\delta'),0)$,  $\sigma:\Delta\to \Delta^+$, $\delta\mapsto (\delta,1)$ and $\rho:\Delta^+\to \Delta$, $ (\delta,n)\mapsto \delta$.

\item \label{defitypeiii} We say that a morphism of diagrams $i:Y/\Delta'\to X/\Delta$ is of \emph{extremal type} if the underlying functor $\Delta'\to \Delta$ is extremal. For $i$ of extremal type, we denote $$X^+=(Y\to X)$$ the diagram of type $\Delta^+$ having $\smash{X^+_{i(\delta'),0}}=\smash{Y_{\delta'}}$ ($\delta'$ in $\Delta'$) and $\smash{X^+_{\delta,1}}=X_\delta$ ($\delta$ in $\Delta$) as vertices and whose edges are the following:

     - $X^+(i(f'),id_0):{X^+_{i(\delta'_1),0}} \to \smash{X^+_{i(\delta'),0}}$ is $Y(f'):Y_{\delta'_1}\to Y_{\delta'}$  for $f':\delta'_1\to \delta'$ in $\Delta'$,

     - $X^+(f,id_1):{X^+_{\delta_1,1}}\to \smash{X^+_{\delta,1}}$ is $X(f):X_{\delta_1}\to X_{\delta}$  for $f:\delta_1\to \delta$ in $\Delta$ and

     - $X^+(i(f'),0\to 1):{X^+_{i(\delta_1'),0}}\to \smash{X^+_{i(\delta'),1}}$ is $i_{\delta'}\circ Y(f'): Y_{\delta_1'}\to X_{\delta'}$ (also equal to $X(i(f'))\circ i_{\delta'}$) for $f':\delta'_1\to \delta'$ in $\Delta'$. \petit

     The functors  $i^+$, $\sigma$ and $\rho$ defined in \ref{defitypeii} uniquely extend to morphisms of diagrams $$\xymatrix{Y\ar @<+2pt> `d[r] `[rr]_-{i} [rr]\ar[r]^-{i^+}&X^+\ar[r]^-\rho&X\ar @/_1pc/[l]_{\sigma}}$$ such that $\rho\circ i^+=i$, $\rho\circ \sigma=id_X$, $i^+_{\delta'}:\smash{Y_{\delta'}\to X_{i(\delta'),0}}$  is the identity for all $\delta'\in\Delta'$ and $\sigma_\delta:X_\delta\to \smash{X^+_{\delta,1}}$ is the identity for all $\delta\in\Delta$.

\finrom
\end{defn}

\begin{lem} \label{lemitype} Consider a ringed variable topos $(\cal T,A)$ on $\cal B$ (Def. \ref{deftopfib} \ref{deftopfibii}). Let $i$ be of extremal type as in  Def. \ref{defitype} \ref{defitypeiii}.
\debrom
\item \label{lemitypei} The direct image functors induced by $i:Y\to X$ and $i^+:Y\to X^+$ have the following description:

- $i_*:Mod(\cal T(Y),A_Y)\to Mod(\cal T(X),A_X)$ satisfies $(i_*M)_\delta={i_{\delta',*}}M_{\delta '}$ if $i(\delta')=\delta$ and $(i_*M)_\delta=0$
if $\delta\notin i(\Delta')$. A module $M$ is acyclic for $i_*$ if and only if each $M_{\delta'}$ is acyclic for $i_{\delta',*}$.

- $i^+_*:Mod(\cal T(Y),A_Y)\to Mod(\cal T(X^+),A_{X^+})$ satisfies $(i^+_*M)_{i(\delta'),0}=M_{\delta'}$,
$(i^+_*M)_{i(\delta),1}={i_{\delta',*}}M_{\delta'}$ and $(i^+_*M)_{\delta,1}=0$ is $\delta\notin i(\Delta')$.
A module $M$ is acyclic for $i^+_*$ if and only if each $M_{\delta'}$ is acyclic for $i_{\delta',*}^+$.

\item \label{lemitypeii} There is a natural isomorphism $\sigma^{-1}A_{X^+}\simeq A_X$. The functors $\sigma^{-1}$ and $\rho_*:$ \break $Mod(\cal T(X^+),A_{X^+})\to Mod(\cal T(X),A_X)$ are naturally isomorphic. In particular, they commute to arbitrary limits. The functors $\sigma_*$ and $\rho^*$ are fully faithful.
\finrom\end{lem}

Proof. \ref{lemitypei} The calculation of $i_*$ and $i^+_*$ is immediate from the projective limit formula of Lem. \ref{lemfibtop} \ref{lemfibtopii}. Since injective modules of $(\cal T(Y),A_Y)$ have flasque components by Lem.-Def. \ref{acycf} \ref{acycfv}, a similar formula holds for $Ri_*$ and $Ri^+_*$. The acyclicity statements follow.

\ref{lemitypeii} Here again, the calculation of $\rho_*$ in \ref{lemitypeii} is straightforward from the projective limit formula. The last statement of \ref{lemitypeii} follows from the isomorphism $\sigma^{-1}\simeq \rho_*$, given that
$\rho\sigma=id_X$.
\begin{flushright}$\square$\end{flushright}

In Lem. \ref{lemitype}, we have used the simplified notation $\cal T(X)$ instead of $\cal T^{codiag}(X)$ or $\cal T_{cof}(X)$. We will continue to do so
from now on. We will sometimes use the suggestive notation $(-)_{|X}$ (resp. $(-)_{|Y}$, resp. $(-)_{|X^+}$) instead of $\sigma^*$ (resp. $i^*$ or
$i^{+,-1}$, resp. $\rho^*$) when no danger of confusion arises.

\para \label{vancoh} We are now ready to define the functors of vanishing sections and study their elementary properties in the general context of fibered topoi.

\begin{defn} \label{defvan} Let $\cal T$, $\cal B$, $A$, $i:Y\to X$ and $X^+$ be as in Lem. \ref{lemitype}.
\debrom
\item \label{defvani} We define the following functors for modules on $(\cal T(X),A_X)$:


 - The functor of \emph{local sections vanishing at $Y$} $$\underline \G^Y(-):Mod(\cal T(X),A_X)\to Mod(\cal T(X),A_X)$$ is defined by $\underline \G^Y(M):=Ker(M\rightarrow i_*i^*M)$.

- Given $f:X'\to X$ in $Diag(\cal B)$, the functor of 
\emph{local sections on $X'$ vanishing at $Y$} $$\underline \G^Y(X',-):Mod(\cal T(X),A_X)\to Mod(\cal T(X'),A_X')$$
is defined as $\underline \G^Y(X',-)=f^*\underline \G^Y(-)$. The functor of 
\emph{sections on $X'$ vanishing at $Y$} $$\G^Y(X',-):Mod(\cal T(X),A_X)\to Mod(\Gamma(X',A_X'))$$
is defined as $\G^Y(X',-)=\G(X',f^*\underline \G^Y(-))$. 

\item \label{defvanii} As a particular case of \ref{defvani}, the following functors are attached to $i^+:Y\to X^+$ and $\sigma:X\to X^+$ (recall that $\sigma^*\simeq \sigma^{-1}\simeq \rho_*$ since $\sigma^{-1}A_{X^+}\simeq \rho_*A_{X^+}\simeq A_X$):

- The functor of \emph{local sections on $X$ vanishing at $Y$} $$\underline \G^Y(X,-):Mod(\cal T(X^+),A_{X^+})\to Mod(\cal T(X),A_X)$$ sends $M$ to $\underline \G^Y(X,M)=Ker(\sigma^{-1}M\to i_*i^{+,-1}M)$, where the arrow is defined using the adjunction morphism $id\to i^+_*i^{+,*}$ and the isomorphisms $\sigma^{-1}i^+_*\simeq \rho_*i^+_*\simeq i_*$ (explicitly the arrow is thus described on component $\delta$ as $M_{\delta,1}\to i_{\delta',*}M_{\delta,0}$ if $\delta=i(\delta')$ and $M_{\delta,1}\to 0$ if $\delta\notin i(\Delta')$).
    
 - The functor of \emph{sections on $X$ vanishing at $Y$} $$\G^Y(X,-):Mod(\cal T(X^+),A_{X^+})\to Mod(\Gamma(X,A_X))$$ sends $M$ to $\G^Y(X,M)=Ker(\G(X,\sigma^{-1}M)\to \G(Y,i^{+,-1}M))$.
\finrom
\end{defn}

\begin{rem} Let us emphasize that the functors $\underline \G^Y$ should not be mistaken for  a functor of the type ``local sections with
support in $Y$'' whose usual notation is $\underline \G_Y$. Let us explain the relation to the latter in a specific situation. Set $A=\Z$ for simplicity. Assume that $i:Y\to X$ induces an open immersion of ringed topoi $\cal T(Y)\to \cal T(X)$ (in practice, this will be the case only when $\cal T$ is one of the usual big topoi). Let then $c:\cal T(Y)^c\to \cal T(X)$ denote the inclusion of the complementary closed subtopos (\cite{SGA4-I} IV, Prop. 9.3.4). Then $\underline \G^Y\simeq  \smash{\underline\G_{\cal T(Y)^c}}:=c_*c^!$ is the functor of local sections with support in $Y^c$. \end{rem}


Let us introduce the following condition, which will be helpful for technical reasons: \petit

\noindent $opimm(f,\cal T,A)$: the morphism of ringed topoi $(\cal T(X_{1,\delta}),A_{X_{1,\delta}})\to (\cal
T(X_{2,f(\delta)}),A_{X_{2,f(\delta)}})$ induced by $f_\delta$ is an open immersion (ie. is isomorphic to the morphism of localization attached to an open of the topos $\cal T(X_{2,f(\delta)})$)  for all $\delta$ in $\Delta_1$). \petit

\petit


Let us now gather some properties of the derived functors $R\underline \G^Y(-)$ and $R\underline \G^Y(X,-)$.

\begin{lem} \label{soritesvan} Keep the notations and assumptions of Def. \ref{defvan}.
\debrom
\item \label{soritesvani} There are canonical isomorphism and distinguished triangle $$\begin{array}{c}R\underline \G^Y(X,M)\simeq MF(R(\sigma^{-1}\to i_*i^{+,-1})(M))\\
    \xymatrix{R\underline \G^Y(X,M)\ar[r]&M_{|X}\ar[r]&Ri_*M_{|Y}\ar[r]^-{+1}&}\end{array}$$
which are functorial with respect to $M$ in $D^+(\cal T(X^+),A_{X^+})$. If $M$ is a complex such that the objects of $M_{|Y}$ are acyclic for $i_*$ (e.g. if the objects of $M$ have flasque components by Lem. \ref{lemitype} \ref{lemitypei}) then the right hand side in the above isomorphism is canonically isomorphic to the mapping fiber of $M_{|X}\to i_*M_{|Y}$.

\item \label{soritesvanii} Assume that $opimm(i,\cal T,A)$ holds. There are canonical isomorphism and distinguished triangle $$\begin{array}{c}R\underline \G^Y(M)\simeq MF(R(id \to i_*i^{-1})(M))\\ \xymatrix{R\underline \G^Y(M)\ar[r]&M\ar[r]&Ri_*M_{|Y}\ar[r]^-{+1}&}\end{array}$$
which are functorial with respect to $M$ in $D^+(\cal T(X),A_X)$. If $M$ is a complex such that the objects of $M_{|Y}$ are acyclic for $i_*$ (e.g. if the objects of $M$ have flasque components by Lem. \ref{lemitype} \ref{lemitypei}) then the right hand side in the above isomorphism is canonically isomorphic to the mapping fiber of $M\to i_*M_{|Y}$.

\item \label{soritesvaniii} Assume that $opimm(i,\cal T,A)$ holds. There is a  canonical isomorphism
$$R\underline \G^Y(X,\rho^{-1}M)\simeq R\underline \G^Y(M)$$
    which is functorial with respect to $M$ in  $D^+(X,A_X)$.

\item \label{soritesvaniv} Consider another morphism $i_1:Y_1/\Delta_1'\to X_1/\Delta_1$ of extremal type. Assume given some morphisms of diagrams $f:X_1\to X$ and $f_Y:Y_1\to Y$ satisfying $fi_1=if_Y$ and denote $f^+:X_1^+\to X^+$ the induced morphism. There is a canonical isomorphism $$R\underline \G^{Y}(X,Rf^{+}_*M)\simeq  Rf_*R\underline \G^{Y_1}(X_1,M)$$ which is functorial with respect to $M$ in $D^+(\smash{\cal T(X_1^+),A_{X_1^+}})$.

%
%

\item \label{soritesvanv} Let $g:(\cal T,A)\to (\cal T',A')$ be a morphism of variable topoi.  There is a canonical isomorphism $$R\underline \G^Y(X,Rg_*M)\simeq Rg_*R\underline \G^Y(X,M)$$ which is functorial with respect to $M$ in $D^+(\cal T(X^+),A_{X^+})$.
\finrom
\end{lem}
Proof. 
\ref{soritesvani}  Let $f:|X^+|\to X^+$ denote the inclusion of the discrete diagram underlying $X^+$.
As in the proof of Lem. \ref{defMV} \ref{defMVii}, the stated isomorphism will follow once checked that $M_{|X}\to i_*M_{|Y}$ is epimorphic when $M$ is of the form $f_*N$ for some $\smash{f^{-1}A_{X^+}}$-module $N$ on $|X^+|$. If $\delta=i(\delta')$, we have natural isomorphisms $$\begin{array}{rcl}(f_*N)_{\delta,1}&\simeq &\prod_{(g':\delta_1'\to \delta')\in\Delta'/\delta'}  (i(g'),0\to 1)_*N_{i(\delta_1'),0} \times (i(g'),id_1)_*N_{i(\delta'_1),1} \\
(f_*N)_{\delta,0}&\simeq &\prod_{(g':\delta_1'\to \delta')\in \Delta'/\delta'} (i(g'),id_0)_*N_{i(\delta'_1),0}\end{array}$$
Using this description, we see that $(f_*N)_{\delta,1}\twoheadrightarrow i_{\delta',*} (f_*N)_{\delta,0}$ is split epimorphic. The distinguished triangle and the last assertion follow from the stated isomorphism since $R(i^+_*i^{+,-1})\simeq (Ri^+_*)i^{+,-1}$ (note that $i^{+,-1}$ sends injectives to componentwise flasques by Lem.-Def. \ref{acycf} \ref{acycfv}, hence is $i^+_*$-acyclic by Lem. \ref{lemitype} \ref{lemitypei}).

\ref{soritesvanii} Let $f:|X|\to X$ denote the inclusion of the discrete diagram underlying $X$. Once again, it suffices to check that  $M\to i_*i^{-1}M$ is
epimorphic for $M=f_*N$ with $N$ injective. Using  \cite{SGA4-II} V, Prop. 4.7, this follows from $opimm(i,\cal T,A)$, thanks to the fact that each  $(f_*N)_\delta$ is flasque (Lem.-Def. \ref{acycf} \ref{acycfv}). The distinguished triangle and the last assertion follow as in \ref{soritesvani} (here $i^{-1}$ preserves injectives thanks to $opimm(i,\cal T,A)$).

\ref{soritesvaniii} Since $id\to i_*i^{-1}$ identifies with $\sigma^{-1}\rho^{-1}\to i_*i^{+,-1}\rho^{-1}$, we have a natural isomorphism $\underline \G^Y(X,\rho^{-1}(-))\simeq \underline \G^Y(-)$. In order to pass to derived functors, we need to know that $\rho^{-1}$ sends injectives to $\underline \G^Y(X,-)$-acyclics. Consider an injective module $M$ over $(X,A_X)$. Then $\rho^{-1}M$ has flasque components and $\sigma^{-1}\rho^{-1}M\to i_*i^{+,-1}\rho^{-1}M$ (ie. $M\to i_*i^{-1}M$) is epimorphic as already observed in the proof of \ref{soritesvani}. It thus follows from \ref{soritesvani} that $\rho^{-1}M$ is
$\underline \G^Y(X,-)$-acyclic as desired.



\ref{soritesvaniv} For any module $M$ over $(\cal T(X_1^+),A_{X_1^+})$, the left commutative diagram below
\begin{eqnarray}\label{MVfl1}\begin{array}{ccc}
\xymatrix{Y_1\ar[d]_{f_Y}\ar[r]^-{i_1^+}\ar @<+2pt> `u[r] `[rr]^-{i_1} [rr]&X_1^+\ar[d]^{f^+}\ar[r]^-{\rho_1}&X_1\ar[d]^f\\
Y\ar[r]^-{i^+}\ar @<+2pt> `d[r] `[rr]_-{i} [rr]&X^+\ar[r]^-\rho&X}
&\hspace{2cm}&
$$\xymatrix{f_*i_{1,*}i_1^{+,-1}M&f_*\rho_{1,*}M\ar[l]\\
i_*i^{+,-1}f^+_*M\ar[u]_\wr&\rho_*f^+_*M\ar[u]_\wr\ar[l]}$$
\end{array}\end{eqnarray}
induces the right commutative square in $Mod(\cal T(X),A_X)$ (use the projective limit formula of Lem. \ref{lemfibtop} \ref{lemfibtopii}) and the fact that $i_1$ is extremal to check the base change isomorphism $\smash{f_{Y,*}i_1^{+,-1}\simeq i^{+,-1}f^+_*}$ underlying the first vertical isomorphism).
Deriving the bottom arrow with
respect to $M$ gives $R(i_*i^{+,-1}\to \rho_*)\circ Rf^+_*(M)$, since $f^+_*$ sends injectives to componentwise flasques (hence $i_*i^{+,-1}$-acyclics thanks to $R(i^+_*i^{+,-1})\simeq (Ri^+_*)i^{+,-1}$ and Lem. \ref{lemitype} \ref{lemitypei}) by Lem.-Def. \ref{acycf} \ref{acycfiv}, \ref{acycfv}. Deriving the top arrow on the other hand, gives
$Rf_*\circ R(\rho_{1,*}\to i_{1,*}i_1^{+,-1})(M)$ since $\rho_{1,*}$ and $i_{1,*}i_1^{+,-1}$ send injectives to $f_*$-acyclics. The latter claim deserves an explanation. The functors $\rho_{1,*}$ and $i_{1,*}$ both preserve d-injectives. Since $i_1$ is of extremal type, it may be checked easily that $i_1^{+,-1}$ preserves d-injectives as well. In particular both $\rho_{1,*}$ and $i_{1,*}i_1^{+,-1}$ send injectives to direct factors of d-injectives (hence $f_*$-acyclics by Lem.-Def. \ref{acycf} \ref{acycfvi}). The
desired isomorphism comes by functorial mapping fibers.

%
%

\ref{soritesvanv} The construction of the morphism is similar to the first part of  \ref{soritesvaniv} (replace $f$ by $g$) using only that $i_*i^{+,-1}$ sends injectives to componentwise flasques (hence $g_*$-acyclic by Lem. \ref{DfuncT} \ref{DfuncTi}).
\begin{flushright}$\square$\end{flushright}

Let us emphasize the following result which justifies the introduction of $X^+$ and Def. \ref{defvan} \ref{defvanii}. In practice, it will be applied when $g$ is the morphism $\epsilon:(-)_{FL}\to (-)_{et}$.

\begin{cor} \label{corsoritesvan} Keep the notations and assumptions of Def. \ref{defvan}. Consider $g:(\cal T,A)\to (\cal T',A')$ as in Lem. \ref{soritesvan} \ref{soritesvanv}. If $opimm(i,\cal T,A)$ holds, there is a canonical isomorphism $$Rg_*R\underline \G^Y(M)\simeq R\underline \G^Y(X,Rg_*M_{|X^+})$$ which is functorial with respect to $M$ in $D^+(\cal T(X),A_X)$.
\end{cor}
Proof. This follows from Lem. \ref{soritesvan} \ref{soritesvaniii} and \ref{soritesvanv}.
\begin{flushright}$\square$\end{flushright}

\begin{lem} \label{resvan} Let $i:X\to Y$, $\cal B$ and  $\cal T$ be as in Def. \ref{defvan} and consider the variable topos of projective systems $\cal T^\N$ together with a ring $A_.$ in $\G(\cal B^{op},\cal T^\N)$. Let $\iota_k:(\cal T,A_k)\to (\cal T,A_.)$ denote the natural morphism (so that $\iota_k^*M_.=M_k$).

\debrom
\item \label{resvani} There is a canonical isomorphism
$$\iota_k^{-1}R\underline \G^Y(X,M_.)\simeq R\underline \G^Y(X,M_k)$$
which is functorial with respect to $M_.$ in $\smash{D^+(\cal T^\N(X^+), A_{.,X^+})}$.
If $opimm(i,\cal T,A)$ holds, there is also a canonical isomorphism $$\iota_k^{-1}R\underline \G^Y(N_.)\simeq R\underline \G^Y(N_k)$$
which is functorial with respect to $N_.$ in $\smash{D^+(\cal T^\N(X), A_{.,X})}$.
\item \label{resvanii} Consider a subcategory $\Delta_1$ of $\Delta$ and set $\Delta'_1=\Delta'\times_\Delta \Delta_1$. Let $X_1/\Delta_1$ (resp. $Y_1/\Delta'_1$) denote the restriction of $X$ (resp. $Y$) to $\Delta_1$ (resp. $\Delta_1'$) and denote $f:X_1\to X$ and $f^+:X_1^+\to X^+$ the inclusion morphisms. There is a canonical isomorphism $$f^{-1}R\underline \G^Y(X,M_.)\buildrel{\sim}\over\to R\underline \G^{Y_1}(X_1,f^{-1}M_.)$$
which is functorial with respect to $M_.$ in $\smash{D^+(\cal T^\N(X^+), A_{.,X^+})}$. If $opimm(i,\cal T,A)$ holds, there is also a canonical isomorphism $$f^{-1}R\underline \G^Y(M_.)\buildrel{\sim}\over\to R\underline \G^{Y_1}(f^{-1}M_.)$$
which is functorial with respect to $N_.$ in $\smash{D^+(\cal T^\N(X), A_{.,X})}$.
\finrom
\end{lem}
Proof. Thanks to Lem. \ref{soritesvan} \ref{soritesvaniii}, it suffices to establish the first isomorphism in \ref{resvani} and \ref{resvanii}. Now both cases follow from Lem. \ref{soritesvan} \ref{soritesvani} using the following remark. If $M_.$ is injective in $Mod(\cal T^\N(Y), A_{.,Y})$, then each $M_{k,\delta'}$ is flasque (use Lem.-Def. \ref{acycf} \ref{acycfv} applied to the variable topos $\cal T^{-}(-)$ on $\cal B\times \N^{op}$). In particular $\iota_k^{-1}M_.$ is acyclic for $i_*: Mod(\cal T(Y),A_{k,Y})\to Mod(\cal T(X),A_{k,X})$ and  $f^{-1}M_.$ is acyclic for  $i_*:Mod(\cal T(Y_1)^\N,A_{.,Y_1})\to Mod(\cal T(X_1)^\N,A_{.,X_1})$
(Lem. \ref{lemprojcomp} \ref{lemprojcompi} and Lem. \ref{lemitype} \ref{lemitypei}).

\begin{flushright}$\square$\end{flushright}

\subsection{Complete Mayer-Vietoris for semi-Abelian schemes over $C$} \label{cmvfnm} ~~ \\

We establish a Mayer-Vietoris isomorphism (and triangle) relatively to the complete neighborhoods of the points of $Z\subset C$. We begin by recalling an acyclicity result with regards to the morphism $$\epsilon:(-)_{FL}\to (-)_{et}$$

\begin{lem} \label{acycsm} Consider a morphism of extremal type $i:Y/\Delta'\to X/\Delta$ between diagrams of schemes and assume that each $i_{\delta'}:Y_{\delta'}\to X_\delta$ is a closed immersion. Consider an Abelian group $M$ in $X_{FL}$ (resp. a module $M_.$ on $(X_{FL}^\N,\Z/p^.)$) and assume that $M=i_*N$ (resp.  $M_.=i_*N_.$) where each $N_{\delta'}$ (resp. each $N_{\delta',k}$) is representable by a smooth group scheme.

\debrom \item \label{acycsmi} The group $M$ (resp. the module $M_.$) is acyclic for $\epsilon_*$. The group $N$ (resp. $N_.$) is acyclic for $\epsilon_*$ and $i_*$.
\item \label{acycsmii} Both $R\underline\G^Y(M)$ and $R\underline\G^Y(X,\epsilon_*M_{|X^+})$ (resp.  $R\underline\G^Y(M_.)$ and $R\underline\G^Y(X,\epsilon_*M_{.,|X^+})$) are zero.
\finrom
\end{lem}
Proof. Using Lem. \ref{DfuncT} \ref{DfuncTi}, Lem. \ref{lemitype} \ref{lemitypei} and Lem. \ref{resvan} \ref{resvani}, \ref{resvanii} it is sufficient to consider the case where $X$ and $Y$ are schemes (ie. punctual diagrams) and $M$ is an Abelian group.

\ref{acycsmi} We know from \cite{Mi1} III, Thm. 3.9 or \cite{Gr1} Brauer III, Thm. 11.7, that $N$ is acyclic for $\epsilon_*$. Since moreover the functor $i_*:Y_{et}\to X_{et}$ is exact, we
find that $N$ is acyclic for the functor $i_*:Y_{FL}\to X_{FL}$ as well (sheafify the isomorphism $\smash{H^q(X'_{FL},Ri'_*N')}\simeq \smash{H^q(X'_{et},i'_*\epsilon_*N')}$ for $X'$ varying in
the big flat site of $X$, $i':Y'\to X'$ the base change of $i$ to $X'$ and $N'=\smash{N_{|Y'}}$).  Finally, the acyclicity of $M=i_*N$ for $\epsilon_*$ results from the isomorphisms $R\epsilon_*M\simeq
R\epsilon_*Ri_*N\simeq Ri_*R\epsilon_*N\simeq i_*\epsilon_*N$.

\ref{acycsmii} We remark that the assumption $opimm(i,(-)^\N_{FL},\Z/p^.)$ holds since arbitrary subschemes give rise to open subtopoi  in the big $top$ topos for as long $top$ is coarser than or equal to $fl$. It thus follows from Lem. \ref{soritesvan} \ref{soritesvanii} that $R\underline \G^Y(M)$ is the mapping fiber of $M\to i_*i^{-1}M$. Now, the latter arrow is invertible since $M$ is in the essential image of the fully faithful functor $i_*$. We have thus proven that $R\underline \G^Y(M)=0$. The case of $\smash{R\underline \G^Y(X,\epsilon_*M_{|X^+})}$ follows by Cor. \ref{corsoritesvan}.
\begin{flushright}$\square$\end{flushright}

We can now establish the following result which states that the diagram $R\epsilon_*A_{|J^+,p^.}$ is sufficient to retrieve the projection of $R\underline \G^Z(A_{p^.})$ to the small \'etale topos of $C$.

\begin{prop}\label{MVfl} Consider a semi-Abelian scheme $A/C$, whose restriction to $U$ is Abelian, and denote $A_{p^.}$ the projective system of its $p$-primary torsion sugroups. Recall the closed immersions of extremal type  $z:Z\to C$ and $z_J:Z_J\to J$ and the morphism $m:J\to C$ occurring in (\ref{diagdiag}). There is a canonical isomorphism $$R\epsilon_*R\underline \G^ZA_{p^.}\simeq Rm_*R\underline \G^{Z_J}(J,R\epsilon_*A_{|J^+,p^.})$$ in  $D^+(C^\N_{et},\Z/p^.)$. Moreover, either side of the isomorphism remains unchanged if $A_{p^.}$ is replaced by $A^0_{p^.}$ or $\Z/p^.\otimes_\Z^LA$.
\end{prop}
Proof. Consider the commutative diagram \begin{eqnarray}\label{diagJplus}\xymatrix{Z_J\ar[r]^-{z_J^+}\ar[d]_-{m_Z} &J^+\ar[d]_-{m^+}\ar[r]^-{\rho_J}&J\ar[d]^-m\ar @/^1pc/[l]^{\sigma_J}\\
Z\ar[r]^-{z^+}&C^+\ar[r]^-\rho&C\ar @/^1pc/[l]^{\sigma}}\end{eqnarray}
For any $M$ in $D^+(C^\N_{FL},\Z/p^.)$, $D^+(C_{FL},\Z/p^k)$ or $D^+(C_{FL})$ we have canonically \begin{eqnarray}\label{MV1} R\epsilon_*R\underline\G^Z(M)&\simeq & R\epsilon_* R\underline \G^Z(C,\rho^{-1}M)\\
\label{MV2}&\simeq &R\underline \G^Z(C,R\epsilon_*\rho^{-1}M)\\
\label{MV3}&\to &R\underline \G^Z(C,R\epsilon_*Rm^+_*m^{+,-1}\rho^{-1}M)\\
\label{MV4}&\simeq  &R\epsilon_*Rm_*R\underline \G^{Z_J}(J,m^{+,-1}\rho^{-1}M)\\
\label{MV5}&\simeq &Rm_*R\epsilon_*R\underline \G^{Z_J}(J,\rho_J^{-1}m^{-1}M)\\
\label{MV6}&\simeq & Rm_*R\underline \G^{Z_J}(J,R\epsilon_*\rho_J^{-1}m^{-1}M)\end{eqnarray} where (\ref{MV1}) is Lem. \ref{soritesvan} \ref{soritesvaniii},
(\ref{MV2}) is Lem. \ref{soritesvan} \ref{soritesvanv}, (\ref{MV3}) is the adjunction morphism for $m^+$, (\ref{MV4}) is Lem. \ref{soritesvan}
\ref{soritesvaniv}, \ref{soritesvanv}, (\ref{MV5}) is clear and (\ref{MV6}) is Lem. \ref{soritesvan} \ref{soritesvanv}. It thus remains to prove that
(\ref{MV3}) is invertible for $M=A_{p^.}$ or equivalently for $M=A_{p^k}$ (Lem. \ref{resvan} \ref{resvani}). \petit

\emph{Claim.} Consider $N=(N_Z\leftarrow N_C)$ in  $Mod(C^+_{FL},\Z)$ (the meaning of this notation is similar to the one of Sect. \ref{diagopi}). If $N_C=A$ or $N_C\simeq z_*z^{-1}N_C$ then the following natural morphism
is invertible: \begin{eqnarray}\label{claimMVfl1} R\epsilon_*N\to R\epsilon_*Rm^+_*m^{+,-1}N\end{eqnarray}

Let us first explain why the claim implies the proposition. The first case of the claim applied to $N=\rho^{-1}M$ shows that (\ref{MV3}) is invertible for $M=A\in
D(C_{FL})$, and thus also (by scalar extension via $\Z\to \Z/p^k$) for $M=\Z/p^k\otimes_\Z^LA\in D^+(C_{FL},\Z/p^k)$.  Since multiplication by $p^k$ is epimorphic on $A^0$ (see Lem. \ref{pdivfunct} \ref{pdivfunctiii}), we have a long exact sequence in the first line below, where $\Phi=A/A^0$ denote the group of components of $A$. The second line is the tautological distinguished triangle of truncation in $D(C_{FL},\Z/p^k)$: $$\begin{array}{c}\xymatrix{0\ar[r]&A^0_{p^k}\ar[r]&A_{p^k}\ar[r]&\Phi_{p^k}\ar[r]&0\ar[r]&A/p^k\ar[r]&\Phi/p^k\ar[r]&0}\\
\xymatrix{A_{p^k}\ar[r]&\Z/p^k\otimes_\Z^LA[-1]\ar[r]&A/p^k[-1]\ar[r]^{+1}&}\end{array}$$
Since $\Phi\simeq z_*z^{-1}\Phi$ (Lem. \ref{components} \ref{componentsii}), the second case of the claim implies that (\ref{MV3}) is
invertible for $M=\Phi/p^k\simeq A/p^k$ or $\Phi_{p^k}$. Using the distinguished triangle, we conclude that  (\ref{MV3}) is invertible, for $M=\smash{A_{p^k}}$, and finally for  $\smash{A^0_{p^k}}$ as well, using the exact sequence.  The last statement of the proposition follows from the above proof together with the vanishing of $R\underline \G^Z(\Phi/p^k)$ and $R\underline \G^Z(\smash{\Phi_{p^k}})$ (Lem. \ref{components} \ref{componentsi} and Lem. \ref{acycsm} \ref{acycsmii}). \petit

Let us now prove the claim, using the conservative couple of functors $z^{+,-1}:D^+(C^+_{et})\to D^+(Z_{et})$, $\sigma^{-1}:D^+(C^+_{et})\to D^+(C_{et})$. Some preliminary observations are in order. First, we note that $z^{+,-1}$ and $\sigma^{-1}$ both commute to $R\epsilon_*$ (because $z^+$ and $\sigma$ are inclusions of diagrams, see Lem. \ref{DfuncT} \ref{DfuncTi}). Next, we observe that the natural base change morphisms $z^{+,-1}Rm^+_*\to m_{Z,*}z_J^{+,-1}$ and $\sigma^{-1}Rm^+_*\to Rm_*\sigma_J^{-1}$ are both invertible (use the projective limit formula of Lem. \ref{lemfibtop} \ref{lemfibtopii} for $m^+_*$, note that $m_{Z,*}$ is exact and that $\sigma^{-1}=\rho_*$ preserves injectives). With these observations in mind, one may check without difficulty that the image of the morphism (\ref{claimMVfl1}) under $z^{+,-1}$ and $\sigma^{-1}$  respectively identifies with \begin{eqnarray}\label{claimMVfl2} &R\epsilon_*z^{+,-1}N\to R\epsilon_*m_{Z,*}m_Z^{-1}z^{+,-1}N\\
\label{claimMVfl3}\hbox{and}& R\epsilon_*\sigma^{-1}N\to R\epsilon_*Rm_*m^{-1}\sigma^{-1}N&\end{eqnarray} First, note that (\ref{claimMVfl2}) is clearly an
isomorphism, since  $m_Z$ is in fact an equivalence of  topoi. It thus remains to prove that the morphism \begin{eqnarray}\label{claimMVfl4}R\epsilon_*N_C\to R\epsilon_*Rm_*m^{-1}N_C\end{eqnarray} is invertible. Let us investigate the target of (\ref{claimMVfl4}). By Lem. \ref{defMV} \ref{defMVii} and Lem. \ref{iotaexact}, we have a distinguished triangle {\small$$\hspace{-.4cm}\xymatrix{R\epsilon_*Rm_*m^{-1}N_C\ar[r] &
(\mathop\prod\limits_{v\in Z}\iota_{C_v,*}R\epsilon_*(N_{C,|C_v}))\oplus Rj_*R\epsilon_*(N_{C,|U}) \ar[r]&\mathop\prod\limits_{v\in
Z}Rj_*\iota_{U_v,*}R\epsilon_*(N_{C,|U_v})\ar[r]^-{+1}&}$$}
In the case $N_C\simeq z_*z^{-1}N_C$, the third term and the second summand of the middle term vanish and the claim follows immediately. Let us now explain the case $N_C=A$. By \cite{Ra4} and the previous case, we may (and will) assume that $A$ is the N\'eron model of its generic fiber (note that the quotient $Q$ of the component group of the N\'eron model by $\Phi$ satisfies $Q\simeq z_*z^{-1}Q$ thanks to Lem. \ref{acycsm} \ref{acycsmi}).
Recall from Lem. \ref{acycsm} that $A$ is acyclic for $\epsilon_*$, as well as its restriction to $U$, $\cal C_v$ or $U_v$. Thus we only have to show that the above distinguished triangle induces an exact sequence
  \begin{eqnarray}\label{compl1}\xymatrix{\ 0\ar[r]&\epsilon_*A\ar[r]&(\mathop\prod\limits_{v\in Z}\iota_{C_v,*}\epsilon_*A_{|C_v})\oplus j_*\epsilon_*A_{|U}\ar[r]&\mathop\prod\limits_{v\in
Z}\iota_{C_v,*}j_{v,*}\epsilon_*(A)_{|U_v}\ar[r]&0}\end{eqnarray} and isomorphisms  \begin{eqnarray}\label{compl2}\xymatrix{R^qj_*\epsilon_*A_{|U} \ar[r]^-\sim&
\mathop\prod\limits_{v\in
Z}\iota_{C_v,*}R^qj_{v,*}\epsilon_*A_{|U_v}}\end{eqnarray} for $q\ge 1$.
The exactness of (\ref{compl1}) follows from the N\'eron extension property which implies that $\epsilon_*A\to j_*\epsilon_*A_{|U}$ and $\epsilon_*A_{|C_v}\to j_{v,*}\epsilon_*A_{|U_v}$ are invertible. Let us now investigate the stalks of the morphism (\ref{compl2}) at a geometric point $\overline x$ of $C$. Both sides give zero unless $\overline x$ lies above a point $v\in Z$. In the latter case we have on the one hand $$\begin{array}{rcl}(R^qj_*\epsilon_*A_{|U})_{\overline x}&\simeq &H^q(K^{sh},\epsilon_*A_{|K^{sh}})\\ & \simeq & \limi_{K'}H^q(K'^{h},\epsilon_*A_{|K'^{h}})\end{array}$$ where both isomorphisms follow from \cite{SGA4-II} VI, Prop. 5.7. Here $K^{sh}$ denotes the fraction field of  the  strictly henselian ring $\cal O_{C,\overline x}$ (with the notations of Lem. \ref{iotaexact}, $K^{sh}=(K^{sep})^{I_{v_{sep}}}$), $K'/K$ runs through finite subextensions of $K^{sh}/K$ and $K'^h=(K^{sep})^{D_{v_{sep}}}$ is the fraction field of the corresponding henselian ring. On the other hand, we have $$\begin{array}{rcl}(\iota_{C_v,*}R^qj_{v,*}\epsilon_*A_{|U_v})_{\overline x}&\simeq &H^q(K_v^{ur},\epsilon_*A_{|K_v^{ur}})\\ & \simeq & \limi_{K'}H^q(K'_{v'},\epsilon_*A_{|K'_{v'}})\end{array}$$
where $K_v^{ur}/K_v$ denotes the maximal unramified subextension of $K_v^{sep}/K_v$.
Here, the first isomorphism follows from Case 2 in Lem. \ref{iotaexact} (note that $(R^qj_{v,*})_{\overline \eta_v}=0$) and the other one follows from \cite{SGA4-II} VI, Prop. 5.7. The morphism (\ref{compl2}) is thus the direct limit of the morphisms
$$H^q(K'^{h},\epsilon_*A_{|K'^{h}})\to H^q(K'_{v'},\epsilon_*A_{|K'_{v'}})$$
for $K'$ as above. We may now conclude using \cite{Mi2} I, Rem. 3.10 (a) if $q=1$ and  III, Thm. 6.10 and Rem. 6.13 of \emph{loc. cit.} if $q\ge 2$ (in that case both sides are in fact trivial), recalling that  the \'etale cohomology of $A_{|K'^{h}}$ (resp. $A_{|K'_{v'}}$) in degree $\ge 1$ identifies with the direct limit of the flat cohomology of its torsion points.
\begin{flushright}$\square$\end{flushright}

\subsection{Rigid uniformization around semi-stable fibers} \label{subsectionrigidunif} \label{ruapoz}

\para For the purpose of this section, let $R$ denote a complete discrete valuation ring, $t$ a uniformizer and $k=R/(t)$ (resp. $K=R[{1\over t}]$) its residual (resp. fraction) field. As explained in \cite{Be2} 0, one has a diagram of functors \begin{eqnarray}\label{diagrigfor}\xymatrix{Sch_{lft}/R\ar[d]_{(-)_{|K}}\ar[r]^{(-)^{for}}&For/R\ar[d]^{(-)_{|K}}\\
Sch_{lft}/K\ar[r]^{(-)^{an}}&Rig/K}\end{eqnarray} where we have used the following notations:

- $Sch_{lft}/R$ (resp. $Sch_{lft}/K$) denotes the category of $Spec(R)$-schemes (resp. $Spec(K)$-schemes) which are locally of finite type,

- $For/R$ denotes the category of $t$-adic formal schemes over $Spf(R)$ which are locally of finite type over $Spf(R)$, and $(-)^{for}$ is the
functor of ($t$-adic) completion defined in \cite{EGA1}, I Sect. 10.8.

- $Rig/K$ denotes the category of rigid analytic spaces over the $t$-adic field $K$ and the functors $(-)_{|K}$, $(-)^{an}$ are respectively defined
in \cite{Be2} Sect. 0.2 and Sect. 0.3.  \petit

\noindent This diagram is not commutative. Instead, for $X/R$ in $Sch_{lft}/R$, there is a functorial morphism $(X^{for})_{|K}\rightarrow (X_{|K})^{an}$ which is
an open immersion if $X/R$ is separated and an isomorphism if $X/R$ is proper [\emph{loc. cit.}, 0.3.5]. \petit

\para \label{Raygp}
Start with a semi-stable Abelian variety $A_K/K$. Let $A/R$ denote the N\'eron model of $A_K/K$ and $A^0/R$ its connected component. There is an essentially unique pair $(G,e)$ where $G/R$ is a smooth algebraic group scheme which
is the extension of a torus $T/R$ by an Abelian scheme $B/R$ and \begin{eqnarray}\label{isoe}e: G^{for}\simeq (A^0)^{for}\end{eqnarray} is an
isomorphism of groups in $For/R$ (note that $e$ is not
algebraic in general). We refer to \cite{SGA7-I} IX, Sect. 7 for the construction of $(G,e)$ and simply refer to $G/R$ as the \emph{Raynaud group} attached to
$A_K/K$.

The next result is originally due to \cite{Ra3}. Here, we slightly reformulate the statement given in \cite{BL} where the reader is referred for a
proof.

\begin{prop} \label{rayunif} \emph{(Raynaud's uniformization)}

\debrom
\item \label{rayunifi} There exists a unique arrow $e'$ making the following diagram commutative:
\begin{eqnarray}\label{diage}\xymatrix{(G^{for})_{|K}\ar[r]^-{e_{|K}}\ar[d]& (A^{0,for})_{|K}\ar[d]\\
(G_{|K})^{an}\ar[r]^{e'}&(A_K)^{an}}\end{eqnarray}

\item \label{rayunifii} The kernel of $e'$ computed in the category of groups in $Rig/K$ is of the form $\smash{\G_{|K}^{an}}$ for some group $\G$ in $Sch_{lft}/R$ which is \'etale locally isomorphic to a constant free Abelian group of rank $dim_K \smash{T_{|K}}$.

\item \label{rayunifiii} If the torus $T$ is split then $\G$ is constant and $e'$ induces a surjection $$\xymatrix{(G_{|K})^{an}(S)\ar @{->>}[r]&(A_K)^{an}(S)}$$ for any $S$ in $Rig/K$ satisfying $H^1(S,\Z)=0$.
\finrom
\end{prop}
\begin{flushright}$\square$\end{flushright}


\begin{cor}  \label{rig}
\debrom
\item \label{rigi} Let $Fin/R$ denote the category of finite $Spec(R)$-schemes. The isomorphism $e$ of (\ref{isoe}) induces an isomorphism of Abelian presheaves on $Fin/R$:
$$\xymatrix{G\ar[r]^{e}_{\sim}& A^0}$$

\item \label{rigii}
Let $Fin/K$ denote the category of finite $Spec(K)$-schemes endowed with any topology which is finer than (or equal to) the \'etale topology. The
morphism $e'$ of Prop. \ref{rayunif} \ref{rayunifi} induces an exact sequence of Abelian sheaves on $Fin/K$:
$$\xymatrix{0\ar[r]&\G_{|K}\ar[r]&G_{|K}\ar[r]^{e'}& A_K\ar[r]& 0}$$
\item \label{rigiii} The morphisms \ref{rigi} and \ref{rigii} are compatible, ie. the following square   \begin{eqnarray}\label{compforrig}\xymatrix{G(S)\ar[r]^{e}\ar[d]& A^0(S)\ar[d]\\ G_{|K}(S_{|K})\ar[r]^{e'}& A_K(S_{|K})}\end{eqnarray}
commutes for any $S$ in $Fin/R$. \finrom
\end{cor}

Proof. \ref{rigi}  it suffices to notice that for any $H$ in $Sch_{lft}/R$ and $S$ in $Fin/R$ the functor $(-)^{for}$ identifies $H(S)$ with
$H^{for}(S^{for})$ (just because finite $R$-algebras are $t$-adically complete).

\ref{rigii} Similarly for any $H$ in $Sch_{lft}/K$ and $S$ in $Fin/K$ the functor $(-)^{an}$ identifies $H(S)$ with $H^{an}(S^{an})$ (because finite
$K$-algebras are complete for the $t$-adic topology).

\ref{rigiii} is a straightforward consequence of the definition of the analytification functor.
\begin{flushright}$\square$\end{flushright}

\para \label{paradefqff} In order to collect the sheaf theoretic consequences of Cor. \ref{rig}, it will be useful to introduce some intermediary sites. We use the general conventions of Sect. \ref{usualtop} for big and small pretopologies.

\begin{defn} \label{defqff} Let $X$ be any scheme. We will consider the small (resp. big) $top$ pretopology of $X$, $top(X)$ (resp. $TOP(X)$) and the associated topos $X_{top}$ (resp. $X_{TOP}$) for $top=qff$, $ff$, or $fet$ defined as follows:

- $qff$ stands for \emph{locally of finite presentation locally quasi-finite and flat},

- $ff$ stands for \emph{finite and flat},

- $fet$ stands for \emph{finite and \'etale}.
\end{defn}

The various inclusions of sites viewed as premorphisms of pretopologies give rise to weak morphisms of topoi fitting in the following pseudo-commutative diagram.
\begin{eqnarray}\label{epsalbet}\xymatrix{X_{FL}\ar[rd]_\gamma\ar[r]^\alpha\ar @<+2pt> `u[r] `[rr]^-{\epsilon} [rr]&X_{qff}\ar[r]^\beta\ar[d]&X_{et}\ar[d]\\ & X_{ff}\ar[r]&X_{fet}}\end{eqnarray}
\noindent It might be useful to point out that the topologies generated by $QFF(X)$ and $FL(X)$ coincide, ie. $X_{QFF}=X_{FL}$. The morphism $\alpha$
is thus nothing but the projection morphism $\pi:X_{FL}\to X_{fl}$. In particular $\alpha_*$ commutes to arbitrary limits. Here is an alternative description of
the topos $X_{qff}$ which will be useful for our purpose.

\begin{lem}\label{lemqff-top} Consider a  locally Noetherian  scheme $X$ and let $top$ stand for \emph{separated, quasi-finite and flat}. Then $X_{qff}$ is equivalent to $X_{top}$.
\end{lem}

\noindent Proof. By \cite{Ar} Lem. 3.1.3 , we are reduced to check that the following properties: 

-  Consider a covering $(U_i\to U)_i$ of the $qff$ pretopology. If $U$ and the $U_i$'s are in $top(X)$, 
then $(U_i\to U)_i$ is a covering for the $top$ pretopology.

- For any $U$ in $qff(X)$, there exists a covering of the $qff$ topology $(U_i\to U)_\alpha$ with $U_i \in top(X)$.


Observe that if $X'$ is a locally Noetherian scheme, then a morphism with target $X'$ is $top$ if and only if it is $qff$, separated and quasi-compact (\cite{EGA1}, Chap. 1, Prop. 6.6.3). The first property follows, by \cite{EGA1} Prop. 5.5.1 (v) and Prop. 6.6.4 (v). The second property follows as well, once noticed that morphisms between affine schemes, as well as open immersions into a locally noetherian scheme, are always separated and quasi-compact (\cite{EGA1}, Chap. 1, Prop. 6.6.4 (i)). 
\begin{flushright}$\square$\end{flushright}


Let us finally point out the equivalences $X_{qff}\simeq X_{ff}$ and $X_{fet}\simeq X_{et}$ in the particular case where $X$ is the spectrum of a
field. \petit

Let us go back to our situation, using the notation $X=Spec(R)$ and
\begin{eqnarray}\label{ijloc}\xymatrix{x=Spec(k)\ar[r]^-s&X&\ar[l]_-j\eta=Spec(K)}\end{eqnarray}

\begin{lem} \label{lemqff} \debrom \item \label{lemqffi} The continuous functor $ff(X)\to qff(X)$ (resp. $fet(X)\to et(X)$) inducing the left (resp. right) vertical premorphism in (\ref{epsalbet}) is cocontinuous as well. Whence a morphism of topoi in the opposite direction $$\begin{array}{cc}\iota:X_{ff}\to X_{qff}&\hspace{1cm} \hbox{(resp. $\iota:X_{fet}\to X_{et}$)}\end{array}$$
\item \label{lemqffii} The morphism $\iota$ defined in \ref{lemqffi} is a closed immersions of topoi. The complementary open immersion is the morphism denoted $j$ in the top (resp. bottom) line of the following  pseudo-commutative diagram of topoi \begin{eqnarray}\label{diagtopqff}\xymatrix{&X_{ff}\ar[d]^-\beta\ar[r]^-{\iota}&X_{qff}\ar[d]^\beta&\ar[d]^\beta
    \ar[l]_-j\eta_{qff}\\
    x_{et}\ar[r]^-\sim & X_{fet}\ar[r]^-\iota&X_{et}&\ar[l]_-j\eta_{et}}
    \end{eqnarray}
\item \label{lemqffiii} Smooth group schemes over $X$ are acyclic for the functor $\beta_*$.
\finrom
\end{lem}
Proof. \ref{lemqffi} In view of \ref{lemqff-top}, it is sufficient to prove that the inclusion of  $ff(X)$ into the category of separated quasi-finite flat $X$-schemes admits a right adjoint which is continuous. If $U/X$ is separated and quasi-finite, there is a disjoint decomposition $U=U^{fin}\sqcup U_2$ where $U^{fin}/X$ is finite
and $x\times_X U_2=\emptyset$ (\cite{Mi1} I, Thm. 4.2 $(c)$). Note that if $U'/X$ is finite, then $Hom_X(U',U_2)=\emptyset$ (since
$Hom_X(U',\eta)=\emptyset$). It follows in particular that $U^{fin}$ varies functorially with respect to $U/X$. Whence the desired functor $qff(X)\to ff(X)$. This functor is easily seen to be continuous (use that the finite part is a direct sum of spectrums of local rings by \cite{Mi1} I, Thm. 4.2 $(c)$) as well as right adjoint to the inclusion.

\ref{lemqffii} Let us explain the first line of (\ref{diagtopqff}) since the second is similar (and well known). The fact that $\iota$ is an
immersion follows from the fact that the inclusion $ff(X)\to qff(X)$ is fully faithful. According to \cite{SGA4-I} IV, Prop. 9.3.4, we have to show that the
following properties $(a)$, $(b)$ are equivalent for any sheaf $F$ on $qff(X)$: $(a)$ $F\to \iota_*\iota^{-1}F$ is an isomorphism, $(b)$ $j^{-1}F$ is a final object. Now this
equivalence is clear from the fact that $\iota_*\iota^{-1}F(U)=F(U^{fin})$.

\ref{lemqffiii} We already know that smooth group schemes are acyclic for $\epsilon_*$ (Lem. \ref{acycsm}) and that $\alpha_*$ is exact. The result
follows formally.


\begin{flushright}$\square$\end{flushright}

In view of the next proposition, let us emphasize that Lem. \ref{lemqff} \ref{lemqffii}  implies an equivalence between $Ab(X_{qff})$ (resp. $Ab(X_{et})$) and the
category of triples $(F_1,F_2,f)$ where $F_1\in Ab(X_{ff})$ (resp. $Ab(X_{fet})$), $F_2\in Ab(\eta_{qff})$ (resp. $Ab(\eta_{et})$) and $f:F_1\to
\iota^{-1}j_*F_2$ (\cite{SGA4-I} IV, Thm. 9.5.4). These equivalences (and thus also the functors $j_!$ of extension by zero) are clearly compatible via the
restriction functor $\beta_*$.

\para Let us come back to the d\'evissage of $A^0_{p^\infty}$. We begin in the topos $qff$.

\begin{cor} \label{rigqff}

\debrom
\item   \label{rigqffi} The projections of Cor. \ref{rig} \ref{rigi} and \ref{rigii} to $X_{ff}$ and $\eta_{ff}\simeq \eta_{qff}$ glue into an exact sequence as follows in $Ab(X_{qff})$: $$\xymatrix{0\ar[r]&j_!\alpha_*\G_{|\eta}\ar[r]&\alpha_*G\ar[r]&\alpha_*A^0\ar[r]&0}$$
\item  \label{rigqffii} Applying $\Z/p^k\otimes^L_\Z(-)$ to \ref{rigqffi} gives an exact sequence in $Mod(X_{qff},\Z/p^k)$:
$$\xymatrix{0\ar[r]&\alpha_*G_{p^k}\ar[r]&\alpha_*A^0_{p^k}\ar[r]&j_!\alpha_*(\G_{|\eta}/p^k)\ar[r]&0}$$
\finrom
\end{cor}
Proof. \ref{rigqffi} We use the interpretation of $Ab(X_{qff})$ in terms of triples recalled after Lem. \ref{lemqff}. The sheaf $\alpha_*G$ corresponds to the triple
$(\iota^{-1}\alpha_*G,j^{-1}\alpha^*G,nat:\iota^{-1}\alpha_*G\to \iota^{-1}j_*j^{-1}\alpha^*G)$. Explicitly for $S/X$ in $ff(X)$ (resp. $S/\eta$ in
$ff(\eta)$)  we have $\iota^{-1}\alpha_*G(S)=G(S)$ (resp.  $j^{-1}\alpha_*G(S)=G_{|K}(S)$), $\iota^{-1}j_*j^{-1}\alpha^*G(S)=G_{|K}(S_{|K})$ and $nat$
is the left vertical arrow in Cor. \ref{rig} \ref{rigiii}. A similar interpretation holds for $A^0$. The desired exact sequence is thus nothing
but a reformulation of Cor. \ref{rig} \ref{rigi}, \ref{rigii}, \ref{rigiii}.

\ref{rigqffii} Exactness on the right follows from the fact that multiplication by $p^k$ is epimorphic on $G$ (Cor. \ref{corpdivtop} \ref{corpdivtopi}) and thus on $\alpha_*G$ too.
\begin{flushright}$\square$\end{flushright}

We may now go back to the big flat site.

\begin{lem} \label{remrigqff}
\debrom
\item \label{remrigqffi} The adjunction morphism $\alpha^{-1}\alpha_*\G_{|\eta}\to \G_{|\eta}$ is invertible. In particular, the first arrow of Cor. \ref{rigqff} \ref{rigqffi} (or equivalently Cor. \ref{rig} \ref{rigii})  defines a morphism \begin{eqnarray} \label{rigetaFL} \G_{|\eta}\to G_{|\eta} \hbox{\hspace{.5cm} in $Ab(\eta_{FL})$}\end{eqnarray}

\item \label{remrigqffii} Similarly, the first arrow of Cor. \ref{rigqff} \ref{rigqffii} defines a morphism \begin{eqnarray} \label{rigXFL} G_{p^.}\to A^0_{p^.}\hbox{\hspace{.5cm} in $Ab(X_{FL}^\N)$}\end{eqnarray}

\finrom
\end{lem}
Proof. \ref{remrigqffi} (resp. \ref{remrigqffii}) follows from the fact that $\G_{|\eta}$ (resp. $G_{p^k}$) is representable by a group of $qff(\eta)$ (resp. $qff(X)$, see Lem. \ref{pdivfunct} \ref{pdivfunctii}).
\begin{flushright}$\square$\end{flushright}

\begin{rem}
\label{remrigqffiii} We may view the exact sequence of Cor. \ref{rig} \ref{rigii} as a quasi-isomorphism  $[\G_{|\eta}\to G_{|\eta}]\to A_{|\eta}$ in the category of complexes of Abelian groups of $\eta_{qff}$, the first complex being placed in degree $[-1,0]$. Applying $\Z/p^.\otimes^L(-)$ and pulling back to the big flat topos, yields an isomorphism of exact sequences in $Ab(\eta_{FL}^\N)$:  \begin{eqnarray}\xymatrix{0\ar[r]&G_{|\eta,p^.}\ar[d]^\wr\ar[r]&\Z/p^.\Ltens [\G_{|\eta}\to G_{|\eta}]\ar[r]\ar[d]^\wr&\G_{|\eta}/p^.\ar[d]^\wr\ar[r]&0\\
    0\ar[r]&G_{|\eta,p^.}\ar[r]&A_{|\eta,p^.}\ar[r]&\G_{|\eta}/p^.\ar[r]&0}
    \end{eqnarray}
    where the bottom left arrow is compatible with (\ref{rigXFL}).
\end{rem}

\begin{prop} \label{riget} The image of the morphism (\ref{rigXFL}) under the functor  $$R\underline\G^x(X,R\epsilon_*\rho^{-1}):D^+(X_{FL},\Z/p^.)\to D^+(X_{et},\Z/p^.)$$ fits into a canonical distinguished triangle $$\xymatrix{R\underline \G^x(X,R\epsilon_*G_{|X^+,p^.})\ar[r]&R\underline \G^x(X,R\epsilon_*A^0_{|X^+,p^.})\ar[r]&j_!\epsilon_*(\G_{|\eta}/p^.)\ar[r]^-{+1}&}$$
\end{prop}
Proof.
%
The restriction of the morphism defined in Lem. \ref{remrigqff} \ref{remrigqffii} to $FL(x)$ (resp. $qff(X)$) clearly coincides with the natural isomorphism of group schemes $\smash{G_{|x,p^k}\simeq A^0_{|x,p^k}}$ induced by (\ref{isoe}) (resp. coincides with the first arrow of Cor. \ref{rigqff} \ref{rigqffii}). If $X^+$ denote the diagram $(x\to X)$ as in Def. \ref{defitype}, the exact sequence of Cor. \ref{rigqff} \ref{rigqffi} can thus be completed into an exact sequence (with notations similar to Sect. \ref{diagopi}) $$\xymatrix{0\ar[r]&(\alpha_*G_{|x,p^.}\leftarrow
\alpha_*G_{p^.})\ar[r]&(\alpha_*A^0_{|x,p^.}\leftarrow \alpha_*A^0_{p^.})\ar[r]&(0\leftarrow j_!\alpha_*\G_{|\eta}/p^k)\ar[r]&0}$$
over  $(\smash{X_{qff}^{+,\N}},\Z/p^.)$.
Apply now the functor $R\beta_*:D^+(X^{+,\N}_{qff},\Z/p^.)\to D^+(X^{+,\N}_{et},\Z/p^.)$ and then   $R\underline \G^x(X,-):D^+(X^{+,\N}_{et},\Z/p^.)\to
D^+(X^{\N}_{et},\Z/p^.)$  to get a distinguished triangle $$\xymatrix{R\underline \G^x(X,R\beta_*\alpha_*G_{|X^+,p^.})\ar[r]&R\underline
\G^x(X,R\beta_*\alpha_*A^0_{|X^+,p^.})\ar[r]&R\beta_*j_!\alpha_*\G_{|\eta}/p^.\ar[r]^-{+1}&}$$ (use e.g.  Lem. \ref{soritesvan} \ref{soritesvani} for the third term). The announced distinguished triangle follows since  $\alpha_*$ is exact and $j_!\alpha_*\G_{|\eta}/p^.$ is $\beta_*$ acyclic (use Lem. \ref{acycsm} to check this).
\begin{flushright}$\square$\end{flushright}


\subsection{D\'evissage in $p$-divisible groups} \label{dipdg}

\para We are now in a position to define a $p$-divisible group $H$ over the diagram $J^+$ which will serve as a replacement for $A_{|J^+,p^\infty}$. By a $p$-divisible group $H$ over $J^+$, we mean an object of $pdiv(J^+)$, where $pdiv$ denotes the cofibered category of $p$-divisible groups over $\cal Sch^{op}$ (here we use the conventions of Sect. \ref{fibcofdiag}); in other terms, $H$ is a group in $J^+_{FL}$, whose components are $p$-divisible groups.


\begin{defn} \label{defH} Let $A/C$ be a semi-Abelian scheme whose restriction to $U$ is Abelian. We define a $p$-divisible group $H$ over $J^+=(Z_v\to C_v\leftarrow U_v\to U)$: $$H:=(G_{v|Z_v,p^\infty}\leftarrow  G_{v,p^\infty}\to A_{|U_v,p^\infty}\leftarrow A_{|U,p^\infty})$$ where $G_v/C_v$ is the Raynaud group attached to $A_{|U_v}/U_v$ as in Sect. \ref{Raygp} and $j_v^{-1}G_{v,p^\infty}\to A_{|U_v,p^\infty}$ is the arrow induced by Lem. \ref{remrigqff} \ref{remrigqffii}. \end{defn}

The following result is the final stage of our d\'evissage. It states that $R\epsilon_*H_{p^.}$ is sufficient to retrieve the projection of $R\underline \G^Z(A_{p^.})$ to the small \'etale topos of $C$.

\begin{prop} \label{devissageFL} Recall the functor $Sma:D(J^\N_{et},\Z/p^.)\to D(J^\N_{et},\Z/p^.)$ from Def. \ref{defsma}. There is a canonical isomorphism $$R\epsilon_*R\underline\G^ZA_{p^.}\simeq Rm_*Sma\, R\underline \G^{Z_J}(J,R\epsilon_*H_{p^.})$$ \end{prop}
Proof. According to Prop. \ref{MVfl}, it suffices to build an isomorphism between the right hand side and $Rm_*R\underline \G^{Z_J}(J,R\epsilon_*
A^0_{|J^+,p^.})$ in $D(C_{et}^\N,\Z/p^.)$. Consider the morphism $H_{p^.}\to A^0_{|J^+,p^.}$ in $Mod(J^{+,\N}_{FL},\Z/p^.)$,  which is induced by Lem. \ref{remrigqff} \ref{remrigqffii} on the vertices $Z_v$  and $C_v$ and which is the identity on the vertices $U$ and $U_v$. This, and the natural transformation $Sma\to id$ (Lem. \ref{lemsma} \ref{lemsmaii}) induce morphisms:
\begin{eqnarray}\label{glufl1}Sma\, R\underline \G^{Z_J}(J,R\epsilon_*H_{p^.})&\to & Sma\, R\underline \G^{Z_J}(J,R\epsilon_*A^0_{|J^+,p^.})\\
\label{glufl2}&\to &  R\underline \G^{Z_J}(J,R\epsilon_*A^0_{|J^+,p^.})\end{eqnarray}
in $D(J_{et}^\N,\Z/p^.)$. We will prove that (\ref{glufl1}) and (\ref{glufl2}) are invertible.

In order to prove that (\ref{glufl1}) is an isomorphism, it suffices to check the conditions of Cor. \ref{corsma} for the morphism \begin{eqnarray}\label{glufl3}R\underline \G^{Z_J}(J,R\epsilon_*H_{p^.})&\to & R\underline \G^{Z_J}(J,R\epsilon_*A^0_{|J^+,p^.})\end{eqnarray}
Let us examine the component $U$. By Lem. \ref{resvan} \ref{resvanii}, there is a functorial isomorphism $$(R\underline\G^{Z_J}(J,M))_U\simeq R\underline\G^\emptyset(U,M_{|U^+})$$ for $M$ in $D^+(\smash{J^{+,\N}_{et}},\Z/p^.)$. Here, $\emptyset$ denote the diagram of empty type and the right hand side is thus naturally isomorphic to $M_{U}$. It follows that the morphism $(\ref{glufl3})_{U}$ identifies with $\smash{(R\epsilon_*H_{p^.})_{U}}\to \smash{(R\epsilon_*A^0_{p^.,|J^+})_{U}}$. It is an isomorphism, since $R\epsilon_*$ can be computed componentwise on $J^+$. The same arguments show that $\smash{(\ref{glufl3})_{U_v}}$ is an isomorphism. Let us now examine the component $C_v$. By Lem. \ref{resvan} \ref{resvanii} again, we have $$(R\underline \G^{Z_J}(J,M))_{C_v}\simeq R\underline \G^{Z_v}(C_v,M_{|C_v^+})$$
The image of the morphism $(\ref{glufl3})_{C_v}$ by $z_v^{-1}$ thus identifies with $$z_v^{-1}R\underline \G^{Z_v}(C_v,R\epsilon_*H_{p^.,|C_v^+})\to z_v^{-1}R\underline \G^{Z_v}(C_v,R\epsilon_*A^0_{p^.,|C_v^+})$$
and is an isomorphism by Prop. \ref{riget}.

In order to prove that (\ref{glufl2}) is an isomorphism, it suffice to check the conditions of Lem. \ref{lemsma} \ref{lemsmaii}, ie. that $$j_v^{-1}(R\G^{Z_J}(J,R\epsilon_*A^0_{|J^+,p^.}))_{C_v}\to  (R\G^{Z_J}(J,R\epsilon_*A^0_{|J^+,p^.}))_{U_v}$$ is an isomorphism. This, in turn, follows immediately from  Lem. \ref{resvan} \ref{resvanii},  Lem. \ref{soritesvan} \ref{soritesvani} and the fact that $j_v^{-1}R\epsilon_*A^0_{|C_v,p^.}\to R\epsilon_*A^0_{|U_v,p^.}$ is invertible (because $j_v$ is \'etale).
\begin{flushright}$\square$\end{flushright}


\section{D\'evissage of twisted syntomic complexes} \label{secdotsc}


In this section, we prove two d\'evissage properties of the twisted syntomic complexes defined in Def. \ref{defSonephi} and Prop. \ref{indepSet} \ref{indepSeti}.

\subsection{Complete Mayer-Vietoris}
\label{cmv} ~~

Recall the natural morphism of diagrams $m^\s:J^\s\to C^\s$ from (\ref{diagdiag}). We have the following diagram of ringed topoi
\begin{eqnarray}\label{gludiagcrys}\xymatrix{((J^\s/\Sigma_.)_{crys,et},\O)\ar[r]^{m^\s}\ar[d]^{u}& ((C^\s/\Sigma_.)_{crys,et},\O)\ar[d]^{u}\\
(J^\N_{et},\widetilde\O^{crys}_.)\ar[r]^{m}& (C^\N_{et},\widetilde\O^{crys}_.)}\end{eqnarray} Here, we have identified the small
\'etale site of a log scheme with the small \'etale site of the underlying scheme. The following proposition expresses that Dieudonn\'e crystals and their
syntomic complexes glue from $J^\s$ to $C^\s$.

\begin{prop} \label{MVsyn} Let us denote $\Delta_J$ the type of the diagram  $J^\s$.

\debrom \item \label{MVsyni} The  pullback functor  $$m^{\s,*}:Mod((C^\s/\Sigma_\infty)_{crys,et},\O)\to Mod((J^\s/\Sigma_\infty)_{crys,et},\O)$$  induces  equivalences of
categories:
$$\begin{array}{c}\xymatrix{\cal Crys_{lfft}((C^\s/\Sigma_\infty)_{crys,et},\O)\ar[r]^-{m^{\s,*}}_-\sim& \G_{cart}(\Delta_J^{op}\times_{J^\s,(\cal Sch^\s/\Sigma_1)^{op}}\cal Crys_{lfft}((-/\Sigma_\infty)_{crys,et},\O)/\Delta_J^{op})} \\
\xymatrix{\DC(C^\s)\ar[r]^-{m^{\s,*}}_-\sim& \G_{cart}(\Delta_J^{\s,op}\times_{J^\s,(\cal Sch^\s/\Sigma_1)^{op}}\DC/\Delta_J^{op})}\end{array}$$
where $\G_{cart}$ denotes the full subcategory of  (\ref{notsectioncof}) formed by cartesian sections.
\item \label{MVsynii} For $M_.\in Mod((C^\s/\Sigma_.)_{crys,et},\O)$,  there is a canonical morphism
$$\xymatrix{R\smash{u}_*M_.\ar[r]&  R\smash{m}_*R\smash{u}_*m^{\s,*}M_.&\hbox{in $\smash{D(C_{et}^\N,\O_.^{crys})}$}}$$
Assume now that $M_.$ is a crystal. Let us choose a morphism $\smash{m^\s_{[.],Y}}$ in  $\smash{HR_F^{\s,et}}$  above $m^\s$ (Lem. \ref{lememb1} \ref{lememb1iv} and \ref{lememb2} \ref{lememb2ii}, \ref{lememb2iii}, \ref{lememb2iv}) and denote
$\smash{m^\s_{T,[.]}}:\smash{T_{J^\s,[.]}\to T_{C^\s,[.]}}$ its logarithmic divided power envelope.  The previous morphism of
$\smash{D(C_{et}^\N,\O_.^{crys})}$ identifies with the natural morphism
$$\xymatrix{Rf_{T_{C^\s,[.],.},*}\Omega^\bullet_{T_{C^\s,[.],.}}(M_.)\ar[r]&
Rm_*Rf_{T_{J^\s,[.],.},*}\Omega^\bullet_{T_{J^\s,[.],.}}(m^{\s,*}M_.)}$$
It is an isomorphism if $M_.$ is a locally free crystal of finite type.

\item \label{MVsyniii} For $D\in \DC(C^\s)$ and $h\in \G(C,M_C/\Gm)$, the canonical base change morphism induces an isomorphism
$$\xymatrix{\S^{1,\varphi}_{et,.,C^\s}(-h)(D)\simeq Rm_*\S^{1,\varphi}_{et,.,J^\s}(-h_{|J})(D_{|J^\s})&\hbox{in $\smash{D(Mod^{1,\varphi}(C^\N_{et},\widetilde\O_.^{crys}))}$}}$$ as well as a similar isomorphism without the superscripts ${}^{1,\varphi}$ in $\smash{D(C^\N_{et},\widetilde\O_.^{crys,F=1})}$.
\finrom
\end{prop}

Proof. \ref{MVsyni} Thanks to $m^*$ being compatible with $F^*$, it suffices to prove that the first functor is fully faithful and essentially
surjective. We will do this by Zariski descent along the coskeleton of an appropriate covering.  Let  $\smash{(C_\lambda^\s\rightarrow C^\s)}$ be an open covering such that each $\smash{C_\lambda^\s}$ is affine, has
$p$-bases of the form $(\emptyset,t)$ and such that each point $v$ of $Z$ is contained in exactly one of the $C_\lambda^\s$'s denoted $C_{\lambda_v}^\s$. Consider $\smash{U_{\lambda}:=C_\lambda\times_{C}U}$,
$\smash{C_{v,\lambda}^\s:=C_\lambda^\s\times_{C^\s}C_v^\s}$, $\smash{U_{v,\lambda}:=C_\lambda\times_{C}U_v}$,
$\smash{J_{\lambda}^\s:=C_\lambda^\s\times_{C^\s}J^\s}$ and form  $\smash{C^\s_{[0]}}:=\smash{\sqcup C^\s_\lambda}$,  $\smash{U_{[0]}}:=\smash{\sqcup
U_\lambda}$, $\smash{C^\s_{v,[0]}}:=\smash{\sqcup C^\s_{v,\lambda}}$, $\smash{U_{v,[0]}}:=\smash{\sqcup U_{v,\lambda}}$ as well as the diagram of type
$\Delta_J$ represented by $(\smash{C^\s_{v,[0]}}\leftarrow \smash{U_{v,[0]}}\to \smash{U_{[0]}})$. Let us furthermore denote
$\smash{C^{\s}_{[\nu]}}$,  $\smash{U_{[\nu]}}$, $\smash{C^{\s}_{v,[\nu]}}$, $\smash{U_{v,[\nu]}}$ and $\smash{J^\s_{[\nu]}}$ the respective
$\nu+1$-th fibered power over $C^\s$, $U$, $C_v^\s$, $U_v$ and $J^\s$ of $\smash{C^{\s}_{[0]}}$,  $\smash{U_{[0]}}$, $\smash{C^{\s}_{v,[0]}}$,
$\smash{U_{v,[0]}}$ and $\smash{J^\s_{[0]}}$. By Zariski descent along $\smash{C^{\s}_{[0]}}\rightarrow C^\s$ and $\smash{J^{\s}_{[0]}}\rightarrow
J^\s$, it is enough to prove that the induced functor $$\xymatrix{\cal Crys_{lfft}((C^{\s}_{[\nu]}/\Sigma_\infty)_{crys,et},\O)\ar[r]^-{m^{\s,*}_{[\nu]}}&
\G_{cart}(\Delta_J^{op}\times_{J^{\s}_{[\nu]},(\cal Sch^\s/\Sigma_1)^{op}}Crys_{lfft}((-/\Sigma_\infty)_{crys,et},\O)/\Delta_J^{op})}$$ is an equivalence, ie that for each $\nu\ge 0$
and multi-index  $\underline \lambda=(\lambda_0,\dots, \lambda_\nu)$ the induced functor
$$\xymatrix{\cal Crys_{lfft}((C_{[\nu],\underline \lambda}^{\s}/\Sigma_\infty)_{crys,et},\O)\ar[r]^-{m^{\s,*}_{[\nu],\underline \lambda}}& \G_{cart}(\Delta_J^{op}\times_{J^\s_{[\nu],\underline \lambda},(\cal Sch^\s/\Sigma_1)^{op}}Crys_{lfft}((-/\Sigma_\infty)_{crys,et},\O)/\Delta_J^{op})}$$ is an equivalence where
$\smash{C_{[\nu],\underline \lambda}^{\s}}:=\smash{C_{\lambda_0}^\s\times_{C^\s}\dots \times_{C^\s}C_{\lambda_\nu}^\s}$ and
$\smash{J^{\s}_{[\nu],\underline \lambda}}:=\smash{J^\s_{\underline \lambda}=J_{\lambda_0}^\s\times_{J^\s}\dots \times_{J^\s}J_{\lambda_\nu}^\s}$.

Since $U_{v,\lambda}=C_{v,\lambda}$ for $\lambda\ne \lambda_v$, the equivalence is trivial if  $\lambda_i\ne \lambda_j$ for some $i\ne j$ (in that case $\smash{C^\s_{[\nu],\underline
\lambda}}=\smash{U_{[\nu],\underline \lambda}}$ and  $\smash{C^\s_{v,[\nu],\underline \lambda}}=\smash{U_{v,[\nu],\underline \lambda}}$ for all $v$).
On the other hand, if  $\underline \lambda=(\lambda,\dots,\lambda)$ then $\smash{C_{[\nu],\underline \lambda}^{\s}}$ and
$\smash{J^{\s}_{[\nu],\underline \lambda}}$ respectively coincide with $\smash{C_{\lambda}^\s}$ and $\smash{J_\lambda^\s}$. We are thus reduced to
the case $\nu=0$, and it suffices to treat the case $\lambda=\lambda_v$ for some $v$ in $Z$ since the other cases are trivial as well. In that case, 
let us choose a $p$-basis $(\emptyset,t)$ for $C_\lambda$ and denote $R=\G(C_\lambda,\O)$. The ring $R$ is thus an \'etale  $\Fp[T]$ algebra via
$T\mapsto t$ and that $t$ is a uniformizer at $v$. Denote furthermore $R_v$ the $t$-adic completion of $R$. Hence with these notations, we have
$\smash{C_\lambda^\s}=(Spec(R),1\mapsto t)$, $U_\lambda=Spec(R[t^{-1}])$, $\smash{C_{v,\lambda}^\s}=(Spec(R_v),1\mapsto t)$ and 
$\smash{U_{v,\lambda}=Spec(R_v[t^{-1}])}$. Let  $\smash{\tilde R_k}$, (resp. $\smash{\tilde R_{v,k}}$) denote the essentially unique \'etale
$\Z/p^{k}[t]$- (resp. $\Z/p^{k}[[t]]$-) algebra lifting $R$ (resp. $R_v$) and endow it with the log structure induced by $t$. We will show that the
functor
$$m_{[0],\lambda_v}^{\s,*}:\nabla\hbox{-}Mod_{lfft}(\tilde R_k,1\mapsto t)\rightarrow  \nabla\hbox{-}Mod_{lfft}(\tilde R_{k}[t^{-1}])\times_{\nabla\hbox{-}Mod_{lfft}(\tilde R_{v,k}[t^{-1}])}\nabla^\s\hbox{-}Mod_{lfft}(\tilde R_{v,k},1\mapsto t)$$ is fully faithful and essentially surjective. The proof will then be finished using
Prop. \ref{crystals} and Rem. \ref{lemlcrys} \ref{lemlcrysi} (the arguments below will show that the property of being normalized when $k$ varies is preserved by the equivalence).

As a preliminary remark, we observe that the following sequence is exact  thanks to the fact that $t$ does not divide zero in the Noetherian ring
$\smash{\tilde R_k}$:  \begin{eqnarray}\label{exR}\xymatrix{0\ar[r]&\tilde R_k\ar[r]&\tilde R_k[t^{-1}]\oplus \tilde R_{v,k}\ar[r]^-{(1,-1)}&\tilde
R_{v,k}[t^{-1}]\ar[r]&0}\end{eqnarray}

Consider the functor sending an object $$((M_1,\nabla_1),(M_2,\nabla_2),\alpha:\smash{\tilde R_{v,k}[t^{-1}]\otimes_{\tilde
R_{k}[t^{-1}]}(M_1,\nabla_1)\simeq \tilde R_{v,k}[t^{-1}]\otimes_{\tilde R_{v,k}}(M_2,\nabla_2)})$$ in the target of category of
$\smash{m_{[0],\lambda_v}^{\s,*}}$ to the $\smash{\tilde R_k}$ module $$M:=Ker(M_1\oplus M_2\mathop\rightarrow\limits^{(\alpha,-1)} M_2[t^{-1}])$$
endowed with the induced connection (which turns out to be well defined by exactness of  $$\xymatrix{0\ar[r]&{\Omega_{\tilde
R_k}}\ar[r]&{\Omega_{\tilde R_k[t^{-1}]}}\oplus{\Omega_{\tilde R_{v,k}}}\ar[r]^{(1,-1)}&\smash{\Omega_{\tilde R_{v,k}[t^{-1}]}}}).$$
This functor is clearly right adjoint to $\smash{m_{[0],\lambda_v}^{\s,*}}$. Note that this functor preserves the property of being normalized when $k$ varies, thanks to the fact that the arrow $M_1\oplus
M_2\to M_2[t^{-1}]$ is in fact onto and its target is flat over $\Z/p^k$.  Let us now check that it is in fact a quasi-inverse to
$\smash{m_{[0],\lambda_v}^{\s,*}}$.

- full faithfulness of $\smash{m_{[0],\lambda_v}^{\s,*}}$. It suffices to check that for a locally  $\smash{\tilde R_k}$-module $M$, the sequence \begin{eqnarray}\label{qsdfg}\xymatrix{0\ar[r]&M\ar[r]&M[t^{-1}] \oplus M_v\ar[r]^-{(1,-1)}& M_v[t^{-1}]}\end{eqnarray} is exact. This is indeed the case thanks to (\ref{exR}).  

- full faithfulness of the right adjoint. If $((M_1,\nabla_1),(M_2,\nabla_2),\alpha)$ is in the target category of $\smash{m_{[0],\lambda_v}^{\s,*}}$, 
the following natural morphisms are invertible: $$\diagram{\tilde R_k[t^{-1}]\otimes_{\tilde R_k} Ker(M_1\oplus
M_2\mathop\rightarrow\limits^{(\alpha,-1)} M_2[t^{-1}])&\hfl{(1,pr_1)}{}&M_1&&\hbox{in $Mod(\tilde R_{k}[t^{-1}])$}\cr \tilde R_{v,k}\otimes_{\tilde
R_k} Ker(M_1\oplus M_2\mathop\rightarrow\limits^{(\alpha,-1)} M_2[t^{-1}])&\hfl{(1,pr_2)}{}&M_2&&\hbox{in $Mod(\tilde R_{v,k})$}}$$ Indeed Nakayama's
lemma reduces us to the case $k=1$ and then, we may use the classification of finitely generated modules over the Dedekind rings $R$ and $R_{v}$ to
show that $(M_1,M_2,\alpha)$ is in fact necessarily of the form $(M[t^{-1}],M_v,can)$ for some $R$-module $M$. Localizing if necessary, we conclude using (\ref{qsdfg}).

\ref{MVsynii} The claimed morphism is simply induced by the natural morphism $M\to \smash{Rm^\s_*m^{\s,*}}M$.
If $M$ is a crystal, the claimed interpretation of the image of this morphism by $Ru_*$ in terms of de Rham complexes using hypercoverings follows easily from the proof of
Prop. \ref{cdOmega} \ref{cdOmegaii}.

Assume now that $M_.$  is a crystal of locally free modules of finite type. In order to prove the claimed isomorphism, we may use the de Rham
interpretation via global embeddings. In fact, we may even choose $\tilde m:\smash{\tilde J^\s}\to \smash{\tilde C^\s}$ in $Emb^{glob}$ above $m$ such that
$\smash{\tilde J^\s}$ and $\smash{\tilde C^\s}$ are respectively liftings of $C^\s$ and $J^\s$. Let us briefly explain this. First, note that since
$C$ is smooth of dimension $1$, the obstruction to the existence of a formally smooth $p$-adic lifting ${\tilde C}$ furnished by deformation theory
vanishes (however, in general it not possible to lift the Frobenius). Next, we may define $\smash{\tilde C^\s}$ by choosing an arbitrary
lifting of the log structure of $C^\s$. The lifting $\smash{\tilde J^\s}$ and the arrow $\tilde m$ are then obtained by relative perfectness of the edges of $J^\s$ over $C^\s$. It
remains to check that the natural morphism $$\Omega^\bullet_{{\tilde C^\s,.}}(M_.)\to R\tilde m_{*}\Omega^\bullet_{{\tilde J^\s,.}}(m^{\s,*}M_.)$$ is
invertible in $\smash{D(\tilde C^{\s,\N}_{et},\O^{crys}_.)}$. Since this question is \'etale local on ${\tilde C}$, we may replace $C^\s$ by
$\smash{C^\s_\lambda}$ as in \ref{MVsyni}. Since quasi-coherent modules are acyclic for the direct image of affine morphisms, we find (using Lem. \ref{defMV} \ref{defMVii}) that the statement  boils down to the exactness of the sequence
$$\hspace{-.5cm}\xymatrix{0\ar[r]&\Omega^\bullet_{{\tilde C_\lambda^\s,.}}(M_.)\ar[r]&\iota_{C_{v,\lambda},*}\Omega^\bullet_{{\tilde
C_{v,\lambda}^\s,.}}(\iota_{C_{v,\lambda}}^*M_.)\oplus j_{\lambda,*}\Omega^\bullet_{{\tilde
U_\lambda,.}}(j_\lambda^*M_.)\ar[r]&j_{\lambda,*}\iota_{U_{v,\lambda},*}\Omega^\bullet_{{\tilde U_{v,\lambda},.}}(\iota_{U_{v,\lambda}}^*j_\lambda^*M_.)\ar[r]&0}$$ where the
notations $j_\lambda$, $j_{v,\lambda}$, $\smash{\iota_{U_{v,\lambda}}}$, $\smash{\iota_{C_{v,\lambda}}}$ are used abusively to denote either the
morphisms induced by $j$, $j_v$, $\smash{\iota_{U_v}}$, $\smash{\iota_{C_v}}$ in the diagram $J_\lambda$ or its lifting $\smash{\tilde J_\lambda}$.
Now, the exactness in question is trivial if $|C_\lambda|\cap |Z|=\emptyset$ and results from (\ref{exR}) if $\lambda=\lambda_v$.

\ref{MVsyniii} It suffices to prove the  isomorphism in $\smash{D(Mod^{1,\varphi}(C_{et}^\N,\widetilde\O_.^{crys}))}$. Consider a morphism
$\smash{Y_{J^\s,[.]}\to Y_{C^\s,[.]}}$ above $m^\s$ in $HR^{\s,et}_{F}$ and denote  $\smash{T_{J^\s,[.]}\to T_{C^\s,[.]}}$ its logarithmic divided power
envelope. We have to show that the induced morphism $$R f_{T_{C^\s,[.],.},*}\cal
S\Omega^{\bullet,1,\varphi}_{et,.,T_{C^\s,[.]}}(-h_{|U_{C^\s,[.]}})(D_{|U_{C^\s,[.]}})\to Rm_* R f_{T_{J^\s,[.],.},*}\cal
S\Omega^{\bullet,1,\varphi}_{et,.,T_{J^\s},[.]}(-h_{|U_{J^\s,[.]}})(D_{|U_{J^\s,[.]}})$$ is invertible. We have the following distinguished triangles in
$D(Mod^{1,\varphi}(T_{C^\s,[.],.,et},\widetilde\O^{crys}_.))$ (see Lem. \ref{1phimod} \ref{1phimodiii} and   Lem. \ref{lemfil} \ref{lemfili}):
$$\begin{array}{ccc}
\xymatrix{(0,\widetilde{Fil}^0\Omega^\bullet_{.,T_{C^\s,[.]}}(-h_{|U_{C^\s,[.]}})(D_{|U_{C^\s,[.]}}),0,0)\ar[d]\\
\cal S\Omega^{\bullet,1,\varphi}_{et,.,T_{C^\s,[.]}}(-h_{|U_{C^\s,[.]}})(D_{|U_{C^\s,[.]}})\ar[d]\\
(\widetilde{Fil}^1\Omega^\bullet_{.,T_{C^\s,[.]}}(-h_{|U_{C^\s,[.]}})(D_{|U_{C^\s,[.]}}),0,0,0)\ar[d]^-{+1}\\
\hbox{$ $}}& \hbox{ and }&
\xymatrix{(\widetilde{Fil}^1\Omega^\bullet_{.,T_{C^\s,[.]}}(-h_{|U_{C^\s,[.]}})(D_{|U_{C^\s,[.]}}),0,0,0)
\ar[d]^-{(1,0)}\\(\widetilde{Fil}^0\Omega^\bullet_{.,T_{C^\s,[.]}}(-h_{|U_{C^\s,[.]}})(D_{|U_{C^\s,[.]}}),0,0,0)
\ar[d]^{(can_D,0)}\\(\widetilde{Lie}{}_{.,T_{C^\s,[.]}}(D_{|U_{C^\s,[.]}})(-h_{|U_{C^\s,[.]}}),0,0,0)\ar[d]^-{+1}\\
\hbox{$ $}}
\end{array}$$
\begin{flushright}$\square$\end{flushright}
These distinguished triangles are compatible with the analogous ones where $C^\s$ is replaced by $J^\s$. Thanks to \ref{MVsynii}, we are thus reduced to prove that $$R f_{T_{C^\s,[.],.},*}\widetilde{Lie}{}_{.,T_{C^\s,[.]}}(D_{|U_{C^\s,[.]}})(-h_{|U_{C^\s,[.]}})\to Rm_* R
f_{T_{J^\s,[.],.},*}\widetilde{Lie}{}_{.,T_{J^\s,[.]}}(D_{|U_{J^\s,[.]}})(-h_{|U_{J^\s,[.]}})$$ is invertible. We know from Prop. \ref{liesurj}
\ref{liesurjii} that $$ch_{Lie}:m^*Lie(D)\to Lie(m^{\s,*}D)$$ is invertible. As before, it follows easily from  (\ref{exR}) and Lem. \ref{defMV} \ref{defMVii} that $$Lie(D)(-h)\to
Rm_*(Lie(m^{\s,*}D)(-h_{|J}))$$ is invertible as desired.  \begin{flushright}$\square$\end{flushright}

Let us write down the \emph{complete Mayer-Vietoris triangles} encoded in the above proposition.

\begin{cor} \label{corMVsyn} \debrom
\item \label{corMVsyni} If $M$ is a locally  free crystal of finite type of $((C^\s/\Sigma_\infty)_{crys,et},\O)$, there is a canonical distinguished triangle in $D(C^{\s,\N}_{et},\O^{crys}_.)$:
$$\xymatrix{Ru_*M_.\ar[r]&(\oplus_{v\in |Z|}\iota_{C_v,*}Ru_*\iota_{C_v^\s}^*M_.)\oplus Rj_*Ru_*j^*M_.\ar[r]&\oplus_{v\in |Z|}Rj_*\iota_{U_v,*}Ru_*\iota_{U_v}^*M_.\ar[r]^-{+1}&}$$

\item \label{corMVsynii} For $D\in \DC(C^\s)$ and $h\in \G(C,M_C/\Gm)$,  there is a canonical distinguished triangle in  $D(Mod^{1,\varphi}(C^{\s,\N}_{et},\widetilde \O^{crys}_.))$:
$$\begin{array}{r}\xymatrix{\cal S_{et,.,C^\s}^{1,\varphi}(-h)(D)\ar[r]&(\oplus_{v\in |Z|}\iota_{C_v,*}\cal S_{et,.,C_v^\s}^{1,\varphi}(-h_{|C_v})(\iota_{C_v^\s}^*D))\oplus Rj_*\cal S_{et,.,U}^{1,\varphi}(j^*D)&&}\\
\xymatrix{\ar[r]&\oplus_{v\in |Z|}Rj_*\iota_{U_v,*}\cal S_{et,.,U_v}^{1,\varphi}(\iota_{U_v}^*j^*D)\ar[r]^-{+1}&}\end{array}$$ and similarly without the superscripts ${}^{(1,\varphi)}$ in $D(C^{\s,\N}_{et},\Z/p^.)$.
\finrom
\end{cor}
Proof. This follows from the previous proposition thanks to Lem. \ref{iotaexact} and the natural distinguished triangle describing  $\smash{Rm_*}$ (Lem. \ref{defMV} \ref{defMVii}).
\begin{flushright}$\square$\end{flushright}

\subsection{A localization triangle} ~~ \\
\label{alt}

We show that for smooth divisors, the twisted syntomic complex of a logarithmic Dieudonn\'e crystal with trivial residue can be recovered from a diagram of syntomic complexes without twists. This can be
viewed as a refined version of the expected localization triangle for compactly supported log crystalline cohomology as defined in \cite{Ts1}.

In order to make a precise statement, we will need the functor $$R\underline \G^{Z_J}(J,-):D^+(J_{et}^{+,\N},\Z/p^.)\longrightarrow
D^+(J_{et}^{\N},\Z/p^.)$$ defined in Def. \ref{defvan} \ref{defvanii}. In order to treat simultaneously the case of
$J^\s$ and $C^\s$, it will be convenient to modify slightly the notations introduced in (\ref{diagdiag}) as follows:
\begin{eqnarray}\label{diagdiagmod}\xymatrix{Z_J\ar[r]^-{z_J}
\ar[d]_{m_Z}&J\ar[d]^m&\ar[l]_-{o_J}J^\s\ar[d]^{m^\s}\\ Z_C
\ar[r]^-{z_C}&C&\ar[l]_-{o_C}C^\s}
\end{eqnarray}

\begin{prop} \label{tdspe} Consider
the morphism of extremal type $z_J:Z_J\to J$
and let us furthermore denote  $(-)_{|J^+}:=\rho^*$ the pullback via the following morphism (Def. \ref{defitype} \ref{defitypeiii})
$$\rho:((J^+/\Sigma_.)_{crys,et},\O)\to ((J/\Sigma_.)_{crys,et},\O).$$

\debrom \item \label{tdspei} If $M_.$ is a crystal of $(J/\Sigma_.)_{crys,et}$,  there is a canonical morphism $$o_{J,*} Ru_*(o_{J}^*M_.(-Z_{J}))\to
R\underline\G^{Z_J}(J,Ru_*(M_{.,|J^+}))$$ in $D(J^\N_{et},\O^{crys}_.)$.  It is an isomorphism if $M_.$ is locally free. \petit

\item \label{tdspeii} If $D$ is in $\DC(J)$, there is a canonical isomorphism
$$o_{J,*}\cal S_{et,.,J^\s}^{1,\varphi}(-Z_{J})(o_{J}^*D)\simeq  R\underline\G^{Z_J}(J,\cal S_{et.,.,J^+}^{1,\varphi}(D_{|J^+}))$$
in $D(Mod^{1,\varphi}(J^\N_{et},\widetilde\O^{crys}_.))$. A similar isomorphism  holds in $D(J^\N_{et},\widetilde \O^{crys}_.) )$ without the
superscripts ${}^{1,\varphi}$.

\item \label{tdspeiii} The statements \ref{tdspei} and \ref{tdspeii} are functorial with respect to $D$, $C$ and $Z$. Moreover, they hold verbatim if the letter $J$ is replaced by the letter $C$.
    \finrom
\end{prop}
Proof. We only prove \ref{tdspeii} since \ref{tdspei} is easier and will be essentially proven along the way. Using Lem. \ref{lememb3}, we choose an object
$\smash{Y_{[.],J^\s}}=(\smash{U_{[.],J^\s}}/J^\s,\smash{Y_{[.],J^\s}},\iota,F)$ in $\smash{HR^{\s,et,ex}_F(J^\s)}$. Let us denote $\smash{U_{[.],J}}$
(resp. $\smash{Y_{[.],J}}$) the diagram obtained from $\smash{U_{[.],J^\s}}$  (resp.  $\smash{Y_{[.],J^\s}}$) by forgetting log structures. Let us
furthermore denote $\smash{U_{[.],Z_J}}$ (resp. $\smash{Y_{[.],Z_J}}$) the diagram obtained from $\smash{U_{[.],X^\s}}$   (resp.
$\smash{Y_{[.],J^\s}}$) as follows.
For each $\nu$ and $\delta$, the vertex $\smash{U_{[\nu],\delta,Z_J}}$ (resp. $\smash{Y_{[\nu],\delta,k,Z_J}}$) is the support of the log structure of $U_{[\nu],\delta,J^\s}$ (resp. $Y_{[\nu],\delta,k,J^\s}$) as defined in Lem.-Def. \ref{suppcenter}  and in particular, a reduced closed subscheme of $U_{[\nu],\delta,J}$ (resp.  $Y_{[\nu],\delta,k,J}$). If $S=J$, $J^\s$, or $Z_J$, we let furthermore
$\smash{U_{[.],S}}\to \smash{T_{[.],S}}$ denote the logarithmic divided power envelope of $\smash{U_{[.],S}}\to \smash{Y_{[.],S}}$.
Summarizing, we
have thus obtained the following commutative diagram of diagrams of simplicial $p$-adic log schemes:
\begin{eqnarray}\label{diagembspe}\xymatrix{Y_{[.],Z_J}\ar[r]^-{z_Y}&Y_{[.],J}&\ar[l]_-{o_J}Y_{[.],J^\s}\\
T_{[.],Z_J}\ar[r]^-{z_T}\ar[u]&T_{[.],J}\ar[u]&\ar[l]_-{o_T}\ar[u]T_{[.],J^\s}\\
U_{[.],Z_J}\ar[r]^-{z_U}\ar[u]&U_{[.],J}\ar[u]&\ar[l]_-{o_U}\ar[u]U_{[.],J^\s}}\end{eqnarray}
We claim that each square of this diagram is cartesian. Let us explain this using that $\smash{Y_{[.],J^\s}}$ is in $\smash{HR^{\s,et,ex}_F(J^\s)}$. On
the right side, this follows from the fact that each $\smash{U_{[\nu],J^\s,\delta}\to Y_{[\nu],J^\s,\delta,k}}$ is an exact closed immersion. Note
that this implies in particular that the morphisms of schemes underlying $o_T$ are isomorphisms and that the morphism of topoi
$\smash{T^\s_{[.],.,et}}\to \smash{T_{[.],.,et}}$ is an equivalence. Next, let us choose a $p$-basis of the form $(\underline s,t)$ for
$\smash{Y_{[\nu],J^\s,\delta}}$. On the left side, the exterior square is cartesian since the ideal of the closed immersions
$\smash{z_{Y,[\nu],\delta,k}}$ and $\smash{z_{U,[\nu],\delta}}$ are both  generated by the image of $t$. Since the image of $t$ does not divide zero
in the structure sheaf of $U$, it follows from \cite{BO} Lem. 3.5 that the top square is cartesian as well.

It is also true that $\smash{U_{[.],J^\s}/J^\s}$ (resp. $\smash{U_{[.],Z_J}/Z_J}$) is the base change of $\smash{U_{[.],J}}/J$ by $o_J:J^\s\to J$
(resp. $z_J:Z_J\to J$). It is thus a hypercovering for the \'etale topology and the previous constructions in particular gives rise to a couple of
arrows \begin{eqnarray}\label{diag1prIJCY}\xymatrix{Y_{[.],Z_J}\ar[r]^-{z}&Y_{[.],J}&\ar[l]_-{o}Y_{[.],J^\s}}
\end{eqnarray} in
$\smash{Diag(HR^{\s,et}_F)}$ above the top line in (\ref{diagdiagmod}). It will be useful to notice that the formation of (\ref{diag1prIJCY}) is
functorial with respect to the object $\smash{Y_{[.],J^\s}}\in  \smash{HR^{\s,et,ex}_F(J^\s)}$ which has been chosen at the beginning of the proof.

Let us denote $h\in \G(J,M_J/\Gm)$, the diagram of effective log divisors corresponding to $-Z_J$. In virtue of Lem. \ref{injtwist2}
\ref{injtwist2ii} and cartesianity properties explained before, we have canonical exact sequences $$\begin{array}{cl}
\xymatrix{0\ar[r]& o_{T,*}(\O(-h_{|U_{[.],J^\s}})_{T_{[.],J^\s,.}})\ar[r]& \O_{T_{[.],J,.}}\ar[r]& z_{T,*}(\O_{T_{[.],Z_J,.}})\ar[r]& 0}&\hbox{in $Mod(T_{[.],J,.,et},\O)$}\\
\xymatrix{\hbox{and }&0\ar[r]&o_{U,*}\O(-h_{|U_{[.],J^\s}})\ar[r]& \O\ar[r]& z_{U,*}\O\ar[r]& 0}&\hbox{in $Mod(U^\N_{[.],J,et},\O)$}\end{array}$$
which  are compatible with each other in the obvious way. Tensoring the first  (resp. second) one over $\O$ with $\smash{D_{T_{[.],J,.}}}$ (resp. with
$\smash{Lie(D_{|U_{[.]},J})}$ then pushing forward to $\smash{T_{[.],J,.}}$) and restricting scalars to $\smash{\O^{crys}_.}$ gives an exact sequence
which may be identified  with the first (resp. second) line of the following commutative diagram  of $\smash{Mod(T_{[.],J,.,et},\O^{crys}_.)}$ by using
that for $M$ locally free $M\otimes f_*N\simeq f_*(f^*M\otimes N)$ (resp. by using Prop. \ref{liesurj} \ref{liesurjii} and Lem. \ref{lemlielfft}).
$$\xymatrix{o_{T,*}((o_J^*D)_{T_{[.],J^\s,.}}(-h_{|T_{[.],J^\s}}))\ar@{^{(}->}[r]\ar@{->>}[d]& D_{T_{[.],J,.}}\ar@{->>}[r]\ar@{->>}[d]& z_{T,*}((z_J^*D)_{T_{[.],Z_J,.}})\ar@{->>}[d]\\
o_{T,*} (Lie_{.,T_{[.],J^\s}}((o_{J}^*D)_{|U_{[.]},J^\s})(-h_{|T_{[.]},J^\s}))\ar@{^{(}->}[r]&
Lie_{.,T_{[.],J}}(D_{|U_{[.]},J})\ar@{->>}[r]& z_{T,*} Lie_{.,T_{[.],Z_J}}((z_J^*D)_{|U_{[.]},Z_J})}$$

Next, let $q\ge 1$ and consider the functoriality morphisms induced by (\ref{diagembspe}):
$$\xymatrix{o_{T,*}\Omega^q_{T_{[.],J^\s,.}}(o_{J}^*D(-h)_{|U_{[.],J^\s}})&\ar[l]_-{o_T^*} \Omega^q_{T_{[.],J,.}}(D)\ar[r]^-{z^*_T}&
z_{T,*}\Omega^q_{T_{[.],Z_J,.}}(z_{J}^*D)}$$ as well as the natural morphism
$$\xymatrix{o_{T,*}\Omega^q_{T_{[.],J^\s,.}}((o_{J}^*D)(-h)_{|U_{[.],J^\s}})\ar[r]^-{nat}&
o_{T,*}\Omega^q_{T_{[.],J^\s,.}}((o_{J}^*D)_{|U_{[.],J^\s}})}$$ Using flatness of the realization of $o_{J}^*D(-h)$ over $(T_{[.],J^\s,.,et},\O)$ and
comparing the local description of logarithmic differentials on $T_{[.],J^\s,.}$ and $T_{[.],J,.}$, we find that $o^*_T$ is monomorphic. The arrow
$nat$ is monomorphic as well, thanks to Lem. \ref{injtwist2} \ref{injtwist2ii} and flatness of the logarithmic differentials on $T_{[.],J^\s,.}$. \petit

\emph{Claim}: The image of $nat$ coincides with $o_T^*(Ker\, z_T^*)$.  \petit

\noindent Let us prove this. We can assume $D=\O$, fix $\nu$, $\delta$ and choose a $p$-basis $((s_1,\dots,s_d),t)$ for $\smash{Y_{[\nu],J^\s,\delta}}$. Then $((s_1,\dots,s_d,t),\emptyset)$ is a $p$-basis for  $\smash{Y_{[\nu],J,\delta}}$ and  $((s_1,\dots,s_d),\emptyset)$ is a $p$-basis for $\smash{Y_{[\nu],Z_J,\delta}}$. We find in particular the following bases for the modules of differential forms of degree $q$:   $$\begin{array}{rcl} \Omega^q_{T_{[\nu],\delta,J^\s,k}}(-h_{|T^\s_{[.],J^\s,.}}):&((tds_{i_1}\wedge\dots \wedge ds_{i_q})_{i_1<\dots<i_q}, (tds_{i_1}\wedge\dots \wedge ds_{i_{q-1}}\wedge d\log t)_{i_1<\dots<i_q-1})\\
\Omega^q_{T_{[\nu],\delta,J,k}}:&((ds_{i_1}\wedge\dots \wedge ds_{i_q})_{i_1<\dots<i_q}, (ds_{i_1}\wedge\dots \wedge dx_{i_{q-1}}\wedge dy)_{i_1<\dots<i_q-1})\\
\Omega^q_{T_{[\nu],\delta,Z_J,k}}:&(ds_{i_1}\wedge\dots \wedge ds_{i_q})_{i_1<\dots<i_q}\end{array}$$ The claim follows immediately , thanks to the
following isomorphisms of rings: $$\O_{T_{[\nu],\delta,J,k}}\simeq o_{T,*}\O_{T_{[\nu],\delta,J^\s,k}} \hspace{1cm} \hbox{ and }\hspace{.5cm}
\O_{T_{[\nu],\delta,J,k}}/(t)\simeq z_{T,*}\O_{T^\s_{[\nu],\delta,Z_J,k}}$$

Putting everything together, we find compatible exact sequences $$\xymatrix{o_{T,*}
Fil^i\Omega^\bullet_{.,T_{[.],J^\s}}(-h_{|U_{[.],J^\s}})((o_J^*D)_{|U_{[.],J^\s}})\ar@{^{(}->}[r] &Fil^i\Omega^\bullet_{.,T_{[.],J}}(D_{|U_{[.],J}})
\ar@{->>}[r] &z_{T,*} Fil^i\Omega^\bullet_{.,T_{[.],Z_J,.}}((z_J^*D)_{|U_{[.],Z_J}})}$$ for $i=0,1$ on $\smash{(T_{[.],J,.,et},\O^{crys}_.)}$. These
sequences are compatible with each other via $1$ and $Fr$.  Let us now restrict scalars to $\smash{\widetilde\O^{crys}_.}$ and apply the functor
$\jj^*$. Since $\smash{o_{T,*}}$ and $\smash{z_{T,*}}$ are exact functors, they commute to $\jj^*$ and this yields exact sequences
$$\xymatrix{o_{T,*} \widetilde{Fil}{}^i\Omega^\bullet_{.,T_{[.],J^\s}}(-h_{|U_{[.],J^\s}})((o_J^*D)_{|U_{[.],J^\s}})\ar@{^{(}->}[r]
&\widetilde{Fil}{}^i\Omega^\bullet_{.,T_{[.],J}}(D_{|U_{[.],J}}) \ar@{->>}[r]
&z_{T,*} \widetilde{Fil}{}^i\Omega^\bullet_{.,T_{[.],Z_J}}((z_J^*D)_{|U_{[.],Z_J}})}$$
(injectivity of the left arrow is checked easily by diagram chasing, using that $\smash{T_{Z_J}}$ is flat over $\Sigma_\infty$ for the case $i=0$ and Lem. \ref{lemfil} \ref{lemfili} together with the fact that the map $$\xymatrix{o_{T,*}H^{-1}\widetilde{Lie}_{.,T_{[.],J^\s}}((o_J^*D)_{|U_{[.],J^\s}})(-h_{|T_{[.],J^\s}}) \ar@{->}[r]&
H^{-1}\widetilde{Lie}_{.,T_{[.],J}}(D_{|U_{[.],J}})}$$ is monomorphic for the case $i=1$). These exact sequences are compatible with each other via $1$ and $\varphi$
(use Lem. \ref{cartgen} \ref{cartgenii}).

We may thus interpret them as exact sequences in the category of $(1,\varphi)$-modules over
$\smash{(T_{[.],J,.,et},\widetilde \O^{crys}_.)}$ or equivalently as an isomorphism  $$\xymatrix{o_{T,*}\cal
S^{1,\varphi}_{et,.,T_{[.],J^\s}}(-h_{|U_{[.],J^\s}})((o_J^*D)_{|U_{[.],J^\s}})\ar[r]^-\sim & MF(\cal S^+_{|J}\to z_{J,*}\cal S^+_{|Z_J})}$$ where $$\cal
S^+:=\cal S^{1,\varphi}_{et,.,T_{[.],J^+}}(D_{|J^+})$$ is the syntomic complex on $J^+$ associated to the object of $\smash{{HR^{et}_F(J^+)}}$ which is
defined by the arrow $\smash{Y_{[.],Z_J}\to Y_{[.],J}}$.  The isomorphism of the proposition follows by applying $\smash{Rf_{T_{[.],J,.,*}}}$ (use
Lem. \ref{soritesvan} \ref{soritesvani}, \ref{soritesvaniv}). Finally, we note that the independence of choices is a formal consequence of the
connectedness of the category $HR_{F}^{\s,et,ex}(J^\s)$ (which, in turn, is ensured by the existence of finite non empty products,  see Lem. \ref{lememb3}).
\begin{flushright}$\square$\end{flushright}

Let us write down the \emph{localization triangles} encoded in the above proposition.

\begin{cor} \label{cortdspe} \debrom
\item \label{remtdspei} If $M_.$ is a locally free crystal of $((X/\Sigma_.)_{crys,et},\O)$, there is a canonical distinguished triangle in $\smash{D(X^\N_{et},\O^{crys}_.)}$:
$$\xymatrix{Ru_*(o_{X}^*M_.(-Z_{|X}))\ar[r]&Ru_*M_.\ar[r]&z_{X,*}Ru_*(z_X^*M_.)\ar[r]^-{+1}&}$$

\item \label{remtdspeii} If $D$ is in $\cal DC(X)$, there is a canonical distinguished triangle in $\smash{D(X^\N_{et},\widetilde \O^{crys}_.)}$:
$$\xymatrix{\cal S_{X^\s,.,et}^{1,\varphi}(-Z_{|X})(o_{X}^*D)\ar[r]&\cal S_{X,.,et}^{1,\varphi}(D)\ar[r]&z_{X,*}\cal S_{Z_X,.,et}^{1,\varphi}(z_X^*D)\ar[r]^-{+1}&}$$
and similarly without the superscripts ${}^{(1,\varphi)}$ in $\smash{D(X^\N_{et},\Z/p^.)}$
\finrom
\end{cor}
Proof. Assertion \ref{remtdspei} (resp. \ref{remtdspeii}) follows from Prop. \ref{tdspe} \ref{tdspei} (resp. Prop. \ref{tdspe} \ref{tdspeii}) by Lem. \ref{soritesvan} \ref{soritesvani}.
\begin{flushright}$\square$\end{flushright}


\section{Dieudonn\'e crystals for semi-stable Abelian varieties} \label{secdcfssav}

The purpose of this chapter is to define a Dieudonn\'e crystal over $(C^\s/\Sigma_\infty)_{crys,et}$ associated to a semi-Abelian scheme $A/C$ whose restriction to $U$ is Abelian.
We begin with a result concerning the structure of log $1$-motives on a complete discrete valuation ring.

\subsection{D\'evissage of log $1$-motives}
\label{dol1m} ~~ \\

\label{subsectionog1mot} Let $X$ denote the spectrum of a complete discrete valuation ring $R$, and let $s:Spec(k)\rightarrow X$ (resp.
$j:Spec(K)\rightarrow X$) the inclusion of its special (resp. generic) point. For any $X'/X$, we let $s'$, $j'$ denote the morphisms deduced from
$s$, $j$ by base change.

\para The following definition is \cite{KT} Sect. 4.6.1.

\begin{defn} \label{defMlog}

\debrom \item \label{defMlogi} The category $\M_{log}(X)$ of \emph{log $1$-motives over $X$} is defined as follows:

- an object is a triple $(\G,G,f)$ where $\G$ is a  twisted constant group (Def. \ref{defmotabtor}, \ref{defmotabtoriv}), while $G$ is the extension of an
Abelian scheme $B$ by a torus $T$ and $f:\G\rightarrow j_*j^{-1}G$ is a morphism in $Ab(X_{FL})$.

- a morphism from $(\G,G,f)$ to $(\G',G',f')$ is a couple of morphisms $\G\rightarrow \G'$, $G\rightarrow G'$, compatible with $f$ and $f'$.

\item \label{defMlogii} Denote $\B_X$ the category $fet(X)$. Define $\M_{log}/\B_X$ as the fibered category corresponding to the contravariant pseudo-functor $(X'/X)\mapsto \M_{log}(X')$, $f\mapsto f^{-1}$ where $f^{-1}$ is the pullback functor deduced from $(-)_{FL}$.
\finrom
\end{defn}

The functor, sending  a usual $1$-motive $f:\G\rightarrow G$ over $X'$
to the log $1$-motive $f:\G\rightarrow j_*j^{-1}G$ deduced from it, is fully faithful. We may thus identify the category of usual $1$-motives over variable bases with a full subcategory of $\M_{log}$. \petit

Just as usual $1$-motives, any log $1$-motive $(\G,G,f)$ over $X$ comes with a functorial increasing \emph{weight filtration}. It is defined as follows: $$\begin{array}{rcl}Fil^0&=&(\G,G,f)\\
Fil^{-1}&=&(0,G,0)\\Fil^{-2}&=&(0,T,0)\\Fil^{-3}&=&(0,0,0)\end{array}$$ $T$ being the maximal subtorus of $G$ (see Lem. \ref{lemG} \ref{lemGiii} for the functoriality). \petit

We consider the following full $\B_X$-
subcategories of $\M_{log}$:

- $\M{}:$ $1$-motives viewed as log $1$-motives, as explained above.

- $\MT_{log}:$ those $(\G,G,f)$ with $Gr^{-1}=0$, ie. such that $G=T$ is a torus.


- $\MT=\M\cap \MT_{log}$. \petit

\begin{rem} \label{remlog1motex}
The $\cal B_X$-categories $\M$, $\M_{log}$, $\MT$ and $\MT_{log}$ are naturally $\B_X$-$e$-exact categories (ie. the fibers are naturally exact and the pullback functors are $e$-exact) and the inclusion functors are $\cal B_X$-$e$-exact. We leave this to the reader, using only the exactness of $j^{-1}$ together with the stability under extensions of the categories of twisted constant groups, tori and extensions of Abelian schemes by tori (see the proof of Lem. \ref{1motex}).
\end{rem}

%
%
\petit


\para Any log $1$-motive $(\G,G,f)$ over $X$ can be seen canonically as object of the category of extensions $\smash{EXT^1_{\M_{log}(X)}((\G,0,0),(0,G,0))}$. A
nice feature of the latter is that it is endowed with the Baer sum $\boxplus$ which is an exact bifunctor underlying the usual addition of the
group $Ext^1_{\M_{log}(X)}$. We will need to extend the bifunctor $\boxplus$ as follows.

\begin{defn} \label{defbaer} Consider an exact category $\cal C$.

\debrom \item \label{defbaeri} We define $\C\square \cal C$ as the category of diagrams of the form

$$\xymatrix{0\ar[r]&A_1\ar[r]^-{i_1}&B_1\ar[r]^-{p_1}&C_1\ar[r]&0\\
0\ar[r]&A_2\ar[u]^{\alpha}\ar[r]^-{i_2}&B_2\ar[r]^-{p_2}&C_2\ar[u]_\gamma\ar[r]&0}$$


\noindent in which both lines are short $e$-exact sequences and both following morphisms are admissible (ie. are part of a short $e$-exact sequence of $\cal C$):
$A_2\hookrightarrow B_1\times B_2$, $a\mapsto (i_1\alpha(a),-i_2(a))$, $B_1\times B_2\twoheadrightarrow C_1$, $(b_1,b_2)\mapsto p_1(b_1)-\gamma
p_2(b_2)$. For simplicity, such an object is usually denoted $(B_1,B_2)$.

\item \label{defbaerii} We define $\boxplus:\C\square \cal C\rightarrow \cal C$ as the functor sending $(B_1,B_2)$ as above to its \emph{Baer sum} $B_1\boxplus B_2:=Ker (B_3\rightarrow
C_1,\, (b_1,b_2)\mapsto p_1(b_1)-\gamma p_2(b_2))$  where $B_3:=Coker(A_2\rightarrow B_1\times B_2,\, a\mapsto (i_1\alpha(a),-i_2(a)))$. \finrom
\end{defn}
The reader is invited to check that this definition makes sense, ie. that the involved kernels and cokernels exist indeed using the axioms of an
exact category (see e.g. \cite{Bu} Def. 2.1). Note that  $\boxplus$ is an $e$-exact functor, if  $\cal C\square\cal C$ is endowed with its obvious exact
structure. The formation of the category $\cal C\square \cal C$ and of the functor $\boxplus$ is moreover canonically pseudo-functorial with respect
to $\cal C$ in the obvious sense. \petit

In our situation, we thus have  a $\B_X$-$e$-exact functor \begin{eqnarray}\label{extbaersum}\boxplus:\M_{log}\square_{\B_X} \M_{log}\rightarrow \M_{log}\end{eqnarray} where
$\smash{\M_{log}\square_{\B_X}\M_{log}}$ denote the natural $\B_X$-$e$-exact fibered category whose fiber at $X'/X$ is $\M_{log}(X')\square\M_{log}(X')$. Let us
indicate a useful computation in a special case. If an object of $\M_{log}(X)\square\M_{log}(X)$ is of the particular form

\begin{eqnarray} \label{calculbaer}\xymatrix{0\ar[r]&(0,G_1,0)\ar[r]^-{(0,id)}&(\G_1,G_1,f_1)\ar[r]^-{(id,0)}&(\G_1,0,0)\ar[r]&0\\
0\ar[r]&(0,G_2,0)\ar[u]^{\alpha}\ar[r]^-{(0,id)}&(\G_2,G_2,f_2)\ar[r]^-{(id,0)}&(\G_2,0,0)
\ar[u]_\gamma\ar[r]&0}\end{eqnarray}
then its Baer sum is naturally isomorphic to $(\G_2,G_1,f_1\circ \gamma +j_*j^{-1}(\alpha)\circ f_2)$.

\para As recalled in Lem. \ref{lemG} \ref{lemGi}, tori and twisted  constant groups are locally trivial for the \'etale topology. Here $X$ is the spectrum of a complete discrete valuation ring and the finite \'etale topology is in fact sufficient (use \cite{Mi1} I., Thm. 4.2 c)). This and Lem. \ref{lemcst} \ref{lemcstiii}
implies in particular that for a twisted constant group $\G$, we have $\epsilon^{-1}\epsilon_*\G\simeq \G$ if $$\epsilon:X_{FL}\to X_{fet}$$ denotes the weak morphism induced by the inclusion $fet(X)\subset FL(X)$. \petit



Consider a fixed uniformizer $t$ of $R$. In $Ab(X_{fet})$, this choice provides a splitting of the valuation exact sequence
$$\xymatrix{0\ar[r]&\epsilon_*\Gm\ar[r]&\epsilon_*j_*j^{-1}\Gm\ar[r]_-v&\ar@/_2pc/[l]_{t}\epsilon_*\Z\ar[r]&0}$$
More generally, consider an extension $G$ of an Abelian scheme $B$ by a torus $T$ and denote $\Gamma^*:=\cal Hom(T,\Gm)$ (resp. $\Gamma^{*\vee}:=\cal Hom(\G^*,\Z)$) the character (resp. cocharacter) group of $T$. Then $t$ induces a split monomorphism $t_T:\epsilon_*\Gamma^{*\vee}\simeq \cal Hom(\epsilon_*\Gamma^*,\epsilon_*\Z)\to  \cal Hom(\epsilon_*\Gamma^*,\epsilon_*j_*j^{-1}\Gm)=\epsilon_*j_*j^{-1}T$. Let us furthermore denote $t_G:\epsilon_*\G^{*\vee}\to \epsilon_*j_*j^{-1}G$ the morphism given by $t_T$ and the inclusion $T\subset G$.

\begin{lem} \label{prelimtdec} The morphism $t_G$ induces the following isomorphisms:
\begin{eqnarray}\label{prelimtdeci}\xymatrix{\epsilon_*G\times\epsilon_*\G^{*\vee}\ar[rr]_{\sim}^{(j^*,t_G)}& & \epsilon_*j_*j^{-1}G}\\
 \label{prelimtdecii} \xymatrix{Hom(\G,G\times \G^{*\vee})\ar[r]_-\sim& Hom(\G,j_*j^{-1}G)}\end{eqnarray}
\end{lem}
Proof. The second isomorphism follows immediately from the first one. Let us thus explain (\ref{prelimtdeci}). First, we note that smooth group schemes are acyclic for $\epsilon_*$ in virtue of Lem. \ref{acycsm} \ref{acycsmi} and Lem. \ref{lemqff} \ref{lemqffi}. Applying this to $T$ and using that $\epsilon_*B\to  \epsilon_*j_*j^{-1}B$ is an isomorphism (this follows from the N\'eron extension property), we get the following commutative diagram with exact lines: $$\xymatrix{0\ar[r]&\epsilon_*T\ar[r]\ar[d]&\epsilon_* G\ar[r]\ar[d]&\epsilon_*B\ar[r]\ar[d]^\wr&0\\
0\ar[r]&\epsilon_*j_*j^{-1}T\ar[r]&\epsilon_* j_*j^{-1}G\ar[r]&\epsilon_*j_*j^{-1}B\ar[r]&0}$$ Diagram chasing reduces us to the case $T=G$ but then, we are done by the construction of $t_T$ explained above.
\begin{flushright}$\square$\end{flushright}
Note that $T$ being functorial with respect to $G$, the isomorphisms of the lemma are functorial as well.


%

%
%

\begin{defn} \label{tdec}

\debrom \item \label{tdeci}  Let $(\G,G,f)$ in $\M_{log}(X)$ and denote $(f_0,v(f))$ the preimage of $f$ under the isomorphism (\ref{prelimtdecii}).  The
resulting couple of log $1$-motives $$((\G,G,f_0),(\G,T,t_T\circ v(f))) \hbox{ in $\M(X)\times \MT_{\log}(X)$}$$ (here we have abusively denoted $t_T$ the morphism $\G^{*\vee}\to j_*j^{-1}T$ obtained using that $\epsilon^{-1}\epsilon_*\G^{*\vee}\simeq \G^{*\vee}$)
coming from this decomposition is called the \emph{$t$-decomposition of $(\G,G,f)$}.

\item \label{tdecii} We say that $(\G,G,f)$ in $\M_{log}(X)$ is a multiple of $t$ if $f_0=0$ (or equivalently if $f=t_G\circ v(f)$). The full subcategory of $\M_{log}(X)$ formed by the multiples of $t$ is denoted $\M^t_{log}(X)$.
\finrom
\end{defn}
Let us emphasize that the $t$-decomposition introduced in \ref{tdeci} is functorial with respect to $(\G,G,f)$ in $\M_{log}(X)$. Let us also notice
that a log $1$-motive $(\G,G,f)$ is in fact in $\M(X)$ if and only $v(f)=0$.


\begin{lem} \label{lemmot1}
\debrom
\item \label{lemmot1i} Let $\M(X)\square \MT^t_{log}(X)$ denote the  full subcategory of $\M_{log}(X)\square\M_{log}(X)$ whose objects are of the form
$$\xymatrix{0\ar[r]&(0,G,0)\ar[r]^{(0,id)}&(\G,G,f)\ar[r]^{(id,0)}&(\G,0,0)\ar[r]&0\\
0\ar[r]&(0,T,0)\ar[u]^{\cup}\ar[r]^{(0,id)}&(\G,T,g)\ar[r]^{(id,0)}&(\G,0,0)\ar[u]_{\parallel}\ar[r]&0}$$
where $(\G,G,f)$ is in $\M(X)$ and $(\G,T,g)$ is in $\MT^t_{log}(X)$ (ie. $g=t_T\circ v(g)$). The functor $(\ref{extbaersum})$ induces an equivalence of categories
$$\xymatrix{\M(X)\square \MT^t_{log}(X)\ar[r]_-\sim^-{\boxplus}& \M_{log}(X)}.$$

\item \label{lemmot1ii} Let $\M(X)\square\MT^0(X)$ denote the full subcategory of $\M(X)\square \MT^t_{log}(X)$ whose objects are those of the above form
satisfying furthermore $g=0$. Then we have an equivalence of categories $$\xymatrix{\M(X)\square \MT^0(X)\ar[r]^-{\boxplus}_-\sim &\M(X)}.$$ \finrom
\end{lem}

Proof. \ref{lemmot1ii} The assertion follows from the functoriality of the weight filtration together with the computation of the Baer sum on $(\ref{extbaersum})$ (alternatively, it is also a consequence of \ref{lemmot1i}).

\ref{lemmot1i} An object of the indicated form is sent by
$\boxplus$ to $(\G,G,f+g)$, where $g$ is seen as a morphism $\G\rightarrow G$. Essential surjectivity thus results from Lem. \ref{prelimtdec}. Full faithfulness on the other hand, boils down to the fact that the $t$-decomposition is functorial with respect to the log $1$-motive.
\begin{flushright}$\square$\end{flushright}

Let $\Z$ and $\Z(1)$ respectively denote the objects $(\Z,0,0)$ and $(0,\Gm,0)$ of $\M(X)$. If $x\in R^\times$ (resp. $K^\times$), we denote
$Kum(x):=(\Z,\Gm,x:\Z\rightarrow \Gm)$ (resp. $Kum(x):=(\Z,\Gm,x:\Z\rightarrow j_*j^{-1}\Gm)$) the $1$-motive (resp. log
$1$-motive) where we consider $x$ as a map sending 1 to $x$ . We call it a \emph{Kummer $1$-motive} (resp. \emph{log
$1$-motive}). Note that sending $x$ to the class of $Kum(x)$ realizes isomorphisms of groups $R^\times\simeq \smash{Ext^1_{\M(X)}(\Z,\Z(1))}$ and $K^\times\simeq \smash{Ext^1_{\M_{log}(X)}(\Z,\Z(1))}$ (compatibility with group laws relies on the computation of $\boxplus$ on $(\ref{extbaersum})$). If $M$ is a twisted constant group, we denote furthermore $M\otimes Kum(x)$ the $1$-motive (resp. log $1$-motive)
$(M,M\otimes \Gm,id_M\otimes x)$.

\begin{lem}  \label{lemmot2} The Baer sum functor induces an equivalence between $\MT^t_{log}(X)$ and the full subcategory of $\M_{log}(X)\square\M_{log}(X)$ whose objects are of the form
$$\xymatrix{0\ar[r]&(0,M\otimes \Gm,0)\ar[r]^-{(0,id)}&M\otimes Kum(t)\ar[r]^-{(id,0)} &(M,0,0)\ar[r]&0 \\
0\ar[r]&(0,M\otimes \Gm,0)\ar[u]^{\parallel}\ar[r]^-{(0,id)} &(\G,M\otimes \Gm,0)\ar[r]^-{(id,0)}&(\G,0,0)\ar[u]_f\ar[r]&0}$$
with $M$ a twisted constant group. \end{lem}

Proof. An object of the above form is sent to $(\G,M\otimes \Gm,g)$, where $g$ is
the composition of the morphism $f:\G\rightarrow M$ with the morphism $id_M\otimes t:M\rightarrow M\otimes j_*j^{-1} \Gm$. Essential surjectivity
results from the definition of $\MT^t_{log}(X)$ (take $M=\G^{*\vee}$) while full faithfulness results from the fact that $id_M\otimes t$ is
monomorphic.
\begin{flushright}$\square$\end{flushright}

\begin{prop} \label{prol} Let $\C/\B_X$ be a stack of exact categories and assume given a $\B_X$-$e$-exact cartesian functor $F:\M\rightarrow \C$ (ie. a collection of $e$-exact functors $F_{X'}:\M(X')\rightarrow \C(X')$ together with isomorphisms $f^{*}F_{X'}\simeq F_{X''}f^{-1}$ for all $f:X''\rightarrow X'$, such that the obvious composition constraint holds). The category of extensions of $F$ to a $\B_X$-$e$-exact functor $$F_{log}:\M_{log}\rightarrow \C$$ is canonically equivalent
to the discrete category whose underlying set is the set of homomorphisms $$|F_{log}|:K^\times\simeq Ext^1_{\M_{log}(X)}(\Z,\Z(1))\rightarrow
Ext^1_{\C(X)}(F(\Z),F(\Z(1)))$$ extending the homomorphism
$$|F|:R^\times\simeq Ext^1_{\M(X)}(\Z,\Z(1))\rightarrow Ext^1_{\C(X)}(F(\Z),F(\Z(1)))$$
induced by $F$. This equivalence is independent  on the choice of $t$ made earlier.
\end{prop}
Proof. Consider the category $\smash{\C artEx^F_{\B_X}(\M_{log},\C)}$ whose objects (resp. morphisms) are the cartesian $\B_X$-$e$-exact functors extending $F$ (resp.
$\B_X$-morphisms between them) and the discrete category $\smash{Hom^{|F|}(K^\times,Ext^{1}_{\C(X)}(F(\Z),F(\Z(1))))}$ whose objects are the homomorphisms
extending  $|F|$.
We need to show that the natural functor $$\begin{array}{rcccc}|-|&:&\C artEx^F_{\B_X}(\M_{log},\C)&\rightarrow & Hom^{|F|}(K^\times,Ext^{1}_{\C(X)}(F(\Z),F(\Z(1))))\\
&& F_{log}&\mapsto& |F_{log}|\end{array}$$
is an equivalence. Let us fix a uniformizer $t$ and consider the following natural factorization (via
restriction to $\M^t_{log}$) of the functor $|-|$: $$\diagram{\C art^F_{\B_X}(\M_{log},\C)&\hfl{res^t}{}&\C art^{res^t
F}_{\B_X}(\M^t_{log},\C)&\hfl{|-|^t}{}&  Hom^{|F|}(K^\times,Ext^{1}_{\C(X)}(F(\Z),F(\Z(1))))}$$ We are going to show that both functors $res^t$ and
$|-|^t$ are equivalences. \petit

The given $\B_X$-$e$-exact functor $F$ induces a pseudo-commutative diagram which corresponds to the solid part of the following diagram where the vertical
equivalences  are given by Lem. \ref{lemmot1}.
$$\xymatrix{\M\square_{\B_X}\M\ar @<+2pt> `u[r] `[rrrr] [rrrr]\ar[rr]^{F\square F}\ar@/_1pc/[dddr]_{\boxplus}&&\C\square_{\B_X}\C\ar'[d][dd]^{\boxplus}&&\M_{log}\square_{\B_X}\M_{log}
\ar @{-->}[ll]_{F_{log}\square F_{log}}\ar@/^1pc/[dddl]^{\boxplus}\\
&
\M\square_{\B_X}\MT^0 \ar[ul]\ar[ur]^{\smash{\tilde F}}\ar[dd]^{\boxplus}_{\wr}\ar[rr]&& \M\square_{\B_X}\MT^t_{log}\ar @{-->}[ul]_{\smash{\tilde F}_{log}}\ar[dd]^{\vspace{.5cm}\boxplus}_{\wr}\ar[ur]\\
&&\C\\
& \M\ar[ur]^{F}\ar[rr]&& \M_{log}\ar @{-->}[ul]_{F_{log}}}$$
This diagram together with pseudo-functoriality of $\square$ and $\boxplus$ shows that it is equivalent to extend  either one of the $\cal B_X$-exact functors $F$,
$F\square F$ or $\smash{\tilde F}$ respectively into a $\cal B_X$-exact functor  $F_{log}$, $F_{log}\square F_{log}$ or $\smash{\tilde F_{log}}$. Now extending $\smash{\tilde
F}$ into $\smash{\tilde F_{log}}$ is also equivalent to extending $res^t(F):\MT^0\rightarrow \C$ into a functor $\smash{\MT^t_{log}}\rightarrow \C$.
This shows that  $res^t$ is an equivalence.

Using Lem. \ref{lemmot2} together with \'etale descent in $\C$, we are reduced to the case  $\G\simeq \Z^a$ and $M\simeq \Z^b$. Thus, we see that extending $res^t(F)$ to $\MT^t_{log}$ simply boils down to the choice of the image of the Kummer log $1$-motive $Kum(t)$. Since the group
$K^\times /R^\times$ is generated by $t$, this shows that $|-|^t$ is an equivalence.
\begin{flushright}$\square$\end{flushright}


\subsection{The semi-stable Dieudonn\'e functor}
\label{tssdf} ~~ \\

The goal of this section is to construct the Dieudonn\'e crystal of a semi-Abelian scheme over $C_v$ 
and then over $C$ by gluing the local constructions for $v$ in $Z$ with the usual Dieudonn\'e crystal of \cite{BBM} over $U$.

\para \label{paradeffv} The Dieudonn\'e functor  to be constructed, will take its values into the exact category of Dieudonn\'e
crystals. In view of cohomological methods, it will be convenient to consider the larger category of $(f,v)$-modules defined just below.


\begin{defn} \label{deffvmod} Consider a ringed topos $(E,A)$ (with $A$ commutative) and an object $X$ of $E$ together with an endomorphism $F:X\to X$. We still denote $F$, the associated localization morphism $F:(E_{/X},A_{|X})\to (E_{/X},A_{|X})$.
\debrom \item \label{deffvmodi} We define the category $Mod^{fv}(E_{/X},A_{|X})$ of $(f,v)$-modules as follows:

\noindent - An object is a triple $(M,f_M,v_M)$ where $M$ is a module of $(E_{/X},A_{|X})$ and  $f_M:F^{-1}M\rightarrow M$, $v_M:M\rightarrow F^{-1}M$
are morphisms in $Mod(E_{/X},A_{|X})$.

\noindent - A morphism $(M,f_M,v_M)\to (N,f_N,v_N)$ is a morphism $a:M\to N$ in $Mod(E_{/X},A_{|X})$ which is compatible with the $f$'s and $v$'s.


\item \label{deffvmodii} We define a bifunctor $$\cal Hom_{A_{|X}}^{fv}:Mod^{fv}(E_{/X},A_{|X})^{op}\times Mod^{fv}(E_{/X},A_{|X})\to Mod^{fv}(E_{/X},A_{|X})$$ by the formula $$\cal Hom_{A_{|X}}^{fv}((M,f_M,v_M),(N,f_N,v_N)):= (\cal
Hom_{A_{|X}}(M,N),\cal Hom_{A_{|X}}(v_{M},f_{N}),\cal Hom_{A_{|X}}(f_{M},v_{N}))$$ where  $\smash{\cal Hom_{A_{|X}}}$ denotes  inner
homomorphisms over $\smash{(E_{/X},A_{|X})}$.
\finrom \end{defn}

Consider the category $\smash{E^{fv}_{/X}}$, of triples $(Y,f_Y,v_Y)$ where $Y$, $f_Y:F^{-1}Y\to Y$ and $v_{Y}:Y\to F^{-1}Y$ are in $\smash{E_{/X}}$. The ring $\smash{(A_{|X},f_A,v_A)}$ where $f_A$ is the canonical identification $F^{-1}A_{|X}\simeq A_{|X}$ and $v_A=f_A^{-1}$ will simply be denoted $\smash{A_{|X}}$.

\begin{lem} \label{lemfv1} Let $(E,A,X,F)$, $E_{/X}^{fv}$ and  $A_{|X}$ as above.

\debrom \item \label{lemfv1i}  The pair $\smash{(E_{/X}^{fv},A_{|X})}$ is a ringed topos. It is pseudo functorial with respect to $(E,A,X,F)$ if a
morphism $(E,A,X,F)\to (E',X',A',F')$ means a morphism of ringed topoi $g:(E,A)\to (E',A')$ together with an isomorphism $g^{-1}X'\simeq X$ which is
compatible with $F$ and $F'$.

\item \label{lemfv1ii} There is a canonical isomorphism $Mod^{fv}(E_{/X},A_{|X})\simeq Mod(\smash{E^{fv}_{/X}},A_{|X})$. Via this identification the adjunction $(g^*,g_*):Mod(\smash{E^{fv}_{/X}},A_{|X})\to Mod(\smash{E'^{fv}_{/X'}},A'_{|X'})$ attached to $(E,A,X,F)\to (E',X',A',F')$ translates as follows: $$\begin{array}{rrcll}
g^*(M',f',v')=(g^*M',f,v)\hbox{, where}&f&:&F^{-1}g^*M'\simeq g^*F'^{-1}M'\mathop{\longrightarrow}\limits^{g^*(f')}g^*M'&\\
\hbox{and }&v&:&g^*M'\mathop{\longrightarrow}\limits^{g^*(v')} g^*F'^{-1}M'\simeq F^{-1}g^*M'&\\
g_*(M,f,v)=(g_*M,f',v')\hbox{, where}&f'&:&F'^{-1}g_*M\simeq g_*F^{-1}M\mathop{\longrightarrow}\limits^{g_*(f)}g_*M&\\
\hbox{ and }&v'&:&g_*M\mathop{\longrightarrow}\limits^{g_*(v)} g_*F^{-1}X\simeq F'^{-1}g_*M& \hbox{ )}\end{array}$$
(here the base change isomorphisms are due to the fact that $F'$ is a localization morphism whose pullback by $g$ is $F$).
\item \label{lemfv1iii} The forgetful functor $\chi^{-1}:Mod^{fv}(E_{/X},A_{|X})\to Mod(E_{/X},A_{|X})$ is the pullback of a morphism of ringed topoi. It has moreover a left adjoint $\chi_!$ which is exact.
\finrom
\end{lem}
Proof. \ref{lemfv1i} and \ref{lemfv1ii} are routine.

\ref{lemfv1iii} Let us give an explicit description of $\chi_!$ in order to check that it is exact. Let $[F^{-1}]^\N*[F_!]^\N$ denote the free
associative monoid with unit on the labels $[F^{-1}]$ and $[F_!]$. Sending $[F^{-1}]$ to $F^{-1}$ and $[F_!]$ to $F_!$, defines an action of   $[F^{-1}]^\N*[F_!]^\N$
on the category $Mod(E_{/X},A_{|X})$ in the strict sense. Explicitly,  given a word $w=[F^{-1}]^{\alpha_1}[F_!]^{\beta_1} \dots
[F^{-1}]^{\alpha_l}[F_!]^{\beta_l}$ of length $2l$ ($l\ge 0$, $\alpha_i\ge 1$, $\beta_i\ge 1$) and an $A$-module $M$, we define  $w.M$ as the $A$-module $(F^{-1})^{\alpha_1}(F_!)^{\beta_1} \dots (F^{-1})^{\alpha_l}(F_!)^{\beta_l}M$.
Set $\chi_!M=(N,f_N,v_N)$, where $$N=\mathop\oplus\limits_{w\in [F^{-1}]^\N*[F_!]^\N} w.M\, ,$$ $f_M:F^{-1}M\rightarrow M$ is induced by multiplying the
indices on the left by $[F^{-1}]$, together with the identity $F^{-1}(w.M)=([F^{-1}]w).M$ in component $w$, and $v_M:M\rightarrow F^{-1}M$ is induced by multiplying the indices on the left by
$[F_!]$, together with the adjunction morphism $w.M\rightarrow F^{-1}F_!(w.M)=F^{-1}([F_!]w.M)$ in component $[F_!]w$. The adjunction morphisms $M\rightarrow \chi^{-1}\chi_!M$
and  $\chi_!\chi^{-1}(M,f_M,v_M)\rightarrow (M,f_M,v_M)$ are respectively defined by sending $M$ into the component $1$ and sending the component
$[F^{-1}]^{\alpha_1}[F_!]^{\beta_1} \dots [F^{-1}]^{\alpha_l}[F_!]^{\beta_l}$ into $M$ by applying successively  the morphism $F_!M\rightarrow M$ deduced
from $v_M$ by adjunction $\beta_l$ times, the morphism $f_M$ $\alpha_l$ times, and so on. We leave it to the reader to check the adjunction property.
This description makes it clear that $\chi_!$ is exact, since it is built from the exact functors $F_!$ and $F^{-1}$ using direct sums. Let us mention that the functor
$\chi_*$ has an analogous description with $[F^{-1}]^\N*[F_!]^\N$ replaced by $[F_*]^\N*[F^{-1}]^\N$.
\begin{flushright}$\square$\end{flushright}

\begin{rem} \label{remfv} The bifunctor $\cal Hom^{fv}_{A_{|X}}$ should not be confused with the bifunctor $\smash{\cal Hom_{E_{/X}^{fv},A_{|X}}}$ of inner homomorphisms in the ringed topos $\smash{(E_{/X}^{fv},A_{|X})}$.
\end{rem}


The bifunctor $\cal Hom_{A_{|X}}^{fv}$ is right derivable and thus induces: \begin{eqnarray}\label{HomfvA}R\cal Hom^{fv}_{A_{|X}}(-,-):D^-(E_{/X}^{fv},A_{|X})^{op}\times
D^+(E_{/X}^{fv},A_{|X})\to D^+(E_{/X}^{fv},A_{|X})\end{eqnarray} We will also need a variant $\cal Hom^{fv}$, where the first argument is in the category of
Abelian groups instead of $A$-modules. Deriving gives: \begin{eqnarray}\label{Homfvab}R\cal Hom^{fv}(-,-):D^-(E_{/X}^{fv})^{op}\times D^+(E^{fv}_{/X},A_{|X})\to D^+(E^{fv}_{/X},A_{|X})\end{eqnarray}
We will also use variants where the first argument runs in the derived category of inductive (resp. projective) systems of $(f,v)$-Abelian groups and the bifunctor takes its values in the derived category of projective (resp. inductive) systems of $(f,v)$-modules. Derived functors in this setting can be computed by taking components of the projective or inductive systems accordingly.

\begin{lem} \label{lemfv2}
\debrom \item \label{lemfv2i} One has a canonical bifunctorial isomorphism  $$\cal Ext^{fv,i}_{A_{|X}}((M,f_M,v_M),(N,f_N,v_N))=(\cal
Ext^i_{A_{|X}}(M,N),f,v)$$ where $f$ is induced by $v_M$ and $f_N$ while $v$ is induced by $f_M$ and $v_N$.

\item \label{lemfv2ii} Let $A_{|X}\otimes(-):Mod^{fv}(E_{/X},\Z)\to Mod^{fv}(E_{/X},A_{|X})$ denote the pullback functor induced by $(E,A,X,F)\to (E,\Z,X,F)$. There is a canonical isomorphism of bifunctors $D^-(E_{/X}^{fv})^{op}\times D^+(E_{/X}^{fv},A_{|X})\to D^+(E_{/X}^{fv},A_{|X})$:
$$R\cal Hom^{fv}_{A_{|X}}(A_{|X}\Ltens(M,f_M,v_M),(N,f_N,v_N))\mathop\simeq\limits^{adj} R\cal Hom^{fv}((M,f_M,v_M),(N,f_N,v_N))$$

\item \label{lemfv2iii} Let $l_*:\smash{Mod^{fv}(E_{/X},A_{|X})^{\N}\to Mod^{fv}(E_{/X},A_{|X})}$ denote the functor taking a projective system of $(f,v)$-modules to its inverse limit. There is a canonical isomorphism of bifunctors $\smash{D^-(E_{/X}^{fv,\N^{op}})^{op}}\times \smash{D^+(E_{/X}^{fv},A_{|X})}\to \smash{D^+(E_{/X}^{fv},A_{|X})}$: $$R\cal Hom^{fv}(\limi (M_.,f_{M_.},v_{M_.}),(N,f_N,v_N))\simeq Rl_* R\cal Hom^{fv}((M_.,f_{M_.},v_{M_.}),(N,f_N,v_N))$$
A similar isomorphism holds with  $R\cal Hom^{fv}_{A_{|X}}$ instead of $R\cal Hom^{fv}$.

\item \label{lemfv2iv}  For each $n\ge 0$, there is a canonical morphism of functors $D^-(E_{/X}^{fv},A_{|X})\to D(E_{/X}^{fv},A_{|X})$ $$id\to R\cal Hom^{fv}_{A_{|X}}(\tau_{\le n}R\cal Hom^{fv}_{A_{|X}}(-,A_{|X}),A_{|X})$$
The morphisms obtained for different values of $n$ are compatible in the obvious way. \finrom \end{lem} Proof. \ref{lemfv2i} This follows from the fact
that the functor $\chi^{-1}$ preserve injectives (Lem. \ref{lemfv1} \ref{lemfv1iii}) and the functor $F^{(-1)}$ as well  (localization).

\ref{lemfv2ii} The isomorphism is clear if $(M,f_M,v_M)$ is flat over $(\smash{E_{/X}^{fv}},\Z)$ and $(N,f_N,v_N)$ is injective over
$(\smash{E_{/X}^{fv}},A_{|X})$. The general case follows by taking resolutions.

\ref{lemfv2iii}  Start with the obvious isomorphism $$\cal Hom^{fv}(\limi (M_.,f_{M_.},v_{M_.}),(N,f_N,v_N))\simeq l_* \cal
Hom^{fv}((M_.,f_{M_.},v_{M_.}),(N,f_N,v_N))$$ Let us explain why the announced isomorphism can be deduced by right derivation. Say that an object
$M_.=(M_.,f_{M_.},v_{M_.})$ in $\smash{Mod^{fv}(E_{/X}^{\N^{op}},\Z)}$ is good if $M_k\hookrightarrow \smash{M_{k'}}$ for $k\le k'$. Every inductive system
$M_.$ is a quotient of a good one e.g. $gd(M_.)\twoheadrightarrow M_.$ where $gd(M_.)_k=\oplus_{j\le k}M_j$. The result will follow if we can prove
that the projective system $\cal Hom^{fv}(M_.,N)$ is $l_*$-acyclic as soon as $M_.$ is good and $N$ is injective. Since
$\chi^{-1}:\smash{Mod^{fv}(E_{/X},A_{|X})^{\N}}\to \smash{Mod(E_{/X},A)^{\N}}$ preserves injectives (Lem. \ref{lemfv1} \ref{lemfv1iii}),
we are thus reduced to prove that for $M_.$ good and $N$ injective $L_.:=\cal Hom(M_.,N)$ is acyclic for the functor $l_*$. In order to check that
$l_*L_.\simeq Rl_*L_.$, it is sufficient to check that $\G(U,l_*L_.)\simeq R\G(U,Rl_*L_.)$ for every $U$ in $E$, or equivalently, $limproj
Hom(M_{.|U},N_{|U})\simeq Rlimproj R\G(U,L_.)$. Since $N$ is injective, we find that $R\G(U,L_.)\simeq Hom(M_{.|U},N_{|U})$ (\cite{SGA4-II} V) and that the transition maps of the latter projective system are surjective. The
result follows (using e.g. the Mittag Leffler criterion for $A(U)$-modules).

\ref{lemfv2iv} The desired morphism may be obtained by using any injective resolution of $A_{|X}$ in $(E_{/X}^{fv},A_{|X})$.
\begin{flushright}$\square$\end{flushright}

\para \label{paranotfv1} In the context of big or $\s$-big crystalline topoi, we will use the following notations:

- For $Y^\s$ in $\cal Sch^\s/\Sigma_1$, we denote $\smash{(Y^\s/\Sigma_\infty)_\s^{fv}}$ the topos denoted $\smash{E_{/X}^{fv}}$ in Lem. \ref{lemfv1} with $E=\smash{(\Sigma_1/\Sigma_\infty)_{CRYS^\s,fl}}$, $X=i_*Y^\s$ and $F$ the endomorphism of $X$ induced by the absolute Frobenius of $Y^\s$.

- For $Y$ in $\cal Sch/\Sigma_1$, we denote $\smash{(Y/\Sigma_\infty)^{fv}}$ the topos denoted $\smash{E_{/X}^{fv}}$ in Lem. \ref{lemfv1} with $E=\smash{(\Sigma_1/\Sigma_\infty)_{CRYS,fl}}$, $X=i_*Y$ and $F$ the endomorphism of $X$ induced by the absolute Frobenius of $Y$.

\begin{lem} \label{DCfv}  Assume that  $X^\s/\Sigma_1$ is locally embeddable.
\debrom \item \label{DCfvi} The inclusion functor $\cal DC_{CRYS^\s,fl}(X^\s)\to Mod(\smash{(X^\s/\Sigma_\infty)^{fv}_\s,\O})$ is $e$-exact and reflects $e$-exactness.

\item \label{DCfvii} For $D_1$, $D_2$ in $\cal DC(X^\s)$, the functor of \ref{DCfvi} induces an injection      $$\xymatrix{Ext^1_{DC_{CRYS^\s,fl}(X^\s)}(D_1,D_2)\ar@{^(->}[r]&Ext^{1}_{(X^\s/\Sigma_\infty)^{fv}_\s,\O}(D_1,D_2)}$$

\item \label{DCfviii} If $X^\s=X$, then \ref{DCfvi} and \ref{DCfvii} hold as well with $(X/\Sigma_\infty)^{fv}_\s$ and $CRYS^\s$ respectively replaced by $(X/\Sigma_\infty)^{fv}$ and $CRYS$.
\finrom
\end{lem}
Proof. \ref{DCfvi} This follows from Rem. \ref{eqDC}, given that the forgetful functor $\chi^{-1}$ is exact. Statement  \ref{DCfvii} follows formally. The case of \ref{DCfviii} is similar.
\begin{flushright}$\square$\end{flushright}

\para \label{paranotfv2} Let us introduce some simplified notations for certain $(f,v)$-modules in crystalline topoi. In what follows, we implicitly use the obvious isomorphism $F^{-1}i_*\simeq i_*F^{-1}$.

 - If $G$ is an Abelian group of $\smash{(\Sigma_1/\Sigma_\infty)_{CRYS^\s,fl}}$, we use the notation
\begin{eqnarray}\label{fvtriv}G=(G,1,1)\end{eqnarray} for the corresponding \emph{trivial $(f,v)$-module of $\smash{((X^\s/\Sigma_\infty)_{CRYS^\s,fl},\Z)}$}, ie. the Abelian group $(G_{|X},f,v)$ of
$\smash{(X^\s/\Sigma_\infty)_\s^{fv}}$,  where $f$ is the trivial identification $G\simeq F^{-1}G$ (not to be
confused with the relative Frobenius) and $v$ is the inverse of $f$. This convention applies for
instance to   $\Gm$, $\mathbb G_m^{log}$, $\Ga$, $\O$, $\cal I$.

-If $G$ is a finite locally free group or a $p$-divisible group over $X$,  then we denote \begin{eqnarray}\label{fvfrel}
G^{fv}=(i_*G,f_G,v_G)\end{eqnarray} the Abelian group of  $\smash{(X^\s/\Sigma_\infty)^{fv}_\s}$, where $v_G$ is induced by the relative Frobenius $\smash{F^{(G/X^\s)}}:G\to
F^{-1}G$ and $f_G$ is induced by $(\smash{F^{(G^*/X^\s)}})^*:F^{-1}G\to G$, where $(-)^*$ denotes the Cartier dual (see section \ref{pdg}). \petit

We use similar notations in the setting of big (instead of $\s$-big) topoi in the case $X^\s=X$.


%

\para \label{paradefbid} Our starting point for the construction of Dieudonn\'e crystals of semi-Abelian schemes is a result from \cite{BBM} for the Dieudonn\'e crystal of a $p$-divisible group.

If $M$ is an Abelian group (resp. $\O$-module) of $\smash{(X/\Sigma_\infty)^{fv}}$, we use the following notation for bidualizing functors $$\begin{array}{rcl}B(M)&=&R\cal Hom^{fv}_{\O_{(X/\Sigma_\infty)}}(\tau_{\le 1}R\cal Hom^{fv}_{(X/\Sigma_\infty)}(M,\O),\O) \hbox{ in $D((X/\Sigma_\infty)^{fv},\O)$}\\
\hbox{(resp. \hspace{.5cm}} B_\O(M)&=&R\cal Hom^{fv}_{\O_{(X/\Sigma_\infty)}}(\tau_{\le 1}R\cal Hom^{fv}_{\O_{(X/\Sigma_\infty)}}(M,\O),\O) \hbox{ in
$D((X/\Sigma_\infty)^{fv},\O)$).}\end{array}$$ Here, and in the following, we denote respectively $\smash{\cal
R\cal Hom^{fv}_{\O_{(X/\Sigma_\infty)}}}$ and $\smash{R\cal Hom^{fv}_{(X/\Sigma_\infty)}}$ the functors (\ref{HomfvA}) and (\ref{Homfvab}).

\begin{rem} \label{remBtrunc} Because the definition of the functor $B$ involves truncation, it does not commute to shifting. Instead we have natural distinguished triangles: $$\begin{array}{ll}\xymatrix{B(M[-1])\ar[r]&B(M)[-1]\ar[r]&R\cal Hom^{fv}_{\O_{(X/\Sigma_\infty)}}(\cal Ext^{fv,2}_{(X/\Sigma_\infty)}(M,\O),\O)[2]\ar[r]^-{+1}&}\\
\xymatrix{B(M)[1]\ar[r]&B(M[1])\ar[r]&R\cal Hom^{fv}_{\O_{(X/\Sigma_\infty)}}(\cal Ext^{fv,1}_{(X/\Sigma_\infty)}(M,\O),\O)[3]\ar[r]^-{+1}&}\end{array}$$
A similar remark holds for $B_\O$.
\end{rem}

\begin{prop} \label{rappelBBM} Let $X$ be in $\cal Sch/\Sigma_1$ and $G$ in $pdiv(X)$.
\debrom
\item \label{rappelBBMi} If $G^*$ denotes the dual $p$-divisible group then $\cal Hom_{(X/\Sigma_\infty)}^{fv}(G^{*,fv},\O)=0$ and
 $$D(G):=\cal Ext^{fv,1}_{(X/\Sigma_\infty)}(G^{*,fv},\O)\hspace{.2cm}\hbox{in $Mod((X/\Sigma_\infty)^{fv},\O)$}$$
is a Dieudonn\'e crystal. The resulting functor $pdiv(X)\to \DC_{CRYS,fl}(X)$ is $e$-exact. \petit

\item \label{rappelBBMii} There are canonical isomorphisms in $Mod((X/\Sigma_\infty)^{fv},\O)$: $$\begin{array}{rcl}D(G)&\simeq &B(G^{fv})[-1]\\
&\simeq & \cal Hom^{fv}_{\O_{(X/\Sigma_\infty)}}(\cal Ext_{(X/\Sigma_\infty)}^{fv,1}(G^{fv},\O),O)\end{array}$$ \finrom
\end{prop}
Proof. \ref{rappelBBMi} follows immediately from Lem. \ref{lemfv2} \ref{lemfv2i} and \cite{BBM} Thm. 3.3.3.

\ref{rappelBBMii} This is proven in \cite{BBM} Prop. 5.3.6 by reduction to the case of finite locally free groups. Let us only recall the definition of the
arrow $B(G^{fv})[-1]\to D(G)$ in question, since it will be needed below. Consider the  exact sequence
$$\xymatrix{0\ar[r]&1+\cal I\ar[r]&\O^\times \ar[r]&\Gm\ar[r]&0}$$
of $\smash{(\Sigma_1/\Sigma_\infty)_{CRYS,fl}}$. Combining with the logarithm $1+\cal I\to \O$ and passing to trivial $(f,v)$-modules of $\smash{(X/\Sigma_\infty)}$, we get a
morphism \begin{eqnarray}\label{deflog}log:\Gm[-1]\to \O\hspace{.2cm}\hbox{ in $D((X/\Sigma_\infty)^{fv})$}\end{eqnarray} According to Lem. \ref{lemfv2} \ref{lemfv2i} and \cite{BBM} Thm. 5.2.7, the morphism (\ref{deflog}) induces an isomorphism
$$R\cal Hom^{fv}_{\O_{(X/\Sigma_\infty)}}(\tau_{\le 1}R\cal Hom^{fv}_{(X/\Sigma_\infty)}(G_{p^.}^{fv},\O),\O)\simeq
 \tau_{\le 0}R\cal Hom^{fv}_{{(X/\Sigma_\infty)}}(\tau_{\le 1}R\cal Hom^{fv}_{(X/\Sigma_\infty)}(G_{p^.}^{fv},\Gm[-1]),\O)$$
in $\smash{D((X/\Sigma_\infty)^{fv,\N},\O)}$. The claimed isomorphisms will follow by computing the effect of $Rl_*$ on both side. We need the following \petit

\noindent \emph{Fact}. There is a natural isomorphism $\smash{limind_k\, \cal Hom^{fv}_{(X/\Sigma_\infty)}(\smash{G_{p^k}^{fv}},\O)}\simeq \smash{\cal Ext^{fv,1}_{(X/\Sigma_\infty)}(G^{fv},\O)}$ and  $\smash{limind_k\, \cal Ext^{fv,1}_{(X/\Sigma_\infty)}}(\smash{G_{p^k}^{fv}},\O)$ vanishes. \petit


\noindent
The proof of this fact is left to the reader (look at realizations of the exact sequence of $Mod(\smash{(X/\Sigma_\infty)^{fv}},\O)$
$$\hspace{-.8cm}\xymatrix{\cal Hom^{fv}_{(X/\Sigma_\infty)}(G_{p^k}^{fv},\O)\ar@{^(->}[r]&\cal E xt^{fv,1}_{(X/\Sigma_\infty)}(G^{fv},\O)\ar[r]^{p^k}&\cal E xt^{fv,1}_{(X/\Sigma_\infty)}(G^{fv},\O)\ar@{->>}[r]&\cal Ext^{fv,1}_{(X/\Sigma_\infty)}(G_{p^k}^{fv},\O)}$$
resulting from \cite{BBM} Thm. 3.3.3 and let $k$ run in $\N^{op}$).

Using this fact, we obtain the following series of isomorphisms:
\begin{eqnarray} \label{isoko1}Rl_*(LHS)
&\simeq & R\cal Hom^{fv}_{\O_{(X/\Sigma_\infty)}}(\limi \tau_{\le 1}R\cal Hom^{fv}_{(X/\Sigma_\infty)}(G_{p^.}^{fv},\O),\O)\\
\label{isoko2} &\simeq &R\cal Hom^{fv}_{\O_{(X/\Sigma_\infty)}}(\cal Ext^{fv,1}_{(X/\Sigma_\infty)}(G^{fv},\O),\O)\\
\label{isoko3} &\simeq & \cal Hom^{fv}_{\O_{(X/\Sigma_\infty)}}(\cal Ext^{fv,1}_{(X/\Sigma_\infty)}(G^{fv},\O),\O) \end{eqnarray}
where (\ref{isoko1}) is by  Lem. \ref{lemfv2} \ref{lemfv2iii} and (\ref{isoko3}) is by local freeness of $\smash{\cal Ext^{1}_{(X/\Sigma_\infty)}(G,\O)}$ (see \ref{rappelBBMi}). Having noticed that $Rl_*(LHS)$ is concentrated in degree $0$, we find that
\begin{eqnarray} \label{isoko4} Rl_*(RHS)&\simeq &R^0l_*(RHS)\\
\label{isoko5}&\simeq &R^0l_*\tau_{\le 0}R\cal Hom^{fv}_{{(X/\Sigma_\infty)}}(\tau_{\le 1}R\cal Hom^{fv}_{(X/\Sigma_\infty)}(G_{p^.}^{fv},\Gm[-1]),\O)\\
\label{isoko6}&\simeq & R^0l_*\tau_{\le 0}R\cal Hom^{fv}_{{(X/\Sigma_\infty)}}(G_{p^.}^{*,fv}[-1],\O)\\
\label{isoko7}&\simeq & R^1l_*\tau_{\le 1}R\cal Hom^{fv}_{{(X/\Sigma_\infty)}}(G_{p^.}^{*,fv},\O)\\
\label{isoko8}&\simeq & R^1l_*R\cal Hom^{fv}_{{(X/\Sigma_\infty)}}(G_{p^.}^{*,fv},\O)\\
\label{isoko9}&\simeq & \cal Ext^{fv,1}_{{(X/\Sigma_\infty)}}(G^{*,fv},\O)\end{eqnarray}
where (\ref{isoko4}) - (\ref{isoko8}) are trivial and (\ref{isoko9}) is by Lem. \ref{lemfv2} \ref{lemfv2iii}.
\begin{flushright}$\square$\end{flushright}

\begin{defn} \label{defD1mot} Let  $X$ in $\cal Sch/\Sigma_1$ and $M$ in $\M(X)$. We set $$D_X(M):=D(M_{p^\infty})$$
where $\smash{M_{p^\infty}}$ is the $p$-divisible group associated to $M$  (Lem. \ref{pdivfunct} \ref{pdivfunctii}). This defines a functor $D_X:\cal M(X)\to \cal DC_{CRYS,fl}(X)$, which is canonically pseudo functorial with respect to $X$ and $e$-exact if $X$ is regular (Lem. \ref{1motex}).
\end{defn}

Another basic input is the following calculation from \cite{BM}.

\begin{prop} \label{propDZ01} Let $X$ in $\cal Sch/\Sigma_1$.
\debrom
\item \label{propDZ01i}
There are canonical isomorphisms in $\DC_{CRYS,fl}(X)$: $$\begin{array}{rcccl}D_X(\Z)&:=&D(\Qp/\Zp)&\simeq &(\O,p,1)\\
D_X(\Z(1))    &:=&D(\mu_{p^\infty})&\simeq &(\O,1,p)\end{array}$$

\item \label{propDZ01ii} Assume that $X$ is regular.
There is a commutative diagram of Abelian groups
{\footnotesize$$\hspace{-1cm}\xymatrix{Ext^1_{\cal
M(X)}(\Z,\Z(1))\ar[d]_{|D_X|}&\ar[l]^\sim_-\delta\Gm(X)\ar[r]^-{nat}&Hom_{(X/\Sigma_\infty)^{fv}}((\Z,p,1),(\Gm,1,p))
\ar[d]^{log}\\
Ext^1_{\DC_{CRYS,fl}(X)}(D_X(\Z),D_X(\Z(1)))\ar[r]^-{\ref{propDZ01i}}_-{\ref{DCfv}}& Ext^{1}_{(X/\Sigma_\infty)^{fv},\O}((\O,p,1),(\O,1,p))&
\ar[l]^-\sim_-{adj}Ext^{1}_{(X/\Sigma_\infty)^{fv}}((\Z,p,1),(\O,1,p)) }$$} \noindent where $\delta$ denotes the obvious isomorphism $x\mapsto [Kum(x)]$ and $|D_X|$ denotes the map induced by the $e$-exact functor $D_X$ defined in
Def. \ref{defD1mot}.
\finrom
\end{prop}

%



Proof. \ref{propDZ01i} Since the structural ring $\O$ of $\smash{(X/\Sigma_\infty)_{CRYS,fl}}$ is $p$-torsion, the natural morphism $(\Qp/\Zp)^{fv}=(\Qp/\Zp,p,1)\to (\Z,p,1)[1]$ of $D((X/\Sigma_\infty)^{fv})$  induces isomorphisms  $$\begin{array}{rclll}\O\Ltens (\Qp/\Zp)^{fv}&\simeq& (\O,p,1)[1]&&\\
R\cal Hom^{fv}_{(X/\Sigma_\infty)}((\Qp/\Zp)^{fv},\O)&\simeq& (\O,1,p)[-1]&&\hbox{(Lem. \ref{lemfv2} \ref{lemfv2ii})}\end{array}$$
in $\smash{D((X/\Sigma_\infty)^{fv},\O)}$.
The announced computation of $D_X(\Z)$ and $D_X(\Z(1))$ follows immediately (using Prop. \ref{rappelBBM} \ref{rappelBBMii} for the first one):

\begin{eqnarray}\label{explicitZ01}\begin{array}{rcl}D_X(\Z)&:=&D(\Qp/\Zp)\\
& \simeq& B((\Qp/\Zp)^{fv})[-1]\\
& \simeq & \smash{R\cal Hom^{fv}_{\O_{(X/\Sigma_\infty)}}((\O,1,p)[-1],\O)}[-1]\\
&\simeq & (\O,p,1)
\\ D_X(\Z(1))&:=& D(\mu_{p^\infty})\\
& := & \smash{\cal Ext^{fv,1}_{(X/\Sigma_\infty)}((\Qp/\Zp)^{fv},\O)}\\
& \simeq & (\O,1,p)\end{array} \end{eqnarray} The following claim will be needed for the proof of \ref{propDZ01ii}. \petit

\emph{Claim}. The above isomorphism $D(\mu_{p^\infty})\simeq (\O,1,p)$ coincides with
$$\xymatrix{D(\mu_{p^\infty})\ar[r]^-{\ref{rappelBBM}\ref{rappelBBMii}} &B((\mu_{p^\infty})^{fv})[-1]\ar[r]^-{log}&B_\O((\O,1,p))& \ar[l]_-\sim (\O,1,p)}$$
where $log$ is obtained by applying $B_\O$ to the morphism $\O\otimes^L(\mu_{p^{\infty}})^{fv}\to (\O,1,p)[1]$ of $D((X/\Sigma_\infty)^{fv},\O)$ deduced from (\ref{deflog}) (apply Rem. \ref{remBtrunc} with $M=(\O,1,p)$). \petit

Let us prove the claim. To begin with, we notice that for any module or Abelian group $M$ we have a natural morphism $$\Qp/\Zp\Ltens \tau_{\le i}R\cal
Hom^{fv}_{(X/\Sigma_\infty)}((\mu_{p^\infty})^{fv},M)\to \limi \tau_{\le i}R\cal Hom^{fv}_{(X/\Sigma_\infty)}((\mu_{p^.})^{fv},M)$$ (since $\smash{\mu_{p^.}\simeq
\Z/p^.\otimes^L\mu_{p^\infty}[-1]}$) which is an isomorphism if and only if $$\cal Tor_1(\Qp/\Zp,\smash{\cal Ext^{fv,i+1}_{(X/\Sigma_\infty)}((\mu_{p^\infty})^{fv},M)})=0$$


\noindent
This is in particular the case if $M=\O$, $i=1$ (\cite{BBM} Thm. 3.3.3 (iii)) or  $M=\Gm$, $i=0$. To see that, use Lem. \ref{lemfv2} \ref{lemfv2i} together with the vanishing of  $\smash{\cal Ext^{1}_{(X/\Sigma_\infty)}(\cal \mu_{p^.},\Gm)}$ (\cite{BBM} Cor. 1.1.8 and  \cite{SGA7-I} VIII, Prop. 3.3.1) and the vanishing of $\Z/p^k\otimes^LR^1l_*\smash{\cal Hom_{(X/\Sigma)}(\cal \mu_{p^.},\Gm)}$ (Example \ref{exRl}).

This observation explains the second commutative square in the following diagram
{\scriptsize$$\hspace{-1.2cm}\xymatrix{B((\mu_{p^\infty})^{fv})[-1]\ar[d]^\wr\ar[rr]_-\sim^-{\ref{rappelBBM} \ref{rappelBBMii}}&& D(\mu_{p^\infty})\ar[d]_\wr\\
R\cal Hom^{fv}_{\O_{(X/\Sigma_\infty)}}(\limi \tau_{\le 1}R\cal Hom^{fv}_{{(X/\Sigma_\infty)}}((\mu_{p^.})^{fv},\O),\O)\ar[rr]^-{log:\Gm[-1]\to \O}&&R\cal Hom^{fv}_{{(X/\Sigma_\infty)}}(\limi \cal Hom^{fv}_{{(X/\Sigma_\infty)}}((\mu_{p^.}^{fv},\Gm),\O)[1]\\
R\cal Hom^{fv}_{\O_{(X/\Sigma_\infty)}}(\Qp/\Zp\Ltens\tau_{\le 1}R\cal Hom^{fv}_{{(X/\Sigma_\infty)}}((\mu_{p^\infty})^{fv},\O),\O)\ar[rr]^-{log:\Gm[-1]\to \O}\ar[u]_\wr\ar[d]^{log:(\mu_{p^\infty})^{fv}\to (\O,1,p)[1]}&&R\cal Hom^{fv}_{{(X/\Sigma_\infty)}}(\Qp/\Zp \Ltens \cal Hom^{fv}_{{(X/\Sigma_\infty)}}((\mu_{p^\infty})^{fv},\Gm),\O)[1]\ar[u]^\wr\ar[d]_\wr\\
R\cal Hom^{fv}_{{(X/\Sigma_\infty)}}(\Qp/\Zp\Ltens \cal Hom^{fv}_{\O_{(X/\Sigma_\infty)}}((\O,1,p),\O))[1]\ar[rr]^-\sim \ar[d]^\wr&& R\cal Hom^{fv}_{{(X/\Sigma_\infty)}}(\Qp/\Zp\Ltens (\Z,p,1),\O)[1]\ar[d]_\wr\\
B_\O((\O,1,p))&&\ar[ll]_\sim (\O,1,p)}$$} Let us explain the other commutative squares. The first commutative square is clear, from the construction of the
isomorphism Prop. \ref{rappelBBM} \ref{rappelBBMii}. The fourth one is obvious. To get the third one, we notice that each corner is in fact concentrated in
degree $0$ so that we may drop the superscripts $fv$ to check commutativity. Then, we are done using the following commutative square of Abelian
groups in $\smash{(X/\Sigma_\infty)_{CRYS,fl}}$:
$$\xymatrix{\cal Ext^1_{(X/\Sigma_\infty)}(\mu_{p^\infty},\O)&&&\ar[lll]_-{log:\Gm\to \O[1]}\cal Hom_{(X/\Sigma_\infty)}(\mu_{p^\infty},\Gm)\\
\cal Hom_{\O_{(X/\Sigma_\infty)}}(\O,\O)\ar[u]_{log:\mu_{p^\infty}[-1]\to \O}&&&\ar[lll]\Z\ar[u]}$$
The claim is now proven. \petit

\ref{propDZ01ii} In the diagram below, we let \ref{propDZ01i} denote the arrow induced by the isomorphisms established in \ref{propDZ01i}.
{\scriptsize$$\hspace{-3.4cm}\xymatrix{&&Ext^1_{\cal M(X)}(\Z,\Z(1))\ar[d]^{(-)_{p^\infty}}\ar[ld]^-{|D|}&\\ &Ext^1_{\DC_{CRYS,fl}(X)}(D_X(\Z),D_X(\Z(1)))\,\, \,\, \,\, \,\, \,\,
\,\, \,\, \ar  @<-0pt> `l[d] `[dddd]^(.75){(i)} [ddddr] \ar[d]&Ext^{1}_{(X/\Sigma_\infty)^{fv}}((\Qp/\Zp)^{fv},(\mu_{p^\infty})^{fv})\ar[ld]^-{|B|}&
\Gm(X)\ar[ul]^-\delta\ar[d]^-{nat}\\
&\,\, \,\, \,\, \,\, \,\,
\,\, \,\, \,\, \,\, \,\,
Ext^{1}_{{(X/\Sigma_\infty)^{fv}},\O}(B((\Qp/\Zp)^{fv}),B((\mu_{p^\infty})^{fv}))\ar[d]^-{log}
\ar@{->}[dr]&&
Hom_{(X/\Sigma_\infty)^{fv}}((\Z,p,1),(\Gm,1,p))\ar[d]^-{log}\ar[dl]_-{|B|}\ar[ul]^-{\Qp/\Zp\Ltens(-)}\\
& \,\, \,\, \,\, \,\, \,\,
\,\, \,\, \,\, \,\, \,\, Ext^{2}_{{(X/\Sigma_\infty)^{fv}},\O}(B((\Qp/\Zp)^{fv}),B_\O(\O,1,p))\ar@{->}[dr]&
Hom_{{(X/\Sigma_\infty)^{fv}},\O}(B(\Z,p,1),B(\Gm,1,p))\ar[d]^-{log}&
Ext^{1}_{(X/\Sigma_\infty)^{fv}}((\Z,p,1),(\O,1,p)) \ar[dl]^{|B|}\ar@<-0pt>
`d[dd]^(.75){adj} [ddl]\\
&&Ext^{1}_{{(X/\Sigma_\infty)^{fv}},\O}(B(\Z,p,1),B_\O(\O,1,p))\ar[d]^\wr&\\
&&Ext^{1}_{{(X/\Sigma_\infty)},\O}((\O,p,1),(\O,1,p))&}$$}

The claim established in \ref{propDZ01i} implies that the pentagon on the left is commutative. The commutativity of the rest of the diagram causes no difficulty and will be left to the reader. Let us only give some explanations concerning the definition of some arrows. The arrow denoted $(-)_{p^\infty}$ is well defined thanks to Lem. \ref{pdivfunct} \ref{pdivfunctii} and Cor. \ref{corpdivtop} \ref{corpdivtopii}. The arrow denoted $nat$ follows from the following observation: for $x:\Z\rightarrow \Gm$ the morphism $F^{-1}x:F^{-1}\Z\to F^{-1}\Gm$ translates as
$1\mapsto x(1)^p$ if one makes the obvious identifications $F^{-1}\Z=\Z$ and $F^{-1}\Gm=\Gm$. The arrows denoted $log$ are deduced from $\smash{\O\otimes^L(\Gm)^{fv}\to (\O,1,p)[1]}$ by applying $B_\O$ and Rem. \ref{remBtrunc}. The first and third arrow denoted $|B|$ are defined using respectively the isomorphisms $B((\Qp/\Zp[-1]))\simeq B(\Qp/\Zp)[-1]$ (by Rem. \ref{remBtrunc} and \cite{BBM} Thm. 3.3.3, (iii)) and $B((\Z,p,1)[-1])\simeq B(\Z,p,1)[-1]$ (by Rem. \ref{remBtrunc} again).
\begin{flushright}$\square$\end{flushright}



\begin{prop} \label{prolD1} Let $X$, $R$, $K$, $k$ as in Sect. \ref{subsectionog1mot}. Consider furthermore $X^\s=(X,Spec(k))$ and the natural morphism $o:X^\s\rightarrow X$. Denote $|D_X|$ and $|D_K|$ the homomorphisms on $Ext^1$'s induced by the functors defined by Prop. \ref{rappelBBM} \ref{rappelBBMi} and  Rem. \ref{eqDC} \ref{eqDCii}  for $X$ and $Spec(K)$ respectively. There exists a canonical homomorphism $\smash{|D_{X^\s}|}$ rendering the following diagram commutative.
{\footnotesize$$\hspace{-2cm}\xymatrix{Ext^1_{\M(X)}(\Z,\Z(1))\ar[r]\ar[d]^{|D_X|}& Ext^1_{\M_{log}(X)}(\Z,\Z(1))\ar[r]\ar @{-->}[d]^{|D_{X^\s}|}&Ext^1_{\M(K)}(\Z,\Z(1))\ar[d]^{|D_K|}\\
Ext^1_{\DC_{CRYS^\s,fl}(X)}(D_X(\Z),D_X(\Z(1)))\ar[r]&Ext^1_{\DC_{CRYS^\s,fl}(X^\s)}(o^*D_X(\Z),o^*D_X(\Z(1))) \ar[r]&Ext^1_{\DC_{CRYS^\s,fl}(K)}(D_K(\Z),D_K(\Z(1)))}$$}
\end{prop}
Proof. Using Prop. \ref{propDZ01} \ref{propDZ01ii}, we find that it is enough to define dotted arrows (the definition of the other ones is analogous and compatible via $\pi^*$ with their counterpart in the context of big crystalline sites) rendering the following
diagram commutative and to check that the composed vertical dotted arrow takes its values into
${Ext^1_{\DC_{CRYS^\s,fl}(X^\s)}((\O,p,1),(\O,1,p))}$ viewed as a subgroup of ${Ext^{1}_{(X^\s/\Sigma_\infty)^{fv}_\s}((\O,p,1),(\O,1,p))}$ (Lem. \ref{DCfv} \ref{DCfvii}).

{\scriptsize$$\hspace{-2cm}\xymatrix{Ext^1_{\M(X)}(\Z,\Z(1))\ar @/_3pc/@<-15ex>[ddddd]_{|D_X|} \ar[r] &Ext^1_{\M_{log}(X)}(\Z,\Z(1))\ar[r]^{\sim}&Ext^1_{\M(K)}(\Z,\Z(1))\ar @/^3pc/@<+15ex>[ddddd]^{|D_K|}\\
\Gm(X)\ar[u]^{\wr}_{\delta_X}\ar[r]\ar[d]^{nat_X}&\G(X^\s,M_{X^\s}^{gp})\ar@{-->}[u]^{\wr}_{\delta^{log}_{X^\s}} \ar[r]^{\sim}
\ar@{-->}[d]^{nat_{X^\s}}&\Gm(K)\ar[u]^{\wr}_{\delta_K}\ar[d]^{nat_K}\\
Ext^{0}_{(X/\Sigma_\infty)^{fv}_\s}((\Z,p,1),(\Gm,1,p))\ar[r]\ar[d]^{log_X}& Ext^{0}_{(X^\s/\Sigma_\infty)^{fv}_\s}((\Z,p,1),(M_X^{gp},1,p))\ar[r]\ar@{-->}[d]^{log_{X^\s}}
&Ext^{0}_{(K/\Sigma_\infty)^{fv}_\s}((\Z,p,1),(\Gm,1,p))\ar[d]^{log_K}\\
Ext^{1}_{(X/\Sigma_\infty)^{fv}_\s}((\Z,p,1),(\O,1,p))\ar[r]\ar[d]^{adj_X}& Ext^{1}_{(X^\s/\Sigma_\infty)^{fv}_\s}((\Z,p,1),(\O,1,p))\ar[r]\ar@{-->}[d]^{adj_{X^\s}}
&Ext^{1}_{(K/\Sigma_\infty)^{fv}_\s}((\Z,p,1),(\O,1,p))\ar[d]^{adj_K}\\
Ext^{1}_{(X/\Sigma_\infty)^{fv}_\s,\O}((\O,p,1),(\O,1,p))\ar[r]\ar[d]_{\wr}& Ext^{1}_{(X^\s/\Sigma_\infty)^{fv}_\s,\O}((\O,p,1),(\O,1,p))\ar[r]
&Ext^{1}_{(K/\Sigma_\infty)^{fv}_\s,\O}((\O,p,1),(\O,1,p))\ar[d]_{\wr}\\
Ext^{1}_{(X/\Sigma_\infty)^{fv}_\s,\O}(D_X(\Z),D_X(\Z(1)))& &Ext^{1}_{(K/\Sigma_\infty)^{fv}_\s,\O}(D_K(\Z),D_K(\Z(1))) }$$}

\noindent Here $\smash{M_{X^\s}^{gp}}$ denotes the sheaf of groups on $\smash{FL^\s(X^\s)}$ defined by the collection of sheaves $\smash{M_{U^\s}^{gp}}$ (the fraction group of the structural monoid on $et(U^\s)$) together with the natural functoriality morphisms for $U^\s$ varying in $\cal Sch^\s/X^\s$. That this is indeed a sheaf, is explained in the proof of \cite{Ka3} Thm. 3.2 and relies essentially on the fact that small \'etale sheaves (and thus in particular morphisms of such) satisfy $fl$ descent (\cite{SGA4-II} VIII, Thm. 9.4). As usual, $\smash{M_{X^\s}^{gp}}$ is viewed as a sheaf on $\smash{CRYS^\s_{fl}(X^\s/\Sigma_\infty)}$ via $i_*$. It should not be confused with $\smash{M_{{X^\s/\Sigma_\infty}}^{gp}}: (U^\s,T^\s)\mapsto \smash{\G(T^\s,M_{T^\s}^{gp})}$ (which is also a sheaf for similar reasons). Both are related by the following exact sequence: \begin{eqnarray}\label{secMcrys}\xymatrix{0\ar[r]&1+\cal I_{X^\s/\Sigma_\infty}\ar[r]&M_{{X^\s/\Sigma_\infty}}^{gp}\ar[r]&M_{X^\s}^{gp}\ar[r]&0}\end{eqnarray}
We define the desired dotted arrows as follows.

%


\noindent - The morphism $\smash{\delta^{log}_{X^\s}}$ is defined so that the top right square  (and thus also the top left one)  is commutative.

\noindent - The morphisms $\smash{nat_{X^\s}}$ and $\smash{adj_{X^\s}}$ are the obvious ones.

\noindent - The arrow $\smash{log_{X^\s}}$ is deduced from (\ref{secMcrys})
and the logarithm $\smash{1+\cal I_{X^\s/\Sigma_\infty}}\rightarrow \smash{\cal O_{X^\s/\Sigma_\infty}}$.

It now remains to check that if an extension class $[E]$ is produced by the dotted vertical arrow from an element of $\smash{\G(X^\s,M_{X^\s}^{gp})}$, then $E$ is in
fact in $\smash{\cal DC_{CRYS^\s,fl}(X^\s)}$. First, we note that $E$ is a locally free crystal since it is an extension of two copies of $\O$. Finally, we notice that the
condition $fv=p$, $vf=p$ can be checked after pulling back to $\smash{((Spec(K)/\Sigma_\infty)_{CRYS^\s,fl},\O)}$ (observe realizations at a lifting).
\begin{flushright}$\square$\end{flushright}

\para \label{deflogDC}
We are now in a position to define the Dieudonn\'e functor for log $1$-motives and semi-Abelian schemes over $X$. We use the following notations $$\xymatrix{Spec(K)\ar[r]_-{j^\s}\ar
@<-0pt> `u[r] `[rr]^-j [rr] &X^\s\ar[r]_-{o}&X}$$ for the obvious morphisms.

\begin{cor} \label{defDMlog} There exists a canonical $e$-exact functor $\smash{D_{X^\s}}$ fitting into a canonically pseudo-commutative diagram as follows: $$\xymatrix{\cal M(X)\ar[r]^-\subset \ar[d]_{D_{X}}&\cal M_{log}(X)\ar[r]^-{j^{-1}}\ar[d]_-{D_{X^\s}}&\cal M(K)\ar[d]_-{D_K}
\\
\DC(X)\ar[r]^-{o^*}&\DC(X^\s)\ar[r]^-{j^{\s,*}}&\DC(K)}$$ This diagram is  canonically pseudo-functorial with respect to finite \'etale
base change.
\end{cor}
Proof. This directly results from Prop. \ref{prol} and Prop. \ref{prolD1}.
\begin{flushright}$\square$\end{flushright}


\begin{defn} \label{defSASloc} Consider a semi-Abelian scheme $A/X$.

\debrom \item \label{defSASloci} The  \emph{log $1$-motive of
$A$} and the \emph{$1$-motive of $A_{|\eta}$} are defined respectively as $$\begin{array}{rcl}M_{log}(A)&:=&(\G,G,f)\hbox{\hspace{.5cm} in $\cal M_{log}(X)$}
\\
M(A_{|\eta})&:=&(\G_{|\eta},G_{|\eta},f)\hbox{\hspace{.5cm} in $\cal M(\eta)$}\end{array}$$ where $\G$, $G$ (resp. $f$) are as in Cor. \ref{rig} (resp. Lem. \ref{remrigqff} \ref{remrigqffi}).

\item \label{defSASlocii} The \emph{Dieudonn\'e crystal of $A/X$} is defined as follows:  $$D_{X^\s}(A):=D_{X^\s}(M_{log}(A/X))$$
\finrom
\end{defn}

Let us write down explicitly how the Dieudonn\'e crystals of $A$, $G$ and $\G$ are related.

\begin{prop} \label{lemDCloc} \debrom \item \label{lemDCloci} The weight filtration on $M_{log}(A)$ induces a canonical exact sequence as follows in $\DC(X^\s)$: $$\xymatrix{0\ar[r]&D_{X^\s}(0,G,0)\ar[r]&D_{X^\s}(M_{log}(A))\ar[r]&D_{X^\s}(\G,0,0)\ar[r]&0}$$

\item \label{lemDClocii} The weight filtration of $M(A_{|\eta})$ and the bottom exact sequence of Rem. \ref{remrigqffiii} induce isomorphic exact sequences as follows in $\DC(\eta)$:
    $$\xymatrix{0\ar[r]&D_{\eta}(0,G_{|\eta},0)\ar[r]\ar[d]^\wr&D_{\eta}(M(A_{|\eta}))\ar[r]\ar[d]^\wr&D_{\eta}(\G_{|\eta},0,0)\ar[r]\ar[d]^\wr&0\\
    0\ar[r]&D_{\eta}(G_{|\eta,p^\infty})\ar[r]&D_{\eta}(A_{|\eta,p^\infty})\ar[r]&D_{\eta}(\Qp/\Zp\otimes \G_{|\eta})\ar[r]&0}$$
The top exact sequence is naturally isomorphic to the one deduced from \ref{lemDCloci} by $j^{\s,*}:\DC(X^\s)\to \DC(\eta)$.
\finrom
\end{prop}
Proof. \ref{lemDCloci} This follows from the fact that $\smash{D_{X^\s}}:\cal M_{log}(X)\to \DC(X^\s)$ is exact.

\ref{lemDClocii} The compatible exact sequences are deduced from Lem. \ref{remrigqff} \ref{remrigqffii} by applying the $e$-exact functor $D_{\eta}:pdiv(\eta)\to \DC(\eta)$. The last assertion is by compatibility of the functors $D_{X^\s}$ and $D_{\eta}$ (Cor. \ref{defDMlog}).
\begin{flushright}$\square$\end{flushright}

\para
We are now in a position to define the Dieudonn\'e crystal of a semi-Abelian scheme  over $C$.
\begin{defn} \label{defDCglob} Let $SAS(C,Z)$ denote the category of semi-Abelian schemes over $C$ which are Abelian above $U$. We define a functor $$D_{C^\s}:SAS(C,Z)\longrightarrow \DC(C^\s)$$ by sending $A/C$ to the Dieudonn\'e crystal over $(C^\s/\Zp)$ corresponding to $$m^*D_{C^\s}(A/C):=(D_{C_v^\s}(A_{|C_v})\rightarrow D_{U_v}(A_{|U_v})\leftarrow D_U(A_{|U}))_{v\in Z}\hbox{ in $\DC(J^\s)$}$$ under the equivalence of Prop. \ref{MVsyn} \ref{MVsyni}.
\end{defn}
Note that this definition uses the isomorphism $\smash{j_v^{\s,*}D_{C_v^\s}(A_{|C_v})}\simeq \smash{D_{U_v}(A_{|U,v})}$ of Prop. \ref{lemDCloc}
\ref{lemDClocii}. The following result explains the relation with the Dieudonn\'e crystal arising from the diagram of $p$-divisible groups
introduced in Def. \ref{defH}.

\begin{lem} \label{lemDCH} Consider the restriction $H_{|J}$ to $J$ of the diagram of $p$-divisible groups $H$ defined in Def. \ref{defH}. We have a canonical exact sequence as follows in $\DC(J^\s)$: $$\xymatrix{0\ar[r]&o_{J}^*D_J(H_{|J})\ar[r]&m^*D_{C^\s}(A/C)\ar[r]&(o_{J}^*D_J((\Qp/\Zp\otimes \G_v\rightarrow 0\leftarrow 0)_{v\in Z})\ar[r]&0}$$
\end{lem}
Proof. This results directly from Prop. \ref{lemDCloc}.
\begin{flushright}$\square$\end{flushright}

%
%
%


\section{The comparison theorem} \label{sectct}

\subsection{Diagrams of $p$-divisible groups} \label{dopdg} ~~ \\


The purpose of this section is to establish the comparison theorem for $p$-divisible groups over a diagram of schemes whose vertices have local finite $p$-bases.
This will be achieved in Lem. \ref{lemcomp2} for $\Qp/\Zp$ and $\mu_{p^\infty}$, and in Thm. \ref{thmcomppdiv} for arbitrary $p$-divisible groups, using Cartier biduality.

\para Let us begin with technical observations regarding the computation of $\cal Ext^q$'s on diagrams.

\begin{lem} \label{lemext1} Consider a ringed variable topos $(\cal T,A)$ on a category $\cal B$ and let $X/\Delta$ in $Diag(\cal B)$.
\debrom \item \label{lemext1i}
Consider an Abelian group $M$ of $\cal T(X)$ and a module $N$ over $(\cal T(X),A_X)$. If $M$ is cocartesian (ie. $f^{-1}M_\delta\simeq M_{\delta'}$ for all $f:\delta'\to \delta$ in $\Delta$) then the following natural morphism of $\smash{D^+(\cal T(X_\delta),A_{X_\delta})}$ is invertible for all $\delta\in \Delta$: $$\xymatrix{R\cal Hom_{\cal T(X)}(M,N)_\delta\ar[r]& R\cal Hom_{\cal T(X_\delta)}(M_\delta,N_\delta)}$$


\item \label{lemext1iii} Consider $\cal Hom_{\cal T(X)}$ as a bifunctor $Ab(\cal T(X))^{\N^{op}}\times Mod(\cal T(X),A_X)\to Mod(\cal T(X),A_X)^{\N}$. This bifunctor is right derivable into $D^-(\cal T(X))^{\N^{op}})\times D^+(\cal T(X),A_X)\to D^+(\cal T(X),A_X)^{\N}$  and there is a bifunctorial isomorphism $$R\cal Hom_{\cal T(X)}(\limi M_k,N)\simeq Rl_*\cal Hom_{\cal T(X)}(M_.,N)$$
\finrom
\end{lem}
Proof. \ref{lemext1i}
Let us compute the $A_X$-module $\cal Hom_{\cal T(X)}(M,N)$. Consider $\delta$ in $\Delta$, $U$ in $\cal T(X_\delta)$ and let $\cal T(X/U)$ denote the topos of sections of the cofibered category over $\Delta/\delta$ whose fiber at $g:\delta'\to \delta$ is $\cal T(X_{\delta'})/g^{-1}U$. Let $h:\cal T(X/U)\to \cal T(X)$ denote the morphism induced by $\Delta/\delta\to \Delta$, $g\mapsto \delta'$ and the localization morphisms $\cal T(X_{\delta'})/g^{-1}U\to \cal T(X_{\delta'})$. Using \cite{SGA4-II} VI, Prop. 7.4.7, we find a natural isomorphism  $$\cal Hom_{\cal T(X)}(M,N)_\delta(U)\simeq Hom_{\cal T(X/U)}(h^{-1}M,h^{-1}N)$$  Via this identification, the image by $H^0$ of the natural morphism in question translates into $$\xymatrix{Hom_{\cal T(X/U)}(h^{-1}M,h^{-1}N)\ar[r] & Hom_{\cal T(X_\delta)}(M_{\delta,|U},N_{\delta,|U})}$$
and is thus an isomorphism, since $M$ is cocartesian. The case of derived functors follows, since $N_\delta$ is flasque if $N$ is injective (Lem.-Def. \ref{acycf} \ref{acycfv}).


\ref{lemext1iii} Right derivability causes no difficulty (use injective resolutions of the second argument). The proof of the claimed isomorphism is similar to (and easier than) the proof of Lem. \ref{lemfv2} \ref{lemfv2iii}.
\begin{flushright}$\square$\end{flushright}

%


Let us now gather some results about the behaviour of $\cal Ext^i$'s while traveling through the topoi which will be involved in our proof. The next statement and later ones implicitly use Cor. \ref{corpdivtop} to switch between the various incarnations of a $p$-divisible group.

\begin{lem}\label{lemext2} Let $X$ in $\cal Sch/\Sigma_1$ and $\Sigma=\Sigma_k$, $1\le k\le \infty$. Let $fl\preceq top\preceq zar$  (e.g. $top=syn$) and let $\epsilon$ denote the weak morphism from $fl$ to $top$.
\debrom
\item \label{lemext2i} Let $G/X$ denote  a quasi-compact quasi-separated group (resp. a $p$-divisible group viewed as a group in $(X/\Sigma)_{CRYS,fl}$). Consider a quasi-coherent crystal $M$ of $((X/\Sigma)_{CRYS,fl},\O)$ (resp. Consider a quasi-coherent crystal $M$ of $((X/\Sigma)_{CRYS,fl},\O)$ and assume that $top\preceq syn$). Then for $0\le i\le 2$ and $n\ge 0$: $$\epsilon_*\cal Ext^i_{(X/\Sigma)_{CRYS,fl}}(G,\cal I^{[n]}M)\simeq \cal Ext^i_{(X/\Sigma)_{CRYS,top}}(\epsilon_*G,\epsilon_*\cal I^{[n]}M)$$
\item \label{lemext2ii} Assume $k<\infty$. Let $G/X$ denote a $p$-divisible group viewed as a group in $X_{SYN}$. If $M$ is a module of $((X/\Sigma)_{CRYS,syn},\O)$ then for $0\le i\le 1$: $$u_*\cal Ext^i_{(X/\Sigma)_{CRYS,syn}}(u^{-1}G,M)\simeq \cal Ext^i_{X_{SYN}}(G,u_*M)$$
    The same is true for $syn$ and $crys$  instead of  $SYN$ and $CRYS$.
\item \label{lemext2iii} Consider a syntomic group scheme (resp. a $p$-divisible group) $G/X$ and let $0\le i\le 2$ (resp. $0\le i\le 1$).

    - If we view $G$ as a group in $X_{SYN}$ and if $M$ is any Abelian group of $X_{SYN}$ (resp. an Abelian group of $X_{SYN}$ which is killed by a power of $p$) then: $$\pi_*\cal Ext^i_{X_{SYN}}(G,M)\simeq \cal Ext^i_{X_{syn}}(\pi_*G,\pi_*M)$$

    - If we view $G$ as a group of $(X/\Sigma)_{CRYS,syn}$ and if $M$ is an Abelian group of $(X/\Sigma)_{CRYS,syn}$  satisfying $M\simeq i_*i^{-1}M$ (resp. an Abelian group of $(X/\Sigma)_{CRYS,syn}$  which is killed by a power of $p$ and satisfies $M\simeq i_*i^{-1}M$) then: $$\pi_*\cal Ext^i_{(X/\Sigma)_{CRYS,syn}}(G,M)\simeq \cal Ext^i_{(X/\Sigma)_{crys,syn}}(\pi_*G,\pi_*M)$$

\item \label{lemext2iv} Let $G/X$ denote a $p$-divisible group viewed as an Abelian group of $(X/\Sigma)_{CRYS,top}$. If $M$ is a quasi-coherent crystal of $((X/\Sigma)_{CRYS,top},\O)$ then for $0\le i\le 2$, $n\ge 0$: $$\cal Ext^i_{(X/\Sigma)_{CRYS,top}}(G,\cal I^{[n]}M)\simeq \limp_j\cal Ext^i_{(X/\Sigma)_{CRYS,top}}(G_{p^j},\cal I^{[n]}M)$$

\item \label{lemext2v} Let $G/X$ denote a $p$-divisible group viewed as an Abelian group of $X_{TOP}$. If $top\preceq syn$ and if $M$ is an Abelian group of $X_{TOP}$ which is killed by a power of $p$ then for $0\le i\le 1$: $$\cal Ext^i_{X_{TOP}}(G,M)\simeq \limp_j \cal Ext^i_{X_{TOP}}(G_{p^j},M)$$ The same is true with $top$ instead of $TOP$.
\finrom
\end{lem}

Proof. \ref{lemext2i} In the first case, the proof of \cite{BBM} Cor. 2.3.11 works here as well (just notice that the $CRYS,top$ analogue of \emph{loc. cit.} Prop. 1.1.19 and Prop. 1.3.6 is true and then use directly \emph{loc. cit.} Prop. 2.3.1 and Cor. 2.3.9).
Similarly, \ref{lemext2iv} is a straightforward adaptation of \cite{BBM} Prop. 2.4.5. The second case of \ref{lemext2i} follows from the first one using \ref{lemext2iv} and Lem. \ref{lemlimind}.

\ref{lemext2v} This follows from Lem. \ref{lemext1} \ref{lemext1iii} using that $Rl_* \cal Hom_{X_{TOP}}(G_{p^.},M)=0$ (thanks to the fact that the transition morphisms $\cal Hom_{X_{TOP}}(G_{p^{i+j}},M)\to \cal Hom_{X_{TOP}}(G_{p^j},M)$ are zero if $M$ is killed by $p^i$).

\ref{lemext2iii} The first assertion is \cite{Ba} Prop. 1.18 if $G$ is a syntomic group scheme. The case of $p$-divisible groups follows by Lem. \ref{flfsyn}, \ref{lemext2v} and Lem. \ref{lemlimind}.  The second assertion follows from the first one since $i_*$ commutes to $\pi_*$ and $\cal Ext^i$.

\ref{lemext2ii} The claimed isomorphism is obtained from  $$Ru_*R\cal Hom_{(X/\Sigma)_{CRYS,syn}}(u^{-1}G,M)\simeq R\cal Hom_{X_{SYN}}(G,Ru_*M)$$ by forming $H^1$. The computation of the left hand term uses the vanishing of \break $\cal Hom_{(X/\Sigma)_{CRYS,syn}}(u^{-1}G,M)$ (which follows from the fact that $p$ is epimorphic on $u^{-1}G$ while $M$ is killed by $p^k$) and the computation of the right hand term uses the vanishing of $\cal Hom_{X_{SYN}}(G,R^1u_*M)=0$ (which follows from the fact that $p$ is epimorphic on $G$ while $R^1u_*M$ is killed by $p^k$). The same arguments apply in the case of small topoi.
\begin{flushright}$\square$\end{flushright}

\begin{rem} Thanks to Lem. \ref{lemext1} \ref{lemext1i}, the isomorphisms of Lem. \ref{lemext2} have obvious variants in the case $\Sigma=\Sigma_.$ (replace $G$ with $l^{-1}G$  in the formulae). \end{rem}

The following lemma relies on the vanishing theorem of \cite{Br} and will be a cornerstone in our proof.

\begin{lem} \label{lemext3} Consider $X$ in $\cal Sch/\Sigma_1$ and a $p$-divisible group $G/X$.
\debrom
\item \label{lemext3i} Let $1\le k<\infty$ If we view $G$ as an Abelian group in $X_{SYN}$ then: $$\begin{array}{rrcl}&\cal Ext^2_{X_{SYN}}(G,\mu_{p^k})&=&0\\
    \hbox{and} & \cal Ext^2_{X_{SYN}}(G,\O)&=&0\end{array}$$
    The same is true with $syn$ instead of $SYN$.
\item \label{lemext3ii} Let $k\le 1\le \infty$, then: $$\cal Ext^2_{(X/\Sigma_k)_{CRYS,syn}}(G,\cal I)\simeq \cal Ext^2_{(X/\Sigma_k)_{CRYS,syn}}(G,\cal O)=0$$\finrom
\end{lem}
Proof. \ref{lemext3i} Let us examine the first case. The starting point is the vanishing of $\cal
Ext^i(G_{p^l},\Gm)$ for $i=1,2$ which is an easy consequence of \cite{SGA7-I} VIII, Prop. 3.3.1 for $i=1$ (use that $\smash{{G_{p^l}^*}}$ is acyclic for $\epsilon_*$, $\epsilon:X_{FL}\to X_{SYN}$) and of the main
result of \cite{Br} for $i=2$ (here we use that $p\ne 2$). In the case of $X_{syn}$, we use furthermore Lem. \ref{lemext2} \ref{lemext2iii}. Using the Kummer exact sequence $0\to \smash{\mu_{p^k}}\to \Gm\to \Gm\to 0$, we obtain the following facts:

- $\cal Ext^2(G_{p^l},\mu_{p^k})=0$;

- for $l'\ge l$ the morphism $$\cal Ext^1(G_{p^{l'}},\mu_{p^k})\to \cal Ext^1(G_{p^{l}},\mu_{p^k})$$
identifies with the morphism $\smash{G_{p^{l'}}^*/p^k}\to \smash{G_{p^{l}}^*/p^k}$ induced by $\smash{p^{l'-l}:G_{p^{l'}}^*}\twoheadrightarrow
\smash{G_{p^l}^*}$. It is thus an isomorphism for  $l\ge k$.

The vanishing of $\cal Ext^2(G,\mu_{p^k})$ will follow from the spectral sequence $$R^il_*\cal
Ext^j(G_{p^.},\mu_{p^k})\Rightarrow \cal Ext^{i+j}(G,\mu_{p^k})$$ (Lem. \ref{lemext1} \ref{lemext1iii}) once checked that $R^{2-j}l_*\cal
Ext^j(G_{p^.},\mu_{p^k})$ vanishes for $j=0,1,2$. Recall that $R^ql_*M_.$ can be computed by sheafifying $U\mapsto Rlimproj_kR\G(U,M_k)$. For $j=1,2$, the desired vanishing follows directly from the above facts. For $j=0$, it results from the fact that  $$\cal Hom(G_{p^{l'}},\mu_{p^k})\to \cal Hom(G_{p^{l}},\mu_{p^k})$$ is zero for $l'\ge l+k$.

The vanishing of $\cal Ext^2(G,\O)$ is proved in the same way using the following facts, which result from \cite{BBM} Prop. 3.3.2 and its proof (thanks to the isomorphisms $\epsilon_*\cal Ext_{X_{FL}}^i(G_{p^l},\O)\simeq \cal Ext_{X_{SYN}}^i(\epsilon_*G_{p^l},\epsilon_*\O)$  (Lem. \ref{acycqcohsch} \ref{acycqcohschi}) and $\pi_*\cal Ext_{X_{SYN}}^i(G_{p^l},\O)\simeq \cal Ext_{X_{syn}}^i(\pi_*G_{p^l},\pi_*\O)$ (Lem. \ref{lemext2} \ref{lemext2iii}):

- $\cal Ext^2(G_{p^l},\O)=0$;

- for $l'\ge l$ the following morphism is invertible: $$\cal Ext^1(G_{p^{l'}},\O)\to \cal Ext^1(G_{p^{l}},\O)$$

- for $l'\ge l+1$ the following morphism is zero:   $$\cal Hom(G_{p^{l'}},\O)\to \cal Hom(G_{p^{l}},\O)$$

\ref{lemext3ii} This is Thm. 3.3.3 and Prop. 3.3.4 of \cite{BBM}  modulo Lem. \ref{lemext2} \ref{lemext2i}.
\begin{flushright}$\square$\end{flushright}

\begin{rem} \label{remext3} \debrom \item \label{remext3i} In Lem. \ref{lemext3} \ref{lemext3i}, $X_{SYN}$ or $X_{syn}$ (resp. $G$)  can be replaced by $X_{SYN}^\N$ or $X_{syn}^\N$ (resp. the constant projective system $l^{-1}G$) thanks to Lem. \ref{lemext1} \ref{lemext1i}.

\item  \label{remext3ii} We do not know whether or not the analogue of Lem. \ref{lemext3} \ref{lemext3ii} holds on the small  crystalline topos. Unpleasantly, this will force us to switch frequently between small and big topoi.
\finrom
\end{rem}

\para Let us come to the techniques of \cite{FM} which are needed for our purpose.
We begin with the description of $\O_k^{crys}$ using divided power envelopes of Witt vectors and syntomic sheafification.

\begin{lem} \label{lemcomp3} Consider the presheaf $\smash{W_k^{dp}}$ on $SYN(\Sigma_1)$ sending $U$ to the divided power envelope of $W_k(U)$ with respect to the kernel $I_k(U)$ of $W_k(U)\to \O(U)$, $(x_0,\dots, x_k)\mapsto \smash{x_0^{p^k}}$. Let $\smash{I_k^{dp}(U)}$ denote the tautological dp-ideal of $\smash{W_k^{dp}(U)}$.
\debrom
\item \label{lemcomp3i} The sheaf associated to $W_k^{dp}$ is canonically isomorphic to $\O^{crys}_k$ in $\Sigma_{1,SYN}$. The ring $\O^{crys}_.$ is in particular normalized.
\item \label{lemcomp3ii} If $U$ is an affine semi-perfect $\Sigma_1$-scheme (ie. $U=Spec(A)$ and the Frobenius of $A$ is surjective), the isomorphism of \ref{lemcomp3i} induces $$W_k^{dp}(U)\simeq \O^{crys}_k(U)$$
\item \label{lemcomp3iii} If $(U_i)_{i\in I}$ is a filtrant projective system of affine $\Sigma_1$-schemes and $U_\infty=limproj U_i$ then $$\limi \O_k^{crys}(U_i)\simeq \O_k^{crys}(U_\infty)$$
\item \label{lemcomp3iv} If $U$ is an affine $\Sigma_1$-scheme, there exists a filtrant projective system $(U_i)_{i\in I}$ of affine syntomic surjective $U$-schemes whose transition morphisms are syntomic surjective as well and such that $U_\infty$ is semi-perfect. If there exists a closed immersion $U\hookrightarrow Y$ where $Y$ has finite $p$-bases over $\Sigma_1$, then we may chose $I=\N$ and $\smash{U_i=U^{(p^i/Y)}}$. If the ideal of the closed immersion is moreover generated by a regular sequence, then $W_{k}^{dp}(U_\infty)$ is flat over $\Z/p^{k}$.
\finrom
\end{lem}

\noindent Proof. Recall that for any $U/\Sigma_1$, we have $\O_k^{crys}(U)=limproj \O(T')$ where the projective limit is indexed by $(U',T',\gamma')$ in $CRYS_{syn}(U/\Sigma_k)$ and does not change if only affine $T'$'s are considered. Thanks to the existence of $\gamma'$, the ring homomorphism \begin{eqnarray}\label{ThetakU}W_k(U)&\to& \O(T')\\
(x_0,x_1,\dots, x_{k-1})&\mapsto &(\tilde x_0)^{p^k}+p(\tilde x_1)^{p^{k-1}}+\dots+ p^{k-1}(\tilde x_k)^p\end{eqnarray}
is well defined if  $\tilde x_i\in \O(T')$ denotes an arbitrary lift of the image of $x_i$ in $\O(U')$. Passing to divided powers and letting $(U',T',\gamma')$ and $U$ vary, this homomorphism defines a morphism of presheaves of rings on $SYN(\Sigma_1)$: \begin{eqnarray}\label{defThetak} \Theta_k: W_k^{dp}\to \O_k^{crys}\end{eqnarray} We will prove that $\Theta_k$ induces the isomorphism announced in \ref{lemcomp3i} after sheafification. The other statements will be proved along the way.

\ref{lemcomp3ii} Since $F_A:x\mapsto x^p$ is surjective, the ring homomorphism $(x_0,\dots, x_k)\mapsto \smash{x_0^{p^k}}$ is surjective and defines a $pd$-thickening $(U,T):=(Spec(A),Spec(W_k^{dp}(A)))$. Let us check that $(U,T)$ is final in $CRYS_{syn}(U/\Sigma_k)$. Consider a pd-ideal $J$ of some $\Z/p^k$-algebra $B$. Any ring homorphism $f:A\to B/J$ lifts uniquely to $\smash{\tilde f}:W_k(A)\to B$ by the formula (\ref{ThetakU}). To see this, we first note that thanks to divided powers on $J$, $\smash{\tilde f([x_0^{p^k}])}=\smash{(\tilde f([x_0]))^{p^k}}$ is uniquely determined by the image $\smash{f(x_0^{p^k})}$ of $\smash{\tilde f([x_0])}$ in $B/J$ and we conclude using the formula $(x_0,\smash{x_1^p},\dots, \smash{x_{k-1}^{p^{k-1}}})=[x_0]+p[x_1]+\dots +p^{k-1}[\smash{x_{k-1}}]$ in $W_k(A)$.

\ref{lemcomp3iii}  Using a construction similar to the proof of Lem. \ref{lememb2} \ref{lememb2iv} (here the finiteness of $\delta/\Delta$ is not required), we find easily a projective system $(U_i,Y_i)_{i\in I}$ where $Y_i=Spec(B_i)$, with $B_i$ is a polynomial algebra (in possibly infinitely many variables) and $B_\infty:=limind_i B_i$ as well (it is in fact an infinite tensor product indexed by $I$ of polynomial algebras).
For $i\in I\sqcup \{\infty\}$ let us denote respectively $T_i$ and $\smash{T^{(1)}_i}$ the divided power envelope of $U_i$ inside $Y_i=Spec(B_i)$  and $Y_i\times Y_i$. Then $\smash{\O_k^{crys}(U_i)}\simeq Ker(\cal O(T_i)\rightrightarrows \smash{\cal O(T^{(1)}_i)})$ and the claim follows.

\ref{lemcomp3iv} Let $A=\O(U)$ and consider a set $I:=\sqcup_{n\ge 1} I_n$ where $I_1$ is the set of finite subsets of $A$, $A_{1,i}:=\otimes_{A,a\in i}A[X]/(X^p-a)$, $A_1=limind_{i\in I_1} A_{1,i}$, and where for $n\ge 2$, $I_n$, $A_{n,i}$ and $A_n$ are inductively defined respectively as follows: $I_n$ is the set of finite subsets in $A_{n-1}$, $A_{n,i}:=\otimes_{A_{n-1},a\in i}A_{n-1}[X]/(X^p-a)$, $A_n=limind_{i\in I_n} A_{n,i}$. The set $I$ is naturally ordered and the inductive system $i\in I\mapsto U_i:=Spec(A_i)$ satisfies the required conditions. In the second case, the choice of a $p$-basis $\underline s=(s_1,\dots,s_d)$ for $Y$ induces the right commutative square below $$\xymatrix{U^{(p^{i+1}/Y)}\ar[r]\ar[d]&Y\ar[r]^{\underline s}\ar[d]^{F}&\mathbb A^d_{\Sigma_1}\ar[d]^{F}\\ U^{(p^{i}/Y)}\ar[r]&Y\ar[r]^{\underline s}&\mathbb A^d_{\Sigma_1}}$$
The right vertical arrow is clearly syntomic surjective. The other vertical ones are thus syntomic surjective as well since both squares are cartesian. Let us explain  the last statement. Set $B=\O(Y)$ and chose a regular sequence $\underline f=(f_1,\dots,f_r)$ generating the ideal of $U$. Then $U_\infty=Spec(A_\infty)$ where $A_\infty=B_\infty/(f_1,\dots,f_r)$ and $B_\infty=limind_F B$. Also, $W_k^{dp}(U_\infty)$ is the divided power envelope (compatible with $(p)$, as always) of the ideal of $W_k(B_\infty)$ generated by the $\smash{[f_i^{p^{-k}}]}$'s. According to \cite{Be3} Prop. 1.5.3 (i), it is sufficient to prove that $W_k(B_\infty)/\smash{([f_1^{p^{-k}}],\dots,[f_r^{p^{-k}}])}$ is flat over $\Z/p^k$. We do this by induction on $r$. For $r=0$, this follows from the fact that $B_\infty$ is perfect. Let us denote $J_i=\smash{([f_1^{p^{-k}}],\dots,[f_i^{p^{-k}}])}$, if $0\le i\le r$. For $r\ge 1$, it is sufficient to prove that $\smash{Tor_1^{\Z/p^k}}(\Z/p,W_k(B_\infty)/J_r)$ vanishes, ie. that the following sequence is exact: $$\xymatrix{0\ar[r]&\Z/p\otimes (J_r/J_{r-1})\ar[r]&\Z/p\otimes (W_k(B_\infty)/J_{r-1})\ar[r]&
\Z/p\otimes (W_k(B_\infty)/J_{r})\ar[r]&0}$$
Now the image of $\smash{[f_r^{p^{-k}}]}$ in $\Z/p\otimes (W_k(B_\infty)/J_{r-1})\simeq B_\infty/\smash{(f_1^{p^{-k}},\dots,f_{r-1}^{p^{-k}})}$ is a non zero divisor and thus $J_r\cap (pW_k(B_\infty)+J_{r-1})=pJ_r+J_{r-1}$ as desired.

\ref{lemcomp3i} Let $F$ denote either the kernel or cokernel presheaf of the morphism $\Theta_k$ and let $U\in SYN(\Sigma_1)$. Thanks to \ref{lemcomp3ii}, \ref{lemcomp3iii} and \ref{lemcomp3iv}, we find that for each $x\in F(U)$, there exists a syntomic surjective family $(U_i\to U)$ such that the restriction of $x$ to each $U_i$ vanishes. It follows in particular that the sheaf associated to $F$ is zero.
\begin{flushright}$\square$\end{flushright}

Let us continue with some technical complements in the setting of small syntomic topoi. We denote $\smash{\cal I_.^{crys}}$ the canonical $dp$-ideal of $\smash{\O_.^{crys}}$ and  $$\smash{\widetilde{\cal I}{}^{crys}_.:=\jj^* \cal I_.^{crys}}$$ where $\jj$ has the same meaning as in Def. \ref{deffilsynmod} (ie. $\smash{\widetilde{\cal I}{}^{crys}_k}=\smash{{\cal I}{}^{crys}_{k+1}}/p^k$, see Def. \ref{defj}).


\begin{lem} \label{lemcomp4} Consider a regular closed immersion of affine schemes $U\to Y$, where $Y$ has finite $p$-bases over $\Sigma_1$.
\debrom
\item \label{lemcomp4i} The modules $\O_{.}^{crys}$ and $\smash{\widetilde{\cal I}{}^{crys}_.}$ of $(U_{syn}^\N,\Z/p^.)$ are flat and normalized.
\item \label{lemcomp4ii} If $Y$ is semi-perfect then $\Theta_k$ induces $$I_{k+1}^{dp}(U)/p^k \simeq {\widetilde{\cal I}}{}^{crys}_k(U)$$
\item \label{lemcomp4iii} If $U_\infty:=limproj U^{(p^i/Y)}$ then $$\limi {\widetilde{\cal I}}{}^{crys}_k(U_i)\simeq {\widetilde{\cal I}}{}^{crys}_k(U_\infty)$$
\finrom
\end{lem}
\noindent Proof. \ref{lemcomp4i} The case of $\O_{.}^{crys}$ follows easily from the flatness statement in Lem. \ref{lemcomp3} \ref{lemcomp3iv} as in the proof of Lem. \ref{lemcomp3} \ref{lemcomp3i}.
Alternatively, one could refer to Lem. \ref{lemsyn} \ref{lemsyniv}. The case of  $\smash{\widetilde{\cal I}{}^{crys}_.}$ follows by Lem. \ref{cartgen} \ref{cartgeniv}.

\ref{lemcomp4ii} Starting with the exact sequence defining $I^{dp}_{k+1}(U)$ (resp.  $\cal I^{crys}_{k+1}$) and applying Lem. \ref{cartgen} \ref{cartgeniv} (resp. applying Lem. \ref{cartgen} \ref{cartgeniv} and evaluating at $U$) gives the top  (resp. bottom) exact sequence in the following commutative diagram:
\begin{eqnarray}\label{diaglemcomp4}\xymatrix{0\ar[r]&\O(U)\ar[r]\ar[d]^\parallel&I^{dp}_{k+1}(U)/p^k\ar[r]\ar[d]&I^{dp}_k(U)\ar[r]\ar[d]&0\\
0\ar[r]&\O(U)\ar[r]&{\widetilde{\cal I}}{}^{crys}_k(U)\ar[r]&\cal I^{crys}_k(U)}\end{eqnarray}
The result follows, since the right vertical arrow is invertible thanks to Lem. \ref{lemcomp3} \ref{lemcomp3ii}.

\ref{lemcomp4iii}  Because $\O$ is ayclic on the $U_i$'s,  we have an exact sequence $0\to \O(U_i)\to \smash{{\widetilde{\cal I}}{}^{crys}_k(U_i)}\to \smash{{{\cal I}}{}^{crys}_k(U_i)}\to 0$. Thanks to Lem. \ref{lemcomp3} \ref{lemcomp3iii}, we may conclude by taking $limind_i$ and comparing with the bottom line of (\ref{diaglemcomp4}).
\begin{flushright}$\square$\end{flushright}

\begin{rem} \label{remThetak} \debrom \item \label{remThetaki} The Cartier morphism $c_1:=C^{-1}:\O\to \O_1^{crys}$ of Lem.-Def. \ref{defCinv} \ref{defCinvii}  generalizes as follows. For any $i\ge 1$ there is a unique $c_i$ fitting into a commutative triangle of rings $$\xymatrix{\O_i^{crys}\ar[rd]_{nat}\ar[rr]^{F^i}&& \O_i^{crys}\\ &\cal O\ar[ru]_{c_i}&}.$$

\item \label{remThetakii} It follows e.g. from Lem. \ref{lemcomp3} that $\smash{\O_.^{crys}}$ is normalized in $\smash{(\Sigma_{1,SYN}^\N,\Z/p^.)}$. As a result, we have a well defined morphism: $p^{k-i}:\O_{i}^{crys}\simeq \O^{crys}_k/p^{i}\to \O_k^{crys}$. Using this, we find that $\Theta_k$ can be described globally as the morphism deduced from $$\begin{array}{rcl}W_k(U)&\to &\O^{crys}_k(U)\\ (x_0,x_1,\dots, x_{k-1})&\mapsto& c_k(x_0)+pc_{k-1}(x_1)+\dots+ p^{k-1}c_1 (x_{k-1})\end{array}$$ by the universal property of divided power envelopes.

    \item \label{remThetakiii} If $U$ is any scheme of characteristic $p$, we have a commutative diagram of rings: \begin{eqnarray}\label{diagTheta}\xymatrix{W_k^{dp}(U)\ar[d]_{F^{k-1}}\ar[r]^{\Theta_k}&\cal O^{crys}_k(U)\ar[d]^{nat}\\
 W_1^{dp}(U)\ar[d]_{F}\ar[r]^{\Theta_1}&\cal O^{crys}_1(U)\ar[d]^{nat}\\
 \O(U)\ar@{=}[r]&\O(U)}\end{eqnarray} This diagram is functorial with respect to $k$ and $U$ in the obvious way.
\finrom
\end{rem}

\begin{lem}  \label{lemcomp2}
Assume that $X$ admits locally a closed regular immersion into a $\Sigma_1$-scheme $Y$ with finite $p$-bases.

\debrom
\item \label{lemcomp2i} The following sequence of $(X_{syn}^\N,\O^{crys}_.)$ is exact:
$$\xymatrix{0\ar[r]&\I^{crys}_.\ar[r]&\O^{crys}_.\ar[r]&\O\ar[r]&0}$$
The modules $\smash{\cal O}$ and $\smash{\cal I}$ of $(X/\Sigma_.,\O)_{crys,syn}$ are acyclic for $u_*$.

\item \label{lemcomp2ii} There exists a unique $\varphi$ making the following square commutative over $(X_{syn}^\N,\O^{crys}_.)$
$$\xymatrix{\widetilde {\cal I}^{crys}_.\ar[r]^-\varphi\ar[d]_-1&\O^{crys}_.\ar[d]^p\\ \O^{crys}_.\ar[r]^-F&\O_.^{cris}}$$

\item \label{lemcomp2iii} There are canonical exact sequences as follows over $(X^\N_{syn},\Z/p^.)$:
$$\begin{array}{c}
\xymatrix{0\ar[r]&\Z/p^.\ar[r]&\O^{crys}_.\ar[r]^{1-F}&\O^{crys}_.\ar[r]&0}\\
\xymatrix{0\ar[r]&\mu_{p^.}\ar[r]&\widetilde{\cal I}^{crys}_.\ar[r]^{1-\varphi}&\O^{crys}_.\ar[r]&0}\end{array}$$

\item \label{lemcomp2iv} The morphism $\varphi$ and the exact sequences of \ref{lemcomp2i}, \ref{lemcomp2iii} are functorial with respect to $X$. All statements thus generalize \emph{verbatim} if $X$ is replaced by a diagram of $Sch/\Sigma_1$ whose vertices have local finite $p$-bases.
\finrom
\end{lem}

Proof. \ref{lemcomp2i} The exact sequence follows either from Lem. \ref{lemsyn} \ref{lemsyniii} or Lem. \ref{lemcomp3} \ref{lemcomp3i}. The acyclicity of $\O$ (resp. $\Ga$) is by Lem. \ref{lemsyn} \ref{lemsyniii} (resp. is immediate) and that of $\cal I$ follows from the exact sequence.

\ref{lemcomp2ii} Thanks to Lem. \ref{lemcomp4} \ref{lemcomp4i}, existence and unicity follow from Lem. \ref{cartgen} \ref{cartgenii}, \ref{cartgeniii} as in the proof of Prop. \ref{defphi}, once noticed that the Frobenius endomorphism of $\smash{\O_1^{crys}}$ vanishes on $\smash{\cal I_1^{crys}}$ ($x^{[p]}=p!x^p=0$).

\ref{lemcomp2iii} The first exact sequence is left to the reader. Let us explain the second one. Starting with the exact sequence of Abelian groups $$\begin{array}{rl}\xymatrix{0\ar[r]&1+\cal I\ar[r]&\O^\times \ar[r]&\Gm\ar[r]&0}&\hbox{in $(\Sigma_1/\Sigma_k)_{CRYS,syn}$}\end{array}$$
then combining with $log:1+\I\to \I$, applying $\Z/p^k\otimes^L(-)$ and $u_*$ we find a morphism $\smash{\mu_{p^k}}\to \cal I_k^{crys}$ in $\Sigma_{1,SYN}$. Since the Frobenius acts like $p$ on $\mu_{p^k}$ this morphism factors via the kernel of $p-F:\smash{\cal I_k^{crys}}\to \smash{\O_k^{crys}}$. Next, we let $k$ vary, restrict to the small syntomic site of $X$ and apply the functor $\jj^*$ in order to get a morphism $\jj^*\mu_{p^.}\to \jj \smash{\widetilde{\cal I}{}_k^{crys}}$ whose composition with $p-F:\smash{\widetilde{\cal I}{}_k^{crys}}\to \smash{\O_k^{crys}}$ is zero. We note that  $\jj^*\mu_{p^.}\simeq \mu_{p^.}$. Since $p\circ (1-\varphi)=p-F$ and $\smash{\O_.^{crys}}$ is flat and normalized over $\Z/p^.$, this yields by Lem. \ref{cartgen} \ref{cartgeniii} a complex $$\xymatrix{0\ar[r]&\mu_{p^.}\ar[r]^-{\cal L}&\widetilde{\cal I}_.^{crys}\ar[r]^-{1-\varphi}&\cal O_.^{crys}\ar[r]&0}$$
over $\smash{(X^\N_{syn},\Z/p^.)}$. It remains to show that this complex is exact. Since $\smash{\mu_{p^.}}$, $\smash{\widetilde{\cal  I}_.^{crys}}$ and $\smash{\O_.^{crys}}$ are $L$-normalized, it is sufficient to check exactness on the component $k=1$ (use Lem. \ref{cartgen} \ref{cartgenii}). Using Lem. \ref{lemcomp3} \ref{lemcomp3ii},  \ref{lemcomp3iii}, \ref{lemcomp3iv} and Lem. \ref{lemcomp4} \ref{lemcomp4ii}, \ref{lemcomp4iii} we are reduced to check the exactness of the sequence $$\xymatrix{0\ar[r]&\mu_{p}(A_\infty)\ar[r]^-{\cal L}&I_2^{dp}(A_\infty)/p\ar[r]^-{1-\varphi} & W_1^{dp}(A_\infty)\ar[r]&0}$$ with $A_\infty=B_\infty/(f_1,\dots,f_r)$, as in the proof of the last statement in Lem. \ref{lemcomp3} \ref{lemcomp3iv}. In order to do this, our first task is to give an explicit description of $W_2^{dp}(A_\infty)/pI_2^{dp}(A_\infty)$, $W_1^{dp}(A_\infty)$ and of the maps $\cal L$ and $1-\varphi$. We will use the following facts.
\petit

\noindent \emph{(Fact 1)}: The natural map $\smash{D(W_k(B_\infty),([f_1^{p^{-k}}],\dots,[f_r^{p^{-k}}]))\to W^{dp}(A_\infty)}$ is an isomorphism.
\petit

\noindent \emph{(Fact 2)}: The $[f_i^{p^{-k}}]$'s form a regular sequence in $W_2(B_\infty)$.

\noindent \emph{(Fact 3)}: Consider the ideal $J$ generated by a regular sequence $(g_1,\dots,g_r)$ in a  $\Z/p^k$-algebra $C$. For $d\ge r$, let $<\underline Y>$ denote the tautological pd-ideal of $C<\underline Y>=C<Y_1,\dots, Y_d>$ and let $I$ denote the ideal generated by the $Y_i-f_i$'s. Then $I\cap <\underline Y>$ coincides with $I<\underline Y>$. It is in particular a sub pd-ideal of $<\underline Y>$ and $D(C,J)=C<\underline Y>/I$.
\petit

Let us explain this briefly. Fact 1: this simply follows from the fact that the image of the $\smash{p^{j}[f_i^{p^{-j}}]}$'s vanish in $\smash{D(W_k(B_\infty),([f_1^{p^{-k}}],\dots,[f_r^{p^{-k}}]))}$, turning the latter into an algebra over $W_k(A_\infty)$. Fact 2: During the proof of Lem. \ref{lemcomp3} \ref{lemcomp3iv} we have established that $\smash{W_k(B_\infty)/([f_1^{p^{-k}}],\dots, [f_i^{p^{-k}}])}$ is $\Z/p^k$-flat for $0\le i\le r$. If $0\le i\le r-1$, this implies that $\Z/p\otimes K=0$, where $K$ denotes the kernel of the multiplication by $[f_{i+1}]$ on $\smash{W_k(A_\infty)/([f_1^{p^{-k}}],\dots, [f_i^{p^{-k}}])}$. It follows that $K=pK=p^kK=0$ as desired. Fact 3: the equality $I\cap <\underline Y>=I<\underline Y>$ can easily be proven by induction on $r$ and the remaining statements follow immediately.
\petit

Using these facts we find the following decomposition into cyclic $A_\infty$-modules  (resp. $W_2(A_\infty)$-modules): $$\begin{array}{rcll}
 W_1^{dp}(A_\infty)&\simeq &D(B_\infty,(f_1^{p^{-1}},\dots, f_r^{p^{-1}}))&
\\ &\simeq& B_\infty<Y_1,\dots,Y_r>/(Y_1-f_1^{p^{-1}},\dots, Y_r-f_r^{p^{-1}})&
\\ &\simeq &\oplus_{\underline n} A_\infty.(\underline f^{p^{-1}})^{[p\underline n]}&
\\  W_2^{dp}(A_\infty)/pI_2^{dp}(A_\infty)&\simeq & D(W_2(B_\infty),([f_1]^{p^{-2}},\dots, [f_r]^{p^{-2}}))/pI_2^{dp}(B_\infty)&
\\ &\simeq & W_2(B_\infty)<Y_1,\dots,Y_r>/p<\underline Y>+(Y_1-[f_1^{p^{-2}}],\dots, Y_r-[f_r^{p^{-2}}])&
\\ &\simeq & W_2(A_\infty/(f_1^{p^{-1}},\dots, f_r^{p^{-1}})).1&
\\ && \oplus (\oplus_{|\underline n|\ge 1} (W_2(A_\infty/(f_1^{p^{-1}},\dots, f_r^{p^{-1}}))/p).[\underline f^{p^{-2}}]^{[p\underline n]})&
\end{array}$$

Next we observe that the following square is commutative: $$\xymatrix{(0,a).[f^{p^{-2}}]^{[p\underline n]}&W_2^{dp}(A_\infty)\ar[r]^-{\Theta_2}_-\sim&\cal O_2^{crys}(A_\infty)\\
a.(\underline f^{p^{-1}})^{[p\underline n]}\ar[u]\ar@{|->}[u]&W_1^{dp}(A_\infty)\ar[u]^V\ar[r]^-{\Theta_1}_-\sim&\cal O_1^{crys}(A_\infty)\ar[u]^-p}$$ where the vertical arrow denoted $V$ is defined by the given formula.  With this in mind, we find that the morphism $\varphi:\smash{I_2^{dp}(A_\infty)}\to W_1^{dp}(A_\infty)$ is described by the formulae (which are valid only since $p\ge 3$) $$\begin{array}{rcl} \varphi(0,a)&=&a^p\\
\varphi((a_0,a_1).[f_i^{p^{-2}}])&=& (p-1)!a_0^{p^2}(f_i^{p^{-1}})^{[p]}\\
\varphi((a_0,a_1).[\underline f^{p^{-2}}]^{[
p\underline n]})&=&0 \hbox{ if $\underline n\ne 0$}\end{array}$$
An arbitrary element of $\smash{I_2^{dp}(A_\infty)}/p$ may be written $x=(0,a)+\sum_i[a_if_i^{p^{-2}}]+\sum_{\underline n\ne 0}[a_{\underline n}\underline f^{p^{-2}}]^{[p\underline n]}$. Then we have: $$\begin{array}{c}(1-\varphi)(x)=(\sum_ia_i^pf_i^{p^{-1}} + \sum_{\underline n\ne 0}(a_{\underline n}^p\underline f^{p^{-1}})^{[p\underline n]})-(a^p-\sum(a_i^pf_i^{p^{-1}})^{[p]})\end{array}$$
This formula shows that $1-\varphi: \smash{I_2^{dp}(A_\infty)}/p\to \smash{W_1^{dp}(A_\infty)}$ is surjective and that $x$ is in its kernel if and only if the following are satisfied in $\smash{A_\infty/(\underline f^{p^{-1}})}$: $$\left\{\begin{array}{l}a=\sum a_if_i^{p^{-2}}\\ a_{e_i}^p=-a_i^p \hbox{ ($i=1,\dots, r$)} \\ a^{p\underline n}_{\underline n}=0 \hbox{ ($|\underline n|\ge 2$)}\end{array}\right.$$
ie. if $x$ is of the form $$\begin{array}{c}x=(0,\sum a_if_i^{p^{-2}})+\sum [a_if_i^{p^{-2}}] - \sum [ a_if_i^{p^{-2}}]^{[p]}\end{array}$$
In other terms, $Ker(1-\varphi)$ is generated by elements of the form $(z,z)-[z]^{[p]}$ with $z$ running in one of the ideals $\smash{(f_i^{p^{-2}})}$ of $A_\infty/(\smash{\underline f^{p^{-1}}})$.

Let us now investigate the map $\cal L:\mu_p(A_\infty)\to \smash{I_2^{dp}(A_\infty)}$. Let $\zeta=\smash{1+z^{p^3}}$ denote an element of $\mu_p(A_\infty)$ with $z\in (\smash{\underline f^{p^{-4}}})$. By definition of $\cal L$, we have $\cal L(\zeta)= \smash{log((1+[z])^{p^2})}$. We have the following equalities in $\smash{I_2^{dp}(A_\infty)}/p$: $$\begin{array}{rcl}\cal L(1+z^{p^3})&=&log(1+p(\sum_{k=1}^{p-1}({1\over p}\left(p^2\atop pk\right))[z]^{pk})+[z^{p^2}])\\
&=&p\sum_{k=1}^{p-1}{(-1)^{k-1}\over k}[z^p]^k+\sum_{k=1}^{p-1}{(-1)^{k-1}\over k}[z^{p^2}]^k-[z^{p^2}]^{[p]}\\
&=&(\sum_{k=1}^{p-1}{(-1)^{k-1}\over k}z^{p^2k},\sum_{k=1}^{p-1}{(-1)^{k-1}\over k}z^{p^2k})-[\sum_{k=1}^{p-1}{(-1)^{k-1}\over k}z^{p^2k}]^{[p]}\end{array}$$
Now $x^p\mapsto \smash{\sum_{k=1}^{p-1}{(-1)^{k-1}\over k}x^k}$ is a bijection of  the ideal $(\smash{\underline f^{p^{-1}}})$ of $A_\infty$ onto the ideal  $(\smash{\underline f^{p^{-2}}})$ of $A_\infty/(\smash{\underline f^{p^{-1}}})$. In particular, the map $\cal L$ is injective. Since its image is   \emph{a priori} contained in $Ker\, (1-\varphi)$ the previous description of the latter yields the desired equality.

\ref{lemcomp2iv} Consider $X'$ with local finite $p$-bases and a morphism $f:X'\to X$. We have to check that the following squares of Abelian groups of $X'_{syn}$ are commutative: $$\xymatrix{\cal O_.^{crys}\ar[r]^-F&\cal O_.^{crys}&\widetilde{\cal I}{}^{crys}_.\ar[r]^-{\varphi}&\cal O^{crys}_.&\mu_{p^.}\ar[r]^-{\cal L}&\widetilde{\cal I}{}^{crys}_.\\
f^{-1}\cal O_.^{crys}\ar[u]\ar[r]^-{f^{-1}(F)}&f^{-1}\cal O_.^{crys}\ar[u]&f^{-1}\widetilde{\cal I}{}^{crys}_.\ar[r]^-{f^{-1}(\varphi)}\ar[u]&f^{-1}\cal O^{crys}_.\ar[u]&f^{-1}\mu_{p^.}\ar[u]\ar[r]^-{\cal L}&\widetilde{\cal I}{}^{crys}_.\ar[u]}$$
The first square is clear. The second follows from the first one using Lem. \ref{cartgen} \ref{cartgeniii}. The third one is obvious since both horizontal arrows come from a morphism in $\smash{\Sigma^\N_{1,SYN}}$ by restriction.

\begin{flushright}$\square$\end{flushright}



\begin{lem} \label{lemcomp33} Assume $X/\Sigma_1$ has local finite $p$-bases.
\debrom  \item \label{lemcomp33i} The Cartier morphism $C^{-1}:\O\to \O_1^{crys}$ is monomorphic in $X_{syn}$. It induces an exact sequence as follows over
$((X/\Sigma_1)_{crys,syn},\O)$: $$\xymatrix{0\ar[r]&\cal I\ar[r]&\cal O\ar[r]&i_*\phi_*\O}$$
\item \label{lemcomp33ii} Let us use the notation $\epsilon$ for either one of the morphisms $(X_{syn},\O)\to (X_{et},\O)$ or $((X/\Sigma_1)_{crys,syn},\O)\to
((X/\Sigma_1)_{crys,et},\O)$. If $M$ is a locally free crystal with trivial $p$-curvature on
$((X/\Sigma_1)_{crys,et},\O)$ then the canonical base change morphism $$\epsilon^*\phi_*M\to \phi_*\epsilon^*M$$ is monomorphic in $Mod(X_{syn},\O)$.
\finrom
\end{lem}
Proof. \ref{lemcomp33i} Let $A_\infty$ be as in the proof of Lem. \ref{lemcomp3} \ref{lemcomp3iv}. We have the following commutative triangle (see (\ref{diagTheta})): $$\xymatrix{W_1^{dp}(A_\infty)\ar[rd]_F\ar[rr]^F&&W_1^{dp}(A_\infty)\\ &A_\infty\ar[ru]_{\Theta_1^{-1}\circ C^{-1}}&}$$
The explicit computations of the proof of Lem. \ref{lemcomp2} \ref{lemcomp2iii} then show that the arrow $\smash{\Theta_1^{-1}\circ C^{-1}}$ is injective and this suffices to conclude.

\ref{lemcomp33ii} Consider the following commutative square: $$\xymatrix{\epsilon^*i^{-1}M\ar[r]_\sim \ar[d]_{\epsilon^*(C^{-1})}&i^{-1}\epsilon^*M\ar[d]^{C^{-1}}\\ \epsilon^*\phi_*M\ar[r]^{ch}&\phi_*\epsilon^*M}$$ We want to check that the bottom horizontal arrow is monomorphic. The left vertical one is an isomorphism by the Cartier equivalence on small \'etale sites (see the proof of Prop. \ref{liesurj} \ref{liesurji} for details).
We claim that the right vertical arrow is monomorphic. Since  $M$ is locally free and has trivial $p$-curvature, \'etale localization reduces us to the case $M=\phi^*\O$, ie. $M=\O$. But then we are done by \ref{lemcomp33i}.

\begin{flushright}$\square$\end{flushright}

\para The following result explains a posteriori the ad hoc definition given in Def. \ref{deflie}.

\begin{lem} \label{synsyn1} Assume that $X/\Sigma_1$ has local finite $p$-bases and consider $G$ in $pdiv(X)$. Let $D$ denote the associated Dieudonn\'e crystal over
$((X/\Sigma_\infty)_{crys,et},\O)$. 
Using the definitions of Def. \ref{deffilsyn}, we have the following.

\debrom \item \label{synsyn1i} There is a canonical isomorphism of exact sequences  over
$((X/\Sigma_\infty)_{crys,syn},\O)$: $$\hspace{-1cm}\xymatrix{0\ar[r]&\pi_*\cal Ext^1_{(X/\Sigma_\infty)_{CRYS,syn}}(G^*,\cal I)\ar[d]_\wr\ar[r]&
\pi_*\cal Ext^1_{(X/\Sigma_\infty)_{CRYS,syn}}(G^*,\cal O)\ar[d]_\parallel\ar[r]&\cal Ext^1_{(X/\Sigma_\infty)_{crys,syn}}(G^*,\Ga)\ar[d]_\wr\ar[r]&0\\
0\ar[r]&Fil^1D^{syn}\ar[r]&Fil^0D^{syn}\ar[r]&i_*Lie^{syn}(D)\ar[r]&0}$$
\item \label{synsyn1ii} There is a canonical isomorphism of exact sequences  over $(\smash{X_{syn}^\N,\O_.^{crys}})$:
$$\hspace{-.1cm}\xymatrix{0\ar[r]&\cal Ext^1_{X_{syn}^\N}(l^{-1}G^*,\cal I^{crys}_.)\ar[r]\ar[d]^\wr&\cal Ext^1_{X_{syn}^\N}(l^{-1}G^*,\cal O^{crys}_.)\ar[r]\ar[d]^\wr&\cal Ext^1_{X_{syn}^\N}(l^{-1}G^*,\O)\ar[r]\ar[d]^\wr&0\\
0\ar[r]&Fil^{1,crys}_.D\ar[r]&Fil^{0,crys}_.D\ar[r]&Lie^{syn}_.(D)\ar[r]&0}$$
\item \label{synsyn1iii} The diagrams of \ref{synsyn1i} and \ref{synsyn1ii} are naturally functorial with respect to $X$.
\finrom
\end{lem}
\noindent Proof. \ref{synsyn1i}  Recall that by definition $D^{syn}=\epsilon^*D$ and $\overline D=\iota_1^{-1}D$. Let us furthermore denote $\overline D{}^{syn}:=\iota_1^{-1}D^{syn}\simeq \epsilon^*\overline D$, $D^{SYN}:=\pi^*D^{syn}$ and $\overline D{}^{SYN}:=\iota_1^{-1}D^{SYN}\simeq \pi^*\overline D{}^{syn}$.
By Lem. \ref{compFext} \ref{compFextii}, we have a commutative square as follows in $((X/\Sigma_1)_{CRYS,syn},\O)$:
(note that $\smash{\iota_1^{-1}\cal Ext^1_{(X/\Sigma_\infty)_{CRYS,syn}}(-,-)}\simeq \smash{\cal Ext^1_{(X/\Sigma_1)_{CRYS,syn}}(\iota_1^{-1}(-),\iota_1^{-1}(-))}$ since  $\smash{\iota^{-1}_1}$ has an exact left adjoint).

\begin{eqnarray} \label{squiobig} \xymatrix{\cal Ext^1_{(X/\Sigma_1)_{CRYS,syn}}(G^*,\O)\ar[d]^\wr \ar[r]^-{C^{-1}}&  \cal Ext^1_{(X/\Sigma_1)_{CRYS,syn}}(G^*,i_*\phi_*\O)\ar[d]^\wr \\
\overline D{}^{SYN}\ar[r]^-{i_*\phi_*(\overline f)\circ C^{-1}}&i_*\phi_*\overline D{}^{SYN}}
\end{eqnarray} Projecting to $((X/\Sigma_1)_{crys,syn},\O)$ we get (see Lem. \ref{lemext2} \ref{lemext2iii})
\begin{eqnarray} \label{squio} \xymatrix{\pi_*\cal Ext^1_{(X/\Sigma_1)_{CRYS,syn}}(G^*,\O)\ar[d]^\wr \ar[r]^-{C^{-1}}&  \cal Ext^1_{(X/\Sigma_1)_{crys,syn}}(G^*,i_*\phi_*\O)\ar[d]^\wr \\
\overline D{}^{syn}\ar[r]^-{i_*\phi_*(\overline f)\circ C^{-1}}&i_*\phi_*\overline D{}^{syn}}
\end{eqnarray}

\petit

\emph{Claim}. The image of the bottom (resp. top) horizontal arrow in (\ref{squio}) identifies canonically with $i_*Lie^{syn}(D)$ (resp. $\smash{\cal Ext^1_{(X/\Sigma_1)_{crys,syn}}(G^*,\Ga)}$).
\petit

\noindent Let us prove the first part of the claim. To begin with, we notice that the monomorphism \begin{eqnarray}\label{monozerzer}\xymatrix{Coker\, \overline v^{syn}\ar @{^(->} [r]^-{\overline f{}^{syn}}&\overline D^{syn}} \hbox{in $\cal Crys((X/\Sigma_1)_{crys,syn},\O)$}\end{eqnarray} induces a monomorphism \begin{eqnarray}\label{monozer}\xymatrix{\phi_* Coker\, \overline v{}^{syn}\ar @{^(->} [rr]^-{\phi_*(\overline f{}^{syn})}&&\phi_*\overline D^{syn}} \hbox{in $Mod(X_{syn},\O)$}\end{eqnarray}
of $\O$-modules on $X_{syn}$.
The claim concerning the image of the top arrow then follows from the commutative diagram $$\xymatrix{\overline D{}^{syn}\ar[d]^{\parallel}\ar[r]^-{C^{-1}}&i_*\phi_*Coker\, \overline v{}^{syn}\ar @{^(->}[rr]^-{i_*\phi_*(\overline f{}^{syn})} &&i_*\phi \overline D{}^{syn}\\
\epsilon^*\overline D\ar@{->>}[r]^-{C^{-1}}&i_*\epsilon^*\phi_*Coker\, \overline v\ar[r]^-=\ar[u]^{ch}&i_*Lie^{syn}(D)}$$ since the arrow denoted $ch$ is monomorphic by Lem. \ref{lemcomp33} \ref{lemcomp33ii}. Note that the bottom arrow denoted $C^{-1}$ uses the isomorphism $\epsilon^*i_*\simeq i_*\epsilon^*$ of Lem. \ref{lemsyn} \ref{lemsyni}.


Let us now prove the second part of the claim. Consider the exact sequences $$\begin{array}{rl} \xymatrix{0\ar[r]&\cal I\ar[r]& \O\ar[r]^{nat}& \Ga\ar[r]& 0}& \hbox{over $(X/\Sigma_1)_{CRYS,syn}$}\\
  \xymatrix{0\ar[r]&\Ga\ar[r]^{C^{-1}}& i_*\phi_*\O\ar[r]& Coker\, C^{-1}\ar[r]& 0}& \hbox{over $(X/\Sigma_1)_{crys,syn}$ (see Lem. \ref{lemcomp33} \ref{lemcomp33i})}\end{array}$$
It suffices to prove that the left (resp. right) vertical  arrow in the following commutative square is epimorphic (resp. monomorphic) $$\xymatrix{\pi_*\cal Ext^1_{(X/\Sigma_1)_{CRYS,syn}}(G^*,\O)\ar[d]^{nat}\ar[r]^-{C^{-1}}&\cal Ext^1_{(X/\Sigma_1)_{crys,syn}}(G^*,i_*\phi_*\O)\\
\pi_* \cal Ext^1_{(X/\Sigma_1)_{CRYS,syn}}(G^*,\Ga)\ar[r]^\sim_{\ref{lemext2} \ref{lemext2iii}}&\cal Ext^1_{(X/\Sigma_1)_{crys,syn}}(G^*,\Ga)\ar[u]_{C^{-1}}}$$
Now, this follows from the vanishing of $\smash{\cal Ext^2_{(X/\Sigma_1)_{CRYS,syn}}(G^*,\cal I)}$, Lem. \ref{lemext3} \ref{lemext3ii} (resp. of  \break $\smash{\cal Hom_{(X/\Sigma_1)_{crys,syn}}(G^*,Coker\, C^{-1})}$).

Using the claim and applying $\iota_{1,*}$ to (\ref{squio}), we find compatible epimorphisms as follows over $\smash{((X/\Sigma_\infty)_{crys,syn},\O)}$:

$$\xymatrix{\pi_*\cal Ext^1_{(X/\Sigma_\infty)_{CRYS,syn}}(G^*,\O)\ar[d]^\wr\ar@{->>}[r]&\iota_{1,*}\pi_*\cal Ext^1_{(X/\Sigma_1)_{CRYS,syn}}(G^*,\O)\ar[d]^\wr\ar@{->>}[r]^-{nat}&\cal Ext^1_{(X/\Sigma_\infty)_{crys,syn}}(G^*,\Ga)\ar[d]^\wr\\
D^{syn}\ar@{->>}[r]&\iota_{1,*}\overline D{}^{syn}\ar@{->>}[r]^-{can_D^{syn}}&i_*Lie^{syn}(D)}$$
The result now follows immediately from the vanishing of $\smash{\cal Hom_{(X/\Sigma_\infty)_{CRYS,syn}}(G^*,\Ga)}$ and the definition of $Fil^1D^{syn}$.

\ref{synsyn1ii} By definition, the second line is obtained from the second line of \ref{synsyn1i} by the functor $u_*l^{-1}$ (Def. \ref{deffilsynmod} \ref{deffilsynmodi}). It is exact by Prop. \ref{tdfilonesyn} \ref{tdfilonesyni}. It is thus sufficient to observe that the first line is also obtained from the first one by applying  $u_*l^{-1}$ (recall that $l^{-1}\pi_*\simeq \pi_*l^{-1}$, $u_*\pi_*\simeq \pi_*u_*$ and use Lem. \ref{lemext1} \ref{lemext1i} and Lem. \ref{lemext2} \ref{lemext2ii}).

\ref{synsyn1iii} We have to check that there are naturally defined base change morphisms for each vertex of the diagrams in question and that those are compatible to each other.
Consider a morphism $f:X'\to X$ where  $X'$ has local finite $p$-bases as well. Let $G'$ denote the base change of $G$ to $X'$. For the purpose of this proof, we view $G$ (resp. $G'$) as an Abelian group in $X_{SYN}$ (resp. $X'_{SYN}$). We thus have $f^{-1}G\simeq G'$. Using the natural isomorphisms $f^{-1}i_*\simeq i_*f^{-1}$,
$\pi f\simeq f\pi$
and the compatibility of $\cal Ext^1$ with localization on the big syntomic crystalline site, we find compatible base change morphisms (these are in fact isomorphisms but we don't need that here) over $((X'/\Sigma_\infty)_{crys,syn},\O)$:
{\small $$\hspace{-1cm}\xymatrix{\pi_*\cal Ext^1_{(X'/\Zp)_{CRYS,syn}}(i_*G'^*,\cal I)\ar[r]&
\pi_*\cal Ext^1_{(X'/\Zp)_{CRYS,syn}}(i_*G'^*,\cal O)\ar[r]&\pi_*\cal Ext^1_{(X'/\Zp)_{crys,syn}}(i_*G'^*,\Ga)\\
f^*\pi_*\cal Ext^1_{(X/\Zp)_{CRYS,syn}}(i_*G^*,\cal I)\ar[r]\ar[u]&
f^*\pi_*\cal Ext^1_{(X/\Zp)_{CRYS,syn}}(i_*G^*,\cal O)\ar[r]\ar[u]&f^*\pi_*\cal Ext^1_{(X/\Zp)_{crys,syn}}(i_*G^*,\Ga)\ar[u]    }$$}
Let us now explain the base change morphisms for the second line of \ref{synsyn1i}. As before we use the following notations: $$\begin{array}{rcll}D^{SYN}&:=&\cal Ext^1_{(X/\Sigma_\infty)_{SYN}}(i_*G,\cal O)&\hbox{in $\cal Crys((X/\Sigma_\infty)_{CRYS,syn},\O)$}\\
\overline D{}^{SYN}&:=&\iota_1^{-1}D^{SYN}&\hbox{in $\cal Crys((X/\Sigma_1)_{CRYS,syn},\O)$}\\
D^{syn}&:=&\pi_*D^{syn}&\hbox{in $\cal Crys((X/\Sigma_\infty)_{crys,syn},\O)$}\\
\overline D{}^{syn}&:=&\iota_1^{-1}\overline D{}^{syn}&\hbox{in $\cal Crys((X/\Sigma_1)_{crys,syn},\O)$}\\
D&:=&\epsilon_*D^{syn}&\hbox{in $\cal Crys((X/\Sigma_\infty)_{crys,et},\O)$}\\
\overline D&:=&\iota_1^{-1}D&\hbox{in $\cal Crys((X/\Sigma_1)_{crys,et},\O)$}\end{array}$$
We define similarly crystals $D'^{SYN}$, $\overline D{}'^{SYN}$, $D'^{syn}$, $\overline D{}'^{syn}$, $D'$, $\overline D{}'$ starting from $G'$ instead of $G$.  A base change morphism $f^*D{}^{syn}\to D'{}^{syn}$ has been described just above. Note that this also induces $f^*\overline D{}^{syn}\to \overline D'{}^{syn}$ (using $\iota_{1}f\simeq f\iota_1$), $f^*D\to D'$ (using $\epsilon f\to f\epsilon$) and $f^*\overline D\to \overline D'$ (using $\iota_1f\simeq f\iota_1$). From the latter and $\phi f\to f\phi$, we deduce the base change morphism $f^*Lie(D)\to Lie(D')$ considered in Prop. \ref{liesurj} \ref{liesurjii}. Recall now that $Lie^{syn}(D):=\epsilon^*Lie(D)$. Using $f\epsilon \to \epsilon f$, we next deduce the base change morphism $f^*Lie^{syn}(D)\to Lie^{syn}(D')$  implicitly hinted in Def. \ref{defSonephisyn} \ref{defSonephisyni}. Finally, a base change morphism $f^*i_*Lie^{syn}(D)\to i_*Lie^{syn}(D')$ for the right term in the second line of \ref{synsyn1i} follows using $if\to fi$. In order to get the expected base change morphisms for the second line of \ref{synsyn1i}, we have to check that the following natural square (where the horizontal arrows are defined using $can$ and $i\epsilon\to \epsilon i$)  is commutative:
\begin{eqnarray}\label{chlierrr}\xymatrix{D'^{syn}\ar[r]&i_*Lie^{syn}(D')\\ f^*D^{syn}\ar[u]\ar[r]&f^*i_*Lie^{syn}(D)\ar[u]}\end{eqnarray}
Even though this was also already implicit in Def. \ref{defSonephisyn} \ref{defSonephisyni}, we give some details here. We have to prove that the exterior square of the following diagram (where the arrows are the obvious ones) commutes:
$$\xymatrix{\epsilon^*D'\ar[r]&\epsilon^*i_*Lie(D')\ar[r]&i_*\epsilon^*Lie(D')\\ \epsilon^*f^*D\ar[u]\ar[r]&\epsilon^*f^*i_*Lie(D)\ar[r]\ar[u]&i_*\epsilon^*f^*Lie(D)\ar[u]\\
f^*\epsilon^*D\ar[u]\ar[r]&f^*\epsilon^*i_*Lie(D)\ar[r]\ar[u]&f^*i_*\epsilon^*Lie(D)\ar[u]}$$
Commutativity of the bottom left (resp. top right) square is tautological (resp. is immediate if one introduces $\epsilon^*i_*f^*Lie(D)$). Commutativity of the top left square is deduced from Prop. \ref{liesurj} \ref{liesurjii} using only pseudo-functoriality of the morphism $i$. Commutativity of the bottom right hand square is a formal consequence of the pseudo-functoriality of the isomorphism $i\epsilon\simeq \epsilon i$.

We have now checked that (\ref{chlierrr}) commutes. As a result, $f^*D^{syn}\to D'^{syn}$ also induces a base change morphism $f^*Fil^1D^{syn}\to Fil^1D'^{syn}$ so that we have described base change morphisms for each vertex of the commutative diagram \ref{synsyn1i}. It remains to check their compatibility with respect to the vertical arrows of \ref{synsyn1i}.
In the case of the middle one there is nothing to prove. The case of the left and right ones follow formally. The tenacious reader may check by himself that the base change morphisms just described above satisfy the composition constraint with respect to morphisms $X''\to X'\to X$.

Functoriality of the diagram in \ref{synsyn1ii} with respect to $X$ follows immediately since it is obtained from \ref{synsyn1i} by applying  $u_*l^{-1}$ while both $u$ and $l$ are pseudo-functorial with respect to $X$.
\begin{flushright}$\square$\end{flushright}





\begin{lem} \label{synsyn2} Assume that $X/\Sigma_1$ has local finite $p$-bases. Consider $G$ in $pdiv(X)$ and view it as group in $X_{syn}$. We use the simplified notation $\cal Ext^q$ for  $\cal Ext^q_{X_{syn}}$ or $\cal Ext^q_{X^\N_{syn}}$.

\debrom \item \label{synsyn2i} The modules $\smash{\cal Ext^1(l^{-1}G^*,\widetilde{\cal I}^{crys}_.)}$ and $\smash{\cal Ext^1(l^{-1}G^*,\cal O^{crys}_.)}$ are normalized and $\Z/p^.$-flat. Moreover Lem. \ref{synsyn1} \ref{synsyn1ii} induces compatible isomorphisms as follows over $(\smash{X_{syn}^\N,\O_.^{crys}})$:
$$\xymatrix{
\widetilde{Fil}{}^{1,crys}_.D\ar[rrr]^-1\ar[d]^{\wr}&&&\widetilde{Fil}{}^{0,crys}_.D\ar[d]^\wr\\
\cal Ext^1(l^{-1}G^*,\widetilde{\cal I}{}^{crys}_.)\ar[rrr]^-{\cal Ext^1(l^{-1}G^*,1)}&&&\cal Ext^1_{X^\N_{syn}}(l^{-1}G^*,\cal O^{crys}_.)}$$
\item \label{synsyn2ii} The following square of $(\smash{X_{syn}^\N,\O_.^{crys}})$ is commutative:
$$\xymatrix{\widetilde{Fil}{}^{1,crys}_.D\ar[d]^\wr\ar[rrr]^-\varphi&&&(F)_*\widetilde{Fil}{}^{0,crys}_.D\ar[d]^\wr\\
\cal Ext^1(l^{-1}G^*,\widetilde{\cal I}^{crys}_.)\ar[rrr]^-{\cal Ext^1(l^{-1}G^*,\varphi)}&&&(F)_*\cal Ext^1(l^{-1}G^*,\cal O^{crys}_.)\\}$$
\item \label{synsyn2iii} There is an isomorphism
$$\begin{array}{rl}\cal S^{1,\varphi}_{syn,.,X}(D(G))\simeq \cal Ext^1(\varpi^{-1}l^{-1}G^*,(\widetilde I{}_.^{crys},\O^{crys}_.,1,\varphi))& \hbox{in $Mod^{1,\varphi}(X^\N_{syn},\O^{crys}_.)$}\end{array}$$ where $\varpi$ denotes the canonical morphism $\smash{X^{\N,1,\varphi}_{syn}}\to
\smash{X^{\N}_{syn}}$ (Lem. \ref{1phitop}).
\finrom
\end{lem}
\noindent Proof. \ref{synsyn2i} By Lem. \ref{lemsyn} \ref{lemsyniv}, we already know that $D^{crys}_.\simeq \smash{\cal Ext^1(l^{-1}G^*,\cal O^{crys}_.)}$ is $\Z/p^.$-flat and normalized. The same holds for $\jj^*\smash{\cal Ext^1(l^{-1}G^*,\cal I^{crys}_.)}$, thanks to Lem. \ref{cartgen} \ref{cartgeniv} applied to the top exact sequence of Lem. \ref{synsyn1} \ref{synsyn1ii}.
In order to get the claimed compatible isomorphisms, it suffices to apply the functor $\jj^*$ to the left square in Lem. \ref{synsyn1} \ref{synsyn1ii} and to prove that the natural arrow
$$\xymatrix{\jj^* \cal Ext^1(l^{-1}G^*,\cal I^{crys}_.)\ar[r]&\cal Ext^1(l^{-1}G^*,\jj \cal I^{crys}_.)}$$ is an isomorphism. This, in turn, will follow from the following commutative diagram by the five lemma, once proven that its lines are exact:
$$\xymatrix{0\ar[r]&\cal Ext^1(G^*,\O)\ar[r]\ar[d]^\parallel&\cal Ext^1(G^*,\I^{crys}_{k+1})/p^k\ar[r]\ar[d]&\cal Ext^1(G^*,\O^{crys}_k)\ar[r]\ar[d]^\parallel&\cal Ext^1(G^*,\O)\ar[d]^\parallel\ar[r]&0
\\
0\ar[r]&\cal Ext^1(G^*,\O)\ar[r]&\cal Ext^1(G^*,\I^{crys}_{k+1}/p^k)\ar[r]&\cal Ext^1(G^*,\O^{crys}_k)\ar[r]&\cal Ext^1(G^*,\O)\ar[r]&0}$$
The first exact line is deduced from the first exact sequence of Lem. \ref{synsyn1} \ref{synsyn1ii} using Lem. \ref{cartgen} \ref{cartgeniv}. The second line on the other hand, is deduced from the exact sequence \begin{eqnarray}\label{uigiuhg} \xymatrix{0\ar[r]&\O\ar[r]&\smash{\cal I^{crys}_{k+1}}/p^k\ar[r]&\smash{\cal O^{crys}_k}\ar[r]&\O\ar[r]&0}\end{eqnarray} (apply Lem. \ref{cartgen} \ref{cartgeniv} to the exact sequence of Lem. \ref{lemcomp2} \ref{lemcomp2i}). Note that (\ref{uigiuhg}) remains exact after applying $\cal Ext^1(G^*,-)$ since $\cal Hom(G^*,\I^{crys}_k)$, $\cal Ext^2(G^*,\O)$ and $\cal Hom(G^*,\O)$ vanish by Cor. \ref{corpdivtop} \ref{corpdivtopi} and Lem. \ref{lemext3} \ref{lemext3i}.

\ref{synsyn2ii} Using Lem. \ref{cartgen} \ref{cartgeniii} and the definition of $\varphi$ (Prop. \ref{defphisyn}),  it is sufficient to prove that the following square commutes: $$\xymatrix{{Fil}{}^{1,crys}_.D\ar[d]^\wr\ar[rrr]^-{Fr}&&&(F)_*\widetilde{Fil}{}^{0,crys}_.D\ar[d]^\wr\\
\cal Ext^1(l^{-1}G^*,{\cal I}^{crys}_.)\ar[rrr]^-{\cal Ext^1(l^{-1}G^*,F)}&&&(F)_*\cal Ext^1(l^{-1}G^*,\cal O^{crys}_.)\\}$$
where $Fr$ is as in Def. \ref{defFrsyn} and $F$ denotes the Frobenius endomorphism of $\smash{\cal O_.^{crys}}$. We conclude by Lem. \ref{compFext} \ref{compFexti}, \ref{compFextiii} thanks, to the fact that the absolute Frobenius of $X$ is syntomic (choose locally a finite $p$-basis and argue as in the proof of \ref{lemcomp3} \ref{lemcomp3iv}).

\ref{synsyn2iii}  In view of Lem. \ref{lemext1} \ref{lemext1i}, this is just a reformulation of \ref{synsyn2i} and \ref{synsyn2ii}.
\begin{flushright}$\square$\end{flushright}



We are now in a position to prove the expected comparison theorem for $p$-divisible groups.

\begin{thm} \label{thmcomppdiv}
Let $X$ denote a diagram of $\cal Sch/\Sigma_1$ whose vertices have local finite $p$-bases and consider $G\in pdiv(X)$. There is a canonical isomorphism $$G_{p^.}\simeq \cal S_{syn,.,X}(D(G))\hbox{\hspace{1cm} in $D(X_{syn}^\N, \Z/p^.)$}$$
\end{thm}

\noindent Proof. Since $\smash{G_{p^.}}$ is concentrated in degree $0$, it suffices to construct the desired isomorphism in the case where $X$ is a single scheme and to check its compatibility with respect to the natural base change morphisms relative to a morphism $X'\to X$. Let thus $X$ denote a scheme with local finite $p$-bases over $\Sigma_1$. In virtue of Lem. \ref{synsyn2} \ref{synsyn2iii} and of the tautological isomorphism $$\cal S_{syn,.,X}(D(G))\simeq R\varpi_*\cal S_{syn,.,X}^{1,\varphi}(D(G))$$ it suffices to establish a canonical isomorphism in $X^\N_{syn}$: $$G_{p^.}\simeq R\varpi_*\cal Ext^1(\varpi^{-1}l^{-1}G^*,(\widetilde I{}_.^{crys},\O^{crys}_.,1,\varphi))$$
Let us compute the right hand side. Recall from Lem. \ref{1phimod} \ref{1phimodiv} that one has a canonical distinguished triangle $$\xymatrix{R\varpi_*\cal Ext^1(\varpi^{-1}l^{-1}G^*,(\widetilde I{}_.^{crys},\O^{crys}_.,1,\varphi))\ar[r]&\cal Ext^1(l^{-1}G^*,\widetilde I{}_.^{crys})\ar[r]^-{1-\varphi}&\cal Ext^1(l^{-1}G^*,\O^{crys}_.)\ar[r]^-{+1}&}$$
Since $\smash{\cal Ext^q(l^{-1}G^*,\mu_{p^.})}$ vanishes for $q=0,2$ (Lem. \ref{lemext1} \ref{lemext1i}, Cor. \ref{corpdivtop} \ref{corpdivtopi} and  Lem. \ref{lemext3} \ref{lemext3i}),
the second exact sequence in Lem. \ref{lemcomp2} \ref{lemcomp2iii} produces an exact sequence $$\xymatrix{0\ar[r]&\cal Ext^1(l^{-1}G^*,\mu_{p^.})\ar[r]&\cal Ext^1(l^{-1}G^*,\widetilde I{}_.^{crys})\ar[r]^-{1-\varphi}&\cal Ext^1(l^{-1}G^*,\O{}_.^{crys})\ar[r]&0&}$$
Comparing with the above, we get \begin{eqnarray}\label{isoza}\begin{array}{rrrrl}&&&R\varpi_*\cal Ext^1(\varpi^{-1}l^{-1}G^*,(\widetilde I{}_.^{crys},\O^{crys}_.,1,\varphi))\simeq \cal Ext^1(l^{-1}G^*,\mu_{p^.})&\hbox{in $Mod(X^\N_{syn},\Z/p^.)$}\end{array}\end{eqnarray}
Let us compute the right hand side. Thanks to Cor. \ref{corpdivtop} \ref{corpdivtopi}, we find that for $k$ fixed, the exact sequence $$\xymatrix{0\ar[r]& G^*_{p^k}\ar[r]& G^*\ar[r]^{p^k}&G^*\ar[r]& 0}$$ gives rise to an isomorphism between $G_{p^k}\simeq \cal Hom(G_{p^k}^*,\mu_{p^k})$ and $\cal Ext^1(G^*,\mu_{p^k})$. Letting $k$ vary, we get \begin{eqnarray}\begin{array}{rl}\label{isozb}G_{p^.}\simeq \cal Ext^1(l^{-1}G^*,\mu_{p^.})&\hbox{in $Mod(X^\N_{syn},\Z/p^.)$}\end{array}\end{eqnarray}
This ends the proof in the case of a single scheme $X$. In order to pass to the case of a diagram $X/\Delta$ of arbitrary type, it suffices to check that the isomorphisms (\ref{isoza}) and (\ref{isozb}) are both naturally functorial with respect to $X$. This causes no difficulty using the first isomorphism of Lem. \ref{lemext2} \ref{lemext2iii} (note that (\ref{isoza}) and (\ref{isozb}) both come from a morphism in the big topos).
\begin{flushright} $\square$ \end{flushright}

\begin{rem} The isomorphism of the theorem may be reformulated as an exact sequence over $(X_{syn}^\N, \Z/p^.)$ $$\xymatrix{0\ar[r]&G_{p^.}\ar[r]&\widetilde{Fil}{}^{1,crys}_.(D(G))\ar[r]^{1-\varphi}&
\widetilde{Fil}{}^{0,crys}_.(D(G))\ar[r]&0}$$
\end{rem}

\begin{cor} \label{corcomppdiv} Keep the assumptions of Thm. \ref{thmcomppdiv} and denote $\epsilon:X_{FL}\to X_{et}$. There is a canonical isomorphism: $$\begin{array}{rl}R\epsilon_*G_{p^.}\simeq \cal S_{et,.,X}(D(G))&\hbox{in $D(X_{et}^\N,\Z/p^.)$}\end{array}$$
\end{cor}
Proof. This follows immediately from Thm. \ref{thmcomppdiv} and Prop. \ref{compsynet} \ref{compsynetii}, given that $G_{p^.}$ is acyclic for the projection functor $X_{FL}\to X_{syn}$ (this follows easily from \cite{Gr1} Thm. 11.7 using \cite{BBM} Thm. 3.1.1).
\begin{flushright}$\square$\end{flushright}

\subsection{Semi-Abelian schemes over $C$} \label{sasoc} ~~ \\

We can finally prove our main result.

\begin{thm} \label{compss} Assume $A$ is a semi-Abelian scheme over $C$ whose restriction to $U$ is Abelian. Let  $\epsilon:C_{FL}\to C_{et}$ denote the canonical morphism. There is a canonical isomorphism $$\begin{array}{rl}R\epsilon_* R\underline \G^ZA_{p^.}\simeq \cal S_{et,.,C^\s}(-Z)(D_{C^\s}(A))&\hbox{in $D(C_{et}^\N,\Z/p^.)$}\end{array}$$
\end{thm}
\noindent Proof. Recall the diagram of $p$-divisible groups $H\in pdiv(J^+)$ defined in Def. \ref{defH}. Denoting $D:=D_{J^+}(H)$ the associated diagram of Dieudonn\'e crystals we have the following series of canonical  isomorphisms in $D(C_{et}^\N,\Z/p^.)$: $$\begin{array}{rcll}R\epsilon_* R\underline \G^ZA_{p^.}&\simeq& Rm_*Sma R\underline \G^{Z_J}(J,R\epsilon_*H_{p^.})&(Prop.\ \ref{devissageFL})\\
&\simeq & Rm_*Sma R\underline \G^{Z_J}(J,\cal S_{et,.,J^+}(D))&(Cor.\ \ref{corcomppdiv})\\
&\simeq & Rm_* Sma\cal S_{et,.,J^\s}(-Z_J)(D_{|J^\s})&(\hbox{Prop.\ \ref{tdspe} \ref{tdspeii} applied to $D_{|J}$})\\
&\simeq &  Rm_* Sma\cal S_{et,.,J^\s}(-Z_J)(D_{C^\s}(A)_{|J^\s}) &(\hbox{see below})\\
&\simeq &Rm_* \cal S_{et,.,J^\s}(-Z_J)(D_{C^\s}(A)_{|J^\s})& (Lem.\ \ref{lemsma} \ref{lemsmaii})\\
&\simeq & \cal S_{et,.,C^\s}(-Z)(D_{C^\s}(A))& (Prop.\ \ref{MVsyn} \ref{MVsyniii})\end{array}$$
Let us explain the fourth isomorphism in details. According to Lem. \ref{lemDCH} and Lem. \ref{Sisexact}, we have a canonical distinguished triangle  \begin{eqnarray}&&\label{triuo}\xymatrix{\cal S_{et,.,J^\s}(-Z_J)(D_{|J^\s})\ar[r]&\cal S_{et,.,J^\s}(-Z_J)(D_{C^\s}(A)_{|J^\s})\ar[r]&\cal S_{et,.,J^\s}(-Z_J)(D'_{|J^\s})\ar[r]^-{+1}&}\end{eqnarray} in $D(J^\N_{et},\Z/p^.)$. Here $D':=D(H')$ is the Dieudonn\'e crystal associated to the following $p$-divisible group on $J^+=(Z_v\to C_v\leftarrow U_v\to U)$: $$H'=(\Qp/\Zp\otimes \G_{v,|Z_v}\leftarrow \Qp/\Zp\otimes \G_{v}\rightarrow 0 \leftarrow 0)$$
We claim that the third term in (\ref{triuo}) vanishes after applying the functor $Sma$ (Def. \ref{defsma} \ref{defsmaii}). Let us explain this. From Prop. \ref{tdspe} \ref{tdspeiii}  and Cor. \ref{corcomppdiv}  we have $$\begin{array}{rcl}\cal S_{et,.,J^\s}(-Z_J)(D'_{|J^\s})&\simeq& R\underline \G^{Z_{J}}(J,\cal S_{et,.,J^+}(D'))\\ &\simeq & R\underline \G^{Z_J}(J,\epsilon_*H'_{p^.})\end{array}$$
By Lem. \ref{soritesvan} \ref{soritesvanii}, \ref{soritesvaniii} we have a distinguished triangle $$\xymatrix{R\underline \G^{Z_J}(J,\epsilon_*H'_{p^.})\ar[r]&\epsilon_*H'_{|J,p^.}\ar[r]&z_{J,*}\epsilon_*H'_{|Z_J,p^.}\ar[r]^-{+1}&}$$
and thus an isomorphism $$R\underline \G^{Z_J}(J,\epsilon_*H'_{p^.})\simeq ( j_{v,!}\epsilon_* \G_{v,|U_v}/p^.\rightarrow 0\leftarrow 0)$$ in $D(J^{\N}_{et},\Z/p^.)$ (recall that $J=(C_v\leftarrow U_v\to U)$). The claim then follows from Lem. \ref{lemsma} \ref{lemsmaiii}.
\begin{flushright}$\square$\end{flushright}


It might be worth to state a down to earth (weaker)  version of this result.

\begin{cor} \label{compssweak} There is a canonical distinguished triangle $$\xymatrix{R\G^Z(C,T_p(A))\ar[r]&R\G(C^\s/\Sigma_\infty,Fil^1 D_{C^\s}(A)(-Z))\ar[r]^-{1-\varphi}& R\G(C^\s/\Sigma_\infty,D_{C^\s}(A)(-Z))\ar[r]^-{+1}&}$$
where $p.\varphi$ is induced by the Frobenius endomorphism of the third term.
\end{cor}
Proof. Apply $Rl_*$ and $R\G(C_{et},-)$  to Thm. \ref{compss} and then apply Lem. \ref{1phimod} \ref{1phimodiv}. 
\begin{flushright}$\square$\end{flushright}



\section{Index of notations}



\noindent \hspace{-1cm} $X/\Delta$ or $(X_\delta)_{\delta\in \Delta}$: a diagram of type $\Delta$, Def. \ref{defdiagto} \ref{defdiagtoi}.

\noindent \hspace{-1cm} $Diag(\C)$: the category of diagrams of a category $\cal C$, Def. \ref{defdiagto} \ref{defdiagtoiii}.

\noindent \hspace{-1cm} $\cal F_{fib}$, resp. $\cal F_{cof}$: the fibered, resp. cofibered, category associated to a contravariant pseudo-functor $\cal F$ on a category $\cal B$, Lem. \ref{eqfib}.

\noindent \hspace{-1cm} $\cal F^{diag}$, $\cal F^{codiag}$: extensions of $\cal F$ to $Diag(\cal B)$, Lem. \ref{lemdiagcodiag}.

\noindent \hspace{-1cm} $\mathfrak Cat$, resp. $\mathfrak Pre$: the $2$-category of categories, resp. pretopologies, Def. \ref{deftopfib}.

\noindent \hspace{-1cm} $Mod(\cal T(-),A)$, $Kom^*(\cal T(-),A)$ and $D^*(\cal T(-),A)$: the category of modules, complexes and its derived category viewed as pseudo-functors, Lem. \ref{lemDfib1}.

\noindent \hspace{-1cm} $\cal T^\N$, $\iota_k:\cal T\to \cal T^\N$ and $l:\cal T^\N\to \cal T$: the topos of projective systems and the associated morphisms, Def. \ref{defprojtop}, Lem. \ref{lemprojtop} \ref{lemprojtopiii}, \ref{lemprojtopiv}.

\noindent \hspace{-1cm} $\cal Sch$: the category of schemes, Sect. \ref{usualtop}.

\noindent \hspace{-1cm} $\cal Sch^\s$: the category of fine log schemes, Sect. \ref{usualtop}.

\noindent \hspace{-1cm} $X^\s=(X,M_X\to O)$: an object of $\cal Sch^\s$, Sect. \ref{usualtop}.

\noindent \hspace{-1cm} $TOP^\s(X^\s)$, $TOP(X^\s)$ and $top(X^\s)$: the $\s$-big, big and small $top$-site of a log scheme $X^\s$, Def. \ref{defusualtop}.

\noindent \hspace{-1cm} $\smash{X^\s_{TOP^\s}}$, $\smash{X^\s_{TOP}}$ and $\smash{X^\s_{top}}$: their associated topoi, Def. \ref{defusualtop}.

\noindent \hspace{-1cm} $\pi$, $r$: the projection weak morphism between the $\s$-big, big and small topos, its right section (\ref{defp}).

\noindent \hspace{-1cm} $\epsilon$: the weak morphism of change of topology, (\ref{defepsilon}).

\noindent \hspace{-1cm} $fl$, $syn$, $et$, $zar$: the usual instances of $top$, Sect. \ref{usualtopusual}.

\noindent \hspace{-1cm} $\cal O$: the structural ring of the usual topoi,  Def.  \ref{defusualO}.

\noindent \hspace{-1cm} $\Sigma_k:=Spec(\Z/p^k)$, $\Sigma_\infty:=Spf(\Zp)$, Sect. \ref{sectioncrys}.

\noindent \hspace{-1cm} $\smash{\cal Sch^\s_{p,nil}}=\smash{\cal Sch^\s/\Sigma_\infty}$, Sect. \ref{sectioncrys}.

\noindent \hspace{-1cm} $\smash{CRYS^\s_{top}(X^\s/\Sigma_k)}$, $\smash{CRYS^\s_{top}(X^\s/\Sigma_k)}$ and $\smash{CRYS^\s_{top}(X^\s/\Sigma_k)}$: the $\s$-big, big and small crystalline sites, Def. \ref{defcrystop} \ref{defcrystopiii}, \ref{defcrystopiv} and \ref{defcrystopv}.

\noindent \hspace{-1cm} $\smash{(X^\s/\Sigma_k)_{CRYS^\s,top}}$,  $\smash{(X^\s/\Sigma_k)_{CRYS,top}}$ and  $\smash{(X^\s/\Sigma_k)_{crys,top}}$: the associated topoi,  Def. \ref{defcrystop} \ref{defcrystopiii}, \ref{defcrystopiv} and \ref{defcrystopv}.

\noindent \hspace{-1cm} $\pi$, $r$: the projection weak morphism between the $\s$-big, big and small crystalline topos, its right section, (\ref{defpcrys}).

\noindent \hspace{-1cm} $\iota_{k,k'}$, $\iota_k:=\iota_{k,\infty}$: the morphisms of change of $k$,  (\ref{defiotacrys}).

\noindent \hspace{-1cm} $\epsilon$: the weak morphism of change of topology, (\ref{defepsiloncrys}).

\noindent \hspace{-1cm} $i$ and $u$: the morphism of immersion of the usual topos into the crystalline one and its right retraction when it exists (\ref{defiu}).

\noindent \hspace{-1cm} $(lift)$: the lifting property for $top$, proof of Lem. \ref{lemiucocont}.

\noindent \hspace{-1cm} $\smash{f_{T^\s}}$, $\smash{\lambda_{T^\s}}$: the localization morphism and the realization weak morphism, $(\ref{deflambda})$.

\noindent \hspace{-1cm} $\smash{F_{|T^\s}}$, $\smash{F_{T^\s}}$: the restriction and realization of a crystalline sheaf, Def. \ref{defreal}.

\noindent \hspace{-1cm} $\cal O=\smash{\cal O_{X^\s/\Sigma_k}}$ and $\mathbb G_a=i_*\O$: the usual rings of the crystalline topos, Def. \ref{defcrysO} \ref{defcrysOi}, \ref{defcrysOii}.

\noindent \hspace{-1cm} $\smash{\cal O_k^{crys}}$, $\cal O^{crys}=\smash{\cal O_\infty^{crys}}$: the crystalline rings of the usual topos, Def. \ref{defcrysO}.

\noindent \hspace{-1cm} $\smash{f_{CRYS^\s}}$, $\smash{f_{CRYS}}$, $\smash{f_{crys}}$: (only in Sect. \ref{scsf}) the crystalline functoriality weak morphisms.

\noindent \hspace{-1cm} $Real(\smash{CRYS^\s_{top}(X^\s/\Sigma)})$: the category of realizations of $top$ crystalline sheaves resp. satisfying $top'$ descent, Lem.  \ref{lemreal} \ref{lemrealiii}.

\noindent \hspace{-1cm} $\cal Crys(\smash{(X^\s/\Sigma)_{CRYS^\s,top}},\O)$: the category of crystals, Def. \ref{defcrystal}.

\noindent \hspace{-1cm} $flf(X)$, $pdiv(X)$: the category of finite locally free groups and $p$-divisible groups over $X$, Lem. \ref{gpflfexact} \ref{gpflfexactii}, \ref{gpflfexacti}.

\noindent \hspace{-1cm} $\cal M(X)$: the category of $1$-motives over $X$, Lem. \ref{1motex}.

\noindent \hspace{-1cm} $\Phi$: the component group of a semi-Abelian scheme, Lem. \ref{components}.

\noindent \hspace{-1cm} $Ab(E):=Mod(E,\Z)$: the category of Abelian groups of a topos $E$, Lem. \ref{corpdivtop}.

\noindent \hspace{-1cm} $flf_{TOP}(X)$, $flf_{top}(X)$, $flf_{CRYS,top}(X/\Sigma)$, $flf_{crys,top}(X/\Sigma)$: incarnations of $flf(X)$, Lem. \ref{corpdivtop}.

\noindent \hspace{-1cm} $pdiv_{TOP}(X)$, $pdiv_{top}(X)$, $pdiv_{CRYS,top}(X/\Sigma)$, $pdiv_{crys,top}(X/\Sigma)$: incarnations of $pdiv(X)$, Lem. \ref{corpdivtop}.

\noindent \hspace{-1cm} $(E_.,A_.)$: the ringed total topos associated to a variable topos  on $\N$, Sect. \ref{ARmod}.

\noindent \hspace{-1cm} $\j$: a certain endomorphism of $(E_.,A_.)$, Def. \ref{defj}.

\noindent \hspace{-1cm} $Mod(E_{.},A_{.})$, $D(E_{.},A_{.})$: the category of modules and its derived category, Def. \ref{defnorm} \ref{defnormi}, \ref{defnormii}.

\noindent \hspace{-1cm} $\cal Sch_p$: the category of $p$-adic schemes, Def. \ref{defpadic} \ref{defpadici}.

\noindent \hspace{-1cm} $\smash{\cal Sch^\s_p}$: the category of $p$-adic log schemes, Def. \ref{defpadic} \ref{defpadicii}.

\noindent \hspace{-1cm} $X_k$ (resp. $\smash{X^\s_k}$): the reduction mod $p^k$ of $X$, (resp. $X^\s$), comment after Def. \ref{defpadic}.

\noindent \hspace{-1cm} $TOP^\s(X^\s)$, $TOP(X^\s)$ and $top(X^\s)$: the $\s$-big, big and small $top$ site of a $p$-adic log scheme $X^\s$, Def. \ref{usualtoppadic}, \ref{usualtoppadici}.

\noindent \hspace{-1cm} $\smash{X^\s_{TOP^\s}}$, $\smash{X^\s_{TOP}}$ and $\smash{X^\s_{top}}$: the associated topoi, Def. \ref{usualtoppadic} \ref{usualtoppadici}.

\noindent \hspace{-1cm} $\smash{X^\s_{.,TOP^\s}}$, $\smash{X^\s_{.,TOP}}$ and $\smash{X^\s_{.,top}}$: the $\s$-big, big and small $top$ topoi of $\smash{X^\s_.}$, Def. \ref{usualtoppadic} \ref{usualtoppadicii}.

\noindent \hspace{-1cm} $\pi$, $r$: the weak morphism of projection between the $\s$-big, big and small topos, its right section, comment below Def. \ref{usualtoppadic}.

\noindent \hspace{-1cm} $\epsilon$: the weak morphism of change of topology, comment below Def. \ref{usualtoppadic}.

\noindent \hspace{-1cm} $\iota_{k,.}$, $\iota_k$, $l$: the weak morphisms relating the $top$ topoi of $\smash{X_k^\s}$, $\smash{X^\s_.}$ and $X^\s$, (\ref{iotalTOP}) and comment below.

\noindent \hspace{-1cm} $\O$: the structural ring of the $top$ topoi of $X^\s$ or $\smash{X^\s_.}$: Def. \ref{defOM}, \ref{defOMi}.

\noindent \hspace{-1cm} $M_X$ and $M_{X,.}$: the monoids of $X^\s$ and $\smash{X^\s_.}$: Def. \ref{defOM}, \ref{defOMii}.

\noindent \hspace{-1cm} $Mod_{prop}(E,\O)$ for $prop\in \{qcoh,lf,lfft\}$ and $E\in \{\smash{X^\s_{TOP^\s}},
\smash{X^\s_{TOP}}, \smash{X^\s_{TOP}}\}$, $X^\s$ a log scheme: Def. \ref{defqcohsch} and Lem.-Def. \ref{lemqcohsch} \ref{lemqcohschi}, \ref{lemqcohschii}.

\noindent \hspace{-1cm} $V_X((f_i)_i)$: the closed subscheme of a scheme $X$ defined by a set of sections of $\O$, Rem. \ref{remqcohsch1} \ref{remqcohsch1iii}.

\noindent \hspace{-1cm} $Mod_{norm,qcoh}(X_{.,et},\O)$, $Mod_{qcoh}(X_{et},\O)$ for $X$ a $p$-adic scheme: Rem. \ref{remqcohpadic} \ref{remqcohpadiciii}.

\noindent \hspace{-1cm} $D^+_{Lnorm,qcoh}(X_{.,et},\O)$, $D^+_{Lqcoh}(X_{et},\O)$ for $X$ a flat $p$-adic scheme: Rem. \ref{remLqcohpadic}.

\noindent \hspace{-1cm} $\smash{(X^\s/\Sigma_.)_{CRYS^\s,top}}$, $\smash{(X^\s/\Sigma_.)_{CRYS,top}}$ and $\smash{(X^\s/\Sigma_.)_{crys,top}}$: the $\s$-big, big and small crystalline $top$ topos of $(X/\Sigma_.)$.  Def. \ref{defprocrystop}.

\noindent \hspace{-1cm} $\iota_{k,.}$, $\iota_k$, $l$: the weak morphisms relating the crystalline $top$
topoi of $(X^\s/\Sigma_k)$, $(X^\s/\Sigma_.)$ and $(X^\s/\Sigma_\infty)$, (\ref{diagiotalcrys}).

\noindent \hspace{-1cm} $\pi$, $r$, $\epsilon$, $i$, $u$ in this context: comments of (\ref{diagiotalcrys}) and (\ref{defproiu}).

\noindent \hspace{-1cm} $\smash{f_{T^\s_.}}$, $\smash{\lambda_{T^\s_.}}$: the localization morphism and the realization weak morphism in this context, (\ref{defprolambda}).

\noindent \hspace{-1cm} $\smash{F_{.,|T^\s_.}}$ and $\smash{F_{.,T^\s_.}}$: the restriction and realization of $F_.$, comment below (\ref{defprolambda}).

\noindent \hspace{-1cm} $Mod_{prop}(E,\O)$ for $prop\subset \{norm,qcoh,lf,lfft\}$ and $E\in \{\smash{(X^\s/\Sigma_.)_{CRYS^\s,top}}, \smash{(X^\s/\Sigma_.)_{CRYS,top}},$ $\smash{(X^\s/\Sigma_.)_{crys,top}}\}$: Lem.-Def. \ref{lemlimcrys} \ref{lemlimcrysii} and Rem.  \ref{lemlcrys} \ref{lemlcrysii}.

\noindent \hspace{-1cm} $\cal Crys_{prop}(E,\O)$ for $prop\subset \{norm,qcoh,lf,lfft\}$ and $E\in \{\smash{(X^\s/\Sigma_.)_{CRYS^\s,top}}, \smash{(X^\s/\Sigma_.)_{CRYS,top}},$ $\smash{(X^\s/\Sigma_.)_{crys,top}}\}$:  Rem.  \ref{lemlcrys} \ref{lemlcrysiv}.

\noindent \hspace{-1cm} $F^{(X/Y)}$: the relative Frobenius, Def. \ref{defpb} \ref{defpbi}.

\noindent \hspace{-1cm} $\mathbb A^d$: the affine space of dimension $d$, Def. \ref{defpb} \ref{defpbiii}.

\noindent \hspace{-1cm} $(\mathbb A^e,\mathbb N^e)$: the affine space of dimension $e$ with its canonical log structure, Def. \ref{defpb} \ref{defpbiv}.

\noindent \hspace{-1cm} $\Omega_{X^\s/Y^\s}$: the module of logarithmic differentials of a morphism of $p$-adic schemes,  (\ref{isopbdiff}).

\noindent \hspace{-1cm} $Emb^\s$, $Emb^{\s,glob}$: the respective category of local and global embeddings, Def. \ref{defemb1}, \ref{defemb1i}.

\noindent \hspace{-1cm} $\cal Sch^{\s,slfpb}/\Sigma_1$, $\cal Sch^{slfpb}/\Sigma_1$: Def. \ref{defemb1var} \ref{defemb1vari}.

\noindent \hspace{-1cm} $Emb^*$ for $*\subset \{\s,glob,lfpb\}$: Def. \ref{defemb1var} \ref{defemb1vari}, \ref{defemb1varii}.

\noindent \hspace{-1cm} $T^\s=D(X^\s,Y^\s)$ (resp. $D^n(X^\s,Y^\s)$): the logarithmic divided power envelope (resp. of order $n$) of $X^\s\to Y^\s$, Sect. \ref{theringP}.

\noindent \hspace{-1cm} $\cal P_T^{\s(i)}$ (resp. $\cal P_T^{\s(i),n}$): the $i$-th sheaf of $dp$-principal parts (resp. of order $n$) of $X^\s$ inside $Y^\s$, viewed as an algebra of $(T_{et},\O)$ via $d_0$, Def. \ref{lemdefP2} \ref{lemdefP2i}.

\noindent \hspace{-1cm} $\Omega_T$: the pullback of $\Omega_{Y^\s/\Sigma_k}$, (\ref{defdiffT}).

\noindent \hspace{-1cm} $\cal D^\s$: the algebra of $dp$-differential operators of $X^\s$ inside $Y^\s$. Def. \label{defDiff}.

\noindent \hspace{-1cm} $K$: the curvature of a connection $\nabla$, Lem. \ref{defnablaint} \ref{defnablainti}.

\noindent \hspace{-1cm} $\theta=(\theta_n)$: the Taylor morphism of an integrable connection $\nabla$. (\ref{deftheta}).

\noindent \hspace{-1cm} $\varepsilon$: the hyper $dp$-stratification of an integrable quasi-nilpotent connection. Prop. \ref{hyperstrat}.

\noindent \hspace{-1cm}   $\nabla$-$Mod_{prop}(T^\s)=\nabla$-$Mod^{pd}_{prop}(X^\s,Y^\s)$ with $prop\subset \{qcoh,lf,lfft\}$: the category of $prop$ modules with quasi-nilpotent integrable connection, Def. \ref{defnablamod}.

\noindent \hspace{-1cm} $Hdp(T^\s)=Hdp^{dp}(X^\s,Y^\s)$: the category of hyper $dp$-differential operators, Def. \ref{defHdp}.

\noindent \hspace{-1cm} $\Omega^\bullet_{T^\s}$, $\Omega^\bullet_{T^\s}(M)$: the de Rham complex, Def. \ref{dR} \ref{dRi}, Def. \ref{dRM}.

\noindent \hspace{-1cm} $L_{T^\s}$: the linearization functor with value in $\nabla\hbox{-}Mod(T^\s)$, Lem. \ref{defLT}.

\noindent \hspace{-1cm} $e_M$: the exact structure of modules on $\cal Crys_{prop}((X^\s/\Sigma)_{crys,et},\O)$, $\Sigma=\Sigma_.$ or $\Sigma_\infty$, Lem. \ref{mexactcrys}.

\noindent \hspace{-1cm} $\nabla\hbox{-}Mod^{pd}_.(X^\s,Y^\s)$: explanations before Lem. \ref{lemflatlift2}.

\noindent \hspace{-1cm} $e$: the exact structure of crystals on $\cal Crys_{prop}((X^\s/\Sigma)_{crys,et},\O)$, $\Sigma=\Sigma_.$ or $\Sigma_\infty$, Prop. \ref{exactcrysk} \ref{exactcryskiii}, Prop. \ref{exactcrysinfty} \ref{exactcrysinftyiii}.

\noindent \hspace{-1cm} $M(-h)$ for $h\in M_E/\smash{M_E^\times}$: Def. \ref{deftwlogtop} \ref{deftwlogtopi}, \ref{deftwlogtopiii}

\noindent \hspace{-1cm} $Supp(M_X)$, $Center(M_X)$: the support and center of a fine log structure on a scheme $X$, Lem.-Def. \ref{suppcenter}.

\noindent \hspace{-1cm} $M_{X/\Sigma}\to O_{X/\Sigma}$: the log structure of the ringed small crystalline topos induced by that of $X^\s$, explanations before Lem. \ref{injtwist2}.

\noindent \hspace{-1cm} $E^{1,\varphi}$, $0$, $1$, $\varpi$: the topos of diagrams of $E$ of type $1\leftleftarrows 0$ and the natural morphisms relating it to the topos $E$, Lem. \ref{1phitop}.

\noindent \hspace{-1cm} $Mod^{1,\varphi}(E,A)=Mod^{1,\varphi}(E,(A,F))$: the category of $(1,\varphi)$-modules, Lem. \ref{1phimod}.

\noindent \hspace{-1cm} $\smash{T^\s_{[.]}}=(\smash{U^\s_{[.]}/X^\s,T^\s_{[.]}})$: a semi-simplicial $p$-adic $dp$-thickening, viewed as a semi-simplicial crystalline sheaf, Sect. \ref{prelimemb}.

\noindent \hspace{-1cm} $\O^{crys}_.$: the crystalline ring of $\smash{X^{\s,\N}_{[.]}}$ or $\smash{U^{\s,\N}_{[.]}}$, also viewed as a ring of $\smash{T^{\s}_{[.],.}}$ in the latter case,  (\ref{descentdiaget}).

\noindent \hspace{-1cm} $\lambda_{T^\s_{[.],.}}$, $f_{T^\s_{[.],.}}$,  $f_{T^\s_{[.],.}/U^\s_{[.],.}}$, $f_{U^\s_{[.],.}}$, $\iota$: the morphisms of topoi associated to $\smash{T^\s_{[.]}}$, (\ref{descentdiaget}).

\noindent \hspace{-1cm} $Emb^*_F$ (resp. $Emb^*_{F,div}$) for $*\subset \{\s,lfpb,glob\}$: categories of embeddings with Frobenius lifts (resp. and effective log divisors), Def. \ref{defemb2} \ref{defemb2i}, \ref{defemb2iii}.

\noindent \hspace{-1cm} $HR^{*,et}_F$ (resp. $HR^{*,et}_{F,div}$) or  $HR^{*,crys}_F$ (resp. $HR^{*,crys}_{F,div}$) for $*\subset \{\s,lfpb\}$: categories of semi-simplicial embeddings with  Frobenius lifts (resp. and effective log divisors) inducing a hypercovering in the small $et$ topos or in the $et$ small crystalline topos, Def. \ref{defemb2} \ref{defemb2ii}, \ref{defemb2iii}.

\noindent \hspace{-1cm} $Emb^{\s,ex}_F$, $HR^{\s,et,ex}_F$ and $HR^{\s,crys,ex}_F$:
some more categories of embeddings or semi-simplicial embeddings with Frobenius lifts, Lem. \ref{lememb3}.

\noindent \hspace{-1cm} $\cal B_0$ (resp. $\cal B^\s_0$, $\cal C^\s_0$): a certain category of diagrams of schemes (resp. log schemes with effective log divisors, resp. semi-simplicial embeddings with Frobenius lifts and effective log divisors),  Def. \ref{defB0} \ref{defB0i}, \ref{defB0ii}.

\noindent \hspace{-1cm} $\cal F^{emb}$, $\cal F^{glob}$, $\cal F^{loc}$, $\cal K^{emb}$, $\cal K^{loc}$, $\cal K^{loc}{}'$, $\cal D^{loc}{}'$: some pseudo-functors occurring in the globalization of a colax morphism, Sect. \ref{paraprelimS}.

\noindent \hspace{-1cm} $\cal S^{emb}$, $\cal S^{loc}$, $\cal S^{glob}$: some variants of the colax morphism in question, Sect. \ref{paraprelimS}.

\noindent \hspace{-1cm} $F=F_{X^\s}$: the absolute Frobenius endomorphism, Sect. \ref{sectintrocartier}.

\noindent \hspace{-1cm} $F_0U^\s$, $F_0(U^\s,T^\s)$ for $U^\s$ in $TOP^\s(X^\s)$, $(U^\s,T^\s)$ in $CRYS^\s_{top}(X/\Sigma_1)$, Sect. \ref{sectintrocartier}.

\noindent \hspace{-1cm} $F^{(-/X^\s)}:F\to id$: the relative Frobenius, viewed as a morphism between endomorphisms of (crystalline or usual) $top$ topoi, Def. \ref{defFrel} \ref{defFreli}, \ref{defFrelii}.

\noindent \hspace{-1cm} $C^{-1}:F\to iu$, $nat:iu\to id$: the morphisms factorizing the relative Frobenius, Def. \ref{defCinv} \ref{defCinvi}, \ref{defCinvii}.

\noindent \hspace{-1cm} $(F)$, $(C^{-1})$, $(nat)$: the induced ring homorphisms, viewed as morphisms of ringed topoi, Def. \ref{defphitop} \ref{defphitopi}, \ref{defphitopiii}.

\noindent \hspace{-1cm} $\phi$: the ringed variant of $u$, using the inverse Cartier ring homomorphism, Def. \ref{defphitop} \ref{defphitopii}, \ref{defphitopiii}.

\noindent \hspace{-1cm} $F^{(-/X^\s)}:F\to (F)$, $C^{-1}:F\to i\phi$: the ringed variants of the relative Frobenius and the inverse Cartier morphism, Lem. \ref{cartierlin} \ref{cartierlini}, \ref{cartierlinii}.

\noindent \hspace{-1cm} $\cal P^{\s(1)}_X$ (resp. $\cal P^{(1)}_X$): the $1$st sheaf of $dp$-principal parts of $X^\s$ (resp. $X$), Sect. \ref{cartiergen}.

\noindent \hspace{-1cm} $\cal D^\s$ (resp. $\cal D$): the algebra of $dp$-differential operators of $X^\s$ (resp. $X$), Sect. \ref{cartiergen}.

\noindent \hspace{-1cm} $\cal P^{(1),F}$ and $\cal D^{(1),F}$: the structural ring of $X\times_{F,X,F}X$ and it dual, Def. \ref{defKprings} \ref{defKpringsi}, \ref{defKpringsii}.

\noindent \hspace{-1cm} $\cal Crys^F((X/\Sigma_1)_{crys,et},\O)$: the category of crystals with trivial $p$-curvature, Lem. \ref{defKp}.

\noindent \hspace{-1cm} $\smash{\phi_{X^\s}}$, $\phi_X$: the morphism of ringed topoi $\phi$ for $X^\s$ and $X$, (\ref{defo}).

\noindent \hspace{-1cm} $o$: the functoriality weak morphism induced by the forgetful morphism $X^\s\to X$, (\ref{defo}).

\noindent \hspace{-1cm} $\smash{Mod_{\underline t-fr}(X_{et},\O)}$, $\smash{\cal Crys_{\underline t-fr}}((X/\Sigma_1)_{crys,et},\O)$, $\smash{\cal Crys_{\underline t-fr}}((X^\s/\Sigma_1)_{crys,et},\O)$: the categories of $\underline t$-torsion free modules and crystals, Def. \ref{deftfr} \ref{deftfri}, \ref{deftfrii}.

\noindent \hspace{-1cm} $\cal DC(X^\s)$: the category of Dieudonn\'e crystals, Def. \ref{defDC} \ref{defDCi}.

\noindent \hspace{-1cm} $\cal DC_1(X^\s)$: the category of truncated Dieudonn\'e crystals of level $1$, Def. \ref{defDC} \ref{defDCii}.

\noindent \hspace{-1cm} $(\overline D,\overline f,\overline v)$: the reduction mod $p$ of $(D,f,v)$ for $X^\s$ locally embeddable, Lem. \ref{exDC} \ref{exDCii}.

\noindent \hspace{-1cm} $\cal DC_{CRYS^\s,top}(X^\s)$, $\cal DC_{CRYS,top}(X^\s)$, $\cal DC_{crys,top}(X^\s)$: various incarnations of $DC(X^\s)$: Rem. \ref{eqDC}.

\noindent \hspace{-1cm} $Lie(D):=Lie(\overline D):=Coker\, \overline v$ for $D$ in $\cal DC(X^\s)$, $X^\s$ with local $p$-bases: Def. \ref{deflie} \ref{defliei}, \ref{deflieii}.

\noindent \hspace{-1cm} $\smash{can_{\overline D}}:D\to i_*Lie(\overline D)$, $can_D:D\to i_*Lie(D)$: the natural morphism built using the inverse Cartier morphism, Def. \ref{deflie} \ref{defliei}, Prop. \ref{liesurj}, \ref{liesurji}.

\noindent \hspace{-1cm} $Fil^h\overline D$, $Fil^hD$, $h=0,1$: the $mod$ $p$ Hodge filtration on the crystalline  modules underlying $\overline D$ and $D$, Def. \ref{deflie} \ref{defliei}, \ref{deflieii}.

\noindent \hspace{-1cm} $ch_f:f^*Lie(D)\to Lie(f^*D)$: the base change isomorphism along $f$, Prop. \ref{liesurj}, \ref{liesurjii}.

\noindent \hspace{-1cm} $Hdp_{norm}(T^\s_.,\O)$: the category of normalized modules and hyper $dp$-differential operators, (\ref{functorLsanspoint}).

\noindent \hspace{-1cm} $L_.$, $L$: the linearization functor with respective values in $\cal Crys_{norm}((X^\s/\Sigma_.)_{crys,et},\O)$, \break $\cal Crys((X^\s/\Sigma_\infty)_{crys,et},\O)$ (\ref{functorLsanspoint}).

\noindent \hspace{-1cm} $M_.$, $\smash{M_{T^\s_.}}$, for $M$ in $Mod((X^\s/\Sigma_\infty),\O)$: explanations below (\ref{functorLsanspoint}).

\noindent \hspace{-1cm} $\smash{\Omega^\bullet_{T^\s_.}}(M)$, for $M$ in $\cal Crys\smash{((X^\s/\Sigma_\infty)_{crys,et},\O)}$:   explanations below (\ref{functorLsanspoint}).

\noindent \hspace{-1cm} $\smash{can_L:L(D_{T^\s_.})}\to i_*Lie(D)$, $\smash{Fil^1L(D_{T^\s_.})}:=Ker\, can_L$: Prop. \ref{complin} \ref{complini}.

\noindent \hspace{-1cm} $aug:D\to \smash{L(D_{T^\s_.})}$:  Prop. \ref{complin} \ref{complini}.

\noindent \hspace{-1cm} $Kom(T_{.,et}^\s,\O_.^{cris})$ (resp. $Fil^{0,1}Kom(T_{.,et}^\s,\O_.^{cris})$): the categories of complexes (resp. endowed with a one step filtration), Def. \ref{deffil}.

\noindent \hspace{-1cm} $\smash{{Fil}^._{.,T^\s}(-h)(D)},\, \smash{{Fil}^.\Omega^\bullet_{.,T^\s}(-h)(D)},\, \smash{{Fil}^.L\Omega^\bullet_{.,T^\s}(-h)(D)}$: the respective filtered twisted realization of $D$, twisted de Rham complex, twisted linearized de Rham complex, associated to a given global embedding with Frobenius lift and effective log divisor, Def. \ref{deffil} \ref{deffili}, \ref{deffilii}, \ref{deffiliii}, or a diagram of such, Rem. \ref{remdeffil}.

\noindent \hspace{-1cm} $\widetilde\O^{crys}_.:=\Z/p^.\otimes l^{-1}\O^{crys}$: Def. \ref{defocrysmod}.

\noindent \hspace{-1cm} $\smash{\widetilde{Fil}^._{.,T^\s}(-h)(D)},\, \smash{\widetilde{Fil}^.\Omega^\bullet_{.,T^\s}(-h)(D)},\, \smash{\widetilde{Fil}^.L\Omega^\bullet_{.,T^\s}(-h)(D)}$: the $L$-normalized variant of the above filtered complexes: Def. \ref{deffilmod}, Rem. \ref{remdeffilmod}, Lem. \ref{lemfil} \ref{lemfilii}.

\noindent \hspace{-1cm} $\smash{{Fil}^._{T^\s}(-h)(D)},\, \smash{{Fil}^.\Omega^\bullet_{T^\s}(-h)(D)},\, \smash{{Fil}^.L\Omega^\bullet_{T^\s}(-h)(D)}$: the limit variant of the above filtered complexes: Def. \ref{deffillim}.

\noindent \hspace{-1cm} $\smash{Lie_{.,T^\s}(-h)(D)}$, $\smash{\widetilde{Lie}_{.,T^\s}(-h)(D)}$,
$\smash{Lie_{T^\s}(-h)(D)}$: Lem. \ref{lemfil} \ref{lemfilii}, Lem. \ref{lemdeffillim} \ref{lemdeffillimi}.

\noindent \hspace{-1cm} $\tilde F^{(-/T^\s_.)}:id\rightarrow \tilde F_{T_.^\s}^{-1}$: the lifted relative Frobenius associated to a Frobenius lift on a $p$-adic log scheme $T^\s$, Lem. \ref{lemfroblift} \ref{lemfroblifti}.

\noindent \hspace{-1cm} $\tilde F$: the lifted Frobenius endomorphism of the structural ring of $\smash{T^\s_{et}}$, Lem. \ref{lemfroblift} \ref{lemfrobliftii}.

\noindent \hspace{-1cm}  $\smash{\tilde F^{(-/T_.^\s)}}:\smash{(\tilde F)^*}\to \smash{\tilde F_{T_.^\s}^*}$ or  $\smash{(F)^*}\to \smash{\tilde F_{T_.^\s}^*}$: ringed variants of the lifted relative Frobenius, Rem. \ref{remfroblift} \ref{remfrobliftiii}.

\noindent \hspace{-1cm} $Fr$: the semi-linear endomorphism of $D(-h)_{T^\s_.}$, $\smash{\Omega^\bullet_{T^\s_.}(D(-h))}$ or $\smash{(L(\Omega^\bullet_{T^\s_.}(D(-h))))_{T^\s_.}}$ built from the lifted relative Frobenius and the Frobenius of the Dieudonn\'e crystal, Def. \ref{defFr}.

\noindent \hspace{-1cm} $\varphi:\smash{\widetilde{Fil}^1_{.,T^\s}(-h)(D)} \to {\widetilde{Fil}^0_{.,T^\s}(-h)(D)}$, $\smash{\widetilde{Fil}^1\Omega^\bullet_{.,T^\s}(-h)(D)}\to \smash{\widetilde{Fil}^0\Omega^\bullet_{.,T^\s}(-h)(D)}$ or \break ${\widetilde{Fil}^1L\Omega^\bullet_{.,T^\s}(-h)(D)}\to \smash{\widetilde{Fil}^0L\Omega^\bullet_{.,T^\s}(-h)(D)}$: the unique morphism such that $p\varphi$ is induced by $Fr$, Prop. \ref{defphi}.

\noindent \hspace{-1cm} $Kom^{1,\varphi}(T^\s_{[.],.,et},\widetilde \O^{crys}_.)$: the category of complexes of $(1,\varphi)$-modules, Def. \ref{defSonephi} and explanations below.

\noindent \hspace{-1cm} $\cal S^{1,\varphi}_{et,.,T^\s}(-h)$, $\cal S\Omega^{\bullet,1,\varphi}_{et,.,T^\s}(-h)$, $\cal SL\Omega^{\bullet,1,\varphi}_{et,.,T^\s}(-h)$: three variants of the twisted syntomic complex functor attached to a given global (or semi-simplicial local) embedding with Frobenius lift and effective log divisor, taking their values in the category of complexes of $(1,\varphi)$-modules of $(\smash{T^{\s}_{.,et},\widetilde \O_.^{crys}})$, Def. \ref{defSonephi}.

\noindent \hspace{-1cm} $\cal S^{1,\varphi}_{et,.,X^\s}(-h)$: the twisted syntomic complex functor attached to $(X^\s,h)$ in $\cal B^\s_0$, taking its values in the derived category of $(1,\varphi)$-modules of $(\smash{X^{\s,\N}_{et},\widetilde \O_.^{crys}})$, Prop.  \ref{indepSet}.

\noindent \hspace{-1cm}  $\cal S_{et,.,X^\s}(-h)$ (resp. $\cal S_{et,X^\s}(-h)$): the twisted syntomic complex functor attached to $(X^\s,h)$ in $\cal B^\s_0$, taking its values in the derived category of modules of $(\smash{X^{\s,\N}_{et},\widetilde \O_.^{crys,F=1}})$ (resp. $(\smash{X^{\s}_{et},\widetilde \O^{crys,F=1}})$), Def.  \ref{defSX}.

\noindent \hspace{-1cm} $Lie^{syn}(D)$: the pullback of $Lie(D)$ to the small syntomic site, Def. \ref{lemsyn} \ref{lemsyni}.

\noindent \hspace{-1cm} $D^{syn}$, $\smash{L^{syn}\Omega^\bullet_{T_.}(D)}$: the pullbacks of $D$ and $\smash{L\Omega^\bullet_{T_.}(D)}$  to the small syntomic crystalline site, Def. \ref{lemsyn} \ref{lemsynii}.

\noindent \hspace{-1cm} $can_D^{syn}$, $can_L^{syn}$:   Def. \ref{lemsyn} \ref{lemsyniii}.

\noindent \hspace{-1cm} $Fil^1D^{syn}$, $\smash{Fil^1L^{syn}\Omega^\bullet_{T_.}(D)}$: Def. \ref{lemsyn} \ref{lemsyniii}.

\noindent \hspace{-1cm} $\smash{Fil^{i,crys}_.D}$ (resp. $\smash{Fil^{i}L^{crys}\Omega^\bullet_{.,T}(D)}$:  the resulting complexes of $(X^{\N}_{syn},\O^{crys}_.)$, Def. \ref{deffilsynmod} \ref{deffilsynmodi}.

\noindent \hspace{-1cm} $\smash{\widetilde{Fil}^{i,crys}_.D}$, $\smash{\widetilde{Fil}^{i}L^{crys}\Omega^\bullet_{.,T}(D)}$ and $\smash{\widetilde{Lie_.^{syn}}(D)}$: $L$-normalized variants, Lem. \ref{tdfilonesyn} \ref{tdfilonesynii}.

\noindent \hspace{-1cm} $\smash{Lie_.^{syn}(D)}$, $\smash{\widetilde{Lie}_.^{syn}(D)}$: Prop. \ref{tdfilonesyn}
\ref{tdfilonesyni}.

\noindent \hspace{-1cm} $Fr$: the semi-linear endomorphism of $Fil^{0,crys}_.D$
(resp. $Fil^0L^{crys}\Omega_{.,T}(D)$)
built from the relative Frobenius (resp. , the Frobenius lift)
and the Frobenius of the Dieudonn\'e crystal, Def. \ref{defFrsyn}.

\noindent \hspace{-1cm} $\varphi:\smash{\widetilde{Fil}{}^{1,crys}_.D}\to (F)_*\smash{\widetilde{Fil}{}^{0,crys}_.D}$
or $\smash{Fil^{1}L^{crys}\Omega^\bullet_{.,T}(D)}\to (F)_*\smash{Fil^{1}L^{crys}\Omega^\bullet_{.,T}(D)}$:
the unique morphism such that $p\varphi$ is induced by $Fr$, Prop. \ref{defphisyn}.

\noindent \hspace{-1cm} $\smash{\cal S^{1,\varphi}_{syn,.,X}}$: the syntomic complex functor attached to a
diagram of separated
schemes with local $p$-bases, taking its values in the category of $(1,\varphi)$-modules of
$\smash{(X^\N_{syn},\O_.^{crys})}$, Def. \ref{defSonephisyn}, \ref{defSonephisyni}.

\noindent \hspace{-1cm} $\smash{\cal SL\Omega^{\bullet,1,\varphi}_{syn,.}}$: the syntomic complex functor attached
to a diagram of global embeddings with Frobenius lift, taking its values in the category of
complexes of $(1,\varphi)$-modules of $\smash{(X^\N_{syn},\O_.^{crys})}$, \ref{defSonephisyn} \ref{defSonephisynii}.

\noindent \hspace{-1cm} $\cal S_{syn,.,X}$:
the syntomic complex functor, taking its values in the derived category
of modules of $\smash{(X^\N_{syn},\O_.^{crys,F=1})}$, Def. \ref{defSsyn}.

\noindent \hspace{-1cm} $Kom(\cal A)$, $D(\cal A)$: the category of complexes and the derived category of $\cal A$.

\noindent \hspace{-1cm} $MF$: the mapping fiber functor, Def. \ref{functoMF}.

\noindent \hspace{-1cm} $C$: an irreducible smooth curve over $\Sigma_1$, (\ref{diagcomplet}).

\noindent \hspace{-1cm} $Z$, $C^\s$, $U$, $Z_v$, $C_v$, $C_v^\s$, $U_v$:
various log schemes related to $C$, (\ref{diagcomplet}).

\noindent \hspace{-1cm} $z,j,o,z_v,j_v,o_v,{\iota_{Z_v}},{\iota_{C_v}},{\iota_{C_v^\s}},{\iota_{U_v}}$:
natural morphisms relating them, (\ref{diagcomplet}).

\noindent \hspace{-1cm} $O_{C,v}$, $O_v$, $K_v$, $k_v$: the local ring at $v$, the complete ring at $v$,
the complete field at $v$, the residue field at $v$, (\ref{diagcomplet}).

\noindent \hspace{-1cm} $Z_J$, $J$, $J^\s$: diagrams of log schemes related to $C$, (\ref{defdiag}), (\ref{defdiagZ}).

\noindent \hspace{-1cm} $z_J$, $o_J$, $m_Z$, $m$, $m^\s$: natural morphisms relating them to each other and to $Z$, $C$, $C^\s$, (\ref{diagdiag}).

\noindent \hspace{-1cm} $Sma_v$, $Sma$: smashing functors, Def. \ref{defsma} \ref{defsmai}, \ref{defsmaii}.

\noindent \hspace{-1cm} $X^+$: the diagram attached to a morphism of diagrams of extremal type $Y\to X$,
Def. \ref{defitype} \ref{defitypeiii}.

\noindent \hspace{-1cm} $i$, $i_+$, $\rho$, $\sigma$: natural morphisms between the diagrams $X$, $Y$ and $X^+$,
Def. \ref{defitype} \ref{defitypeiii}.

\noindent \hspace{-1cm} $\Gamma^Y(X,-)$: the functor taking a sheaf on $X$ or $X^+$ to its sections
on $X$ vanishing on $Y$, Def. \ref{defvan} \ref{defvani}, \ref{defvanii}.

\noindent \hspace{-1cm} $\underline\Gamma^Y(-)$ (resp. $\underline\Gamma^Y(X,-)$) the functor taking
a sheaf on $X$ (resp. on $X^+$) to the subsheaf (resp. to the subsheaf of its restriction to $X$)
formed by the sections which vanish on $Y$, Def. \ref{defvan} \ref{defvani}, \ref{defvanii}.

\noindent \hspace{-1cm} $A$: a semi-Abelian scheme over $C$, Prop. \ref{MVfl}.

\noindent \hspace{-1cm} $A_{p^.}$: the projective system of its $p$-primary torsion subgroups, Prop. \ref{MVfl}.

\noindent \hspace{-1cm} $A^0$: its connected component, Prop. \ref{MVfl}.

\noindent \hspace{-1cm} $\Phi=A/A^0$: its component group, proof of Prop. \ref{MVfl}.

\noindent \hspace{-1cm} $Sch_{lft}/R$, $Sch_{lft}/K$, $For/R$ and $Rig/K$ for $R$ a complete discrete
valuation ring with fraction field $K=R[{1\over t}]$: the categories of schemes which are locally of finite type over $R$, schemes which are locally of finite type over $K$, $t$-adic formal schemes over $R$ and rigid analytic spaces over $K$, (\ref{diagrigfor}).

\noindent \hspace{-1cm} $(-)^{for}$, $(-)^{an}$, $(-)_{|K}$: functors relating these categories, (\ref{diagrigfor}).

\noindent \hspace{-1cm} $Fin/R$, $Fin/K$: the category of finite schemes over $R$, $K$, Cor. \ref{rig} \ref{rigi}, \ref{rigii}.

\noindent \hspace{-1cm} $qff$, $ff$, $fet$: other instances of $top$, Def. \ref{defqff}.

\noindent \hspace{-1cm} $G_v/C_v$: the Raynaud group of $A_{|U_v}/U_v$, Def. \ref{defH}.

\noindent \hspace{-1cm} $H/J^+$: the diagram of $p$-divisible groups associated to $A/C$, Def. \ref{defH}.

\noindent \hspace{-1cm} $\G_{cart}(-)$: the category of cartesian sections of a cofibered category, Prop.
\ref{MVsyn} \ref{MVsyni}.

\noindent \hspace{-1cm} $\cal M_{log}(X)$: the category of log $1$-motives over $X$, the spectrum of a complete discrete valuation ring $R$, Def. \ref{defMlog} \ref{defMlogi}.

\noindent \hspace{-1cm} $\B_X$: the category of finite \'etale $X$-schemes,   Def. \ref{defMlog} \ref{defMlogii}.

\noindent \hspace{-1cm} $\cal M(X)$, $\cal M\cal T_{log}(X)$, $\cal M\cal T(X)$: subcategories of $\cal M_{log}(X)$, comments below Def. \ref{defMlog}.

\noindent \hspace{-1cm} $\cal C\square \cal C$: a certain category of diagrams of an exact category $\cal C$, Def. \ref{defbaer} \ref{defbaeri}.

\noindent \hspace{-1cm} $\boxplus:\C\square \cal C\rightarrow \cal C$: the Baer sum functor, Def. \ref{defbaer} \ref{defbaerii}.

\noindent \hspace{-1cm} $\cal M^t_{log}(X)$: another subcategory of $\cal M(X)$, Def. \ref{tdec}, \ref{tdecii}.

\noindent \hspace{-1cm} $\Z(1)$: the Tate $1$-motive, comment before Lem. \ref{lemmot2}.

\noindent \hspace{-1cm} $Kum(x)$: the Kummer log $1$-motive associated to some $x\in K^\times$, comment before Lem. \ref{lemmot2}.

\noindent \hspace{-1cm} $Mod^{fv}(E_{/X},A_{|X})$: the category of $(f,v)$-modules, Def. \ref{deffvmod} \ref{deffvmodi}.

\noindent \hspace{-1cm} $\cal Hom_{A_{|X}}^{fv}$: the bifunctor of inner homomorphisms of $A_{|X}$-modules enriched with an $(f,v)$-module structure, Def. \ref{deffvmod} \ref{deffvmodii}.

\noindent \hspace{-1cm} $E_{/X}^{fv}$: the topos of $(f,v)$-objects, comment before Lem. \ref{lemfv1}.

\noindent \hspace{-1cm} $\chi:(E_{/X},A_{|X})\to (E_{/X}^{fv},A_{|X})$: the natural morphism of ringed topoi, Lem. \ref{lemfv1} \ref{lemfv1iii}.

\noindent \hspace{-1cm} $G^{fv}$: the $\Z$-$(f,v)$-module of the crystalline topos defined by locally free group or $p$-divisible group $G$ together with its Verschiebung and relative Frobenius, Sect. \ref{paranotfv2}.

\noindent \hspace{-1cm} $B(-)$, $B_{\cal O}(-)$: bidualizing functors for $(f,v)$-modules, Sect. \ref{paradefbid}.

\noindent \hspace{-1cm} $D(G)$: \cite{BBM}'s covariant Dieudonn\'e crystal of some $p$-divisible group $G$, Prop. \ref{rappelBBMi}.

\noindent \hspace{-1cm} $D_X:\cal M(X)\to \cal DC_{CRYS,fl}(X)$: the Dieudonn\'e functor for $1$-motives over a $\Sigma_1$-scheme $X$, Def. \ref{defD1mot}.

\noindent \hspace{-1cm} $|D_{X^\s}|:Ext^1_{\M_{log}(X)}(\Z,\Z(1))\to Ext^1_{\cal DC_{CRYS^\s,fl}}(X^\s)(o^*D_X(\Z), o^*D_X(\Z(1)))$ for $X^\s=(Spec(R),Spec(k))$: a canonical map which is compatible with $D_X$ and $D_K$, Prop. \ref{prolD1}.

\noindent \hspace{-1cm} $D_{X^\s}:\M_{log}(X)\to \cal DC(X^\s)$: a canonical functor which is compatible with $D_X$ and $D_K$ and induces the above map, Cor. \ref{defDMlog}.

\noindent \hspace{-1cm} $M_{log}(A)$: the log $1$-motive associated to a semi-Abelian scheme over $Spec(R)$, Def. \ref{defSASloc} \ref{defSASloci}.

\noindent \hspace{-1cm} $D_{X^\s}(A)=D_{X^\s}(M_{log}(A))$: the Dieudonn\'e crystal associated to a semi-Abelian scheme over $Spec(R)$, Def. \ref{defSASloc} \ref{defSASlocii}.

\noindent \hspace{-1cm} $D_{C^\s}:SAS(C,Z)\to \cal DC(C,Z)$: the Dieudonn\'e functor for semi-Abelian schemes over $C$ which are Abelian outside of $Z$, Def. \ref{defDCglob}.

\noindent \hspace{-1cm} $W_k(U)$, $I_k(U)$: the ring of Witt vectors of some scheme $U$ of characteristic $p$, a certain ideal inside it, Lem.  \ref{lemcomp3}.

\noindent \hspace{-1cm} $W_k^{dp}(U)$, $I_k^{dp}(U)$: the divided power envelope of $W_k(U)$ with respect to $I_k(U)$, the tautological $dp$-ideal, Lem. \ref{lemcomp3} \ref{lemcomp3i}.

\noindent \hspace{-1cm} $\Theta_k:W_k^{dp}\to \O_k^{crys}$: the canonical morphism, (\ref{defThetak}).

\noindent \hspace{-1cm} $\cal I_.^{crys}$, $\widetilde{\cal I}_.^{crys}$: the tautological $dp$-ideal of $\cal O_.^{crys}$, its normalization on the small syntomic site of certain schemes, Lem. \ref{lemcomp4} and comment above.

\noindent \hspace{-1cm} $\varphi$: the Frobenius divided by $p$ on $\widetilde{\cal I}_.^{crys}$, Lem. \ref{lemcomp2} \ref{lemcomp2ii}.

\section{Index of terminology}


\noindent \hspace{-1cm} \emph{functor descending along a change of the base category}, comment after Lem. \ref{CBdescent1}.

\noindent \hspace{-1cm} \emph{pretopology, premorphism of pretopologies}: Def. \ref{defweakmor} \ref{defweakmori}, \ref{defweakmorii}.

\noindent \hspace{-1cm} \emph{weak morphism, morphism of topoi}: Def.  \ref{defweakmor} \ref{defweakmoriv}.

\noindent \hspace{-1cm} \emph{prevariable, pretopology}: Def. \ref{deftopfib} \ref{deftopfibi}.

\noindent \hspace{-1cm} \emph{weakly variable, variable topos}: Def. \ref{deftopfib} \ref{deftopfibii}.

\noindent \hspace{-1cm} \emph{projective limit formula for the direct image of a morphism of diagrams}: Lem. \ref{lemfibtop} \ref{lemfibtopii}.

\noindent \hspace{-1cm} \emph{componentwise injective, componentwise flasque, componentwise $\cal P$-acyclic}: Lem. + Def. \ref{acycf} \ref{acycfiii}.

\noindent \hspace{-1cm} \emph{d-injective, d-flasque, d-$\cal P$-acyclic}: Lem. + Def. \ref{acycf} \ref{acycfiv}.

\noindent \hspace{-1cm} \emph{premorphism of prevariable pretopologies}: Def. \ref{defmorT} \ref{defmorTi}.

\noindent \hspace{-1cm} \emph{weak morphism, morphism of weakly variable topoi}: Def. \ref{defmorT} \ref{defmorTii}.

\noindent \hspace{-1cm} \emph{crystal condition}: Lem. \ref{lemcrystal1}.

\noindent \hspace{-1cm} \emph{$e$-exact, reflects $e$-exactness, fully $e$-exact, fully Abelian}: Sect. \ref{deffullyexact}.

\noindent \hspace{-1cm} \emph{finite locally free group, $p$-divisible group}: Def. \ref{defflfgppdiv} \ref{defflfgppdivi}, \ref{defflfgppdivii}.

\noindent \hspace{-1cm} \emph{Cartier dual}: (\ref{defcartierduality1}), (\ref{defcartierduality2}).

\noindent \hspace{-1cm} \emph{Abelian scheme, torus, semi-Abelian scheme, twisted constant group, $1$-motive}: Def. \ref{defmotabtor} \ref{defmotabtori}, \ref{defmotabtoriii}, \ref{defmotabtorii}, \ref{defmotabtoriv}, \ref{defmotabtorv}.

\noindent \hspace{-1cm} \emph{normalized module, $L$-normalized complex}, Def. \ref{defnorm} \ref{defnormi}, \ref{defnormii}.

\noindent \hspace{-1cm} \emph{$p$-adic scheme, $p$-adic log scheme}, Def. \ref{defpadic} \ref{defpadici}, \ref{defpadicii}.

\noindent \hspace{-1cm} \emph{quasi-coherent (or locally free or locally free of finite type) module on a scheme}, Def. \ref{defqcohsch}, Lem. \ref{lemqcohsch}  \ref{lemqcohschi},  \ref{lemqcohschii}.

\noindent \hspace{-1cm} \emph{quasi-coherent projective system of modules}, Def. \ref{defqcohpadic} \ref{defqcohpadici}.

\noindent \hspace{-1cm} \emph{quasi-coherent module on a $p$-adic scheme}, Def. \ref{defqcohpadic} \ref{defqcohpadicii}.

\noindent \hspace{-1cm} \emph{quasi-coherent complex of projective systems of modules}, Def. \ref{Lqcohpadic} \ref{defLqcohpadici}.

\noindent \hspace{-1cm} \emph{quasi-coherent complex on a $p$-adic scheme}, Def. \ref{Lqcohpadic} \ref{defLqcohpadicii}.

\noindent \hspace{-1cm} \emph{module with quasi-coherent (or locally free or locally free of finite type) realizations on the crystalline site}, Rem. \ref{lemlcrys} \ref{lemlcrysii}.

\noindent \hspace{-1cm} \emph{quasi-coherent (or locally free or locally free of finite type) crystal}, Rem. \ref{lemlcrys} \ref{lemlcrysiv}.

\noindent \hspace{-1cm} \emph{relatively perfect morphism}, Def. \ref{defpb} \ref{defpbi}, \ref{defpbii}.

\noindent \hspace{-1cm} \emph{finite $p$-basis, morphism with local finite $p$-bases}, Def. \ref{defpb} \ref{defpbii}, \ref{defpbiii}, \ref{defpbiv}.

\noindent \hspace{-1cm} \emph{local embedding, global embedding}, Def. \ref{defemb1} \ref{defemb1i}.

\noindent \hspace{-1cm} \emph{morphism of embeddings extending or lifting a morphism of schemes}, Def. \ref{defemb1} \ref{defemb1iii}.

\noindent \hspace{-1cm} \emph{global embedding lifting a scheme}, Lem. \ref{lememb1} \ref{lememb1iii}.

\noindent \hspace{-1cm} \emph{algebra of $dp$-differential operators}, Def. \ref{defDiff}.

\noindent \hspace{-1cm} \emph{module with connection}, Lem. \ref{defnabla}.

\noindent \hspace{-1cm} \emph{integrable connection}, Lem. \ref{defnablaint}.

\noindent \hspace{-1cm} \emph{quasi-nilpotent connection}, Lem. \ref{defqnil}.

\noindent \hspace{-1cm} \emph{hyper $dp$-stratification}, Prop. \ref{hyperstrat}.

\noindent \hspace{-1cm} \emph{hyper $dp$-differential operators}, Def. \ref{tensnablaHdp}.

\noindent \hspace{-1cm} \emph{de Rham complex}, Lem. \ref{dR} \ref{dRi}, Def. \ref{dRM}.

\noindent \hspace{-1cm} \emph{linearization functor}, Lem. \ref{defLT}.

\noindent \hspace{-1cm} \emph{local crystal}, Lem. \ref{lemft}.

\noindent \hspace{-1cm} \emph{morphism admitting flat liftings to local embeddings}, Def. \ref{defflatlift}.

\noindent \hspace{-1cm} \emph{integral log structure on a ringed topos}, Sect. \ref{twlogtop}.

\noindent \hspace{-1cm} \emph{effective log divisor}, Def. \ref{deftwlogtop} \ref{deftwlogtopi}.

\noindent \hspace{-1cm} \emph{twisted module},  Def. \ref{deftwlogtop} \ref{deftwlogtopiii}.

\noindent \hspace{-1cm} \emph{$top$ effective Cartier divisor}, Def. \ref{defcartierdiv}.

\noindent \hspace{-1cm} \emph{flat fine chart, local flat fine charts}, Def. \ref{flatcharts} \ref{flatchartsi}, \ref{flatchartsii}.

\noindent \hspace{-1cm} \emph{relative Frobenius}, Def. \ref{defFrel} \ref{defFreli}, \ref{defFrelii}.

\noindent \hspace{-1cm} \emph{inverse Cartier morphism}, Def, \ref{defCinv} \ref{defCinvi}, \ref{defCinvii}.

\noindent \hspace{-1cm} \emph{crystal with trivial $p$-curvature}, Lem. \ref{defKp}.

\noindent \hspace{-1cm} \emph{$\underline t$-torsion free module or crystal}, Def. \ref{deftfr} \ref{deftfri}, \ref{deftfrii}.

\noindent \hspace{-1cm} \emph{Dieudonn\'e crystal, truncated Dieudonn\'e crystal of level $1$}, Def. \ref{defDC} \ref{defDCi}, \ref{defDCii}.

\noindent \hspace{-1cm} \emph{$mod$ $p$ Hodge filtration}, Def. \ref{deflie} \ref{defliei}, \ref{deflieii}.

\noindent \hspace{-1cm} \emph{lifted relative Frobenius}, Lem. \ref{lemfroblift} \ref{lemfroblifti}, Rem. \ref{remfroblift} \ref{remfrobliftiii}.

\noindent \hspace{-1cm} \emph{twisted syntomic complex functor on the \'etale site}, Def.  \ref{defSX}.

\noindent \hspace{-1cm} \emph{global situation, local situation}, Sect. \ref{scotss}.

\noindent \hspace{-1cm} \emph{syntomic complex functor on the syntomic site}, Def. \ref{defSsyn}.

\noindent \hspace{-1cm} \emph{category of true arrows}, Sect. \ref{prelimMF1}.

\noindent \hspace{-1cm} \emph{functorial mapping fiber}, Def. \ref{functoMF}.

\noindent \hspace{-1cm} \emph{smashing functor}, Def. \ref{defsma} \ref{defsmai}, \ref{defsmaii}.

\noindent \hspace{-1cm} \emph{morphism of extremal type}, Def. \ref{defitype} \ref{defitypeiii}.

\noindent \hspace{-1cm} \emph{functor of vanishing sections},  Def. \ref{defvan} \ref{defvani}, \ref{defvanii}.

\noindent \hspace{-1cm} \emph{Raynaud group of a semi-stable Abelian variety}, Sect. \ref{Raygp}.

\noindent \hspace{-1cm} \emph{complete Mayer-Vietoris triangle for the twisted syntomic complex on the \'etale site}, Cor. \ref{corMVsyn}, \ref{corMVsynii}.

\noindent \hspace{-1cm} \emph{localization triangle for the twisted syntomic complex on the \'etale site},
Cor. \ref{cortdspe}, \ref{remtdspeii}.

\noindent \hspace{-1cm} \emph{log $1$-motive over a complete discrete valuation ring}, Def. \ref{defMlog} \ref{defMlogi}.

\noindent \hspace{-1cm} \emph{the weight filtration of a log $1$-motive}, comments below Def. \ref{defMlog}.

\noindent \hspace{-1cm} \emph{Baer sum}, Def \ref{defbaer} \ref{defbaerii} and comment above.

\noindent \hspace{-1cm} \emph{$t$-decomposition of a log $1$-motive}, Def. \ref{tdec} \ref{tdeci}.

\noindent \hspace{-1cm} \emph{Kummer log $1$-motive}, comment before Lem. \ref{lemmot2}.

\noindent \hspace{-1cm} \emph{log $1$-motive of a semi-Abelian scheme}. Def. \ref{defSASloc} \ref{defSASloci}

\end{document}